\documentclass{article}

\RequirePackage{amsthm,amsmath,amsfonts,amssymb}

\RequirePackage{hyperref}%
\hypersetup{
    colorlinks,
    linkcolor={red!50!black},
    citecolor={blue!50!black},
    urlcolor={blue!80!black}
}
\RequirePackage{graphicx}%
\usepackage{setspace}
\usepackage{booktabs}
\usepackage{mathtools}

\usepackage[capitalise,nameinlink]{cleveref}

\crefname{supp}{Supplement}{Supplements}
\crefname{appendix}{Supplement}{Supplements}
\Crefname{appendix}{Supplement}{Supplements}
\crefname{subappendix}{Supplement}{Supplements}
\Crefname{subappendix}{Supplement}{Supplements}
\crefname{subsubappendix}{Supplement}{Supplements}
\Crefname{subsubappendix}{Supplement}{Supplements}
\crefname{enumi}{part}{parts}
\Crefname{enumi}{Part}{Parts}
\crefname{enumii}{part}{parts}
\Crefname{enumii}{Part}{Parts}

\usepackage[margin=1in]{geometry}

\usepackage[shortlabels]{enumitem}

\usepackage{mathrsfs}

\usepackage{algorithm}
\usepackage{algpseudocode}
\usepackage{tikz}
\usepackage{subcaption}
\usetikzlibrary{shapes,arrows,fadings}
\usepackage{circuitikz}
\usepackage{apptools}
\usepackage{lmodern}
\usepackage{placeins}
\usepackage{authblk}

\usepackage[
natbib=true,
backend=biber,
sorting=nyt,
useprefix=true,
style=authoryear,
url=false,
doi=true,
eprint=false,
giveninits=true,
uniquename=false]{biblatex}
\AtEveryBibitem{%
  \clearfield{note}%
}
\AtEveryBibitem{\clearlist{language}}
\AtEveryBibitem{\clearfield{month}}
\AtEveryCitekey{\clearfield{month}}

\usepackage[toc, title, page,header]{appendix}

\usepackage{etoc}
\makeatletter
\patchcmd{\l@section}{1.5em}{2.5em}{}{}
\makeatother

\makeatletter
\newcommand{\supplementnumbering}[1]{%
  \setcounter{#1}{0}%
  \@namedef{the#1}{S.\arabic{#1}}%
  \@namedef{theH#1}{S.\arabic{#1}}%
}
\makeatother

\theoremstyle{plain}

\newtheorem{theorem}{Theorem}

\newtheorem{lemma}{Lemma}
\newtheorem{proposition}{Proposition}
\newtheorem{remark}{Remark}

\theoremstyle{remark}
\newtheorem{assumption}{Assumption}
\newtheorem{definition}{Definition}
\newtheorem{example}{Example}

\crefname{assumption}{Assumption}{Assumptions}
\Crefname{assumption}{Assumption}{Assumptions}
\crefname{counterexample}{Counterexample}{Counterexamples}
\Crefname{counterexample}{Counterexample}{Counterexamples}
\crefname{examplerepeat}{Example}{Examples}
\Crefname{examplerepeat}{Example}{Examples}
\crefname{examplerepeattwo}{Example}{Examples}
\Crefname{examplerepeattwo}{Example}{Examples}
\crefname{algorithm}{Algorithm}{Algorithms}
\Crefname{algorithm}{Algorithm}{Algorithms}

\newcommand{\mc}[1]{\mathcal{#1}}

\DeclareMathOperator*{\esssup}{ess\,sup}
\DeclareMathOperator*{\essinf}{ess\,inf}

\newcounter{examplecontid}

\newenvironment{examplecont}[2][continued]
  {\stepcounter{examplecontid}%
   \renewcommand{\theexample}{\ref*{#2}}%
   \addtocounter{example}{-1}%
   \example[#1]}
  {\endexample}

\setlist[enumerate,1]{label=(\alph*)}

\title{On Local Inference for Functionals of Ill-Posed Inverse Problems}
\date{October 2026}
\author[1]{Ellen Graham\thanks{Corresponding author: \texttt{Ellen.Graham@PennMedicine.upenn.edu}}}
\author[2]{Andrea Rotnitzky}
\affil[1]{Department of Biostatistics, Epidemiology, and Informatics, University of Pennsylvania}
\affil[2]{Department of Biostatistics, University of Washington}

\begin{document}

\maketitle

\begin{abstract}
    We study inference for continuous linear functionals of solutions to
ill-posed inverse problems, including Euclidean parameters arising in
nonparametric instrumental variables, proximal causal inference, and
shadow-variable models. Although pointwise root-$n$ inference is
possible under standard conditions, we show that this need not translate
into locally reliable inference. Under a model imposing only minimal
regularity conditions, the target is Hellinger discontinuous when both
inverse problems are ill-posed: arbitrarily small changes in the
data-generating law can produce arbitrarily large changes in the target,
ruling out informative locally uniform inference. Imposing a uniform
bound on a minimum-norm solution restores Hellinger continuity under regularity conditions, but does
not by itself ensure stable root-$n$ inference. If, in addition, both inverse problems are ill-posed and the solution whose
norm is uniformly bounded is nonunique at a law $P_0$, then under regularity conditions the target fails
to be pathwise differentiable; consequently, asymptotically linear
estimators exhibit nonvanishing root-$n$ bias under some local
alternatives, causing Wald confidence intervals to undercover
asymptotically. Under regularity conditions, pathwise differentiability is restored if the uniformly
bounded solution is unique at $P_0$. It is also restored, regardless of
uniqueness, under the stronger requirement that the corresponding
solution satisfy a uniformly bounded $\beta$-source condition for some
$\beta>0$. We relate these regimes to very weak, weak, and strong
identification in instrumental variables models.
\end{abstract}

\begin{refsection}
\etocdepthtag.toc{mtmain}

 \section{Introduction}
\label{sec:intro}

Many parameters of interest can be expressed as real-valued continuous linear functionals of solutions to a pair of integral operator equations, including estimands in nonparametric instrumental variables (NPIV) regression \parencite{newey_instrumental_2003-2, darolles_nonparametric_2011-1, blundell_semi-nonparametric_2007, horowitz_applied_2011, severiniEfficiencyBoundsEstimating2012, bennettSourceConditionDouble2023b, bennett_inference_2023, olivas-martinez_source-condition_2025}, in proximal causal inference \parencite{cui_semiparametric_2024}, and in inference under nonignorable missingness using shadow variables \parencite{miao_identification_2018, liNonparametricInferenceMean2022}. Such equations are often ill-posed in the sense of Hadamard \parencite{hadamard_lectures_1952}: they may fail to admit a unique solution \parencite{newey_instrumental_2003-2, santos_inference_2012, kallusCausalInferenceUnmeasured2022}, and even when a unique solution exists, it need not depend continuously on the inputs to the equation \parencite{carrasco_chapter_2007}. When the variables of interest are continuously distributed, this behavior is typical rather than exceptional \parencite{carrasco_chapter_2007}. Nevertheless, pointwise root-$n$ asymptotically linear estimators of the target functional remain attainable through Neyman-orthogonal estimating equations whose bias is insensitive to first-order errors in the ill-posed solutions \parencite{chernozhukovDoubleDebiasedMachine2018b, bennettSourceConditionDouble2023b, bennett_inference_2023, olivas-martinez_source-condition_2025}.

A pointwise asymptotic guarantee at a fixed law $P_0$ describes the
behavior of an inferential procedure when the data are generated from
exactly the same law $P_0$ for every sample size as $n\to\infty$.
However, at any given sample size $n$, there may be laws $P_n$ whose
Hellinger distance from $P_0$ is sufficiently small that data generated
under $P_n$ cannot be reliably distinguished from data generated under
$P_0$ on the basis of $n$ observations. This makes it important to ask
whether the target parameter, denoted by $\varphi(P)$, is itself stable
as $P$ varies near $P_0$. If arbitrarily small changes in the
data-generating law, as measured by Hellinger distance, can produce
nonvanishing, or even arbitrarily large, changes in $\varphi(P)$, then
no inferential procedure can perform uniformly well over sufficiently
small neighborhoods of $P_0$. As established in a large body of work
\parencite{bahadur_nonexistence_1956,donoho_one-sided_1988,
donoho_geometrizing_1991,low_nonparametric_1997,cai_adaptation_2004,
pfanzagl_nonexistence_1998,potscher_lower_2002,dufour_impossibility_1997},
failure of Hellinger continuity has direct inferential consequences:
point estimators must incur substantial error at some nearby laws, while
confidence intervals that maintain coverage throughout such
neighborhoods must remain correspondingly wide; in the most severe cases
considered below, the expected diameter of any locally honest confidence
interval is necessarily infinite.

Even when $P\mapsto\varphi(P)$ is Hellinger continuous at $P_0$, continuity
alone does not determine whether the parameter varies at a scale
compatible with ordinary root-$n$ inference. The standard way to study
this question is to examine the behaviour of the target functional along each possible quadratic mean differentiable (QMD) submodel
$\{P_t:t\in[0,\delta)\}$ through $P_0=P_{t=0}$ \parencite[Chapter 25]{vaart_asymptotic_1998-2}. Along any
such submodel with non-zero score, the Hellinger distance between $P_t$ and
$P_0$ is of order $t$. Consequently, alternatives of the form
$P_{h/\sqrt n}$, for fixed $h > 0$, lie at the critical statistical
scale: with $n$ observations they remain distinguishable from $P_0$,
but not consistently distinguishable, in the sense that
level-$\alpha$ Neyman--Pearson test of $P_0$ against
$P_{h/\sqrt n}$ has limiting power strictly between $\alpha$ and one.
To study root-$n$ inference under such alternatives, it is typical to assess the pathwise differentiability of the parameter $\varphi$. Pathwise differentiability requires that along each QMD submodel, $t \mapsto \varphi(P_t)$ is differentiable at $t=0$, and that these derivatives depend on the submodel only through a bounded linear functional of the submodel's score, called the local parameter,  this bounded linear functional being common to all QMD submodels. For pathwise differentiable parameters, under some conditions \parencite{klaassen_consistent_1987-2}, it is possible to construct asymptotically linear estimators whose normal limiting distributions are invariant to the local submodel, facilitating asymptotic inference based on this common limiting distribution that is asymptotically valid under sampling from any such local alternative. Conversely, when pathwise differentiability fails, for any asymptotically linear estimator, there are local alternatives under which the estimator has non-vanishing root-$n$ bias, and the corresponding Wald intervals centered at asymptotically linear estimators undercover asymptotically under those local alternatives.

For the target linear functionals that are the focus of the present paper, recent
work has shown that Wald-type confidence intervals centered at estimators
obtained from Neyman-orthogonal estimating equations are not globally
uniformly valid over a minimal model that suffices for the existence of
such estimating equations \parencite{smuclerAsymptoticValidityConfidence2025}
when the dependence between the variables defining the operator equation
is allowed to become arbitrarily weak and the range of the target
parameter is unbounded. However, such global guarantees are strictly stronger than the local
guarantees considered here. Global validity requires an inferential
procedure to perform well uniformly over the entire model, including at
laws that may be sufficiently far from $P_0$ that, when $P_0$ is the
true law, they can be distinguished from $P_0$ with probability tending
to one as the sample size increases. Local validity, by contrast, only
requires good behavior over laws that remain close to $P_0$ and therefore
cannot be so easily ruled out from the observed data. There are many
settings in which locally valid inference is possible even though
globally valid inference is necessarily degenerate. Examples include
inference for the Wald estimand under arbitrarily weak instruments, for
a ratio of means when the denominator is not uniformly bounded away
from zero, and for the average treatment effect under arbitrarily poor
overlap
\parencite{gleser_nonexistence_1987,smuclerNotePropertiesConfidence2026,
fiellerProblemsIntervalEstimation1954,khanIrregularIdentificationSupport2010}.
By contrast, the question of locally valid inference for continuous
linear functionals of ill-posed inverse problems has received
comparatively little attention.

In the present work, we investigate Hellinger continuity and pathwise differentiability of parameters defined as continuous linear functionals of solutions to ill-posed inverse problems under a baseline model that suffices for the existence of a Neyman-orthogonal estimating equation and three of its submodels that make stronger assumptions about the minimum norm solutions to the inverse problems. Below, we give an overview of our results, and briefly discuss their implications for inference.

We begin by showing that, under the baseline model and some regularity conditions on $P_0$, the target parameter is not Hellinger continuous over this model at points $P_0$ where both inverse problems are ill-posed, which, as noted above, is typical when $X$ and $Z$ are continuous. The value of the parameter can diverge to infinity along a sequence of laws that converge in Hellinger distance to $P_0$. Failure of continuity further implies failure of pathwise differentiability. Standard two-point arguments then imply that no locally uniformly consistent estimators exist, and that the expected diameter of all locally uniformly honest confidence intervals must be infinite at every sample size $n$ \parencite{potscher_lower_2002}. In summary, even though pointwise inference is possible within this baseline model, local inference is impossible.

We additionally show that the maximal tangent space of this statistical model is equal to that of a locally non-parametric model under regularity conditions when the integral operators map between spaces of the same dimension. This result stands in contrast to related results that implicitly or explicitly compute tangent spaces only relative to the set of QMD at $P_0$ submodels along which the minimum norm solution is differentiable \parencite{severiniEfficiencyBoundsEstimating2012, chen_overidentification_2018-1}.

We then examine the target parameter in a smaller statistical model that additionally enforces that the minimum norm solution to one of the two inverse problems is bounded with a bound that holds uniformly over the model. This assumption is justified in several applications, such as the well-specified causal-NPIV model with bounded counterfactual outcomes and the shadow variable model under an overlap condition. We show that under a uniform norm bound, and under regularity conditions, Hellinger continuity of the target parameter is restored.

We then investigate pathwise differentiability within a model where we assume the minimum norm solutions to one or both inverse problems are uniformly bounded. We show that even within this more restrictive model, the target parameter fails to be pathwise differentiable relative to the set of all one-dimensional QMD at $P_0$ submodels, under regularity conditions, when there exist multiple solutions at $P_0$ to whichever inverse problems are known to admit bounded solutions uniformly and both inverse problems are ill-posed. When only one minimum norm solution is bounded and the other inverse problem admits a unique solution at $P_0$, our result requires $(X, Z)$ to be continuous; for $(X, Z)$ finite discrete, the parameter is then pathwise differentiable. Consequently, at laws where pathwise differentiability fails, as discussed above, all asymptotically linear estimators will exhibit non-vanishing root-$n$ scaled bias under repeated sampling from some local submodel, and Wald confidence intervals centered at them will undercover asymptotically.

In proving failure of pathwise differentiability in this model, we construct a submodel for which one singular value of the operators defining the inverse problems tends to 0 at a root-$n$ rate. We discuss the relationship between this submodel and the analysis of weak instruments in the linear instrumental variables problem in \cref{sec:weak-instruments}. We further argue that this submodel gives intuition for the properties of estimators in the NPIV, proximal causal inference, and shadow variables models when instruments, proxies, and shadow variables, respectively, are only weakly informative relative to the sample size (see \cref{subsec:weak-instrument-discrete}).

In contrast, if the solution that admits a uniform norm bound is unique at $P_0$, the parameter is pathwise differentiable relative to the set of all one-dimensional QMD at $P_0$ submodels under regularity conditions. When no restrictions are placed on the solution to the other inverse problem, pathwise differentiability holds regardless of whether the other inverse problem admits a unique solution at $P_0$. We also prove pathwise differentiability under the stronger model that assumes a uniform bound on the norm of the $\beta$-source representer of the primal (or symmetrically dual) solution for $\beta_h > 0$ (respectively, $\beta_q > 0$), a technical condition that is related to a strong-instrument assumption (see \cref{sec:weak-instruments}). Here, differentiability holds whether or not the primal (or symmetrically, dual) inverse problem admits a unique solution at $P_0$. We further show that, in both cases, the maximal tangent space is no larger than the ones derived by \cite{severiniEfficiencyBoundsEstimating2012} and \cite{chen_overidentification_2018-1} for related models. When there exist multiple solutions to the other inverse problem, the maximal tangent space is strictly smaller than that of a locally non-parametric model under mild regularity conditions, implying that some regular asymptotically linear estimators are locally inefficient at $P_0$.

To give intuition for the above results, we conduct a simple simulation study demonstrating the poor finite-sample performance of the Tikhonov-regularized DML estimator in a fully discrete NPIV problem when data are generated from local perturbations of a law $P_0$ at which both inverse problems admit multiple solutions. We compare the performance under perturbations where the minimum norm solutions diverge to infinity, where they are uniformly bounded, and under repeated sampling from the fixed law $P_0$.

The remainder of the paper is organized as follows. In \cref{sec:target-of-inference}, we fix notation, define the target functional together with the primal and dual inverse problems, review the relevant semiparametric theory, and review the implications of Hellinger discontinuity and of failure of pathwise differentiability for inference. In \cref{sec:stat-models}, we define the baseline statistical model, show that the target parameter is Hellinger discontinuous over this model when both inverse problems are ill-posed, apply the results of \cref{sec:target-of-inference} to conclude that locally uniform inference is impossible, and compute the maximal tangent space. In \cref{sec:norm-bound-models}, we consider submodels that uniformly bound the norms of the minimum norm solutions. We show that such bounds restore Hellinger continuity, give conditions under which pathwise differentiability fails and under which it holds, including under a uniform bound on the $\beta$-source representer, and bound the corresponding tangent spaces. In \cref{sec:weak-instruments}, we interpret these results by analogy with very weak, weak, and strong instrument settings in the linear instrumental variables model. In \cref{sec:relation-to-neyman-orthogonality}, we reconcile our results with existing pointwise results based on Neyman-orthogonal estimating equations, and distinguish Neyman orthogonality from pathwise differentiability. In \cref{sec:simulation}, we illustrate failure of continuity and of pathwise differentiability in a finite discrete NPIV example with a small simulation study. We conclude in \cref{sec:discussion}. All proofs are deferred to the supplement.

\section{Preliminaries}
\label{sec:target-of-inference}

\subsection{Notation}
\label{subsec:notation}
We first set notation used throughout. Let $O = (W, X, Z, D)$ denote a $\mathbb{R}^\ell$ valued random variable defined over the Borel $\sigma$-algebra $\mathcal{B}(\mathbb{R}^\ell)$, with $W, X, Z, D$ taking values in $\mathbb{R}^{\ell_W}, \mathbb{R}^{\ell_X}, \mathbb{R}^{\ell_Z}, \mathbb{N}^{\ell_D}$, respectively, and $\ell \coloneqq \ell_W + \ell_X + \ell_Z + \ell_D$. Throughout, we assume that $D$ is finite-discrete. Let $\mc{P}$ denote our statistical model, i.e. a collection of probability measures defined on the measurable space $(\mathbb{R}^\ell, \mathcal{B}(\mathbb{R}^\ell))$. We suppose that the true data generating process $P_0$ is restricted only to belong to $\mc{P}$, and that we observe $n$ i.i.d. draws of $O$ distributed according to $P_0 \in \mc{P}$. { Throughout we will assume all distributions $P$ are supported on a bounded rectangle in $\mathbb{R}^\ell$ and absolutely continuous with respect to a product measure $\lambda = \lambda_W\otimes\lambda_X\otimes\lambda_Z\otimes\lambda_D$. We let $p \coloneqq \frac{dP}{d\lambda}$ be a version of the density of $O$ with respect to $\lambda$, $p_X$ a version of the density of $X$ with respect to $\lambda_X$, $p_{X|Z}$ a version of the conditional density of $X$ given $Z$ with respect to $\lambda_X$ (e.g. $p_{X,Z}/p_Z$), and $p_{X|Z, d}$ a version of the conditional density of $X$ given $Z, D=d$ with respect to $\lambda_X$. We similarly define these densities for all other combinations of $W, X, Z, D$.}

For $P \in \mc{P}$, and $q \in [1, \infty)$, we let $L^q(P)$ denote the Banach space of (equivalence classes of $P$-a.e. equal) measurable functions of $O$ with finite $q$-th moment under $P$ equipped with the $L^q(P)$ norm  $\|f\|_{L^q(P)} \coloneqq E_{P}[|f(O)|^{q}]^{1/q}$. We let $L^q_{0}(P) \coloneqq \{f \in L^q(P): E_{P}[f(O)] = 0\}$. For $Y= f_Y(O)$ a measurable function of $O$, we let $L^q(P_Y)$ denote the Banach space of measurable functions of $Y$ equipped again with the $L^q(P_Y)$ norm $\|f\|_{L^q(P_Y)} \coloneqq E_{P}[|f(Y)|^{q}]^{1/q}$, and again let $L^q_0(P_Y) \coloneqq \{f \in L^q(P_Y): E_{P}[f(Y)] = 0\}$. For $d$ a point in the support of $D$ under $P$, we let $L^q(P_{Y|d})$ denote the collection of functions of $Y$ that have finite $q$-th moment conditional on $D=d$, similarly let $\|f\|_{L^q(P_{Y|d})} \coloneqq E_{P}[|f(Y)|^{q}|D=d]^{1/q}$. For $q = \infty$, we let $L^\infty(P)$ denote the Banach space of $P$-essentially bounded functions of $O$, and similarly define $L^\infty(P_Y)$ and $L^\infty(P_{Y|d})$.

We note that for $q = 2$, $L^q(P) = L^2(P)$ is a Hilbert space equipped with the inner product $\langle f, g \rangle_{L^2(P)} = E_P[f(O)g(O)]$. We similarly define the $L^2(P_Y)$ inner product via $\langle f, g \rangle_{L^2(P_Y)} = E_P[f(Y)g(Y)]$ and  $L^2(P_{Y|d})$ inner product via $\langle f, g \rangle_{L^2(P_{Y|d})} = E_P[f(Y)g(Y)|D=d]$.

For $\mc{T}$ a bounded linear operator between Banach spaces $\mc{B}_1 \rightarrow \mc{B}_2$, we let $\mathscr{N}(\mc{T}) \coloneqq \{b_1 \in \mc{B}_1: \mc{T}b_1 = 0\}$ denote the null space of $\mc{T}$, and $\mathscr{R}(\mc{T}) \coloneqq \{b_2 \in \mc{B}_2: b_2 = \mc{T}b_1 \text{ for some } b_1 \in \mc{B}_1\}$ denote the range of $\mc{T}$. When $\mc{B}_1$, $\mc{B}_2$ are both Hilbert spaces, we let $\mc{T}^*$ denote the Hermitian adjoint of $\mc{T}$. In this case, for $b_2 \in \mathscr{R}(\mc{T})$ the set $\{b_1 \in \mc{B}_1: \mc{T}b_1 = b_2\}$ contains a unique element of smallest norm, which we call the minimum norm solution to $\mc{T}b_1 = b_2$ and denote with $b_1^\dagger = \mc{T}^\dagger b_2$; we refer to $\mc{T}^\dagger$ as the Moore-Penrose pseudo-inverse of $\mc{T}$. Throughout, a superscript $\dagger$ on a solution to an operator equation indicates the minimum norm solution. For $\mc{A}$ a linear subspace of a Hilbert space $\mc{H}$, we let $\overline{\mc{A}}$ denote the closure of $\mc{A}$.

For $Q$ a distribution, we let $\textsf{Supp}[Q]$ denote its support, and for $Y$ a measurable function of $O$ we write $\textsf{Supp}[Y;P]$ for the support of the distribution of $Y$ under $P$. Finally, for random elements $A_n$ and $A$, we write $A_n \rightsquigarrow A$ to denote convergence in distribution, and write the law under which $A_n$ is generated above the arrow when we wish to emphasize it.

Throughout, we will make use of the Hellinger distance between two distributions for a random variable, which we define below.
\begin{definition}
    \label{def:hellinger-dist}
    For $P, P'$ two probability measures on a measurable space $(\mc{X}, \mc{F})$, the (squared) Hellinger distance $H^2(P, P')$ is defined as
    \begin{align}
        \label{eq:hellinger-dist}
        H^2(P, P') \coloneqq \frac{1}{2}\int_{\mathcal{X}}\left\{\sqrt{p(x)}- \sqrt{p'(x)} \right\}^2d\lambda(x)
    \end{align}
    for $\lambda$ any dominating measure, for example $ \lambda \coloneqq P + P'$.
\end{definition}

\subsection{The Target of Inference}
\label{subsec:target-of-inference}

We now define our target of inference $\varphi$. Let the functions
\begin{align*}
    \widetilde{m}&: \mathscr{M}(X, D) \rightarrow \mathscr{M}(W, X, D),  \quad h \mapsto \widetilde{m}(\cdot;h)\\
    m&: \mathscr{M}(Z, D) \rightarrow \mathscr{M}(W, Z, D),\quad  q \mapsto {m}(\cdot;q)
\end{align*}
be known linear maps, where $\mathscr{M}(Y)$ denotes the set of measurable functions of $Y$ for any $Y$ a measurable function of $O$. We will work at laws $P$ such that the functionals
\begin{align}
    \label{eq:mtilde-m-def}
    L^2(P_{X,D}) \rightarrow \mathbb{R}, h \mapsto E_{P}[\widetilde{m}(W, X, D;h)] && L^2(P_{Z,D}) \rightarrow \mathbb{R}, q \mapsto E_{P}[{m}(W, Z, D;q)]
\end{align}
are bounded with Riesz representers $\alpha_P \in L^2(P_{X,D})$ and $r_P \in L^2(P_{Z,D})$ respectively, i.e.
\begin{align*}
    E_{P}[\widetilde{m}(W, X, D;h)] = E_{P}[h(X, D)\alpha_P(X, D)], && E_{P}[{m}(W, Z, D;q)] = E_{P}[q(Z, D)r_P(Z, D)]
\end{align*}
where the equalities hold for all $h \in L^2(P_{X, D})$ and $q \in L^2(P_{Z, D})$, respectively.
We define the conditional expectation operators
\begin{align*}
    &\mc{E}_{P}: L^2(P_{X,D}) \rightarrow L^2(P_{Z,D}), \mc{E}_{P}h = E_{P}[h(X, D)|Z= \cdot, D = \cdot],\\
    & \mc{E}_{P}^*: L^2(P_{Z,D}) \rightarrow L^2(P_{X,D}), \mc{E}^*_{P}q = E_{P}[q(Z, D)|X= \cdot, D=\cdot].
\end{align*}
We also define the conditional expectation operators operating within fixed levels of $d$, i.e.
\begin{align*}
    &\mc{E}_{P;d}: L^2(P_{X|d}) \rightarrow L^2(P_{Z|d}), \mc{E}_{P;d}h = E_{P}[h(X)|Z= \cdot, D = d],\\
    & \mc{E}_{P;d}^*: L^2(P_{Z|d}) \rightarrow L^2(P_{X|d}), \mc{E}^*_{P;d}q = E_{P}[q(Z)|X= \cdot, D=d]
\end{align*}
where we recall that $D$ is finite-discrete.

We then define our target parameter as
\begin{align*}
    \varphi(P) = E_{P}[\widetilde{m}(W, X, D;h_P)] = E_{P}[{m}(W, Z, D;q_P)] = E_{P}[h_P(X, D)q_P(Z, D)]
\end{align*}
for $h_P \in L^2(P_{X,D})$, $q_P \in L^2(P_{Z,D})$ solving the primal and dual inverse problems
\begin{align}
    \label{eq:primal-dual-inverse-problem}
    \mc{E}_{P}h = r_{P}, && \mc{E}_{P}^*q = \alpha_{P},
\end{align}
respectively, where these solutions are assumed to exist over $\mc{P}$. A slightly larger class of parameters was studied in \cite{bennettSourceConditionDouble2023b, bennett_inference_2023, olivas-martinez_source-condition_2025}, among others. These authors allow $\widetilde{m}, m$ to depend on $O$ instead of $(W, X, D)$ and $(W, Z, D)$, respectively. We consider only a smaller class of parameters, but conjecture the forthcoming results hold over the larger class. The parameter $\varphi(P)$ is well defined even when $h_P$ and $q_P$ are not uniquely defined, i.e. when there exist multiple solutions to the primal or dual inverse problem. Prior work (\cite{bennettSourceConditionDouble2023b, olivas-martinez_source-condition_2025, bennett_inference_2023}) has also shown that, at $P$ satisfying the existence of a primal and dual solution,
\begin{align}
    \varphi^1(O;h, q, \varphi) \coloneqq \widetilde{m}(W, X, D;h) + m(W,Z, D;q) - h(X, D)q(Z, D) - \varphi \label{eq:neyman-orthogonal-ee}
\end{align}
is a Neyman-orthogonal estimating equation, in that its mean is Gateaux differentiable in the nuisance functions $h$, $q$ at $(h_P,q_P)$ with Gateaux derivative 0. This estimating function admits a doubly robust expansion, i.e.
\begin{align*}
    E_{P}[\varphi^1(O;h, q, \varphi(P))] = -E_{P}\left[\{h(X, D) - h_{P}(X, D)\}\{q(Z, D) - q_{P}(Z, D)\}\right],
\end{align*}
for $h_P$, $q_P$ any solutions to \eqref{eq:primal-dual-inverse-problem}. As such, under certain conditions on the estimation rates of $h_{P}$ and $q_{P}$, the estimator
\begin{align}
    \label{eq:dml-estimator}
    \widehat{\varphi}_n \coloneqq \frac{1}{n}\sum_{i=1}^n\widetilde{m}(W_i, X_i, D_i;\widehat{h}_n) + m(W_i, Z_i, D_i;\widehat{q}_n) - \widehat{h}_n(X_i, D_i)\widehat{q}_n(Z_i, D_i)
\end{align}
(or rather a cross-fit version of it) is asymptotically linear with influence function $\varphi^1(O; h_P, q_P, \varphi(P))$. Remarkably, such estimators exist even when the operators $\mathcal{E}_P, \mathcal{E}_P^*$ are ill-posed, in the sense that they do not have continuous inverses, leading to instability in estimating $h_{P}$, $q_{P}$. For $X, Z$ finite discrete, this ill-posedness arises when there exist multiple solutions to the inverse problems. For $X, Z$ continuously valued, this ill-posedness arises under mild conditions on $P$, for example when the joint density of $(X, Z)|D=d$ with respect to the Lebesgue measure is bounded above and below by a constant (See the discussion in \cite{carrasco_chapter_2007} following Assumption A.1) for each $d$.

{
This class of estimands includes several target parameters central in causal inference, including a weighted average of the dose response curve under an additive non-parametric instrumental variable model (NPIV), the average treatment effect under proximal causal inference, and the mean outcome among nonrespondents under a shadow variable model, as we now discuss. }

{ \begin{example}[Additive Nonparametric Instrumental Variables]
    \label{example:npiv}

    Suppose that \(W, X, Z\) are each one-dimensional random variables, where
    \(W\) is an outcome, \(X\) is a treatment, and \(Z\) is an instrumental
    variable. For simplicity, suppose that \(D\) is absent, though finite-discrete
    covariates could be incorporated through \(D\). Consider the additive
    structural causal model
    \begin{equation}
    \label{eq:npiv-model}
        W(x,z)
        =
        h_{P;\mathrm{cdr}}(x)
        +
        \epsilon_{W;P},
        \qquad
        \text{with}
        \qquad
        E_P(\epsilon_{W;P}\mid Z)=0,
    \end{equation}
    where \(h_{P;\mathrm{cdr}}\) is unspecified and \(W(x,z)\) denotes the counterfactual outcome that would be observed if
    \(X\) and \(Z\) were set to \(x\) and \(z\). This model imposes that \(Z\)
    has no direct effect on the counterfactual outcome and does not modify the
    effect of \(X\). Hence \(W(x,z)\) does not depend on \(z\), and we write
    \(W(x)\) for this common counterfactual value. Since
    \(E_P(\epsilon_{W;P})=0\), then
    \begin{equation}
    \label{eq:dose-response-curve}
        h_{P;\mathrm{cdr}}(x)
        =
        E_P\{W(x)\}.
    \end{equation}

Thus, the function \( x \mapsto h_{P;\mathrm{cdr}}(x)\) is interpreted as the causal
    dose-response curve of treatment $X$.
    By consistency, the observed outcome satisfies
    \(W=W(X,Z)=h_{P;\mathrm{cdr}}(X)+\epsilon_{W;P}\). Therefore,
    \begin{equation}
    \label{eq:primal-npiv-eq}
        E_P(W\mid Z)
        =
        E_P\{h_{P;\mathrm{cdr}}(X)\mid Z\}.
    \end{equation}
    Thus the causal
    dose-response curve is a solution to the primal inverse problem
    \begin{equation*}
        \mathcal E_P h
        =
        r_P
        \qquad \text{with} \qquad
        r_P(Z):=E_P(W\mid Z).
    \end{equation*}

    Given a known weight function \(\psi\in L^2(P_X)\), the target parameter is
    the weighted mean of the causal dose-response curve,
    \begin{equation*}
        \varphi(P)
        :=
        E_P\{\psi(X)h_{P;\mathrm{cdr}}(X)\}.
    \end{equation*}
    If the primal inverse problem has multiple solutions, this weighted mean is
    identified under the further condition that
    \(
        \psi\in \mathscr R(\mathcal E_P^*),
    \)
    equivalently, there exists \(q_P\in L^2(P_Z)\) solving the dual inverse
    problem
    \begin{equation*}
        \psi(X)
        =
        E_P\{q_P(Z)\mid X\}.
    \end{equation*}
    Under this condition, the value of \(E_P\{\psi(X)h(X)\}\) is the same for
    all solutions \(h\) of the primal inverse problem, and hence
    \(\varphi(P)\) is identified.

    This example fits the general framework with
    \begin{align*}
        \widetilde m(W,X;h)
        &:= \psi(X)h(X),
        &m(W,Z;q)
        &\coloneqq Wq(Z),\\
        \alpha_P(X)
        &\coloneqq \psi(X),
        &r_P(Z)
        &:= E_P(W\mid Z).
    \end{align*}
\end{example}}

{ \begin{example}[Proximal Causal Inference]
    \label{example:proximal}
    Suppose that  $W, X, Z$ are 1-dimensional random variables and $D$ is binary. $W$ is an outcome of interest, $X$ is an outcome proxy, $Z$ a treatment proxy, and $D$ is a binary treatment. For simplicity again we do not include covariates, but $D$ could be extended to accommodate finite discrete covariates. Suppose that there exist solutions $h_P$, $q_P$, to the primal and dual inverse problems
    \begin{align*}
        (\mc{E}_{P}h)(Z,D) = E_{P}[W|Z, D], && (\mc{E}_{P}^*q)(X,D) = \frac{2D-1}{P(D|X)}.
    \end{align*}
    The parameter of interest is
    \begin{align*}
        \varphi(P) \coloneqq E_P[h_P(X, 1) - h_P(X, 0)].
    \end{align*}
    The statistical parameter $\varphi$ fits into the above framework with
    \begin{align*}
        \widetilde{m}(W, X, D;h) \coloneqq& h(X, 1)- h(X, 0),  & m(W, Z, D;q) \coloneqq& Wq(Z, D)\\
        \alpha_P(X, D) \coloneqq&\frac{2D-1}{P(D|X)} , & r_{P}(Z, D) \coloneqq& E_{P}[W|Z, D].
    \end{align*}
    Under certain causal assumptions \cite{cui_semiparametric_2024} showed that  $\varphi(P)$ equals the average treatment effect (ATE) of $D$ on $W$, i.e. $E[W(d=1) - W(d=0)]$ where $W(d)$ is the counterfactual outcome when $D$ is set to $d$.
\end{example}

}

{

\begin{example}[Shadow variable model]
    \label{example:shadow}
Suppose that \(X\) and \(Z\) are one-dimensional random variables and that
\(D\) is binary. In this example, \(W\) is absent. Suppose that there exist
solutions \(h_P\) and \(q_P\) to the primal and dual inverse problems
\begin{align*}
    (\mathcal E_P h)(Z,D)
    &=
    DZ,
    &
    (\mathcal E_P^* q)(X,D)
    &=
    D\frac{P(D=0\mid X)}{P(D=1\mid X)}.
\end{align*}
The parameter of interest is
\begin{equation*}
    \varphi(P)
    :=
    E_P\{(1-D)h_P(X,1)\}.
\end{equation*}
This parameter fits into the general framework with
\begin{align*}
    \widetilde m(X,D;h)
    &:=
    (1-D)h(X,1),
    &
    m(Z,D;q)
    &:=
    DZq(Z,D),
    \\
    \alpha_P(X,D)
    &:=
    D\frac{P(D=0\mid X)}{P(D=1\mid X)},
    &
    r_P(Z,D)
    &:=
    DZ.
\end{align*}

Interest in the parameter \(\varphi(P)\) arises because, as shown by \cite{miao_identification_2018} and
\cite{liNonparametricInferenceMean2022}, \(\varphi(P)/P(D=0)\) equals the mean of \(Z\) among nonrespondents under
the following shadow-variable model. In this
model, \(Z\) is the outcome of interest, observed only in a subsample; \(X\) is
an outcome proxy that is always observed; and \(D\) is a response indicator,
with \(D=1\) if and only if \(Z\) is observed. For simplicity, we again omit
covariates, although the setup could be extended to accommodate finite-discrete
covariates.

The shadow-variable model can be described as follows. For nonrespondents, that
is, conditional on \(D=0\), the pair \((Z,X)\) is drawn from an unknown law
\(Q\), but only \(X\) is observed. For respondents, that is, conditional on
\(D=1\), the pair \((Z,X)\) is observed; \(Z\) is drawn from an unknown law
\(U\), and, given \(Z=z\), \(X\) is drawn from the conditional law
\(Q_{X\mid Z=z}\) induced by the joint law \(Q\). The model implicitly assumes
that \(U\) is absolutely continuous with respect to the marginal law \(Q_Z\) of
\(Z\) under \(Q\).

Under this model, it can be shown that
\begin{equation*}
    q_P(Z,D)
    =
    D
    \frac{P(D=0)}{P(D=1)}
    \frac{dQ_Z}{dU}(Z)
\end{equation*}
solves the dual inverse problem
\begin{equation}
\label{eq:dual2}
    E_P\{q_P(Z,D)\mid X,D\}
    =
    \alpha_P(X,D),
\end{equation}
and that
\begin{equation*}
      E_Q(Z)
    =
    \frac{E_P\{m(Z,D;q_P)\}}{P(D=0)}.
\end{equation*}
\end{example}

}

\subsection{Review of Semiparametric Theory}
\label{subsec:semiparametrics}
Before giving our results, we give a brief review of semiparametric efficiency theory, including the definition of a tangent space, regular estimators, and pathwise differentiability. See \cite{bickelEfficientAdaptiveEstimation1998a} for a textbook treatment of this material.

Suppose that distributions in a model $\mc{P}$ are dominated by the measure $\lambda$, as will be assumed throughout this work. Then, we say a one-dimensional parametric submodel $\{P_t:t\in[0, \delta)\} \subseteq \mc{P}$ passing through $P_0$ with $P_{t=0}  = P_0 $ is quadratic mean differentiable (QMD) at $P_0$ if there exists $s \in L^2_0(P_0)$ with
\begin{align*}
    \lim_{t\rightarrow 0} \int \left\{\frac{\sqrt{p_t(o)} - \sqrt{p_0(o)}}{t} - \frac{1}{2}s(o)\sqrt{p_{0}(o)}\right\}^2 d\lambda(o) = 0.
\end{align*}
The function $s$ is the so called score, and generalizes the classical pointwise derivative of the log-likelihood function.

We let $\mathscr{P}(P_0, \mc{P})$ denote the set of all one dimensional QMD at $P_0$ submodels of $\mc{P}$. We define the maximal tangent set as the set of scores of submodels $\mathscr{P}(P_0, \mc{P})$ and the maximal tangent space at $P_0$ relative to $\mc{P}$ as the closed linear span of the tangent set. We denote the maximal tangent set with $\dot {\mc{T}}(P_0, \mc{P})$ and the maximal tangent space with $\mc{T}(P_0, \mc{P})$. Note that the latter is a subspace of $L^2_0(P_0)$. When the model consists of the set of all distributions for $O$, the maximal tangent space is equal to $ L^2_0(P_0)$. Consequently, whenever $\mc{T}(P_0, \mc{P}) = L^2_0(P_0)$, the model is said to be (locally) non-parametric at $P_0$. When $\mc{T}(P_0, \mc{P})$ is finite dimensional, the model is said to be (locally) parametric at $P_0$. When $\mc{T}(P_0, \mc{P})$ is neither equal to $L^2_0(P_0)$ nor finite dimensional, the model is said to be (locally) semiparametric at $P_0$.

The set of submodels $\mathscr{P}(P_0, \mc{P})$ can be used to state an efficiency theory for regular estimators of parameters that are pathwise differentiable relative to $\mathscr{P}(P_0, \mc{P})$. We first introduce pathwise differentiability for parameters that take values in $\mathbb{R}$, then extend the definition to parameters that take values in a real separable Hilbert space. A real valued parameter $P \mapsto \phi(P) \in \mathbb{R}$ is pathwise differentiable relative to $\mathscr{P}(P_0, \mc{P})$ if there exists a bounded linear functional $\dot \phi_{P_0}: \mc{T}(P_0, \mc{P}) \rightarrow \mathbb{R}$ such that, for all $\{P_t:t\in [0, \delta)\} \in \mathscr{P}(P_0, \mc{P})$,
\begin{align}
    \label{eq:local-parameter-def}
    \dot\phi_{P_0}[s] = \left.\frac{d}{dt}\phi(P_t)\right\vert_{t=0}
\end{align}
where $s$ is the score of the submodel. The functional $\dot\phi_{P_0}$ is called the local parameter. By the Riesz representation theorem, existence of the local parameter is, for real valued parameters, equivalent to the existence of a gradient $\phi^1_{P_0} \in L^2_0(P_0)$ satisfying
\begin{align*}
    \langle s, \phi^1_{P_0}\rangle_{L^2(P_0)}  = \left.\frac{d}{dt}\phi(P_t)\right\vert_{t=0}.
\end{align*}

Recall that an estimator $\widehat{\phi}_n$  is any measurable function of $(O_i)_{i=1}^n$. We say that an estimator $\widehat{\phi}_n$ of $P \mapsto \phi(P)$ is regular relative to $\mathscr{P}(P_0, \mc{P})$ if for any sequence $t_n \rightarrow t$ with $t_n, t \geq 0$,
\begin{align*}
    \sqrt{n}\left\{\widehat{\phi}_n - \phi(P_{t_n / \sqrt{n}})\right\}  \overset{P_{t_n / \sqrt{n}}}{\rightsquigarrow} \mathbb{L}
\end{align*}
for every $\{P_t:t \in [0, \delta)\} \in \mathscr{P}(P_0, \mc{P})$ where the limit distribution $\mathbb{L}$ does not depend on the submodel $\{P_t:t \in [0, \delta)\}$ (See \cite{bickelEfficientAdaptiveEstimation1998a} Section 2.2 Definition 7 and Section 3.1). We often drop "relative to $\mathscr{P}(P_0, \mc{P})$" when discussing tangent spaces, pathwise differentiability, and regularity. In this work, these should all be taken relative to $\mathscr{P}(P_0, \mc{P})$.

Regular estimators achieve their limits in a locally uniform sense, where uniformity is with respect to the submodels in $\mathscr{P}(P_0, \mc{P})$. The existence of regular estimators is intimately related to the pathwise differentiability of the parameter $\phi$. If there exists a regular and asymptotically linear estimator of $\phi$ at $P_0$, then $\phi$ is pathwise differentiable at $P_0$ \parencite{vaart_differentiable_1991}. On the other hand, if $\phi$ is pathwise differentiable at $P_0$, under some conditions and in particular if a consistent estimator of the gradient exists, there exist regular and asymptotically linear estimators of $\phi$ at $P_0$ \parencite{klaassen_consistent_1987-2}. Recall that an estimator is asymptotically linear at $P_0$ if there exists a function $\phi^1_{P_0} \in L^2_0(P_0)$ such that
\begin{align*}
    \widehat{\phi}_n - \phi(P_0) = \frac{1}{n}\sum_{i=1}^n \phi^1_{P_0}(O_i) + o_p(n^{-1/2}).
\end{align*}
We call $\phi^1_{P_0}$ the influence function of $\widehat{\phi}_n$. We use the same notation for gradients and influence functions because, if $\phi$ is pathwise differentiable at $P_0$ and $\widehat{\phi}_n$ is regular, $\phi^1_{P_0}$ is a gradient.

The tangent space $\mc{T}(P_0, \mc{P})$ provides lower bounds on the asymptotic variance of any regular estimator sequence $\widehat{\phi}_n$ of pathwise differentiable parameters $P \mapsto \phi(P)$ at $P_0$. These lower bounds are formalized in the Convolution and Local Asymptotic Minimax Theorems (Theorems 25.20 and 25.21 of \cite{vaart_asymptotic_1998-2}). In short, the minimum asymptotic variance of any regular estimator sequence at $P_0$ is $E_{P_0}[\{\phi^1_{P_0;eff}(O)\}^2]$ where $\phi^1_{P_0;eff}$ is the unique gradient of $\phi$ that belongs to the tangent space $\mc{T}(P_0, \mc{P})$. If an estimator is asymptotically linear with influence function equal to $\phi^1_{P_0;eff}(O)$, i.e.
\begin{align*}
    \widehat{\phi}_n - \phi(P_0) = \frac{1}{n}\sum_{i=1}^n\phi^1_{P_0;eff}(O_i) + o_p(n^{-1/2}),
\end{align*}
it is regular and achieves the asymptotic variance lower bound. We call such estimators efficient. When the model tangent space is $L^2_0(P_0)$, all regular and asymptotically linear estimators share a single influence function, and hence are asymptotically equivalent; no regular estimator can exploit restrictions in $\mc{P}$ to gain asymptotic efficiency. When the tangent space is strictly included in $L^2_0(P_0)$, some regular and asymptotically linear estimators will be more efficient than others, if they exist. The efficient estimators are exactly those that best exploit the modeling restrictions to estimate the parameter of interest.

The above definitions apply to real-valued parameters. Many of the definitions extend readily to real separable Hilbert space-valued parameters. See \cite{luedtke_one-step_2024} for a detailed discussion of efficiency theory and estimation for such parameters. In particular, let $\mc{H}$ denote a fixed real separable Hilbert space, for example $L^2(P_0)$. Suppose now that $\phi: \mc{P} \rightarrow \mc{H}$ takes values in $\mc{H}$. Then, we say that $\phi$ is pathwise differentiable relative to $\mathscr{P}(P_0, \mc{P})$ if there exists a bounded linear operator $\dot \phi_{P_0}: \mc{T}(P_0, \mc{P}) \rightarrow \mc{H}$ such that, for all QMD at $P_0$ submodels $\{P_t:t\in [0, \delta)\} \in \mathscr{P}(P_0, \mc{P})$, \eqref{eq:local-parameter-def} holds where $s$ is again the score of the submodel. As in the real valued setting, we call $\dot \phi_{P_0}$ the local parameter of $\phi$ at $P_0$. In contrast to the real valued setting, when $\mc{H}$ is infinite dimensional, there does not exist an equivalent formulation in terms of gradients $\phi^1_{P_0} \in L^2_0(P_0)$ as the Riesz-representation Theorem does not apply.

\subsection{Parameter Smoothness and Implications for Inference}

A large body of classical asymptotic theory is devoted to the connection between the smoothness of a parameter $\phi$ at $P_0$ in a statistical model $\mc{M}$, and the possibility or impossibility of inference about $\phi$ within the model $\mc{M}$ \parencite{bahadur_nonexistence_1956, donoho_one-sided_1988, donoho_geometrizing_1991, donoho_statistical_1994, gleser_nonexistence_1987, low_nonparametric_1997, cai_adaptation_2004, pfanzagl_nonexistence_1998, dufour_impossibility_1997, potscher_lower_2002}. Here, we briefly review the implications of two measures of smoothness: Hellinger continuity, and pathwise differentiability. Briefly, Hellinger discontinuity implies that no estimator is locally uniformly consistent, and that no confidence interval can be locally uniformly honest while having diameter shrinking to zero. When the local modulus of continuity is infinite, as it is for our target parameter under the conditions of \cref{subsec:failure-of-continuity}, every locally uniformly honest confidence interval must have infinite expected diameter at every sample size. Failure of pathwise differentiability implies that all asymptotically linear estimators are irregular, exhibiting non-vanishing local bias when data are generated from some local alternatives, leading to under-coverage of Wald-type confidence intervals under those alternatives. We leave a more in depth discussion to \cref{sup-sec:smoothness-implications}. The forthcoming results should be taken as standard; our contribution is not their derivation, but their application to the parameter defined in \cref{subsec:target-of-inference}.

\subsubsection{Implications of Hellinger Discontinuity}
\label{subsubsec:implications-discontinuity-main}

To fix ideas, Hellinger continuity is defined as follows.
\begin{definition}
    \label{def:hellinger-continuous}
    We say that a real valued functional $\phi:\mc{M} \rightarrow \mathbb{R}$ is Hellinger continuous at $P \in \mc{M}$ if, for all $\epsilon > 0$, there exists $\delta > 0$ such that for all $P' \in \mc{M}$ such that $H(P', P) < \delta $, it holds that $|\phi(P) - \phi(P')| < \epsilon$.
\end{definition}

We begin by giving general finite-sample consequences of Hellinger discontinuity. The following finite-sample bounds follow from elementary two-point testing arguments. They are in the same spirit as the lower
bounds for confidence intervals developed by \cite{low_nonparametric_1997} and
\cite{cai_adaptation_2004}, and are closely related to the
modulus-of-continuity approach of \cite{donoho_geometrizing_1991,donoho_statistical_1994}.

Given a collection $\mathcal{M}$ of probability laws defined on a common measurable space, the local modulus of continuity at $P_0 \in \mathcal{M}$ of a parameter $\phi$  over $\mathcal{M}$ for a Hellinger radius $r > 0$ is defined as:
\begin{equation*}
\omega_{P_0}(r) = \sup \{ |\phi(P) - \phi(P_0)| : H(P, P_0) < r, \, P \in \mathcal{M} \}.
\end{equation*}
The local modulus of continuity quantifies the local sensitivity, i.e. instability to small changes in the data generating distribution. A rapidly vanishing modulus (as $r$ goes to 0) indicates local stability. On the other hand, a modulus that fails to vanish as $r$ goes to 0 indicates discontinuity at $P_0$. Indeed, it is immediate from the definition of continuity that at $P_0$ such that $\phi$ is Hellinger discontinuous, there exists $\epsilon > 0$ such that $\omega_{P_0}(r) > \epsilon$ for every $r > 0$ (see the definition of Hellinger continuity in \cref{def:hellinger-continuous}). This fact combined with a
standard two-point argument leads to the conclusion that, at every fixed sample size, any
confidence interval with honest coverage over a sufficiently small Hellinger
neighborhood of \(P_0\) must have length bounded away from zero with positive probability
under \(P_0\), and that the corresponding local uniform estimation error for point
estimation is also bounded away from 0. Throughout, \(P^n\) denotes the law of \(n\) i.i.d. draws from \(P\), \(C_n\) denotes a random interval whose endpoints are measurable functions of \(n\) i.i.d. draws, and \(|C_n|\) denotes its diameter, allowing the value \(+\infty\). We say that \(C_n\) has coverage at least \(1-\alpha\) over a set \(\mathcal{A} \subseteq \mc{M}\) if \(\inf_{P \in \mathcal{A}} P^n\{\phi(P) \in C_n\} \geq 1-\alpha\). For the following result, define
\begin{equation}
\label{eq:rn-tau}
    r_{n,\tau}
    :=
    \left\{
        1-(1-\tau^2)^{1/(2n)}
    \right\}^{1/2}.
\end{equation}
\begin{lemma}
    \label{lemma:finite-sample-discontinuity}
    Suppose that $P \mapsto \phi(P)$ is Hellinger discontinuous at $P_0$. Then, there exists $\epsilon > 0$ such that \(\omega_{P_0}(r)> \epsilon\) for every \(r>0\). Fix \(0<\alpha<1/2\), \(n\in\mathbb N\), and \(\tau\in(0,1-2\alpha)\). If \(C_n\) has coverage at least \(1-\alpha\) over
    the ball \(\{P\in\mathcal M:H(P,P_0)<r_{n,\tau}\}\), then,
    \begin{equation}
    \label{eq:tail-any-M}
        P_0^n\{|C_n|\ge \epsilon\}
        \ge
        1-2\alpha-\tau.
    \end{equation}
    If \(\omega_{P_0}(r) = \infty\) for every \(r>0\), then
    \begin{equation}
    \label{eq:infinite-length-positive-prob}
        P_0^n\{|C_n| = \infty\}
        \ge
        1-2\alpha-\tau,
        \qquad
        E_{P_0^n}\{|C_n|\}
        =  \infty.
    \end{equation}
    Furthermore, for every estimator \(\widehat\phi_n\) of \(\phi\), every \(n \in \mathbb{N}\), and every \(\tau\in(0,1)\),
    \begin{equation}
    \label{eq:finite-sample-estimator-bad-somewhere}
        \sup_{P\in\mathcal M:H(P,P_0)<r_{n,\tau}}
        P^n\{|\widehat\phi_n-\phi(P)|\ge \epsilon/2\}
        \ge
        \frac{1-\tau}{2},
    \end{equation}
    and if \(\omega_{P_0}(r) = \infty\) for every \(r>0\), \eqref{eq:finite-sample-estimator-bad-somewhere} holds with $\epsilon/2$ replaced by any $M > 0$.
\end{lemma}
It follows from the preceding lemma that, under Hellinger discontinuity, any
finite-sample interval that is honest over a sufficiently small Hellinger
neighborhood of a discontinuity point \(P_0\) must assign positive probability, under \(P_0\), to
having length greater than $\epsilon$, with $\epsilon$ independent of $n$. In particular, this immediately rules out Wald confidence intervals as having local coverage greater than or equal to their stated nominal level for sufficiently large $n$, as the diameter of Wald intervals converges to 0 as $n\rightarrow \infty$. When \(\omega_{P_0}(r) = \infty\) for every \(r>0\), Wald intervals will locally undercover for every $n$.

Similarly, by \eqref{eq:finite-sample-estimator-bad-somewhere}, no estimator can have uniformly controlled finite-sample error less than $\epsilon/2$ over the
local Hellinger neighborhoods $\left\{ P:H(P,P_0)<r_{n,\tau}\right\}$, and when \(\omega_{P_0}(r) = \infty\) for every \(r>0\), no estimator can have uniformly controlled finite-sample error at all.

We now study the implications of these finite-sample statements for locally uniform asymptotic inference at points $P_0$ where $\phi$ is Hellinger discontinuous. Let \(\delta_n\downarrow0\) be any
sequence of positive numbers, and
define the shrinking Hellinger neighborhoods
\begin{equation*}
    \mathcal U_n(P_0)
    :=
    \{P\in\mathcal M:H(P,P_0)<\delta_n\}.
\end{equation*}
Applying \cref{lemma:finite-sample-discontinuity} inside the ball of radius \(\min\{\delta_n,r_{n,\tau}\}\) shows that, for every estimator
\(\widehat\phi_n\), and every \(\tau\in(0,1)\),
\begin{equation}
\label{eq:no-local-uniform-estimation}
    \liminf_{n\to\infty}
    \sup_{P\in\mathcal U_n(P_0)}
    P^n\{|\widehat\phi_n-\phi(P)|\ge \epsilon/2\}
    \ge
    \frac{1-\tau}{2}.
\end{equation}
Since \(\tau\) can be chosen arbitrarily small, the lower bound can be made
arbitrarily close to \(1/2\).  In particular, no sequence of estimators is
uniformly consistent over \(\mathcal U_n(P_0)\).

Similarly, suppose that \(C_n\) is a sequence of confidence intervals satisfying
locally uniform asymptotic coverage over \(\mathcal U_n(P_0)\):
\begin{equation}
\label{eq:local-uniform-asymptotic-coverage}
    \liminf_{n\to\infty}
    \inf_{P\in\mathcal U_n(P_0)}
    P^n\{\phi(P)\in C_n\}
    \ge
    1-\alpha,
    \qquad
    0<\alpha<1/2.
\end{equation}
Then, for every \(\tau\in(0,1-2\alpha)\),
\begin{equation}
\label{eq:no-local-uniform-ci-tightness}
    \liminf_{n\to\infty}
    P_0^n\{|C_n|\ge \epsilon\}
    \ge
    1-2\alpha-\tau.
\end{equation}
Thus the intervals cannot be locally uniformly valid and at the same time have length converging
to zero. Since \(\tau\) can be taken arbitrarily small, the right-hand side can
be made arbitrarily close to \(1-2\alpha\).

In summary, the failure of Hellinger continuity has three consequences.  First,
for every fixed sample size, honest confidence intervals over sufficiently small
Hellinger neighborhoods must have length bounded away from zero with positive probability
under \(P_0\).  Second, point estimators have error bounded away from zero somewhere
in those same neighborhoods.  Third, as a result, locally uniform asymptotic
estimation and locally uniform asymptotic confidence intervals with diameter that shrinks to zero are impossible,
even over neighborhoods whose radii shrink to zero arbitrarily fast.

\subsubsection{Implications of Failure of Pathwise Differentiability}
\label{subsubsec:implications-non-diff-main}

Pathwise differentiability, formally defined in \cref{subsec:semiparametrics}, is a strictly stronger smoothness condition than Hellinger continuity, see \cref{remark:failure-of-continuity-fails-diff}. While failure of pathwise differentiability has less severe consequences for inference than failure of Hellinger continuity, it nevertheless implies poor asymptotic coverage of Wald-type intervals for asymptotically linear estimators under sampling from some local alternative distributions. We explore these implications below.

In what follows, let $\mathscr{P}_{\Theta} \coloneqq \{ \{P_{t, \theta} \in \mc{M}: t\in [0, \delta_{\theta})\}: \theta \in \Theta\}$ be a collection of QMD at $P_0$ submodels of $\mc{M}$ where $\Theta$ is a (possibly infinite dimensional) index set. Let $s_{\theta}$ be the score of the submodel $\{P_{t, \theta}:t\in[0, \delta_\theta)\}$ for each $\theta \in \Theta$.

To analyze $\phi$ under sampling from local alternatives of the form $P_{1/\sqrt{n}, \theta}$, one decomposes the root-$n$ error as
\[
\sqrt n\{\widehat\phi_n-\phi(P_{1/\sqrt n, \theta})\}
=
\sqrt n\{\widehat\phi_n-\phi(P_0)\}
-
\sqrt n\{\phi(P_{1/\sqrt n, \theta})-\phi(P_0)\}.
\]
If $\widehat\phi_n$ is root-$n$ consistent at $P_0$, the first term is
$O_{P_0}(1)$ and, by contiguity, remains $O_P(1)$ under
$P_{1/\sqrt n, \theta}$. Hence the root-$n$ error can remain $O_P(1)$ only if
the second term is $O(1)$. If, moreover, $\widehat\phi_n$ is
asymptotically linear at $P_0$ with influence function
$\phi^1_{P_0}$, then Le Cam's third lemma (\cite{vaart_asymptotic_1998-2}, Theorem 6.6 and Example 6.7; \cite{lecamAsymptoticsStatistics2000}, Chapter 3) implies that, under
$P_{1/\sqrt n, \theta}$, the first term has an asymptotic normal distribution with mean
$E_{P_0}\{\phi^1_{P_0}(O)s_{\theta}(O)\}$.
If $t\mapsto\phi(P_t)$ is differentiable at zero, then
$\sqrt n\{\phi(P_{1/\sqrt n, \theta})-\phi(P_0)\}$ converges to
$\,\frac{d}{dt}\phi(P_{t, \theta})|_{t=0}$. Thus the root-$n$ estimation error
under $P_{1/\sqrt n, \theta}$ is asymptotically normal with mean, i.e. asymptotic $\sqrt{n}$-scaled bias,
\[
\left[
E_{P_0}\{\phi^1_{P_0}(O)s_{\theta}(O)\}
-
\left.\frac{d}{dt}\phi(P_{t, \theta})\right|_{t=0}
\right].
\]
For Wald intervals centered at $\widehat\phi_n$ to retain their nominal
asymptotic coverage along every such local alternative, it is therefore necessary that this mean
vanish for every QMD submodel. This requires the derivatives of
$\phi(P_t)$ to be representable as a single bounded linear functional
of the corresponding scores, namely $s \mapsto E_{P_0}\{\phi^1_{P_0}(O)s(O)\}$. In particular, it requires that the parameter $\phi(P)$ be pathwise differentiable at $P_0$. When pathwise differentiability fails, there are local alternatives along which an asymptotically linear estimator has
nonvanishing root-$n$ bias, and the corresponding Wald intervals centered at asymptotically linear estimators
undercover asymptotically. Failure of pathwise differentiability relative to $\mathscr{P}_\Theta$ can occur in two ways: 1) the derivatives $\left.\frac{d}{dt}\phi(P_{t, \theta}) \right\vert_{t=0}$ exist for every $\theta \in \Theta$, but are not given by a bounded linear functional of the scores $s_{\theta}$, for example because they are not a function of $s_\theta$, are a non-linear function of $s_\theta$, or are linear but unbounded in $s_\theta$, or 2) the derivative $\left.\frac{d}{dt}\phi(P_{t, \theta}) \right\vert_{t=0}$ does not exist for some $\theta \in \Theta$.

We formalize the above discussion in the following lemma, where in the lemma $\varsigma^2_{P_0} \coloneqq E_{P_0}[\{\phi^1_{P_0}(O)\}^2]$:
\begin{lemma}
    \label{lemma:failure-differentiability-irregular}
    Suppose that $\widehat{\phi}_n$ is an asymptotically linear estimator of $\phi$ at $P_0$ with influence function $\phi^1_{P_0} \in L^2_0(P_0)$, and that $P \mapsto \phi(P)$ is not pathwise differentiable at $P_0$ relative to $\mathscr{P}_{\Theta}$. Then, there exist a submodel $\{P_{t, \theta}: t\in [0, \delta_\theta)\} \in \mathscr{P}_{\Theta}$ and sequences $n_k \rightarrow \infty$ and $h_k \rightarrow 1$ such that, under sampling from $P_{h_k/\sqrt{n_k}, \theta}$, either
    \begin{align*}
        \sqrt{n_k}(\widehat{\phi}_{n_k} - \phi(P_{h_k/\sqrt{n_k}, \theta})) {\rightsquigarrow} N\left(\Delta_{\theta}, \varsigma^2_{P_0}\right)
    \end{align*}
    for some $\Delta_\theta \in \mathbb{R}\setminus\{0\}$, or $\sqrt{n_k}(\widehat{\phi}_{n_k} - \phi(P_{h_k/\sqrt{n_k}, \theta}))$ diverges in probability. Consequently, for each $\alpha \in (0, 1)$, letting $Q_\alpha$ denote the $1 - \alpha/2$ quantile of the standard normal distribution, the corresponding Wald-type confidence interval
    \begin{align*}
        C_n \coloneqq \left(\widehat{\phi}_{n} - Q_{\alpha}\varsigma_{P_0}/\sqrt{n},  \widehat{\phi}_{n} + Q_{\alpha}\varsigma_{P_0}/\sqrt{n}\right)
    \end{align*}
    undercovers asymptotically along this sequence of local alternatives:
    \begin{align*}
        \lim_{k\rightarrow\infty} P_{h_k/\sqrt{n_k}, \theta}^{n_k}\left(\phi(P_{h_k/\sqrt{n_k}, \theta}) \in C_{n_k}\right) = \Phi\left(Q_\alpha - \frac{\Delta_\theta}{\varsigma_{P_0}}\right) - \Phi\left(-Q_\alpha - \frac{\Delta_\theta}{\varsigma_{P_0}}\right) < 1 - \alpha,
    \end{align*}
    where $\Phi$ is the standard normal distribution function, with the convention that the limit is $0$ when $\sqrt{n_k}(\widehat{\phi}_{n_k} - \phi(P_{h_k/\sqrt{n_k}, \theta}))$ diverges in probability or $\varsigma_{P_0} = 0$.
\end{lemma}
The above lemma shows that for any asymptotically linear estimator $\widehat{\phi}_n$, there exists a QMD at $P_0$ submodel that defines a local alternative data generating process for which $\widehat{\phi}_n$ exhibits non-zero $\sqrt{n}$-scaled bias, and for which Wald-type confidence intervals undercover asymptotically. The necessity of the subsequence $n_k$ is to handle the case where $P \mapsto \phi(P)$ does not even admit directional derivatives along the submodels in $\mathscr{P}_{\Theta}$. For related results, see \cite{vaart_differentiable_1991} and \cite{hirano_impossibility_2012}, who study the implications of non-differentiability for regular and locally unbiased estimation, and \cite{kaji_theory_2021}, who develops an efficiency theory for weakly identified parameters, which are not pathwise differentiable.

\section{The Baseline Statistical Model}

\label{sec:stat-models}

{
We now define $\mc{P}$, the largest statistical model we consider over which we desire to perform statistical inference, of which all subsequent models are submodels. Our choice of $\mc{P}$ suffices for the existence of a Neyman orthogonal estimating equation of the form \eqref{eq:neyman-orthogonal-ee}. However, this model does not encode any additional restrictions that may be implied by certain structural causal models, as in the case of NPIV and the shadow variable model. We return to this point in \cref{subsec:recover-continuity} below.

As we will show in \cref{subsec:failure-of-continuity}, the mapping $P \mapsto \varphi(P)$ is Hellinger discontinuous over $\mc{P}$ at distributions $P_0$ satisfying mild conditions. Combined with \cref{subsubsec:implications-discontinuity-main}, such discontinuity implies that local uniform inference is impossible over $\mc{P}$ at such $P_0$. In \cref{subsec:tangent-spaces} we derive the tangent space within this model, and discuss the relationship between the present results and the related works \cite{severiniEfficiencyBoundsEstimating2012, chen_overidentification_2018-1}.

}
Our first assumption concerns the smoothness of the functionals in \eqref{eq:mtilde-m-def}.
\begin{assumption}
    \label{as:continuous-linear-maps}
    At $P$, the functionals \eqref{eq:mtilde-m-def} are bounded and linear.
\end{assumption}
By the Riesz representation theorem, \cref{as:continuous-linear-maps} suffices for the existence of the Riesz representers $\alpha_P$ and $r_{P}$. This assumption is satisfied under mild conditions on $P$ for each of \cref{example:npiv}, \cref{example:shadow}, and \cref{example:proximal}, as can be seen from application of Example 5 in \cite{luedtke_one-step_2024}.

The next assumption concerns the subspaces to which these Riesz representers belong.
\begin{assumption}
    \label{as:solutions-exist}
    At $P$, $r_{P} \in \mathscr{R}(\mc{E}_{P})$ and $\alpha_{P} \in \mathscr{R}(\mc{E}_{P}^*)$, equivalently there exist solutions to \eqref{eq:primal-dual-inverse-problem}.
\end{assumption}
\cite{bennett_inference_2023} show that \cref{as:solutions-exist} suffices for identification of $\varphi(P)$ even when there exist multiple solutions to \eqref{eq:primal-dual-inverse-problem}, and for the validity of the doubly robust expansion of the estimating equation \eqref{eq:neyman-orthogonal-ee} from which DML estimators are constructed.  \cref{as:continuous-linear-maps,as:solutions-exist} are the minimal assumptions required to define the parameters considered in \cref{subsec:target-of-inference} via the dual formulation we give. Nevertheless, they should not be taken as a given. When $X, Z$ are continuously valued, \cref{as:solutions-exist} places intricate assumptions on the structure of $r_P$, $\alpha_P$, and the joint distribution of $(X, Z, D)$: Informally, $r_P$ and $\alpha_P$ must not be too strongly supported on the tailing singular vectors of the singular value decomposition of the operator $\mc{E}_P$ for each $P$.

For inference based on \eqref{eq:neyman-orthogonal-ee}, we require that $h_P$ and $q_P$ can be estimated sufficiently fast as the sample size grows. A typical regularity condition to facilitate fast estimation rates is the so-called $\beta$-source conditions \parencite{carrasco_chapter_2007, bennettSourceConditionDouble2023b, olivas-martinez_source-condition_2025}:
\begin{assumption}
    \label{as:beta-source}
    Let $\beta_h, \beta_q \in [0, \infty)$. At $P$,
    \begin{align}
        \label{eq:beta-source-representers}
        r_{P} \in \mathscr{R}(\{\mc{E}_{P}\mc{E}_{P}^*\}^{(\beta_h+1)/2}),  && \alpha_{P} \in \mathscr{R}(\{\mc{E}_{P}^*\mc{E}_{P}\}^{(\beta_q+1)/2}).
    \end{align}
\end{assumption}
\begin{remark}
    For $\beta_h, \beta_q \in (0, \infty)$, \eqref{eq:beta-source-representers} holds if and only if there exist solutions $h_{P}$, $q_{P}$ to the primal and dual inverse problems \eqref{eq:primal-dual-inverse-problem} satisfying
    \begin{align}
        \label{eq:beta-source-solutions}
        h_{P} \in \mathscr{R}(\{\mc{E}_{P}^*\mc{E}_{P}\}^{\beta_h/2}),  && q_{P} \in \mathscr{R}(\{\mc{E}_{P}\mc{E}_{P}^*\}^{\beta_q/2}).
    \end{align}
\end{remark}
For $\beta_h = \beta_q = 0$, \cref{as:beta-source} is equivalent to \cref{as:solutions-exist}. When $X, Z$ are finite discretely valued, these assumptions are equivalent for any choice of $\beta_h, \beta_q$. \cite{bennett_inference_2023} refer to \cref{as:beta-source} with $\beta_h = 0$, $\beta_q = 1$ as Strong Identification. Higher values of $\beta_h$, $\beta_q$ imply stronger smoothness of the solutions, as they ensure that $h_{P}$ and $q_P$ concentrate less mass on the tailing singular vectors of $\mc{E}_{P}$. Higher values of $\beta_h, \beta_q$ allow more flexible function classes to be used to estimate $h_P$, $q_P$, as measured by the critical radii of these classes. \cite{bennettSourceConditionDouble2023b} and \cite{olivas-martinez_source-condition_2025} analyze the behavior of DML estimators based on \eqref{eq:neyman-orthogonal-ee} under various values of $\beta_h$, $\beta_q$.

We also require that every distribution in $\mc{P}$ be mutually absolutely continuous relative to a reference distribution $P_0$, whose role is made precise following \eqref{eq:model-definition} below.
\begin{assumption}
    \label{as:mutual-absolute-continuity}
    $P$ is mutually absolutely continuous relative to the reference distribution $P_0$.
\end{assumption}
In addition to mutual absolute continuity, we will require that the support of $O$ be a bounded rectangle, and that for each $P$, the density of $P$ with respect to the appropriate product of Lebesgue and counting measures is essentially bounded above and below by a constant that does not depend on $P$.
\begin{assumption}
    \label{as:dominating-measure-support}
    At $P$,
    \begin{enumerate}
        \item The joint distribution of $(W, X, Z, D)$ is dominated by the product measure {$ \lambda = \lambda_W\otimes\lambda_X\otimes\lambda_Z\otimes\lambda_D$. $\lambda_W$ is equal to (possibly a normalized version of) the product measure $\bigotimes_{j=1}^{\ell_W}\mu_{j}$ where $\mu_j$ is either the Lebesgue or counting measure, and similarly for $\lambda_X, \lambda_Z$. $\lambda_D$ is the (possibly normalized) counting measure.}  In addition, the joint density $p$ satisfies that $L_{bounded}^{-1} < \essinf p \leq \esssup p < L_{bounded}$ for $L_{bounded}$ a universal constant not depending on $P$ where the essential infimum and supremum are with respect to $\lambda$ over the support of $P$.
        \item  $(W, X, Z, D)$ are supported on a bounded rectangle in $\prod_{j=1}^{\ell_W + \ell_X + \ell_Z} \mathbb{A}_j \times \prod_{j=1}^{\ell_D}\mathbb{B}_j$, where $\mathbb{A}_j$ is either $\mathbb{R}$ or $\mathbb{N}$ for each $j \in \{1, \dots, \ell_W + \ell_X + \ell_Z\}$ and $\mathbb{B}_j = \mathbb{N}$ for each $j \in \{1, \dots, \ell_D\}$.
    \end{enumerate}
\end{assumption}
These conditions, among other things, ensure the forthcoming negative results do not arise from any pathologies in the support of $P$ or due to extreme values of the density of $P$. For example, without uniform upper and lower bounds on the density of $O$ under $P$, there may exist functions $g \in L^2(P_0)$ such that even the mean functional $P \mapsto E_P[g(O)]$ fails to be Hellinger continuous \parencite{bahadur_nonexistence_1956}.

\cref{as:dominating-measure-support} allows for discrete shared covariates $D$ between the primal and dual inverse problems in \eqref{eq:primal-dual-inverse-problem}. { However, because the joint distribution of $(X, Z)$ under each $P$ is absolutely continuous with respect to a product measure, supported on a rectangle, and with density bounded below by a constant, shared continuous covariates are not allowed. In other words, we do not allow $X = (X', R), Z = (Z', R)$ with $R$ a random variable. Because $D$ is only able to accommodate discrete covariates, exclusion of such $R$ prevents inclusion of continuous covariates in \cref{example:npiv}, \cref{example:proximal}, and \cref{example:shadow}}. We leave an extension of the forthcoming results to settings with shared continuous covariates to future work.

We are now able to define the statistical model we consider inference over.
\begin{align}
    \label{eq:model-definition}
    \mc{P} \coloneqq \left\{P \text{ a distribution for }O \text{ such that Assumptions \ref{as:continuous-linear-maps}, \ref{as:solutions-exist}, \ref{as:beta-source}, \ref{as:mutual-absolute-continuity}, and \ref{as:dominating-measure-support} hold} \right\}
\end{align}
In the above definition, \cref{as:beta-source} is assumed to hold for a fixed but arbitrary choice of $\beta_h, \beta_q \in [0, \infty)$, and \cref{as:mutual-absolute-continuity} relative to the reference distribution $P_0$. As the forthcoming results hold for any choice of $\beta_h, \beta_q \in [0, \infty)$ and $P_0 \in \mc{P}$, we suppress this dependence in the notation.

Of these assumptions, \cref{as:continuous-linear-maps}, \cref{as:solutions-exist}, and \cref{as:beta-source} are standard in the literature on inference for these parameters, while \cref{as:mutual-absolute-continuity} and \cref{as:dominating-measure-support} are regularity conditions of our own. In restricting the laws available to perturb toward $P_0$, they only strengthen the impossibility results that follow: Hellinger discontinuity of $\varphi$ at $P_0$ over $\mc{P}$ implies discontinuity at $P_0$ over any model containing $\mc{P}$, and the maximal tangent space at $P_0$ can only be larger over a larger model.

In the remainder of this section, we give the main results of the paper. Under the model $\mc{P}$, the parameter $\varphi(P)$ is not Hellinger continuous at $P_0$ satisfying mild regularity conditions. We will further argue that the maximal tangent space of the model $\mc{P}$ at $P_0$ is equal to $L^2_0(P_0)$ under mild regularity conditions, implying these assumptions impose no first order restriction on the model. { We then consider continuity and pathwise differentiability under various restrictions of the model $\mc{P}$.}

\subsection{Failure of Continuity over the Baseline Model}
\label{subsec:failure-of-continuity}
We now demonstrate that the parameter $\varphi(P)$ is not Hellinger continuous over $\mc{P}$ defined in \eqref{eq:model-definition} at $P_0$ satisfying mild regularity conditions. Recall the definition of Hellinger continuity in \cref{def:hellinger-continuous}.
In fact, we will prove a stronger result than lack of continuity at $P$: that there exists a sequence $\{P_j\}\subseteq \mc{P}$ such that $\lim_{j\rightarrow \infty}H^2(P_j, P) = 0$ but where $\lim_{j\rightarrow \infty}|\varphi(P_j) - \varphi(P)| = \infty$.

We split our results for this subsection into two pieces. First, we consider the case where $X, Z$ are finite-discretely valued, then consider the case where $X, Z$ are continuously valued. For $X, Z$ discrete, we have the following result.
\begin{lemma}
    \label{lemma:discrete-necessary-sufficient-continuity}
    Suppose that $X, Z$ are supported on $\{1, \dots, K_X\}$, $\{1, \dots, K_Z\}$, respectively at $P_0$. Suppose $\mc{P}$ is as in \eqref{eq:model-definition} with $P_0 \in \mc{P}$. We have the following:
    \begin{enumerate}
        \item \label{lemma:discrete-necessary-sufficient-continuity:discontinuity}If there exists $d_e \in \textsf{Supp}[D;P_0]$ such that both $r_{P_0}(\cdot, {d}_e)$ and $\alpha_{P_0}(\cdot, {d}_e)$ are not constant functions and both $\mathscr{N}(\mc{E}_{P_0;d_e}) \not= \{0\}$ and $\mathscr{N}(\mc{E}_{P_0;d_e}^*) \not= \{0\}$, then
        $\varphi$ is not Hellinger continuous at $P_0$ over the model $\mc{P}$. Furthermore, there exists a sequence $\{P_j\}_{j=1}^\infty\subseteq\mc{P}$ such that $\lim_{j \rightarrow \infty}H^2(P_j, P_0) = 0$ but where $\lim_{j\rightarrow \infty}|\varphi(P_j) - \varphi(P_0)| = \infty$.
        \item \label{lemma:discrete-necessary-sufficient-continuity:continuity} If $W$ is supported on $\{1, \dots, K_W\}$ and for each $d \in \textsf{Supp}[D;P_0]$, either $\mathscr{N}(\mc{E}_{P_0;d}) = \{0\}$ or $\mathscr{N}(\mc{E}_{P_0;d}^*) = \{0\}$, $\varphi$ is Hellinger continuous at $P_0$ over the model $\mc{P}$.
    \end{enumerate}
\end{lemma}
The condition that $\mathscr{N}(\mc{E}_{P_0;d_e}) \not= \{0\}$ and $\mathscr{N}(\mc{E}_{P_0;d_e}^*) \not= \{0\}$ is equivalent to the statement that both the primal and dual inverse problems admit multiple solutions, i.e. that they are both ill-posed. When one of $\mathscr{N}(\mc{E}_{P_0;d_e})$ or $\mathscr{N}(\mc{E}^*_{P_0;d_e})$ is trivial, $X, Z$ finite discrete ensures the corresponding inverse problem is well-posed.

We now give some intuition for the above lemma. We begin with \cref{lemma:discrete-necessary-sufficient-continuity:continuity}. When either $\mc{E}_{P_0}$ or its adjoint has a trivial null space ($\mathscr{N}(\mc{E}_{P_0}) = \{0\}$ or $\mathscr{N}(\mc{E}^*_{P_0}) = \{0\}$), the operators $\mc{E}_{P_0}$ and $\mc{E}_{P_0}^*$ are full rank (i.e. $\text{rank}(\mc{E}_{P_0;d}) = \min(K_X, K_Z)$ for each $d$). In such cases, the Moore-Penrose pseudo-inverses of these operators are continuous at $P_0$. By this we mean that for every sequence of laws $\{P_j\}_{j =1 }^\infty\subseteq \mc{P}$ with $\lim_{j\rightarrow \infty}H^2(P_j, P_0) = 0$, the Moore-Penrose pseudo-inverses $\mc{E}_{P_j;d}^\dagger, \{\mc{E}_{P_j;d}^*\}^\dagger$ of the operators $\mc{E}_{P_j;d}, \mc{E}_{P_j;d}^*$ converge to those pseudo-inverses at $P_0$ (with respect to any norm equivalent with the Frobenius norm). Furthermore, under the assumptions of the Lemma, and in particular because $W$ is finite discrete, the Riesz representers $\alpha_{P_j}, r_{P_j}$ converge to the Riesz representers at $P_0$.  Consequently, the minimum norm solutions $h_{P_j}^\dagger \coloneqq \mc{E}_{P_j;d}^\dagger r_{P_j}$, $q_{P_j}^\dagger \coloneqq \{\mc{E}_{P_j;d}^*\}^\dagger \alpha_{P_j}$ converge to the minimum norm solutions at $P_0$, and so $\lim_{j\rightarrow \infty}\varphi(P_j) = \varphi(P_0)$.

Convergence of the Moore-Penrose pseudo-inverses holds if and only if the operators $\mc{E}_{P}, \mc{E}_{P}^*$ have constant rank locally at $P_0$ over the model $\mc{P}$. In finite dimensions, full rank is preserved under small perturbations. However, when $\mc{E}_{P}, \mc{E}_{P}^*$ are not full rank, i.e. they both have non-trivial null spaces, an arbitrarily small perturbation to $P_0$ can change the rank of the operators. Under such perturbations, the Moore-Penrose pseudo-inverse is discontinuous at $P_0$, and the minimum norm solutions to the inverse problems \eqref{eq:primal-dual-inverse-problem} can diverge to infinity, which can in turn cause $\varphi(P)$ to diverge. The sequence $\{P_j\}_{j=1}^\infty$ we construct along which $|\varphi(P_j)|$ diverges to infinity is one such sequence where the rank of $\mc{E}_{P_j;d_e}$ is greater than the rank of $\mc{E}_{P_0;d_e}$ for every $j$. Informally, the rank of the operator $\mc{E}_{P_j;d_e}$ encodes how much information $X$ and $Z$ carry about each other. When the rank decreases, $X$ and $Z$ provide an order of magnitude less information about each other. Hence, discontinuity of the parameter occurs when information relating $X$ and $Z$ can be lost at $P_0$ relative to laws near $P_0$. In \cref{sec:weak-instruments}, we relate this phenomenon to the very weak instrument setting in linear instrumental variables models.

\cref{lemma:discrete-necessary-sufficient-continuity:discontinuity} of the above lemma further requires that $r_{P_0}(\cdot, d_e)$ and $\alpha_{P_0}(\cdot, d_e)$ are not constant functions for some $d_e$ for which both inverse problems are ill-posed. Without this assumption, the parameter at $P_0$ degenerates to the $P_0$-mean of a known function.

We note that when $W, X, Z$ are discrete, the laws $P$ in $\mc{P}$ at which $\mc{E}_{P;d}$ fails to have full rank for some $d$ form a set of measure zero. To make this precise, identify $\mc{E}_{P;d}$ with the $K_Z \times K_X$ row-stochastic matrix whose $(z, x)$ entry is $P(X = x|Z = z, D = d)$, so that each row lies in the simplex of probability vectors on $\{1, \dots, K_X\}$. With respect to the uniform measure on the resulting product of simplices, the rank deficient matrices form a set of measure zero, since failure of full rank requires the vanishing of a polynomial in the entries that does not vanish identically.
On the other hand, when $X$ and $Z$ are continuous, failure of continuity is the typical case, as demonstrated in \cref{thrm:continuous-discontinuity} below. Before giving the Theorem statement, we introduce a further regularity condition for the functions $m, \widetilde{m}$ to satisfy at $P_0$.

\begin{assumption}
    \label{as:m-m-tilde-bounded}
    Suppose that $\widetilde{m}, m, P_0$ satisfy that there exists $0 < L < \infty$ such that for all $h \in L^1(P_{X,D;0})$, $q \in L^1(P_{Z,D;0})$,
    \begin{align*}
        \|\widetilde{m}(W, X, D;h)\|_{L^1(P_{W,X,D;0})} \leq L \|h\|_{L^1(P_{X,D;0})}, && \|m(W, Z, D;q)\|_{L^1(P_{W,Z,D;0})} \leq L \|q\|_{L^1(P_{Z,D;0})}.
    \end{align*}
\end{assumption}

\cref{as:m-m-tilde-bounded} is a technical condition used in the proofs of the results below. In particular, for any $P$ mutually absolutely continuous with $P_0$ and with density bounded above and below by a constant, it implies that $P$ satisfies \cref{as:continuous-linear-maps} and that the Riesz representers $\alpha_P$ and $r_P$ are essentially bounded (\cref{lemma:l1-bounded-functional}). It further implies that if such a sequence $\{P_j\}_{j=1}^\infty$ converges to $P_0$ in Hellinger distance, with densities bounded above and below by a constant not depending on $j$, then $\alpha_{P_j}$ and $r_{P_j}$ converge weakly to $\alpha_{P_0}$ and $r_{P_0}$ (\cref{lemma:weak-convergence-of-representers}).

Under mild conditions, this condition is satisfied for NPIV, for proximal causal inference, and for the shadow variable model.
\begin{examplecont}{example:npiv}
    Suppose that  $W, X, Z$ are 1-dimensional continuous random variables jointly supported on $[0, 1]^{3}$ at $P_0$. Suppose that at $P_0$, there exists $L_{P_0} \geq 1$ such that the density $p_0$ with respect to the Lebesgue measure satisfies $L_{P_0}^{-1} < p_0 < L_{P_0}$.
    Let $x \mapsto \psi(x)$ be a known uniformly bounded function on $[0, 1]$.
    Then, $\widetilde{m}, m, P_0$ satisfy \cref{as:m-m-tilde-bounded}:
    \begin{align*}
        \|\widetilde{m}(W, X;h) \|_{L^1(P_{W,X;0})} =&  \|\psi h\|_{L^1(P_{W,X;0})}\leq \|\psi\|_{L^\infty(P_{X;0})}\|h\|_{L^1(P_{X;0})}\\
         \|{m}(W, Z;q) \|_{L^1(P_{W,Z;0})} =&  \|W q\|_{L^1(P_{W,Z;0})}\leq \|W\|_{L^\infty(P_{W;0})}\|q\|_{L^1(P_{Z;0})}.
    \end{align*}
\end{examplecont}

\begin{examplecont}{example:proximal}
    Suppose that  $W, X, Z$ are 1-dimensional continuous random variables jointly supported on $[0, 1]^3$ at $P_0$. Recall $D$ is binary.  Suppose that at $P_0$, there exists $L_{P_0} \geq 1$ such that  $L_{P_0}^{-1} <  p_{0}(O)  < L_{P_0}$ $P_0$-a.e.
    Then, $\widetilde{m}, m, P_0$ satisfy \cref{as:m-m-tilde-bounded}:
    \begin{align*}
        \|{m}(W, Z, D;q) \|_{L^1(P_{W,Z,D;0})} =&  \|W q\|_{L^1(P_{W,Z,D;0})}\leq \|W\|_{L^\infty(P_{W;0})}\|q\|_{L^1(P_{Z;0})}\\
        \|\widetilde{m}(W, X;h) \|_{L^1(P_{W,X;0})} =& E_{P_0}[| h(X, 1) - h(X, 0)|] \\
        \leq & \sum_{d =0}^1 E_{P_0}\left[\frac{I(D=d)}{p_0(D=d|X)}|h(X, D)|\right]\\
        \leq & 2 L_{P_0}^2\|h\|_{L^1(P_{X,D;0})}.
    \end{align*}
\end{examplecont}

\begin{examplecont}{example:shadow}
    Suppose that $X, Z$ are $1$-dimensional continuous random variables jointly supported on $[0, 1]^2$ at $P_0$. Recall that $D$ is binary. Suppose that at $P_0$, there exists $L_{P_0} \geq 1$ such that $L_{P_0}^{-1} < p_0(O) < L_{P_0}$ $P_0$-a.e. Then, $\widetilde{m}, m,  P_0$ satisfy \cref{as:m-m-tilde-bounded}:
    \begin{align*}
        \|\widetilde{m}(X, D;h)\|_{L^1(P_{X,D;0})} =& E_{P_0}[(1-D)|h(X,1)|]\\
        \leq& E_{P_0}[|h(X, 1)|]\\
        \leq& E_{P_0}\left[\frac{D}{p_0(D|X)}|h(X, D)|\right]\\
        \leq& 2L_{P_0}^{2} \|h\|_{L^1(P_{X,D;0})}\\
        \|m(Z,D;q)\|_{L^1(P_{Z,D;0})}
    =&
    E_{P_0}[|DZq(Z,D)|]\\
    \leq& E_{P_0}[|q(Z,D)|]=\|q\|_{L^1(P_{Z,D;0})}.
    \end{align*}
\end{examplecont}

We are now ready to state our discontinuity theorem in the continuous case.
\begin{theorem}
    \label{thrm:continuous-discontinuity}
    Suppose $\mc{P}$ is as in \eqref{eq:model-definition} and $P_0 \in \mc{P}$. Suppose $m, \widetilde{m}$, $P_0$ satisfy \cref{as:m-m-tilde-bounded}. Suppose the joint distribution of $(X, Z)|D=d$ under $P_0$ is dominated by the Lebesgue product measure for each $d \in \textsf{Supp}(D;P_0)$. Suppose $r_{P_0}(\cdot, d)$ and $\alpha_{P_0}(\cdot, d)$ are not constant functions $P_0$-a.e for some $d \in \textsf{Supp}(D;P_0)$. Then, $\varphi$ is not Hellinger continuous at $P_0$ over the model $\mc{P}$, and there exists a sequence of distributions $\{P_j\}_{j=1}^\infty \subseteq \mc{P}$ such that $\lim_{j\rightarrow \infty}H^2(P_j, P_0) = 0$ but where $\lim_{j\rightarrow \infty}|\varphi(P_j) - \varphi(P_0)| =\infty$.
\end{theorem}
{ Under \cref{as:dominating-measure-support}, when $(X, Z)|D=d$ are continuous random variables for each $d$, the operators $\mc{E}_{P_0}$ and $\mc{E}_{P_0}^*$ are both compact. As such, the inverse problems \eqref{eq:primal-dual-inverse-problem} are ill-posed: when $\mc{E}_{P_0}$ has infinite rank, holding the operators fixed, arbitrarily small perturbations of the right hand sides $r_{P_0}, \alpha_{P_0}$ can produce arbitrarily large changes in the solutions $h_{P_0}$ and $q_{P_0}$.

Failure of continuity of $\varphi$ arises instead from discontinuity of the operator valued mappings $P \mapsto \mc{E}_{P}^\dagger$, $P \mapsto \{\mc{E}_{P}^*\}^\dagger$. Even though the operator valued mappings $P \mapsto \mc{E}_{P}$, $P \mapsto  \mc{E}_{P}^*$ are continuous at $P_0$, when the inverse problems are ill-posed, the mapping taking each law to the corresponding pseudo-inverses $\mc{E}_{P}^\dagger, \{\mc{E}_{P}^*\}^\dagger$ may be discontinuous at $P_0$. Such discontinuity can lead to discontinuity of the parameter even when the Riesz representers $r_P, \alpha_P$ do not depend on $P$, as with $\alpha_P = \psi$ in \cref{example:npiv} and $r_P(Z, D) = DZ$ in \cref{example:shadow}.

The discrete and continuous settings differ in how typical the discontinuity is. When $X, Z$ are finite discrete, the maps $P \mapsto \mc{E}_{P}^\dagger, P \mapsto \{\mc{E}_{P}^*\}^\dagger$ are discontinuous  at $P_0$ only when the rank of $\mc{E}_{P_0}$ can change under an arbitrarily small perturbation of $P_0$, and so arises only at the measure zero set of laws described above. When $X, Z$ are continuous and the joint density is bounded above and below by a constant, the maps $P \mapsto \mc{E}_{P}^\dagger, P \mapsto \{\mc{E}_{P}^*\}^\dagger$ are guaranteed to be discontinuous at $P_0$. No restriction imposed by the model $\mc{P}$ prevents this discontinuity, including the $\beta$-source conditions. In our construction of the sequence $\{P_j\}_{j=1}^\infty\subseteq \mc{P}$ along which $\varphi(P_j) \rightarrow \infty$, the operators $\mc{E}_{P_j}$, $\mc{E}^*_{P_j}$ are finite rank for each $j$, and hence satisfy \cref{as:beta-source} for every choice of $\beta_h, \beta_q \in [0, \infty)$. As in the discrete case, we require that $\alpha_{P_0}(\cdot, d)$ and $r_{P_0}(\cdot, d)$ are not constant functions so that $\varphi(P_0)$ does not degenerate to the mean of a fixed function at $P_0$.

}

{The proof of ~\cref{thrm:continuous-discontinuity} admits an interpretation in terms of weak instruments, which we develop in \cref{sec:weak-instruments}. Informally, the singular functions of $\mc{E}_P$ act as canonical treatment and instrument contrasts, and compactness of $\mc{E}_P$ produces infinitely many such contrasts that are arbitrarily weakly instrumented. The approximating laws in the proof place signal on a contrast whose strength tends to zero faster than the target-relevant mass placed on it, so that a Hellinger-small perturbation of the observed law, though hard to detect statistically, produces an unbounded change in the estimand. In this sense, \cref{thrm:continuous-discontinuity} can be read as a continuum analogue of a many-weak-instruments phenomenon.}

\cref{lemma:discrete-necessary-sufficient-continuity} and \cref{thrm:continuous-discontinuity} have severely negative implications for inference, as they imply that under the conditions of those results, the local modulus of continuity $\omega_{P_0}(r)$ of $\varphi$ over $\mc{P}$ is infinite for every $r > 0$. By \cref{lemma:finite-sample-discontinuity}, at such points $P_0$, locally uniform inference over $\mc{P}$ is impossible for each $n$: the expected diameter of every locally uniformly honest confidence interval is infinite, and for every $M > 0$, every estimator has error exceeding $M$ with probability at least $(1-\tau)/2$ at some law in a Hellinger neighborhood of $P_0$ of radius $r_{n, \tau}$. Asymptotically, no estimator is locally uniformly consistent, and every confidence interval with locally uniform asymptotic coverage has diameter that is not bounded in probability under $P_0$.

\subsection{The Model Tangent Spaces}
\label{subsec:tangent-spaces}

We conclude this section by deriving the maximal tangent space of the model $\mc{P}$ at selected laws $P_0$. As we will see, the tangent space is unrestricted when $L^2(P_{X, D;0})$ and $L^2(P_{Z, D;0})$ have the same dimension, meaning that the model $\mc{P}$ is unrestricted in a local sense at such $P_0$. Note that when $X$ and $Z$ are continuous, $\dim(L^2(P_{X, D;0})) = \dim(L^2(P_{Z, D;0}))$ automatically. This gives some further intuition on why the assumptions on the model $\mc{P}$ do not suffice for locally uniformly valid inference. The present results stand in contrast to related work that implicitly or explicitly derives the tangent space relative only to particularly smooth submodels of $\mc{P}$ \parencite{severiniEfficiencyBoundsEstimating2012, chen_overidentification_2018-1}.  We still consider it instructive to derive this tangent space even though $P \mapsto \varphi(P)$ is Hellinger discontinuous over this model, and hence not pathwise differentiable (\cref{remark:failure-of-continuity-fails-diff}), and so the tangent space cannot be used to gain insights into the behavior of regular estimators of $\varphi(P)$.

As in the discussion of continuity, we begin with the discrete case. We only give our first result for $K_X = K_Z$.
\begin{lemma}
    \label{lemma:discrete-tangent-space}
    Suppose $\mc{P}$ is as in \eqref{eq:model-definition} and $P_0 \in \mc{P}$. { Suppose that $X$ is supported on $\{1, \dots, K_X\}$ and $Z$ is supported on $\{1, \dots, K_Z\}$ with $K_X = K_Z$. Then, $\mc{T}(P_0, \mc{P}) = L^2_0(P_0)$.}
\end{lemma}
The above lemma reveals the model $\mc{P}$ is locally non-parametric at any $P_0$ in the model. In particular, \cref{as:solutions-exist}, the assumption that there exist solutions to the primal and dual inverse problems, does not locally restrict the model, even when there exist multiple solutions to both the primal and dual inverse problems at $P_0$. Intuitively, this result follows from identifying $\mc{E}_{P}$ and $\mc{E}_{P}^*$ with $K_X \times K_X$ matrices. Solutions to the primal and dual inverse problems exist whenever $\alpha_P, r_P$, represented as vectors, are in the range of these matrices. Such solutions trivially exist if these matrices are invertible. But, every $K_X \times K_X$ matrix is arbitrarily close in Frobenius norm to an invertible $K_X \times K_X$ matrix. Upon establishing a relation between Frobenius norm and Hellinger distance, it becomes clear that every law at which solutions do not exist is arbitrarily close in Hellinger distance to a law where solutions do exist. One can show, as we do in the proof of \cref{lemma:discrete-tangent-space}, that such result implies a locally non-parametric tangent space at every $P_0 \in \mc{P}$.

The above intuition persists in the case where $X$, $Z$ are both continuous. In such cases, the spaces $L^2(P_{Z, D})$ and $L^2(P_{X, D})$ are the same dimension. Again, every law $P$ at which there do not exist solutions to the primal and/or dual inverse problems is arbitrarily close in Hellinger distance to a law $P'$ for which solutions exist. The technical argument is much more involved in this case, because \cref{as:dominating-measure-support} precludes $\mc{E}_{P}$ and $\mc{E}_{P}^*$ from admitting a bounded inverse. Nevertheless, an appropriate adaptation can be made, from which we arrive at the following theorem.
\begin{theorem}
    \label{thrm:continuous-tangent-space}
    Suppose $\mc{P}$ is as in \eqref{eq:model-definition} and $P_0 \in \mc{P}$. Suppose $m$, $\widetilde{m}, P_0$ satisfy \cref{as:m-m-tilde-bounded}. Suppose the joint distribution of $(X, Z)|D=d$ under $P_0$ is dominated by the Lebesgue product measure for each $d \in \textsf{Supp}[P_{D;0}]$. Then, $\mc{T}(P_0, \mc{P}) = L^2_0(P_0)$.
\end{theorem}

\cref{lemma:discrete-tangent-space} and \cref{thrm:continuous-tangent-space} diverge from existing results computing tangent spaces within models similar to $\mc{P}$, for example \cite{severiniEfficiencyBoundsEstimating2012} and \cite{chen_overidentification_2018-1}. Such divergence arises because those authors only derived tangent spaces relative to the set of submodels along which the derivatives of particular solutions $h_P$ to the primal inverse problem exist. For concreteness, consider the following set of submodels:
\begin{align}
    \label{eq:differentiable-solution-paths}
    \mathscr{P}_{h_P^\dagger \text{ diff}}(P_0, \mc{P}) \coloneqq \{\{P_t:t\in [0, \delta)\}: \{P_t:t\in [0, \delta)\} \in \mathscr{P}(P_0, \mc{P}), \frac{d}{dt}h_{P_t}^\dagger\vert_{t=0} \text{ exists }\}
\end{align}
where $\mathscr{P}(P_0, \mc{P})$ is the set of all one-dimensional QMD at $P_0$ submodels of $\mc{P}$ and the derivative is taken in $L^2(P_{X, D;0})$. The proofs of \cite{severiniEfficiencyBoundsEstimating2012} only consider scores of submodels in a related set of submodels, but with different assumptions on $\mc{P}$.  \cite{chen_overidentification_2018-1} also consider differentiability relative to a related set of submodels, see Condition A in that work.

Restricting to such submodels, those authors show the tangent space is typically restricted when $\mathscr{N}(\mc{E}_{P_0}^*) \not = \{0\}$, i.e. the dual inverse problem admits multiple solutions. We say typically because the exact nature of the restrictions depend on $m, \widetilde{m}$, and $P_0$. The restrictions on the submodels in \eqref{eq:differentiable-solution-paths} are strong because, when the primal inverse problem is ill-posed, the map $t \mapsto h_{P_t}^\dagger$ may fail to even be continuous, let alone differentiable, even when $\{P_t: t\in [0, \delta)\} \subseteq \mc{P}$ is QMD at $P_0$. See for example the sequences $P_j$ along which $\varphi(P_j)$ diverges in the proofs of \cref{lemma:discrete-necessary-sufficient-continuity} and \cref{thrm:continuous-discontinuity}. However, such paths contribute to the maximal tangent space, as seen in the discrepancy. It is exactly those submodels along which $t \mapsto h_{P_t}^\dagger$ is not differentiable that give rise to the scores appearing in the maximal tangent space that do not appear as scores for submodels in $\mathscr{P}_{h_P^\dagger \text{ diff}}(P_0, \mc{P})$.

When $X, Z$ are finite discrete, there exists a simple sufficient condition under which every one dimensional QMD at $P_0$ submodel belongs to the set in \eqref{eq:differentiable-solution-paths}, namely that the rank of the operators $\mc{E}_{P;d}$ and $\mc{E}_{P;d}^*$ is constant for each $d$, and that the map $P \mapsto r_{P}$ is pathwise differentiable as a Hilbert-valued parameter. We call the submodel of $\mc{P}$ satisfying this constant rank condition $\mc{P}_{rank}$. In such cases, the map $t \mapsto \mc{E}_{P_t}^\dagger$ is itself differentiable in the appropriate metric at $t=0$, as is $t \mapsto r_{P_t}$ at $t=0$. Then, by application of the chain rule, $t \mapsto h_{P_t}^\dagger = \mc{E}_{P_t}^\dagger r_{P_t}$ is differentiable at $t=0$. Under this condition, the argument of \cite{severiniEfficiencyBoundsEstimating2012} applies to every QMD at $P_0$ submodel of $\mc{P}_{rank}$, giving the following lemma. We note that we require both $\alpha_{P}$ and $r_P$ to be pathwise differentiable, since one can then make a symmetric argument for the dual inverse problem. We require this differentiability relative to the model
\begin{align}
    \label{eq:riesz-model}
    \mc{P}_{\text{Riesz}} \coloneqq \{P: P \text{ satisfies \cref{as:continuous-linear-maps,as:mutual-absolute-continuity,as:dominating-measure-support}}\},
\end{align}
obtained by dropping \cref{as:solutions-exist} and \cref{as:beta-source} from the definition of $\mc{P}$, so that $\mc{P} \subseteq \mc{P}_{\text{Riesz}}$. The Riesz representers $\alpha_P$ and $r_P$ are defined at every $P \in \mc{P}_{\text{Riesz}}$. If $m, \widetilde{m}, P_0$ satisfy \cref{as:m-m-tilde-bounded}, then $\mc{P}_{\text{Riesz}}$ contains the linear tilts with densities $(1 + ts)p_0$ for each bounded $s \in L^2_0(P_0)$ and $t$ small enough, so that the local parameters $\dot\alpha_{P_0}$ and $\dot r_{P_0}$ are defined on all of $L^2_0(P_0)$.
\begin{lemma}
    \label{lemma:constant-rank-tangent-space}
    Suppose $\mc{P}$ is as in \eqref{eq:model-definition} and $P_0 \in \mc{P}$. { Suppose that $X$ is supported on $\{1, \dots, K_X\}$ and $Z$ is supported on $\{1, \dots, K_Z\}$. }
    \begin{enumerate}
        \item \label{lemma:discrete-tangent-space:part:constant-rank} Define
            \begin{align*}
                \mc{P}_{\text{rank}} = \{P: P \in \mc{P},  \text{rank}(\mc{E}_{P;d}) = \text{rank}(\mc{E}_{P_0;d}) \text{ for each $d \in \textsf{Supp}[P_{D;0}]$}\},
            \end{align*}
        Suppose $m, \widetilde{m}, P_0$ satisfy \cref{as:m-m-tilde-bounded}, and that $P \mapsto \alpha_{P}$ and $P \mapsto r_{P}$ are pathwise differentiable at $P_0$ relative to the set of all QMD at $P_0$ submodels of $\mc{P}_{\text{Riesz}}$, and let $\dot \alpha_{P_0}: L^2_0(P_{0}) \rightarrow L^2(P_{X,D;0})$, $\dot r_{P_0}: L^2_0(P_{0}) \rightarrow L^2(P_{Z,D;0})$ denote their local parameters.

        Define
        \begin{align*}
            \mc{C}_{rank}(P_0) \coloneqq \Big\{ &s \in L^2_0(P_0): \text{For each $h_{P_0}, q_{P_0}$ solving } \mc{E}_{P_0}h = r_{P_0}, \mc{E}_{P_0}^*q = \alpha_{P_0}, \\
                \nonumber &\dot\alpha_{P_0}[s](X, D) - E_{P_0}[q_{P_0}(Z, D)s_{Z|X, D}(X, Z, D)|X, D]  \in \mathscr{N}(\mc{E}_{P_0})^\perp \text{ and } \\
                &\dot r_{P_0}[s](Z, D) - E_{P_0}[h_{P_0}(X, D)s_{X|Z, D}(X, Z, D)|Z, D]  \in \mathscr{N}(\mc{E}_{P_0}^*)^\perp\Big\}
        \end{align*}
        where
            \begin{align*}
                s_{Z|X, D}(X, Z, D) &\coloneqq E_{P_0}[s(O)|X, Z, D] - E_{P_0}[s(O)|X, D] \quad \text{and} \\
                s_{X|Z, D}(X, Z, D) &\coloneqq E_{P_0}[s(O)|X, Z, D] - E_{P_0}[s(O)|Z, D].
            \end{align*}

        Then,
            \begin{align}
                \label{eq:discrete-tangent-space-restriction}
                \mc{T}(P_0, \mc{P}_{\text{rank}}) \subseteq \overline{\mc{C}_{rank}(P_0)}.
            \end{align}

            If ${m}(W, Z, D;\eta) \in L^2(P_0)$ and $\widetilde{m}(W, X, D;\widetilde{\eta}) \in L^2(P_0)$ for each $\eta \in \mathscr{N}(\mc{E}_{P_0}^*)$, $\widetilde{\eta} \in \mathscr{N}(\mc{E}_{P_0})$, respectively, the right hand side of \eqref{eq:discrete-tangent-space-restriction} admits an alternative expression as
            \begin{align*}
                \overline{\mc{C}_{rank}(P_0)} = cl\Big(\Big\{ &s \in L^2_0(P_0): \text{For each $h_{P_0}, q_{P_0}$ solving } \mc{E}_{P_0}h = r_{P_0}, \mc{E}_{P_0}^*q = \alpha_{P_0}, \\
                \nonumber &\langle s, m(\cdot;\eta) - h_{P_0}\eta\rangle_{L^2(P_0)}=0 \text{ for all }  \eta\in \mathscr{N}(\mc{E}_{P_0}^*) \text{ and } \\
                \nonumber &\langle s, \widetilde{m}(\cdot;\widetilde{\eta}) - q_{P_0}\widetilde{\eta}\rangle_{L^2(P_0)}=0 \text{ for all }  \widetilde{\eta}\in \mathscr{N}(\mc{E}_{P_0})\Big\}\Big).
            \end{align*}
            Under these same conditions, if for some $d \in \textsf{Supp}[P_{D;0}]$ either $\mathscr{N}(\mc{E}_{P_0;d}^*)\not= \{0\}$ and $r_{P_0}(\cdot, d)$ is a non-constant function of $Z$, or $\mathscr{N}(\mc{E}_{P_0;d})\not= \{0\}$ and $\alpha_{P_0}(\cdot, d)$ is a non-constant function of $X$, then $\mc{T}(P_0, \mc{P}_{rank}) \subsetneq L^2_0(P_0)$.
        \item \label{lemma:discrete-tangent-space:part:full-rank} If $\mathscr{N}(\mc{E}_{P_0}) = \{0\}$ or $\mathscr{N}(\mc{E}_{P_0}^*) = \{0\}$, then $\mc{T}(P_0, \mc{P}) = \mc{T}(P_0, \mc{P}_{\text{rank}})$.
    \end{enumerate}
\end{lemma}

\begin{remark}
    \label{remark:simplified-restrictions}
    If $K_Z > K_X$ and $\mathscr{N}(\mc{E}_{P_0}) = \{0\}$, then $\mathscr{R}(\mc{E}_{P_0}^*)  = L^2(P_{X,D;0})$ and there is a unique solution $h_{P_0}$ to $\mc{E}_{P_0}h = r_{P_0}$. Furthermore, \cref{lemma:discrete-tangent-space:part:full-rank} applies, so that $\mc{T}(P_0, \mc{P}) = \mc{T}(P_0, \mc{P}_{\text{rank}})$.

    In this setting, the right hand side of \eqref{eq:discrete-tangent-space-restriction} simplifies significantly, leading to the inclusion
    \begin{align*}
        \mc{T}(P_0, \mc{P})\subseteq \Big\{s \in L^2_0(P_0):& \text{For $h_{P_0}$ the unique solution to } \mc{E}_{P_0}h = r_{P_0}, \\
                &\dot r_{P_0}[s](Z, D) - E_{P_0}[h_{P_0}(X, D)s_{X|Z, D}(X, Z, D)|Z, D]  \in \mathscr{R}(\mc{E}_{P_0})\Big\}.
    \end{align*}
    A similar simplification holds if $K_X > K_Z$ and $\mathscr{N}(\mc{E}_{P_0}^*) = \{0\}$.
\end{remark}

In the above lemma, we additionally require that $P \mapsto \alpha_P$, $P \mapsto r_P$ be pathwise differentiable in a non-parametric model containing $\mc{P}_{rank}$ so that their local parameters are defined with domain $L^2_0(P_0)$. \cref{lemma:discrete-tangent-space:part:full-rank} of the above lemma reveals that the constant rank model is locally equivalent to the model $\mc{P}$ when the operators $\mc{E}_{P_0}$ and $\mc{E}_{P_0}^*$ are full rank. This result follows from the stability of the rank of $\mc{E}_{P}$ and $\mc{E}_{P}^*$ locally around $P_0$ under such full rank condition. We again note that in the discrete case, the set of laws $P$ at which $\mc{E}_{P}$ and $\mc{E}_{P}^*$ are not full rank is in an appropriate sense measure zero, and so \cref{lemma:discrete-tangent-space:part:full-rank} of the above lemma is in some sense the typical case for the discrete problem.

For $X, Z$ continuous, such a constant rank assumption becomes tenuous, as the rank of the operators $\mc{E}_{P}$ and $\mc{E}_{P}^*$ is typically infinite. Instead, we derive in \cref{thrm:injective-tangent-space} a bound on the maximal tangent space in a model where the norm of $h_P^\dagger$ is uniformly bounded over the model, at $P_0$ where $\mc{E}_{P_0}$ is injective, as well as in a model where the norm of the $\beta_h$-source representer of $h_P^\dagger$ is uniformly bounded over the model, whether or not $\mc{E}_{P_0}$ is injective (see \cref{subsec:recover-differentiability}). This bound recovers the bound derived by \cite{severiniEfficiencyBoundsEstimating2012} and \cite{chen_overidentification_2018-1} relative only to submodels in \eqref{eq:differentiable-solution-paths}, even though in the model we consider, the maps $t \mapsto h_{P_t}^\dagger$ may not be differentiable.

{
\section{Models with Uniform Norm Bounds}
\label{sec:norm-bound-models}
We now consider properties of $\varphi$ relative to submodels of the model $\mc{P}$ that uniformly bound the minimum norm solutions to either the primal, the dual, or both the primal and dual inverse problems.

\subsection{Recovery of Continuity under Uniform Norm Bounds}
\label{subsec:recover-continuity}

As shown in \cref{lemma:discrete-necessary-sufficient-continuity}, a sufficient condition for Hellinger continuity of $P \mapsto \varphi(P)$ is that at $P_0$, $\mc{E}_{P_0;d}$ is full rank for each $d$. More generally, for $(X, Z)$ discrete or continuous, uniformly bounding the norm of the minimum norm solutions of the inverse problems \eqref{eq:primal-dual-inverse-problem}, together with further regularity conditions, suffices to recover continuity.

We first give our assumption bounding the norm of the minimum norm solutions to the inverse problems \eqref{eq:primal-dual-inverse-problem}. We consider two versions, one in which $h_P^\dagger$ is uniformly bounded and one in which $q_P^\dagger$ is uniformly bounded.
\begin{assumption}
    \label{as:uniform-norm-bound}
    Let $C > 0$ be a constant not depending on $P$. Consider the following two conditions at $P$:
    \begin{enumerate}
        \item \label{as:uniform-norm-bound:part:h}

        $\|h^\dagger_P(\cdot, d)\|_{L^2(P_{X|d})} < C$ for each $d \in \textsf{Supp}[D;P]$ with $h^\dagger_P(\cdot, d)$ the minimum norm solution to $\mc{E}_{P;d}h = r_P(\cdot, d)$.
        \item \label{as:uniform-norm-bound:part:q}

        $\|q^\dagger_P(\cdot, d)\|_{L^2(P_{Z|d})} < C$ for each $d \in \textsf{Supp}[D;P]$
    with $q^\dagger_P(\cdot, d)$ the minimum norm solution to $\mc{E}_{P;d}^*q = \alpha_P(\cdot, d)$.
    \end{enumerate}
\end{assumption}
The above conditions restrict each $d$ separately. Because $\mc{E}_P$ acts blockwise in $d$, the function $\sum_{d}I(D=d)h_P^\dagger(\cdot, d)$ is the minimum norm solution to the overall inverse problem $\mc{E}_{P}h = r_P$, and $\|h_P^\dagger\|^2_{L^2(P_{X,D})} = \sum_d P(D=d)\|h_P^\dagger(\cdot, d)\|^2_{L^2(P_{X|d})}$. Bounding the by-$d$ solutions uniformly is therefore equivalent, up to the factors $P(D=d)$, to bounding the solution to the overall inverse problem, and analogously for $q_P^\dagger$. We discuss the plausibility of this assumption after giving the main results of this subsection. We state many of the following results under the assumption that \cref{as:uniform-norm-bound:part:h} of \cref{as:uniform-norm-bound} holds, although by symmetry they could be stated similarly under \cref{as:uniform-norm-bound:part:q}. To this end, define the following model.
\begin{align}
    \label{eq:h_P-bounded-model}
    \mc{P}_{h_{P}^\dagger \text{ bound}} \coloneqq \{P \in \mc{P}: \cref{as:uniform-norm-bound} \cref{as:uniform-norm-bound:part:h} \text{ holds}\}.
\end{align}
We also define the smaller model
\begin{align}
    \label{eq:h_P-q_P-bounded-model}
    \mc{P}_{h_{P}^\dagger, q_{P}^\dagger \text{ bound}} \coloneqq \{P \in \mc{P}: \cref{as:uniform-norm-bound} \cref{as:uniform-norm-bound:part:h} \text{ and \cref{as:uniform-norm-bound:part:q} hold}\}.
\end{align}

}

To demonstrate continuity of $P \mapsto \varphi(P)$ at $P_0$ over any submodel $\mc{P}_{cont}$ of $\mc{P}_{h_{P}^\dagger \text{ bound}}$, we further require that the Riesz representer maps are Hellinger continuous at $P_0$ over $\mc{P}_{cont}$. We introduce the submodel $\mc{P}_{cont}$ as it may happen that the Riesz representers are not continuous at $P_0$ over the larger model $\mc{P}_{h_{P}^\dagger \text{ bound}}$. If they are continuous at $P_0$ over $\mc{P}_{h_{P}^\dagger \text{ bound}}$, then $\mc{P}_{cont}$ can be taken equal to $\mc{P}_{h_{P}^\dagger \text{ bound}}$.

As $r_P$ and $\alpha_P$ take values in $L^2(P_{Z, D;0})$ and $L^2(P_{X, D;0})$, respectively, Hellinger continuity of $P \mapsto r_P$ at $P_0$ over $\mc{P}_{cont}$ means that for all $\epsilon > 0$, there exists $\delta > 0$ such that for all $P' \in \mc{P}_{cont}$ such that $H(P', P_0) < \delta$, $\|r_{P'} - r_{P_0}\|_{L^2(P_{Z, D;0})} < \epsilon$. A similar definition holds for $P \mapsto \alpha_P$. In \cref{example:npiv}, \cref{example:proximal}, and \cref{example:shadow}, under mild conditions on $P_0$, the Riesz representer maps are Lipschitz with respect to the Hellinger distance over $\mc{P}_{\text{Riesz}}$, and hence Hellinger continuous at $P_0$ over any submodel of $\mc{P}_{\text{Riesz}}$, as we show in \cref{sup-sec:examples} using Example 5 of \cite{luedtke_one-step_2024}. In such cases, $\mc{P}_{cont}$ can be taken equal to $\mc{P}_{h_{P}^\dagger \text{ bound}}$.
We are now able to state the following theorem.
\begin{theorem}
    \label{thrm:bounded-means-continuous}
    Let $\mc{P}_{cont}$ be any submodel of $\mc{P}_{h_{P}^\dagger \text{ bound}}$, with $P_0 \in \mc{P}_{cont}$, and suppose that $P \mapsto r_P$ and $P \mapsto \alpha_P$ are Hellinger continuous at $P_0$ over $\mc{P}_{cont}$. Then, $P \mapsto \varphi(P)$ is Hellinger continuous at $P_0$ over $\mc{P}_{cont}$.
\end{theorem}

{

~\cref{thrm:bounded-means-continuous} implies that the impossibility results of \cref{subsec:failure-of-continuity} do not apply
over models on which the minimum-norm solution to at least one of the two
inverse problems is uniformly bounded at points $P_0$ where the Riesz representer maps are continuous.  \cref{as:uniform-norm-bound} should therefore be read
as a substantive stability condition on at least one of the bridge functions that appear in the definition of
\(\varphi(P)\).  It says that, uniformly over the model, the primal bridge (\cref{as:uniform-norm-bound:part:h}) or dual bridge (\cref{as:uniform-norm-bound:part:q}) can be chosen to have moderate \(L^2\) size.  Equivalently,
the integral equation is not allowed to be solved only by bridge functions
with very large positive and negative values that cancel after conditioning on
the proxy variables. In the language of \cref{sec:weak-instruments}, the bound on the primal bridge plays the role of a uniform bound on the structural coefficient in a linear instrumental variables model. We next discuss the plausibility of \cref{as:uniform-norm-bound} in \cref{example:npiv}, then in a further example involving a misspecified NPIV model, and finally in \cref{example:proximal} and \cref{example:shadow}.

\begin{examplecont}{example:npiv}  As pointed out in \cref{example:npiv}, under the additive structural NPIV model \eqref{eq:npiv-model}, the causal dose-response function $h_{P,cdr}$ solves the primal inverse problem, i.e. it satisfies \eqref{eq:primal-npiv-eq}. Thus, because  equality \eqref{eq:dose-response-curve} holds, when the absolute value of the counterfactual outcome $W(x)$ is uniformly bounded over $x$ and over the model by a constant $B$, we have that, under any $P$ that obeys the model, $\|h_{P;\mathrm{cdr}}\|_{L^\infty(P_X)}< B$. Since \(h_P^\dagger\) is the
minimum-norm solution, then
\(
    \|h_P^\dagger\|_{L^2(P_X)}
    \le
    \|h_{P;\mathrm{cdr}}\|_{L^2(P_X)}.
\)
Thus \cref{as:uniform-norm-bound} \cref{as:uniform-norm-bound:part:h} follows from a transparent substantive condition:
the causal dose-response curve is uniformly bounded.  In this setting,
\cref{as:uniform-norm-bound} is not merely a technical restriction; it is implied by the
structural interpretation of the bridge under the additive NPIV model and boundedness of the counterfactual outcomes.
\end{examplecont}

\begin{example}
[Misspecified NPIV model]
\label{example:wrongNPIV}
Suppose now that the structural model is allowed to include a direct effect of
the instrument:
\begin{equation*}
    W(x,z)
    =
    h_{P;\mathrm{cdr}}(x)
    +
    u_P(z)
    +
    \epsilon_{W;P},
    \qquad
    E_P(\epsilon_{W;P}\mid Z)=0,
\end{equation*}
where \(h_{P;\mathrm{cdr}}\) is unrestricted and \(u_P\) is normalized so that
\(E_P\{u_P(Z)\}=0\). The term \(u_P(z)\) represents a violation of the
exclusion restriction. Indeed, for any fixed \(x\), the contrast
\(W(x,z)-W(x,z')=u_P(z)-u_P(z')\) is the individual direct effect of changing
the instrument from \(z'\) to \(z\) while holding the treatment fixed at \(x\).
Because \(u_P\) depends on the counterfactual outcomes and is not identified
from the observed data without the exclusion restriction, throughout this
example we let \(P\) index the law of the full data, including the
counterfactual outcomes \(W(x,z)\).

Let \(W(x,\cdot):=W(x,Z)\) denote the counterfactual outcome when \(X\) is set
to \(x\) but \(Z\) is left at its natural value. Since \(E_P\{u_P(Z)\}=0\) and
\(E_P(\epsilon_{W;P})=0\), the structural model implies that
\(E_P\{W(x,\cdot)\}=h_{P;\mathrm{cdr}}(x)\). Thus \(x\mapsto
h_{P;\mathrm{cdr}}(x)\) retains the interpretation of a causal dose-response
curve, now averaging over the natural distribution of \(Z\). In particular, if
the counterfactual outcome is uniformly bounded over the model, then
\(h_{P;\mathrm{cdr}}\) is uniformly bounded as well.

Suppose further that the direct effect of the instrument is mild, in the sense
that \(|u_P(z)-u_P(z')|\le \delta\) for all \(z,z'\). Since \(u_P(Z)\) has mean
zero and takes values in an interval of length at most \(\delta\), it follows
that \(\|u_P\|_{L^2(P_Z)}\le \delta/2\). Thus this condition gives a clear sense
in which the violation of the exclusion restriction is small on the \(Z\)-scale.

The observed regression of \(W\) on \(Z\) satisfies
\(E_P(W\mid Z)=E_P\{h_{P;\mathrm{cdr}}(X)\mid Z\}+u_P(Z)\). Equivalently,
\(h_{P;\mathrm{cdr}}\) solves the correctly specified inverse problem with
right-hand side \(E_P(W\mid Z)-u_P(Z)\). Now suppose that the parameter of
interest is the weighted average
\[
\varphi(P)=E_P\{\psi(X)
h_{P;\mathrm{cdr}}(X)\},
\]
where \(\psi\) belongs to the range of
\(\mathcal E_P^*\). Suppose also that \(u_P=\mathcal E_P s_P\) for some
\(s_P\in L^2(P_X)\). Then \(h_{P;\mathrm{cdr}}+s_P\) solves the misspecified
inverse problem \((\mathcal E_P h)(Z)=E_P(W\mid Z)\), which therefore admits a
solution. If, incorrectly, we impose the exclusion restriction
\(u_P=0\), then we instead target
\[
\varphi_{\mathrm{bias}}(P):=E_P\{\psi(X)h_{P;\mathrm{bias}}(X)\},
\]
where
\(h_{P;\mathrm{bias}}\) is any solution of the misspecified inverse problem
\((\mathcal E_P h)(Z)=E_P(W\mid Z)\). The range condition on \(\psi\) ensures that
this functional is well-defined, in the sense that it does not depend on which
solution of the misspecified inverse problem is chosen.

A crucial point is that the solutions to the misspecified inverse problem may fail to be uniformly bounded in norm, and thus may fail to satisfy \cref{as:uniform-norm-bound}.
To see this, recall that boundedness of the counterfactual outcome implies that
\(h_{P;\mathrm{cdr}}\) is uniformly bounded over the model, and the mild-direct
effect condition implies that \(u_P\) is uniformly small in \(L^2(P_Z)\).
However, neither fact implies that the solutions of
\((\mathcal E_P h)(Z)=E_P(W\mid Z)\) have uniformly bounded \(L^2(P_X)\) norm. Specifically, recall that
\(h_{P;\mathrm{cdr}}+s_P\) solves the misspecified equation, where \(u_P=\mathcal E_P s_P\). Because
\(\mathcal E_P\) is ill-posed, \(s_P\) may have very large \(L^2(P_X)\) norm
even when \(u_P\) has very small \(L^2(P_Z)\) norm. Consequently, the minimal norm solution of the misspecified inverse problem
need not satisfy \cref{as:uniform-norm-bound} \cref{as:uniform-norm-bound:part:h}, even when the causal
dose-response curve \(h_{P;\mathrm{cdr}}\) is uniformly bounded.

For a fixed law \(P\), the difference between the misspecified and correctly
specified targets admits an easy to establish bound. Indeed, if \(q_P\) is the minimum norm solution to \(\mathcal E_P^*q_P=\psi\), then
\begin{equation*}
    \varphi_{\mathrm{bias}}(P)-\varphi(P)
    =
    E_P\{q_P(Z)u_P(Z)\}.
\end{equation*}
Therefore, by Cauchy--Schwarz,
\begin{equation}
\label{eq:rep_bias}
    |\varphi_{\mathrm{bias}}(P)-\varphi(P)|
    \le
    \|q_P\|_{L^2(P_Z)}\|u_P\|_{L^2(P_Z)}.
\end{equation}
Consequently, at a fixed law \(P\), as the direct effect \(u_P\) tends to zero
in \(L^2(P_Z)\), then the bias
\(|\varphi_{\mathrm{bias}}(P)-\varphi(P)|\) also goes to zero. That is, at each fixed law $P$,
small violations of the exclusion restriction produce small representation bias. This point was raised by
\cite{deaner_nonparametric_2022}.  Theorem 1.4 of that paper is a consequence of the inequality \eqref{eq:rep_bias}.

This pointwise sensitivity statement is not, however, a local Hellinger
continuity statement. To obtain a bound that is uniform over laws in a
Hellinger neighborhood of \(P\) following the argument of \cite{deaner_nonparametric_2022}, one would need a constant \(C<\infty\) such
that \(\|q_{P'}\|_{L^2(P'_Z)}\le C\) for all nearby laws \(P'\). This is a local version of \cref{as:uniform-norm-bound} \cref{as:uniform-norm-bound:part:q}. However, this is an assumption on the solution of the dual
inverse-problem with no apparent substantive interpretation in this setting.
In fact, without \cref{as:uniform-norm-bound}, the discontinuity results of \cref{subsec:failure-of-continuity} may apply and one may construct laws \(P_j\) converging to \(P\)
in Hellinger distance for which the misspecified parameter \(\varphi_{\mathrm{bias}}(P_j)\) diverges to $\infty$. On
the other hand, if the correctly specified target \(\varphi(P)\) is Hellinger
continuous, then \(\varphi(P_j)\to\varphi(P)\). Along such a sequence it must be
that
\begin{equation*}
    |\varphi_{\mathrm{bias}}(P_j)-\varphi(P_j)|
    \to
    \infty.
\end{equation*}
This does not contradict the inequality \eqref{eq:rep_bias}. Rather, it shows that \(\|q_{P_j}\|_{L^2(P_{Z;j})}\) must be diverging.

\end{example}

}

\begin{examplecont}{example:proximal}
\cref{as:uniform-norm-bound} is difficult to assess in the proximal causal inference application of \cref{example:proximal}.
There, the bridge functions are generally not themselves causal dose-response
curves and we are not aware of any substantive interpretation for them.    Even when the
observed outcome is bounded, an outcome bridge need not be uniformly bounded:
a conditional expectation equation can sometimes be solved only by a function
with large positive and negative values that cancel after averaging over the
proxy distribution. { On the other hand, under a correctly specified proximal causal model, if the outcome $W$ is bounded, the ATE $\varphi(P) = E_P[W(d=1) - W(d=0)]$ is bounded. Consequently, our constructions in \cref{lemma:discrete-necessary-sufficient-continuity} \cref{lemma:discrete-necessary-sufficient-continuity:discontinuity} and \cref{thrm:continuous-discontinuity} along which $\varphi(P_j) \rightarrow \infty$ must violate the proximal causal model. The specifics of this violation are important to investigate. If one or both bridge functions are truly uniformly bounded under the causal model, the parameter $\varphi$ will be Hellinger continuous. Alternatively, it could be possible that the causal model only restricts $\varphi$ to belong to $[\min(\Omega_W) - \max(\Omega_W), \max(\Omega_W) - \min(\Omega_W)]$, where $\Omega_W$ denotes the support of $W$, without sufficiently restricting the model to give Hellinger continuity of $\varphi$. In any case, as in the setting of misspecified NPIV, it seems plausible that even mild deviations from the proximal causal model could return us to the setting of \cref{subsec:failure-of-continuity}.}
\end{examplecont}

\begin{examplecont}{example:shadow}
 As noted in \cref{example:shadow}, the likelihood ratio
\(q_P(Z):=\{P(D=0)/P(D=1)\}\{dQ_Z/dU\}(Z)\) solves the dual integral equation
\eqref{eq:dual2}. Thus, in this example, \cref{as:uniform-norm-bound} \cref{as:uniform-norm-bound:part:q} on the
dual bridge is implied by a positivity-like overlap condition between
the respondent and nonrespondent outcome distributions, namely a uniform bound on
\(\|q_P\|_{L^2(U)}\). This overlap condition rules out, in an
\(L^2(U)\) sense, the possibility that values of \(Z\) that occur only rarely
among respondents can occur with high probability among nonrespondents.
\end{examplecont}

\subsection{Failure of Pathwise Differentiability under Uniform Norm Bounds}
\label{subsec:failure-differentiability}

Hellinger continuity of $\varphi$ was recovered in \cref{subsec:recover-continuity} by bounding uniformly over the model the minimum norm solution of the primal inverse problem, or symmetrically by bounding uniformly the minimum norm solution of the dual inverse problem. We now show that, under regularity conditions, these uniform norm bounds alone do not suffice to ensure pathwise differentiability. Under the model $\mc{P}_{h_{P}^\dagger \text{ bound}}$, pathwise differentiability fails at laws $P_0$ where both inverse problems are ill-posed if the primal inverse problem admits multiple solutions, as we show below. Even under the stronger model $\mc{P}_{h_{P}^\dagger, q_{P}^\dagger \text{ bound}}$, pathwise differentiability fails at laws where both inverse problems are ill-posed if both inverse problems admit multiple solutions. On the other hand, in \cref{subsec:recover-differentiability} below, we consider conditions under which pathwise differentiability holds when both inverse problems are ill-posed. We show that at $P_0$ where the primal inverse problem admits a unique solution, $\varphi$ is pathwise differentiable at $P_0$ relative to $\mc{P}_{h_{P}^\dagger \text{ bound}}$, under regularity conditions. If we instead assume the stronger model that uniformly bounds the norm of the $\beta$-source representer of the primal solution, we show that pathwise differentiability holds whether or not the primal solution is unique. In both cases, we provide a bound on the tangent spaces. Analogous results can be stated in terms of the dual inverse problem.

We provide the following table summarizing the forthcoming results at laws $P_0$ such that both the primal and dual inverse problems are ill-posed:
\begin{center}
    \begin{tabular}{llcl}
        \toprule
        Norm bounded over the model & Unique solution at $P_0$ & $\varphi$ differentiable at $P_0$ & Result\\
        \midrule
        $h^\dagger_P$ & primal & yes & \cref{thrm:sufficient-pathwise-diff}\\
        $h^\dagger_P$ & dual only & no & \cref{thrm:mismatched-injective-not-differentiable}\\
        $h^\dagger_P$ and $q^\dagger_P$ & neither & no & \cref{thrm:not-injective-not-differentiable}\\
        $\omega_{P;\beta_h}^\dagger$ & neither & yes & \cref{thrm:sufficient-pathwise-diff}\\
        \bottomrule
    \end{tabular}
\end{center}
In the above, $\omega_{P;\beta_h}^\dagger$ is the $\beta$-source representer of $h_P^\dagger$, satisfying $(\mc{E}_{P}^*\mc{E}_{P})^{\beta_h/2}\omega_{P;\beta_h}^\dagger = h_{P}^\dagger$, which we define formally in \cref{subsec:recover-differentiability}.

In the statistical literature, assumptions that the solutions to the primal and/or dual inverse problems are unique are often referred to as completeness assumptions \cite{newey_instrumental_2003-2,darolles_nonparametric_2011-1,horowitz_applied_2011}. When $(X, Z)$ are finite discrete, such completeness assumptions hold under mild assumptions for $X, Z$ taking the same number of levels. However, for $(X, Z)$ continuous, such completeness assumptions are difficult to reason about. See \cite{huNONPARAMETRICIDENTIFICATIONUSING2018} for some sufficient conditions for completeness. Difficulty arises in part because, within the class of distributions considered by \cite{canayTestabilityIdentificationNonparametric2013a}, at any distribution $P$ for which the completeness assumptions hold, there exists another distribution $P'$, arbitrarily close to $P$ in Hellinger distance, for which the completeness assumptions do not hold. Consequently, in light of Le Cam's two point method, for continuous $(X, Z)$ it is impossible to construct a valid test of the null hypothesis that completeness fails with power exceeding its size against any alternative at which completeness holds \parencite{canayTestabilityIdentificationNonparametric2013a}. Such difficulty has been extensively discussed in the setting of NPIVs, where the existence of a unique solution to the primal inverse problem is equivalent to identification of the structural function of interest, instead of simply certain linear functionals of it \parencite{santos_inference_2012, severiniEfficiencyBoundsEstimating2012}.

Our present results show that even if uniqueness of solutions is not required for $\sqrt{n}$-asymptotically linear inference, uniqueness of the solution whose norm is uniformly bounded is required for regular asymptotically linear inference under the uniform norm bounds of \cref{as:uniform-norm-bound} alone. As a consequence, Wald-type confidence intervals based on the standard DML estimator under-cover asymptotically under some local alternatives (\cref{lemma:failure-differentiability-irregular}), a fact which is hidden by a pointwise asymptotic analysis of the estimator. We show this poor behavior in a simple simulation in \cref{sec:simulation} below.

In the following results, we assume the Riesz representers $P \mapsto \alpha_P$, $P \mapsto r_P$ are pathwise differentiable at $P_0$ relative to a submodel of $\mc{P}_{h_{P}^\dagger \text{ bound}}$ in \cref{thrm:sufficient-pathwise-diff}, or $\mc{P}_{\text{Riesz}}$ in \cref{thrm:injective-tangent-space}. Recall from \cref{subsec:semiparametrics} that pathwise differentiability of Hilbert valued parameters is equivalent to the existence of the local parameter. Application of Example 5 of \cite{luedtke_one-step_2024} shows that this assumption holds relative to the model $\mc{P}_{\text{Riesz}}$ under mild conditions on $P_0$ in \cref{example:npiv}, \cref{example:proximal}, and \cref{example:shadow}, as we show in \cref{sup-sec:examples}. As $\mc{P}_{h_{P}^\dagger \text{ bound}} \subseteq \mc{P}_{\text{Riesz}}$, these conditions also suffice for differentiability in any submodel of $\mc{P}_{h_{P}^\dagger \text{ bound}}$. We also recall that the primal solution is unique at $P_0$ if and only if $\mathscr{N}(\mc{E}_{P_0;d}) = \{0\}$, and that the dual solution is unique if and only if $\mathscr{N}(\mc{E}^*_{P_0;d}) = \{0\}$ for each $d$.

\begin{theorem}
    \label{thrm:mismatched-injective-not-differentiable}
    Suppose that $P_0 \in \mc{P}_{h_{P}^\dagger \text{ bound}}$ and that $(X, Z)$ are continuous. Suppose there exists $d_e \in \textsf{Supp}[P_{D;0}]$ such that $\alpha_{P_0}(\cdot, d_e), r_{P_0}(\cdot, d_e)$ are not a.e. constant functions, that
    \begin{align*}
        \mathscr{N}(\mc{E}_{P_0;d_e}) \cap L^\infty(P_{X|d_e;0}) \not = \{0\},
    \end{align*}
    that there exists a singular system $(\sigma_{i;d_e}, u_{i;d_e}, v_{i;d_e})_{i=1}^\infty$ for $\mc{E}_{P_0;d_e}$, consisting of the non-zero singular values arranged in non-increasing order together with associated singular functions, satisfying
    \begin{align}
        \label{eq:mismatched-injective-not-differentiable:uniform-summability}
        \sum_{i=1}^\infty \sigma_{i;d_e}\|u_{i;d_e}\|_{L^\infty(P_{X|d_e;0})}\|v_{i;d_e}\|_{L^\infty(P_{Z|d_e;0})} < \infty, \quad \sum_{i=1}^\infty \sigma_{i;d_e}\|u_{i;d_e}\|_{L^\infty(P_{X|d_e;0})}\|v_{i+1;d_e}\|_{L^\infty(P_{Z|d_e;0})} < \infty,
    \end{align}
    and that $\mc{E}_{P_0;d_e}^*$ is injective. Then, $P \mapsto \varphi(P)$ is not pathwise differentiable at $P_0$ relative to the set of all one-dimensional submodels of $\mc{P}_{h_{P}^\dagger \text{ bound}}$ that are QMD at $P_0$.
\end{theorem}
\begin{remark}
    We prove this theorem by constructing a collection of submodels for which $\varphi$ does not admit a local parameter satisfying \eqref{eq:local-parameter-def}. \eqref{eq:mismatched-injective-not-differentiable:uniform-summability} and the assumption that $\mathscr{N}(\mc{E}_{P_0;d_e})$ contains uniformly bounded elements are technical conditions that ensure our proposed submodels remain valid probability distributions. There are probability distributions $P_0 \in \mc{P}$ that satisfy these restrictions and probability distributions $P_0 \in \mc{P}$ that do not. We conjecture that a more sophisticated construction would prove this theorem without these regularity conditions.
\end{remark}
The above result applies even if we further restricted ourselves to laws $P_0$ at which the Riesz representers are pathwise differentiable. In the proof of the above theorem, we construct a collection of submodels along which pathwise differentiability fails but for which $\alpha_P, r_P$ are both constant. Hence, along those submodels, they are pathwise differentiable with local parameter equal to 0. Indeed, this demonstrates failure of pathwise differentiability is intrinsic to the inverse problem, and not due to any pathologies in the Riesz representers $\alpha_P$ and $r_P$.

The above result does not apply to finite discrete random variables because, in such cases, injectivity of $\mc{E}_{P_0}^*$ ensures the pseudo-inverse operator $t \mapsto \mc{E}_{P_t}^\dagger$ is differentiable. Consequently, application of the chain rule proves that $t \mapsto \varphi(P_t) = E_{P_t}[\alpha_{P_t}(X, D) (\mc{E}_{P_t}^\dagger r_{P_t})(X, D)]$ is pathwise differentiable at $P_0$ over the model $\mc{P}_{h_{P}^\dagger \text{ bound}}$ whenever both $\alpha_{P_t}$ and $r_{P_t}$ are pathwise differentiable at $P_0$ over the model $\mc{P}_{h_{P}^\dagger \text{ bound}}$.

We next give a complementary theorem that proves pathwise differentiability can fail even when both $h_{P}^\dagger$ and $q_{P}^\dagger$ are uniformly bounded in norm. We now require that multiple solutions exist to the primal and dual inverse problems at $P_0$, as \cref{thrm:sufficient-pathwise-diff} below shows that pathwise differentiability holds were either solution unique.
\begin{theorem}
    \label{thrm:not-injective-not-differentiable}
    Let $P_0 \in \mc{P}_{h_{P}^\dagger, q_{P}^\dagger \text{ bound}}$. Suppose that there exists $d_e \in \textsf{Supp}[D;P_0]$ such that $\alpha_{P_0}(\cdot, d_e)$ and $r_{P_0}(\cdot, d_e)$ are both not constant functions and
    \begin{align*}
        \mathscr{N}(\mc{E}_{P_0;d_e})\cap L^\infty(P_{X, D;0}) \not= \{0\}, && \mathscr{N}(\mc{E}_{P_0;d_e}^*) \cap L^\infty(P_{Z, D;0}) \not= \{0\}.
    \end{align*}
    Then $P \mapsto \varphi(P)$ is not pathwise differentiable at $P_0$ relative to the set of all one-dimensional submodels of $\mc{P}_{h_{P}^\dagger, q_{P}^\dagger \text{ bound}}$ that are QMD at $P_0$.
\end{theorem}
As with \cref{thrm:mismatched-injective-not-differentiable}, the above result applies even at laws $P_0$ for which $P \mapsto \alpha_P$ and $P \mapsto r_P$ are both pathwise differentiable. We again conjecture that the assumption that $\mathscr{N}(\mc{E}_{P_0;d_e})$ and $\mathscr{N}(\mc{E}_{P_0;d_e}^*)$ contain uniformly bounded elements could be relaxed with more sophisticated arguments. We note that this assumption is trivially true when $X, Z$ are finite discrete and those null spaces are non-trivial.

\cref{thrm:mismatched-injective-not-differentiable} and \cref{thrm:not-injective-not-differentiable} have implications for inference. Under the conditions of those theorems, all asymptotically linear estimators, and in particular the estimator \eqref{eq:dml-estimator}, are irregular: their limiting distributions are sensitive to vanishing perturbations of the true data generating process. In particular, by \cref{lemma:failure-differentiability-irregular}, under the conditions of either theorem, every asymptotically linear estimator exhibits non-vanishing local asymptotic bias along some sequence of local alternatives, so that Wald-type confidence intervals under-cover asymptotically. See also the references on irregular estimation in \cref{subsubsec:implications-non-diff-main}.

We illustrate \cref{lemma:failure-differentiability-irregular} using the proof of \cref{thrm:not-injective-not-differentiable}, which exhibits a collection of QMD at $P_0$ submodels $\{\{P_{t, \gamma}:t\in [0, \delta)\}: \gamma \in (a_1, a_2)\}$ with scores $s_{\gamma} = s^* + \gamma \bar{s}$ affine in $\gamma$ (where $s^*, \bar{s} \in L^2(P_0)\setminus\{0\}$), but with derivatives $\left.\frac{d}{dt}\varphi(P_{t, \gamma})\right\vert_{t=0} = \frac{\mu}{\gamma}$ for some $\mu \not=0$. For any asymptotically linear estimator $\widehat{\varphi}_n$ with influence function $\varphi^1_{P_0}$, such as the DML estimator \eqref{eq:dml-estimator}, the local asymptotic bias under sampling from $P_{1/\sqrt{n}, \gamma}$ is
\begin{align*}
    \Delta_{\gamma} = \langle \varphi^1_{P_0}, s^*\rangle_{L^2(P_0)} + \gamma\langle \varphi^1_{P_0}, \bar{s}\rangle_{L^2(P_0)} - \frac{\mu}{\gamma}.
\end{align*}
Because $\mu \not= 0$, $\gamma\Delta_\gamma$ is a non-zero polynomial of degree at most two in $\gamma$, so $\Delta_\gamma = 0$ for at most two $\gamma \in (a_1, a_2)$. Hence every asymptotically linear estimator is locally asymptotically biased along all but at most two of these submodels, and the corresponding Wald-type confidence intervals under-cover asymptotically along them. We demonstrate poor coverage in a small simulation in \cref{sec:simulation} under a related data generating process, where the $\sqrt{n}$-scaled bias is essentially unchanged between $n = 200$ and $n = 200{,}000$ while the coverage of a nominal $95\%$ Wald interval falls from $65\%$ to $55\%$ over those sample sizes. In \cref{sec:weak-instruments}, we relate this phenomenon to weak-instrument asymptotics in linear instrumental variables models.

\subsection{Recovering Pathwise Differentiability}
\label{subsec:recover-differentiability}

As mentioned in \cref{subsec:failure-differentiability}, uniqueness of the primal solution suffices for pathwise differentiability at $P_0$ over the model $\mc{P}_{h_{P}^\dagger \text{ bound}}$, and a uniform norm bound on the $\beta$-source representer of $h_P^\dagger$ suffices even without uniqueness, under sufficient regularity. Before stating our results, we define the $\beta$-source representer of $h_P^\dagger$ and the associated model. Recall that \cref{as:beta-source} holds over the model $\mc{P}$. If it holds for some $\beta_h > 0$, then there exists a unique $\omega_{P;\beta_h}^\dagger \in \mathscr{N}(\mc{E}_{P})^\perp$ such that
\begin{align*}
    (\mc{E}_{P}^*\mc{E}_{P})^{\beta_h/2}\omega_{P;\beta_h}^\dagger = h_P^\dagger.
\end{align*}
We call such $\omega_{P;\beta_h}^\dagger$ the $\beta$-source representer of $h_P^\dagger$.
We now give the assumption bounding the norm of $\omega_{P;\beta_h}^\dagger$.
\begin{assumption}
    \label{as:beta-source-bound}
    Let $\beta_h > 0, C_{\beta_h} > 0$ be constants not depending on $P$. At $P$,
    \begin{align*}
        \|\omega_{P;\beta_h}^\dagger\|_{L^2(P_{X, D})} < C_{\beta_h}.
    \end{align*}
\end{assumption}
Then, for fixed $\beta_h > 0, C_{\beta_h} > 0$, we define the model
\begin{align}
    \label{eq:beta-source-bound-model}
    \mc{P}_{\beta_h \text{ bound}} \coloneqq \{P \in \mc{P}_{h_{P}^\dagger \text{ bound}}: \|\omega_{P;\beta_h}^\dagger \|_{L^2(P_{X, D})} < C_{\beta_h}\}.
\end{align}
Because the operator $(\mc{E}_{P}^*\mc{E}_{P})^{\beta_h / 2}$ has operator norm at most one, the equality $(\mc{E}_{P}^*\mc{E}_{P})^{\beta_h/2}\omega_{P;\beta_h}^\dagger = h_P^\dagger$ means that \cref{as:beta-source-bound} implies \cref{as:uniform-norm-bound} \cref{as:uniform-norm-bound:part:h}, with a possibly different constant. After presenting our results, we briefly discuss the interpretation of this assumption.

We now give the formal statement of our pathwise differentiability results.
\begin{theorem}
    \label{thrm:sufficient-pathwise-diff}
    Let $\mc{P}_{diff}$ and $P_0 \in \mc{P}_{diff}$ satisfy that $P \mapsto \alpha_P$ and $P \mapsto r_P$ are pathwise differentiable at $P_0$ relative to $\mc{P}_{diff}$ as $L^2(P_{X, D;0})$ and $L^2(P_{Z, D;0})$-valued parameters, respectively, and either:
    \begin{enumerate}
        \item \label{thrm:sufficient-pathwise-diff:injective} $\mc{P}_{diff} \subseteq \mc{P}_{h_{P}^\dagger \text{ bound}}$ and $\mc{E}_{P_0}$ is injective or
        \item \label{thrm:sufficient-pathwise-diff:beta-bound} $\mc{P}_{diff} \subseteq \mc{P}_{\beta_h \text{ bound}} $.
    \end{enumerate}
    Then, $P \mapsto \varphi(P)$ is pathwise differentiable at $P_0$ relative to $\mc{P}_{diff}$. If additionally $\widetilde{m}(W, X, D;h^\dagger_{P_0}) \in L^2(P_0)$, then the following is a gradient for any solution $q_{P_0} \in L^2(P_{Z, D;0})$ to the dual inverse problem that satisfies ${m}(W, Z, D;q_{P_0}) \in L^2(P_0)$:
    \begin{align}
        \label{eq:injective-gradient} \varphi^1(O;h_{P_0}^\dagger, q_{P_0}, \varphi(P_0)) = \widetilde{m}(W, X, D;h_{P_0}^\dagger) + m(W, Z, D;q_{P_0}) - h^\dagger_{P_0}(X, D)q_{P_0}(Z, D) - \varphi(P_0)
    \end{align}
\end{theorem}
A symmetric statement holds in terms of the dual inverse problem, which we do not restate here. If there exist multiple solutions to the dual inverse problem, the expression \eqref{eq:injective-gradient} reveals there may also exist distinct gradients. Because gradients are unique when the tangent space is non-parametric, the above theorem reveals that the assumptions in $\mc{P}_{diff}$ can restrict the observed data tangent space. The next theorem shows that the tangent space $\mc{T}(P_0, \mc{P}_{diff})$ is contained within a space that is a strict subspace of $L^2_0(P_0)$ when there exist multiple solutions to the dual inverse problem, under regularity conditions. As in \cref{lemma:constant-rank-tangent-space}, we require differentiability of $P \mapsto r_P$ relative to the larger model $\mc{P}_{\text{Riesz}}$ defined in \eqref{eq:riesz-model}. Under the conditions of the theorem, this model is locally non-parametric at $P_0$, so that the local parameter $\dot r_{P_0}$ has domain $L^2_0(P_0)$.
\begin{theorem}
    \label{thrm:injective-tangent-space}
    Suppose $m, \widetilde{m}, P_0$ satisfy \cref{as:m-m-tilde-bounded}, that at $P_0$, $r_P$ is pathwise differentiable relative to $\mc{P}_{\text{Riesz}}$, and that $\mc{P}_{diff}$ and $P_0 \in \mc{P}_{diff}$ satisfy either of the two enumerated conditions of \cref{thrm:sufficient-pathwise-diff}.
    Define
    \begin{align*}
        \mc{C}_{diff}(P_0) \coloneqq \left\{s \in L^2_0(P_0): \dot r_{P_0}[s](Z, D) - E_{P_0}[h_{P_0}^\dagger(X, D)s_{X|Z, D}(X, Z, D)|Z, D] \in \mathscr{N}(\mc{E}_{P_0}^*)^\perp\right\},
    \end{align*}
    where $s_{X|Z, D}(X, Z, D) \coloneqq E_{P_0}[s(O)|X, Z, D] - E_{P_0}[s(O)|Z, D]$. Then,
    \begin{align*}
        \mc{T}(P_0, \mc{P}_{diff}) \subseteq \overline{\mc{C}_{diff}(P_0)}.
    \end{align*}
    If ${m}(W, Z, D;\eta) \in L^2(P_0)$ for every $\eta \in \mathscr{N}(\mc{E}_{P_0}^*)$, then $\mc{C}_{diff}(P_0)$ can be equivalently expressed as
    \begin{align*}
        \mc{C}_{diff}(P_0) = \left\{s \in L^2_0(P_0): \langle s, m(\cdot;\eta) - h_{P_0}^\dagger \eta \rangle_{L^2(P_0)} = 0 \text{ for all }\eta \in \mathscr{N}(\mc{E}_{P_0}^*)\right\}
    \end{align*}

    If in addition there exists $d_e \in \textsf{Supp}[P_{D;0}]$ such that $\mathscr{N}(\mc{E}_{P_0;d_e}^*) \not= \{0\}$ and $r_{P_0}(\cdot, d_e)$ is a non-constant function of $Z$,
    \begin{align*}
        \overline{\mc{C}_{diff}(P_0)} \subsetneq L^2_0(P_0).
    \end{align*}
\end{theorem}
We note the similarity between the present result and the results in \cite{severiniEfficiencyBoundsEstimating2012} and \cite{chen_overidentification_2018-1}. As discussed in \cref{subsec:tangent-spaces}, those authors computed a tangent space only relative to the submodels in a set similar to \eqref{eq:differentiable-solution-paths}. Our results do not restrict attention to the set of submodels in \eqref{eq:differentiable-solution-paths}. As before, not all QMD at $P_0$ submodels of $\mc{P}_{diff}$ necessarily belong to the set in \eqref{eq:differentiable-solution-paths} when $X, Z$ are continuous or when both $\mc{E}_{P_0;d}, \mc{E}_{P_0;d}^*$ are not injective for some $d \in \textsf{Supp}[P_{D;0}]$. Instead, we show that a similar argument to bound the tangent space to the one used by \cite{severiniEfficiencyBoundsEstimating2012} and \cite{chen_overidentification_2018-1} applies even though $t \mapsto h_{P_t}^\dagger$ need not be differentiable. As such, we arrive at the same restrictions on the score $s$ given by those authors. In contrast to those authors, we do not claim that $\mc{T}(P_0, \mc{P}_{diff}) = \overline{\mc{C}_{diff}(P_0)}$, as we do not assume there exists a set of QMD at $P_0$ submodels whose scores are dense in $\overline{\mc{C}_{diff}(P_0)}$. We leave the construction of such sets, if they exist, to future work.

The assumption that \cref{as:beta-source-bound} holds uniformly over the model is considerably stronger than \cref{as:uniform-norm-bound}. Indeed, in the proof of \cref{lemma:weak-convergence-of-solutions} \cref{lemma:weak-convergence-of-solutions:beta-bound-result}, we show that it suffices to ensure that $h_P^\dagger$ converges in $L^2(P_{X, D;0})$ to $h_{P_0}^\dagger$. This result is related to \cite{babiiHonestConfidenceSets2020}, who show that under this and other assumptions, it is possible to construct uniformly honest confidence sets for the structural function $h_{P_0}^\dagger$.  Unlike \cref{as:uniform-norm-bound}, which follows from structural assumptions in \cref{example:npiv} and \cref{example:shadow}, we are not aware of any structural reason why this stronger bound should hold in any of \cref{example:npiv}, \cref{example:proximal}, or \cref{example:shadow}. While it may be possible to test for violations of \cref{as:beta-source-bound}, a similar phenomenon to the impossibility of testing completeness \parencite{canayTestabilityIdentificationNonparametric2013a} likely applies for testing the converse when $X, Z$ are continuous. Informally, it should be possible to perturb the tail of the singular value decomposition of $\mc{E}_{P}$ to violate \cref{as:beta-source-bound} while changing the law $P$ only negligibly in Hellinger distance, at least under certain regularity conditions.

The assumption is easiest to assess in the finite discrete case. Recall the interpretation of the singular functions $u_{P;i}$ and $v_{P;i}$ of $\mc{E}_{P}$ as canonical treatment and instrument contrasts, and the singular values $\sigma_{P;i}$ as measuring the strength of these canonical contrasts. Then, \cref{as:beta-source-bound} requires that
\begin{align*}
    \max_{i: \sigma_{P;i} > 0} \left|\frac{\langle v_{P;i}, r_P\rangle_{L^2(P_{Z, D})}}{\sigma^{1 + \beta_h}_{P;i}} \right|
\end{align*}
be uniformly bounded.
In light of this fact, \cref{as:beta-source-bound} admits an interpretation as a type of strong instrument condition: uniformly over the model, the Riesz representer $r_{P}$ may be supported on a canonical instrument $v_{P;i}$ only to a degree that vanishes faster than the strength $\sigma_{P;i}$ of that instrument. In particular, if the strength $\sigma_{P;i}$ of the canonical treatment-instrument pair $u_{P;i}$, $v_{P;i}$ tends to zero, then $\langle v_{P;i}, r_P\rangle_{L^2(P_{Z, D})}$ must tend to zero at least as fast as $\sigma_{P;i}^{1 + \beta_h}$. We develop this interpretation further in \cref{sec:weak-instruments}.

\section{Relation to Weak Instruments}
\label{sec:weak-instruments}

A useful way to interpret the preceding results is through their
connection with weak-instrument phenomena in linear instrumental
variables models. We begin by recalling the familiar one-treatment,
one-instrument linear IV problem, and then show that the singular-value
decomposition of the NPIV conditional expectation operator transforms
the NPIV problem into an infinite collection of canonical linear IV
problems. This representation also clarifies an important distinction
between the present setting and ordinary scalar IV: the structural
function may contain unidentified components, while the target
functional remains identified because it places zero weight on those
components.

\subsection{The Scalar Linear IV Model}

Consider first the homoskedastic linear IV model with one endogenous
treatment $X$ and one instrument $Z$,
\begin{equation}
    W = X\theta_P + \epsilon_P,
    \qquad
    X = Z\sigma_P + \xi_P,
    \qquad
    E_P(\epsilon_P\mid Z)=E_P(\xi_P\mid Z)=0,
    \label{eq:scalar_lin_iv}
\end{equation}
where $\theta_P$ is the structural coefficient and $\sigma_P$ is the
first-stage coefficient. We allow $E_P(\epsilon_P\mid X)\neq 0$, so
that ordinary least squares need not identify $\theta_P$.

Substituting the first-stage equation into the structural equation gives
a reduced-form regression of $W$ on $Z$ with coefficient
$\delta_P=\sigma_P\theta_P$. Thus, whenever $\sigma_P\neq 0$,
\begin{equation}
    \theta_P=\frac{\delta_P}{\sigma_P}.
    \label{eq:scalar_iv_ratio}
\end{equation}
Population two-stage least squares recovers precisely this ratio, and
sample two-stage least squares replaces $\delta_P$ and $\sigma_P$ by their
estimated reduced-form and first-stage coefficients.

The weak-instrument problem is therefore fundamentally a ratio problem.
When $\sigma_P$ is bounded away from zero, the map
$(\delta_P,\sigma_P)\mapsto\delta_P/\sigma_P$ is smooth. When $\sigma_P$ is
close to zero, small changes in the first-stage and reduced-form
coefficients can produce large changes in their ratio.

Different asymptotic regimes formalize different forms of this
instability. In the very-weak-instrument setting, the model permits
$\sigma_P$ to be arbitrarily close to, or equal to, zero without a
restriction preventing $\delta_P/\sigma_P$ from becoming arbitrarily
large. \cite{gleser_nonexistence_1987} and \cite{dufour_impossibility_1997} show that this has
severe consequences for uniformly valid inference. In the
weak-instrument asymptotics of \cite{staiger_instrumental_1997}, the
data-generating law is allowed to depend on $n$ and the first-stage
coefficient is local to zero, for example
$\sigma_{P^{(n)}}=\eta/\sqrt{n}$. The reduced-form coefficient is then
also of order $1/\sqrt{n}$ when the structural coefficient remains of
order one. Standard Wald inference based on two-stage least squares is
not uniformly approximated by the usual strong-IV limit theory in this
regime. In the classical strong-instrument regime, the relevant
first-stage coefficient remains locally separated from zero and the
ordinary smooth asymptotic theory applies.

We now show that these phenomena have direct analogues in the inverse
problems studied in this paper.

\subsection{NPIV as a Collection of Canonical Linear IV Problems}

Consider the NPIV model \eqref{eq:npiv-model} of \cref{example:npiv}.
The structural function $h_{P;\mathrm{cdr}}$ satisfies the conditional moment restriction \eqref{eq:primal-npiv-eq}. Under the conditions imposed throughout this paper, the linear operator
\[
    \mc{E}_P:L^2(P_X)\longrightarrow L^2(P_Z),
    \qquad
    \mc{E}_Ph=E_P\{h(X)\mid Z=\cdot\}
\]
is compact. Thus, there exists a singular value decomposition $\{(\sigma_{P;i},u_{P;i},v_{P;i}): i \in I_+(P), \sigma_{P;i}>0, u_{P;i}\in L^2(P_X), v_{P;i}\in L^2(P_Z)\}$ where the set $I_+(P)$ is a countable set of indexes, $\{u_{P;i}:i\in I_+(P)\}$ is an
orthonormal basis of $\mathscr{N}(\mc{E}_P)^\perp$, and
$\{v_{P;i}:i\in I_+(P)\}$ is an orthonormal basis of
$\mathscr{N}(\mc{E}_P^*)^\perp$ satisfying
\begin{equation}
    \mc{E}_Pu_{P;i}=\sigma_{P;i}v_{P;i},
    \qquad
    \mc{E}_P^*v_{P;i}=\sigma_{P;i}u_{P;i}.
    \label{eq:positive_singular_system}
\end{equation}

The singular functions $u_{P;i}$ have an immediate IV interpretation.
Specifically, define the $i$th canonical endogenous regressor by
$U_{P;i}:=u_{P;i}(X)$ and the corresponding canonical instrument by
$V_{P;i}:=v_{P;i}(Z)$. Equation~\eqref{eq:positive_singular_system}
implies
$E_P(U_{P;i}\mid Z)=\sigma_{P;i}V_{P;i}$. Moreover, orthonormality
gives
$E_P\{V_{P;i}U_{P;j}\}=\sigma_{P;j}I(i=j)$.
Thus, in the singular coordinates, the population first stage is
diagonal. The quantity $\sigma_{P;i}$ is the first-stage strength of
the canonical endogenous regressor $U_{P;i}$ relative to its
canonical instrument $V_{P;i}$. A small $\sigma_{P;i}$ therefore
corresponds to a weakly instrumented direction.

Now the structural function $h_{P;\mathrm{cdr}}$ can be uniquely decomposed as
\[
    h_{P;\mathrm{cdr}}=h_P^\dagger+h_{P;\mathrm{cdr};\perp}
\]
where  $h_P^\dagger$ is the minimum $L^2(P_X)$-norm solution to the primal inverse problem \eqref{eq:primal-npiv-eq} and $h_{P;\mathrm{cdr};\perp}\in \mathscr{N}(\mc{E}_P)$. Letting $\{u^0_{P;j},j \in I_0(P)\}$ denote an orthonormal basis of $\mathscr{N}(\mc{E}_P)$, we arrive at the expression
\[
    h_{P;\mathrm{cdr}}=\sum_{i\in I_+(P)}
    \theta_{P;i}u_{P;i} + \sum_{j\in I_0(P)}
    \theta^0_{P;j}u^0_{P;j}.
\]
 The coefficients $\theta_{P;i}, i \in I_+(P)$, are identified because taking the $L^2(P_Z)$ inner product of both sides of
$\mc{E}_P h_P^\dagger=r_P$
with $v_{P;i}$  gives
\begin{equation}
    \delta_{P;i}
    :=
    \langle r_P,v_{P;i}\rangle_{L^2(P_Z)}
    =
    \sigma_{P;i}\theta_{P;i},
    \text{   and consequently   }
    \theta_{P;i}
    =
    \frac{\delta_{P;i}}{\sigma_{P;i}}.
    \label{eq:canonical_iv}
\end{equation}
Because $r_P(Z)=E_P(W\mid Z)$,
$\delta_{P;i}=E_P\{Wv_{P;i}(Z)\}$. Thus
$\delta_{P;i}$ is the reduced-form coefficient of the outcome on the
$i$th canonical instrument. We thus see that Equation~\eqref{eq:canonical_iv} is exactly the scalar IV identity
\eqref{eq:scalar_iv_ratio} applied separately to every canonical
direction $u_{P;i}(X), i \in I_+(P)$.

On the other hand, whenever $\mathscr{N}(\mc{E}_P)\neq\{0\}$, the coefficients
$\theta^0_{P;j}$, $j\in I_0(P)$, are not identified by the conditional
moment restriction. Indeed, the corresponding canonical endogenous
regressors
\(
    U^0_{P;j}:=u^0_{P;j}(X), j\in I_0(P),
\)
are completely uninstrumented, since
\(
    E_P(U^0_{P;j}\mid Z)
    =
    \mc{E}_Pu^0_{P;j}
    =
    0.
\)
Consequently, changing any of the coefficients $\theta^0_{P;j}$ leaves
$\mc{E}_P h_{P;\mathrm{cdr}}$ unchanged, and hence leaves the conditional
moment restriction \eqref{eq:primal-npiv-eq} unchanged.

For the identified linear component of the structural function, and its contribution to the target $\varphi(P) = \langle \psi, h_{P;\mathrm{cdr}}\rangle_{L^2(P_X)}$ of \cref{example:npiv}, the correspondence with the linear IV problem is
\[
\begin{array}{c|c}
\text{NPIV canonical coordinate} & \text{linear IV interpretation}
\\ \hline
u_{P;i}(X) & \text{endogenous regressor} \\
v_{P;i}(Z) & \text{instrument} \\
\sigma_{P;i} & \text{first-stage coefficient} \\
\delta_{P;i}=\langle r_P,v_{P;i}\rangle_{L^2(P_Z)} & \text{reduced-form coefficient} \\
\theta_{P;i}=\delta_{P;i}/\sigma_{P;i} & \text{structural coefficient} \\
\psi_{P;i}=\langle \psi, u_{P;i}\rangle_{L^2(P_X)} & \text{target weight on the structural coefficient} \\
\psi_{P;i}\theta_{P;i} & \text{contribution to the target}.
\end{array}
\]
The last two rows have no counterpart in the scalar model \eqref{eq:scalar_lin_iv}, in which the structural coefficient is itself the target. As we discuss below, the behavior of the target weights is what distinguishes the consequences of weak instruments in the present setting from those in the scalar model.

When $X$ and $Z$ are continuously distributed, and the conditional
expectation operator $\mc{E}_P$ has infinite rank, then its
positive singular values can be arranged so that
$\sigma_{P;i}\downarrow 0$ as $i\to\infty$. Hence, the NPIV problem contains an infinite sequence of canonical
endogenous regressors whose corresponding instruments become
arbitrarily weak. If instead $\mc{E}_P$ has finite rank, there are only
finitely many positive singular values; any remaining directions in
$L^2(P_X)$ belong to $\mathscr{N}(\mc{E}_P)$ and are therefore completely
uninstrumented rather than weakly instrumented.

\subsection{Identification in Weak Instruments versus NPIV}

The preceding representation also helps clarify the distinction between
identification of the structural function and identification of the
target functional. Specifically, as argued above, whenever $\mathscr{N}(\mc{E}_P)\neq\{0\}$, the coefficients
$\theta^0_{P;j}$, $j\in I_0(P)$, are not identified, and consequently neither is the structural function $h_{P;\mathrm{cdr}}$.

The target of \cref{example:npiv} is nevertheless identified. Recall that
$\alpha_P=\psi$ and that \cref{as:solutions-exist} requires
$\psi\in \mathscr{R}(\mc{E}_{P}^*)$. Therefore
$\psi\in \mathscr{N}(\mc{E}_{P})^\perp$ and
$\langle\psi,u^0_{P;j}\rangle_{L^2(P_X)}=0$ for every null direction.
Consequently the target can be expressed as a linear combination of  the identified canonical
coefficients only:
\begin{align}
    \varphi(P)
    &= \sum_{i\in I_+(P)} \psi_{P;i}\theta_{P;i}
    \nonumber\\
    &=
    \sum_{i\in I_+(P)}
    \psi_{P;i}\frac{\delta_{P;i}}{\sigma_{P;i}},
    \qquad
    \psi_{P;i}
    :=
    \langle\psi,u_{P;i}\rangle_{L^2(P_X)}.
    \label{eq:canonical_target}
\end{align}
Because the target parameter places zero weight on coefficients
associated with null directions, $\varphi(P)$ is identified even when the structural function is not.

There is an additional useful consequence of the range condition $\psi\in \mathscr{R}(\mc{E}_{P}^*)$.
Fix $P$ and let $q_P$ be any solution of $\mc{E}_P^*q_P=\psi$. For
$i\in I_+(P)$, define
$q_{P;i}:=\langle q_P,v_{P;i}\rangle_{L^2(P_Z)}$. Then
$\psi_{P;i}
=\langle\psi,u_{P;i}\rangle_{L^2(P_X)}
=\sigma_{P;i}q_{P;i}$.
If $\mc{E}_P$ has infinite rank, index the positive singular values so that
$I_+(P)=\mathbb{N}$ and $\sigma_{P;i}\downarrow0$ as $i\to\infty$.
Because $\{v_{P;i}:i\geq1\}$ is an orthonormal sequence and
$q_P\in L^2(P_Z)$, the scalar sequence
$(q_{P;i})_{i\geq1}$ is square summable. Hence
$q_{P;i}\to0$ as $i\to\infty$. It follows that
$\psi_{P;i}/\sigma_{P;i}=q_{P;i}\to0$ as $i\to\infty$.
Equivalently,
\begin{equation*}
    \psi_{P;i}=o(\sigma_{P;i})
    \qquad\text{as }i\to\infty.
\end{equation*}

Thus, at a fixed law $P$, the range condition forces the target weights
to vanish faster than the corresponding singular values along the
sequence of increasingly weak canonical directions. If $\mc{E}_P$ has
finite rank, there are only finitely many positive singular values, so
there is no corresponding $i\to\infty$ statement.

This canonical-IV representation helps interpret the
different statistical models studied in \cref{sec:stat-models} and \cref{sec:norm-bound-models} and the consequences for inference about the target parameter, as we shall see next.

\subsection{The Baseline Model and the Analogue of Very Weak Instruments}

Consider first the baseline model $\mc{P}$. As we argue below, for local inference this
model shares important features with the scalar very weak instrument model
\eqref{eq:scalar_lin_iv}, which places no restrictions on either the
structural coefficient $\theta_P=\delta_P/\sigma_P$ or the first-stage
coefficient $\sigma_P$. The precise nature of this analogy differs,
however, between the finite-discrete and continuous settings, as we now discuss.

\paragraph{Finite-discrete $X$ and $Z$.}

When $X$ and $Z$ are finite discrete, there are only finitely many
canonical directions at any fixed law $P$.

If $\mc{E}_{P_0}$ has full rank, none of its singular values is zero, so its
smallest singular value is strictly positive. Moreover, since singular values vary continuously
with the finite-dimensional matrix representing $\mc{E}_P$, this smallest
singular value remains bounded away from zero for all laws sufficiently
close to $P_0$. Hence a canonical first stage cannot become
arbitrarily weak in a sufficiently small neighborhood of such a law. Thus, locally at $P_0$, the parameter $\varphi(P)$ behaves as in the very-strong IV model that places a uniform lower bound on $|\sigma_P|$.

In contrast, the discontinuity result in \cref{lemma:discrete-necessary-sufficient-continuity} concerns a rank-deficient
law $P_0$. Suppose, as in that lemma, that both $\mc{E}_{P_0}$ and
$\mc{E}_{P_0}^*$ have nontrivial null spaces. The proof selects null
directions $u^0_{P_0}\in\mathscr{N}(\mc{E}_{P_0})$ and
$v^0_{P_0}\in\mathscr N(\mc{E}_{P_0}^*)$. At $P_0$, these directions
correspond to a zero first-stage strength. In addition, the range
conditions imply that the corresponding target weight and reduced-form
coefficient are both zero: $\langle\psi,u^0_{P_0}\rangle_{L^2(P_{X;0})}=0$ because
$\psi\in \mathscr{R}(\mc{E}_{P_0}^*)$, and
$\langle r_{P_0},v^0_{P_0}\rangle_{L^2(P_{Z;0})}=0$ because
$r_{P_0}\in \mathscr{R}(\mc{E}_{P_0})$.
The proof then constructs nearby laws $P_j\to P_0$ for which a
canonical direction approaching these null directions acquires a
strictly positive singular value $\sigma_{P_j;i}$ satisfying
$\sigma_{P_j;i}\to0$ for a fixed index $i$.
Along the constructed sequence, the corresponding reduced-form coefficient and
target weight also converge to zero. However, the reduced-form
coefficient can converge to zero more slowly than the singular value, so
that the structural coefficient
$\theta_{P_j;i}=\delta_{P_j;i}/\sigma_{P_j;i}$ diverges. Moreover, it
can diverge faster than the target weight $\psi_{P_j;i}$ converges to
zero, causing the contribution
$\psi_{P_j;i}\theta_{P_j;i}$, and hence the target
itself, to diverge.

The analogy with the scalar very-weak-IV problem is therefore quite
direct in the finite-discrete case. The null direction at $P_0$
corresponds to $\sigma_{P_0}=0$, and nearby laws $P_j$ correspond to
scalar IV laws for which $\sigma_{P_j}>0$ but $\sigma_{P_j}\to0$ in such a way that $\delta_{P_j} / \sigma_{P_j}$ converges to $\infty$. There is,
however, a subtle difference. In the scalar problem,
$\sigma_{P_0}=0$ leaves the structural coefficient $\theta_{P_0}$
unidentified. In the NPIV problem, the structural coefficient along
the null direction is likewise unidentified, but the target
$\varphi(P_0)$ remains identified because it places zero weight on that
direction.

\paragraph{Continuously distributed $X$ and $Z$.}

The situation is different when $X$ and $Z$ are continuously
distributed because $L^2(P_X)$ and $L^2(P_Z)$ are infinite
dimensional, while $\mc{E}_P$ is compact under the conditions imposed in
this paper.

Suppose first that $\mc{E}_{P_0}$ has infinite rank. Its positive singular
values can then be indexed so that
$\sigma_{P_0;i}\downarrow0$ as $i\to\infty$. Thus, unlike in the
finite-discrete case, arbitrarily weak canonical instruments already
exist at the single fixed law $P_0$. This remains true even when
$\mc{E}_{P_0}$ is injective, in which case every canonical structural
coefficient is identified and there are no null directions: there are
nevertheless identified directions whose first-stage strengths are
arbitrarily close to zero.

The proof of the discontinuity result of \cref{thrm:continuous-discontinuity} constructs laws $P_j\to P_0$ by perturbing a singular value lying
increasingly far into the singular tail as $j$ grows. Thus the canonical direction itself
changes with $j$: one may think of the construction as selecting
indices $i_j\to\infty$, for which the associated first-stage
strengths are increasingly small, and then perturbing the law in
those increasingly weak directions. The perturbations can be made
arbitrarily small in Hellinger distance because they occur farther and
farther into the tail, while the corresponding singular value can be
made sufficiently small relative to the reduced-form coefficients and target
weights that the contribution to $\varphi(P_j)$ diverges.

Hence, in the continuous infinite-rank case, the counterpart of the
scalar sequence $\sigma_{P_j}\to0$ is not necessarily a single canonical
first-stage coefficient becoming weak as $P$ changes. Rather, it is
the existence of a sequence of law-direction pairs $(P_j,i_j)$ for
which $\sigma_{P_j;i_j}\to0$, with increasingly weak directions being
selected farther into the singular tail.

If $\mc{E}_{P_0}$ instead has finite rank, there are only finitely many
positive singular values, but there are null directions in the
infinite-dimensional domain. In that case, the source of arbitrarily
weak directions is closer to the finite-discrete rank-deficient
case: nearby laws can turn directions that are null, and hence
completely uninstrumented, at or near the limiting law into directions
with arbitrarily small positive singular values.

Thus the continuous case does not require an actual zero singular
value at $P_0$ for the discontinuity phenomenon. When $\mc{E}_{P_0}$ has
infinite rank, compactness itself supplies positive singular values
arbitrarily close to zero; when it has finite rank, null directions
are available instead. This is why \cref{thrm:continuous-discontinuity} does not require the
nontrivial-null-space conditions that appear in \cref{lemma:discrete-necessary-sufficient-continuity}.

Whether in the rank-deficient discrete case or in the continuous case, the baseline model $\mc{P}$ permits sequences of law-direction pairs for
which the canonical first-stage strength tends to zero without a
uniform restriction preventing the associated structural coefficient
$\delta_{P;i}/\sigma_{P;i}$ from becoming arbitrarily large. This is
the precise sense in which the discontinuity results are analogous to
the very-weak-instrument phenomenon in the scalar IV model.

\subsection{Bounded Structural Coefficients and Continuity}

Now consider the model
$\mc{P}_{h_{P}^\dagger \text{ bound}}$ defined through
\cref{as:uniform-norm-bound} \cref{as:uniform-norm-bound:part:h}. By orthonormality,
\[
    \|h_P^\dagger\|_{L^2(P_X)}^2
    =
    \sum_{i\in I_+(P)}\theta_{P;i}^2.
\]
Thus a uniform bound on $\|h_P^\dagger\|_{L^2(P_X)}$ implies a
uniform bound on every canonical structural coefficient
$\theta_{P;i}$. Importantly, this assumption does not require the instruments to be
strong: singular values may still approach zero. The restriction is
instead that the reduced-form coefficient must vanish at least
proportionately to the first-stage strength, so that
$\delta_{P;i}/\sigma_{P;i}$ cannot diverge. In the language of the scalar model \eqref{eq:scalar_lin_iv}, the model $\mc{P}_{h_{P}^\dagger \text{ bound}}$ bounds the structural coefficient $|\theta_P|$, but does not bound the first-stage coefficient $\sigma_P$ away from zero.

This explains why the uniform bound restores continuity when $P \mapsto \alpha_P$ and $P \mapsto r_P$ are also continuous at $P_0$. Consider again a canonical direction $u_{P;i}$ that converges to a null direction at $P_0$ as $P \rightarrow P_0$. Because $\psi \in \mathscr{N}(\mc{E}_{P_0})^\perp$, its target weight $\psi_{P;i}$ converges to zero, while its structural coefficient $\theta_{P;i}$ remains bounded. Their product, which is the contribution of this direction to $\varphi(P)$, therefore converges to zero. In the baseline model, by contrast, the structural coefficient may diverge quickly enough that this contribution diverges. Canonical directions whose singular values are bounded away from zero near $P_0$ are locally well behaved, and contribute continuously to $\varphi(P)$.

The same intuition extends to continuous $X, Z$, where the uniform norm bound additionally controls the contribution of the infinitely many weak directions in the singular tail; see the proof of \cref{thrm:bounded-means-continuous}.

\subsection{Weak-Instrument Asymptotics and Failure of Pathwise Differentiability}

Although bounded structural coefficients rule out the analogue of very weak instruments, they do not rule out the analogue of weak-instrument asymptotics. In the scalar model, weak-instrument asymptotics produce a structural coefficient that is a ratio of two quantities local to zero at rate $1/\sqrt{n}$, and which therefore cannot be consistently estimated. In NPIV, the same ratio can arise along a canonical direction, but it is multiplied by a target weight that is itself local to zero at rate $1/\sqrt{n}$. Consequently, the contribution of such a direction to the target vanishes, so that the target remains consistently estimable along these local alternatives, and the weak-instrument instability instead manifests as failure of pathwise differentiability.

The construction in the proof of \cref{thrm:not-injective-not-differentiable} exhibits this phenomenon at laws $P_0$ where both inverse problems admit multiple solutions, whether $X, Z$ are finite discrete or continuous. As in the proof of \cref{lemma:discrete-necessary-sufficient-continuity}, it selects null directions $u_{P_0}^0 \in \mathscr{N}(\mc{E}_{P_0})$ and $v_{P_0}^0 \in \mathscr{N}(\mc{E}^*_{P_0})$, which at $P_0$ have zero first-stage strength, zero target weight, and zero reduced-form coefficient. It then constructs a family of QMD at $P_0$ submodels $\{P_{t, \eta}:t\in[0, \delta_{\eta})\}$ of $\mc{P}_{h_{P}^\dagger \text{ bound}}$, indexed by a scalar $\eta$, along which a pair of canonical directions approaching $(u_{P_0}^0, v_{P_0}^0)$ acquires a strictly positive first-stage strength, with the orders of the relevant quantities given by
\begin{center}
    \begin{tabular}{lcc}
        \toprule
        Canonical quantity & Order along $P_{t,\eta}$ & Order under $P_{1/\sqrt{n},\eta}$ \\
        \midrule
        First-stage coefficient & $\eta t$ & $\eta/\sqrt{n}$ \\
        Reduced-form coefficient & $t$ & $1/\sqrt{n}$ \\
        Structural coefficient & $1/\eta$ & $1/\eta$ \\
        Target weight & $t$ & $1/\sqrt{n}$ \\
        Contribution to $\varphi$ & $t/\eta$ & $1/(\eta\sqrt{n})$ \\
        \bottomrule
    \end{tabular}
\end{center}
Under the local alternative $P_{1/\sqrt{n},\eta}$, the first-stage and reduced-form coefficients are both local to zero at rate $1/\sqrt{n}$, while the structural coefficient remains of constant order $1/\eta$: this is precisely the configuration of the weak-instrument asymptotics of \cite{staiger_instrumental_1997}. Because the structural coefficient remains bounded, the construction respects the uniform norm bound.

The contribution of this direction to the target vanishes, so that $\varphi(P_{t,\eta}) \approx \varphi(P_0) + t/\eta$ is continuous along each submodel. However, its derivative at $t = 0$ is approximately $1/\eta$, whereas the score of the submodel, $s_\eta = s^* + \eta\bar{s}$ for $s^*, \bar{s} \in L^2_0(P_0)$ not depending on $\eta$, is affine in $\eta$. Were $\varphi$ pathwise differentiable, its derivative along this family would be a bounded linear functional of $s_\eta$, and hence also affine in $\eta$. The term $1/\eta$ therefore rules out pathwise differentiability. The consequence for inference is the local asymptotic bias described following \cref{thrm:not-injective-not-differentiable}: under $P_{1/\sqrt{n},\eta}$, no asymptotically linear estimator is $\sqrt{n}$-unbiased for all $\eta$, and Wald-type confidence intervals under-cover (\cref{lemma:failure-differentiability-irregular}). Thus, the same $1/\eta$ term that renders the structural coefficient inconsistently estimable in the scalar weak-instrument problem appears here as non-vanishing $\sqrt{n}$-scaled bias of every asymptotically linear estimator of the target.

\paragraph{Pathwise differentiability when the primal solution is unique.}

The existence of a fixed null direction $u_{P_0}^0$ at $P_0$ is essential to the above construction. When $\mc{E}_{P_0}$ is injective, there is no null direction to turn on, and $P \mapsto \varphi(P)$ is pathwise differentiable over $\mc{P}_{h_{P}^\dagger \text{ bound}}$ under mild regularity conditions (\cref{thrm:sufficient-pathwise-diff} \cref{thrm:sufficient-pathwise-diff:injective}). When $X, Z$ are finite discrete, injectivity of $\mc{E}_{P_0}$ also bounds the smallest singular value away from zero locally, so that $P_0$ is locally analogous to a strong-instrument setting. When $X, Z$ are continuous, arbitrarily weak directions remain, but each has strictly positive first-stage strength $\sigma_{P_0;i}$ at $P_0$. Perturbing such a direction as above, with first-stage strength $\sigma_{P_0;i} + \eta t$ and reduced-form coefficient and target weight perturbed by terms of order $t$, changes its contribution to $\varphi$ by a term whose quadratic part, of order $t^2/(\sigma_{P_0;i} + \eta t)$, has zero derivative at $t = 0$. The weak-instrument term $1/\eta$ therefore does not arise; the proof of \cref{thrm:sufficient-pathwise-diff} controls the infinitely many such directions simultaneously.

\subsection{The \texorpdfstring{$\beta_h$}{beta}-Source Condition and Effective Strong Identification of the Target}
We finally consider the stronger model $\mc{P}_{\beta_h \text{ bound}}$ with $\beta_h > 0$. Recall that the $\beta_h$-source representer $\omega^\dagger_{P;\beta_h}$ is defined as the unique solution in $\mathscr{N}(\mc{E}_{P})^\perp$ satisfying
\begin{align*}
    (\mc{E}_{P}^*\mc{E}_{P})^{\beta_h/2}\omega^\dagger_{P;\beta_h} = h^\dagger_P.
\end{align*}
The $\beta_h$-source representer $\omega^\dagger_{P;\beta_h}$ can alternatively be represented in terms of the canonical treatment-instrument contrasts as
\begin{align*}
    \omega^\dagger_{P;\beta_h} = \sum_{i\in I_+(P)}\frac{\theta_{P;i}}{\sigma_{P;i}^{\beta_h}}u_{P;i},
\end{align*}
so that it is coordinate wise in $u_{P;i}$ equal to the structural coefficient inversely weighted by $\sigma_{P;i}^{\beta_h}$. Because $\mc{E}_P$ is a conditional expectation operator, it has operator norm at most one, and hence $\sigma_{P;i} \leq 1$, so that $|\theta_{P;i}|/\sigma_{P;i}^{\beta_h} \geq |\theta_{P;i}|$, with strict inequality whenever $\sigma_{P;i} < 1$ and $\theta_{P;i} \neq 0$. A uniform bound on $\|\omega^\dagger_{P;\beta_h}\|_{L^2(P_{X})}$ implies
\begin{align*}
    \|\omega^\dagger_{P;\beta_h}\|^2_{L^2(P_X)} = \sum_{i\in I_+(P)}\left(\frac{\theta_{P;i}}{\sigma_{P;i}^{\beta_h}}\right)^2 < C_{\beta_h}^2
\end{align*}
and in particular that $|\theta_{P;i}| \leq C_{\beta_h} \sigma_{P;i}^{\beta_h}$. Because $C_{\beta_h}$ does not depend on $P$ or $i$, the structural coefficient must therefore tend to zero along any sequence of laws and canonical directions $(P_j, i_j)$ for which $\sigma_{P_j;i_j} \rightarrow 0$. Thus, while the uniform $\beta_h$-source condition is not a literal strong instrument condition, as it does not require singular values to be bounded away from zero, it means that the structural coefficient for the $i$-th coordinate is no larger than a constant times the instrument strength $\sigma_{P;i}$ to the positive power ${\beta_h}$.

This condition breaks the counterexample used to show failure of pathwise differentiability at $P_0$ where $\mc{E}_{P_0}$ is not injective under the model $\mc{P}_{h_P^\dagger \text{ bound}}$. In that construction, a submodel $\{P_{t, \eta}: t\in[0, \delta_{\eta})\}$ is constructed so that for one coordinate, the reduced-form coefficient and the instrument strength are of the same order, both decaying to zero as $t \rightarrow 0$, and where the weight $\langle \psi, u_{P_t;i}\rangle_{L^2(P_{X;t})}$ is of order $t$. Consequently, the structural coefficient is of constant order. Such construction violates the uniform bound on the norm of the $\beta_h$-source representer, as for $t$ small enough $|\theta_{P_t;i}| \leq C_{\beta_h} \sigma_{P_t;i}^{\beta_h}$ will not hold. Under the $\beta_h$-source bound, the structural coefficient must instead converge to zero, so that the contribution of the $i$-th canonical direction to $\varphi$ behaves like
\begin{align*}
    \theta_{P_t;i} \langle \psi, u_{P_t;i}\rangle_{L^2(P_{X;t})} = o(t)
\end{align*}
which will not contribute to the derivative of $\varphi(P_t)$.

In light of the above discussion, the uniform $\beta_h$-source condition is intermediate between the framing of weak-instrument asymptotics and traditional strong instrument asymptotics. While the strength of instruments can be local to zero, the reduced-form coefficient on any weak canonical instrument is always of smaller order than the strength of that instrument, since $|\delta_{P;i}| = \sigma_{P;i}|\theta_{P;i}| \leq C_{\beta_h}\sigma_{P;i}^{1+\beta_h}$. The uniform norm bound on the $\beta_h$-source representer may then be interpreted as a type of strong instrument condition relative to the target of inference $\varphi$, insofar as it ensures that weak instruments do not contribute in first order to the behavior of the parameter.

\subsection{Implications of the Weak-Instrument Connection for \texorpdfstring{$X$, $Z$}{X, Z} Finite Discrete}
\label{subsec:weak-instrument-discrete}

As discussed following \cref{lemma:discrete-necessary-sufficient-continuity}, when $X$, $Z$ are discrete, the laws at which there exist multiple solutions to both the primal and dual inverse problems form a set of measure zero. Nevertheless, the connection between \cref{thrm:not-injective-not-differentiable} and the setting of weak instrument asymptotics gives intuition for why the present results are still actionable. In the setting of linear IV, the laws at which pathological behavior occurs ($\sigma_P = 0$) are Lebesgue measure zero. Nevertheless, when $\sigma_P$ is small relative to the sample size, i.e. of order $1/\sqrt{n}$, weak-instrument asymptotics has been demonstrated to give a better finite sample approximation \parencite{andrews_weak_2019}. Taking the same view, in the setting of proximal causal inference, discrete proxies may be weak for the corresponding discrete unmeasured confounder, and in the setting of shadow variables, the shadow variable may be only weakly associated with the outcome subject to missingness, even when it is discrete. If each of $\sigma_{P;i}, \langle \alpha_P, u_{P;i}\rangle_{L^2(P_X)}, \langle r_P, v_{P;i}\rangle_{L^2(P_Z)}$ is of order $1/\sqrt{n}$ (i.e. each is upper and lower bounded by approximately $1/\sqrt{n}$), Wald-type confidence intervals are likely to undercover in finite samples (see \cref{sec:simulation} for an illustration), and alternative interval estimates that account for the weak instruments/proxies/shadow variables may outperform them. See \cite{smuclerNotePropertiesConfidence2026} for a brief discussion of a related point, who argue that inverting a score test provides better finite sample coverage under weak-instrument asymptotics when $X, Z$ are binary. Similarly, it may prove fruitful to develop alternative confidence interval estimators that are locally uniformly honest under the model $\mc{P}_{h_P^\dagger \text{ bound}}$ at laws $P_0$ with multiple solutions to both inverse problems, even when $X, Z$ are discrete. We leave such constructions to future work.

\section{Relation to Results Showing Neyman Orthogonality}
\label{sec:relation-to-neyman-orthogonality}

{
The present results are compatible with standard results deriving asymptotically linear estimators for $\varphi$ using the estimating equation \eqref{eq:neyman-orthogonal-ee} at $P_0$ where the parameter $P \mapsto \varphi(P)$ is Hellinger discontinuous \parencite{bennettSourceConditionDouble2023b, olivas-martinez_source-condition_2025}. As discussed in \cref{subsubsec:implications-discontinuity-main}, pointwise asymptotic consistency of estimators and validity of confidence intervals are strictly weaker than locally uniform consistency of estimators and validity of confidence intervals. Inference can still be performed using standard pointwise asymptotic theory in this setting. However, in light of the finite sample results in \cref{lemma:finite-sample-discontinuity} applied at the laws described in \cref{subsec:failure-of-continuity}, the pointwise asymptotic approximation may be poor for every $n$: at each sample size there are laws in $\mc{P}$ that the data cannot distinguish from $P_0$ at which the approximation is arbitrarily poor. Additionally, because the parameter is not pathwise differentiable, all asymptotically linear estimators are irregular (\cref{lemma:failure-differentiability-irregular}), and semiparametric efficiency theory cannot be used to select between candidate asymptotically linear estimators.

This work emphasizes an important distinction between Neyman-orthogonality, the basis under which the estimating equation \eqref{eq:neyman-orthogonal-ee} is justified, and pathwise differentiability. The two concern different perturbations: Neyman-orthogonality controls sensitivity to perturbations of the nuisance functions, while pathwise differentiability concerns sensitivity to perturbations of the law itself. Viewing the estimating function \eqref{eq:neyman-orthogonal-ee} as the mapping
\begin{align*}
    \mathbb{R} \times  L^2(P_{X,D}) \times L^2(P_{Z,D}) \rightarrow L^2(P), (\varphi, h, q) \mapsto \varphi^1_P(\cdot;h, q, \varphi),
\end{align*}
Neyman orthogonality is a statement about the composition of this map with the expectation $E_{P}$: it states that $(h, q) \mapsto E_{P}[\varphi^1_P(O;h, q, \varphi(P))]$ has vanishing Gateaux derivative at $(h_P, q_P)$, where the nuisance space carries its $L^2$ topology. On the other hand, pathwise differentiability considers differentiability relative to the Hellinger topology, or rather relative to the Hellinger differentiable one dimensional paths in $\mathscr{P}(P_0, \mc{P})$. The latter topology is in some sense more natural for deriving statistical results, at least those related to asymptotic optimality, and is the one on which semiparametric efficiency theory is built \parencite{bickelEfficientAdaptiveEstimation1998a, vaart_asymptotic_1998-2}. Recent work has attempted to unify these related concepts; see for example \cite{chen_equivalence_2026}, who show that under regularity conditions and correct choice of nuisance parameters, pathwise differentiability and Neyman-orthogonality imply each other. The setting of this paper shows that this equivalence does not always hold.

The hypothesis that fails here is a smoothness requirement on the nuisance parameters. Letting $\eta(P)$ denote the nuisance parameters and $\beta(P)$ the target of inference, that work requires the mapping $t \mapsto \eta(P_t)$ be differentiable at $t=0$ whenever $\{P_t:t\in [0, \delta)\} \subseteq \mc{P}$ is a QMD at $P_0$ submodel. In the setting of this paper, $\eta(P) \coloneqq (h_P^\dagger, q_P^\dagger)$ and $\beta(P) = \varphi(P)$. As discussed in \cref{subsec:tangent-spaces}, the mapping $t \mapsto (h_{P_t}^\dagger, q_{P_t}^\dagger)$ fails to be differentiable along some $\{P_t:t\in [0, \delta)\}  \in \mathscr{P}(P_0, \mc{P})$. Because the requirement is imposed at every such submodel, a single submodel along which differentiability fails suffices to place the present setting outside the scope of that equivalence.

}

\section{Counterexamples and Simulations}
\label{sec:simulation}

We now provide a simple example demonstrating lack of continuity of the parameter $\varphi(P)$ and implications for inference in the setting of a finite discrete NPIV problem. We first describe a law $P_0$ and submodel $\{P_t: t\in [0, \delta)\}$ over which the target parameter fails to be continuous, and demonstrate the poor behavior of a non-parametric estimator of $\varphi(P)$ under selected local alternative data generating processes. One local alternative leaves every model with uniform norm bounds, illustrating the results of \cref{sec:stat-models} for the baseline model $\mc{P}$; another remains within $\mc{P}_{h_{P}^\dagger \text{ bound}}$, illustrating the results of \cref{sec:norm-bound-models}. We note that the construction used in this counterexample is different from that used in the proof of \cref{lemma:discrete-necessary-sufficient-continuity} to simplify presentation.
\begin{example}
    \label{counterexample:npiv}
    We return to the setting of NPIV with $D$ absent. In what follows, we do not endow $\mc{P}$ with further restrictions implied by the additive causal structure model, remaining agnostic about the practical interpretation of the parameter $\varphi$, as this example is purely for demonstration purposes. In particular, the model $\mc{P}$ does not impose a uniform bound on the norm of $h^\dagger_P$; below, we consider one path along which the minimum norm solutions diverge, and another along which they remain bounded. Suppose that $W, X, Z$ each take values in $\{1, 2, 3\}$. In the following discussion, we treat univariate functions of $W, X$, or $Z$ as vectors in $\mathbb{R}^3$, and the conditional expectation operators $\mc{E}_{P}, \mc{E}_{P}^*$ as $3\times 3$ matrices, where the $(i, j)$-th entry of $\mc{E}_{P}, \mc{E}_{P}^*$, is equal to $P(X = j|Z=i)$ and $P(Z = j|X=i)$, respectively. Suppose that
    \begin{align*}
        \mc{E}_{P_0} \coloneqq & \begin{pmatrix}
            0.4 & 0.2 & 0.4\\
            0.2 & 0.4 & 0.4\\
            0.2 & 0.4 & 0.4
        \end{pmatrix},
    \end{align*}
    and that $P_0(W=w|Z=z) = P_0(X=w|Z=z)$ for $w, z \in \{1, 2, 3\}$. Suppose that marginally, $P_0(Z=z) = (1/3, 1/3, 1/3)^T$. Let $\tau(w) \coloneqq w$ for $w \in \{1, 2, 3\}$, so that $\tau = (1, 2, 3)^T$ as a vector, and let $\psi(x) \coloneqq E_{P_0}[Z|X=x]$, so that $\psi = \mc{E}_{P_0}^*(1, 2, 3)^T = (7/4, 11/5, 2)^T$. We leave the full joint distribution of $(W, X, Z)$ unspecified for now, noting it could be explicitly given by choosing any conditional copula relating $X|Z$ and $W|Z$ for each $z$.

    Then, $\alpha_{P_0} = \psi \in \mathscr{R}(\mc{E}_{P_0}^*)$ and because $P_0(W=w|Z=z) = P_0(X=w|Z=z)$, $r_{P_0} = E_{P_0}[W|Z] = \mc{E}_{P_0}\tau \in \mathscr{R}(\mc{E}_{P_0})$. Hence, \cref{as:solutions-exist} is satisfied.  Note that $\psi$ is fixed, i.e. it does not change with $P$. \cref{as:continuous-linear-maps}, \cref{as:dominating-measure-support}, \cref{as:beta-source}, and \cref{as:mutual-absolute-continuity} are satisfied at $P_0$ by construction, and so $P_0 \in \mc{P}$ defined as in \eqref{eq:model-definition}.

    Note that $\mathscr{N}(\mc{E}_{P_0}) = \text{span}\{(2, 2, -3)^T\}$. Because $K_X = K_Z$ and $\mathscr{N}(\mc{E}_{P_0}) \not= \{0\}$, the matrix representation of $\mc{E}_{P_0}$ is square and rank deficient, and hence $\mathscr{N}(\mc{E}_{P_0}^*) \not= \{0\}$ as well. Then, the conditions of \cref{lemma:discrete-necessary-sufficient-continuity} \cref{lemma:discrete-necessary-sufficient-continuity:discontinuity} hold. We now give a path $\{P_t:t \in [0, \delta)\} \in \mathscr{P}(P_0, \mc{P})$ such that $\varphi(P_t) \rightarrow \infty$ as $t \rightarrow 0$.

    Let $P_t$ be such that the joint distribution of $(W, Z)$ under $P_t$ equals that distribution under $P_0$. Let
    \begin{align}
        \label{eq:counterexample-discontinuity}
        \mc{E}_{P_t} \coloneqq \mc{E}_{P_0} + B(t), && B(t) \coloneqq t\begin{pmatrix}
            0 & 0 & 0\\
            -1 & 1 - \beta(t) & \beta(t) \\
            0 & 0 & 0
        \end{pmatrix}
    \end{align}
    with $\beta(t) \coloneqq t^2$. Then, $P_t(X=x|Z=z) = P_{0}(X=x|Z=z) + B(t)_{z, x}$ for $x, z \in \{1, 2, 3\}$.

    Now, the determinant of $\mc{E}_{P_t}$ equals $-t^3/5$, and hence $\mc{E}_{P_t}$ is invertible for each $t\not=0$. As such, $r_{P_t} = r_{P_0} \in \mathscr{R}(\mc{E}_{P_t})$ and $\alpha_{P_t} = \psi \in \mathscr{R}(\mc{E}_{P_t}^*)$ for $t\in [0, \delta)$. Furthermore, if $\delta < 0.2$, the entries of the matrix representation $\mc{E}_{P_t}$ are non-negative with row-sums equal to one for each $t \in [0, \delta)$, meaning $P_t(X=x|Z=z)$ is a valid conditional probability distribution. Again we leave the full joint distribution unspecified, noting a QMD $P_t$ could be explicitly given by choosing the same conditional copula relating $X|Z$ and $W|Z$ for each $z$ at $P_0$.

    After simplifying terms, we arrive at
    \begin{align*}
        \varphi(P_t) =& \varphi(P_0) + \frac{1}{15t^2}\left\{\frac{9}{10}t + \frac{15}{4}t^3  -\frac{12}{5}t^5\right\}.
    \end{align*}
    The second term on the right-hand side equals $\frac{3}{50}t^{-1} + O(t)$, and hence $\varphi(P_t) \rightarrow \infty$ as $t\rightarrow 0$. \cref{fig:diverging-fig} plots the value of $\varphi(P_t)$ for $t \in [0, 0.1)$.
    \begin{figure}[htb]
        \centering
        \includegraphics[width=0.75\linewidth]{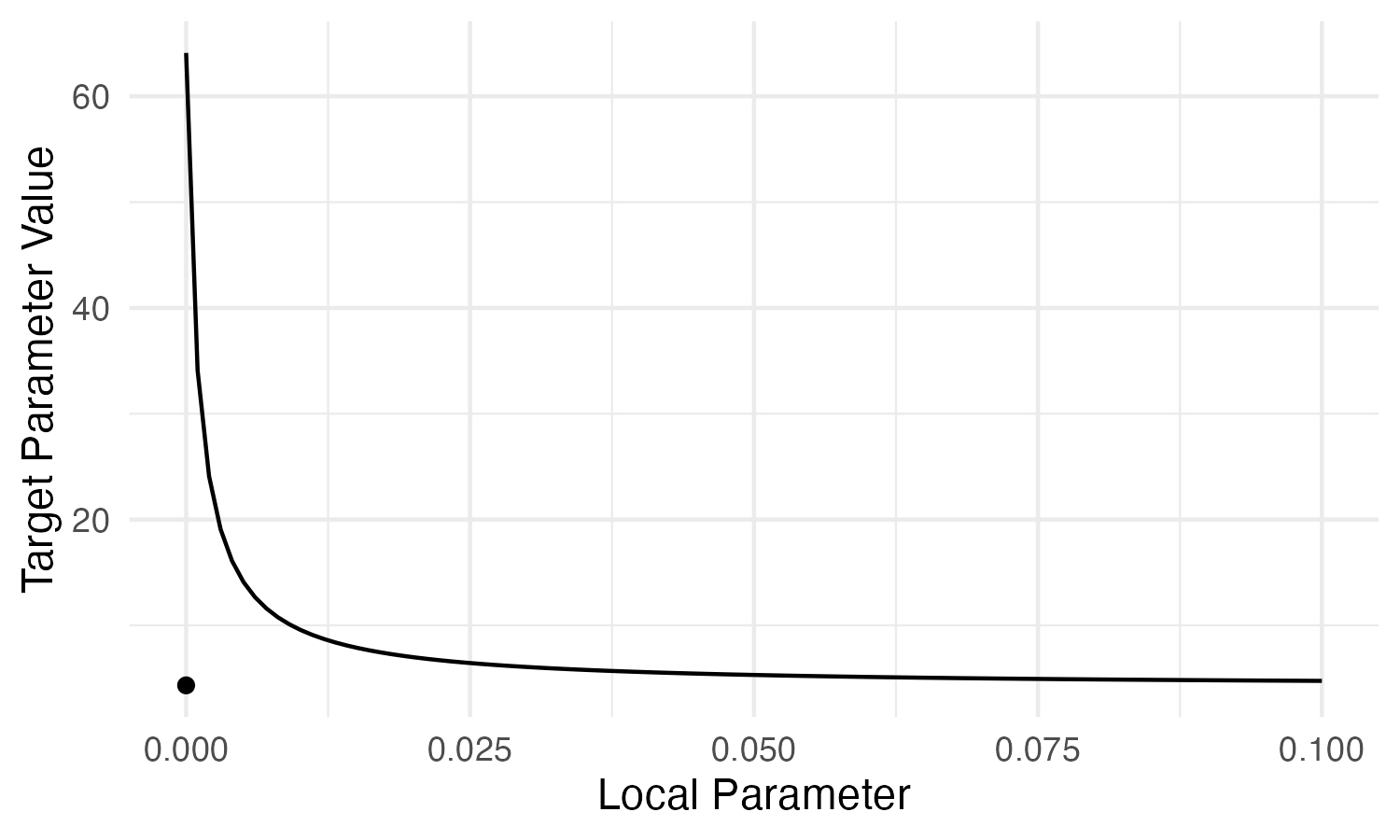}
        \caption[Failure of continuity of the target parameter in the discrete NPIV counterexample]{%
            Failure of continuity of the target parameter in the discrete NPIV counterexample.
            The curve traces the target parameter value $\varphi(P_t)$ against the local
            parameter $t \in [0, 0.1)$ along the path $\{P_t\}$ defined in
            \eqref{eq:counterexample-discontinuity}, and the isolated point marks the value
            $\varphi(P_0)$ at the limit law. Although $P_t \rightarrow P_0$ as $t \rightarrow 0$,
            $\varphi(P_t)$ diverges to infinity rather than approaching $\varphi(P_0)$, so
            $\varphi$ is discontinuous at $P_0$.}
        \label{fig:diverging-fig}
    \end{figure}
\end{example}

We now perform a small simulation study based on the above example. We use R \parencite{rcoreteamLanguageEnvironmentStatistical2026} and the R package SimEngine \parencite{kenny_simengine_2024} to run our simulations. We generate data at sample sizes $n = 200, 2000, 20000, 200000$. For each $n$ and each data generating process, we simulate 1000 replicates.

{ We consider sampling under three local data generating processes, each corresponding to a QMD at $P_0$ submodel $\{P_t: t\in [0, \delta)\} \in \mathscr{P}(P_0, \mc{P})$. Given one such submodel, we sample data according to the local alternative $P_{0.1/\sqrt{n}}$, and compute bias and coverage for the parameter $\varphi(P_{0.1/\sqrt{n}})$. We describe our full data generating process in \cref{sup-sec:simulation}, but give a brief explanation of the expected results under each process here.
\begin{enumerate}
    \item \textbf{The Fixed Alternative:} Define $P_{t} \coloneqq P_0$ for all $t \in [0, \delta)$. This is the standard i.i.d. sampling process. We expect our estimator to be consistent and display  nominal coverage for $n$ large enough.
    \item \textbf{The Irregular Alternative:} Let $\mc{E}_{P_t}$ be defined as in \eqref{eq:counterexample-discontinuity} with $\beta(t) = 0.01$ for all $t$, and let the joint distribution of $(W, Z)$ under $P_t$ be equal to the joint distribution under $P_0$. Algebra shows that the unique solution $h_{P_t}$ to $\mc{E}_{P_t}h_{P_t} = r_{P_t}$ equals $(207/5, 212/5, -288/5)^T$ for all $t > 0$, and is not equal to $h^\dagger_{P_0}$. Consequently, for $C$ sufficiently large, this submodel lies within $\mc{P}_{h_{P}^\dagger \text{ bound}}$. Because $\alpha_{P_t} = \psi$ and $r_{P_t} = r_{P_0}$ are fixed along the submodel, $\varphi$ is Hellinger continuous at $P_0$ along this submodel by \cref{thrm:bounded-means-continuous}; indeed, because $h_{P_t}$ is constant and $P_t(X)$ is linear in $t$, direct calculation gives $\varphi(P_t) = \varphi(P_0) + \frac{15621}{2500}t$. At the same time, both inverse problems admit multiple solutions at $P_0$, so that $\varphi$ is not pathwise differentiable at $P_0$ relative to $\mc{P}_{h_{P}^\dagger \text{ bound}}$ by \cref{thrm:not-injective-not-differentiable}. In particular, along this path, $\frac{d}{dt}\varphi(P_t)\vert_{t=0} \not= \langle s, \varphi^1_{P_0} \rangle_{L^2(P_0)} $ with $s$ the score of the submodel and $\varphi^1_{P_0}$ the influence function of the Tikhonov regularized estimator described below. By Le Cam's third lemma, as in the proof of \cref{lemma:failure-differentiability-irregular} (see \cref{subsubsec:implications-non-diff-main}), the Tikhonov regularized estimator should then display non-vanishing $\sqrt{n}$-scaled bias under sampling from $P_{0.1/\sqrt{n}}$, and the corresponding Wald-type confidence intervals should under-cover asymptotically, even though the estimator should be consistent.
    \item \textbf{The Discontinuous Alternative:}
    Let $\mc{E}_{P_t}$ be defined as in \eqref{eq:counterexample-discontinuity} with $\beta(t) = t^2$
    and let the joint distribution of $(W, Z)$ under $P_t$ be equal to the joint distribution under $P_0$. Then, by the previous discussion, $\varphi(P_t) \rightarrow \infty$ as $t \rightarrow 0$. Along this submodel, both minimum norm solutions $h^\dagger_{P_t}$ and $q^\dagger_{P_t}$ diverge in norm, so that it leaves every model with uniform norm bounds considered in \cref{sec:norm-bound-models}, and illustrates the discontinuity of \cref{lemma:discrete-necessary-sufficient-continuity} over the baseline model $\mc{P}$. Because the submodel is QMD at $P_0$, the laws $P_{0.1/\sqrt{n}}^n$ are contiguous to $P_0^n$, so that our estimator, which is consistent for $\varphi(P_0)$ under $P_0$, remains consistent for $\varphi(P_0)$ under $P_{0.1/\sqrt{n}}$ by Le Cam's first lemma. As $\varphi(P_{0.1/\sqrt{n}}) \rightarrow \infty$, we therefore expect our estimator to be inconsistent for $\varphi(P_{0.1/\sqrt{n}})$ and the coverage to converge to zero, illustrating the impossibility results of \cref{subsubsec:implications-discontinuity-main}.
\end{enumerate}
}

Given a random sample from $P$, we estimate the conditional expectation operators $\mc{E}_{P}$ and $\mc{E}_{P}^*$ and $r_{P}$ using non-parametric MLE, and recall that $\alpha_P = \psi$ is fixed. We then estimate the minimum norm solutions $h^\dagger_{P}$ and $q^\dagger_{P}$ using the Tikhonov regularized estimators $\widehat{h}_n, \widehat{q}_n$ with penalty $\lambda_n = 0.01n^{-0.6}$. This choice of penalty suffices to ensure $E_{P_0}[\{\widehat{h}_n(X) - h^\dagger_{P_0}(X)\}\{\widehat{q}_n(Z) - q^\dagger_{P_0}(Z)\}] = o_p(n^{-1/2})$ under repeated sampling from $P_0$, as shown in \cref{prop:tikhonov-rates}. We construct $\widehat{\varphi}_n$ using the Neyman orthogonal estimating equation \eqref{eq:neyman-orthogonal-ee}, as in \eqref{eq:dml-estimator}. We do not use cross-fitting as all parameters are finite dimensional.

It follows that under sampling from $P_0$, the estimator $\widehat{\varphi}_n$ is asymptotically linear with influence function $\varphi^1_{P_0} \coloneqq \varphi^1(O;h_{P_0}^\dagger, q_{P_0}^\dagger, \varphi(P_0))$ where $\varphi^1$ is defined in \eqref{eq:neyman-orthogonal-ee}. We use this influence function to construct Wald-type confidence intervals, where the asymptotic variance is estimated via the empirical second moment of the estimated influence function. Details of the estimation procedure can be found in \cref{sup-sec:simulation}. The estimator $\widehat{\varphi}_n$ and its associated confidence interval are computed identically under all three sampling processes; no aspect of their construction depends on which process generated the data.

\cref{tab:sim-results} reports the bias, the $\sqrt{n}$-scaled bias, and the coverage of the Wald-type confidence intervals for the three data generating processes. Throughout, bias and coverage are computed for the estimand $\varphi(P_{0.1/\sqrt{n}})$, which varies with both $n$ and the sampling process, and which we evaluate exactly rather than by simulation. For reference, $\varphi(P_0) = 13/3$, and under the fixed alternative the estimand equals $\varphi(P_0)$ at every $n$. With 1000 replicates, the Monte Carlo standard error of the coverage estimates is at most 1.6 percentage points. Under the fixed alternative, the estimator displays low bias and nominal coverage, as standard pointwise asymptotic theory predicts. Under the irregular alternative, the estimator is consistent but has a $\sqrt{n}$-scaled bias that is essentially constant across sample sizes, as predicted by \cref{lemma:failure-differentiability-irregular}. This reflects the fact that the limit distributions of irregular estimators depend on the score of the local alternative path. Here, this dependence prevents the limit distribution under the irregular alternative from having mean zero. The large $\sqrt{n}$-scaled bias in turn leads to poor coverage of the Wald-type confidence intervals under the irregular alternative. As expected from the contiguity argument above, under the discontinuous alternative, the bias, and hence the $\sqrt{n}$-scaled bias, diverges to infinity, and the coverage is 0\% for all sample sizes. The poor asymptotic behavior under the two local alternatives cannot be captured in a pointwise asymptotic analysis, as seen by the low bias and nominal coverage of the estimator under the fixed alternative data generating process.
\begin{table}[t]
\caption{Bias, scaled bias, 95\% confidence interval coverage, and mean confidence interval width of the Tikhonov-regularized IV estimator under each data-generating process.}
\begin{tabular*}{\linewidth}{@{\extracolsep{\fill}}rrrrr}
\toprule
\(n\) & Bias & \(\sqrt{n}\cdot\) Bias & Coverage & Mean CI Width \\
\midrule\addlinespace[2.5pt]
\multicolumn{5}{l}{Fixed Alternative} \\[2.5pt]
\midrule\addlinespace[2.5pt]
200 & -0.00511 & -0.0723 & 0.956 & 0.563 \\
2,000 & -0.00156 & -0.0699 & 0.959 & 0.163 \\
20,000 & 0.000449 & 0.0635 & 0.947 & 0.0490 \\
200,000 & -0.000310 & -0.139 & 0.948 & 0.0152 \\
\midrule\addlinespace[2.5pt]
\multicolumn{5}{l}{Irregular Alternative} \\[2.5pt]
\midrule\addlinespace[2.5pt]
200 & -0.217 & -3.07 & 0.646 & 0.553 \\
2,000 & -0.0688 & -3.08 & 0.618 & 0.162 \\
20,000 & -0.0221 & -3.13 & 0.567 & 0.0490 \\
200,000 & -0.00704 & -3.15 & 0.547 & 0.0152 \\
\midrule\addlinespace[2.5pt]
\multicolumn{5}{l}{Discontinuous Alternative} \\[2.5pt]
\midrule\addlinespace[2.5pt]
200 & -8.49 & -120 & 0.000 & 0.564 \\
2,000 & -26.8 & -1,200 & 0.000 & 0.161 \\
20,000 & -84.9 & -12,000 & 0.000 & 0.0490 \\
200,000 & -268 & -120,000 & 0.000 & 0.0152 \\
\bottomrule
\end{tabular*}
\label{tab:sim-results}
\end{table}

\section{Discussion}
\label{sec:discussion}

We have demonstrated that, under minimal conditions that suffice for the existence of a Neyman orthogonal estimating equation and further regularity, parameters defined as continuous linear functionals of solutions to ill-posed inverse problems are not Hellinger continuous at laws where both inverse problems are ill-posed, which is typical when $X$ and $Z$ are continuous. Locally uniform inference is impossible over this model, even in finite samples: for every $n$, no estimator attains any finite accuracy uniformly over Hellinger neighborhoods of $P_0$, and every confidence interval with nominal coverage locally uniformly has infinite expected diameter. Asymptotically, there do not exist locally uniformly consistent estimators, and all asymptotically locally uniformly honest confidence intervals have diameters that are not bounded in probability under $P_0$. These results do not conflict with recent work constructing pointwise asymptotically linear estimators, as pointwise inference is strictly weaker than locally uniform inference. { Under stronger assumptions, including a uniform bound over the model on the magnitude of the minimum norm solution to one of the inverse problems, Hellinger continuity is restored. However, the parameter may still fail to be pathwise differentiable over this model. Without additional assumptions, and under regularity conditions, pathwise differentiability fails at laws $P_0$ such that the inverse problem with uniformly bounded solutions does not admit a unique solution when both inverse problems are ill-posed. Consequently, at such laws, all asymptotically linear estimators are irregular, and exhibit non-vanishing $\sqrt{n}$-scaled bias and asymptotic undercoverage of Wald-type intervals under some local alternative data generating processes. On the other hand, if this solution is unique, and further regularity holds, the parameter is pathwise differentiable. Pathwise differentiability also holds under the stronger assumption that the $\beta$-source representer of one of the two solutions admits a bound uniformly over the model, whether or not the corresponding inverse problem admits a unique solution at $P_0$.

Whether a uniform norm bound on the minimum norm solution to one of the inverse problems is available depends on the substantive problem. In the shadow variable model, a uniform bound on the dual solution follows from a uniform bound on the density ratio between the outcome distribution in the nonrespondents and the respondents, and so functions as a positivity-like overlap condition. In the correctly specified additive NPIV model it follows from the structural interpretation of the primal solution: the causal dose response curve solves the primal inverse problem, so uniformly bounded counterfactual outcomes bound the minimum norm solution directly. Under mild violations of the exclusion restriction this argument breaks down, and the assumption appears to be a technical stability condition not justified based on any counterfactual interpretation.  In proximal causal inference the situation is unresolved. A correctly specified proximal model with bounded outcomes does bound the target parameter, so our discontinuity constructions must violate that model, but whether the model restricts the bridge functions enough to restore Hellinger continuity remains an open question. We are unaware of any general substantive assumption that leads to a uniform norm bound on the $\beta$-source representer.

We have also investigated the maximal tangent space of various statistical models. We have shown that the assumptions encoded in the baseline model $\mc{P}$ do not restrict the maximal tangent space of the statistical model when the dimensions of the spaces $L^2(P_{X,D})$ and $L^2(P_{Z,D})$ are equal. When the inverse problems are finite dimensional and these dimensions differ, the existence of solutions to both equations is potentially substantive whenever the operators defining them are full rank, a condition that fails only on a set of laws of  measure zero. In contrast to the baseline model, under the conditions that suffice for pathwise differentiability, we are able to derive a superset of the tangent space, proving it is no larger than the space derived by \cite{severiniEfficiencyBoundsEstimating2012} and \cite{chen_overidentification_2018-1} under stronger assumptions. When the primal solution is uniformly norm bounded and admits a unique solution at $P_0$, or when its $\beta$-source representer is uniformly norm bounded, this bound is substantive under mild conditions when there exist multiple solutions to the dual inverse problem, in that the tangent space is strictly smaller than that of a locally non-parametric model.
 }

Our results mirror in many respects the classical results on weak instruments in the linear instrumental variables problem. The results under the baseline model resemble impossibility results in the unrestricted linear IV problem \parencite{gleser_nonexistence_1987, dufour_impossibility_1997}, the results under a uniform norm bound on the primal solution when that solution is not unique at $P_0$ are related to the framing of weak-instrument asymptotics \parencite{staiger_instrumental_1997}, and the uniform bound on the norm of the $\beta$-source representer can be viewed as a type of strong instrument assumption. Given this relation, it may be fruitful to adapt results for performing valid inference in the linear IV problem \parencite{andrews_weak_2019} to the present setting, a task we leave to future work.

These results also separate two notions of differentiability that recent work has sought to unify. Neyman-orthogonality of the estimating equation \eqref{eq:neyman-orthogonal-ee} concerns differentiability with respect to the $L^2$ topology on the nuisance functions, while pathwise differentiability concerns differentiability along Hellinger differentiable paths through the model. \cite{chen_equivalence_2026} show that the two coincide under conditions requiring the nuisance functions to vary smoothly along quadratic mean differentiable submodels. The parameters studied here satisfy the former without satisfying the latter, precisely because that smoothness fails. The practical consequence is visible in \cref{sec:simulation}: Wald intervals constructed from the influence function of \eqref{eq:neyman-orthogonal-ee} undercover under local alternatives, potentially quite severely, as predicted by \cref{lemma:failure-differentiability-irregular}.

\section*{Acknowledgements}

Andrea Rotnitzky's work has been supported by the National Institute of
Allergy and Infectious Diseases grants UM1 AI068635 and R37-AI029168. This material is based upon work supported by the National Science Foundation Graduate Research Fellowship Program under Grant No. GE-2140004. The content is solely the responsibility of the authors. Any opinions, findings, and conclusions or
recommendations expressed in this material are those of the authors and do not necessarily reflect the views
of the National Science Foundation.

\section*{Use of AI Tools}

Claude (Anthropic) and ChatGPT (OpenAI) were used throughout the preparation of this manuscript and its supplementary material for editing and review. In addition, Claude Opus 4.8 developed the main argument underlying the proof of \cref{lemma:uniform-convergence-of-kernel}, and Claude Opus 5.5 developed the main argument underlying the proof of strong convergence of $h_{P}^\dagger \rightarrow h_{P_0}^\dagger$ within the model $\mc{P}_{\beta_h \text{ bound}}$ in \cref{lemma:weak-convergence-of-solutions}; those arguments were subsequently rewritten by the authors and checked thoroughly for correctness. The authors have verified all mathematical content and take full responsibility for the results presented here.

\printbibliography[heading=subbibintoc, title={References}]
\end{refsection}

\appendix

\crefalias{section}{appendix}
\crefalias{subsection}{subappendix}
\crefalias{subsubsection}{subsubappendix}

\supplementnumbering{section}
\supplementnumbering{equation}
\supplementnumbering{theorem}
\supplementnumbering{corollary}
\supplementnumbering{lemma}
\supplementnumbering{proposition}
\supplementnumbering{remark}
\supplementnumbering{assumption}
\supplementnumbering{definition}
\supplementnumbering{example}
\supplementnumbering{examplerepeat}
\supplementnumbering{examplerepeattwo}
\supplementnumbering{counterexample}
\supplementnumbering{figure}
\supplementnumbering{table}
\supplementnumbering{algorithm}

\clearpage

\appendixpage

\etocdepthtag.toc{mtsupplement}
\etocsettagdepth{mtmain}{none}
\etocsettagdepth{mtsupplement}{subsection}
\etocsettocstyle{\section*{Contents of the Supplementary Material}}{}
{\singlespacing\tableofcontents}

\begin{refsection}

To distinguish the results in the supplement from the main text, we preface section numbers, lemma numbers, equation numbers, etc. that appear for the first time in the supplement by ``S.''.

\section{Parameter Smoothness and Implications for Inference}
\label{sup-sec:smoothness-implications}
\subsection{Implications of Hellinger Discontinuity}
\label{subsubsec:implications-discontinuity-supplement}

Here, we give an extended discussion of the results in \cref{subsubsec:implications-discontinuity-main}. As described there, the finite-sample impossibility results that we shall derive below are in the same spirit of classical results on
nonexistence of "honest" inference about locally unstable
functionals.  The earliest example is due to \cite{bahadur_nonexistence_1956}, who
showed that, over large nonparametric classes, effective distribution-free
inference for the mean is impossible.  Subsequent work developed this idea in
several directions: \cite{donoho_one-sided_1988} and
\cite{donoho_geometrizing_1991} related nonparametric inference to moduli of
continuity, defined in \cref{subsubsec:implications-discontinuity-main}; \cite{gleser_nonexistence_1987} showed that uniform confidence sets in
errors-in-variables and related models cannot have finite expected diameter;
\cite{low_nonparametric_1997} and \cite{cai_adaptation_2004} gave lower bounds on
the expected length of honest nonparametric confidence intervals in terms of
moduli of continuity; and \cite{pfanzagl_nonexistence_1998} provided a general
formulation of nonexistence results for confidence procedures for
discontinuous functionals.  Closely related impossibility results also appear
in econometrics in the work of \cite{dufour_impossibility_1997} and \cite{potscher_lower_2002}.

We now prove \cref{lemma:finite-sample-discontinuity}, using the local modulus of continuity $\omega_{P_0}$, the notation $P^n$, $C_n$, and $|C_n|$, and the radius $r_{n,\tau}$ defined in \eqref{eq:rn-tau}, all introduced in \cref{subsubsec:implications-discontinuity-main}. We additionally let \(d_{\mathrm{TV}}\) denote the usual total variation distance
\(d_{\mathrm{TV}}(P,Q):=\sup_A |P(A)-Q(A)|\). We build on the following finite-sample bounds, which follow from elementary two-point testing arguments.

\begin{proposition}[A finite-sample two-point lower bound]
\label{prop:two-point-ci-bound}
Fix two laws \(P_0\) and \(P_1\) in  $\mathcal{M}$. Write
\(\Delta:=|\phi(P_1)-\phi(P_0)|\).  Suppose that \(C_n\) has coverage at least
\(1-\alpha\) at both \(P_0\) and \(P_1\), that is,
\begin{equation*}
    P_0^n\{\phi(P_0)\in C_n\}\ge 1-\alpha,
    \qquad
    P_1^n\{\phi(P_1)\in C_n\}\ge 1-\alpha.
\end{equation*}
Then
\begin{equation}
\label{eq:two-point-ci-tail}
    P_0^n\{|C_n|\ge \Delta\}
    \ge
    1-2\alpha-d_{\mathrm{TV}}(P_0^n,P_1^n),
\end{equation}
and consequently
\begin{equation}
\label{eq:two-point-ci-expected-length}
    E_{P_0^n}\{|C_n|\}
    \ge
    \Delta
    \left[
        1-2\alpha-d_{\mathrm{TV}}(P_0^n,P_1^n)
    \right]_+ .
\end{equation}
\end{proposition}

The following proposition as well as results that we state subsequently invoke the following well known bound between total
variation and Hellinger distance of product measures \parencite{cam_basic_1986},
\begin{equation}
\label{eq:tv-product-hellinger-bound}
    d_{\mathrm{TV}}(P^n,Q^n)
    \le
    \left[
        1-\{1-H^2(P,Q)\}^{2n}
    \right]^{1/2}.
\end{equation}

\begin{proposition}[Finite-sample expected length and the local modulus]
\label{prop:local-modulus-ci-bound}
Fix \(r\in(0,1)\), and suppose that \(C_n\) has coverage at least
\(1-\alpha\) over the Hellinger ball
\(\{P\in\mathcal M:H(P,P_0)<r\}\).  Then
\begin{equation}
\label{eq:local-modulus-expected-length}
    E_{P_0^n}\{|C_n|\}
    \ge
    \omega_{P_0}(r)
    \left[
        1-2\alpha-
        \left\{1-(1-r^2)^{2n}\right\}^{1/2}
    \right]_+ .
\end{equation}
\end{proposition}

The same two-point logic also yields finite-sample lower bounds for point
estimators.  Let \(\widehat\phi_n\) be any estimator of \(\phi(P)\).  For fixed
\(P_0,P_1\in\mathcal M\), write again
\(\Delta:=|\phi(P_1)-\phi(P_0)|\).  Then we have that
\begin{equation}
\label{eq:two-point-estimator-lower-bound}
    \max_{i=0,1}
    P_i^n
    \left\{
        |\widehat\phi_n-\phi(P_i)|\ge \frac{\Delta}{2}
    \right\}
    \ge
    \frac{1-d_{\mathrm{TV}}(P_0^n,P_1^n)}{2}.
\end{equation}
To show \eqref{eq:two-point-estimator-lower-bound}, let \(A:=\{|\widehat\phi_n-\phi(P_0)|<\Delta/2\}\).  On \(A\), the
estimator is closer to \(\phi(P_0)\) than to \(\phi(P_1)\), so
the event \(A\) is contained in
\(\{|\widehat\phi_n-\phi(P_1)|\ge\Delta/2\}\).  Therefore
\begin{align*}
    &P_0^n
    \left\{
        |\widehat\phi_n-\phi(P_0)|\ge \frac{\Delta}{2}
    \right\}
    +
    P_1^n
    \left\{
        |\widehat\phi_n-\phi(P_1)|\ge \frac{\Delta}{2}
    \right\} \\
    &\qquad\ge
    P_0^n(A^c)+P_1^n(A)
    \ge
    1-d_{\mathrm{TV}}(P_0^n,P_1^n),
\end{align*}
which implies \eqref{eq:two-point-estimator-lower-bound}.

\begin{proof}[Proof of \cref{lemma:finite-sample-discontinuity}]
Because $\phi$ is Hellinger discontinuous at $P_0$, there exist $\epsilon > 0$ and, for every $r > 0$, a law $P \in \mc{M}$ with $H(P, P_0) < r$ and $|\phi(P) - \phi(P_0)| > \epsilon$, so that $\omega_{P_0}(r) > \epsilon$ for every $r > 0$. By \eqref{eq:tv-product-hellinger-bound} and the definition of $r_{n,\tau}$ in \eqref{eq:rn-tau}, every $P_1 \in \mc{M}$ with $H(P_1, P_0) < r_{n,\tau}$ satisfies
\begin{equation*}
    d_{\mathrm{TV}}(P_0^n, P_1^n) \leq \left\{1 - (1 - r_{n,\tau}^2)^{2n}\right\}^{1/2} = \tau.
\end{equation*}
Fix such a $P_1$ with $\Delta \coloneqq |\phi(P_1) - \phi(P_0)| > \epsilon$. If $C_n$ has coverage at least $1 - \alpha$ over $\{P \in \mc{M}: H(P, P_0) < r_{n,\tau}\}$, then it has coverage at least $1-\alpha$ at $P_0$ and $P_1$, and \cref{prop:two-point-ci-bound} gives $P_0^n\{|C_n| \geq \epsilon\} \geq P_0^n\{|C_n| \geq \Delta\} \geq 1 - 2\alpha - \tau$, proving \eqref{eq:tail-any-M}. Similarly, for any estimator $\widehat{\phi}_n$, \eqref{eq:two-point-estimator-lower-bound} gives
\begin{equation*}
    \max_{i = 0, 1} P_i^n\{|\widehat\phi_n-\phi(P_i)|\ge \epsilon/2\} \geq \max_{i = 0, 1} P_i^n\{|\widehat\phi_n-\phi(P_i)|\ge \Delta/2\} \geq \frac{1-\tau}{2},
\end{equation*}
proving \eqref{eq:finite-sample-estimator-bad-somewhere}, as both $P_0$ and $P_1$ belong to the ball of radius $r_{n,\tau}$.

If $\omega_{P_0}(r) = \infty$ for every $r > 0$, then for every $M < \infty$ we may instead choose $P_1$ with $|\phi(P_1) - \phi(P_0)| > M$. The preceding arguments then give $P_0^n\{|C_n| \geq M\} \geq 1 - 2\alpha - \tau$ for every $M < \infty$, and \eqref{eq:finite-sample-estimator-bad-somewhere} with $\epsilon/2$ replaced by $M/2$, for every $M$. By continuity from above,
\begin{equation*}
    P_0^n\{|C_n|=\infty\}
    =
    \lim_{M\to\infty}P_0^n\{|C_n|\ge M\}
    \ge
    1-2\alpha-\tau,
\end{equation*}
and the infinite-expectation claim follows immediately because $1 - 2\alpha - \tau > 0$.
\end{proof}

The asymptotic statements \eqref{eq:no-local-uniform-estimation} and \eqref{eq:no-local-uniform-ci-tightness} in \cref{subsubsec:implications-discontinuity-main} follow because \(\omega_{P_0}(r)> \epsilon\) for every \(r>0\), so that the proof of \cref{lemma:finite-sample-discontinuity} applies verbatim with the ball of radius \(r_{n,\tau}\) replaced by the smaller ball of radius \(\min\{\delta_n,r_{n,\tau}\}\). For \eqref{eq:no-local-uniform-ci-tightness}, we additionally use that, under \eqref{eq:local-uniform-asymptotic-coverage}, for every \(\eta > 0\) with \(\alpha + \eta < 1/2\) and \(\tau < 1 - 2(\alpha + \eta)\), \(C_n\) has coverage at least \(1 - \alpha - \eta\) over \(\mathcal U_n(P_0)\) for all \(n\) sufficiently large, and then let \(\eta \downarrow 0\).

\subsection{Implications of Failure of Pathwise Differentiability}
\label{subsubsec:implications-non-diff-supplement}

Here we give a more in depth discussion of \cref{subsubsec:implications-non-diff-main}.
Pathwise differentiability, formally defined in \cref{subsec:semiparametrics}, is a strictly stronger smoothness condition than Hellinger continuity.
\begin{remark}
    \label{remark:failure-of-continuity-fails-diff}
    To see that failure of Hellinger continuity implies failure of pathwise differentiability, note that failure of Hellinger continuity implies the existence of a sequence $\{P_j\}_{j=1}^\infty \subseteq \mc{M}$ and constant $\epsilon > 0$ for which $H(P_j, P_0) \rightarrow 0$ but where $|\phi(P_j) - \phi(P_0)| > \epsilon$ for each $j$. Passing to a subsequence $\{j_i\}_{i=1}^\infty$ with $H(P_{j_i}, P_0) \leq i^{-2}$, we have $H(P_{j_{\lfloor 1/t\rfloor}}, P_0) \leq \lfloor 1/t\rfloor^{-2} = o(t)$ as $t \downarrow 0$. Then, $\{P_{j_{\lfloor 1/t\rfloor}}:t\in [0, 1)\} \subseteq \mc{M}$, where we take the convention that $P_{j_{\lfloor 1/0\rfloor}} = P_0$, is QMD at $P_0$ with score $0$, but $|\phi(P_{j_{\lfloor 1/t\rfloor}}) - \phi(P_0)|/t > \epsilon/t$ diverges to infinity as $t \downarrow 0$, so that $\phi$ is not pathwise differentiable at $P_0$. The converse fails in general; see \cref{subsec:failure-differentiability} for laws and models at which $\varphi$ is Hellinger continuous but not pathwise differentiable.
\end{remark}

We now prove \cref{lemma:failure-differentiability-irregular}, using the notation and setting of \cref{subsubsec:implications-non-diff-main}.

\begin{proof}[Proof of \cref{lemma:failure-differentiability-irregular}]
We first analyze the behavior of $\widehat{\phi}_n$ in case 1). Define the map $D: \Theta \mapsto \mathbb{R}$ via
\begin{align*}
    D(\theta) \coloneqq \left.\frac{d}{dt}\phi(P_{t, \theta}) \right\vert_{t=0}.
\end{align*}
Then, for each $\theta \in \Theta$, application of Le Cam's third Lemma shows that
\begin{align*}
    \sqrt{n}(\widehat{\phi}_n - \phi(P_{1/\sqrt{n}, \theta})) \overset{P_{1/\sqrt{n}, \theta}}{\rightsquigarrow} N\left(\Delta_{\theta}, \varsigma^2_{P_0}\right)
\end{align*}
where
\begin{align*}
    \Delta_{\theta} \coloneqq \langle s_{\theta}, \phi^1_{P_0}\rangle_{L^2(P_0)} - D(\theta).
\end{align*}
Suppose that $\Delta_\theta = 0$ for every $\theta \in \Theta$. Then $D(\theta) = \langle s_{\theta}, \phi^1_{P_0}\rangle_{L^2(P_0)}$ for every $\theta \in \Theta$, so that the bounded linear functional $s \mapsto \langle s, \phi^1_{P_0}\rangle_{L^2(P_0)}$, restricted to the tangent space of $\mathscr{P}_\Theta$, is a local parameter, a contradiction. Hence, there exists $\theta \in \Theta$ for which $\Delta_{\theta} \not=0$, and $\widehat{\phi}_n$ is locally asymptotically biased with local asymptotic bias $\Delta_\theta$ under sampling from $P_{1/\sqrt{n}, \theta}$.

For case 2), fix $\theta \in \Theta$ for which $\frac{d}{dt}\phi(P_{t, \theta})\vert_{t=0}$ does not exist. Then, as $t \downarrow 0$, the difference quotient $\{\phi(P_{t, \theta}) - \phi(P_0)\}/t$ either has at least two cluster points in the extended reals $\overline{\mathbb{R}}$, or has $+\infty$ or $-\infty$ as its only cluster point. In either case, application of Bolzano-Weierstrass in $\overline{\mathbb{R}}$ implies there exist $t_k \downarrow 0$ and $c \in \overline{\mathbb{R}}$ with $c \not = \langle s_{\theta}, \phi^1_{P_0}\rangle_{L^2(P_0)}$ such that
\begin{align*}
    \lim_{k\rightarrow \infty} \frac{\phi(P_{t_k, \theta}) - \phi(P_0)}{t_k} = c.
\end{align*}
Let $n_k \coloneqq \lceil t_k^{-2}\rceil$, passing to a further subsequence so that $n_k$ is strictly increasing, and let $h_k \coloneqq \sqrt{n_k}t_k$, so that $h_k \rightarrow 1$ and $t_k = h_k/\sqrt{n_k}$. Application of Le Cam's third lemma along local alternatives $h_k/\sqrt{n_k}$ with $h_k \rightarrow 1$ (\cite{vaart_asymptotic_1998-2} Chapters 6 and 7) gives, when $c$ is finite,
\begin{align*}
    \sqrt{n_k}(\widehat{\phi}_{n_k} - \phi(P_{h_k/\sqrt{n_k}, \theta})) \overset{P_{h_k/\sqrt{n_k}, \theta}}{\rightsquigarrow} N\left(\Delta_{\theta}, \varsigma^2_{P_0}\right)
\end{align*}
where $\Delta_{\theta} \coloneqq \langle s_{\theta}, \phi^1_{P_0}\rangle_{L^2(P_0)} - c \not= 0$. When $c = \pm\infty$, $\sqrt{n_k}(\widehat{\phi}_{n_k} - \phi(P_{h_k/\sqrt{n_k}, \theta}))$ instead diverges to $\mp\infty$ in probability under $P_{h_k/\sqrt{n_k}, \theta}$. In case 1), the claimed sequences are $n_k = k$ and $h_k = 1$. The coverage claim follows directly from the limiting distribution.
\end{proof}

\section{Further Details for \cref{sec:simulation}}
\label{sup-sec:simulation}

\subsection{Simulation Details}

We first give a full description of our simulation. In what follows, we use matrix and vector notation to describe functions and conditional and marginal probabilities. For all settings, we fix the following constants:
\begin{align*}
    \mc{E}_{P_0} \coloneqq& \begin{pmatrix}
        .4 & .2 & .4\\
        .2 & .4 & .4 \\
        .2 & .4 & .4
    \end{pmatrix},\\
    P_0(W=w|Z=z) \coloneqq& P_0(X=w|Z=z) = \mc{E}_{P_0;z, w} \text{ for } w, z \in \{1, 2, 3\},\\
    P_{Z;0} \coloneqq & \begin{pmatrix}
        1/3 &
        1/3 &
        1/3
    \end{pmatrix}^T,\\
    \psi \coloneqq& P_{Z|X;0} \begin{pmatrix}
        1 &
        2 &
        3
    \end{pmatrix}^T,\\
    \tau \coloneqq& \begin{pmatrix}
        1 &
        2 &
        3
    \end{pmatrix}^T,\\
    \rho(z) \coloneqq& I(z=1)0.5 + I(z=2)0.4 + I(z=3)0.6.
\end{align*}
Noting that $P_{Z|X;0} = \mc{E}_{P_0}^*$, we see that $\psi \in \mathscr{R}(\mc{E}_{P_0}^*)$. Furthermore, because $P_{W|Z;0} = P_{X|Z;0}$, $P_{W|Z;0}\tau \in \mathscr{R}(\mc{E}_{P_0})$.

Now, let $C_z:[0, 1]\times[0, 1] \rightarrow [0, 1]$ for $z \in \{1, 2, 3\}$ denote a copula function implicitly defined by the sampling procedure we define below. Then, we let $P_0$ denote the joint distribution for $(W, X, Z)$ with margin for $Z$ given above and joint CDF for $(W, X)|Z=z$ equal to
\begin{align*}
    P_0(W \leq w, X \leq x|Z=z) =& C_z(P_0(W \leq w|Z=z), P_0(X \leq x|Z=z)).
\end{align*}

Recall that we simulate data under three data generating processes, each corresponding to a submodel of $\mc{P}$. In the fixed alternative submodel, which we denote with $\{P^{Fixed}_{t}:t \in [0, \delta)\}$, we let $P^{Fixed}_{t} \coloneqq P_0$ for each $t$. For the irregular alternative submodel, which we denote with $\{P^{Irreg}_{t}: t\in[0, \delta)\}$, we let
\begin{align*}
    \mc{E}_{P^{Irreg}_t} \coloneqq \mc{E}_{P_0} + t\begin{pmatrix}
            0 & 0 & 0\\
            -1 & 0.99 & 0.01 \\
            0 & 0 & 0
        \end{pmatrix}
\end{align*}
so that $P_t^{Irreg}(X=x|Z=z)$ equals the $(z, x)$ entry of $\mc{E}_{P^{Irreg}_t}$. Then, we let $P_t^{Irreg}$ be the joint distribution with marginal for $Z$ equal to $P_0$ and conditional CDF of $(W, X)$ given $Z=z$ given by
\begin{align*}
    P_t^{Irreg}(W \leq w, X \leq x|Z=z) = C_z(P_{0}(W \leq w|Z=z), P_t^{Irreg}(X \leq x|Z=z)).
\end{align*}
Finally, for the discontinuous submodel, which we denote with $\{P_t^{Discont}: t\in[0, \delta)\}$, we let
\begin{align*}
    \mc{E}_{P^{Discont}_t} \coloneqq \mc{E}_{P_0} + t\begin{pmatrix}
            0 & 0 & 0\\
            -1 & 1 - t^2 & t^2 \\
            0 & 0 & 0
        \end{pmatrix}
\end{align*}
so that $P_t^{Discont}(X=x|Z=z)$ equals the $(z, x)$ entry of $\mc{E}_{P^{Discont}_t}$. Now, we let $P_t^{Discont}$ be the joint distribution with marginal for $Z$ equal to $P_0$ and conditional CDF of $(W, X)$ given $Z=z$ given by
\begin{align*}
    P_t^{Discont}(W \leq w, X \leq x|Z=z) = C_z(P_{0}(W \leq w|Z=z), P_t^{Discont}(X \leq x|Z=z)).
\end{align*}
Since $C_z$ is fixed and infinitely differentiable on $(0, 1)^2$, and the conditional margins are polynomial in $t$, the cell probabilities of $(W, X)|Z=z$ are continuously differentiable in $t$ for each $z$. Furthermore, because $|\rho(z)| < 1$, the Gaussian copula assigns positive mass to every cell, so these cell probabilities are bounded away from zero for $t$ small enough. Hence, $\{P^{Irreg}_t: t\in[0, \delta)\}$ and $\{P^{Discont}_t: t\in[0, \delta)\}$ are both QMD submodels of $\mc{P}$ at $P_0$.

Our choice of copula $C_z$ is implicitly given in steps $2$ and $3$ of the sampling procedure we now outline. For generic (conditional) CDFs $P_Z$ of $Z$, $P_{X|Z}$ of $X|Z$, and $P_{W|Z}$ of $W|Z$ generate an observation as follows:
\begin{enumerate}
    \item Sample $Z \sim P_{Z}$.
    \item Sample
    \begin{align*}
        W', X' \sim N\left(
        \begin{pmatrix}
            0 \\
            0
        \end{pmatrix},
        \begin{pmatrix}
            1 & \rho(Z) \\
            \rho(Z) & 1
        \end{pmatrix}
        \right)
    \end{align*}
    \item Let
    \begin{align*}
        W'' \coloneqq \Phi(W'), X'' \coloneqq \Phi(X')
    \end{align*}
    for $\Phi$ the CDF of the standard normal distribution, and set
    \begin{align*}
        W \coloneqq P^{-1}_{{W|Z}}(W''|Z), && X \coloneqq P^{-1}_{{X|Z}}(X''|Z)
    \end{align*}
    where $P^{-1}_{{W|Z}}(w|z)$ is the quantile function (generalized inverse CDF) of the conditional distribution of $W$ given $Z=z$ under $P$, and similarly for $P^{-1}_{{X|Z}}(x|z)$.
\end{enumerate}
Steps 2 and 3 generate $W, X$ that are correlated given $Z$ using a conditional Gaussian copula, discretized to respect the supports of $W, X$. The parameter $\rho(Z)$ determines the correlation between $W$ and $X$ for each $Z$.

For $n = 200, 2000, 20000, 200000$, we generate $n$ i.i.d. observations from each of the three data generating processes described above according to the sampling procedure described above. For each $n$, we sample 1000 replicates from $P^{Fixed}_{0.1/\sqrt{n}},P^{Irreg}_{0.1/\sqrt{n}}$, and $P^{Discont}_{0.1/\sqrt{n}}$. We note that $0.1/\sqrt{200}  < 0.2$, so that each of $P^{Fixed}_{0.1/\sqrt{n}}, P^{Irreg}_{0.1/\sqrt{n}}, P^{Discont}_{0.1/\sqrt{n}}$ are valid distribution functions with cell probabilities bounded above and below by constants.

For each replicate, we estimate the conditional expectation operators $\mc{E}_{P}$ and $\mc{E}_{P}^*$ by their non-parametric MLEs $\widehat{\mc{E}}_n$ and $\widehat{\mc{E}}_n^*$, which are simply the empirical conditional probability matrices associated with the operators. We estimate $r_P$ by
\begin{align*}
    \widehat{r}_n \coloneqq \widehat{E}_{W|Z} \tau,
\end{align*}
with $\widehat{E}_{W|Z}$ the empirical conditional probability matrix of $W$ given $Z$, and $\alpha_P$ by $\widehat{\alpha}_n \coloneqq \psi$.

We then estimate $h_P$ and $q_P$ by Tikhonov regularization \parencite{carrasco_chapter_2007}. Letting $\lambda_n \coloneqq 0.01 n^{-0.6}$, we define
\begin{align*}
    \widehat{h}_{n} \coloneqq (\widehat{\mc{E}}_{n}^*\widehat{\mc{E}}_{n}+ \lambda_n I)^{-1}\widehat{\mc{E}}_{n}^*\widehat{r}_n, && \widehat{q}_{n} \coloneqq (\widehat{\mc{E}}_{n}\widehat{\mc{E}}_{n}^* + \lambda_n I)^{-1}\widehat{\mc{E}}_{n}\widehat{\alpha}_n.
\end{align*}
By our choice of $\lambda_n$, \cref{prop:tikhonov-rates} below gives that
\begin{align*}
        \|\widehat{h}_n - h_{P_0}^\dagger\|_{L^2(P_{X;0})} =& O_p\left(n^{-2/5}\right), &  \|\widehat{q}_n - q_{P_0}^\dagger\|_{L^2(P_{Z;0})} =& O_p\left(n^{-2/5}\right)\\
        \|\mc{E}_{P_0}\{\widehat{h}_n - h_{P_0}^\dagger\}\|_{L^2(P_{Z;0})} =&O_p\left(n^{-1/2}\right) , & \|\mc{E}_{P_0}^*\{\widehat{q}_n - q_{P_0}^\dagger\}\|_{L^2(P_{X;0})} =&O_p\left(n^{-1/2}\right).
    \end{align*}
Consequently,
\begin{align*}
    \left|E_{P_0}\left[\{\widehat{h}_{n}(X) - h^\dagger_{P_0}(X)\}\{\widehat{q}_{n}(Z) - q^\dagger_{P_0}(Z)\}\right]\right| \leq& \|\mc{E}_{P_0}\{\widehat{h}_{n}(X) - h^\dagger_{P_0}(X)\}\|_{L^2(P_{Z;0})}\\
    &\times\|\widehat{q}_{n}(Z) - q^\dagger_{P_0}(Z)\|_{L^2(P_{Z;0})}\\
    =&  o_p(n^{-1/2}).
\end{align*}
We do not claim our choice of $\lambda_n$ is optimal, merely that it suffices for demonstration.

Then, standard arguments \parencite{bennettSourceConditionDouble2023b, olivas-martinez_source-condition_2025} show that
\begin{align*}
    \widehat{\varphi}_n \coloneqq \frac{1}{n}\sum_{i=1}^n \varphi^1(O_i; \widehat{h}_n, \widehat{q}_n, 0),
\end{align*}
with $\varphi^1$ as in \eqref{eq:neyman-orthogonal-ee}, is asymptotically linear under repeated sampling from $P_0$ with influence function equal to $\varphi^1_{P_0} \coloneqq \varphi^1(O;h_{P_0}^\dagger, q_{P_0}^\dagger, \varphi(P_0))$. Because $W, X, Z$ are each finite dimensional, cross-fitting is not needed. Consequently,
\begin{align*}
    \sqrt{n}(\widehat{\varphi}_n - \varphi(P_0)) \rightsquigarrow N(0, \varsigma_{P_0}^2), && \varsigma_{P_0}^2 \coloneqq E_{P_0}[\{\varphi^1_{P_0}(O)\}^2].
\end{align*}
We estimate $\varsigma_{P_0}^2$ with
\begin{align*}
    \widehat{\varsigma}^2_n \coloneqq \frac{1}{n}\sum_{i=1}^n \{\varphi^1(O_i; \widehat{h}_n, \widehat{q}_n, \widehat{\varphi}_n)\}^2,
\end{align*}
and construct Wald-type $1-\alpha$ confidence intervals with
\begin{align*}
    C_n \coloneqq (\widehat{\varphi}_n - \Phi^{-1}(1-\alpha/2)\widehat{\varsigma}_n/\sqrt{n}, \widehat{\varphi}_n + \Phi^{-1}(1-\alpha/2)\widehat{\varsigma}_n/\sqrt{n}).
\end{align*}

The following proposition gives rates of convergence for our estimators of $h_{P_0}^\dagger, q_{P_0}^\dagger$. The proof uses standard arguments, but we repeat them here to keep the results self contained.
\begin{proposition}
    \label{prop:tikhonov-rates}
    For $\widehat{h}_n$, $\widehat{q}_n$ estimated as described above, under sampling from $P_0$, and with $\delta_n = O_p(n^{-1/2})$ defined in \eqref{eq:delta-n-def-tikhonov} below,
    \begin{align*}
        \|\widehat{h}_n - h_{P_0}^\dagger\|_{L^2(P_{X;0})} =& O_p\left(\delta_n + \lambda_n + \frac{\delta_n^2}{\lambda_n}\right), &  \|\widehat{q}_n - q_{P_0}^\dagger\|_{L^2(P_{Z;0})} =& O_p\left(\delta_n + \lambda_n + \frac{\delta_n^2}{\lambda_n}\right)\\
        \|\mc{E}_{P_0}\{\widehat{h}_n - h_{P_0}^\dagger\}\|_{L^2(P_{Z;0})} =&O_p\left(\delta_n + \lambda_n + \frac{\delta_n^3}{\lambda_n}\right) , & \|\mc{E}_{P_0}^*\{\widehat{q}_n - q_{P_0}^\dagger\}\|_{L^2(P_{X;0})} =&O_p\left(\delta_n + \lambda_n + \frac{\delta_n^3}{\lambda_n}\right).
    \end{align*}
\end{proposition}

\begin{proof}
    Note that $\widehat{\mc{E}}_n = \mc{E}_{P_n}$ and $\widehat{\mc{E}}^*_n = \mc{E}^*_{P_n}$ with $P_n$ the empirical distribution of $(O_i)_{i=1}^n$ and that $\mc{E}_{P_n}: L^2(P_{X;n}) \rightarrow L^2(P_{Z;n})$ and $\mc{E}^*_{P_n}: L^2(P_{Z;n}) \rightarrow L^2(P_{X;n})$ are $P_n$-adjoints.

    Now, let
    \begin{align*}
        \gamma_n \coloneqq& \max\left(\max_{x}\left|\frac{\widehat{P}_n(X=x)}{P_0(X=x)}\right|, \max_{z}\left|\frac{\widehat{P}_n(Z=z)}{P_0(Z=z)}\right|\right), \\
        \underline{\gamma}_n \coloneqq& \min\left(\min_{x}\left|\frac{\widehat{P}_n(X=x)}{P_0(X=x)}\right|, \min_{z}\left|\frac{\widehat{P}_n(Z=z)}{P_0(Z=z)}\right|\right).
    \end{align*}
    Then, because $P_0(X=x), P_0(Z=z)$ are bounded away from zero, $ \gamma_n, \underline{\gamma}_n \overset{p} \rightarrow 1$. Consequently,
    \begin{align*}
        \|\mc{E}_{P_n} - \mc{E}_{P_0}\|_{L^2(P_{X;n})\rightarrow L^2(P_{Z;n})} \leq \sqrt{\gamma_n/\underline{\gamma}_n} \|\mc{E}_{P_n} - \mc{E}_{P_0}\|_{L^2(P_{X;0})\rightarrow L^2(P_{Z;0})} = O_p(n^{-1/2})
    \end{align*}
    where $\|\mc{E}_{P_n} - \mc{E}_{P_0}\|_{L^2(P_{X;0})\rightarrow L^2(P_{Z;0})} = O_p(n^{-1/2})$ because $\mc{E}_{P_n}$ is a matrix of empirical conditional probabilities for a finite-discrete distribution. Furthermore, we have that
    \begin{align}
        \label{eq:norm-eqivalence-in-prob}
        \underline{\gamma}_{n}^{1/2}\|h\|_{L^2(P_{X;0})} < \|h\|_{L^2(P_{X;n})} < \gamma_{n}^{1/2}\|h\|_{L^2(P_{X;0})}
    \end{align}
    for any $h \in L^2(P_{X;0})$, and similarly for $L^2(P_{Z;0})$.

    Define
    \begin{align}
        \label{eq:delta-n-def-tikhonov}
        \delta_n \coloneqq \max\left(\|\mc{E}_{P_n} - \mc{E}_{P_0}\|_{L^2(P_{X;n})\rightarrow L^2(P_{Z;n})}, \|\mc{E}^*_{P_n} - \mc{E}^*_{P_0}\|_{L^2(P_{Z;n})\rightarrow L^2(P_{X;n})}, \|\widehat{r}_n - r_{P_0}\|_{L^2(P_{Z;n})}, \|\widehat{\alpha}_n - \alpha_{P_0}\|_{L^2(P_{X;n})}\right)
    \end{align}
    and note that because $\mc{E}_{P_n}$ consist of empirical conditional probabilities and $\widehat{r}_n$ empirical means for a finite-discrete distribution, the above display is $O_p(n^{-1/2})$.

    Let $({\sigma}_{i;P_n}, {u}_{i;P_n}, {v}_{i;P_n})$ denote a singular system for $\mc{E}_{P_n}$,  $({\sigma}_{i;P_0}^n, {u}^n_{i;P_0}, {v}_{i;P_0}^n)$ denote a singular system for $\mc{E}_{P_0}$ computed with respect to the empirical distribution inner products, and $({\sigma}_{i;P_0}, {u}_{i;P_0}, {v}_{i;P_0})$ a singular system for $\mc{E}_{P_0}$ computed with respect to the $P_0$ inner products. Now, because $P_n(X=x) \overset{p}\rightarrow P_0(X=x), P_n(Z=z) \overset{p}\rightarrow P_0(Z=z)$ for each $x, z$, $\sigma_{i;P_0}^n \rightarrow \sigma_{i;P_0}$ for each $i$. We let $k$ denote the index of the smallest non-zero singular value of $\mc{E}_{P_0}$, which is the same if the singular system is computed with respect to the empirical distribution or $P_0$.

    It follows from Weyl's inequality that
    \begin{align*}
        \sigma_{i;P_n} > \sigma_{i;P_0}^n - \delta_n, i \leq  k &&\sigma_{i;P_n} < \delta_n, i > k.
    \end{align*}
    Now, in light of $r_{P_0} = \mc{E}_{P_0}h^\dagger_{P_0}$,
    \begin{align*}
        |\langle v_{i;P_n}, r_{P_0}\rangle_{L^2(P_{Z;n})}| =& |\langle v_{i;P_n}, \{\mc{E}_{P_0} - \mc{E}_{P_n}\}h^\dagger_{P_0}\rangle_{L^2(P_{Z;n})} + \langle v_{i;P_n}, \mc{E}_{P_n}h^\dagger_{P_0}\rangle_{L^2(P_{Z;n})}|\\
        =& |\langle v_{i;P_n}, \{\mc{E}_{P_0} - \mc{E}_{P_n}\}h^\dagger_{P_0}\rangle_{L^2(P_{Z;n})} + \langle \mc{E}_{P_n}^*v_{i;P_n}, h^\dagger_{P_0}\rangle_{L^2(P_{X;n})}| \\
        =& |\langle v_{i;P_n}, \{\mc{E}_{P_0} - \mc{E}_{P_n}\}h^\dagger_{P_0}\rangle_{L^2(P_{Z;n})} + \sigma_{i;P_n}\langle u_{i;P_n}, h^\dagger_{P_0}\rangle_{L^2(P_{X;n})}| \\
        \leq&\delta_n \|h^\dagger_{P_0}\|_{L^2(P_{X;n})} + \sigma_{i;P_n}\|h^\dagger_{P_0}\|_{L^2(P_{X;n})}.
    \end{align*}

    Then, for $i > k$,
    \begin{align*}
        |\langle v_{i;P_n}, \widehat{r}_n\rangle_{L^2(P_{Z;n})}|=&  |\langle v_{i;P_n}, \widehat{r}_n - {r}_{P_0}\rangle_{L^2(P_{Z;n})} + \langle v_{i;P_n}, {r}_{P_0}\rangle_{L^2(P_{Z;n})}|\\
        \leq & \|\widehat{r}_n - {r}_{P_0}\|_{L^2(P_{Z;n})} + \delta_n \|h^\dagger_{P_0}\|_{L^2(P_{X;n})} + \sigma_{i;P_n}\|h^\dagger_{P_0}\|_{L^2(P_{X;n})}\\
        \leq& \delta_n( 1 + 2\|h^\dagger_{P_0}\|_{L^2(P_{X;n})})
    \end{align*}
    where we have used that $\sigma_{i;P_n} < \delta_n$ for $i > k$.

    Now, note that
    \begin{align*}
        \widehat{h}_n \coloneqq \sum_{i=1}^{\min(K_X, K_Z)}\frac{\sigma_{i;P_n}}{\sigma_{i;P_n}^2 + \lambda_n} \langle v_{i;P_n}, \widehat{r}_n\rangle_{L^2(P_{Z;n})}u_{i;P_n}.
    \end{align*}
    Let
    \begin{align}
        \label{eq:h_2-def-convergence-rates}
        \widehat{h}^{(2)}_n \coloneqq \sum_{i> k}\frac{\sigma_{i;P_n}}{\sigma_{i;P_n}^2 + \lambda_n} \langle v_{i;P_n}, \widehat{r}_n\rangle_{L^2(P_{Z;n})}u_{i;P_n}
    \end{align}
    Then, for $i > k$, $\frac{\sigma_{i;P_n}}{\sigma_{i;P_n}^2 + \lambda_n} \leq \frac{\sigma_{i;P_n}}{\lambda_n} \leq \frac{\delta_n}{\lambda_n}$ and so
    \begin{align*}
        \|\widehat{h}^{(2)}_n\|_{L^2(P_{X;n})} \leq& \sqrt{\min(K_X, K_Z)}\frac{\delta_n}{\lambda_n} \delta_n( 1 + 2\|h^\dagger_{P_0}\|_{L^2(P_{X;n})})
    \end{align*}

    On the other hand, let
    \begin{align*}
        \widehat{h}^{(1)}_n \coloneqq \sum_{i\leq k}\frac{\sigma_{i;P_n}}{\sigma_{i;P_n}^2 + \lambda_n} \langle v_{i;P_n}, \widehat{r}_n\rangle_{L^2(P_{Z;n})}u_{i;P_n}
    \end{align*}
    and let
    \begin{align*}
        \Pi_{n}: L^2(P_{X;n}) \rightarrow L^2(P_{X;n}), \Pi_{n}h \coloneqq \sum_{i \leq k} u_{i;P_n}\langle h, u_{i;P_n}\rangle_{L^2(P_{X;n})}
    \end{align*}
    denote the projection onto the first $k$ right singular functions of $\mc{E}_{P_n}$. Then, $\widehat{h}^{(1)}_n -  \Pi_{n}h^\dagger_{P_0} = \sum_{i \leq k}c_{i;n}u_{i;P_n}$, where, using $r_{P_0} = \mc{E}_{P_0}h^\dagger_{P_0}$ and $\langle v_{i;P_n}, \mc{E}_{P_n}h^\dagger_{P_0}\rangle_{L^2(P_{Z;n})} = \sigma_{i;P_n}\langle u_{i;P_n}, h^\dagger_{P_0}\rangle_{L^2(P_{X;n})}$,
    \begin{align*}
        c_{i;n} \coloneqq& \frac{\sigma_{i;P_n}}{\sigma_{i;P_n}^2 + \lambda_n} \langle v_{i;P_n}, \widehat{r}_n\rangle_{L^2(P_{Z;n})} - \langle u_{i;P_n}, h^\dagger_{P_0}\rangle_{L^2(P_{X;n})}\\
        =&\frac{\sigma_{i;P_n}}{\sigma_{i;P_n}^2 + \lambda_n} \left\{\langle v_{i;P_n}, \widehat{r}_n - r_{P_0}\rangle_{L^2(P_{Z;n})} + \langle v_{i;P_n}, \{\mc{E}_{P_0} - \mc{E}_{P_n}\}h^\dagger_{P_0}\rangle_{L^2(P_{Z;n})}\right\} \\
        & + \left\{\frac{\sigma^2_{i;P_n}}{\sigma_{i;P_n}^2 + \lambda_n}-1\right\}\langle u_{i;P_n}, h^\dagger_{P_0}\rangle_{L^2(P_{X;n})}.
    \end{align*}
    Consequently, by the definition of $\delta_n$ and Cauchy--Schwarz,
    \begin{align*}
        |c_{i;n}| \leq \frac{\sigma_{i;P_n}}{\sigma_{i;P_n}^2 + \lambda_n} \delta_n(1 + \|h^\dagger_{P_0}\|_{L^2(P_{X;n})}) + \left|\frac{\sigma^2_{i;P_n}}{\sigma_{i;P_n}^2 + \lambda_n}-1\right|\|h^\dagger_{P_0}\|_{L^2(P_{X;n})}.
    \end{align*}

    Now, conditioning on the event that $\sigma_{k;P_n} > \sigma_{k;P_0}^n/2$, which happens with probability tending to one in light of Weyl's inequality, for each $i \leq k$,
    \begin{align*}
        \frac{\sigma_{i;P_n}}{\sigma_{i;P_n}^2 + \lambda_n} \delta_n(1 + \|h^\dagger_{P_0}\|_{L^2(P_{X;n})}) \leq \sigma_{i;P_n}^{-1} \delta_n(1 + \|h^\dagger_{P_0}\|_{L^2(P_{X;n})}) \leq \{\sigma^n_{k;P_0}/2\}^{-1} \delta_n(1 + \|h^\dagger_{P_0}\|_{L^2(P_{X;n})}).
    \end{align*}
    On the other hand,
    \begin{align*}
        \left|\frac{\sigma^2_{i;P_n}}{\sigma_{i;P_n}^2 + \lambda_n}-1\right| = \frac{\lambda_n}{\sigma_{i;P_n}^2 + \lambda_n} \leq \frac{\lambda_n}{\sigma^2_{i;P_n}} \leq \frac{\lambda_n}{\{\sigma^n_{k;P_0}/2\}^2}.
    \end{align*}
    Consequently, with probability tending to 1,
    \begin{align*}
        &\|\widehat{h}^{(1)}_n -  \Pi_{n}h^\dagger_{P_0}\|_{L^2(P_{X;n})} \\
        &\leq \sqrt{\min(K_X, K_Z)}\left\{\frac{\lambda_n}{\{\sigma^n_{k;P_0}/2\}^2}\|h^\dagger_{P_0}\|_{L^2(P_{X;n})} +  \{\sigma^n_{k;P_0}/2\}^{-1} \delta_n(1 + \|h^\dagger_{P_0}\|_{L^2(P_{X;n})})\right\}.
    \end{align*}

    We turn now to the term $h^\dagger - \Pi_{n}h^\dagger$. Now, because $h^\dagger \in \mathscr{N}(\mc{E}_{P_0})^\perp$, $h^\dagger \in \mathscr{R}(\mc{E}_{P_0}^*)$. Let $\omega_{P_0} \in L^2(P_{Z;0}) $ be such that $h^\dagger_{P_0} = \mc{E}_{P_0}^* \omega_{P_0}$. By a similar argument as before,
    \begin{align*}
        |\langle h^\dagger, u_{i;n}\rangle_{L^2(P_{X;n})}| =& |\langle \mc{E}_{P_0}^* \omega_{P_0}, u_{i;n}\rangle_{L^2(P_{X;n})}|\\
        =& |\langle \{\mc{E}_{P_0}^* - \mc{E}_{P_n}^*\}  \omega_{P_0} , u_{i;n}\rangle_{L^2(P_{X;n})} +  \langle \mc{E}_{P_n}^*  \omega_{P_0} , u_{i;n}\rangle_{L^2(P_{X;n})}|\\
        \leq& \delta_n\|\omega_{P_0}\|_{L^2(P_{Z;n})} + \sigma_{i;P_n} \|\omega_{P_0}\|_{L^2(P_{Z;n})}.
    \end{align*}
    Then, for $i > k$,
    \begin{align*}
        |\langle h^\dagger, u_{i;P_n}\rangle| \leq \delta_n\|\omega_{P_0}\|_{L^2(P_{Z;n})} + \sigma_{i;P_n} \|\omega_{P_0}\|_{L^2(P_{Z;n})} \leq 2\delta_n \|\omega_{P_0}\|_{L^2(P_{Z;n})}.
    \end{align*}
    Consequently,
    \begin{align*}
        \|h^\dagger - \Pi_{n}h^\dagger\|_{_{L^2(P_{X;n})}} \leq& \left\{\sum_{i> k} \langle h^\dagger, u_{i;n}\rangle^2_{L^2(P_{X;n})}\right\}^{1/2}\\
        \leq & \sqrt{K_X} 2\delta_n \|\omega_{P_0}\|_{L^2(P_{Z;n})}.
    \end{align*}

    Then, by the above arguments,
    \begin{align*}
        \|\widehat{h}_n - h^\dagger_{P_0}\|_{L^2(P_{X;n})} \leq& \|h^\dagger - \Pi_{n}h^\dagger\|_{_{L^2(P_{X;n})}} +  \|\widehat{h}^{(1)}_n -  \Pi_{n}h^\dagger_{P_0}\|_{L^2(P_{X;n})} + \|\widehat{h}^{(2)}_n\|_{L^2(P_{X;n})} \\
        \leq & \sqrt{K_X} 2\delta_n \|\omega_{P_0}\|_{L^2(P_{Z;n})}\\
        &+ \sqrt{\min(K_X, K_Z)}\bigg\{\frac{\lambda_n}{\{\sigma^n_{k;P_0}/2\}^2}\|h^\dagger_{P_0}\|_{L^2(P_{X;n})}\\
        &\qquad+  \{\sigma^n_{k;P_0}/2\}^{-1} \delta_n(1 + \|h^\dagger_{P_0}\|_{L^2(P_{X;n})})\bigg\}\\
        &+ \sqrt{\min(K_X, K_Z)}\frac{\delta_n}{\lambda_n} \delta_n( 1 + 2\|h^\dagger_{P_0}\|_{L^2(P_{X;n})})
    \end{align*}
    with probability tending to 1. Recalling that for any functions $h \in L^2(P_{X;n}), q \in L^2(P_{Z;n})$, $\|h\|_{L^2(P_{X;n})} \overset{p}{\rightarrow} \|h\|_{L^2(P_{X;0})}$ and  $\|q\|_{L^2(P_{Z;n})} \overset{p}{\rightarrow} \|q\|_{L^2(P_{Z;0})}$,  and that $\sigma^n_{k;P_0} \overset{p}\rightarrow \sigma_{k;P_0}$, we arrive at
    \begin{align*}
        \|\widehat{h}_n - h^\dagger_{P_0}\|_{L^2(P_{X;n})} = O_p\left(\delta_n + \lambda_n + \frac{\delta_n^2}{\lambda_n}\right).
    \end{align*}
    Then, by \eqref{eq:norm-eqivalence-in-prob}, we have that
    \begin{align*}
        \|\widehat{h}_n - h^\dagger_{P_0}\|_{L^2(P_{X;0})} = O_p\left(\delta_n + \lambda_n + \frac{\delta_n^2}{\lambda_n}\right).
    \end{align*}

    The above gives the bound on the strong norm. For the weak norm bound, note that
    \begin{align*}
        \mc{E}_{P_0}u_{i;P_n} =& \{\mc{E}_{P_0} - \mc{E}_{P_n}\}u_{i;P_n} + \mc{E}_{P_n}u_{i;P_n}\\
        =& \{\mc{E}_{P_0} - \mc{E}_{P_n}\}u_{i;P_n} + \sigma_{i;P_n}v_{i;P_n}.
    \end{align*}
    Then, for $i > k$,
    \begin{align*}
        \|\mc{E}_{P_0}u_{i;P_n}\|_{L^2(P_{Z;n})} =& \|\{\mc{E}_{P_0} - \mc{E}_{P_n}\}u_{i;P_n} + \mc{E}_{P_n}u_{i;P_n}\|_{L^2(P_{Z;n})}\\
        \leq & \delta_n + \sigma_{i;P_n}\\
        \leq&  2\delta_n.
    \end{align*}
    Then, recalling the definition in \eqref{eq:h_2-def-convergence-rates} and subsequent arguments, we have that
    \begin{align*}
        \|\mc{E}_{P_0}\widehat{h}^{(2)}_n\|_{L^2(P_{Z;n})} \leq & \sum_{i> k}\frac{\sigma_{i;P_n}}{\sigma_{i;P_n}^2 + \lambda_n} \langle v_{i;P_n}, \widehat{r}_n\rangle_{L^2(P_{Z;n})}\|\mc{E}_{P_0}u_{i;P_n}\|_{L^2(P_{Z;n})}\\
        \leq& \sqrt{\min(K_X, K_Z)}2\frac{\delta_n^3}{\lambda_n}( 1 + 2\|h^\dagger_{P_0}\|_{L^2(P_{X;n})}).
    \end{align*}
    On the other hand,
    \begin{align*}
        \|\mc{E}_{P_0}\{\widehat{h}^{(1)}_n -  \Pi_{n}h^\dagger_{P_0}\}\|_{L^2(P_{Z;n})} &\leq \|\mc{E}_{P_0}\|_{L^2(P_{X;n}) \rightarrow L^2(P_{Z;n})}\|\widehat{h}^{(1)}_n -  \Pi_{n}h^\dagger_{P_0}\|_{L^2(P_{X;n})} \quad \text{and} \\
        \|\mc{E}_{P_0}\{h^\dagger - \Pi_{n}h^\dagger\}\|_{_{L^2(P_{X;n})}} &\leq \|\mc{E}_{P_0}\|_{L^2(P_{X;n}) \rightarrow L^2(P_{Z;n})}\|h^\dagger - \Pi_{n}h^\dagger\|_{_{L^2(P_{X;n})}}.
    \end{align*}
    Then, because $\|\mc{E}_{P_0}\|_{L^2(P_{X;n}) \rightarrow L^2(P_{Z;n})} \overset{p}{\rightarrow} \|\mc{E}_{P_0}\|_{L^2(P_{X;0}) \rightarrow L^2(P_{Z;0})}$, we have that
    \begin{align*}
        \|\mc{E}_{P_0}\{\widehat{h}_n - h^\dagger_{P_0}\}\|_{L^2(P_{Z;n})} = O_p\left(\delta_n + \lambda_n + \frac{\delta_n^3}{\lambda_n}\right).
    \end{align*}
    Then, by \eqref{eq:norm-eqivalence-in-prob}, we have that
    \begin{align*}
        \|\mc{E}_{P_0}\{\widehat{h}_n - h^\dagger_{P_0}\}\|_{L^2(P_{Z;0})} = O_p\left(\delta_n + \lambda_n + \frac{\delta_n^3}{\lambda_n}\right).
    \end{align*}
    The bounds for $\widehat{q}_n$ follow from the symmetric argument, interchanging the roles of $\mc{E}_{P_n}, \mc{E}_{P_0}, \widehat{r}_n, r_{P_0}, h^\dagger_{P_0}$ with $\mc{E}^*_{P_n}, \mc{E}^*_{P_0}, \widehat{\alpha}_n, \alpha_{P_0}, q^\dagger_{P_0}$. This completes the proof.
\end{proof}

\subsection{Derivations in \cref{sec:simulation}}

Here we give derivations showing the claimed divergence in the discussion of \cref{counterexample:npiv} in \cref{sec:simulation}.

Direct calculation shows that $\det(\mc{E}_{P_t}) = -t^3/5$. Now, by construction $r_{P_t} = r_{P_0} = (2, 2.2, 2.2)^T$. Algebra then gives that for $t> 0$,
\begin{align*}
    h_{P_t} = \mc{E}_{P_t}^{-1}r_{P_t} = \mc{E}_{P_t}^{-1}(2, 2.2, 2.2)^T = \frac{1}{5t^2}(7t^2 + 2, 2(6t^2 + 1), -3(1-4t^2))^T.
\end{align*}
Now, one can check that $\mc{E}_{P_0}h_{P_t} = r_{P_0}$. Hence, it also solves the primal inverse equation at $P_0$.

On the other hand,
\begin{align*}
    P_t(X) = (P_0(Z=1),P_0(Z=2), P_0(Z=3) )^T\mc{E}_{P_t} = P_0(X) + (- t, t(1-t^2), t^3)/3.
\end{align*}
Consequently,
\begin{align*}
    \varphi(P_t) =& E_{P_t}[\alpha_{P_t}(X)h_{P_t}(X)]\\
    =& E_{P_t}[\psi(X)h_{P_t}(X)]\\
    =& E_{P_0}[\psi(X)h_{P_t}(X)]\\
    &+ \frac{1}{3}(- t, t(1-t^2), t^3)^T \text{diag}(\psi(X))\frac{1}{5t^2}(7t^2 + 2, 2(6t^2 + 1), -3(1-4t^2))^T\\
    =& \varphi(P_0) + \frac{1}{15t^2}\left\{\frac{9}{10}t + \frac{15}{4}t^3  -\frac{12}{5}t^5\right\}.
\end{align*}
The final equality holds because $h_{P_t}$ also solves the primal inverse problem at $P_0$, so that $h_{P_t} - h^\dagger_{P_0} \in \mathscr{N}(\mc{E}_{P_0})$. Since $\psi = \alpha_{P_0} \in \mathscr{R}(\mc{E}_{P_0}^*) = \mathscr{N}(\mc{E}_{P_0})^{\perp}$, with the orthogonal complement taken in $L^2(P_{X;0})$, it follows that $E_{P_0}[\psi(X)\{h_{P_t}(X) - h^\dagger_{P_0}(X)\}] = 0$, and hence $E_{P_0}[\psi(X)h_{P_t}(X)] = E_{P_0}[\psi(X)h^\dagger_{P_0}(X)] = \varphi(P_0)$.

\section{Further Details for the Examples}
\label{sup-sec:examples}

Here we argue that under mild conditions on $P_0$, the Riesz representers $P \mapsto \alpha_P$, $P \mapsto r_P$ are both Hellinger continuous and pathwise differentiable at $P_0$ in the model $\mc{P}_{\text{Riesz}}$ defined in \eqref{eq:riesz-model}, and hence also relative to any submodel of $\mc{P}_{\text{Riesz}}$, including $\mc{P}$.

\begin{examplecont}{example:npiv}
    Recall that $P \mapsto \alpha_P = \psi$ which does not depend on $P$, and so is pathwise differentiable with local parameter $\dot\alpha[s] = 0$. On the other hand, recall that $r_P \coloneqq E_{P}[W|Z = \cdot]$, and that the support of $W$ is bounded. Then, $\sup_{P \in \mc{P}_{\text{Riesz}}}\esssup_zE_{P}[W^2|Z=z] < \infty$. Consequently, Example 5 of \cite{luedtke_one-step_2024} shows that $P \mapsto r_P$ is pathwise differentiable at $P$ in model $\mc{P}_{\text{Riesz}}$, and so also for any submodel of $\mc{P}_{\text{Riesz}}$. Example 5 of \cite{luedtke_one-step_2024} further shows that $P \mapsto r_P$ is Lipschitz with respect to the Hellinger distance over $\mc{P}_{\text{Riesz}}$, and so is Hellinger continuous at $P_0$ over $\mc{P}_{\text{Riesz}}$ and any submodel of it. As $\alpha_P = \psi$ does not depend on $P$, $P \mapsto \alpha_P$ is trivially Hellinger continuous.
\end{examplecont}

\begin{examplecont}{example:proximal}
    Recall that $r_{P} \coloneqq E_{P}[W|Z=\cdot, D = \cdot]$, which is pathwise differentiable by the same reasoning as $r_P$ in \cref{example:npiv}. On the other hand, $\alpha_{P}(X, D) \coloneqq \frac{2D -1}{P(D|X)}$. $P(D|X) = DE_{P}[D|X] + (1-D)E_{P}[1-D|X]$, which is a linear combination of two regression functions, each pathwise differentiable, and so $P(D|X)$ is pathwise differentiable at $P$ in $\mc{P}_{\text{Riesz}}$. If we additionally assume that $P_0(D|X) > \delta$ for some $\delta > 0$, then by \cref{as:dominating-measure-support} $P_0(D|X) > \delta'$ for some $\delta'> 0$. Consequently, the chain rule shows that $P \mapsto \frac{2D -1}{P(D|X)}$ is pathwise differentiable at $P$ in the model $\mc{P}_{\text{Riesz}}$.

    For Hellinger continuity, $P \mapsto r_P$ is Lipschitz with respect to the Hellinger distance over $\mc{P}_{\text{Riesz}}$ by the same reasoning as in \cref{example:npiv}, and similarly for the regression functions $E_P[D|X]$ and $E_P[1-D|X]$, so that $P \mapsto P(D|X)$ is also Lipschitz with respect to the Hellinger distance. Furthermore, because $D$ is binary and the density of every $P \in \mc{P}_{\text{Riesz}}$ is bounded above and below on a common rectangle by \cref{as:dominating-measure-support}, there exists $\delta' > 0$ such that $P(D=d|X) > \delta'$ $P$-a.s. for every $P \in \mc{P}_{\text{Riesz}}$ and $d \in \{0, 1\}$. As $x \mapsto 1/x$ is Lipschitz on $[\delta', 1]$ and $|2D - 1| = 1$, $P \mapsto \alpha_P$ is Lipschitz with respect to the Hellinger distance over $\mc{P}_{\text{Riesz}}$, and so is Hellinger continuous at $P_0$ over $\mc{P}_{\text{Riesz}}$ and any submodel of it.
\end{examplecont}

\begin{examplecont}{example:shadow}
    Recall that $r_P = DZ$, which does not depend on $P$ and as such is pathwise differentiable with local parameter $\dot r[s] = 0$. On the other hand, if $P_0(D=1|X) > \delta > 0$ for some $\delta>0$, the same reasoning as $\alpha_P$ in \cref{example:proximal} shows that $P \mapsto \alpha_P$ is pathwise differentiable at $P$ in the  model  $\mc{P}_{\text{Riesz}}$. The same reasoning as in \cref{example:proximal} also shows that $P \mapsto \alpha_P$ is Lipschitz with respect to the Hellinger distance over $\mc{P}_{\text{Riesz}}$, and as $r_P$ does not depend on $P$, both Riesz representer maps are Hellinger continuous at $P_0$ over $\mc{P}_{\text{Riesz}}$ and any submodel of it.
\end{examplecont}

\section{Proofs of Main Text Results}

{
\begin{remark}
    \label{remark:marginal-conditional-density-bounds}
    In the following proofs, we always take $P_0 \in \mc{P}$. Hence, $P_0$ has bounded rectangular support. Let $\Omega_W$ denote the support of $W$ under $P_0$, and similarly for $\Omega_X, \Omega_Z$, $\Omega_D$. Throughout, we will typically assume that the reference measures $\lambda_W$, $\lambda_X, \lambda_Z, \lambda_D$ are normalized to satisfy that $\lambda_W(\Omega_W) = \lambda_X(\Omega_X) = \lambda_Z(\Omega_Z) = \lambda_D(\Omega_D) = 1$ to avoid tracking unnecessary constants. This is easily done by taking $\lambda_W$ as the normalized Lebesgue measure over $\Omega_W$, when $W$ is continuous, the normalized counting measure when $W$ is discrete, or the product of the normalized measures when $W$ contains both continuous and discrete components, and similarly for $\lambda_X, \lambda_Z, \lambda_D$.

    Consequently, if $P$, mutually absolutely continuous with $P_0$, has density $p$ satisfying
    \begin{align*}
        L_{P}^{-1} < p(O) < L_P,
    \end{align*}
    it holds that
    \begin{align*}
        L_{P}^{-1} <p_{X}(X) < L_P, && L_{P}^{-2} <p_{X|Z}(X|Z) < L_P^2,
    \end{align*}
    and similarly for any other combination of marginal and conditional density functions. To see this, note that
    \begin{align*}
        \int p(w, x, z, d) d(\lambda_W\otimes\lambda_Z\otimes\lambda_D)(w, z, d) > \int_{\Omega_W \times \Omega_Z \times \Omega_D} L_{P}^{-1} d(\lambda_W\otimes\lambda_Z\otimes\lambda_D)(w, z, d) = L_{P}^{-1},
    \end{align*}
    and
    \begin{align*}
        p_{X|Z}(X|Z) = \frac{p(X, Z)}{p(Z)} > \frac{L_{P}^{-1}}{L_P} = L_{P}^{-2}.
    \end{align*}
    A similar calculation gives the upper bound.
\end{remark}
}

\subsection{Technical Lemmas}
We begin by proving some technical lemmas. In what follows, given a law $P$ for $(W,X,Z, D)$ with density $p$ with respect to some dominating measure, we write $P_V$ and $p_V$ for the corresponding law and density of a subvector $V$ of $(W,X,Z, D)$. Because $D$ is finite discrete, the operators $\mc{E}_P$ and $\mc{E}_P^*$ act blockwise in $d$, i.e. $(\mc{E}_Ph)(\cdot, d) = \mc{E}_{P;d}h(\cdot, d)$. When we describe such operators as integral operators with a kernel in $(x, z, d)$, or refer to their Hilbert--Schmidt norms, we mean that for each $d$ the block is an integral operator with the corresponding kernel in $(x, z)$. As $D$ takes finitely many values, a property of the Hilbert--Schmidt norms that holds for every block, such as finiteness or convergence, holds for the operator as a whole.
\begin{lemma}
\label{lemma:constant-riesz-representers}
Let \(P\) satisfy \cref{as:continuous-linear-maps} and
\cref{as:dominating-measure-support}. { Suppose that \(P'_{X,Z, D}\) is a law
on the measurable space of \((X,Z, D)\), with density \(p'_{X,Z, D}\) with respect to
\(\lambda_X\otimes\lambda_Z\otimes\lambda_D\), such that
\(
p'_{X, D}=p_{X, D}, \,
p'_{Z, D}=p_{Z, D} \,
\text{a.e.}(P)
\), such that \(p'_{X,Z,D}\) vanishes outside the support of \(P_{X,Z,D}\), and additionally satisfying
\[
\|p'_{X,Z, D}-p_{X,Z, D}\|_{L^\infty(\lambda_X\otimes\lambda_Z\otimes\lambda_D)}
<
L_{bounded}^{-1}\min(\essinf p - L_{bounded}^{-1}, L_{bounded} - \esssup p)\{\lambda_X(\mathcal X)\lambda_Z(\mathcal Z)\lambda_D(\mc{D})\}^{-1}.
\]
Then there exists a probability law \(P^\dagger\) on the measurable
space of \((W,X,Z, D)\), with the same support as \(P\), satisfying
\cref{as:continuous-linear-maps} and \cref{as:dominating-measure-support},
such that
\[
p^\dagger_{X,Z, D}=p'_{X,Z, D},
\qquad
r_{P^\dagger}=r_P,
\qquad
\alpha_{P^\dagger}=\alpha_P.
\]
Moreover,
\begin{equation}
\label{eq:firstbound}
H^2(P^\dagger,P)
\le
\sqrt{2H^2(P'_{X,Z, D},P_{X,Z, D})},
\end{equation}
\begin{equation}
\label{eq:secondbound}
H^2(P^\dagger,P)
\le
\frac12
\|p'_{X,Z, D}-p_{X,Z, D}\|_{L^1(\lambda_X\otimes\lambda_Z\otimes\lambda_D)}.
\end{equation}}
\begin{equation}
    \label{eq:constant-riesz-representers:supbound}
    \|p^\dagger - p\|_{\infty} \leq L_{bounded}\|p'_{X, Z, D} - p_{X, Z, D}\|_{L^\infty(\lambda_X \otimes \lambda_Z \otimes \lambda_D)}\{\lambda_X(\mathcal X)\lambda_Z(\mathcal Z)\lambda_D(\mc{D})\}.
\end{equation}
\end{lemma}

\begin{proof}
    Because the maps
    \begin{align*}
        h \mapsto E_{P}[\widetilde{m}(W, X, D;h)], &&q \mapsto E_{P}[{m}(W, Z, D;q)]
    \end{align*}
    depend on $P$ only through the joint distributions of $(W, X, D)$ and $(W, Z, D)$, respectively, $\alpha_{P^\dagger} = \alpha_{P}$ and $r_{P^\dagger} = r_{P}$ if $p^\dagger(W, X, D) = p(W, X, D)$ and $p^\dagger(W, Z, D) = p(W, Z, D)$. In that case, the two functionals and the spaces $L^2(P_{X, D})$ and $L^2(P_{Z, D})$ are unchanged, so that $P^\dagger$ also satisfies \cref{as:continuous-linear-maps}. Hence, it suffices to find a joint density $p^\dagger(w, x, z, d)$ with $L_{bounded}^{-1}< \essinf p^\dagger \leq \esssup p^\dagger < L_{bounded}$ such that
    \begin{align}
        \label{eq:marginals-hold}
        \nonumber \int p^\dagger(w, x, z, d) d\lambda_W(w)  = p'(x, z, d)\\
        \nonumber \int p^\dagger(w, x, z, d) d\lambda_X(x)  = p(w, z, d)\\
        \int p^\dagger(w, x, z, d) d\lambda_Z(z)  = p(w, x, d)
    \end{align}
    and with { corresponding law $P^\dagger$ satisfying \eqref{eq:firstbound} and \eqref{eq:secondbound}.

    Let
\(
\Delta(x,z, d):=p'_{X,Z, D}(x,z, d)-p_{X,Z, D}(x,z, d)
\)
and define
\[
p^\dagger(w,x,z, d)
:=
p(w,x,z, d)+p_W(w)\Delta(x,z, d).
\]
Then, each of the equalities in \eqref{eq:marginals-hold} hold. Consequently, $\int p^\dagger(w, x,z, d)\allowbreak d\lambda_W(w)\lambda_X(x)\lambda_Z(z)\lambda_D(d) = 1$.

Set
\(
M_{XZD}:=\lambda_X(\mathcal X)\lambda_Z(\mathcal Z)\lambda_D(\mathcal{D}).
\)
By \cref{as:dominating-measure-support}, \(p(w,x,z, d)\le L_{bounded}\) on the support of $P$, and since
\((W,X,Z, D)\) is supported on the bounded rectangle
\(\mathcal W\times\mathcal X\times\mathcal Z \times \mathcal{D}\), we have
\[
p_{W}(w)
=
\int_{\mathcal X\times\mathcal Z\times \mathcal{D}}
p(w,x,z, d)\,d\lambda_X(x)d\lambda_Z(z)d\lambda_D(d)
\le
L_{bounded} M_{XZD}.
\]
Hence,
\begin{align*}
    \|p_W\Delta\|_{L^\infty(\lambda_W\otimes\lambda_X\otimes\lambda_Z\otimes\lambda_D)}
\le&
\|p_W\|_{L^\infty(\lambda_W)} \|\Delta\|_{L^\infty(\lambda_X\otimes\lambda_Z\otimes\lambda_D)}\\
<&
L_{bounded} M_{XZD}\,
L_{bounded}^{-1}\min(\essinf p - L_{bounded}^{-1}, L_{bounded} - \esssup p)M_{XZD}^{-1}
\\
=&
\min(\essinf p - L_{bounded}^{-1}, L_{bounded} - \esssup p).
\end{align*}
Because \(p^\dagger - p = p_W\Delta\), the first inequality in the preceding display, together with \(\|p_W\|_{L^\infty(\lambda_W)} \le L_{bounded}M_{XZD}\), gives \eqref{eq:constant-riesz-representers:supbound}.
It follows that
\begin{gather*}
    \essinf p^\dagger(w,x,z, d)
=
\essinf \{p(w,x,z, d)+p_W(w)\Delta(x,z, d)\}
>
L_{bounded}^{-1}\\
\esssup p^\dagger(w,x,z, d)
=
\esssup \{p(w,x,z, d)+p_W(w)\Delta(x,z, d)\}
<
L_{bounded}
\end{gather*}
where the essinf is taken with respect to $\lambda$ on the support of $P$. Moreover, \(p^\dagger(w,x,z, d)=0\) outside the support of
\(P\), so \(P^\dagger\) has the same support as \(P\).

It remains to prove the Hellinger bounds. Using the elementary inequality
\[
(\sqrt a-\sqrt b)^2\le |a-b|,
\qquad a,b\ge0,
\]
we obtain
\begin{align*}
H^2(P^\dagger,P)
&=
\frac12
\int
\left\{
\sqrt{p^\dagger(w,x,z, d)}
-
\sqrt{p(w,x,z, d)}
\right\}^2
\,d\lambda_W(w)d\lambda_X(x)d\lambda_Z(z)d\lambda_D(d)
\\
&\le
\frac12
\int
\left|
p^\dagger(w,x,z, d)-p(w,x,z, d)
\right|
\,d\lambda_W(w)d\lambda_X(x)d\lambda_Z(z)d\lambda_D(d)
\\
&=
\frac12
\int
\left|
p_W(w)\{p'_{X,Z, D}(x,z, d)-p_{X,Z, D}(x,z, d)\}
\right|
\,d\lambda_W(w)d\lambda_X(x)d\lambda_Z(z)d\lambda_D(d)
\\
&=
\frac12
\int
\left|
p'_{X,Z, D}(x,z, d)-p_{X,Z, D}(x,z, d)
\right|
\,d\lambda_X(x)d\lambda_Z(z)d\lambda_D(d)
\\
&=
\operatorname{TV}(P'_{X,Z, D},P_{X,Z, D}).
\end{align*}
By the standard relation between total variation and Hellinger distance,
\[
\operatorname{TV}(P'_{X,Z, D},P_{X,Z, D})
\le
\sqrt{2}\,H(P'_{X,Z, D},P_{X,Z, D})
=
\sqrt{2H^2(P'_{X,Z, D},P_{X,Z, D})}.
\]
Therefore,
\[
H^2(P^\dagger,P)
\le
\sqrt{2H^2(P'_{X,Z, D},P_{X,Z, D})}.
\]
The preceding display also gives
\[
H^2(P^\dagger,P)
\le
\frac12
\|p'_{X,Z, D}-p_{X,Z, D}\|_{L^1(\lambda_X\otimes\lambda_Z\otimes\lambda_D)}.
\]}

\end{proof}

{ In what follows  $C^\infty(\mathbb{R}^d)$ is the class of real valued infinitely differentiable functions on $\mathbb{R}^d$ and $C_c^\infty(\mathbb{R}^d)$ is the subclass of $C^\infty(\mathbb{R}^d)$ of functions that take value 0 outside a compact set. For $\mc{A} \subseteq \mathbb{R}^d$ open and bounded, we define $C^\infty(\mc{A})$ the space of infinitely differentiable functions with domain $\mc{A}$, and $C^\infty(\overline{\mc{A}})$ the space of infinitely differentiable functions whose derivatives extend continuously to $\overline{\mc{A}}$. The class $L_{loc}^1(\lambda)$ is the class of locally integrable functions, that is, the class of real valued functions $f$ on $\mathbb{R}^d$ such that $\int_K |f(x)|d\lambda(x) <\infty$ for any compact set $K$, where $\lambda$ is the Lebesgue measure on $\mathbb{R}^d$ and the dimension $d$ is understood according to the context. Furthermore, for $g \in C_c^\infty(\mathbb{R}^d)$ and $f \in L_{loc}^1(\lambda)$, $g*f$ is the convolution of $g$ with $f$, i.e.
\[
(g*f) (y) :=\int g(y-y') f(y') dy'
\]
and it holds that $g*f \in C_c^\infty(\mathbb{R}^d).$

The next lemma will invoke the following results:
\begin{enumerate}

\item \textit{(Special case of the Approximate Identity Theorem (\cite{folland_real_2013}, Theorem 8.14.))} Suppose that $\rho$ is a non-negative real-valued function on $\mathbb{R}^d$ supported on $B(0,1):=\{y \in \mathbb{R}^d:\|y\|\le 1\}$ with $\int \rho(x)d\lambda(x)=1$, and for each $t>0$ let $\rho_t(y) \coloneqq t^{-d}\rho(y/t)$, so that $\rho_t$ is supported on $B(0,t)$ and $\int \rho_t(x)d\lambda(x)=1$, where $\lambda$ is the Lebesgue measure on $\mathbb{R}^d$. Let $q \in [1,\infty)$. Then, if $f \in L^q(\lambda)$ it holds that $\|\rho_t * f-f\|_{L^q(\lambda)} \rightarrow 0$ as $t \rightarrow 0$.

\item \textit{(Young's convolution inequality).}
Let \(1\le q\le\infty\). Suppose
\(f\in L^1(\mathbb R^d)\) and
\(g\in L^q(\mathbb R^d)\).
Then the convolution \(f*g\) belongs to \(L^q(\mathbb R^d)\), and
\(
\|f*g\|_{L^q}
\le
\|f\|_{L^1}
\|g\|_{L^q}.
\)

\end{enumerate}

We note that this lemma will be invoked in the proofs of \cref{thrm:continuous-discontinuity}, \cref{thrm:continuous-tangent-space} to construct a sequence of laws $P_j \rightarrow P$ that each admit continuously differentiable versions of their densities. These well behaved $P_j$'s will then be further perturbed to construct the sequences and submodels used in those results. A similar construction appears in \cite{canayTestabilityIdentificationNonparametric2013a} to prove that completeness is untestable when $X, Z$ are continuous.

\begin{lemma}
\label{lemma:dense-continuous-densities}
Let \(\Omega=\prod_{\ell=1}^{d_Y}(a_\ell,b_\ell)\subset\mathbb R^{d_Y}\) be a bounded open rectangle.
Let \(p\) be a density for \(Y\) with respect to the Lebesgue measure $\lambda$, supported on \(\Omega\), and suppose that there exist constants \(0<a<b<\infty\) such that
\(
a<p(y)<b \,\,
\text{for Lebesgue-a.e. }y\in\Omega .
\)
Then, for every \(\varepsilon>0\), there exists a law on the measurable space of $Y$ with density \(p_\varepsilon\) supported on \(\Omega\) such that:
\begin{enumerate}
\item \(p_\varepsilon\in C^\infty(\overline{\Omega})\), where \(\overline{\Omega} :=\prod_{\ell=1}^{d_Y}[a_\ell,b_\ell]\);
\item \(p_\varepsilon\) and \(p\) are mutually absolutely continuous;
\item there exists \(\gamma(\varepsilon)>0\) such that
\(
\gamma(\varepsilon)\to1
\, \text{as }\varepsilon\downarrow0,
\)

\[
\frac{a}{\gamma(\varepsilon)}
<
p_\varepsilon(y)
<
\frac{b}{\gamma(\varepsilon)}
\qquad
\text{for Lebesgue-a.e. }y\in\Omega;
\]
\end{enumerate}
Furthermore, for every \(1\le q<\infty\),
\[
\|p_\varepsilon-p\|_{L^q(\lambda)}\to0
\qquad
\text{as }\varepsilon\downarrow0 .
\]

\end{lemma}

\begin{proof}
Let
\begin{align*}
        \rho(y) \coloneqq \begin{cases}
            k\exp\left(\frac{-1}{1-\|y\|^2}\right) &\text{if } \|y\| < 1 \\
            0 &\text{if } \|y\| \geq 1
        \end{cases}
    \end{align*}
where $k:=\big\{\int_{y:\|y\|<1} \exp\left(\frac{-1}{1-\|y\|^2}\right) d\lambda(y)\big\}^{-1} $ and for \(t>0\), define
\[
\rho_t(y)=t^{-d_Y}\rho(y/t).
\] Then
\(
\rho_t:\mathbb R^{d_Y}\mapsto[0,\infty)
\)
satisfies \(\rho_t\in C_c^\infty(\mathbb R^{d_Y})\), it is supported on $B(0,t)$
and $\int \rho_t(x)dx=1$. Consequently,
for any $q \in [1,\infty)$, if $f \in L^q(\lambda) \cap L^1_{loc}(\lambda)$ it holds that $\rho_t * f\in C^\infty(\mathbb{R}^{d_Y})$ and $\|\rho_t * f-f\|_{L^q(\lambda)} \rightarrow 0$ as $t \rightarrow 0$.

For \(y\in\mathbb R^{d_Y}\), define the distance from \(y\) to \(\Omega\) by
\(
\operatorname{dist}(y,\Omega)
:=
\inf_{z\in\Omega}\|y-z\|,
\)
where \(\|\cdot\|\) denotes the Euclidean norm.
Let
\(
\Omega_t:=\{y\in\mathbb R^{d_Y}:\operatorname{dist}(y,\Omega)<t\}
\) and define
\[
f_t(y):=p(y)+a\mathbf I\{y\in\Omega_t\setminus\Omega\}.
\]
Then \(0\le f_t\le b\), and \(f_t\to p\) in \(L^q(\lambda)\) for every \(1\le q<\infty\), because
\(
\lambda(\Omega_t\setminus\Omega)\to0
\,
\text{as }t\downarrow0.
\)

Now define
\[
p_t^\ast(y):=(\rho_t*f_t)(y)
=
\int \rho_t(y-y')f_t(y')\,dy' .
\]
Since \(f_t\in L_{loc}(\lambda)\) (by virtue of being bounded) and \(\rho_t\in C_c^\infty(\mathbb{R}^{d_Y})\), the function \(p_t^\ast\) is in \(C^\infty(\mathbb{R}^{d_Y})\).

Moreover, for \(y\in\Omega\) and \(y'\in\operatorname{supp}(\rho_t(y-\cdot))\), we have \(y'\in\Omega_t\).
On \(\Omega_t\), the function \(f_t\) is bounded below by \(a\) and above by \(b\). Therefore,
\begin{equation}
\label{eq:p*tBound}
a<p_t^\ast(y)<b
\qquad
\text{for all }y\in\Omega .
\end{equation}
Now,
\begin{align*}
\|p_t^\ast-p\|_{L^1(\lambda)}
&=
\|\rho_t*f_t-p\|_{L^1(\lambda)}
\\
&\le
\|\rho_t*p-p\|_{L^1(\lambda)}
+
\|\rho_t*(f_t-p)\|_{L^1(\lambda)}
\\
&\le
\|\rho_t*p-p\|_{L^1(\lambda)}
+
\|\rho_t\|_{L^1(\lambda)}\|f_t-p\|_{L^1(\lambda)}
\\
&=
\|\rho_t*p-p\|_{L^1(\lambda)}
+
\|f_t-p\|_{L^1(\lambda)}
\\
&\to0.
\end{align*}

Let
\[
\gamma(t):=\int_\Omega p_t^\ast(y)\,dy.
\]
Then,
\[
|\gamma(t)-1| \le
\int_\Omega |p_t^\ast(y)-p(y)|\,dy \le \|p_t^\ast-p\|_{L^1(\lambda)} \rightarrow 0 \qquad
\text{as } t\downarrow0
\]

Moreover, \(\gamma(t)\ge a\lambda(\Omega)>0\) for every \(t>0\) by \eqref{eq:p*tBound}. Define
\[
p_t(y):=
\begin{cases}
p_t^\ast(y)/\gamma(t), & y\in\Omega,\\
0, & y\notin\Omega.
\end{cases}
\]
Then \(p_t\) is a density supported on \(\Omega\), and
\[
\frac{a}{\gamma(t)}
<
p_t(y)
<
\frac{b}{\gamma(t)}
\qquad
\text{for Lebesgue-a.e. }y\in\Omega.
\]
Since both \(p\) and \(p_t\) are strictly positive and finite on \(\Omega\), and both vanish outside \(\Omega\), we have
\(
p_t\) and $p$ are mutually absolutely continuous.

Finally, we show that
\(
\|p_t-p\|_{L^q(\lambda)}\longrightarrow0.
\)
Define for any $g$, $\|g\|_{L^q(\Omega)} := \{\int_{\Omega}|g(x)|^qd\lambda(x)\}^{1/q}$. Since both \(p_t\) and \(p\) vanish outside \(\Omega\), we have that
\(
\|p_t-p\|_{L^q(\lambda)}
=
\|p_t-p\|_{L^q(\Omega)}.
\)
Moreover,
\[
p_t-p
=
\left(\frac{1}{\gamma(t)}-1\right)p_t^\ast
+
(p_t^\ast-p),
\]
so that, by the triangle inequality,
\[
\|p_t-p\|_{L^q(\lambda)}
\le
\left|
\frac{1}{\gamma(t)}-1
\right|
\|p_t^\ast\|_{L^q(\Omega)}
+
\|p_t^\ast-p\|_{L^q(\lambda)}.
\]
The bounds in~\eqref{eq:p*tBound} yield
\(
\|p_t^\ast\|_{L^q(\Omega)}
\le
b\,\lambda(\Omega)^{1/q},
\)
and since
\(\gamma(t)\to1\), we obtain
\(
\left|
\frac{1}{\gamma(t)}-1
\right|
\|p_t^\ast\|_{L^q(\Omega)}
\rightarrow0.
\)

It therefore remains to prove that
\(
\|p_t^\ast-p\|_{L^q(\lambda)}
\rightarrow0.
\)
Write
\[
p_t^\ast-p
=
(\rho_t*p-p)
+
\rho_t*(f_t-p),
\]
Then, the triangle inequality gives
\[
\|p_t^\ast-p\|_{L^q(\lambda)}
\le
\|\rho_t*p-p\|_{L^q(\lambda)}
+
\|\rho_t*(f_t-p)\|_{L^q(\lambda)}.
\]
The first term converges to zero by the Approximate Identity Theorem (see fact 1 above), since
\(p\in L^q(\lambda)\). For the second term, Young's convolution inequality (fact 2 above) yields
\[
\|\rho_t*(f_t-p)\|_{L^q(\lambda)}
\le
\|\rho_t\|_{L^1(\lambda)}
\|f_t-p\|_{L^q(\lambda)}
=
\|f_t-p\|_{L^q(\lambda)}.
\]
Finally,
\(
f_t-p
=
a\mathbf 1_{\Omega_t\setminus\Omega},
\)
so
\[
\|f_t-p\|_{L^q(\lambda)}
=
a
\lambda(\Omega_t\setminus\Omega)^{1/q}
\rightarrow0,
\]
because \(\Omega\) is a bounded rectangle. Hence
\(
\|p_t-p\|_{L^q(\lambda)}\to0.
\)
Taking \(p_\varepsilon:=p_t\) with \(t=\varepsilon\) proves the lemma.
\end{proof}}

\begin{lemma}
    \label{lemma:uniform-convergence-of-kernel}
    Let \(\Omega=\Omega_X\times\Omega_Z\subset\mathbb R^{\ell_X+\ell_Z}\) be a bounded
open rectangle. Let \(P\) be a law for \((X,Z)\), dominated by Lebesgue measure on
\(\Omega\), with density \(p_{X,Z}\) satisfying
\[
C^{-1}<p_{X,Z}(x,z)<C
\qquad
\text{for Lebesgue-a.e. }(x,z)\in\Omega
\]
for some \(C>1\). Define
\[
k(x,z)
:=
\frac{p_{X,Z}(x,z)}{p_X(x)p_Z(z)}.
\]
Suppose that
\(
k\in W^{s,2}(\Omega)
\)
for some integer \(s>2(\ell_X+\ell_Z)+4\).

Let
\(
(\sigma_i,u_i,v_i)_{i\ge1}
\)
be a singular system of \(\mathcal{E}_P\), ordered so that
\(
\sigma_1\ge\sigma_2\ge\cdots>0,
\)
with singular values repeated according to multiplicity, and write \(d\coloneqq \ell_X+\ell_Z\). Then the following hold.
\begin{enumerate}
    \item \label{lemma:uniform-convergence-of-kernel:decay}
    The singular values satisfy
    \(
    \sigma_i=O(i^{-s/d}).
    \)

    \item \label{lemma:uniform-convergence-of-kernel:sup-norms}
    For every integer \(\ell\) with \(d/2<\ell<s/4\), there is a finite constant \(A\), not depending on \(i\), such that
    \[
    \|u_i\|_{L^\infty(P_X)}
    \le
    A\sigma_i^{-\ell/s}
    \qquad\text{and}\qquad
    \|v_i\|_{L^\infty(P_Z)}
    \le
    A\sigma_i^{-\ell/s}
    \]
    for every \(i\ge1\). In particular \(\ell/s<1/4\), so that
    \(
    \sigma_i\|u_i\|_{L^\infty(P_X)}\|v_i\|_{L^\infty(P_Z)}
    \le
    A^2\sigma_i^{1-2\ell/s}
    \)
    with \(1-2\ell/s>1/2\).

    \item \label{lemma:uniform-convergence-of-kernel:summability}
    For every such \(\ell\),
    \[
    \sum_{i=1}^\infty\sigma_i^{1-2\ell/s}<\infty,
    \qquad\text{and consequently}\qquad
    \sum_{i=1}^\infty\sigma_i\|u_i\|_{L^\infty(P_X)}\|v_i\|_{L^\infty(P_Z)}<\infty.
    \]

    \item \label{lemma:uniform-convergence-of-kernel:uniform-convergence}
    The truncated expansions of \(k\) converge uniformly, i.e.
    \[
    \left\|
    k-\sum_{i=1}^j\sigma_i u_i v_i
    \right\|_{L^\infty(P_X\otimes P_Z)}
    \longrightarrow0
    \]
    as \(j\rightarrow\infty\).

    \item \label{lemma:uniform-convergence-of-kernel:shifted-summability}
    Write \(M_i\coloneqq\max\{\|u_i\|_{L^\infty(P_X)},\|v_i\|_{L^\infty(P_Z)}\}\). If, in addition, either
    \(
    \sup_{i\ge1}\sigma_i/\sigma_{i+1}<\infty
    \)
    or
    \(
    \sup_{i\ge1}M_i<\infty,
    \)
    then the index-shifted sum is also finite:
    \[
    \sum_{i=1}^\infty\sigma_i\|u_i\|_{L^\infty(P_X)}\|v_{i+1}\|_{L^\infty(P_Z)}<\infty.
    \]
\end{enumerate}
Moreover, since \(p_X\) and \(p_Z\) are bounded above and below on their supports,
the convergence in \cref{lemma:uniform-convergence-of-kernel:uniform-convergence} also holds in \(L^\infty(\Omega)\), after choosing appropriate versions
of the singular functions.

\end{lemma}

\begin{proof}
We establish \cref{lemma:uniform-convergence-of-kernel:sup-norms} first, then \cref{lemma:uniform-convergence-of-kernel:decay}, and deduce \cref{lemma:uniform-convergence-of-kernel:summability,lemma:uniform-convergence-of-kernel:uniform-convergence,lemma:uniform-convergence-of-kernel:shifted-summability} from the two. Recall that \(d=\ell_X+\ell_Z\). Since \(p_{X,Z}\) is bounded above and below on the bounded
rectangle \(\Omega_X\times\Omega_Z\), the marginal densities \(p_X\) and \(p_Z\) are also
bounded above and below on their supports. Indeed, for Lebesgue-a.e. \(x\in\Omega_X\),
\[
C^{-1}\lambda_Z(\Omega_Z)
\le
p_X(x)
=
\int_{\Omega_Z}p_{X,Z}(x,z)\,d\lambda_Z(z)
\le
C\lambda_Z(\Omega_Z),
\]
and similarly, for Lebesgue-a.e. \(z\in\Omega_Z\),
\[
C^{-1}\lambda_X(\Omega_X)
\le
p_Z(z)
\le
C\lambda_X(\Omega_X).
\]
Since \(p_X\) and \(p_Z\) are bounded above and below by positive constants on their
supports, the weighted \(L^2\) norms induced by \(P_X\), \(P_Z\), and
\(P_X\otimes P_Z\) are equivalent to the corresponding Lebesgue \(L^2\) norms. More
precisely, there exist constants \(0<c'<C'<\infty\),\(0<c''<C''<\infty\) and \(0<c'''<C'''<\infty\) such that, for every measurable
\(f:\Omega_X\to\mathbb R\), \(g:\Omega_Z\to\mathbb R\), and \(a:\Omega_X\times\Omega_Z\to\mathbb R\),
\[
c'\|f\|_{L^2(\lambda_X)}
\le
\|f\|_{L^2(P_X)}
\le
C'\|f\|_{L^2(\lambda_X)},
\]
\[
c''\|g\|_{L^2(\lambda_Z)}
\le
\|g\|_{L^2(P_Z)}
\le
C''\|g\|_{L^2(\lambda_Z)},
\]
and
\[
c'''\|a\|_{L^2(\lambda_X\otimes\lambda_Z)}
\le
\|a\|_{L^2(P_X\otimes P_Z)}
\le
C'''\|a\|_{L^2(\lambda_X\otimes\lambda_Z)}.
\]

Thus the \(L^2(P_X)\), \(L^2(P_Z)\), and Lebesgue \(L^2\) norms on the corresponding
bounded supports are equivalent.

The operator \(\mathcal{E}_P\) has kernel \(k\), i.e.
\[
(\mathcal{E}_P f)(z)
=
\int_{\Omega_X}
k(x,z)f(x)p_X(x)\,d\lambda_X(x).
\]
Since \(k\in L^2(P_X\otimes P_Z)\), \(\mathcal{E}_P\) is Hilbert--Schmidt and hence compact.
Therefore it admits a singular system \(
(\sigma_i,u_i,v_i)_{i\ge1}
\) satisfying
\[
\mathcal{E}_P u_i=\sigma_i v_i,
\qquad
\mathcal{E}_P^\ast v_i=\sigma_i u_i.
\]
The Schmidt expansion gives
\[
\sum_{i=1}^\infty \sigma_i u_i(x)v_i(z)
=
k(x,z)
\]
with convergence in \(L^2(P_X\otimes P_Z)\).

We next obtain uniform bounds on the singular functions. For a function \(a\) of \(Z\),
write \(\widetilde a(x,z):=a(z)\). Since
\(
v_i=\sigma_i^{-1}\mathcal{E}_Pu_i,
\)
we have
\[
\|\widetilde v_i\|_{W^{s,2}(\Omega)}
=
\sigma_i^{-1}\|\widetilde{\mathcal{E}_Pu_i}\|_{W^{s,2}(\Omega)}.
\]
We claim that
\begin{equation}
\label{eq:bound_derivative}
   \|\widetilde{\mathcal{E}_Pu_i}\|_{W^{s,2}(\Omega)}
\le
M
\end{equation}
for a constant \(M<\infty\) not depending on \(i\). To see this, note first that
\(\widetilde{\mathcal{E}_Pu_i}\) depends only on \(z\). Hence its weak derivatives in the \(x\)
coordinates vanish. For a multi-index \(\alpha=(\alpha_X,\alpha_Z)\), only the terms
with \(\alpha_X=0\) contribute to the computation of the weak derivative of \(\widetilde{\mathcal{E}_Pu_i}\).
For such \(\alpha\), we claim that the weak derivative of
\(\widetilde{\mathcal{E}_Pu_i}\) is given by
\[
D^\alpha \widetilde{\mathcal{E}_Pu_i}(x,z)
=
\int_{\Omega_X}
D_z^{\alpha_Z}k(x',z)u_i(x')p_X(x')\,d\lambda_X(x').
\]
Indeed, let \(\varphi\in C_c^\infty(\Omega)\) be an arbitrary test function. Since
\(\alpha_X=0\), the derivative \(D^\alpha\) acts only on the
\(z\)-coordinates, i.e. \(D^\alpha\varphi = D_z^{\alpha_Z}\varphi\). Furthermore,
\(D_z^{\alpha_Z}k\in L^2(\Omega)\) because
\(k\in W^{s,2}(\Omega)\). Write \(g_i\coloneqq u_ip_X\), which belongs to \(L^2(\lambda_X)\) because \(u_i \in L^2(P_X)\) and \(p_X\) is bounded, and let \(\{g_{i,n}\}_{n\ge1}\subseteq C_c^\infty(\Omega_X)\) satisfy \(g_{i,n}\rightarrow g_i\) in \(L^2(\lambda_X)\). For each fixed \(x\in\Omega_X\), the function \((x',z)\mapsto g_{i,n}(x')\varphi(x,z)\) belongs to \(C_c^\infty(\Omega)\), and its weak derivative of order \((0,\alpha_Z)\) in \((x',z)\) is \(g_{i,n}(x')D_z^{\alpha_Z}\varphi(x,z)\). Applying the definition of the weak derivative \(D_z^{\alpha_Z}k\) of \(k\) with this test function gives
\[
\int_\Omega
k(x',z)g_{i,n}(x')
D^\alpha\varphi(x,z)
\,d\lambda(x',z)
=
(-1)^{|\alpha|}
\int_\Omega
D_z^{\alpha_Z}k(x',z)g_{i,n}(x')
\varphi(x,z)
\,d\lambda(x',z)
\]
for each \(x \in \Omega_X\). We integrate both sides over \(x\in\Omega_X\) and let \(n\rightarrow\infty\). By the Cauchy--Schwarz inequality, for \(a\in\{k, D_z^{\alpha_Z}k\}\) and \(b \in \{D^\alpha\varphi, \varphi\}\),
\[
\left|\int_{\Omega_X}\int_\Omega a(x',z)\{g_{i,n}(x')-g_i(x')\}b(x,z)\,d\lambda(x',z)\,d\lambda_X(x)\right|
\le
\lambda_X(\Omega_X)\lambda_Z(\Omega_Z)^{1/2}\|b\|_{L^\infty(\Omega)}\|a\|_{L^2(\Omega)}\|g_{i,n}-g_i\|_{L^2(\lambda_X)}
\rightarrow0,
\]
so that the identity holds with \(g_{i,n}\) replaced by \(g_i\). Applying Fubini's theorem to each side then gives
\begin{align*}
\int_\Omega
\widetilde{\mathcal{E}_Pu_i}(x,z)
D^\alpha\varphi(x,z)\,
d\lambda(x,z)
&=
(-1)^{|\alpha|}
\int_\Omega
\left\{
\int_{\Omega_X}
D_z^{\alpha_Z}k(x',z)
u_i(x')p_X(x')
\,d\lambda_X(x')
\right\}\\
&\qquad\times
\varphi(x,z)
\,d\lambda(x,z),
\end{align*}
which is precisely the defining identity of the weak derivative. Hence
\[
D^\alpha \widetilde{\mathcal{E}_Pu_i}(x,z)
=
\int_{\Omega_X}
D_z^{\alpha_Z}k(x',z)
u_i(x')p_X(x')
\,d\lambda_X(x')
\]
in the weak sense.

By Cauchy--Schwarz,
\begin{align*}
\left|
\int_{\Omega_X}
D_z^{\alpha_Z}k(x',z)u_i(x')p_X(x')\,d\lambda_X(x')
\right|^2
&\le
\left\{
\int_{\Omega_X}
|D_z^{\alpha_Z}k(x',z)|^2p_X(x')\,d\lambda_X(x')
\right\} \\
&\qquad
\left\{
\int_{\Omega_X}
u_i^2(x')p_X(x')\,d\lambda_X(x')
\right\}.
\end{align*}
Because \(\|u_i\|_{L^2(P_X)}=1\), and because \(p_X\) is bounded above, integration
over \(z\) gives
\(
\|D^\alpha \widetilde{\mathcal{E}_Pu_i}\|_{L^2(\Omega)}
\le
M_\alpha
\|D_z^{\alpha_Z}k\|_{L^2(\Omega)}
\)
for a finite constant \(M_\alpha\) not depending on \(i\). Summing over
\(|\alpha|\le s\), we obtain the bound ~\eqref{eq:bound_derivative}
for some finite \(M\) independent of \(i\). Hence
\[
\|\widetilde v_i\|_{W^{s,2}(\Omega)}
\le
M\sigma_i^{-1}.
\]

We now combine Sobolev embedding with interpolation. Since \(\Omega\) is a bounded
rectangle, it satisfies the cone condition (\cite{adams_sobolev_2003-1} Definition 4.6). Therefore, for any integer \(\ell>d/2\),
the Sobolev embedding theorem (Theorem 4.12 of \cite{adams_sobolev_2003-1}) gives
\[
\|\widetilde v_i\|_{L^\infty(\Omega)}
\le
A
\|\widetilde v_i\|_{W^{\ell,2}(\Omega)}.
\]
Choose an integer \(\ell\) such that
\[
\frac d2<\ell<\frac s4.
\]
This is possible because \(s>2d+4\). By the interpolation inequality for Sobolev norms (Theorem 5.2 of \cite{adams_sobolev_2003-1}), there exists a constant $A'$ such that
\[
\|\widetilde v_i\|_{W^{\ell,2}(\Omega)}
\le
A'
\|\widetilde v_i\|_{W^{s,2}(\Omega)}^{\ell/s}
\|\widetilde v_i\|_{L^2(\Omega)}^{1-\ell/s}.
\]
The \(L^2(\Omega)\)-norm of \(\widetilde v_i\) is uniformly bounded in \(i\), because
\[
\int_{\Omega_Z}v_i^2(z)\,d\lambda_Z(z)
\le
C \,\lambda_X(\Omega_X)^{-1}
\int_{\Omega_Z}v_i^2(z)p_Z(z)\,d\lambda_Z(z)
=
C \,\lambda_X(\Omega_X)^{-1}.
\]
Thus, there exists $A''>0$ such that
\(
\|\widetilde v_i\|_{L^\infty(\Omega)}
\le
A''
\sigma_i^{-\ell/s}.
\)
Since \(\widetilde v_i(x,z)=v_i(z)\), this gives
\[
\|v_i\|_{L^\infty(P_Z)}
\le
A''
\sigma_i^{-\ell/s}.
\]
Applying the same argument to the adjoint operator \(\mathcal{E}_P^\ast\) gives that for some constant $A'''$
\[
\|u_i\|_{L^\infty(P_X)}
\le
A'''
\sigma_i^{-\ell/s}.
\]

Therefore, for some constant $A''''$
\[
\sigma_i
\|u_i\|_{L^\infty(P_X)}
\|v_i\|_{L^\infty(P_Z)}
\le
A''''
\sigma_i^{1-2\ell/s}.
\]
Since \(\ell<s/4\), we have
\[
1-\frac{2\ell}{s}>\frac12.
\]

We next invoke a classical result of Birman and Solomyak on the decay of singular
values of integral operators with Sobolev kernels. Recall that
\(\Omega=\Omega_X\times\Omega_Z\subset\mathbb R^d\), where
\(d=\ell_X+\ell_Z\), is a bounded rectangle, and in particular a bounded Lipschitz domain.

The Birman--Solomyak theorem states (\cite{griebel_decay_2018}, \cite{birman_estimates_1977} Section 2, Proposition 2.1) that if the kernel of an integral operator belongs
to \(W^{s,2}(\Omega)\), then the singular values of the operator decay
polynomially at a rate determined by the smoothness \(s\) and the ambient dimension
\(d\). More precisely, if
\[
(Tf)(z)
=
\int_{\Omega_X}
k(x,z)f(x)\,d\lambda_X(x),
\]
with
\(k\in W^{s,2}(\Omega)\),
then the singular values of $T$ decay as
\(
\sigma_i(T)
=
O(i^{-s/d}).
\)
This result applies to \(T\) with kernel \(k\), which acts between the Lebesgue spaces \(L^2(\lambda_X)\) and \(L^2(\lambda_Z)\). The conditional expectation operator factors as
\(
\mathcal{E}_P = J_Z\,T\,M_{p_X},
\)
where \(M_{p_X}:L^2(P_X)\rightarrow L^2(\lambda_X)\) is multiplication by \(p_X\) and \(J_Z:L^2(\lambda_Z)\rightarrow L^2(P_Z)\) is the identity map. Both are bounded because \(p_X\) and \(p_Z\) are bounded above. Since \(\sigma_i(J_ZTM_{p_X})\le\|J_Z\|\,\|M_{p_X}\|\,\sigma_i(T)\) for every \(i\), this yields
\[
\sigma_i
=
O(i^{-s/d}).
\]
Since \(s>2d+4\), we have \(s/d>2\), and therefore
\(
\sigma_i=o(i^{-2}).
\)
Consequently,
\[
\sigma_i^{\,1-2\ell/s}
=
O\!\left(
i^{-(s/d)(1-2\ell/s)}
\right).
\]
Because \(\ell<s/4\), we have
\(
1-\frac{2\ell}{s}>\frac12,
\)
and hence
\[
\frac{s}{d}\left(1-\frac{2\ell}{s}\right)
>
\frac{s}{2d}
>
1.
\]
It follows that
\[
\sum_{i=1}^\infty
\sigma_i^{\,1-2\ell/s}
<
\infty.
\]
Therefore,
\[
\sum_{i=1}^\infty
\sigma_i
\|u_i\|_{L^\infty(P_X)}
\|v_i\|_{L^\infty(P_Z)}
<
\infty.
\]

By the bounds established above,
\[
\left|
\sigma_i u_i(x)v_i(z)
\right|
\le
\sigma_i
\|u_i\|_{L^\infty(P_X)}
\|v_i\|_{L^\infty(P_Z)}
=:\kappa_i,
\]
where
\(
\sum_{i=1}^\infty \kappa_i<\infty.
\)
Therefore, the Weierstrass \(M\)-test implies that the series
\(
\sum_{i=1}^\infty
\sigma_i u_i(x)v_i(z)
\)
converges absolutely and uniformly on \(\Omega\) to a bounded function \(k'\). In
particular,
\[
\left\|
k'
-
\sum_{i=1}^j
\sigma_i u_i v_i
\right\|_{L^\infty(P_X\otimes P_Z)}
\rightarrow0 \,\, \text{as} \,\, j \rightarrow \infty.
\]
Since \(\Omega\) has finite Lebesgue measure, convergence in \(L^\infty\) implies
convergence in \(L^2(P_X\otimes P_Z)\). Thus,
\[
\left\|
k'
-
\sum_{i=1}^j
\sigma_i u_i v_i
\right\|_{L^2(P_X\otimes P_Z)}
\rightarrow0 \,\, \text{as} \,\, j \rightarrow \infty.
\]

On the other hand, the Schmidt expansion of the Hilbert--Schmidt kernel \(k\) states
that the same sequence of partial sums converges to \(k\) in
\(L^2(P_X\otimes P_Z)\). Since limits in the Hilbert space
\(L^2(P_X\otimes P_Z)\) are unique, we have
\[
k=k'
\qquad
P_X\otimes P_Z\text{-almost everywhere}.
\]
Replacing \(k\) by this equivalent version, we conclude that
\[
\left\|
k-
\sum_{i=1}^j
\sigma_i u_i v_i
\right\|_{L^\infty(P_X\otimes P_Z)}
\longrightarrow0,
\]
as claimed. This proves \cref{lemma:uniform-convergence-of-kernel:uniform-convergence}.

It remains to prove \cref{lemma:uniform-convergence-of-kernel:shifted-summability}. Suppose first that \(R\coloneqq\sup_{i\ge1}\sigma_i/\sigma_{i+1}<\infty\). By the bounds of \cref{lemma:uniform-convergence-of-kernel:sup-norms},
\[
\sigma_i\|u_i\|_{L^\infty(P_X)}\|v_{i+1}\|_{L^\infty(P_Z)}
\le
A^2\sigma_i^{1-2\ell/s}
\left(\frac{\sigma_i}{\sigma_{i+1}}\right)^{\ell/s}
\le
A^2R^{\ell/s}\sigma_i^{1-2\ell/s},
\]
which is summable by \cref{lemma:uniform-convergence-of-kernel:summability}. Suppose instead that \(\bar M\coloneqq\sup_{i\ge1}M_i<\infty\). Then
\(
\sigma_i\|u_i\|_{L^\infty(P_X)}\|v_{i+1}\|_{L^\infty(P_Z)}
\le
\bar M^2\sigma_i,
\)
and \(\sum_{i=1}^\infty\sigma_i<\infty\) because \(\sigma_i=O(i^{-s/d})\) by \cref{lemma:uniform-convergence-of-kernel:decay} with \(s/d>2\). In either case the shifted sum converges.
This proves the lemma.
\end{proof}

\begin{lemma}
    \label{lemma:in-range-uniformly-bounded}
    Suppose that \(P\) is a distribution for \(O=(W,X,Z, D)\) such that the joint
distribution of \((X,Z)\) conditional on $D =d$ is dominated by Lebesgue measure on a bounded rectangle
\(\Omega_X\times\Omega_Z\) for each $d \in \textsf{Supp}[D;P]$, with density \(p_{X,Z|d}\). Suppose that, for some
constant \(C>1\),
\[
C^{-1}<p_{X,Z|d}(x,z)<C
\qquad
\text{Lebesgue-a.e. on }\Omega_X\times\Omega_Z .
\]
If
\(
g\in \mathscr{R}(\mathcal{E}_{P;d})
\) and
\(
g^\ast\in \mathscr{R}(\mathcal{E}_{P;d}^\ast),
\)
then \(g\) has an essentially bounded version as a function of \(Z\), and
\(g^\ast\) has an essentially bounded version as a function of \(X\). That is,
\[
\|g\|_{L^\infty(P_{Z|d})}<\infty,
\qquad
\|g^\ast\|_{L^\infty(P_{X|d})}<\infty.
\]
\end{lemma}

\begin{proof}
Because \(g\in \mathscr{R}(\mathcal{E}_P)\), there exists \(f\in L^2(P_X)\) such that
\(
g=\mathcal{E}_P f.
\)
For \(P_Z\)-a.e. \(z\),
\[
g(z)
=
\int_{\Omega_X}
\frac{p_{X,Z}(x,z)}{p_X(x)p_Z(z)}
f(x)p_X(x)\,d\lambda_X(x).
\]
By Cauchy--Schwarz,
\[
|g(z)|
\le
\left[
\int_{\Omega_X}
\left\{
\frac{p_{X,Z}(x,z)}{p_X(x)p_Z(z)}
\right\}^2
p_X(x)\,d\lambda_X(x)
\right]^{1/2}
\|f\|_{L^2(P_X)}.
\]
Since \(p_{X,Z}\) is bounded above and below on the bounded rectangle
\(\Omega_X\times\Omega_Z\), the marginal densities \(p_X\) and \(p_Z\) are also
bounded above and below by positive constants on their supports. Hence the kernel
\[
k(x,z):=\frac{p_{X,Z}(x,z)}{p_X(x)p_Z(z)}
\]
is essentially bounded. Therefore,
\[
\sup_z
\int_{\Omega_X} k(x,z)^2p_X(x)\,d\lambda_X(x)
<\infty,
\]
and consequently
\(
\|g\|_{L^\infty(P_Z)}<\infty.
\)

A similar argument shows $\|g^*\|_{L^\infty(P_X)} < \infty$.
This proves the lemma.
\end{proof}

\begin{lemma}
    \label{lemma:product-in-l2}
    Let $P$ be a distribution for $O = (W, X, Z, D)$ satisfying \cref{as:dominating-measure-support}, and suppose that $L^{-1} < p(o) < L$ $P$-a.e. for some constant $L \geq 1$. Then, for every $h \in L^2(P_{X, D})$ and every $q \in L^2(P_{Z, D})$, the product $h(X, D)q(Z, D)$ belongs to $L^2(P)$, with
    \begin{align*}
        \|hq\|_{L^2(P)} \leq L^{7/2}\left|\textsf{Supp}[P_D]\right|^{1/2}\|h\|_{L^2(P_{X, D})}\|q\|_{L^2(P_{Z, D})}.
    \end{align*}
    In particular, this constant depends on $P$ only through $L$ and the number of levels of $D$, and so may be taken common to any family of such laws sharing these.
\end{lemma}

\begin{proof}[Proof of \cref{lemma:product-in-l2}]
    We take the convention of \cref{remark:marginal-conditional-density-bounds} to avoid tracking unnecessary constants. By that remark, the conditional density $p_{X, Z|d}$ is bounded above by $L^2$, and the conditional densities $p_{X|d}$ and $p_{Z|d}$ are bounded below by $L^{-2}$, uniformly over $d \in \textsf{Supp}[P_D]$. Hence
    \begin{align*}
        p_{X, Z|d}(x, z) \leq L^{6}p_{X|d}(x)p_{Z|d}(z)
    \end{align*}
    for every $d \in \textsf{Supp}[P_D]$ and almost every $(x, z)$. Consequently,
    \begin{align*}
        \|hq\|^2_{L^2(P)} =& \sum_{d \in \textsf{Supp}[P_D]}P(D=d)\int \int h(x, d)^2q(z, d)^2p_{X, Z|d}(x, z)d\lambda_X(x)d\lambda_Z(z)\\
        \leq& L^{6}\sum_{d \in \textsf{Supp}[P_D]}P(D=d)\|h(\cdot, d)\|^2_{L^2(P_{X|d})}\|q(\cdot, d)\|^2_{L^2(P_{Z|d})}\\
        \leq& L^{6}\left\{\max_{d \in \textsf{Supp}[P_D]}\|q(\cdot, d)\|^2_{L^2(P_{Z|d})}\right\}\sum_{d \in \textsf{Supp}[P_D]}P(D=d)\|h(\cdot, d)\|^2_{L^2(P_{X|d})}\\
        =& L^{6}\left\{\max_{d \in \textsf{Supp}[P_D]}\|q(\cdot, d)\|^2_{L^2(P_{Z|d})}\right\}\|h\|^2_{L^2(P_{X, D})},
    \end{align*}
    where the maximum is attained because $D$ is finite discrete. Now
    \begin{align*}
        \max_{d \in \textsf{Supp}[P_D]}\|q(\cdot, d)\|^2_{L^2(P_{Z|d})} \leq \left\{\min_{d \in \textsf{Supp}[P_D]}P(D=d)\right\}^{-1}\|q\|^2_{L^2(P_{Z, D})},
    \end{align*}
    and under the normalization of \cref{remark:marginal-conditional-density-bounds} the reference measure $\lambda_D$ is the counting measure on $\textsf{Supp}[P_D]$ divided by $\left|\textsf{Supp}[P_D]\right|$. Since $p_D > L^{-1}$ by that same remark, it follows that
    \begin{align*}
        P(D=d) = p_D(d)\lambda_D(\{d\}) > \frac{L^{-1}}{\left|\textsf{Supp}[P_D]\right|}
    \end{align*}
    for every $d \in \textsf{Supp}[P_D]$. Combining the three displays,
    \begin{align*}
        \|hq\|^2_{L^2(P)} \leq L^{7}\left|\textsf{Supp}[P_D]\right|\|h\|^2_{L^2(P_{X, D})}\|q\|^2_{L^2(P_{Z, D})},
    \end{align*}
    which is the claim.
\end{proof}

{ \begin{lemma}
\label{lemma:l1-bounded-functional}
Suppose that \(P_0\) satisfies \cref{as:continuous-linear-maps},  \cref{as:dominating-measure-support}, and that
\(m,\widetilde m,P_0\) satisfy \cref{as:m-m-tilde-bounded}. Let \(P\) be a probability law for
\(O=(W,X,Z,D)\) dominated by the same product measure \(\lambda\) as \(P_0\), and
suppose that \(P\) and \(P_0\) are supported on the same bounded rectangle \(\Omega\). Let \(p\) and \(p_0\) denote the corresponding densities
with respect to \(\lambda\). Assume that for some constant \(C>1\),
\[
C^{-1}<p(o)<C
\qquad
\lambda\text{-a.e. on the common support.}
\]
Then \(P\) and \(P_0\) are mutually absolutely continuous. Moreover, \(P\) satisfies \cref{as:continuous-linear-maps}. Moreover, the Riesz representers at \(P\) admit
essentially bounded versions:
\[
\|\alpha_P\|_{L^\infty(P_{X,D})}<\infty,
\qquad
\|r_P\|_{L^\infty(P_{Z,D})}<\infty .
\]
\end{lemma}

\begin{proof}[Proof of \cref{lemma:l1-bounded-functional}]
The density bounds on $p_0$ and $p$ imply that $P$ and $P_0$ are mutually absolutely continuous since $p$ and $p_0$ are strictly positive and finite  $\lambda$-a.e. on the common support.
We prove the claim for \(\alpha_P\). The argument for \(r_P\) is identical.

By~\cref{as:m-m-tilde-bounded}, there exists \(L_P<\infty\) such that, for every
\(h\in L^1(P_{X,D;0})\),
\[
\|\widetilde m(W,X,D;h)\|_{L^1(P_{W,X,D;0})}
\le
L_P\|h\|_{L^1(P_{X,D;0})} .
\]
Since \(p_0\) and \(p\) are both bounded above and below by positive constants on the
same bounded rectangle, the marginal densities \(p_{0,X,D}\) and \(p_{X,D}\) are also
bounded above and below by positive constants on their common support. Hence the
\(L^1(P_{X,D;0})\) and \(L^1(P_{X,D})\) norms are equivalent: there exists
\(A<\infty\) such that
\[
\|h\|_{L^1(P_{X,D;0})}
\le
A\|h\|_{L^1(P_{X,D})}
\]
for all measurable \(h \in L^1(P_{X,D})\).

Similarly, the marginal densities of \((W,X,D)\) under \(P_0\) and \(P\) are bounded
above and below by positive constants on their common support. Therefore, for some
\(B<\infty\),
\[
\frac{p_{W,X,D}(w,x,d)}{p_{0,W,X,D}(w,x,d)}
\le B
\qquad
P_0\text{-a.e.}
\]
Consequently, for every \(h\in L^1(P_{X,D})\),
\begin{align*}
\left|
E_P\{\widetilde m(W,X,D;h)\}
\right|
&\le
E_P\{|\widetilde m(W,X,D;h)|\}
\\
&=
E_{P_0}\left[
\frac{p_{W,X,D}(W,X,D)}{p_{0,W,X,D}(W,X,D)}
|\widetilde m(W,X,D;h)|
\right]
\\
&\le
B
E_{P_0}\{|\widetilde m(W,X,D;h)|\}
\\
&\le
BL_P\|h\|_{L^1(P_{X,D;0})}
\\
&\le
BL_PA\|h\|_{L^1(P_{X,D})} .
\end{align*}
Thus the map
\[
h\mapsto E_P\{\widetilde m(W,X,D;h)\}
\]
is a bounded linear functional on \(L^1(P_{X,D})\).

By the dual representation of bounded linear functionals on \(L^1\), there exists
\(\alpha_P\in L^\infty(P_{X,D})\) such that
\[
E_P\{\widetilde m(W,X,D;h)\}
=
E_P\{h(X,D)\alpha_P(X,D)\}
\]
for every \(h\in L^1(P_{X,D})\). Since \(P\) is a probability measure,
\(L^2(P_{X,D})\subset L^1(P_{X,D})\), and the same identity holds for every
\(h\in L^2(P_{X,D})\). Moreover,
\[
\left|
E_P\{\widetilde m(W,X,D;h)\}
\right|
\le
\|\alpha_P\|_{L^\infty(P_{X,D})}
\|h\|_{L^1(P_{X,D})}
\le
\|\alpha_P\|_{L^\infty(P_{X,D})}
\|h\|_{L^2(P_{X,D})} .
\]
Hence the functional is also bounded on \(L^2(P_{X,D})\). Since
\(\alpha_P\in L^\infty(P_{X,D})\subset L^2(P_{X,D})\), this is precisely the Riesz
representer required in \cref{as:continuous-linear-maps}.

The same argument, replacing \(\widetilde m,\alpha_P,(X,D)\) by
\(m,r_P,(Z,D)\), shows that
\(
\|r_P\|_{L^\infty(P_{Z,D})}<\infty
\)
and that the second functional in \cref{as:continuous-linear-maps} is bounded on \(L^2(P_{Z,D})\). This
proves the lemma.
\end{proof}
}

\begin{lemma}
    \label{lemma:weak-convergence-of-representers}
    Suppose $P_0$ satisfies \cref{as:continuous-linear-maps},  \cref{as:dominating-measure-support}, and { $\widetilde{m}, m, P_0$ satisfy that for every $h \in L^2(P_{X,D;0}), q \in L^2(P_{Z,D;0})$, $\widetilde{m}(W, X, D;h) \in L^1(P_{W,X,D;0}), m(W, Z, D;q) \in L^1(P_{W,Z,D;0})$, which holds for example under \cref{as:m-m-tilde-bounded}}. Let $P_j$ be any sequence of distributions for $O$ satisfying \cref{as:continuous-linear-maps}, $P_j$ and $P_0$ are mutually absolutely continuous, $\lim_{j\rightarrow \infty}\|p_j-  p_0\|_{L^1(\lambda)} = 0$, and there exists $C \geq 1$ such that $C^{-1} < p_j < C $ a.e. $P_j$. Then, for every $h \in L^2(P_{X,D;0})$ and $q \in L^2(P_{Z,D;0})$,
    \begin{align*}
        E_{P_j}[\widetilde{m}(W, X, D;h)] \rightarrow E_{P_0}[\widetilde{m}(W, X, D;h)], && E_{P_j}[m(W, Z, D;q)] \rightarrow E_{P_0}[m(W, Z, D;q)],
    \end{align*}
    and
    \begin{align*}
        \alpha_{P_j} \rightharpoonup \alpha_{P_0} \text{ weakly in $L^2(P_{X,D;0})$}, &&  r_{P_j} \rightharpoonup r_{P_0} \text{ weakly in $L^2(P_{Z,D;0})$}.
    \end{align*}
\end{lemma}

\begin{proof}[Proof of \cref{lemma:weak-convergence-of-representers}]
    We argue for $\alpha_{P_j}$, and by symmetry the result holds for $r_{P_j}$.

    Take $h \in L^2(P_{X,D;0})$. { Then, by assumption, $\widetilde{m}(W, X, D;h) \in L^1(P_{W,X,D;0})$}. We now argue that
    \begin{align*}
        \lim_{j\rightarrow \infty}|E_{P_j}[\widetilde{m}(W, X, D;h)] - E_{P_0}[\widetilde{m}(W, X, D;h)]| = 0.
    \end{align*}

    Now,
    \begin{align*}
        E_{P_j}[\widetilde{m}(W, X, D;h)] - E_{P_0}[\widetilde{m}(W, X, D;h)]= E_{P_0}\left[\left\{\frac{p_j(W, X, D)}{p_0(W, X, D)} - 1\right\}\widetilde{m}(W, X, D;h)\right].
    \end{align*}
    Because $p_0 > L_{P_0}^{-1}$ by \cref{as:dominating-measure-support},
    \begin{align*}
        \lim_{j\rightarrow \infty}\left\|p_{W, X, D;j}/p_{W, X, D;0} - 1\right\|_{L^1(\lambda)} = 0,
    \end{align*}
    and hence $p_{W, X, D;j}/p_{W, X, D;0} - 1$ converges to 0 in probability. Consequently $\left\{\frac{p_j(W, X, D)}{p_0(W, X, D)} - 1\right\}\widetilde{m}(W, X, D;h)$ converges in probability to $0$. Furthermore, because $C^{-1} < p_j < C$ for all $j$,  $|p_{W, X, D;j}/p_{W, X, D;0} - 1| < C{L_{P_0}} + 1$ for all $j$, the sequence $\left|\left\{\frac{p_j(W, X, D)}{p_0(W, X, D)} - 1\right\}\widetilde{m}(W, X, D;h)\right|$, $j \geq 1$ is uniformly $P_0$-integrable. Then, by Vitali's Convergence Theorem, we have that
    $\left\{\frac{p_j(W, X, D)}{p_0(W, X, D)} - 1\right\}\widetilde{m}(W, X, D;h)$ converges in $L^1(P_0)$ to $0$. But,
    \begin{align*}
        0 &\leq \limsup_{j\rightarrow \infty }|E_{P_j}[\widetilde{m}(W, X, D;h)] - E_{P_0}[\widetilde{m}(W, X, D;h)]| \\
        &\leq \limsup_{j \rightarrow \infty }E_{P_0}\left[\left|\left\{\frac{p_j(W, X, D)}{p_0(W, X, D)} - 1\right\}\widetilde{m}(W, X, D;h)\right|\right] = 0.
    \end{align*}

    Consequently, for each $h \in L^2(P_{X,D;0})$,
    \begin{align*}
        0 =&\lim_{j\rightarrow \infty}|E_{P_j}[\widetilde{m}(W, X, D;h)] - E_{P_0}[\widetilde{m}(W, X, D;h)]|\\
        =& \lim_{j\rightarrow \infty}|E_{P_j}[\alpha_{P_j}(X, D)h(X, D)] - E_{P_0}[\alpha_{P_0}(X, D)h(X, D)]|\\
        =&\lim_{j\rightarrow \infty}\left|E_{P_0}\left[\frac{p_j(X, D)}{p_0(X, D)}\alpha_{P_j}(X, D)h(X, D)\right] - E_{P_0}[\alpha_{P_0}(X, D)h(X, D)]\right|\\
        =& \lim_{j\rightarrow \infty}\left|\left\langle\frac{p_{X, D;j}}{p_{X, D;0}}\alpha_{P_j}, h \right\rangle_{L^2(P_{X,D;0})} - \left\langle\alpha_{P_0}, h \right\rangle_{L^2(P_{X,D;0})}\right|
    \end{align*}
    i.e. $\frac{p_{X, D;j}}{p_{X, D;0}}\alpha_{P_j}\rightharpoonup \alpha_{P_0}$ in $L^2(P_{X,D;0})$.

    Because weakly convergent sequences are bounded, $\sup_{j}\left\|\frac{p_{X, D;j}}{p_{X, D;0}}\alpha_{P_j}\right\|_{L^2(P_{X,D;0})} < \infty$. Moreover, because $C^{-1}L_{P_0}^{-1}< \frac{p_{X, D;j}(X,D)}{p_{X, D;0}(X, D)} < CL_{P_0}$ and $\frac{p_{X, D;j}}{p_{X, D;0}} \rightarrow 1$, Vitali's Convergence Theorem gives that $\frac{p_{X, D;0}(X,D)}{p_{X, D;j}(X, D)}h \rightarrow h $ for each fixed $h \in L^2(P_{X,D;0})$. Hence,
    \begin{align*}
        \left\langle \alpha_{P_j}, h \right\rangle_{L^2(P_{X,D;0})} = \left\langle\frac{p_{X, D;j}}{p_{X, D;0}}\alpha_{P_j}, \frac{p_{X, D;0}}{p_{X, D;j}}h \right\rangle_{L^2(P_{X,D;0})} \rightarrow \left\langle \alpha_{P_0}, h \right\rangle_{L^2(P_{X,D;0})},
    \end{align*}
    i.e. $\alpha_{P_j} \rightharpoonup \alpha_{P_0}$.
\end{proof}

The next lemma demonstrates that injectivity of $\mc{E}_{P_0}$, a uniform bound on the norms of solutions $h_{P_t}$, and additional regularity suffice to ensure $h_{P_t}^\dagger$ converges weakly to $h_{P_0}^\dagger$ when $p_t$ converges in $L^2$ to $p_0$, and that a uniform bound on the norms of the $\beta_h$-source representers of $h_{P_t}^\dagger$ further ensures strong convergence, with or without injectivity.
\begin{lemma}
    \label{lemma:weak-convergence-of-solutions}
    Let $P_0 \in \mc{P}$ and let $\{P_t: t \in [0, \delta)\} \subseteq \mc{P}$ be a family of distributions for $O$ with $P_t = P_0$ at $t=0$, all satisfying \cref{as:dominating-measure-support} with a common constant $L_{bounded}$. For each $t \in [0, \delta)$, let $h_{P_t}$ denote any solution to the primal inverse problem $\mc{E}_{P_t}h = r_{P_t}$, which exists by \cref{as:solutions-exist}. Suppose that
    \begin{enumerate}
        \item \label{lemma:weak-convergence-of-solutions:densities} $\lim_{t\downarrow 0}\|p_t - p_0\|_{L^2(\lambda)} = 0$, equivalently, by \cref{as:dominating-measure-support}, $\lim_{t\downarrow 0}\|p_t - p_0\|_{L^2(P_0)} = 0$;
        \item \label{lemma:weak-convergence-of-solutions:bounded} there exists $C < \infty$ such that $\|h_{P_t}\|_{L^2(P_{X, D;t})} < C$ for all $t \in [0, \delta)$;
        \item \label{lemma:weak-convergence-of-solutions:representers} $r_{P_t} \rightharpoonup r_{P_0}$ in $L^2(P_{Z, D;0})$ as $t \downarrow 0$
    \end{enumerate}
    Then,
    \begin{enumerate}
        \item \label{lemma:weak-convergence-of-solutions:injective-result} If $\mc{E}_{P_0}$ is injective, then $h^\dagger_{P_t} \rightharpoonup h_{P_0}^\dagger$ in $L^2(P_{X, D;0})$ as $t \downarrow 0$.
        \item \label{lemma:weak-convergence-of-solutions:beta-bound-result} If there exist $\beta_h > 0$ and $C_{\beta_h} < \infty$ such that, for all $t \in [0, \delta)$, the minimum norm solution $\omega_{P_t}^\dagger$ to $(\mc{E}_{P_t}^*\mc{E}_{P_t})^{\beta_h/2}\omega = h_{P_t}^\dagger$ exists and satisfies $\|\omega^\dagger_{P_t}\|_{L^2(P_{X, D;t})} < C_{\beta_h}$, where $\mc{E}_{P_t}^*$ is the adjoint with respect to the inner products of $L^2(P_{X, D;t})$ and $L^2(P_{Z, D;t})$, then $h^\dagger_{P_t} \rightarrow h_{P_0}^\dagger$ in $L^2(P_{X, D;0})$ as $t \downarrow 0$.
    \end{enumerate}
\end{lemma}
\begin{remark}
    The symmetric statement, obtained by replacing $\mc{E}_{P_t}, r_{P_t}, h_{P_t}$ by $\mc{E}^*_{P_t}, \alpha_{P_t}, q_{P_t}$ throughout, holds by the same argument, because $\mc{E}^*_{P_t}$ is likewise an integral operator with the roles of $x$ and $z$ interchanged.
\end{remark}

\begin{proof}[Proof of \cref{lemma:weak-convergence-of-solutions}]
    We begin with preliminaries that apply to both parts of the Lemma.

    Throughout, we use that \cref{as:dominating-measure-support} holds at every $P_t$ with the common constant $L_{bounded}$, so that the $L^2(P_t)$ and $L^2(P_0)$ norms are equivalent, with constants not depending on $t$. Consequently the spaces $L^2(P_{X, D;t})$ coincide as sets for all $t$, differing only in their inner products, and we regard each $\mc{E}_{P_t}$ as acting on these common sets. All norms, operator norms, and convergences below are those of $L^2(P_{X, D;0})$ and $L^2(P_{Z, D;0})$. In particular, because $h^\dagger_{P_t}$ is the minimum norm solution, $\|h^\dagger_{P_t}\|_{L^2(P_{X, D;t})} \leq \|h_{P_t}\|_{L^2(P_{X, D;t})}$, so condition \cref{lemma:weak-convergence-of-solutions:bounded} supplies a constant $C' < \infty$ with $\|h^\dagger_{P_t}\|_{L^2(P_{X, D;0})} < C'$ for all $t \in [0, \delta)$, and similarly for the bound in \cref{lemma:weak-convergence-of-solutions:beta-bound-result}.

    \textbf{Convergence of $\mc{E}_{P_t}$ to $\mc{E}_{P_0}$ in operator norm.} Integrating out the remaining coordinates and applying Cauchy-Schwarz over the finite measures $\lambda_W$ and $\lambda_X$, condition \cref{lemma:weak-convergence-of-solutions:densities} implies that $p_{X, Z, D;t}$ and $p_{Z, D;t}$ converge to $p_{X, Z, D;0}$ and $p_{Z, D;0}$, respectively, in $L^2(\lambda)$ and hence, by \cref{as:dominating-measure-support}, in $L^2(P_0)$. Now, recall that
    \begin{align*}
        ({\mc{E}}_{P_t}h)(z, d) =& \int h(x, d)p_{x|z, d;t}(x|z, d)d\lambda_{X}(x)\\
        =&\int h(x, d)\frac{p_{x, z, d;t}(x, z, d)}{p_{x, d;0}(x, d)p_{z, d;t}(z, d)}p_{x, d;0}(x, d)d\lambda_X(x),
    \end{align*}
    i.e. $\mc{E}_{P_t}$ is an integral operator with kernel $\frac{p_{x, z, d;t}(x, z, d)}{p_{x, d;0}(x, d)p_{z, d;t}(z, d)}$. The above convergences then give that $\frac{p_{x, z, d;t}(x, z, d)}{p_{x, d;0}(x, d)p_{z, d;t}(z, d)}$ converges in $L^2(P_{X, D;0})\otimes L^2(P_{Z, D;0})$ to $\frac{p_{x, z, d;0}(x, z, d)}{p_{x, d;0}(x, d)p_{z, d;0}(z, d)}$, the kernel of the integral operator $\mc{E}_{P_0}$, because $p_{x, d;0}, p_{z, d;t}$ are both uniformly bounded below. Convergence in $L^2(P_{X, D;0})\otimes L^2(P_{Z, D;0})$ of kernels implies convergence in Hilbert-Schmidt norm, which in turn implies that $\mc{E}_{P_t}$ converges in operator norm to $\mc{E}_{P_0}$.

    \textbf{Proof of \cref{lemma:weak-convergence-of-solutions:injective-result}} Let $\{t_i\}_{i=1}^\infty$ be any sequence in $[0, \delta)$ with $\lim_{i\rightarrow \infty}t_i = 0$. Because $\|h^\dagger_{P_{t_i}}\|_{L^2(P_{X, D;0})} < C'$ for every $i$, every subsequence of $\{h^\dagger_{P_{t_i}}\}_{i=1}^\infty$ admits a further subsequence converging weakly in $L^2(P_{X, D;0})$. Let $\{h^\dagger_{P_{t_{i_j}}}\}_{j=1}^\infty$ be any such weakly convergent subsequence, and let $h^*_{P_0}$ denote its limit.

    By \cref{as:dominating-measure-support}, the operator $\mc{E}_{P_{0}}$ is compact (see the discussion in \cite{carrasco_chapter_2007} after Assumption A.1). Compact operators map weakly convergent sequences to strongly convergent ones, so $\mc{E}_{P_0}h^\dagger_{P_{t_{i_j}}}$ converges in $L^2(P_{Z, D;0})$ to $\mc{E}_{P_0}h_{P_0}^*$. Combined with the operator norm convergence of $\mc{E}_{P_t}$ to $\mc{E}_{P_0}$ shown above, the expansion
    \begin{align*}
        \mc{E}_{P_{t_{i_j}}}h^\dagger_{P_{t_{i_j}}} - \mc{E}_{P_0}h^*_{P_0} =  \{\mc{E}_{P_{t_{i_j}}} - {\mc{E}}_{P_{0}}\}h^\dagger_{P_{t_{i_j}}} + \mc{E}_{P_0}\{ h^\dagger_{P_{t_{i_j}}} - h_{P_0}^*\}
    \end{align*}
    then gives that $\mc{E}_{P_{t_{i_j}}}h^\dagger_{P_{t_{i_j}}}$ converges in $L^2(P_{Z, D;0})$ to $\mc{E}_{P_0}h_{P_0}^*$, the first term on the right hand side being bounded in norm by $C'\|\mc{E}_{P_{t_{i_j}}} - {\mc{E}}_{P_{0}}\|$.

    On the other hand, $\mc{E}_{P_{t_{i_j}}}h^\dagger_{P_{t_{i_j}}} = r_{P_{t_{i_j}}}$ for each $j$, and $r_{P_{t_{i_j}}} \rightharpoonup r_{P_0}$ by \cref{lemma:weak-convergence-of-solutions:representers}. Because strongly convergent sequences converge weakly to the same limit, and weak limits are unique,
    \begin{align*}
        \mc{E}_{P_{0}}h^*_{P_{0}} = r_{P_0}.
    \end{align*}

    Now, if $\mc{E}_{P_0}$ is injective, $h^\dagger_{P_0}$ is the unique solution to $\mc{E}_{P_0}h = r_{P_0}$. Consequently,  $h^*_{P_0} = h^\dagger_{P_0}$.

    The above shows that for every subsequence of $\{h^\dagger_{P_{t_i}}\}_{i=1}^\infty$, any weakly convergent subsequence converges weakly to $h^\dagger_{P_0}$, and hence that $h^\dagger_{P_{t_i}} \rightharpoonup h^\dagger_{P_0}$. Because the sequence $\{t_i\}_{i=1}^\infty$ in $[0, \delta)$ with $\lim_{i\rightarrow \infty}t_i = 0$ was arbitrary, $h^\dagger_{P_t} \rightharpoonup h^\dagger_{P_0}$ as $t \downarrow 0$.

    \textbf{Proof of \cref{lemma:weak-convergence-of-solutions:beta-bound-result}}
    Now suppose that the conditions of \cref{lemma:weak-convergence-of-solutions:beta-bound-result} hold. By equivalence of the norms of $L^2(P_{X, D;t})$ and $L^2(P_{X, D;0})$, there exists $C'_{\beta_h} < \infty$ such that $\|\omega^\dagger_{P_t}\|_{L^2(P_{X, D;0})} < C'_{\beta_h}$ for all $t \in [0, \delta)$. The adjoint $\mc{E}_{P_t}^*$, taken with respect to the inner products of $L^2(P_{X, D;t})$ and $L^2(P_{Z, D;t})$, is the conditional expectation operator $q \mapsto E_{P_t}[q(Z, D)|X, D]$, an integral operator of the same form as $\mc{E}_{P_t}$ with the roles of $x$ and $z$ interchanged. The argument given above for $\mc{E}_{P_t}$ therefore also gives that $\mc{E}_{P_t}^*$ converges in operator norm to $\mc{E}_{P_0}^*$, and hence that $\mc{E}_{P_t}^*\mc{E}_{P_t}$ converges in operator norm to $\mc{E}_{P_0}^*\mc{E}_{P_0}$. It remains to argue that $(\mc{E}_{P_t}^*\mc{E}_{P_t})^{\beta_h/2}$ converges in operator norm to $(\mc{E}_{P_0}^*\mc{E}_{P_0})^{\beta_h/2}$. In what follows, define
    \begin{align*}
        B_t \coloneqq \mc{E}_{P_t}^*\mc{E}_{P_t}, \quad B_0 \coloneqq \mc{E}_{P_0}^*\mc{E}_{P_0},
    \end{align*}
    so that $B_t$ is self-adjoint and nonnegative in $L^2(P_{X, D;t})$, and $B_0$ is compact because $\mc{E}_{P_0}$ is compact. Let $\rho_t \coloneqq p_{X, D;t}/p_{X, D;0}$ and let $M_t$ denote multiplication by $\sqrt{\rho_t}$, i.e. $M_tf = \sqrt{\rho_t}f$. Because $\|M_tf\|_{L^2(P_{X, D;0})} = \|f\|_{L^2(P_{X, D;t})}$, $M_t$ is unitary as a map from $L^2(P_{X, D;t})$ to $L^2(P_{X, D;0})$. Furthermore, as an operator on $L^2(P_{X, D;0})$, $M_t$ is self-adjoint, and because $\rho_t$ is bounded above and below by constants not depending on $t$ by \cref{as:dominating-measure-support}, $M_t$ and $M_t^{-1}$ are bounded in operator norm uniformly over $t$. Define
    \begin{align*}
        \widetilde{B}_t \coloneqq M_tB_tM_t^{-1}
    \end{align*}
    which, because $B_t$ is self-adjoint in $L^2(P_{X, D;t})$ and $M_t$ is unitary from $L^2(P_{X, D;t})$ to $L^2(P_{X, D;0})$, is self-adjoint and nonnegative in $L^2(P_{X, D;0})$. Because the functional calculus commutes with unitary conjugation,
    \begin{align*}
        M_t^{-1}\widetilde{B}_t^{\beta_h/2}M_t = B_t^{\beta_h/2}.
    \end{align*}
    Now, because $p_{X, D;t}\rightarrow p_{X, D;0}$ in $L^2(P_{X, D;0})$ and $p_{X, D;0}$ is bounded below, $\rho_t \rightarrow 1$ in $P_{X, D;0}$-probability. As $\rho_t$ is bounded above and below uniformly in $t$, dominated convergence along subsequences, together with the Urysohn subsequence principle, gives that $M_t \rightarrow I$ and $M_t^{-1} \rightarrow I$ strongly in $L^2(P_{X, D;0})$. We additionally have that $B_t$ converges in operator norm to $B_0$, and that $B_0$ is compact. Recall that if $\{S_t\}$ is a family of operators that is uniformly bounded in operator norm and converges strongly to $0$, and $K$ is compact, then $S_tK \rightarrow 0$ in operator norm. Now write
    \begin{align*}
        \widetilde{B}_t - B_0 = M_t(B_t - B_0)M_t^{-1} + (M_t - I)B_0M_t^{-1} + B_0(M_t^{-1} - I).
    \end{align*}
    The first term converges to zero in operator norm because $M_t$ and $M_t^{-1}$ are uniformly bounded and $B_t \rightarrow B_0$ in operator norm. The second converges to zero in operator norm because $(M_t - I)B_0 \rightarrow 0$ in operator norm by the above fact and $M_t^{-1}$ is uniformly bounded. For the third, because $M_t^{-1}$ and $B_0$ are self-adjoint in $L^2(P_{X, D;0})$, $\|B_0(M_t^{-1} - I)\| = \|(M_t^{-1} - I)B_0\| \rightarrow 0$ by the same fact. Hence, $\widetilde{B}_t$ converges in operator norm to $B_0$. Because $\widetilde{B}_t$ and $B_0$ are nonnegative self-adjoint operators on the same Hilbert space $L^2(P_{X, D;0})$, we may apply the continuous functional calculus (\cite{conway_course_2007-1} Theorem VIII.2.6): if $A$ is a bounded self-adjoint operator with spectrum $\sigma(A)$ and $f$ is continuous on $\sigma(A)$, then $f(A)$ is well defined, the map $f \mapsto f(A)$ is linear and multiplicative, $p(A)$ is the usual polynomial in $A$ whenever $p$ is a polynomial, and $\|f(A)\| = \sup_{x \in \sigma(A)}|f(x)|$. Because $\mc{E}_{P_t}$ is a contraction from $L^2(P_{X, D;t})$ to $L^2(P_{Z, D;t})$ and $M_t$ is unitary, $\|\widetilde{B}_t\| = \|B_t\|_{L^2(P_{X, D;t})} = \|\mc{E}_{P_t}\|^2 \leq 1$, so that $\sigma(\widetilde{B}_t) \subseteq [0, 1]$, and similarly $\sigma(B_0) \subseteq [0, 1]$. Fix $\epsilon > 0$. By the Weierstrass approximation theorem, there exists a polynomial $p$ such that $\sup_{x \in [0, 1]}|x^{\beta_h/2} - p(x)| < \epsilon$. Applying the norm identity above to the continuous function $x \mapsto x^{\beta_h/2} - p(x)$ then gives
    \begin{align*}
        \|\widetilde{B}_t^{\beta_h/2} - B_0^{\beta_h/2}\| \leq \|\widetilde{B}_t^{\beta_h/2} - p(\widetilde{B}_t)\| + \|p(\widetilde{B}_t) - p(B_0)\| + \|p(B_0) - B_0^{\beta_h/2}\| \leq 2\epsilon + \|p(\widetilde{B}_t) - p(B_0)\|.
    \end{align*}
    Because $p$ is a polynomial and $\widetilde{B}_t \rightarrow B_0$ in operator norm, $\|p(\widetilde{B}_t) - p(B_0)\| \rightarrow 0$. As $\epsilon > 0$ was arbitrary, $\widetilde{B}_t^{\beta_h/2} \rightarrow B_0^{\beta_h/2}$ in operator norm. Finally,
    \begin{align*}
        B_t^{\beta_h/2} - B_0^{\beta_h/2} = M_{t}^{-1}(\widetilde{B}_t^{\beta_h/2} - B_0^{\beta_h/2})M_t + (M_t^{-1}B_0^{\beta_h/2}M_t - B_0^{\beta_h/2}) \rightarrow 0
    \end{align*}
    in operator norm, where the first term vanishes because $M_t^{-1}$ and $M_t$ are uniformly bounded, and the second vanishes by the same argument used for $\widetilde{B}_t$, because $B_0^{\beta_h/2}$ is self-adjoint and, as $B_0$ is compact and $\beta_h > 0$, compact. Consequently, $(\mc{E}_{P_t}^*\mc{E}_{P_t})^{\beta_h/2}$ converges in operator norm to $(\mc{E}_{P_0}^*\mc{E}_{P_0})^{\beta_h/2}$.

    Let $\{t_i\}_{i=1}^\infty$ be any sequence in $[0, \delta)$ with $\lim_{i\rightarrow \infty}t_i = 0$. Because $\|\omega^\dagger_{P_{t_i}}\|_{L^2(P_{X, D;0})} < C_{\beta_h}'$ for every $i$, every subsequence of $\{\omega^\dagger_{P_{t_i}}\}_{i=1}^\infty$ admits a further subsequence converging weakly in $L^2(P_{X, D;0})$. Let $\{\omega^\dagger_{P_{t_{i_j}}}\}_{j=1}^\infty$ be any such weakly convergent subsequence, and let $\omega^*_{P_0}$ denote its limit. Because $(\mc{E}_{P_{t_{i_j}}}^*\mc{E}_{P_{t_{i_j}}})^{\beta_h/2}$ converges to $(\mc{E}_{P_0}^*\mc{E}_{P_0})^{\beta_h/2}$ in operator norm, the first term in the decomposition below is bounded in norm by $C'_{\beta_h}\|(\mc{E}_{P_{t_{i_j}}}^*\mc{E}_{P_{t_{i_j}}})^{\beta_h/2} - (\mc{E}_{P_{0}}^*\mc{E}_{P_{0}})^{\beta_h/2}\| \rightarrow 0$. Because $(\mc{E}_{P_0}^*\mc{E}_{P_0})^{\beta_h/2}$ is compact and $\omega^\dagger_{P_{t_{i_j}}} \rightharpoonup \omega^*_{P_0}$, the second term converges in norm to $(\mc{E}_{P_{0}}^*\mc{E}_{P_{0}})^{\beta_h/2}\omega^*_{P_0}$. Hence,
    \begin{align*}
        h_{{P_{t_{i_j}}}}^\dagger = (\mc{E}_{P_{t_{i_j}}}^*\mc{E}_{P_{t_{i_j}}})^{\beta_h/2}\omega^\dagger_{P_{t_{i_j}}} =& \{(\mc{E}_{P_{t_{i_j}}}^*\mc{E}_{P_{t_{i_j}}})^{\beta_h/2} - (\mc{E}_{P_{0}}^*\mc{E}_{P_{0}})^{\beta_h/2}\}\omega^\dagger_{P_{t_{i_j}}} + (\mc{E}_{P_{0}}^*\mc{E}_{P_{0}})^{\beta_h/2}\omega^\dagger_{P_{t_{i_j}}}\\
        \rightarrow&(\mc{E}_{P_{0}}^*\mc{E}_{P_{0}})^{\beta_h/2}\omega^*_{P_0} \eqqcolon h_{P_0}^*.
    \end{align*}
    On the other hand,
    \begin{align*}
        \mc{E}_{P_{t_{i_j}}}h_{{P_{t_{i_j}}}}^\dagger = r_{P_{t_{i_j}}}
    \end{align*}
    with $r_{P_{t_{i_j}}} \rightharpoonup r_{P_{0}}$ by \cref{lemma:weak-convergence-of-solutions:representers}. Because $\mc{E}_{P_{t_{i_j}}} \rightarrow \mc{E}_{P_0}$ in operator norm and $h^\dagger_{P_{t_{i_j}}} \rightarrow h^*_{P_0}$ in norm, and because weak limits are unique,
    \begin{align*}
        \mc{E}_{P_{t_{i_j}}}h_{{P_{t_{i_j}}}}^\dagger \rightarrow \mc{E}_{P_{0}}h_{P_0}^* = r_{P_0}.
    \end{align*}
    But, $ h_{P_0}^* = (\mc{E}_{P_{0}}^*\mc{E}_{P_{0}})^{\beta_h/2}\omega^*_{P_0} \in \mathscr{R}((\mc{E}_{P_{0}}^*\mc{E}_{P_{0}})^{\beta_h/2}) \subseteq \mathscr{N}((\mc{E}_{P_{0}}^*\mc{E}_{P_{0}})^{\beta_h/2})^\perp = \mathscr{N}(\mc{E}_{P_0})^\perp$, where the inclusion holds because $(\mc{E}_{P_{0}}^*\mc{E}_{P_{0}})^{\beta_h/2}$ is self-adjoint, and the equality because $\beta_h > 0$. Hence, because $h_{P_0}^*$ solves
    \begin{align*}
        \mc{E}_{P_0}h_{P_0}^* = r_{P_0}
    \end{align*}
    and belongs to $\mathscr{N}(\mc{E}_{P_0})^\perp$, $h_{P_0}^* = h^\dagger_{P_0}$.

    We have shown that every subsequence of $\{h^\dagger_{P_{t_i}}\}_{i=1}^\infty$ has a further subsequence converging in $L^2(P_{X, D;0})$ to $h^\dagger_{P_0}$. By the Urysohn subsequence principle, $h^\dagger_{P_{t_i}} \rightarrow h^\dagger_{P_0}$ in $L^2(P_{X, D;0})$. Because the sequence $\{t_i\}_{i=1}^\infty$ in $[0, \delta)$ with $\lim_{i\rightarrow \infty}t_i = 0$ was arbitrary, $h^\dagger_{P_t} \rightarrow h^\dagger_{P_0}$ in $L^2(P_{X, D;0})$ as $t \downarrow 0$.
\end{proof}

\begin{lemma}
    \label{lemma:differentiability-of-operators}
    Suppose that $\{P_t:t \in [0, \delta)\} \subseteq \mc{P}$ is QMD at $P_0$, with $P_t = P_0$ at $t=0$, and that each $P_t$ satisfies \cref{as:dominating-measure-support} with a common constant $L_{bounded}$. Then, $t \mapsto \mc{E}_{P_t}$ and $t \mapsto \mc{E}_{P_t}^*$ are differentiable in  $\|\cdot\|_{L^2(P_{X, D;0}) \rightarrow L^2(P_{Z, D;0})}$, $\|\cdot\|_{L^2(P_{Z, D;0}) \rightarrow L^2(P_{X, D;0})}$, norm, respectively, at $t=0$ with derivatives
    \begin{align*}
        \dot {\mc{E}}_{P_0}[s]:& h \mapsto E_{P_0}[s_{X|Z, D}(X, Z, D)h(X, D)|Z, D]\\
        \dot {\mc{E}}^*_{P_0}[s]:& q \mapsto E_{P_0}[s_{Z|X, D}(X, Z, D)q(Z, D)|X, D]
    \end{align*}
    where $s$ is the score of $\{P_t:t \in [0, \delta)\}$ at $t=0$ and
    \begin{align*}
        s_{Z|X, D}(X, Z, D) \coloneqq& E_{P_0}[s(O)|X, Z, D] - E_{P_0}[s(O)|X, D] \quad \text{and} \\
        s_{X|Z, D}(X, Z, D) \coloneqq& E_{P_0}[s(O)|X, Z, D] - E_{P_0}[s(O)|Z, D].
    \end{align*}
    Both $\dot{\mc{E}}_{P_0}[s]$ and $\dot{\mc{E}}^*_{P_0}[s]$ are Hilbert-Schmidt, and hence bounded and compact. Furthermore, the map $t \mapsto p_t$ is differentiable in $L^2(P_0)$ at $t=0$ with derivative $sp_0$. More generally, for any subvector $V$ of $(W, X, Z, D)$, the map $t \mapsto p_{V;t}$ is differentiable in $L^2(P_0)$ at $t=0$ with derivative $E_{P_0}[s(O)|V]p_{V;0}$, so that
    \begin{align*}
        \left.\frac{d}{dt}\frac{p_{V;t}}{p_{V;0}}\right\vert_{t=0} = E_{P_0}[s(O)|V],
    \end{align*}
    and for any disjoint subvectors $U, V$ of $(W, X, Z, D)$, the map $t \mapsto p_{U|V;t}$ is differentiable in $L^2(P_0)$ at $t=0$ with derivative $\{E_{P_0}[s(O)|U, V] - E_{P_0}[s(O)|V]\}p_{U|V;0}$.
\end{lemma}

\begin{proof}[Proof of \cref{lemma:differentiability-of-operators}]
    In what follows, we prove the result only for $\mc{E}_{P_0}^*$. The result for $\mc{E}_{P_0}$ follows symmetrically. Let $s \in L^2_0(P_0)$ denote the score of $\{P_t:t\in[0, \delta)\}$. We let
    \begin{align*}
        s_{Z|X, D}(X, Z, D) \coloneqq& E_{P_0}[s(O)|X, Z, D] - E_{P_0}[s(O)|X, D].
    \end{align*}
    We also take the convention of \cref{remark:marginal-conditional-density-bounds} to avoid tracking unnecessary constants.

    We first argue that $t \mapsto p_t$ is differentiable in $L^2(P_0)$ at $t=0$. We have
    \begin{align*}
        \frac{p_{t} - p_0}{t} - sp_0 =& \frac{\sqrt{p}_{t} - \sqrt{p}_0}{t}\{\sqrt{p}_{t} + \sqrt{p}_0\} - s\sqrt{p_0}\sqrt{p_0}\\
        =& \left\{\frac{\sqrt{p}_{t} - \sqrt{p}_0}{t} - \frac{1}{2}s\sqrt{p_0}\right\}\{\sqrt{p}_{t} +\sqrt{p_0}\} + \frac{1}{2}s\sqrt{p_0}\{\sqrt{p_t}  - \sqrt{p_0}\}.
    \end{align*}
    Now, QMD implies $\|\sqrt{p_t} - \sqrt{p_0}\|_{L^2(\lambda)} = O(t)$, so that $\sqrt{p_t} \rightarrow \sqrt{p_0}$ in $\lambda$-measure. Because, in addition, $|\sqrt{p_t}  - \sqrt{p_0}| < 2\sqrt{L_{bounded}}$ for each $t$ and $s^2p_0 \in L^1(P_0)$, dominated convergence along subsequences gives $\lim_{t\downarrow0}\|\frac{1}{2}s\sqrt{p_0}\{\sqrt{p_t}  - \sqrt{p_0}\}\|_{L^2(P_0)} = 0$, and
    \begin{align*}
         \lim_{t\downarrow0}\left\|\left\{\frac{\sqrt{p}_{t} - \sqrt{p}_0}{t} - \frac{1}{2}s\sqrt{p_0}\right\}\{\sqrt{p}_{t} +\sqrt{p_0}\}\right\|_{L^2(P_0)}\leq 2\sqrt{L_{bounded}}\lim_{t\downarrow0}\left\|\left\{\frac{\sqrt{p}_{t} - \sqrt{p}_0}{t} - \frac{1}{2}s\sqrt{p_0}\right\}\right\|_{L^2(P_0)}=0
    \end{align*}
    where the final equality follows because $\{P_t:t \in [0, \delta)\}$ is QMD at $P_0$ with score $s$. Now let $V$ be any subvector of $(W, X, Z, D)$, and let $\lambda_{-V}$ denote the product of the reference measures of the coordinates not in $V$, which is finite. Integrating out those coordinates and applying Cauchy-Schwarz over $\lambda_{-V}$,
    \begin{align*}
        \left|\frac{p_{V;t} - p_{V;0}}{t} - E_{P_0}[s(O)|V]p_{V;0}\right|^2 \leq \int \left\{\frac{p_{t} - p_{0}}{t} - sp_0\right\}^2 d\lambda_{-V},
    \end{align*}
    so that, using \cref{as:dominating-measure-support} to pass between $\lambda$ and $P_0$, the map $t \mapsto p_{V;t}$ is differentiable in $L^2(P_0)$ at $t=0$ with derivative $E_{P_0}[s(O)|V]p_{V;0}$. Because $p_{V;0}$ does not depend on $t$ and is bounded away from zero, dividing through by it gives the stated derivative of $t \mapsto \frac{p_{V;t}}{p_{V;0}}$. Applying the quotient rule to $p_{U|V;t} = \frac{p_{U, V;t}}{p_{V;t}}$, using the boundedness of these densities above and below uniformly over $t$ given by \cref{as:dominating-measure-support}, gives the stated derivative of $t \mapsto p_{U|V;t}$ for $U, V$ disjoint subvectors. In particular,
    \begin{align*}
        \left.\frac{d}{dt}p_{Z|X, D;t}\right\vert_{t=0} = s_{Z|X, D}p_{Z|X, D;0}
    \end{align*}
    where the derivative is in $L^2(P_0)$.

    We now argue that $\mc{E}_{P_t}^*$ is differentiable at $t=0$. Let
    \begin{align*}
        R(t) \coloneqq \frac{p_{Z|X, D;t} - p_{Z|X, D;0}}{t} - s_{Z|X, D}p_{Z|X, D;0}
    \end{align*}
    which, by the above discussion, satisfies $\lim_{t\downarrow 0}\|R(t)\|_{L^2(P_0)}= 0$.
    Now, for any $q \in L^2(P_{Z, D;0})$, we have that
    \begin{align*}
        &\left|\frac{\{\mc{E}^*_{P_t} - \mc{E}^*_{P_0}\}}{t}q - E_{P_0}[s_{Z|X, D}(X, Z, D)q(Z, D)|X = \cdot,D = \cdot]\right|^2\\
        =&\left|\int \frac{p_{Z|X, D;t}(z|x, d) - p_{Z|X, D;0}(z|x, d)}{t}q(z, d)d\lambda_Z(z) -\int s_{Z|X, D}(x, z, d)q(z, d)p_{Z|X, D;0}(z|x, d)d\lambda_Z(z)\right|^2\\
        =& \left|\int R(t)(z, x, d)q(z, d)d\lambda_Z(z)\right|^2\\
        \leq& \int R(t)(z, x, d)^2d\lambda_Z(z) \int q^2(z, d)d\lambda_Z(z).
    \end{align*}
    Now, using \cref{as:dominating-measure-support}, $\int \int R(t)(z, x, d)^2d\lambda_Z(z)dP_{X, D;0}(x, d) \rightarrow 0$ as $t\downarrow 0$. Consequently,
    \begin{align*}
        &\sup_{\|q\|_{L^2(P_{Z, D;0})} \leq 1}\left\|\frac{\{\mc{E}^*_{P_t} - \mc{E}^*_{P_0}\}}{t}q - E_{P_0}[s_{Z|X, D}(X, Z, D)q(Z, D)|X = \cdot,D = \cdot]\right\|^2_{L^2(P_{X, D;0})}\\
        \leq&\sup_{\|q\|_{L^2(P_{Z, D;0})} \leq 1}\sum_{d}P_0(D=d)\int \int R(t)(z, x, d)^2d\lambda_Z(z)dP_{X|d;0}(x) \int q^2(z, d)d\lambda_Z(z)\\
        \leq &\left\{\max_{d' \in \textsf{Supp}[P_{D;0}]}\int \int R(t)(z, x, d')^2d\lambda_Z(z)dP_{X|d';0}(x)\right\}\sup_{\|q\|_{L^2(P_{Z, D;0})} \leq 1}\sum_{d}P_0(D=d)\int q^2(z, d)d\lambda_Z(z) \\
        \leq &L_{bounded}^{2}\left\{\min_{d \in \textsf{Supp}[P_{D;0}]} P_0(D=d)\right\}^{-1}\| R(t)\|_{L^2(P_0)}^2 \cdot L_{bounded}^{2}\sup_{\|q\|_{L^2(P_{Z, D;0})} \leq 1}\|q\|^2_{L^2(P_{Z, D;0})}\\
        =&L_{bounded}^{4}\left\{\min_{d \in \textsf{Supp}[P_{D;0}]} P_0(D=d)\right\}^{-1}\| R(t)\|_{L^2(P_0)}^2 = o(1)
    \end{align*}
    where the third inequality bounds the $R(t)$ factor by its maximum over $d$, which is attained because $\textsf{Supp}[P_{D;0}]$ is finite, and the fourth uses \cref{as:dominating-measure-support} together with \cref{remark:marginal-conditional-density-bounds}, which give that $p_{Z|d;0}$ and $p_{Z|X, D;0}$ exceed $L_{bounded}^{-2}$ and that $\min_{d \in \textsf{Supp}[P_{D;0}]}P_0(D=d) > 0$.

    We conclude by arguing that $\dot{ \mc{E}}^*_{P_0}[s]$ is Hilbert-Schmidt. This follows from observing that $\dot{ \mc{E}}^*_{P_0}[s]$ is an integral operator with kernel $s_{Z|X, D}(x, z, d)\frac{p_{Z|X, D;0}(z|x, d)}{p_{Z, D;0}(z, d)}$. This kernel belongs to $L^2(P_{Z, D;0})\otimes L^2(P_{X, D;0})$, because the ratio $\frac{p_{Z|X, D;0}}{p_{Z, D;0}}$ is bounded uniformly by \cref{as:dominating-measure-support} and $s_{Z|X, D} \in L^2(P_0)$ as a difference of conditional expectations of $s \in L^2_0(P_0)$. As such it is Hilbert-Schmidt, and so bounded and compact, as a map from $L^2(P_{Z, D;0})$ to $L^2(P_{X, D;0})$. This completes the proof.
\end{proof}

{
We conclude this subsection with a differentiability property of the Moore-Penrose pseudo-inverse. Standard statements of this result assume that the matrix valued mapping is differentiable on a neighborhood of the point in question. Quadratic mean differentiability of a submodel at $t=0$ delivers differentiability of the induced matrix valued mappings at the single point $t=0$ only, and so we record the pointwise version we require.
\begin{lemma}
    \label{lemma:pseudo-inverse-differentiability}
    Let $\delta > 0$ and let $t \mapsto A(t)$ map $[0, \delta)$ into the set of $m\times n$ real matrices. Suppose that
    \begin{enumerate}
        \item \label{lemma:pseudo-inverse-differentiability:diff} $A'(0) \coloneqq \lim_{t \rightarrow 0}t^{-1}\{A(t) - A(0)\}$ exists, and
        \item \label{lemma:pseudo-inverse-differentiability:rank} $\text{rank}(A(t)) = \text{rank}(A(0))$ for each $t \in [0, \delta)$.
    \end{enumerate}
    Then $t \mapsto A(t)^\dagger$ and $t \mapsto I - A(t)^\dagger A(t)$ are differentiable at $t=0$, with
    \begin{align*}
        \left.\frac{d}{dt}A(t)^\dagger\right\vert_{t=0} =& -A(0)^\dagger A'(0)A(0)^\dagger + A(0)^\dagger \{A(0)^\dagger\}^{T}A'(0)^{T}\{I - A(0)A(0)^\dagger\}\\
        &+ \{I - A(0)^\dagger A(0)\}A'(0)^{T}\{A(0)^\dagger\}^{T}A(0)^\dagger,
    \end{align*}
    where $A^\dagger$ denotes the Moore-Penrose pseudo-inverse of a matrix $A$. Here derivatives and limits at $t=0$ are understood as one-sided.
\end{lemma}

\begin{proof}[Proof of \cref{lemma:pseudo-inverse-differentiability}]
    For $A$ and $B$ two $m \times n$ real matrices of equal rank, writing $E \coloneqq B - A$, it holds that
    \begin{align}
        \label{eq:pseudo-inverse-identity}
        B^\dagger - A^\dagger = -B^\dagger EA^\dagger + B^\dagger\{B^\dagger\}^{T}E^{T}(I - AA^\dagger) + (I - B^\dagger B)E^{T}\{A^\dagger\}^{T}A^\dagger
    \end{align}
    (\cite{stewart_matrix_1990} Chapter 3, Lemma 3.7).

    We apply \eqref{eq:pseudo-inverse-identity} with $B \coloneqq A(t)$ and $A \coloneqq A(0)$, which have equal rank by condition \ref{lemma:pseudo-inverse-differentiability:rank}, so that $E = A(t) - A(0)$. Condition \ref{lemma:pseudo-inverse-differentiability:diff} implies $\lim_{t\rightarrow 0}A(t) = A(0)$, which together with condition \ref{lemma:pseudo-inverse-differentiability:rank} gives $\lim_{t \rightarrow 0}A(t)^\dagger = A(0)^\dagger$ (\cite{stewart_matrix_1990} Chapter 3, Corollary 3.10). Dividing \eqref{eq:pseudo-inverse-identity} by $t$ and letting $t \rightarrow 0$, condition \ref{lemma:pseudo-inverse-differentiability:diff} gives $t^{-1}E \rightarrow A'(0)$, and each of the three terms on the right hand side converges. This gives the claimed derivative. Differentiability of $t \mapsto I - A(t)^\dagger A(t)$ at $t=0$ then follows because both $t \mapsto A(t)^\dagger$ and $t \mapsto A(t)$ are differentiable at $t=0$.
\end{proof}

\begin{lemma}
    \label{lemma:re-expressing-tangent-space}
    Let $P_0 \in \mc{P}$ and let $\mc{P}_{\text{Riesz}}$ be as in \eqref{eq:riesz-model}. Suppose $P \mapsto r_P$ is pathwise differentiable at $P_0$ relative to $\mc{P}_{\text{Riesz}}$ with local parameter $\dot r_{P_0}$, and let $h_{P_0}$ be any solution to $\mc{E}_{P_0}h = r_{P_0}$. Suppose that for each $s \in L^\infty(P_0)\cap L^2_0(P_0)$, the law with density $p_t \coloneqq (1 + ts)p_0$ is in $\mc{P}_{\text{Riesz}}$ for $t$ small enough. If $m(W, Z, D;\eta) \in L^2(P_0)$ for each $\eta \in {\mathscr{N}}(\mc{E}_{P_0}^*)$, then the set
    \begin{align*}
        \mathscr{S} \coloneqq \{s \in L^2_0(P_0): \dot r_{P_0}[s] - E_{P_0}[h_{P_0}(X, D)s_{X|Z, D}(X, Z, D)|Z, D] \in \mathscr{N}(\mc{E}_{P_0}^*)^\perp\}
    \end{align*}
    is closed and can be re-expressed as
    \begin{align*}
        \mathscr{S} = \{s \in L^2_0(P_0): \langle s, m(\cdot;\eta) - h_{P_0}\eta\rangle_{L^2(P_0)} = 0 \text{ for all } \eta \in {\mathscr{N}}(\mc{E}_{P_0}^*)\}.
    \end{align*}
    Furthermore, if $r_{P_0}(\cdot, d)$ is a non-constant function of $Z$ and $\mathscr{N}(\mc{E}_{P_0;d}^*) \not= \{0\}$ for some $d \in \textsf{Supp}[P_{D;0}]$, then $\mathscr{S} \subsetneq L^2_0(P_0)$.
\end{lemma}

\begin{proof}[Proof of \cref{lemma:re-expressing-tangent-space}]
    We begin by showing the alternative representation, then that the space $\mathscr{S}$ is strictly included in $L^2_0(P_0)$ under the additional conditions of the lemma.

    For $\eta \in \mathscr{N}(\mc{E}_{P_0}^*)$, define the linear functional
    \begin{align*}
        \psi_\eta[s] \coloneqq \left\langle \dot r_{P_0}[s] - E_{P_0}[h_{P_0}(X, D)s_{X|Z, D}(X, Z, D)|Z = \cdot, D =\cdot], \eta \right\rangle_{L^2(P_{Z, D;0})}, \qquad s \in L^2_0(P_0).
    \end{align*}
    By definition of $\mathscr{N}(\mc{E}_{P_0}^*)^\perp$, $\mathscr{S} = \bigcap_{\eta \in \mathscr{N}(\mc{E}_{P_0}^*)}\mathscr{N}(\psi_\eta)$. We first show that each $\psi_\eta$ is bounded, so that $\mathscr{S}$ is an intersection of closed sets and hence closed. We then show that $\psi_\eta[s] = \langle s, m(\cdot;\eta) - h_{P_0}\eta\rangle_{L^2(P_0)}$ for every $s \in L^2_0(P_0)$, which gives the alternative representation.

    \textbf{Boundedness of $\psi_\eta$.} The map $s \mapsto \dot r_{P_0}[s]$ is bounded, being the local parameter of $P \mapsto r_P$. Furthermore, $h_{P_0}\eta \in L^2(P_0)$ by \cref{lemma:product-in-l2} applied at $P_0$ with $q = \eta$, whose hypotheses hold by \cref{as:dominating-measure-support}. Hence, by the tower property and the Cauchy--Schwarz inequality,
    \begin{align*}
        \left|\left\langle E_{P_0}[h_{P_0}(X, D)s_{X|Z, D}(X, Z, D)|Z = \cdot, D = \cdot], \eta \right\rangle_{L^2(P_{Z, D;0})}\right| = \left|E_{P_0}[h_{P_0}(X, D)\eta(Z, D)s_{X|Z, D}(X, Z, D)]\right| \leq 2\|h_{P_0}\eta\|_{L^2(P_0)}\|s\|_{L^2(P_0)},
    \end{align*}
    where we used that $\|s_{X|Z, D}\|_{L^2(P_0)} \leq 2\|s\|_{L^2(P_0)}$ because conditional expectations are contractions on $L^2(P_0)$.

    \textbf{Reexpressing $\psi_\eta[s]$.} We now derive a more convenient expression for $\psi_\eta[s]$. We first derive the form of $\langle \dot r_{P_0}[s], \eta \rangle_{L^2(P_{Z, D;0})}$.

    Take $s \in L^2_0(P_0) \cap L^\infty(P_0)$ and let $\delta_s > 0$ be such that $p_t \coloneqq (1 + ts)p_0 \in \mc{P}_{\text{Riesz}}$ for $t < \delta_s$. Note that $p_{Z, D;t} = (1 + ts_{Z, D})p_{Z, D;0}$ where $s_{Z, D} \coloneqq E_{P_0}[s(O)|Z, D]$.
    By the defining property of the Riesz representer $r_{P_t}$,
    \begin{align*}
        E_{P_t}[m(W, Z, D;\eta)] = \langle r_{P_t}, \eta\rangle_{L^2(P_{Z, D;t})} = \left\langle \frac{p_{Z, D;t}}{p_{Z, D;0}}r_{P_t}, \eta\right\rangle_{L^2(P_{Z, D;0})}
    \end{align*}
    for each $t \in [0, \delta_s)$, and we differentiate each side at $t=0$. For the left side, because $p_t = (1 + ts)p_0$,
    \begin{align*}
        E_{P_t}[m(W, Z, D;\eta)] = E_{P_0}[m(W, Z, D;\eta)] + tE_{P_0}[m(W, Z, D;\eta)s(O)]
    \end{align*}
    for each $t \in [0, \delta_s)$, where the second expectation is finite because $m(W, Z, D;\eta) \in L^2(P_0)$ by hypothesis and $s \in L^\infty(P_0)$. Hence,
    \begin{align*}
        \left.\frac{d}{dt}E_{P_t}[m(W, Z, D;\eta)]\right\vert_{t=0} = \langle s, m(\cdot;\eta)\rangle_{L^2(P_0)}.
    \end{align*}
    For the right side, we first note that $\langle r_{P_0}, \eta\rangle_{L^2(P_{Z, D;0})} = 0$, because $r_{P_0} \in \mathscr{R}(\mc{E}_{P_0})$ by \cref{as:solutions-exist} and $\eta \in \mathscr{N}(\mc{E}_{P_0}^*) = \mathscr{R}(\mc{E}_{P_0})^\perp$. Hence,
    \begin{align*}
        \lim_{t\downarrow 0}\frac{1}{t}\left\langle \frac{p_{Z, D;t}}{p_{Z, D;0}}r_{P_t}, \eta\right\rangle_{L^2(P_{Z, D;0})}
        =&\lim_{t\downarrow 0}\frac{1}{t}\left\langle \frac{p_{Z, D;0}(1 + ts_{Z, D})}{p_{Z, D;0}}r_{P_t}, \eta\right\rangle_{L^2(P_{Z, D;0})}\\
        =& \lim_{t\downarrow 0}\frac{1}{t}\left\langle r_{P_t}, \eta\right\rangle_{L^2(P_{Z, D;0})} +  \lim_{t\downarrow 0} \left\langle r_{P_t}s_{Z, D}, \eta\right\rangle_{L^2(P_{Z, D;0})}\\
        =& \left\langle \dot r_{P_0}[s], \eta\right\rangle_{L^2(P_{Z, D;0})} +  \left\langle r_{P_0} s_{Z, D},  \eta \right\rangle_{L^2(P_{Z, D;0})}.
    \end{align*}
    Here, the first limit in the third line follows from $\langle r_{P_0}, \eta\rangle_{L^2(P_{Z, D;0})} = 0$ and pathwise differentiability of $P \mapsto r_P$ relative to $\mc{P}_{\text{Riesz}}$, as $\{P_t: t \in [0, \delta_s)\} \subseteq \mc{P}_{\text{Riesz}}$ is QMD at $P_0$ with score $s$. The second limit follows because $r_{P_t} \rightarrow r_{P_0}$ in $L^2(P_{Z, D;0})$ and $s_{Z, D}\eta \in L^2(P_{Z, D;0})$, as $s$ is bounded. Because $\langle s_{Z, D}r_{P_0}, \eta\rangle_{L^2(P_{Z, D;0})} = E_{P_0}[s_{Z, D}(Z, D)r_{P_0}(Z, D)\eta(Z, D)] = \langle r_{P_0}\eta, s\rangle_{L^2(P_0)}$ by the tower property, equating the two derivatives gives
    \begin{align*}
        \langle \dot r_{P_0}[s], \eta \rangle_{L^2(P_{Z, D;0})} = \langle s,  m(\cdot;\eta)\rangle_{L^2(P_0)} - \langle r_{P_0}\eta, s\rangle_{L^2(P_0)}.
    \end{align*}

    Furthermore, writing $s_{X, Z, D} \coloneqq E_{P_0}[s(O)|X = \cdot, Z = \cdot, D = \cdot]$ and using that $E_{P_0}[h_{P_0}(X, D)|Z, D] = r_{P_0}(Z, D)$,
    \begin{align*}
        &\langle E_{P_0}[h_{P_0}(X, D)s_{X|Z, D}(X, Z, D)|Z = \cdot, D = \cdot], \eta \rangle_{L^2(P_{Z, D;0})}\\
        =& E_{P_0}[h_{P_0}(X, D)\eta(Z, D)\{s_{X, Z, D}(X, Z, D) - s_{Z, D}(Z, D)\}]\\
        =& E_{P_0}[h_{P_0}(X, D)\eta(Z, D)s_{X, Z, D}(X, Z, D)] - E_{P_0}[r_{P_0}(Z, D)\eta(Z, D)s_{Z, D}(Z, D)]\\
        =& -\langle r_{P_0}\eta, s \rangle_{L^2(P_0)} + \langle h_{P_0} \eta, s \rangle_{L^2(P_0)}.
    \end{align*}
    Subtracting, the two occurrences of $\langle r_{P_0}\eta, s\rangle_{L^2(P_0)}$ cancel and
    \begin{align*}
        \psi_\eta[s] =& \langle s,  m(\cdot;\eta)\rangle_{L^2(P_0)} - \langle h_{P_0} \eta, s \rangle_{L^2(P_0)} = \left\langle s,  m(\cdot;\eta) - h_{P_0}\eta\right\rangle_{L^2(P_0)}.
    \end{align*}
    The above display was derived for $s \in L^\infty(P_0)\cap L^2_0(P_0)$. Both sides are bounded linear functionals of $s$ on $L^2_0(P_0)$: the left by the boundedness of $\psi_\eta$ shown above, and the right because $m(W, Z, D;\eta) \in L^2(P_0)$ by hypothesis and $h_{P_0}\eta \in L^2(P_0)$. As $L^\infty(P_0) \cap L^2_0(P_0)$ is dense in $L^2_0(P_0)$, the display holds for every $s \in L^2_0(P_0)$. Consequently, $\mathscr{S} = \bigcap_{\eta \in \mathscr{N}(\mc{E}_{P_0}^*)}\mathscr{N}(\psi_\eta)$ coincides with the set in the alternative representation.

    \textbf{Exhibiting $\eta$ with $\psi_\eta$ not identically zero.} Because $\psi_\eta$ is defined on $L^2_0(P_0)$, which is orthogonal to the constants, $\psi_\eta$ vanishes identically if and only if
    \begin{align*}
        h_{P_0}(X, D) \eta(Z, D) - m(W, Z, D;\eta) = c \quad P_0\text{-a.s.}
    \end{align*}
    for some constant $c \in \mathbb{R}$. Now let $d \in \textsf{Supp}[P_{D;0}]$ be as in the statement of the lemma, so that $\mathscr{N}(\mc{E}_{P_0;d}^*) \not= \{0\}$ and $r_{P_0}(\cdot, d)$ is a non-constant function of $Z$. Then $h_{P_0}(\cdot, d)$ is a non-constant function of $X$, as otherwise $r_{P_0}(\cdot, d) = \mc{E}_{P_0;d}h_{P_0}(\cdot, d)$ would be constant. Take $\eta$ to vanish off that level with $\eta(\cdot, d) \in \mathscr{N}(\mc{E}_{P_0;d}^*)$ not identically zero, which is possible because $\mc{E}_{P_0}$ acts blockwise in $d$. By Fubini's theorem, for almost every $z$ the previous display holds for almost every $(w, x)$ given $(Z, D) = (z, d)$, and $\eta(\cdot, d) \not= 0$ on a set of positive measure, so we may fix such a $z$ with $\eta(z, d) \not=0$. The only dependence on $x$ in the previous display is through $h_{P_0}(x, d)\eta(z, d)$, so that display forces
    \begin{align*}
        h_{P_0}(x, d) = \frac{c + m(w, z, d;\eta)}{\eta(z, d)}
    \end{align*}
    for almost every $(w, x)$ in the support of $(W, X)$ given $(Z, D) = (z, d)$. By \cref{as:dominating-measure-support}, that support is a rectangle on which the density is bounded below, so $x$ ranges over the support of $X$ given $D = d$ for almost every such $w$. The right hand side does not depend on $x$, so $h_{P_0}(\cdot, d)$ is constant, a contradiction. Hence $\psi_\eta$ is not identically zero for this $\eta$, and so $\mathscr{S} \subseteq \mathscr{N}(\psi_\eta) \subsetneq L^2_0(P_0)$.
\end{proof}
}

\subsection{Proofs for \cref{subsec:failure-of-continuity}}

\begin{proof}[Proof of \cref{lemma:discrete-necessary-sufficient-continuity}]

    { We begin by giving an informal proof sketch. The proof of \cref{lemma:discrete-necessary-sufficient-continuity:continuity} relies on the following facts. First, when $\mathscr{N}(\mc{E}_{P_0;d}) = \{0\}$ or $\mathscr{N}(\mc{E}_{P_0;d}^*) = \{0\}$, the minimum singular value of $\mc{E}_{P_0;d}$ is separated from zero by a positive constant. Second, with $t$ real-vector valued, given a continuous at $t=0$ matrix valued function $t \mapsto A(t)$, the singular values are continuous in $t$, and if $A(t)$ has rank equal to $A(0)$ for $t$ near zero, the Moore-Penrose pseudo-inverse of $A(t)$ is also continuous at $t=0$. Because $W, X, Z, D$ are discrete, $P$ can be parameterized by a real vector, namely the cell probabilities $P(W=w, X=x, Z=z, D=d)$. Consequently, for $P$ close enough to $P_0$, continuity of the singular values implies $\mc{E}_{P;d}$ is full rank. Then, the mapping $P \mapsto \mc{E}_{P;d}^\dagger$ is continuous at $P_0$, with $\mc{E}_{P;d}^\dagger$ the Moore-Penrose pseudo-inverse. $W, X, Z, D$ discrete and the restrictions on $\mc{P}$ further guarantee the map $P \mapsto \alpha_{P}, r_{P}$ is continuous at $P_0$, and so the solution to the primal inverse problem given by the matrix pseudo-inverse is continuous at $P_0$. Consequently, because $\varphi(P)$ takes the same value at every solution, $P \mapsto \varphi(P)$ is continuous at $P_0$. Note that this proof relies on $X, Z$ finite discrete. When $X, Z$ are continuous, the minimum singular value of $\mc{E}_{P_0;d}$ is not separated from zero, and the Moore-Penrose pseudo inverse is discontinuous.

    On the other hand, when neither $\mathscr{N}(\mc{E}_{P_0;d}) = \{0\}$ nor $\mathscr{N}(\mc{E}_{P_0;d}^*) = \{0\}$ for some $d$, the minimum singular value of $\mc{E}_{P_0;d}$ is equal to zero. This guarantees that the mapping $P \mapsto \mc{E}_{P;d}^\dagger$ is discontinuous. The proof of \cref{lemma:discrete-necessary-sufficient-continuity:discontinuity} proceeds by constructing a sequence $\{P_j\}_{j=1}^\infty$ that exploits this discontinuity by causing the minimum norm solution to diverge to infinity as $j \rightarrow \infty$. We do so by perturbing the singular system of $\mc{E}_{P_0;d}$ so that $\alpha_{P_0}(\cdot, d), r_{P_0}(\cdot, d)$ place non-zero mass on a pair of singular vectors with singular value equal to zero, then perturbing that singular value to be small but non-zero. These perturbations cause small changes in $P$, but large changes to the value of $h^\dagger_P$, and consequently to $\varphi(P)$. The majority of the argument is devoted to ensuring our perturbations result in a valid probability distribution.
    }

    In the following, we note that $h$ solves the equation $\mc{E}_{P}h = f$ if and only if $h_d \coloneqq h(\cdot, d)$ solves $\mc{E}_{P;d}h_d = f(\cdot, d)$ for each $d$, and similarly for $\mc{E}_{P}^*, \mc{E}_{P;d}^*$.

    \textbf{Proof of \cref{lemma:discrete-necessary-sufficient-continuity:continuity}:}

    Let $\{P_j\}_{j=1}^{\infty} \subseteq \mc{P}$ with $\lim_{j \rightarrow \infty }H^2(P_j, P_0) = 0$. We further define the matrices
    \begin{align*}
        {E}_{P;d} \coloneqq& \begin{pmatrix}
            P(X=1|Z=1, D=d) & \cdots & P(X=K_X|Z=1, D=d)\\
            \vdots\\
            P(X=1|Z=K_Z, D=d) & \cdots & P(X=K_X|Z=K_Z, D=d)
        \end{pmatrix}, \\
        {E}_{P;d}^* \coloneqq& \begin{pmatrix}
            P(Z=1|X=1, D=d) & \cdots & P(Z=K_Z|X=1, D=d)\\
            \vdots\\
            P(Z=1|X=K_X, D=d) & \cdots & P(Z=K_Z|X=K_X, D=d).
        \end{pmatrix}
    \end{align*}
    Treating functions of $X$, $Z$ as vectors taking values in $\mathbb{R}^{K_X}$, $\mathbb{R}^{K_Z}$, respectively, it follows that $\mc{E}_{P;d}h = E_{P;d}h$.

    By \cref{as:dominating-measure-support}, and because $W, X, Z$ are finite discretely supported, $\lim_{j \rightarrow \infty}H^2(P_j, P_0) \rightarrow 0$ implies $\lim_{j\rightarrow \infty}\|p_j - p_0\|_{L^q(\lambda)}=0$ for any $1 \leq q \leq \infty$, and consequently there exists $C$ such that for $j$ large enough, $C^{-1} < p_j(O) < C$.

    { Since \(W,X,Z\) are finite-discrete at \(P_0\), and \(D\) is finite-discrete
throughout, the spaces \(L^2(P_{X,D;0})\) and \(L^2(P_{Z,D;0})\) are finite
dimensional.  Hence, for every \(h\in L^2(P_{X,D;0})\) and
\(q\in L^2(P_{Z,D;0})\), the functions
\(\widetilde m(W,X,D;h)\) and \(m(W,Z,D;q)\) are integrable under \(P_0\).
Thus the integrability condition required in \cref{lemma:weak-convergence-of-representers} is automatically
satisfied in the present finite-discrete setting.}
{ Then, because $P_0$ satisfies \cref{as:continuous-linear-maps}, \cref{as:dominating-measure-support},  and $P_j$ satisfy \cref{as:continuous-linear-maps}, by \cref{lemma:weak-convergence-of-representers} }
    \begin{align*}
        \alpha_{P_j} \rightharpoonup \alpha_{P_0},  && r_{P_j} \rightharpoonup r_{P_0}.
    \end{align*}

    Because $L^2(P_{X,D;0})$ and $L^2(P_{Z,D;0})$ are finite dimensional, weak convergence implies strong convergence. Hence,
    \begin{align*}
        \lim_{j\rightarrow \infty}\|\alpha_{P_j} - \alpha_{P_0}\|_{L^2(P_{X,D;0})} = 0 && \lim_{j\rightarrow \infty}\|r_{P_j} - r_{P_0}\|_{L^2(P_{Z,D;0})}=0.
    \end{align*}

    Now, because $\mathscr{N}(\mc{E}_{P_0;d}) = \{0\}$ or $\mathscr{N}(\mc{E}_{P_0;d}^*) = \{0\}$, the matrices ${E}_{P_0;d}$ and ${E}_{P_0;d}^*$ are both rank $\min(K_X, K_Z)$ for each $d$. Hence, the minimum singular value of ${E}_{P_0;d}$, ${E}_{P_0;d}^*$ is non-zero.

    Now, because $\lim_{j \rightarrow \infty}H^2(P_j, P_0)=0$, and because $P_j$, $P_0$ are supported on a rectangle for each $j$, we have that $\lim_{j\rightarrow \infty}|P_j(X=x|Z=z, D=d) - P_0(X=x|Z=z, D=d)| = 0$ and $\lim_{j\rightarrow \infty}|P_j(Z=z|X=x, D=d) - P_0(Z=z|X=x, D=d)| = 0$ for each $z,x,d$. Consequently, ${E}_{P_j;d} \rightarrow {E}_{P_0;d}$ and ${E}^*_{P_j;d} \rightarrow {E}^*_{P_0;d}$ in Frobenius norm.
    Then, for all $j$ large enough, the minimum singular value of ${E}_{P_j;d}, {E}^*_{P_j;d}$ is bounded away from zero (\cite{stewart_matrix_1990} Chapter IV Theorem 1.3). Consequently, for all $j$ large enough, the matrices ${E}_{P_j;d}$ and ${E}_{P_j;d}^*$ are both rank $\min(K_X, K_Z)$. Then, letting $A^\dagger$ denote the Moore-Penrose pseudo-inverse of a matrix $A$, we have that
    \begin{align*}
        \lim_{j \rightarrow \infty }\|{E}_{P_j;d}^\dagger - {E}_{P_0;d}^\dagger\|_{F} = 0, && \lim_{j \rightarrow \infty }\|\{{E}_{P_j;d}^*\}^\dagger - \{{E}_{P_0;d}^*\}^\dagger \|_{F} = 0
    \end{align*}
    by Chapter 3, Corollary 3.10 of \cite{stewart_matrix_1990}. Here $A^\dagger$ is the Euclidean pseudo-inverse of the matrix $A$, so that, because $r_{P}(\cdot, d) \in \mathscr{R}({E}_{P;d})$, ${E}_{P;d}^\dagger r_{P}(\cdot, d)$ solves ${E}_{P;d}h = r_{P}(\cdot, d)$, although it need not be the minimum $L^2(P_{X|d})$-norm solution. This suffices because, as $\alpha_P \in \mathscr{R}(\mc{E}_P^*)$, $\varphi(P)$ takes the same value at every solution of the primal inverse problem. As $K_X, K_Z< \infty$, and because $\alpha_{P_j} \rightarrow \alpha_{P_0}, r_{P_j} \rightarrow r_{P_0}$, we have that
    \begin{align*}
        \lim_{j\rightarrow\infty}\|{E}_{P_j;d}^\dagger r_{P_j}(\cdot, d) - {E}_{P_0;d}^\dagger r_{P_0}(\cdot, d)\|_{L^2(P_{X;0})} &= 0, \\
        \lim_{j\rightarrow\infty}\|\{{E}_{P_j;d}^*\}^\dagger \alpha_{P_j}(\cdot, d) - \{{E}_{P_0;d}^*\}^\dagger \alpha_{P_0}(\cdot, d)\|_{L^2(P_{Z;0})} &= 0.
    \end{align*}
    Then,
    \begin{align*}
        &\lim_{j \rightarrow \infty}|\varphi(P_j) - \varphi(P_0)| \\
        =& \lim_{j \rightarrow \infty}\Bigg|\sum_{d \in \textsf{Supp}[P_{D;0}]}\Big\{P_j(D=d)E_{P_j}\left[\alpha_{P_j}(X, d)({E}_{P_j;d}^\dagger r_{P_j}(\cdot, d))(X)|D=d\right] \\
        &-P_0(D=d)E_{P_0}\left[\alpha_{P_0}(X, d)({E}_{P_0;d}^\dagger r_{P_0}(\cdot, d))(X)|D=d\right]\Big\}\Bigg|\\
        =&\lim_{j \rightarrow \infty} \Bigg|\sum_{d \in \textsf{Supp}[P_{D;0}]}\Big\{P_j(D=d)E_{P_0}\left[\frac{P_j(X|D=d)}{P_0(X|D=d)}\alpha_{P_j}(X, d)({E}_{P_j;d}^\dagger r_{P_j}(\cdot, d))(X)|D=d\right] \\
        &-P_0(D=d)E_{P_0}\left[\alpha_{P_0}(X, d)({E}_{P_0;d}^\dagger r_{P_0}(\cdot, d))(X)|D=d\right]\Big\}\Bigg|\\
        =& 0.
    \end{align*}
    The final equality follows in light of $P_j(D=d) \rightarrow P_0(D=d)$, $\frac{P_j(X|D=d)}{P_0(X|D=d)} \rightarrow 1$, $\alpha_{P_j} \rightarrow \alpha_{P_0}$, and ${E}_{P_j;d}^\dagger r_{P_j}(\cdot, d) \rightarrow {E}_{P_0;d}^\dagger r_{P_0}(\cdot, d)$. { This completes the proof of this part.

    \textbf{Proof of \cref{lemma:discrete-necessary-sufficient-continuity:discontinuity}}
    }
    Let ${d_e} \in \textsf{Supp}[P_{D;0}]$ be such that both
 $\mathscr{N}(\mc{E}_{P_0;{d_e}}) \not= \{0\}$ and $\mathscr{N}(\mc{E}_{P_0;{d_e}}^*) \not= \{0\}$. We define
    \begin{align*}
        L_{gap} \coloneqq \min( \essinf p_0 - L_{bounded}^{-1}, L_{bounded} - \esssup p_0).
    \end{align*}

    We proceed by constructing a sequence $\{P_j\}_{j=1}^\infty \subseteq \mc{P}$ such that $\lim_{j\rightarrow\infty}H^2(P_j, P_0) = 0$ but where $\lim_{j\rightarrow \infty}|\varphi(P_j) - \varphi(P_0)| = \infty$. We let $K \coloneqq \min(K_X, K_Z)$.

    Let $(\sigma_{i;d}, u_{i;d}, v_{i;d})_{i=1}^{K}$ denote a singular system of $\mc{E}_{P_0;d}$, sorted according to decreasing singular values, including zero singular values. $\sigma_{1;d} \leq 1$ because $\mc{E}_{P_0;d}$ is a contraction. Furthermore, we can order this singular system so that $u_{1;d} = 1$,  $v_{1;d} = 1$, and $\sigma_{1;d} = 1$ because  $\mc{E}_{P_0;d}^*\mc{E}_{P_0;d}1 = 1$ and $\mc{E}_{P_0;d}\mc{E}_{P_0;d}^*1 = 1$.

    Let $i_{\alpha}$ denote the index of the first non-constant singular vector $u_{i;{d_e}}$ for which $\langle \alpha_{P_0}(\cdot, {d_e}),  u_{i_{\alpha};{d_e}}\rangle \not=0$ and $i_{r}$ denote the index of the first non-constant singular vector $v_{i;{d_e}}$ for which $\langle r_{P_0}(\cdot, {d_e}),  v_{i_{r};{d_e}}\rangle \not=0$. We note that  $\mc{E}_{P_0;{d_e}}u_{K;{d_e}} = 0$, $\mc{E}_{P_0;{d_e}}^*v_{K;{d_e}}= 0$
    because both $\mathscr{N}(\mc{E}_{P_0;{d_e}}) \not= \{0\}$ and $\mathscr{N}(\mc{E}_{P_0;{d_e}}^*) \not= \{0\}$.

    In what follows, we define
    \begin{align*}
        k_{0;d}(x, z) \coloneqq \sum_{i=1}^K \sigma_{i;d}u_{i;d}(x)v_{i;d}(z)
    \end{align*}
    and let $L_{k_0} > 0$ such that $L_{k_0}^{-1} < k_{0;d}(X, Z) < L_{k_0}$ $P_0$-a.e. for each $d$. Such $L_{k_0}$ exists because $p_{X, Z|d;0}(X, Z|D=d) = k_{0;d}(X, Z)p_{X|d;0}(X|D=d)p_{Z|d;0}(Z|D=d)$ for each $d$, and $L_{bounded}^{-1}< p_{X, Z, D;0}(X, Z, D) < L_{bounded}$.

    We define
    \begin{align}
        \label{eq:rotation-discrete-alpha}
        \nonumber \widetilde{u}_{i_{\alpha};{d_e}}(\delta) \coloneqq&  {u}_{i_{\alpha} ;{d_e}}\cos(\delta) - {u}_{K;{d_e}}\sin(\delta)\\
        \widetilde{u}_{K;{d_e}}(\delta) \coloneqq&  {u}_{i_{\alpha};{d_e}}\sin(\delta) + {u}_{K;{d_e}}\cos(\delta).
    \end{align}
    and
    \begin{align}
        \label{eq:rotation-discrete-r}
        \nonumber \widetilde{v}_{i_r;{d_e}}(\delta) \coloneqq&  {v}_{i_r;{d_e}}\cos(\delta) - {v}_{K;{d_e}}\sin(\delta)\\
        \widetilde{v}_{K;{d_e}}(\delta) \coloneqq&  {v}_{i_r ;{d_e}}\sin(\delta) + {v}_{K;{d_e}}\cos(\delta).
    \end{align}
    Then, $\lim_{\delta \rightarrow 0}\|\widetilde{u}_{i_{\alpha};{d_e}}(\delta) - {u}_{i_{\alpha};{d_e}} \|_{L^\infty(P_{X;0})} = 0$ and similarly for $\widetilde{u}_{K;{d_e}}(\delta)$. Furthermore, because
    \begin{align*}
        \langle {u}_{i_{\alpha} ;{d_e}},\alpha_{P_0}(\cdot, {d_e}) \rangle_{L^2(P_{X|d_e;0})} &\not=0 \quad \text{and} \\
        \langle {v}_{i_{r} ;{d_e}},r_{P_0}(\cdot, {d_e}) \rangle_{L^2(P_{Z|d_e;0})} &\not=0,
    \end{align*}
    the rotations ensure that for all $\delta \not= 0$ small enough,
    \begin{align}
        \label{eq:rotation-not-zero-discrete}
        \langle \widetilde{u}_{K ;{d_e}}(\delta),\alpha_{P_0}(\cdot, {d_e}) \rangle_{L^2(P_{X|d_e;0})} \not=0,  &&\langle \widetilde{v}_{K ;{d_e}}(\delta),r_{P_0}(\cdot, {d_e}) \rangle_{L^2(P_{Z|d_e;0})} \not=0.
    \end{align}

    For each $i \not\in \{i_{\alpha}, K\}$ we let $\widetilde{u}_{i; {d_e}}(\delta) \coloneqq {u}_{i; {d_e}}$ and for each $i \not\in \{i_{r}, K\}$ we let $\widetilde{v}_{i; {d_e}}(\delta) \coloneqq {v}_{i; {d_e}}$.

    We then choose a sequence of non-zero constants $\{\delta_j\}_{j=1}^\infty$ with $\lim_{j\rightarrow\infty}\delta_j = 0$ such that for each $j$, \eqref{eq:rotation-not-zero-discrete} holds with $\delta_j$ in the place of $\delta$ and
    \begin{align}
        \label{eq:i_alpha_j-bound-discrete}
        \sigma_{i_{\alpha};{d_e}}\|\widetilde{u}_{i_{\alpha};{d_e}}(\delta_j)\widetilde{v}_{i_{\alpha};{d_e}}(\delta_j) - {u}_{i_{\alpha};{d_e}}{v}_{i_{\alpha};{d_e}}\|_{L^\infty(P_{X,Z;0})}
        <\frac{1}{3}\min\left({L_{bounded}^{-4}}\min\left(\frac{1}{j}, \frac{1}{2}L_{bounded}^{-2}L_{gap}\right), L_{k_0}^{-1} \right)
    \end{align}
    and
    \begin{align}
        \label{eq:i_r_j-bound-discrete}
        &\sigma_{i_{r};{d_e}}\|\widetilde{u}_{i_{r};{d_e}}(\delta_j)\widetilde{v}_{i_{r};{d_e}}(\delta_j) - {u}_{i_{r};{d_e}}{v}_{i_{r};{d_e}}\|_{L^\infty(P_{X,Z;0})}
        < \frac{1}{3}\min\left({L_{bounded}^{-4}}\min\left(\frac{1}{j}, \frac{1}{2}L_{bounded}^{-2}L_{gap}\right), L_{k_0}^{-1} \right).
    \end{align}

    Having chosen $\delta_j$, we define
    \begin{align*}
    \widetilde{u}_{i_{\alpha};{d_e};j} \coloneqq& \widetilde{u}_{i_{\alpha};{d_e}}(\delta_j), &\widetilde{v}_{i_{r};{d_e};j} \coloneqq& \widetilde{v}_{i_{r};{d_e}}(\delta_j)\\
        \widetilde{u}_{K;{d_e};j} \coloneqq& \widetilde{u}_{K;{d_e}}(\delta_j), &\widetilde{v}_{K;{d_e};j} \coloneqq&\widetilde{v}_{K;{d_e}}(\delta_j),
    \end{align*}
    and let $\widetilde{u}_{i; {d_e};j} \coloneqq u_{i;{d_e}}$ for $i \in \{1, \dots, K\}\setminus \{i_{\alpha}, K\}$, $\widetilde{v}_{i; {d_e};j} \coloneqq v_{i;{d_e}}$ for $i \in \{1, \dots, K\}\setminus \{i_{r}, K\}$.

    Now, for $i \not=K$, we let $\widetilde{\sigma}_{i;{d_e};j} \coloneqq \sigma_{i;{d_e}}$. We choose $\widetilde{\sigma}_{K;{d_e};j} > 0$ such that both
    \begin{align}
        \label{eq:sigma-bound-discrete}\widetilde{\sigma}_{K;{d_e};j}\|\widetilde{u}_{K;{d_e};j}\widetilde{v}_{K;{d_e};j}\|_{L^\infty(P_{X,Z;0})}  < \frac{1}{3}\min\left({L_{bounded}^{-4}}\min\left(\frac{1}{j}, \frac{1}{2}L_{bounded}^{-2}L_{gap}\right), L_{k_0}^{-1} \right)
    \end{align}
    and
    \begin{align}
        \label{eq:sigma_discontinuity-def-discrete}
        \widetilde{\sigma}_{K;{d_e};j} < \frac{P_0(D={d_e})}{j + |\kappa_j|}|\langle \widetilde{u}_{K;{d_e};j}, \alpha_{P_0}(\cdot, {d_e})\rangle_{L^2(P_{X|d_e;0})}\langle \widetilde{v}_{K;{d_e};j}, r_{P_0}(\cdot, {d_e})\rangle_{L^2(P_{Z|d_e;0})}|
    \end{align}
    where
    \begin{align*}
        \kappa_j \coloneqq& P_0(D={d_e})\sum_{i< K: \widetilde{\sigma}_{i;{d_e};j}>0}\widetilde{\sigma}_{i;{d_e};j}^{-1}\langle \widetilde{u}_{i; {d_e};j}, \alpha_{P_0}(\cdot, {d_e})\rangle_{L^2(P_{X|d_e;0})}\langle \widetilde{v}_{i; {d_e};j}, r_{P_0}(\cdot, {d_e})\rangle_{L^2(P_{Z|d_e;0})}\\
        &+ \sum_{d \in \textsf{Supp}[P_{D;0}]\setminus \{{d_e}\}}P_0(D=d)\sum_{i\leq K: \sigma_{i;{d}}>0}{\sigma}_{i;{d}}^{-1}\langle {u}_{i; {d}}, \alpha_{P_0}(\cdot, {d})\rangle_{L^2(P_{X|d;0})}\langle {v}_{i; {d}}, r_{P_0}(\cdot, {d})\rangle_{L^2(P_{Z|d;0})}.
    \end{align*}

    We define
    \begin{align*}
        k_{j;d}(X, Z) \coloneqq& I(d = {d_e})\sum_{i=1}^K \widetilde{\sigma}_{i;{d_e};j}\widetilde{u}_{i;{d_e};j}(X)\widetilde{v}_{i;{d_e};j}(Z)\\
        &+ (1 - I(d = {d_e}))k_{0;d}(X, Z).
    \end{align*}

    We then have that $\|k_{j;d}(X, Z) - k_{0;d}(X, Z)\|_{L^\infty(P_{X,Z;0})} = 0$ for $d \not= {d_e}$ and
    \begin{align*}
        \|k_{j;{d_e}}(X, Z) - k_{0;{d_e}}(X, Z)\|_{L^\infty(P_{X,Z;0})}\leq& \sigma_{i_{\alpha};{d_e}}\|\widetilde{u}_{i_{\alpha};{d_e};j}\widetilde{v}_{i_{\alpha};{d_e};j} - {u}_{i_{\alpha};{d_e}}{v}_{i_{\alpha};{d_e}}\|_{L^\infty(P_{X,Z;0})}\\
         \nonumber &+\sigma_{i_{r};{d_e}}\|\widetilde{u}_{i_{r};{d_e};j}\widetilde{v}_{i_{r};{d_e};j} - {u}_{i_{r};{d_e}}{v}_{i_{r};{d_e}}\|_{L^\infty(P_{X,Z;0})}\\
        \nonumber&+\widetilde{\sigma}_{K;{d_e};j}\|\widetilde{u}_{K;{d_e};j}\widetilde{v}_{K;{d_e};j}\|_{L^\infty(P_{X,Z;0})}.
    \end{align*}
    The bounds \eqref{eq:i_alpha_j-bound-discrete},\eqref{eq:i_r_j-bound-discrete}, and \eqref{eq:sigma-bound-discrete} ensure that
    \begin{align*}
        \|k_{j;{d_e}}(X, Z) - k_{0;{d_e}}(X, Z)\|_{L^\infty(P_{X,Z;0})}<\min\left({L_{bounded}^{-4}}\min\left(\frac{1}{j}, \frac{1}{2}L_{bounded}^{-2}L_{gap}\right), L_{k_0}^{-1} \right).
    \end{align*}
    Because $k_{0;{d_e}}(X, Z) > L_{k_0}^{-1}$, we have that $k_{j;{d_e}}(X, Z) > 0$ $P_0$-a.e. Recall from \cref{remark:marginal-conditional-density-bounds} that, because $p_0 < L_{bounded}$, we have $p_{X|d_e;0} < L_{bounded}^2$, $p_{Z|d_e;0} < L_{bounded}^2$, and $p_{D;0}(d_e) < L_{bounded}$, where the reference measures are normalized as in that remark. Defining
    \begin{align*}
        \widetilde{p}_{X, Z|{d_e};j}(X, Z) \coloneqq k_{j;{d_e}}(X, Z) p_{X|{d_e};0}(X)p_{Z|{d_e};0}(Z)
    \end{align*}
    then gives that
    \begin{align}
        \label{eq:density-bound-discrete}
        \nonumber&\|\widetilde{p}_{X, Z|{d_e};j}(X, Z) - p_{X, Z|{d_e};0}(X, Z)\|_{L^\infty(P_{X,Z;0})}\\
        \nonumber\leq&\|k_{j;{d_e}}(X, Z) - k_{0;{d_e}}(X, Z)\|_{L^\infty(P_{X,Z;0})} \|p_{X|{d_e};0}\|_{L^\infty(P_{X;0})}\|p_{Z|{d_e};0}\|_{L^\infty(P_{Z;0})} \\
        <& \min\left(\frac{1}{j}, \frac{1}{2}L_{bounded}^{-2}L_{gap}\right).
    \end{align}

    { Because, as noted before, \(k_{j;d_e}(X,Z)>0\) \(P_0\)-a.e. and
\(p_{X|d_e;0}\) and \(p_{Z|d_e;0}\) are strictly positive on
their supports, it follows that
\(\widetilde p_{X,Z|d_e;j}(X,Z)>0\)  \(P_0\)-a.e.}

    We further argue that $\widetilde{p}_{X, Z|{d_e};j}(X, Z)$ is a density for $(X, Z)$. Indeed,
    \begin{align*}
        \int \int \widetilde{p}_{X, Z|{d_e};j}(x, z)d\lambda_X(x)d\lambda_Z(z) =& \int\int k_{j;{d_e}}(x, z)p_{X|{d_e};0}(x)p_{Z|{d_e};0}(z)d\lambda_X(x)d\lambda_Z(z)\\
        =& \int\int \sum_{i=1}^K\widetilde{\sigma}_{i;{d_e};j}\widetilde{u}_{i;{d_e};j}(x)\widetilde{v}_{i;{d_e};j}(z)\\
        &\qquad\times p_{X|{d_e};0}(x)p_{Z|{d_e};0}(z)d\lambda_X(x)d\lambda_Z(z)\\
        =& \sum_{i=1}^K \widetilde{\sigma}_{i;{d_e};j} \langle \widetilde{u}_{i;{d_e};j}, 1\rangle_{L^2(P_{X|d_e;0})}\langle \widetilde{v}_{i;{d_e};j}, 1\rangle_{L^2(P_{Z|d_e;0})}\\
        =& 1
    \end{align*}
    where the last equality follows because $\widetilde{\sigma}_{1;{d_e};j} = 1$, $\widetilde{u}_{1;{d_e};j} = 1$, and $\widetilde{v}_{1;{d_e};j} = 1$. Furthermore, $\widetilde{p}_{Z|{d_e};j} = {p}_{Z|{d_e};0}$ as
    \begin{align*}
        \int \widetilde{p}_{X, Z|{d_e};j}(x, z)d\lambda_X(x) =& \int k_{j;{d_e}}(x, z)p_{X|{d_e};0}(x)p_{Z|{d_e};0}(z)d\lambda_X(x)\\
        =& p_{Z|{d_e};0}(z)\int \sum_{i=1}^K\widetilde{\sigma}_{i;{d_e};j}\widetilde{u}_{i;{d_e};j}(x)\widetilde{v}_{i;{d_e};j}(z)p_{X|{d_e};0}(x)d\lambda_X(x)\\
       =&  p_{Z|{d_e};0}(z)\sum_{i=1}^K \widetilde{\sigma}_{i;{d_e};j} \langle \widetilde{u}_{i;{d_e};j}, 1\rangle_{L^2(P_{X|d_e;0})}  \widetilde{v}_{i;{d_e};j}(z)\\
        =& p_{Z|{d_e};0}(z)
    \end{align*}
    where the final equality holds because $\widetilde{u}_{1;{d_e};j} = 1$, $\widetilde{v}_{1;{d_e};j} = 1$, $\widetilde{\sigma}_{1;{d_e};j} = 1$, and $\langle \widetilde{u}_{i;{d_e};j}, 1\rangle_{L^2(P_{X|d_e;0})} = 0$ for $i \geq 2$, as each such $\widetilde{u}_{i;{d_e};j}$ is a linear combination of singular functions orthogonal to $u_{1;{d_e}} = 1$. Similarly, $\widetilde{p}_{X|{d_e};j} = {p}_{X|{d_e};0}$.

    Let $\widetilde{p}_{X, Z|d;j} \coloneqq p_{X, Z|d;0}$ for $d \not= {d_e}$ and $\widetilde{P}_j(D=d) = P_0(D=d)$, and let $\widetilde{p}_{X, Z, D;j}(x, z, d) \coloneqq \widetilde{p}_{X, Z|d;j}(x, z)p_{D;0}(d)$ denote the corresponding joint density of $(X, Z, D)$. Because $\widetilde{p}_{X, Z, D;j}$ and $p_{X, Z, D;0}$ agree for $d \not= d_e$ and both vanish off the support of $P_0$, \eqref{eq:density-bound-discrete} and $p_{D;0}(d_e) < L_{bounded}$ give
    \begin{align*}
        \|\widetilde{p}_{X, Z, D;j} - p_{X, Z, D;0}\|_{L^\infty(\lambda_X\otimes\lambda_Z\otimes\lambda_D)} = p_{D;0}(d_e)\|\widetilde{p}_{X, Z|{d_e};j} - p_{X, Z|{d_e};0}\|_{L^\infty(P_{X,Z;0})} < \frac{1}{2}L_{bounded}^{-1}L_{gap}.
    \end{align*}
    With $P \coloneqq P_0$, the threshold in \cref{lemma:constant-riesz-representers} is $L_{bounded}^{-1}L_{gap}\{\lambda_X(\mathcal X)\lambda_Z(\mathcal Z)\lambda_D(\mc{D})\}^{-1} = L_{bounded}^{-1}L_{gap}$ under the normalization of \cref{remark:marginal-conditional-density-bounds}. Hence the conditions of \cref{lemma:constant-riesz-representers} are satisfied with $p'_{X, Z, D}\coloneqq\widetilde{p}_{X, Z, D;j}$ and $P \coloneqq P_0$ in that Lemma. We let $P_j$ be $P^\dagger$ of the Lemma. Then, by that Lemma,
    \begin{enumerate}
        \item $P_j$ and $P_0$ are mutually absolutely continuous
        \item $L_{bounded}^{-1} < \essinf p_j \leq \esssup p_j < L_{bounded}$
        \item $P_j$ satisfies \cref{as:continuous-linear-maps} with $\alpha_{P_j} = \alpha_{P_0}$, $r_{P_j} = r_{P_0}$.
        \item $H^2(P_j, P_0) \leq \frac{1}{2}\|\widetilde{p}_{X, Z, D;j} - p_{X, Z, D;0}\|_{L^1(\lambda_X\otimes\lambda_Z\otimes\lambda_D)}$.
    \end{enumerate}

    Now,
    \begin{align*}
        \frac{1}{2}\|\widetilde{p}_{X, Z, D;j} - p_{X, Z, D;0}\|_{L^1(\lambda_X\otimes\lambda_Z\otimes\lambda_D)} <& L' \|\widetilde{p}_{X, Z, D;j} - p_{X, Z, D;0}\|_{L^\infty(\lambda_X\otimes\lambda_Z\otimes\lambda_D)}\\
         =&  L'\|\widetilde{p}_{X, Z, D;j} - p_{X, Z, D;0}\|_{L^\infty(P_{X,Z,D;0})}\\
        &\leq \sum_{d\in \textsf{Supp}[P_{D;0}]}L'\|\widetilde{p}_{X, Z|d;j} - p_{X, Z|d;0}\|_{L^\infty(P_{X,Z;0})}P_0(D=d)\\
        &\leq L'/j
    \end{align*}
    for some constant $L'$ not depending on $j$. The final inequality follows from \eqref{eq:density-bound-discrete}. Consequently, $\lim_{j \rightarrow \infty}H^2(P_j, P_0) = 0$.

    We now argue $P_j$ belongs to $\mc{P}$ for each $j$, i.e. Assumptions \ref{as:continuous-linear-maps}, \ref{as:solutions-exist}, \ref{as:beta-source} (for all $\beta_h, \beta_q \in [0, \infty)$), \ref{as:mutual-absolute-continuity} (relative to $P_0$), and \ref{as:dominating-measure-support} hold.  \cref{as:continuous-linear-maps} holds by the construction in \cref{lemma:constant-riesz-representers} and because it holds at $P_0$. \cref{as:solutions-exist} holds because $\alpha_{P_0}(\cdot, d) \in \mathscr{R}(\mc{E}^*_{P_0;{d}})$, $r_{P_0}(\cdot, d) \in \mathscr{R}(\mc{E}_{P_0;{d}})$ for each $d$ because \cref{as:solutions-exist} holds at $P_0$,  $\mc{E}_{P_j;d} = \mc{E}_{P_0;d}$, $\mc{E}_{P_j;d}^* = \mc{E}_{P_0;d}^*$, for each $d \not = {d_e}$, $\mathscr{R}(\mc{E}_{P_0;{d_e}}) \subset  \mathscr{R}(\mc{E}_{P_j;{d_e}})$, $\mathscr{R}(\mc{E}^*_{P_0;{d_e}}) \subset  \mathscr{R}(\mc{E}^*_{P_j;{d_e}})$ by ~\eqref{eq:rotation-discrete-alpha}, ~\eqref{eq:rotation-discrete-r}. For the last two inclusions, note that $v_{i_r;d_e}, v_{K;d_e} \in \text{span}\{\widetilde{v}_{i_r;d_e;j}, \widetilde{v}_{K;d_e;j}\}$ and $u_{i_\alpha;d_e}, u_{K;d_e} \in \text{span}\{\widetilde{u}_{i_\alpha;d_e;j}, \widetilde{u}_{K;d_e;j}\}$, where each of these rotated singular functions has a strictly positive singular value at $P_j$, while the remaining singular functions with positive singular values at $P_0$ are unchanged.
    Because $X, Z$ are finite discretely valued, \cref{as:solutions-exist} implies \cref{as:beta-source}.  \cref{as:mutual-absolute-continuity} holds because ${P}_j$ is mutually absolutely continuous with $P_0$ by \cref{lemma:constant-riesz-representers}. Finally, \cref{as:dominating-measure-support} holds by the construction of ${P}_j$ in \cref{lemma:constant-riesz-representers}.

    Hence, we have constructed a sequence of distributions $\{P_j\}_{j=1}^\infty$ such that, for each $j$, $P_j$ belongs to the statistical model $\mc{P}$ defined in \eqref{eq:model-definition}  and where $\lim_{j \rightarrow \infty}H^2(P_j, P_0) = 0$. We conclude by demonstrating that $\lim_{j\rightarrow \infty}|\varphi(P_j) - \varphi(P_0)|= \infty$. It suffices to show that $\lim_{j\rightarrow \infty}|\varphi(P_j)|= \infty$.

    Let $h_{{P}_j}^* \coloneqq \sum_{d \in \textsf{Supp}[P_{D;0}]}I(D=d)h_{{P}_j}^*(\cdot, d)$ where for each $d$, $h_{{P}_j}^*(\cdot, d)$ is the  minimum norm solution to the inverse problem
    \begin{align*}
        \mc{E}_{{P}_j;d}h(\cdot, d) = r_{{P}_j}(\cdot, d) = r_{{P}_0}(\cdot, d).
    \end{align*}
    Then,
    \begin{align*}
        (\mc{E}_{{P}_j}h_{{P}_j}^*)(Z, D) =& E_{{P}_j}[h_{{P}_j}^*(X, D)|Z, D]\\
        =& \sum_{d \in \textsf{Supp}[P_{D;0}]}E_{{P}_j}[h_{{P}_j}^*(X, d)|Z, D=d]I(D=d)\\
        =& \sum_{d \in \textsf{Supp}[P_{D;0}]}(\mc{E}_{{P}_j;d}h_{{P}_j}^*(\cdot, d))(Z)I(D=d)\\
        =& \sum_{d \in \textsf{Supp}[P_{D;0}]}r_{P_j}(\cdot, d)I(D=d)\\
        =& r_{P_j},
    \end{align*}
    i.e. $h_{{P}_j}^*$ solves the primal inverse problem. Then, we have that,
    \begin{align*}
        |\varphi({P}_j)| =& \left|E_{{P}_j}[\alpha_{{P}_j}(X, D)h_{{P}_j}^*(X, D)]\right|\\
        =& \left|\sum_{d \in \textsf{Supp}[P_{D;0}]}P_j(D=d)E_{{P}_j}[\alpha_{{P}_j}(X, d)h_{{P}_j}^*(X, d)|D=d]\right|\\
        =& \Bigg|\sum_{d \in \textsf{Supp}[P_{D;0}]\setminus\{{d_e}\}}P_0(D=d)\sum_{i\leq K: \sigma_{i;d}>0}^{K} {\sigma}_{i;d}^{-1}\left\langle \alpha_{{P}_0}(\cdot, d), {u}_{i;d} \right\rangle_{L^2(P_{X|d;0})}\left\langle r_{{P}_0}(\cdot, d), {v}_{i;d} \right\rangle_{L^2(P_{Z|d;0})}\\
        &+P_0(D={d_e})\sum_{i\leq K: \widetilde{\sigma}_{i;{d_e};j}>0}^{K} \widetilde{\sigma}_{i;{d_e};j}^{-1}\left\langle \alpha_{{P}_0}(\cdot, {d_e}), \widetilde{u}_{i;{d_e};j} \right\rangle_{L^2(P_{X|d_e;0})}\left\langle r_{{P}_0}(\cdot, {d_e}), \widetilde{v}_{i;{d_e};j} \right\rangle_{L^2(P_{Z|d_e;0})}\Bigg|\\
        \geq&P_0(D={d_e})\widetilde{\sigma}_{K;{d_e};j}^{-1}\left|\left\langle \alpha_{{P}_0}(\cdot, {d_e}), \widetilde{u}_{K;{d_e};j} \right\rangle_{L^2(P_{X|d_e;0})}\left\langle r_{{P}_0}(\cdot,{d_e}), \widetilde{v}_{K;{d_e};j} \right\rangle_{L^2(P_{Z|d_e;0})}\right| - |\kappa_j| \\
        > & j,
    \end{align*}
    where the first inequality is the reverse triangle inequality, and the second follows from \eqref{eq:sigma_discontinuity-def-discrete}.
    Hence, $\lim_{j\rightarrow\infty}|\varphi({P}_j) - \varphi(P_0)| = \infty$.
\end{proof}

\begin{proof}[Proof of \cref{thrm:continuous-discontinuity}]

    { We begin by providing a sketch of the proof. The core argument is similar to that of the discrete case. Because the minimum eigenvalue of $\mc{E}_{P;d}$ is not separated from zero, the mapping $P \mapsto \mc{E}_{P;d}^\dagger$ with $\mc{E}_{P;d}^\dagger$ the Moore-Penrose pseudo-inverse is discontinuous. We exploit this discontinuity by perturbing $P_0$ to $P_j$ where $\alpha_{P_j}, r_{P_j}$ are supported on a pair of singular vectors with small singular value. Perturbing this singular value closer to zero then causes the minimum norm solution $h^\dagger_{P_j}$, and so $\varphi(P_j)$, to diverge. As $j$ increases, we do this perturbation further into the tail of the singular system, so that the perturbation has a vanishingly small effect on the Hellinger distance between $P_j$ and $P_0$.

    As in the discrete case, the majority of the argument is devoted to ensuring our perturbation $P_j$ is a valid distribution. In contrast to the discrete case, this is complicated by the fact that, in general, the uniform norm of the singular functions can grow to infinity as we move further into the tail of the singular system. Then, a small perturbation to the tail of the singular system may cause the kernel of the conditional expectation operator to be negative for some values of $x, z$. However, as the kernel is a conditional density for the appropriate choice of reference measure, it must be non-negative. We solve this by controlling the uniform norm of these singular vectors by constructing  $x, z \mapsto p_{X, Z|d;j}$ infinitely differentiable for each $j$, $d$. \cref{lemma:uniform-convergence-of-kernel} then shows that a truncated singular system of an infinitely differentiable density will still be a valid density, as the tail contributes negligibly in uniform norm. To ensure $\alpha_{P_j}$, $r_{P_j}$ are in range of this truncated operator, we add the residuals of $\alpha_{P_j}$ and $r_{P_j}$ as singular vectors, with singular value small enough to ensure the kernel of the new operator is a valid conditional density. The $\beta$-source conditions are also trivially satisfied for every $\beta_h, \beta_q$ as the operator $\mc{E}_{P_j;d}$ is finite dimensional.}

    We now begin our formal arguments. Write $W = (W_C, W_D)$, where $W_C$ collects the components of $W$ whose reference measure in \cref{as:dominating-measure-support} is the Lebesgue measure and $W_D$ collects those whose reference measure is the counting measure; either may be empty. Correspondingly, $\lambda_W = \lambda_{W_C}\otimes\lambda_{W_D}$, where $\lambda_{W_C}$, $\lambda_X$, $\lambda_Z$ are Lebesgue product measures and $\lambda_{W_D}$, $\lambda_D$ are counting measures. All densities are taken with respect to $\lambda$, or the appropriate product measure for functions of $W, X, Z, D$ and combinations thereof. We normalize these reference measures as in \cref{remark:marginal-conditional-density-bounds}, so that $\lambda_W(\textsf{Supp}[W])= \lambda_X(\textsf{Supp}[X])= \lambda_Z(\textsf{Supp}[Z])= \lambda_D(\textsf{Supp}[D])=1$.

    Let $\widetilde{d}\in \textsf{Supp}[D;P_0]$ such that $\alpha_{P_0}(\cdot, \widetilde{d}), r_{P_0}(\cdot, \widetilde{d})$ are not $P_0$-a.e. constant functions.

    We first construct a sequence of laws $(P_j)_{j=1}^\infty$ whose densities are smooth in $(W_C, X, Z)$ and lie strictly inside the density bounds of \cref{as:dominating-measure-support}. Because $(W_D, D)$ is supported on a bounded rectangle in a product of copies of $\mathbb{N}$, the set $\mc{S} \coloneqq \textsf{Supp}[(W_D, D);P_0]$ is finite. Let $\Omega_C$ denote the open bounded rectangle supporting $(W_C, X, Z)$ under $P_0$. Because the support of $P_0$ is a rectangle, $\Omega_C$ also supports the conditional law of $(W_C, X, Z)$ given $(W_D, D) = (w_D, d)$ for each $(w_D, d) \in \mc{S}$. By \cref{as:dominating-measure-support}, we may fix constants $c_{lo}, c_{hi}$ with
    \begin{align*}
        L_{bounded}^{-1} < c_{lo} < \essinf p_0 \leq \esssup p_0 < c_{hi} < L_{bounded}.
    \end{align*}
    For $(w_D, d) \in \mc{S}$, let $q(w_D, d) \coloneqq p_{W_D, D;0}(w_D, d) > 0$, and let $p_{W_C, X, Z|w_D, d;0} \coloneqq p_0(\cdot, w_D, \cdot, \cdot, d)/q(w_D, d)$ denote the conditional density of $(W_C, X, Z)$ given $(W_D, D) = (w_D, d)$ with respect to $\lambda_{W_C}\otimes\lambda_X\otimes\lambda_Z$. Then
    \begin{align*}
        \frac{c_{lo}}{q(w_D, d)} < p_{W_C, X, Z|w_D, d;0} < \frac{c_{hi}}{q(w_D, d)} \qquad \text{a.e. on }\Omega_C.
    \end{align*}
    We apply \cref{lemma:dense-continuous-densities} with $Y = (W_C, X, Z)$, $\Omega = \Omega_C$, $p = p_{W_C, X, Z|w_D, d;0}$, $a = c_{lo}/q(w_D, d)$, and $b = c_{hi}/q(w_D, d)$, separately for each of the finitely many $(w_D, d) \in \mc{S}$. Because $\gamma(\varepsilon) \rightarrow 1$ as $\varepsilon \downarrow 0$ in that Lemma and $\mc{S}$ is finite, for each $j$ we may choose $\varepsilon = \varepsilon_j$ small enough that the resulting densities $p_{W_C, X, Z|w_D, d;j}$ supported on $\Omega_C$, with corresponding constants $\gamma_j(w_D, d) > 0$, satisfy for every $(w_D, d) \in \mc{S}$:
    \begin{enumerate}
        \item $p_{W_C, X, Z|w_D, d;j} \in C^\infty(\overline{\Omega}_C)$;
        \item $c_{lo}/\{q(w_D, d)\gamma_j(w_D, d)\} < p_{W_C, X, Z|w_D, d;j} < c_{hi}/\{q(w_D, d)\gamma_j(w_D, d)\}$ a.e. on $\Omega_C$, where $L_{bounded}^{-1} < c_{lo}/\gamma_j(w_D, d)$ and $c_{hi}/\gamma_j(w_D, d) < L_{bounded}$;
        \item $\|p_{W_C, X, Z|w_D, d;j} - p_{W_C, X, Z|w_D, d;0}\|_{L^1(\lambda_{W_C}\otimes\lambda_X\otimes\lambda_Z)} < 1/j$.
    \end{enumerate}
    We then let $P_j$ be the law with density
    \begin{align*}
        p_j(w_C, w_D, x, z, d) \coloneqq q(w_D, d)\,p_{W_C, X, Z|w_D, d;j}(w_C, x, z)
    \end{align*}
    for $(w_D, d) \in \mc{S}$, and $p_j \coloneqq 0$ otherwise. Then $p_j$ is a probability density with $p_{W_D, D;j} = p_{W_D, D;0}$, and in particular $P_j(D=d) = P_0(D=d)$ for each $d$. Multiplying the second item through by $q(w_D, d)$, and taking the minimum and maximum over the finitely many $(w_D, d) \in \mc{S}$, gives
    \begin{align*}
        L_{bounded}^{-1} < \essinf p_j \leq \esssup p_j < L_{bounded},
    \end{align*}
    where the essential infimum and supremum are taken over the support of $P_0$, which is also the support of $P_j$. Hence $P_j$ and $P_0$ are mutually absolutely continuous, and $P_j$ satisfies \cref{as:dominating-measure-support}. By the third item,
    \begin{align*}
        \|p_j - p_0\|_{L^1(\lambda)} = \sum_{(w_D, d) \in \mc{S}} P_0(W_D = w_D, D=d)\|p_{W_C, X, Z|w_D, d;j} - p_{W_C, X, Z|w_D, d;0}\|_{L^1(\lambda_{W_C}\otimes\lambda_X\otimes\lambda_Z)} < 1/j,
    \end{align*}
    so that $H^2(P_j, P_0) \leq \frac{1}{2}\|p_j - p_0\|_{L^1(\lambda)} < 1/j$. Because $p_j$ and $p_0$ lie between $L_{bounded}^{-1}$ and $L_{bounded}$ on their common support, $L^1(\lambda)$ convergence of $p_j$ to $p_0$ implies $L^q$ convergence for every $1 \leq q < \infty$, both of $p_j$ and of each of its marginal and conditional densities. Finally, for each $d$, $p_{X, Z|d;j} \in C^\infty(\overline{\Omega}_{X, Z})$: it is a finite positive linear combination of the functions $(x, z) \mapsto \int p_{W_C, X, Z|w_D, d;j}(w_C, x, z)\,d\lambda_{W_C}(w_C)$ over $w_D$ with $(w_D, d) \in \mc{S}$ (with the integral omitted if $W_C$ is empty), and differentiation under the integral sign is justified because every derivative of $p_{W_C, X, Z|w_D, d;j}$ is continuous on the compact set $\overline{\Omega}_C$, hence bounded.

    Notice that at $P_0$,  $\alpha_{P_0} \in \mathscr{R}(\mc{E}_{P_0}^*)$, $r_{P_0} \in \mathscr{R}(\mc{E}_{P_0})$ by \cref{as:solutions-exist}. Combined with \cref{as:dominating-measure-support}, \cref{lemma:in-range-uniformly-bounded} shows that the Riesz representers $\alpha_{P_0}, r_{P_0}$ are essentially uniformly bounded.

    Because each $P_j$ is mutually absolutely continuous with $P_0$ and has density bounded above and below by constants, \cref{lemma:l1-bounded-functional} shows that each $P_j$ satisfies \cref{as:continuous-linear-maps}, with $\alpha_{P_j} \in L^\infty(P_{X,D;j})$ and $r_{P_j} \in L^\infty(P_{Z,D;j})$. We define
    \begin{align*}
        L_{gap;j} \coloneqq \min(\essinf p_j - L_{bounded}^{-1}, L_{bounded} - \esssup p_j),
    \end{align*}
    which is strictly positive for each $j$ by the construction of $P_j$.

    In what follows, for $P$ a law on $O$, we let
    \begin{align*}
        \mc{E}_{P;d}: L^2(P_{X|d}) \rightarrow L^2(P_{Z|d}), \mc{E}_{P;d}h_d \coloneqq E_{P}[h_d(X)|Z, D=d],\\
        \mc{E}_{P;d}^*: L^2(P_{Z|d}) \rightarrow L^2(P_{X|d}), \mc{E}_{P;d}^*q_d \coloneqq E_{P}[q_d(Z)|X, D=d],
    \end{align*}
    and note that for $h \in L^2(P_{X|d})$, $q \in  L^2(P_{Z|d})$
    \begin{align*}
        (\mc{E}_{P}h)(Z, D) =& \sum_{d \in \textsf{Supp}[D;P]} I(D=d)(\mc{E}_{P;d}h(\cdot, d))(Z),\\
        (\mc{E}_{P}^*q)(X, D) =& \sum_{d \in \textsf{Supp}[D;P]} I(D=d)(\mc{E}^*_{P;d}q(\cdot, d))(X).
    \end{align*}

    Let $\Omega$ denote the open bounded rectangle upon which $P_0$ is supported, i.e. $P_0(O \in \Omega) = 1$, and let $\Omega_{W ,X, Z}, \Omega_{X, Z}$, $\Omega_X, \Omega_Z$ be the open bounded rectangles supporting $(W, X, Z), (X, Z)$, $X$, and $Z$,  respectively. Let $(\sigma_{i;d;j}, u_{i;d;j}, v_{i;d;j})_{i=1}^\infty$ denote the singular system of $\mc{E}_{{P}_j;{d}}$ for each $d, j$ with singular values decreasing in magnitude as $i \rightarrow \infty$ and repeated according to multiplicity. This singular system exists because $L_{bounded}^{-1} < p_j(O) < L_{bounded}$ $P_j$-a.e. implies that $\mc{E}_{P_j;d}$ is compact (\cite{carrasco_chapter_2007} Assumption A.1), and compact operators have countable singular systems with singular values accumulating at $0$ (\cite{kress_linear_2014-2}, Section 15.4). Because $\mc{E}_{P_j;{d}}$ is a contraction, we have $\sigma_{1;d;j} \leq 1$. Furthermore, because $\mc{E}_{P_j;{d}}^*\mc{E}_{P_j;{d}}1 = 1$ and $\mc{E}_{P_j;{d}}\mc{E}_{P_j;{d}}^*1 = 1$, $\mc{E}_{P_j;{d}}$ has $1$ as a singular value and $1$ as both left and right singular functions. We order $(\sigma_{i;d;j}, u_{i;d;j}, v_{i;d;j})_{i=1}^\infty$ such that $\sigma_{1;d;j}=1$, $u_{1;d;j}=1$, $v_{1;d;j}=1$ for each $j$.

    By construction, $p_{X,Z|d;j} \in C^{\infty}(\overline{\Omega}_{X, Z})$. Consequently, $p_{X|{d};j} \in C^{\infty}(\overline{\Omega}_{X})$, and $p_{Z|{d};j} \in C^{\infty}(\overline{\Omega}_{Z})$, respectively.  Furthermore, because $p_{X|{d};j}, p_{Z|{d};j} > L_{bounded}^{-2} > 0$ by \cref{remark:marginal-conditional-density-bounds}, we have that $k_{d;j} \coloneqq p_{X,Z|{d};j} / (p_{X|{d};j}p_{Z|{d};j}) \in C^{\infty}(\overline{\Omega}_{X, Z})$. Because $\overline{\Omega}_{X, Z}$ is compact, every derivative of $k_{d;j}$ is uniformly bounded over $\overline{\Omega}_{X, Z}$. It then follows that $k_{d;j}\in W_{s,2}(\Omega_{X, Z})$, the $(s,2)$-Sobolev space, for any $s \in \mathbb{N}$. In particular, $p_{X, Z|d;j}(\cdot|{d}), k_{d;j}$ satisfy the conditions of \cref{lemma:uniform-convergence-of-kernel} for each $j$. Consequently,
    \begin{align*}
        k_{d;j}^l \coloneqq \sum_{i=1}^l \sigma_{i;d;j}u_{i;d;j}v_{i;d;j}
    \end{align*}
    converges essentially uniformly to $k_{d;j}$ as $l \rightarrow \infty$. Because $L_{bounded}^{-2} < p_{X, Z|{d};j}, p_{X|{d};j}, p_{Z|{d};j} < L_{bounded}^{2}$ by \cref{remark:marginal-conditional-density-bounds}, we have that $L_{bounded}^{-6}< k_{d;j}  < L_{bounded}^{6}$. Consequently, there exists some $l'$ such that for all $l > l'$, $k_{d;j}^l> 0$  $P_j$-a.e for all $d$. For each $l > l'$, let $L_{j}^{l} \geq 1$ be such that $\{L_{j}^{l}\}^{-1} < k_{d;j}^l(X, Z) < L_j^l$ $P_j$-a.e. for all $d$. $k_{d;j}^l$ is upper bounded because each $u_{i;d;j}$, $v_{i;d;j}$ are essentially uniformly bounded by \cref{lemma:in-range-uniformly-bounded}.

    We let
    \begin{align*}
        p^l_{X,Z|{d};j} \coloneqq k_{d;j}^lp_{X|{d};j}p_{Z|{d};j},
    \end{align*}
    which we claim is a density with respect to $\lambda_X\otimes\lambda_Z$. As $(x, z) \mapsto p^l_{X,Z|{d};j}(x, z)$ is a non-negative measurable function, it suffices to show that $\int\int p^l_{X,Z|{d};j}(x, z|{d})d\lambda(x)d\lambda(z) = 1$. Now,
    \begin{align*}
        &\int\int k_{d;j}^l(x, z)p_{X|{d};j}(x)p_{Z|{d};j}(z)d\lambda_X(x)d\lambda_Z(z)\\
        =&\int\int \sum_{i=1}^l\sigma_{i;d;j}u_{i;d;j}(x)v_{i;d;j}(z)p_{X|{d};j}(x)p_{Z|{d};j}(z)d\lambda_X(x)d\lambda_Z(z)\\
        =& \sum_{i=1}^l\sigma_{i;d;j}\langle u_{i;d;j}, 1\rangle_{L^2(P_{X|d;j})}\langle v_{i;d;j}, 1\rangle_{L^2(P_{Z|d;j})}\\
        =& 1
    \end{align*}
    where the final equality holds in light of $\sigma_{1;d;j}=1$, $u_{1;d;j}=1$, $v_{1;d;j}=1$.

    Because $k_{d;j}^l$ converges uniformly to $k_{d;j}$, we may choose $l $ large enough, such that
    \begin{align}
        \label{eq:p-l_j-j-infinity-bound}
        \|p_{X,Z|{d};j} -p^l_{X,Z|{d};j}\|_{L^\infty(P_{X,Z;j})} < \min\left(\frac{1}{j}, \frac{1}{2}L_{bounded}^{-2}L_{gap;j}\right)
    \end{align}
    We also require that $l$ exceed a threshold $L_j$, whose existence for $j$ large enough we argue below, such that for some $1 < i, i' \leq  l$, both $\langle \alpha_{P_j}(\cdot, \widetilde{d}), u_{i;\widetilde{d};j}\rangle_{L^2(P_{X|\widetilde{d};j})} \not= 0$ and $\langle r_{P_j}(\cdot, \widetilde{d}), v_{i';\widetilde{d};j}\rangle_{L^2(P_{Z|\widetilde{d};j})} \not= 0$.

    We now argue such $L_j$ exists. We argue that for $l, j$ large enough $\langle \alpha_{P_j}(\cdot, \widetilde{d}), u_{i;\widetilde{d};j}\rangle_{L^2(P_{X|\widetilde{d};j})} \not= 0$, and by symmetry the same holds for  $r_{P_j}(\cdot, \widetilde{d})$. It suffices to argue that
    \begin{align*}
        g_j \coloneqq \mc{E}_{P_j;\widetilde{d}}\alpha_{P_j}(\cdot, \widetilde{d}) = \sum_{i=1}^\infty \sigma_{i;\widetilde{d};j}\langle \alpha_{P_j}(\cdot, \widetilde{d}), u_{i;\widetilde{d};j}\rangle_{L^2(P_{X|\widetilde{d};j})}v_{i;\widetilde{d};j}
    \end{align*}
    is not a constant function for $j$ large enough.

    Now, with $(\sigma_{i;\widetilde{d};0}, u_{i;\widetilde{d};0}, v_{i;\widetilde{d};0})_{i=1}^\infty$ the singular system for $\mc{E}_{P_0;\widetilde{d}}$, ordered as above so that $u_{1;\widetilde{d};0} = v_{1;\widetilde{d};0} = 1$, let $g_0 \coloneqq \mc{E}_{P_0;\widetilde{d}}\alpha_{P_0}(\cdot, \widetilde{d}) = \sum_{i=1}^\infty \sigma_{i;\widetilde{d};0}\langle \alpha_{P_0}(\cdot, \widetilde{d}), u_{i;\widetilde{d};0}\rangle_{L^2(P_{X|\widetilde{d};0})}v_{i;\widetilde{d};0}$. We claim $g_0$ is not constant. Indeed, $\alpha_{P_0}(\cdot, \widetilde{d}) \in \mathscr{R}(\mc{E}_{P_0;\widetilde{d}}^*) \subseteq \mathscr{N}(\mc{E}_{P_0;\widetilde{d}})^\perp = \overline{\text{span}}\{u_{i;\widetilde{d};0}: i \geq 1\}$, and $\alpha_{P_0}(\cdot, \widetilde{d})$ is not constant by assumption. Hence $\langle \alpha_{P_0}(\cdot, \widetilde{d}), u_{i;\widetilde{d};0}\rangle_{L^2(P_{X|\widetilde{d};0})} \not= 0$ for some $i \geq 2$, and as $\sigma_{i;\widetilde{d};0} > 0$ and $\{v_{i;\widetilde{d};0}\}_{i \geq 1}$ is orthonormal with $v_{1;\widetilde{d};0} = 1$, the component of $g_0$ orthogonal to the constants is non-zero. We argue $g_j \rightarrow g_0$. We have that
    \begin{align*}
        \lim_{j\rightarrow \infty} \|g_j - g_0\|_{L^2(P_{Z|\widetilde{d};0})} =& \lim_{j\rightarrow \infty} \|\mc{E}_{P_j;\widetilde{d}}\alpha_{P_j}(\cdot, \widetilde{d}) -\mc{E}_{P_0;\widetilde{d}}\alpha_{P_0}(\cdot, \widetilde{d})\|_{L^2(P_{Z|\widetilde{d};0})}\\
        \leq& \lim_{j\rightarrow \infty} \|\mc{E}_{P_0;\widetilde{d}}\{\alpha_{P_j}(\cdot, \widetilde{d}) - \alpha_{P_0}(\cdot, \widetilde{d})\}\|_{L^2(P_{Z|\widetilde{d};0})}\\
        &+\lim_{j\rightarrow \infty} \|\{\mc{E}_{P_j;\widetilde{d}} - \mc{E}_{P_0;\widetilde{d}}\}\alpha_{P_j}(\cdot, \widetilde{d})\|_{L^2(P_{Z|\widetilde{d};0})}.
    \end{align*}
    Now, $\alpha_{P_j}(\cdot, \widetilde{d}) - \alpha_{P_0}(\cdot, \widetilde{d}) \rightharpoonup 0$ by \cref{lemma:weak-convergence-of-representers}. Then, $\lim_{j\rightarrow \infty} \|\mc{E}_{P_0;\widetilde{d}}\{\alpha_{P_j}(\cdot, \widetilde{d}) - \alpha_{P_0}(\cdot, \widetilde{d})\}\|_{L^2(P_{Z|\widetilde{d};0})} = 0$ because compact operators map weakly convergent sequences to strongly convergent sequences.

    On the other hand, $\mc{E}_{P_j;\widetilde{d}}$ converges to $\mc{E}_{P_0;\widetilde{d}}$ in operator norm. To see this,
    \begin{align*}
        &\lim_{j \rightarrow \infty}\sup_{h: \|h\|=1}\|\mc{E}_{P_j;\widetilde{d}}h - \mc{E}_{P_0;\widetilde{d}}h\|\\
        =& \lim_{j \rightarrow \infty}\sup_{h: \|h\|=1}\left\|\int h(x) \left\{p_{X|Z, \widetilde{d};j}(x, \cdot) -  p_{X|Z, \widetilde{d};0}(x, \cdot)\right\}d\lambda(x)\right\|_{L^2(P_{Z|\widetilde{d};0})}\\
        \leq & \lim_{j \rightarrow \infty}\sup_{h: \|h\|=1}\left\{\int h(x)^2d\lambda(x)\right\}^{1/2}\left\|\left\{\int \left\{p_{X|Z, \widetilde{d};j}(x, \cdot) -  p_{X|Z, \widetilde{d};0}(x, \cdot)\right\}^2\right\}^{1/2}\right\|_{L^2(P_{Z|\widetilde{d};0})}\\
        \leq & L_{bounded}^2\lim_{j \rightarrow \infty}\sup_{h: \|h\|=1}\left\{\int h(x)^2p_{X|\widetilde{d};0}(x)d\lambda(x)\right\}^{1/2}\left\|p_{X|Z, \widetilde{d};j} -  p_{X|Z, \widetilde{d};0}\right\|_{L^2(\lambda_X\otimes\lambda_Z)}\\
        =& 0.
    \end{align*}
    Furthermore, because weakly convergent sequences are bounded, operator convergence of $\mc{E}_{P_j;\widetilde{d}}$ to $ \mc{E}_{P_0;\widetilde{d}}$ implies
    \begin{align*}
        \lim_{j\rightarrow \infty} \|\{\mc{E}_{P_j;\widetilde{d}} - \mc{E}_{P_0;\widetilde{d}}\}\alpha_{P_j}(\cdot, \widetilde{d})\|_{L^2(P_{Z|\widetilde{d};0})} \leq \lim_{j\rightarrow \infty} \|\{\mc{E}_{P_j;\widetilde{d}} - \mc{E}_{P_0;\widetilde{d}}\}\|\sup_{j}\|\alpha_{P_j}(\cdot, \widetilde{d})\|_{L^2(P_{X|\widetilde{d};0})} = 0
    \end{align*}
    Then, $g_j \rightarrow g_0$, and $g_0$ is not a constant function $P_0$-a.e. Because the $L^2(P_{Z|\widetilde{d};j})$ and $L^2(P_{Z|\widetilde{d};0})$ norms are equivalent uniformly in $j$, there exists $J_\alpha$ such that $g_j$ is not $P_j$-a.e. constant for all $j > J_\alpha$. For such $j$, as $v_{1;\widetilde{d};j} = 1$ and $\{v_{i;\widetilde{d};j}\}_{i \geq 1}$ is orthonormal, the expansion of $g_j$ above forces $\langle \alpha_{P_j}(\cdot, \widetilde{d}), u_{i;\widetilde{d};j}\rangle_{L^2(P_{X|\widetilde{d};j})} \not= 0$ for some finite $i \geq 2$. The same argument applied to $\mc{E}_{P_j;\widetilde{d}}^*r_{P_j}(\cdot, \widetilde{d})$ gives $J_r$ such that, for all $j > J_r$, $\langle r_{P_j}(\cdot, \widetilde{d}), v_{i';\widetilde{d};j}\rangle_{L^2(P_{Z|\widetilde{d};j})} \not= 0$ for some finite $i' \geq 2$. Let $J \coloneqq \max(J_\alpha, J_r)$. For each $j > J$, we let $L_j$ be the larger of the smallest such $i$ and the smallest such $i'$, and let $l_j$ be the first $l \geq L_j$ satisfying \eqref{eq:p-l_j-j-infinity-bound}. In what follows we consider only $j > J$, re-indexing the sequence if necessary.

    For $d \not = \widetilde{d}$ and $i \leq l_j$ we let $\widetilde{u}_{i;d;j} \coloneqq {u}_{i;d;j}$, $\widetilde{v}_{i;d;j} \coloneqq {v}_{i;d;j}$, $\widetilde{\sigma}_{i;d;j} \coloneqq {\sigma}_{i;d;j}$. If
    \begin{align*}
        \alpha_{P_j}(\cdot, {d}) - \sum_{i\leq l_j:\widetilde{\sigma}_{i;{d};j}>0}\langle \widetilde{u}_{i;d;j}, \alpha_{P_j}(\cdot, {d})\rangle_{L^2(P_{X|d;j})}\widetilde{u}_{i;d;j} \not = 0,
    \end{align*}
    we let
    \begin{align}
        \label{eq:residual-projection-generic-d}
        \widetilde{u}_{l_j+1;d;j} \coloneqq& \frac{\alpha_{P_j}(\cdot, {d}) - \sum_{i\leq l_j:\widetilde{\sigma}_{i;{d};j}>0}\langle \widetilde{u}_{i;d;j}, \alpha_{P_j}(\cdot, {d})\rangle_{L^2(P_{X|d;j})}\widetilde{u}_{i;d;j}}{\|\alpha_{P_j}(\cdot, {d}) - \sum_{i\leq l_j:\widetilde{\sigma}_{i;{d};j}>0}\langle \widetilde{u}_{i;d;j}, \alpha_{P_j}(\cdot, {d})\rangle_{L^2(P_{X|d;j})}\widetilde{u}_{i;d;j}\|_{L^2(P_{X|d;j})}},
    \end{align}
    $\widetilde{u}_{l_j+1;{d};j}$ equal to an arbitrary unit vector essentially uniformly bounded orthogonal to $\widetilde{u}_{i;{d};j}$ for $i \leq l_j$ otherwise. We do similarly for $\widetilde{v}_{l_j+1;d;j}$. Recall that $\alpha_{P_j}(\cdot, {d})$, $r_{P_j}(\cdot,{d})$ are essentially uniformly bounded by \cref{lemma:l1-bounded-functional}, and hence $\|\widetilde{u}_{l_j+1;d;j} \|_{L^\infty(P_{X;j})}, \|\widetilde{v}_{l_j+1;d;j} \|_{L^\infty(P_{Z;j})} < \infty$.

    We then choose $\widetilde{\sigma}_{l_j+1;d;j}  > 0$ sufficiently small such that
    \begin{align}
        \label{eq:sigma-1-d-bound}
        \widetilde{\sigma}_{l_j + 1;d;j}\|\widetilde{u}_{l_j+1;d;j}\widetilde{v}_{l_j+1;d;j}\|_{L^\infty(P_{X,Z;j})} < \min\left(L_{bounded}^{-4}\min\left(\frac{1}{j}, \frac{1}{2}L_{bounded}^{-2}L_{gap;j}\right), \{L_{j}^{l_j}\}^{-1}\right).
    \end{align}
    Then, letting
    \begin{align*}
        \widetilde{k}_{d;j} \coloneqq \sum_{i=1}^{l_j+1} \widetilde{\sigma}_{i;d;j}\widetilde{u}_{i;d;j}\widetilde{v}_{i;d;j}
    \end{align*}
    \begin{align*}
        \nonumber \|k_{d;j}^{l_j} - \widetilde{k}_{d;j}\|_{L^{\infty}(P_{X,Z;j})} < \min\left(L_{bounded}^{-4}\min\left(\frac{1}{j}, \frac{1}{2}L_{bounded}^{-2}L_{gap;j}\right), \{L_{j}^{l_j}\}^{-1}\right).
    \end{align*}
    Because ${k}^{l_j}_{d;j} > \{L_{j}^{l_j}\}^{-1} > 0$, $\|k_{d;j}^{l_j} - \widetilde{k}_{d;j}\|_{L^{\infty}(P_{X,Z;j})} < \{L_{j}^{l_j}\}^{-1}$ implies that $\widetilde{k}_{d;j} > 0$.
    Defining
    \begin{align*}
        \widetilde{p}_{X,Z|{d};j}(x, z) \coloneqq \widetilde{k}_{d;j}(x, z) p_{X|{d};j}(x)p_{Z|{d};j}(z)
    \end{align*}
    then gives that
    \begin{align*}
        \|p_{X, Z|{d};j}^{l_j} - \widetilde{p}_{X, Z|{d};j}\|_{L^{\infty}(P_{X,Z;j})}
        &\leq\|k_{d;j}^{l_j} - \widetilde{k}_{d;j}\|_{L^{\infty}(P_{X,Z;j})}\|p_{X|{d};j}\|_{L^\infty(P_{X;j})}\|p_{Z|{d};j}\|_{L^\infty(P_{Z;j})} \\
        &< \min\left(\frac{1}{j}, \frac{1}{2}L_{bounded}^{-2}L_{gap;j}\right)
    \end{align*}
    where we have used that $\|p_{X|{d};j}\|_{L^\infty(P_{X;j})}, \|p_{Z|{d};j} \|_{L^\infty(P_{Z;j})} < L_{bounded}^2$ by application of \cref{remark:marginal-conditional-density-bounds}.

    We also have that $\widetilde{p}_{X,Z|{d};j} > 0$. We further argue that $\widetilde{p}_{X,Z|{d};j}$ is a density for $(X, Z)$. Indeed,
    \begin{align*}
        &\int\int \widetilde{k}_{d;j}(x, z)p_{X|{d};j}(x, z)p_{Z|{d};j}(z)d\lambda_X(x)d\lambda_Z(z)\\
        =&\int\int \sum_{i=1}^{l_j+1}\widetilde{\sigma}_{i;d;j}\widetilde{u}_{i;d;j}(x)\widetilde{v}_{i;d;j}(z)p_{X|{d};j}(x)p_{Z|{d};j}(z)d\lambda_X(x)d\lambda_Z(z)\\
        =& \sum_{i=1}^{l_j+1}\widetilde{\sigma}_{i;d;j}\langle \widetilde{u}_{i;d;j}, 1\rangle_{L^2(P_{X|d;j})}\langle \widetilde{v}_{i;d;j}, 1\rangle_{L^2(P_{Z|d;j})}\\
        =& 1
    \end{align*}
    which again follows because $\widetilde{\sigma}_{1;d;j} = 1$, $\widetilde{u}_{1;d;j} = 1$, $\widetilde{v}_{1;d;j}=1$. Moreover, $\widetilde{p}_{X|d;j} = p_{X|d;j}$ and $\widetilde{p}_{Z|d;j} = p_{Z|d;j}$, by the same calculation given below for $d = \widetilde{d}$, because $\widetilde{u}_{i;d;j}$ and $\widetilde{v}_{i;d;j}$ are orthogonal to $\widetilde{u}_{1;d;j} = 1$ and $\widetilde{v}_{1;d;j} = 1$, respectively, for each $2 \leq i \leq l_j + 1$.

    Now, for $d = \widetilde{d}$, let $i_{\alpha;j} , i_{r;j} > 1$ be the indexes of the first non-constant singular functions $u_{i_{\alpha;j};\widetilde{d};j}$, $v_{i_{r;j};\widetilde{d};j}$ such that
    \begin{align}
        \label{eq:singular-functions-not-orthogonal-continuous}
        \langle u_{i_{\alpha;j};\widetilde{d};j}, \alpha_{P_j}(\cdot, \widetilde{d})\rangle_{L^2(P_{X|\widetilde{d};j})} \not = 0 && \text{ and } &&  \langle v_{i_{r;j};\widetilde{d};j}, r_{P_j}(\cdot,\widetilde{d})\rangle_{L^2(P_{Z|\widetilde{d};j})} \not = 0.
    \end{align}
    We define
    \begin{align}
        \label{eq:rotation-continuous-alpha}
        \nonumber \widetilde{u}_{i_{\alpha;j};\widetilde{d};j}(\delta) \coloneqq&  {u}_{i_{\alpha;j} ;\widetilde{d};j}\cos(\delta) - {u}_{l_j+1;\widetilde{d};j}\sin(\delta)\\
        \widetilde{u}_{l_j+1;\widetilde{d};j}(\delta) \coloneqq&  {u}_{i_{\alpha;j};\widetilde{d};j}\sin(\delta) + {u}_{l_j+1;\widetilde{d};j}\cos(\delta).
    \end{align}
    and
    \begin{align}
        \label{eq:rotation-continuous-r}
        \nonumber \widetilde{v}_{i_{r;j};\widetilde{d};j}(\delta) \coloneqq&  {v}_{i_{r;j};\widetilde{d};j}\cos(\delta) - {v}_{l_j+1;\widetilde{d};j}\sin(\delta)\\
        \widetilde{v}_{l_j+1;\widetilde{d};j}(\delta) \coloneqq&  {v}_{i_{r;j} ;\widetilde{d};j}\sin(\delta) + {v}_{l_j+1;\widetilde{d};j}\cos(\delta).
    \end{align}
    If $\sigma_{l_j+1}=0$, we let ${u}_{l_j+1;\widetilde{d};j}$, ${v}_{l_j+1;\widetilde{d};j}$ in the above be arbitrary unit length essentially uniformly bounded vectors orthogonal to ${u}_{i;\widetilde{d};j}, {v}_{i;\widetilde{d};j}$ for each $i \leq l_{j}$, respectively.

    Then, $\lim_{\delta \rightarrow 0}\|\widetilde{u}_{i_{\alpha;j};\widetilde{d};j}(\delta) - {u}_{i_{\alpha;j};\widetilde{d};j} \|_{L^\infty(P_{X;j})} = 0$ and similarly for $\widetilde{u}_{l_j+1;\widetilde{d};j}(\delta)$. Furthermore, by \eqref{eq:singular-functions-not-orthogonal-continuous},
    the rotations ensure that for all $\delta \not= 0$ small enough,
    \begin{align}
        \label{eq:rotation-not-zero-continuous}
        \langle \widetilde{u}_{l_j+1 ;\widetilde{d};j}(\delta),\alpha_{P_j}(\cdot, \widetilde{d}) \rangle_{L^2(P_{X|\widetilde{d};j})} \not=0,  &&\langle \widetilde{v}_{l_j+1 ;\widetilde{d};j}(\delta),r_{P_j}(\cdot, \widetilde{d}) \rangle_{L^2(P_{Z|\widetilde{d};j})} \not=0.
    \end{align}
    We, for each $j$, choose $\delta_j \not= 0$ such that \eqref{eq:rotation-not-zero-continuous} holds and
    \begin{align}
        \label{eq:i_alpha_j-bound-continuous}
        \sigma_{i_{\alpha;j};\widetilde{d};j}\|\widetilde{u}_{i_{\alpha;j};\widetilde{d};j}(\delta_j)\widetilde{v}_{i_{\alpha;j};\widetilde{d};j}(\delta_j) - {u}_{i_{\alpha;j};\widetilde{d};j}{v}_{i_{\alpha;j};\widetilde{d};j}\|_{L^\infty(P_{X,Z;j})}
        <\frac{1}{4}\min\left(L_{bounded}^{-4}\min\left(\frac{1}{j}, \frac{1}{2}L_{bounded}^{-2}L_{gap;j}\right), \{L_{j}^{l_j}\}^{-1}\right)
    \end{align}
    and
    \begin{align}
        \label{eq:i_r_j-bound-continuous}
        &\sigma_{i_{r;j};\widetilde{d};j}\|\widetilde{u}_{i_{r;j};\widetilde{d};j}(\delta_j)\widetilde{v}_{i_{r;j};\widetilde{d};j}(\delta_j) - {u}_{i_{r;j};\widetilde{d};j}{v}_{i_{r;j};\widetilde{d};j}\|_{L^\infty(P_{X,Z;j})}
        <\frac{1}{4}\min\left(L_{bounded}^{-4}\min\left(\frac{1}{j}, \frac{1}{2}L_{bounded}^{-2}L_{gap;j}\right), \{L_{j}^{l_j}\}^{-1}\right),
    \end{align}
    where $\widetilde{v}_{i_{\alpha;j};\widetilde{d};j}(\delta_j) \coloneqq {v}_{i_{\alpha;j};\widetilde{d};j}$ and $\widetilde{u}_{i_{r;j};\widetilde{d};j}(\delta_j) \coloneqq {u}_{i_{r;j};\widetilde{d};j}$ if $i_{\alpha;j} \not= i_{r;j}$.
    Note that
    \begin{align*}
        \|u_{i_{\alpha;j};\widetilde{d};j}\|_{L^\infty(P_{X;j})}, \|v_{i_{r;j};\widetilde{d};j}\|_{L^\infty(P_{Z;j})}, \|u_{l_j+1;\widetilde{d};j}\|_{L^\infty(P_{X;j})}, \|v_{l_j+1;\widetilde{d};j}\|_{L^\infty(P_{Z;j})} < \infty
    \end{align*}
    in light of the fact that $u_{i;\widetilde{d};j} \in \mathscr{R}(\mc{E}^*_{P_j;\widetilde{d}}) $, $v_{i;\widetilde{d};j} \in \mathscr{R}(\mc{E}_{P_j;\widetilde{d}}) $ for each $i$ and \cref{lemma:in-range-uniformly-bounded}.

    Having chosen $\delta_j$, we define
    \begin{align*}
        \widetilde{u}_{i_{\alpha;j};\widetilde{d};j} \coloneqq& \widetilde{u}_{i_{\alpha;j};\widetilde{d};j}(\delta_j), &\widetilde{v}_{i_{r;j};\widetilde{d};j} \coloneqq& \widetilde{v}_{i_{r;j};\widetilde{d};j}(\delta_j)\\
        \widetilde{u}_{l_j + 1;\widetilde{d};j} \coloneqq& \widetilde{u}_{l_j + 1;\widetilde{d};j}(\delta_j), &\widetilde{v}_{l_j + 1;\widetilde{d};j} \coloneqq&\widetilde{v}_{l_j + 1;\widetilde{d};j}(\delta_j),
    \end{align*}
    and let $\widetilde{u}_{i; \widetilde{d};j} \coloneqq u_{i;\widetilde{d};j}$ for $i \in \{1, \dots, l_j\}\setminus \{i_{\alpha;j}\}$, $\widetilde{v}_{i; \widetilde{d};j} \coloneqq v_{i;\widetilde{d};j}$ for $i \in \{1, \dots, l_j\}\setminus \{i_{r;j}\}$.

    Now, if $\alpha_{P_j}(\cdot, \widetilde{d}) - \sum_{i\leq l_j+1:\widetilde{\sigma}_{i;\widetilde{d};j}>0}\langle \widetilde{u}_{i; \widetilde{d};j}, \alpha_{P_j}(\cdot, \widetilde{d})\rangle_{L^2(P_{X|\widetilde{d};j})}\widetilde{u}_{i; \widetilde{d};j} \not = 0$, we let
    \begin{align}
        \label{eq:residual-projection-1}
        \widetilde{u}_{l_j+2;\widetilde{d};j} \coloneqq& \frac{\alpha_{P_j}(\cdot, \widetilde{d}) - \sum_{i\leq l_j+1:\widetilde{\sigma}_{i;\widetilde{d};j}>0}\langle \widetilde{u}_{i; \widetilde{d};j}, \alpha_{P_j}(\cdot, \widetilde{d})\rangle_{L^2(P_{X|\widetilde{d};j})}\widetilde{u}_{i; \widetilde{d};j}}{\|\alpha_{P_j}(\cdot, \widetilde{d}) - \sum_{i\leq l_j+1:\widetilde{\sigma}_{i;\widetilde{d};j}>0}\langle \widetilde{u}_{i; \widetilde{d};j}, \alpha_{P_j}(\cdot, \widetilde{d})\rangle_{L^2(P_{X|\widetilde{d};j})}\widetilde{u}_{i; \widetilde{d};j}\|_{L^2(P_{X|\widetilde{d};j})}},
    \end{align}
    and set $\widetilde{u}_{l_j+2;\widetilde{d};j}$ equal to an arbitrary unit vector essentially uniformly bounded orthogonal to $\widetilde{u}_{i;\widetilde{d};j}$ for $i \leq l_j+1$, otherwise. We do similarly for $\widetilde{v}_{l_j+2;\widetilde{d};j}$. Recall that $\alpha_{P_j}(\cdot, \widetilde{d})$, $r_{P_j}(\cdot,\widetilde{d})$ are essentially uniformly bounded by \cref{lemma:l1-bounded-functional}, and hence $\|\widetilde{u}_{l_j+2;\widetilde{d};j} \|_{L^\infty(P_{X;j})}, \|\widetilde{v}_{l_j+2;\widetilde{d};j} \|_{L^\infty(P_{Z;j})} < \infty$.

    We let $\widetilde{\sigma}_{i;\widetilde{d};j} \coloneqq {\sigma}_{i;\widetilde{d};j}$ for $i \leq l_j$ and let $\widetilde{\sigma}_{l_j + 1;\widetilde{d};j},\widetilde{\sigma}_{l_j + 2;\widetilde{d};j}  > 0$ to be defined below. We define
    \begin{align*}
        \widetilde{k}_{\widetilde{d};j}  \coloneqq \sum_{i=1}^{l_j + 2} \widetilde{\sigma}_{i;\widetilde{d};j}\widetilde{u}_{i; \widetilde{d};j}\widetilde{v}_{i; \widetilde{d};j}.
    \end{align*}
    It follows that
    \begin{align}
        \label{eq:bound-uniform-norm-rotated-kernel}
        \nonumber \|k_{\widetilde{d};j}^{l_j} - \widetilde{k}_{\widetilde{d};j}\|_{L^{\infty}(P_{X,Z;j})} \leq& \sigma_{i_{\alpha;j};\widetilde{d};j}\|\widetilde{u}_{i_{\alpha;j};\widetilde{d};j}\widetilde{v}_{i_{\alpha;j};\widetilde{d};j} - {u}_{i_{\alpha;j};\widetilde{d};j}{v}_{i_{\alpha;j};\widetilde{d};j}\|_{L^\infty(P_{X,Z;j})}\\
         \nonumber &+\sigma_{i_{r;j};\widetilde{d};j}\|\widetilde{u}_{i_{r;j};\widetilde{d};j}\widetilde{v}_{i_{r;j};\widetilde{d};j} - {u}_{i_{r;j};\widetilde{d};j}{v}_{i_{r;j};\widetilde{d};j}\|_{L^\infty(P_{X,Z;j})}\\
        \nonumber&+\widetilde{\sigma}_{l_j + 1;\widetilde{d};j}\|\widetilde{u}_{l_j+1;\widetilde{d};j}\widetilde{v}_{l_j+1;\widetilde{d};j}\|_{L^\infty(P_{X,Z;j})}\\
        &+\widetilde{\sigma}_{l_j + 2;\widetilde{d};j}\|\widetilde{u}_{l_j+2;\widetilde{d};j}\widetilde{v}_{l_j+2;\widetilde{d};j}\|_{L^\infty(P_{X,Z;j})}.
    \end{align}

    Now,
    \begin{align*}
        \|\widetilde{u}_{l_j+1;\widetilde{d};j}\widetilde{v}_{l_j+1;\widetilde{d};j}\|_{L^\infty(P_{X,Z;j})} \leq& \|\widetilde{u}_{l_j+1;\widetilde{d};j}\|_{L^\infty(P_{X;j})}\|\widetilde{v}_{l_j+1;\widetilde{d};j}\|_{L^\infty(P_{Z;j})} < \infty\\
        \|\widetilde{u}_{l_j+2;\widetilde{d};j}\widetilde{v}_{l_j+2;\widetilde{d};j}\|_{L^\infty(P_{X,Z;j})} \leq& \|\widetilde{u}_{l_j+2;\widetilde{d};j}\|_{L^\infty(P_{X;j})}\|\widetilde{v}_{l_j+2;\widetilde{d};j}\|_{L^\infty(P_{Z;j})} < \infty.
    \end{align*}
    We then choose $\widetilde{\sigma}_{l_j+2;\widetilde{d};j}  > 0$ sufficiently small such that
    \begin{align}
        \label{eq:sigma-2-bound}
        \widetilde{\sigma}_{l_j + 2;\widetilde{d};j}\|\widetilde{u}_{l_j+2;\widetilde{d};j}\widetilde{v}_{l_j+2;\widetilde{d};j}\|_{L^\infty(P_{X,Z;j})} < \frac{1}{4}\min\left(L_{bounded}^{-4}\min\left(\frac{1}{j}, \frac{1}{2}L_{bounded}^{-2}L_{gap;j}\right), \{L_{j}^{l_j}\}^{-1}\right),
    \end{align}
    and choose $ \widetilde{\sigma}_{l_j + 1;\widetilde{d};j} > 0$ to satisfy both that
    \begin{align}
        \label{eq:sigma_1_bound}
        \widetilde{\sigma}_{l_j + 1;\widetilde{d};j}\|\widetilde{u}_{l_j+1;\widetilde{d};j}\widetilde{v}_{l_j+1;\widetilde{d};j}\|_{L^\infty(P_{X,Z;j})} < \frac{1}{4}\min\left(L_{bounded}^{-4}\min\left(\frac{1}{j}, \frac{1}{2}L_{bounded}^{-2}L_{gap;j}\right), \{L_{j}^{l_j}\}^{-1}\right)
    \end{align}
    and
    \begin{align}
        \label{eq:sigma_discontinuity-def}
        \widetilde{\sigma}_{l_j + 1;\widetilde{d};j} < \frac{P_j(D=\widetilde{d})}{j + |\kappa_j|}|\langle \widetilde{u}_{l_j + 1;\widetilde{d};j}, \alpha_{P_j}(\cdot, \widetilde{d})\rangle_{L^2(P_{X|\widetilde{d};j})}\langle \widetilde{v}_{l_j + 1;\widetilde{d};j}, r_{P_j}(\cdot, \widetilde{d})\rangle_{L^2(P_{Z|\widetilde{d};j})}|
    \end{align}
    where
    \begin{align*}
        \kappa_j \coloneqq& P_j(D=\widetilde{d})\sum_{i\leq l_j:\widetilde{\sigma}_{i;\widetilde{d};j}>0}\widetilde{\sigma}_{i;\widetilde{d};j}^{-1}\langle \widetilde{u}_{i; \widetilde{d};j}, \alpha_{P_j}(\cdot, \widetilde{d})\rangle_{L^2(P_{X|\widetilde{d};j})}\langle \widetilde{v}_{i; \widetilde{d};j}, r_{P_j}(\cdot, \widetilde{d})\rangle_{L^2(P_{Z|\widetilde{d};j})}\\
        &+P_j(D=\widetilde{d})\widetilde{\sigma}_{l_j+2;\widetilde{d};j}^{-1}\langle \widetilde{u}_{l_j+2;\widetilde{d};j}, \alpha_{P_j}(\cdot,\widetilde{d})\rangle_{L^2(P_{X|\widetilde{d};j})}\langle \widetilde{v}_{l_j+2;\widetilde{d};j}, r_{P_j}(\cdot, \widetilde{d})\rangle_{L^2(P_{Z|\widetilde{d};j})}\\
        &+ \sum_{d \in \textsf{Supp}[D;P_j]\setminus \{\widetilde{d}\}}P_j(D=d)\sum_{i\leq l_j+1: \widetilde{\sigma}_{i;{d};j}>0}\widetilde{\sigma}_{i;{d};j}^{-1}\langle \widetilde{u}_{i; {d};j}, \alpha_{P_j}(\cdot, {d})\rangle_{L^2(P_{X|d;j})}\\
        &\qquad\times\langle \widetilde{v}_{i; {d};j}, r_{P_j}(\cdot, {d})\rangle_{L^2(P_{Z|d;j})}
    \end{align*}
    The bounds in \eqref{eq:i_alpha_j-bound-continuous}, \eqref{eq:i_r_j-bound-continuous},  \eqref{eq:sigma-2-bound}, \eqref{eq:sigma_1_bound} combined with \eqref{eq:bound-uniform-norm-rotated-kernel} ensure that
    \begin{align*}
        \nonumber \|k_{\widetilde{d};j}^{l_j} - \widetilde{k}_{\widetilde{d};j}\|_{L^{\infty}(P_{X,Z;j})} < \min\left(L_{bounded}^{-4}\min\left(\frac{1}{j}, \frac{1}{2}L_{bounded}^{-2}L_{gap;j}\right), \{L_{j}^{l_j}\}^{-1}\right).
    \end{align*}
    Because ${k}^{l_j}_{\widetilde{d};j} > \{L_{j}^{l_j}\}^{-1} > 0$, $\|k_{\widetilde{d};j}^{l_j} - \widetilde{k}_{\widetilde{d};j}\|_{L^{\infty}(P_{X,Z;j})} < \{L_{j}^{l_j}\}^{-1}$ implies that $\widetilde{k}_{\widetilde{d};j} > 0$.
    Defining
    \begin{align*}
        \widetilde{p}_{X,Z|\widetilde{d};j}(x, z) \coloneqq \widetilde{k}_{\widetilde{d};j}(x, z) p_{X|\widetilde{d};j}(x)p_{Z|\widetilde{d};j}(z)
    \end{align*}
    then gives that
    \begin{align}
        \label{eq:p-l_j-tilde-infinity-bound}
        \begin{split}
        \|p_{X, Z|\widetilde{d};j}^{l_j} - \widetilde{p}_{X, Z|\widetilde{d};j}\|_{L^{\infty}(P_{X,Z;j})}
        &\leq\|k_{\widetilde{d};j}^{l_j} - \widetilde{k}_{\widetilde{d};j}\|_{L^{\infty}(P_{X,Z;j})}\|p_{X|\widetilde{d};j}\|_{L^\infty(P_{X;j})}\|p_{Z|\widetilde{d};j}\|_{L^\infty(P_{Z;j})} \\
        &< \min\left(\frac{1}{j}, \frac{1}{2}L_{bounded}^{-2}L_{gap;j}\right)
        \end{split}
    \end{align}
    where we have again used that $\|p_{X|{d};j}\|_{L^\infty(P_{X;j})}, \|p_{Z|{d};j} \|_{L^\infty(P_{Z;j})} < L_{bounded}^2$ by application of \cref{remark:marginal-conditional-density-bounds}.

    We also have that $\widetilde{p}_{X,Z|\widetilde{d};j} > 0$. We further argue that $\widetilde{p}_{X,Z|\widetilde{d};j}$ is a density for $(X, Z)$. Indeed,
    \begin{align*}
        &\int\int \widetilde{k}_{\widetilde{d};j}(x, z)p_{X|\widetilde{d};j}(x)p_{Z|\widetilde{d};j}(z)d\lambda_X(x)d\lambda_Z(z)\\
        =&\int\int \sum_{i=1}^{l_j+2}\widetilde{\sigma}_{i;\widetilde{d};j}\widetilde{u}_{i; \widetilde{d};j}(x)\widetilde{v}_{i; \widetilde{d};j}(z)p_{X|\widetilde{d};j}(x)p_{Z|\widetilde{d};j}(z)d\lambda_X(x)d\lambda_Z(z)\\
        =& \sum_{i=1}^{l_j+2}\widetilde{\sigma}_{i;\widetilde{d};j}\langle \widetilde{u}_{i; \widetilde{d};j}, 1\rangle_{L^2(P_{X|\widetilde{d};j})}\langle \widetilde{v}_{i; \widetilde{d};j}, 1\rangle_{L^2(P_{Z|\widetilde{d};j})}\\
        =& 1
    \end{align*}
    which again follows because $\widetilde{\sigma}_{1;\widetilde{d};j} = 1$, $\widetilde{u}_{1;\widetilde{d};j} = 1$, $\widetilde{v}_{1;\widetilde{d};j}=1$. Furthermore,
    \begin{align*}
        \widetilde{p}_{X|\widetilde{d};j}(x) = &\int \widetilde{k}_{\widetilde{d};j}(x, z)p_{X|\widetilde{d};j}(x)p_{Z|\widetilde{d};j}(z)d\lambda_Z(z)\\
        =&p_{X|\widetilde{d};j}(x)\int \sum_{i=1}^{l_j+2}\widetilde{\sigma}_{i;\widetilde{d};j}\widetilde{u}_{i; \widetilde{d};j}(x)\widetilde{v}_{i; \widetilde{d};j}(z)p_{Z|\widetilde{d};j}(z)d\lambda_Z(z)\\
        =& p_{X|\widetilde{d};j}(x)\sum_{i=1}^{l_j+2}\widetilde{\sigma}_{i;\widetilde{d};j}\widetilde{u}_{i; \widetilde{d};j}(x)\langle \widetilde{v}_{i; \widetilde{d};j}, 1\rangle_{L^2(P_{Z|\widetilde{d};j})}\\
        =& p_{X|\widetilde{d};j}(x)
    \end{align*}
    and similarly $\widetilde{p}_{Z|\widetilde{d};j} = {p}_{Z|\widetilde{d};j}$.

    Now, for each $d$,
    \begin{align*}
        \|\widetilde{p}_{X,Z|{d};j} - p_{X, Z|{d};j}\|_{L^\infty(P_{X,Z;j})}
        &\leq \|\widetilde{p}_{X,Z|{d};j} - p^{l_j}_{X, Z|{d};j}\|_{L^\infty(P_{X,Z;j})} + \|p^{l_j}_{X, Z|{d};j} - p_{X, Z|{d};j}\|_{L^\infty(P_{X,Z;j})} \\
        &< L_{bounded}^{-2}L_{gap;j}
    \end{align*}
    by the inequalities in \eqref{eq:p-l_j-tilde-infinity-bound} and \eqref{eq:p-l_j-j-infinity-bound}.

    Let $\widetilde{p}_{X, Z, D;j}(x, z, d) \coloneqq \widetilde{p}_{X,Z|{d};j}(x, z)p_{D;j}(d)$ denote the corresponding joint density of $(X, Z, D)$. Because $p_{D;j}(d) < L_{bounded}$ for each $d$ by \cref{remark:marginal-conditional-density-bounds}, and because $\widetilde{p}_{X, Z, D;j}$ and $p_{X, Z, D;j}$ both vanish off the support of $P_j$, the above display gives
    \begin{align*}
        \|\widetilde{p}_{X, Z, D;j} - p_{X, Z, D;j}\|_{L^\infty(\lambda_X\otimes\lambda_Z\otimes\lambda_D)} = \max_{d \in \textsf{Supp}[D;P_j]} p_{D;j}(d)\|\widetilde{p}_{X,Z|{d};j} - p_{X, Z|{d};j}\|_{L^\infty(P_{X,Z;j})} < L_{bounded}^{-1}L_{gap;j}.
    \end{align*}
    With $P \coloneqq P_j$, the threshold in \cref{lemma:constant-riesz-representers} is $L_{bounded}^{-1}L_{gap;j}\{\lambda_X(\mathcal X)\lambda_Z(\mathcal Z)\lambda_D(\mc{D})\}^{-1} = L_{bounded}^{-1}L_{gap;j}$, because the reference measures are normalized. Hence the conditions of \cref{lemma:constant-riesz-representers} are satisfied with $\widetilde{p}_{X, Z, D;j}$ the density $p'$ and $P_j$ the reference distribution $P$ in that Lemma. We let $\widetilde{P}_j$ be $P^\dagger$ of that Lemma. By that Lemma,
    \begin{enumerate}
        \item $P_j$ and $\widetilde{P}_j$ are mutually absolutely continuous
        \item $L_{bounded}^{-1} < \essinf \widetilde{p}_j \leq \esssup \widetilde{p}_j < L_{bounded}$
        \item $\widetilde{P}_j$ satisfies \cref{as:continuous-linear-maps} and $\alpha_{\widetilde{P}_j}  = \alpha_{{P}_j} $, $r_{\widetilde{P}_j}  = r_{{P}_j} $
        \item $H^2(\widetilde{P}_j, P_j) \leq \frac{1}{2}\|\widetilde{p}_{X, Z, D;j} - p_{X, Z, D;j}\|_{L^1(\lambda_X\otimes\lambda_Z\otimes\lambda_D)}$.
    \end{enumerate}

    Now,
    \begin{align*}
        &\frac{1}{2}\|\widetilde{p}_{X, Z, D;j} - p_{X, Z, D;j}\|_{L^1(\lambda_X\otimes\lambda_Z\otimes\lambda_D)} \\
        \leq& L'\|\widetilde{p}_{X, Z, D;j} - p_{X, Z, D;j}\|_{L^\infty(\lambda_X\otimes\lambda_Z\otimes\lambda_D)}\\
        =&L'\|\widetilde{p}_{X, Z, D;j} - p_{X, Z, D;j}\|_{L^\infty(P_{X,Z,D;j})}\\
        \leq&\sum_{d \in \textsf{Supp}[D;P_j]}L'\|\widetilde{p}_{X, Z|{d};j} - p_{X, Z|{d};j}\|_{L^\infty(P_{X,Z,D;j})}P_j(D=d)\\
        \leq & \sum_{d \in \textsf{Supp}[D;P_j]}P_j(D=d)L' \{\|\widetilde{p}_{X, Z|{d};j} - p^{l_j}_{X, Z|{d};j}\|_{L^\infty(P_{X,Z;j})} + \|{p}^{l_j}_{X, Z|{d};j} - p_{X, Z|{d};j}\|_{L^\infty(P_{X,Z;j})}\}\\
        \leq& 2L'/j
    \end{align*}
    for some constant $L'$ not depending on $j$.  The final inequality follows from \eqref{eq:p-l_j-j-infinity-bound} and \eqref{eq:p-l_j-tilde-infinity-bound}. Consequently,
    \begin{align*}
        H^2(\widetilde{P}_j, P_0) \leq H^2(\widetilde{P}_j, P_j) +2  H(\widetilde{P}_j, P_j)H({P}_j, P_0) + H^2({P}_j, P_0) <  2L'/j + 2\sqrt{2L'}/j + 1/j.
    \end{align*}
    and so $\lim_{j\rightarrow \infty}H^2(\widetilde{P}_j, P_0) = 0$.

    We now argue $\widetilde{P}_j$ belongs to $\mc{P}$ for each $j$, i.e. Assumptions \ref{as:continuous-linear-maps}, \ref{as:solutions-exist}, \ref{as:beta-source} (for all $\beta_h, \beta_q \in [0, \infty)$), \ref{as:mutual-absolute-continuity} (relative to $P_0$), and \ref{as:dominating-measure-support} hold.  \cref{as:continuous-linear-maps} holds by the construction in \cref{lemma:constant-riesz-representers} and because it holds at $P_j$. \cref{as:solutions-exist} holds by the choice of $\widetilde{u}_{l_j + 1;d;j}$ and $\widetilde{v}_{l_j + 1;{d};j}$ for $d \not=\widetilde{d}$,  $\widetilde{u}_{l_j + 2;\widetilde{d};j}$ and $\widetilde{v}_{l_j + 2;\widetilde{d};j}$ in \eqref{eq:residual-projection-generic-d}, \eqref{eq:residual-projection-1}. Furthermore, because $\mc{E}_{\widetilde{P}_j;d}$, $\mc{E}_{\widetilde{P}_j;d}^*$ only have a finite number of singular values greater than 0, \cref{as:solutions-exist} implies \cref{as:beta-source} for any choice of $\beta_h, \beta_q \in [0, \infty)$. $\widetilde{P}_j$ is mutually absolutely continuous with $P_j$ by \cref{lemma:constant-riesz-representers}, and $P_j$ is mutually absolutely continuous with $P_0$, so \cref{as:mutual-absolute-continuity} holds. Finally, \cref{as:dominating-measure-support} also holds by the construction of $\widetilde{P}_j$ in \cref{lemma:constant-riesz-representers}.

    Hence, we have constructed a sequence of distributions $\{\widetilde{P}_j\}_{j=1}^{\infty}$ such that $\lim_{j\rightarrow \infty}H^2(\widetilde{P}_j, P_0) = 0$ and where $\widetilde{P}_j$ belongs to the statistical model $\mc{P}$ of \eqref{eq:model-definition} for each $j$. We conclude by demonstrating that $|\varphi(\widetilde{P}_j) - \varphi(P_0)| \rightarrow \infty$. It suffices to show that $|\varphi(\widetilde{P}_j)| \rightarrow \infty$.

    Let $h_{\widetilde{P}_j}^* \coloneqq \sum_{d \in \textsf{Supp}[P_{D;0}]}I(D=d)h_{\widetilde{P}_j}^*(\cdot, d)$ where for each $d$, $h_{\widetilde{P}_j}^*(\cdot, d)$ is the  minimum norm solution to the inverse problem
    \begin{align*}
        \mc{E}_{\widetilde{P}_j;d}h(\cdot, d) = r_{\widetilde{P}_j}(\cdot, d) = r_{{P}_j}(\cdot, d).
    \end{align*}
    Then,
    \begin{align*}
        (\mc{E}_{\widetilde{P}_j}h_{\widetilde{P}_j}^*)(Z, D) =& E_{\widetilde{P}_j}[h_{\widetilde{P}_j}^*(X, D)|Z, D]\\
        =& \sum_{d \in \textsf{Supp}[P_{D;0}]}E_{\widetilde{P}_j}[h_{\widetilde{P}_j}^*(X, d)|Z, D=d]I(D=d)\\
        =& \sum_{d \in \textsf{Supp}[P_{D;0}]}(\mc{E}_{\widetilde{P}_j;d}h_{\widetilde{P}_j}^*(\cdot, d))(Z)I(D=d)\\
        =& \sum_{d \in \textsf{Supp}[P_{D;0}]}r_{P_j}(\cdot, d)I(D=d)\\
        =& r_{P_j},
    \end{align*}
    i.e. $h_{\widetilde{P}_j}^*$ solves the primal inverse problem. Then, we have that,
    \begin{align*}
        |\varphi(\widetilde{P}_j)| =& \left|E_{\widetilde{P}_j}[\alpha_{\widetilde{P}_j}(X, D)h_{\widetilde{P}_j}^*(X, D)]\right|\\
        =& \left|E_{{P}_j}[\alpha_{{P}_j}(X, D)h_{\widetilde{P}_j}^*(X, D)]\right|\\
        =& \left|\sum_{d \in \textsf{Supp}[D;P_j]}P_j(D=d)E_{{P}_j}[\alpha_{{P}_j}(X, d)h_{\widetilde{P}_j}^*(X, d)|D=d]\right|\\
        =& \Bigg|P_j(D=\widetilde{d})\sum_{i\leq l_j+2:\widetilde{\sigma}_{i;\widetilde{d};j}>0}\widetilde{\sigma}_{i;\widetilde{d};j}^{-1}\langle \widetilde{u}_{i; \widetilde{d};j}, \alpha_{P_j}(\cdot, \widetilde{d})\rangle_{L^2(P_{X|\widetilde{d};j})}\langle \widetilde{v}_{i; \widetilde{d};j}, r_{P_j}(\cdot, \widetilde{d})\rangle_{L^2(P_{Z|\widetilde{d};j})}\\
        &+ \sum_{d \in \textsf{Supp}[D;P_j]\setminus \{\widetilde{d}\}}P_j(D=d)\sum_{i\leq l_j+1: \widetilde{\sigma}_{i;{d};j}>0}\Big\{\widetilde{\sigma}_{i;{d};j}^{-1}\langle \widetilde{u}_{i; {d};j}, \alpha_{P_j}(\cdot, {d})\rangle_{L^2(P_{X|d;j})} \\
        &\qquad\qquad \langle \widetilde{v}_{i; {d};j}, r_{P_j}(\cdot, {d})\rangle_{L^2(P_{Z|d;j})}\Big\}\Bigg|\\
        \geq&P_j(D=\widetilde{d})\widetilde{\sigma}_{l_j+1;\widetilde{d};j}^{-1}\bigg|\left\langle \alpha_{{P}_j}(\cdot, \widetilde{d}), \widetilde{u}_{l_j+1;\widetilde{d};j} \right\rangle_{L^2(P_{X|\widetilde{d};j})} \\
        &\qquad\qquad\left\langle r_{{P}_j}(\cdot,\widetilde{d}), \widetilde{v}_{l_j+1;\widetilde{d};j} \right\rangle_{L^2(P_{Z|\widetilde{d};j})}\bigg| - |\kappa_j| \\
        > & j.
    \end{align*}
    The second equality is justified because the $(W, X, D)$ marginals of $P_j$ and $\widetilde{P}_j$ are the same. The first inequality is the reverse triangle inequality, and the second follows from \eqref{eq:sigma_discontinuity-def}.
    Hence, $\lim_{j\rightarrow\infty}|\varphi(\widetilde{P}_j) - \varphi(P_0)| = \infty$.
\end{proof}

\subsection{Proofs for \cref{sup-sec:smoothness-implications}}

\begin{proof}[Proof of \cref{prop:two-point-ci-bound}]
Let
\begin{equation*}
    A_0:=\{\phi(P_0)\in C_n\},
    \qquad
    A_1:=\{\phi(P_1)\in C_n\}.
\end{equation*}
By coverage at \(P_0\), \(P_0^n(A_0)\ge 1-\alpha\).  By coverage at
\(P_1\), \(P_1^n(A_1)\ge 1-\alpha\).  Transferring this latter probability from
\(P_1^n\) to \(P_0^n\) gives
\begin{align*}
    P_0^n(A_1)
    &\ge
    P_1^n(A_1)-d_{\mathrm{TV}}(P_0^n,P_1^n) \\
    &\ge
    1-\alpha-d_{\mathrm{TV}}(P_0^n,P_1^n).
\end{align*}
Therefore, by the union bound,
\begin{align*}
    P_0^n(A_0\cap A_1)
    &\ge
    1-P_0^n(A_0^c)-P_0^n(A_1^c) \\
    &\ge
    1-\alpha-\{\alpha+d_{\mathrm{TV}}(P_0^n,P_1^n)\} \\
    &=
    1-2\alpha-d_{\mathrm{TV}}(P_0^n,P_1^n).
\end{align*}
On the event \(A_0\cap A_1\), the interval \(C_n\) contains both
\(\phi(P_0)\) and \(\phi(P_1)\).  Hence \(|C_n|\ge \Delta\) on this event.
This proves \eqref{eq:two-point-ci-tail}.  Multiplying the indicator of the event \(|C_n|\ge \Delta\) by \(\Delta\) and taking
expectations proves \eqref{eq:two-point-ci-expected-length}.
\end{proof}

\begin{proof}[Proof of \cref{prop:local-modulus-ci-bound}]
Suppose first that \(\omega_{P_0}(r)<\infty\), and let \(\eta>0\).  By the definition of \(\omega_{P_0}(r)\), there exists
\(P_1\in\mathcal M\) with \(H(P_1,P_0)<r\) such that
\begin{equation*}
    |\phi(P_1)-\phi(P_0)|
    \ge
    \omega_{P_0}(r)-\eta.
\end{equation*}  Since \(C_n\)
has coverage at \(P_0\) and \(P_1\), \cref{prop:two-point-ci-bound} gives
\begin{equation*}
    E_{P_0^n}\{|C_n|\}
    \ge
    \{\omega_{P_0}(r)-\eta\}
    \left[
        1-2\alpha-d_{\mathrm{TV}}(P_0^n,P_1^n)
    \right]_+ .
\end{equation*}
Invoking the inequality \eqref{eq:tv-product-hellinger-bound},
the condition \(H(P_1,P_0)<r\) then yields
\begin{equation*}
    d_{\mathrm{TV}}(P_0^n,P_1^n)
    \le
    \left\{1-(1-r^2)^{2n}\right\}^{1/2}.
\end{equation*}
Letting \(\eta\downarrow0\) proves the claim. If instead \(\omega_{P_0}(r)=\infty\), then for every \(M>0\) there exists \(P_1\in\mathcal M\) with \(H(P_1,P_0)<r\) and \(|\phi(P_1)-\phi(P_0)|\ge M\), and the same argument gives
\begin{equation*}
    E_{P_0^n}\{|C_n|\}
    \ge
    M\left[
        1-2\alpha-
        \left\{1-(1-r^2)^{2n}\right\}^{1/2}
    \right]_+ .
\end{equation*}
Letting \(M\rightarrow\infty\) proves the claim, with the convention \(\infty\cdot 0 = 0\) when the bracket vanishes.
\end{proof}

\subsection{Proofs for \cref{subsec:tangent-spaces}}

\begin{proof}[Proof of \cref{lemma:discrete-tangent-space}]

    { We first provide a sketch of the proof. We use that every $K_X \times K_X$ matrix with degenerate rank is arbitrarily close to an invertible matrix. Because $X, Z$ are finite-discrete with $K_X = K_Z$, the operator $\mc{E}_{P;d}$ has an equivalent representation as a $K_X \times K_X$ matrix. Then, given any submodel $\{P_t:t \in [0, \delta)\}$, we may perturb $P_t$ so that the perturbed matrix has full rank, but where this perturbation is $o(t)$ in Hellinger distance, and so does not contribute to the score of the path. But, $\alpha_{P_t},r_{P_t}$ are certainly in the range of an invertible matrix, satisfying the conditions of the model.

    }

    Let $s \in L^\infty(P_0)\cap L^2_0(P_0)$ and for $t < \|s\|^{-1}_{L^\infty}$ define $\widetilde{p}_t(o) \coloneqq (1 + ts(o))p_0(o)$ with $\widetilde{P}_t$ the corresponding distribution for $O$. Then, $\widetilde{p}_t$ is a density for each $t<\|s\|^{-1}_{L^\infty}$ and $\{\widetilde{P}_t: t\in [0, \|s\|^{-1}_{L^\infty})\}$ is QMD at $t=0$ with score $s$, and all such $\widetilde{P}_t$ are mutually absolutely continuous with $P_0$. Now, let
    \begin{align*}
        L_{gap} \coloneqq \min( \essinf {p}_0  - L_{bounded}^{-1},L_{bounded} -  \esssup p_0)
    \end{align*}
    which is strictly positive by \cref{as:dominating-measure-support}. In what follows, we restrict to $t < t_s$, where $t_s < \min(1, \|s\|^{-1}_{L^\infty})$ is such that for all $t < t_s$,
    \begin{align*}
        L_{gap}/2 < \min( \essinf \widetilde{p}_t  - L_{bounded}^{-1},L_{bounded} -  \esssup \widetilde{p}_t).
    \end{align*}
    Such choice of $t$ is possible because $\widetilde{p}_t \rightarrow p_0$ in $L^\infty(P_0)$. This is without loss of generality, because quadratic mean differentiability at $t=0$, and the score of a submodel, depend on that submodel only through its behavior at $t$ arbitrarily close to $0$.

    Now, let $(\widetilde{\sigma}_{i;d;t}, \widetilde{u}_{i;d;t}, \widetilde{v}_{i;d;t})_{i=1}^{K_X}$ denote a singular system of $\mc{E}_{\widetilde{P}_t;d}$, with $\widetilde{\sigma}_{1;d;t} \geq \widetilde{\sigma}_{2;d;t} \geq \dots \geq \widetilde{\sigma}_{K_X;d;t} \geq 0$. In this proof we use the full singular system, i.e. allow for singular values equal to 0, such that $(\widetilde{u}_{i;d;t})_{i=1}^{K_X}, (\widetilde{v}_{i;d;t})_{i=1}^{K_X}$ are orthonormal basis of $L^2(\widetilde{P}_{X|d;t}), L^2(\widetilde{P}_{Z|d;t})$, respectively. Because $\mc{E}_{\widetilde{P}_t;d}$ is a contraction, $\widetilde{\sigma}_{1;d;t} \leq 1$. Furthermore because $\mc{E}_{\widetilde{P}_t;d}^*\mc{E}_{\widetilde{P}_t;d}1 = 1, \mc{E}_{\widetilde{P}_t;d}\mc{E}_{\widetilde{P}_t;d}^*1 = 1$, we may order the singular system so that $\widetilde{u}_{1;d;t} =1$ and $\widetilde{v}_{1;d;t} =1$.

    Define
    \begin{align*}
        \widetilde{k}_{d;t}(x, z) \coloneqq \sum_{i=1}^{K_X} \widetilde{\sigma}_{i;d;t}\widetilde{u}_{i;d;t}(x) \widetilde{v}_{i;d;t}(z)
    \end{align*}
    and let
    \begin{align*}
        C_t &\coloneqq \min_{d \in \textsf{Supp}[P_{D;0}]}\min_{x , z \in \{1, \dots, K_X\}} \widetilde{k}_{d;t}(x, z).\\
        C'_t &\coloneqq  \max_{d \in \textsf{Supp}[P_{D;0}]}\max_{x, z \in \{1, \dots, K_X\}}\left|\sum_{i=2}^{K_X} \widetilde{u}_{i;d;t}(x) \widetilde{v}_{i;d;t}(z)\right|
    \end{align*}
    and for $1 < i \leq K_X$, let
    \begin{align*}
        {\sigma}_{i;d;t} \coloneqq \widetilde{\sigma}_{i;d;t} + t^3\{C'_t\}^{-1}\min(C_t, L_{bounded}^{-6}L_{gap}/2).
    \end{align*}
    We note that $C_t \in (0, \infty)$: because $L_{bounded}^{-1} < \widetilde{p}_t < L_{bounded}$, the density $\widetilde{p}_{X, Z|d;t}$ and the marginal densities $\widetilde{p}_{X|d;t}, \widetilde{p}_{Z|d;t}$ are each bounded above and below by positive constants, and $\widetilde{k}_{d;t} = \widetilde{p}_{X, Z|d;t}/(\widetilde{p}_{X|d;t}\widetilde{p}_{Z|d;t})$. We further note that $C'_{t} \in (0, \infty)$ whenever $K_X \geq 2$: the $K_X\times K_X$ matrix with $(x, z)$ entry $\sum_{i=2}^{K_X} \widetilde{u}_{i;d;t}(x) \widetilde{v}_{i;d;t}(z)$ has rank $K_X - 1 \geq 1$, because $(\widetilde{u}_{i;d;t})_{i=2}^{K_X}$ and $(\widetilde{v}_{i;d;t})_{i=2}^{K_X}$ are each linearly independent, and so is not the zero matrix. When $K_X = 1$, $\mc{E}_{P;d}$ is invertible for every $P$, so that we may take $P_t \coloneqq \widetilde{P}_t$ in what follows and no perturbation is required; we therefore assume $K_X \geq 2$.

    Let $\sigma_{1;d;t} \coloneqq 1$, and ${u}_{i;d;t} \coloneqq \widetilde{u}_{i;d;t}, {v}_{i;d;t} \coloneqq \widetilde{v}_{i;d;t}$. Define
    \begin{align*}
        {k}_{d;t}(x, z) \coloneqq \sum_{i=1}^{K_X} {\sigma}_{i;d;t}{u}_{i;d;t}(x) {v}_{i;d;t}(z).
    \end{align*}
    It follows that
    \begin{align*}
        \max_{x , z \in \{1, \dots, K_X\}} |\widetilde{k}_{d;t}(x, z) - k_{d;t}(x, z)| =&  \max_{x , z \in \{1, \dots, K_X\}} \left|\sum_{i=2}^{K_X} \{\widetilde{\sigma}_{i;d;t} - {\sigma}_{i;d;t}\}\widetilde{u}_{i;d;t}(x) \widetilde{v}_{i;d;t}(z)\right|\\
        =&t^3\{C'_t\}^{-1} \min(C_t, L_{bounded}^{-6}L_{gap}/2) \max_{x , z \in \{1, \dots, K_X\}} \left|\sum_{i=2}^{K_X} \widetilde{u}_{i;d;t}(x) \widetilde{v}_{i;d;t}(z)\right|\\
        \leq& t^3\min(C_t, L_{bounded}^{-6}L_{gap}/2).
    \end{align*}
    Then, because $\widetilde{k}_{d;t}(x, z) \geq C_t$ for each $x, z$ by the definition of $C_t$, and because $\min(C_t, L_{bounded}^{-6}L_{gap}/2) \leq C_t$, the above display gives $k_{d;t}(x, z) \geq (1 - t^3)C_t$. Hence, for $t < 1$, $k_{d;t}(x, z) > 0$ for each $x, z \in \{1, \dots, K_X\}$.
    We further claim that ${p}_{X, Z|d;t}(x, z) \coloneqq k_{d;t}(x, z)\widetilde{p}_{X|d;t}(x)\widetilde{p}_{Z|d;t}(z)$ is a joint density for $(X, Z)$. As it is a non-negative function, it suffices to show it integrates to 1.
    \begin{align*}
        \int\int  {p}_{X, Z|d;t}(x, z)d\lambda_X(x)d\lambda_Z(z) =& \int \int k_{d;t}(x, z)\widetilde{p}_{X|d;t}(x)\widetilde{p}_{Z|d;t}(z)d\lambda_X(x)d\lambda_Z(z)\\
        =& \int \int \sum_{i=1}^{K_X} {\sigma}_{i;d;t}{u}_{i;d;t}(x) {v}_{i;d;t}(z)\widetilde{p}_{X|d;t}(x)\widetilde{p}_{Z|d;t}(z)d\lambda_X(x)d\lambda_Z(z)\\
        =& \sum_{i=1}^{K_X}{\sigma}_{i;d;t} \langle u_{i;d;t}, 1\rangle_{L^2(\widetilde{P}_{X|d;t})}\langle v_{i;d;t}, 1\rangle_{L^2(\widetilde{P}_{Z|d;t})}\\
        =& 1
    \end{align*}
    where the final equality follows from $\sigma_{1;d;t}= 1$, $u_{1;d;t} = 1$, $v_{1;d;t} =1$. Furthermore, we claim ${p}_{X|d;t}(x) \coloneqq \int {p}_{X, Z|d;t}(x, z)d\lambda_Z(z) = \widetilde{p}_{X|d;t}(x)$. To see this,
    \begin{align*}
        \int  {p}_{X, Z|d;t}(x, z)d\lambda_Z(z) =& \int k_{d;t}(x, z)\widetilde{p}_{X|d;t}(x)\widetilde{p}_{Z|d;t}(z)d\lambda_Z(z)\\
        =& \int \sum_{i=1}^{K_X} {\sigma}_{i;d;t}{u}_{i;d;t}(x) {v}_{i;d;t}(z)\widetilde{p}_{X|d;t}(x)\widetilde{p}_{Z|d;t}(z)d\lambda_Z(z)\\
        =& \widetilde{p}_{X|d;t}(x)\sum_{i=1}^{K_X}{\sigma}_{i;d;t}  u_{i;d;t}(x)\langle v_{i;d;t}, 1\rangle_{L^2(\widetilde{P}_{Z|d;t})}\\
        =& \widetilde{p}_{X|d;t}(x).
    \end{align*}
    The same argument shows ${p}_{Z|d;t}(z) \coloneqq \int {p}_{X, Z|d;t}(x, z)d\lambda_X(x) = \widetilde{p}_{Z|d;t}(z)$. We let $p_t(d) = \widetilde{p}_t(d)$ for each $d$.

    Now, because $\widetilde{p}_{X|d;t}, \widetilde{p}_{Z|d;t} < L_{bounded}^2$ by \cref{remark:marginal-conditional-density-bounds}, it holds for each $t < t_s$ that
    \begin{align*}
        &\max_{x, z\in \{1, \dots, K_X\}}|\widetilde{p}_{X, Z|d;t}(x, z) - p_{X, Z|d;t}(x, z)|  \\
        \leq&\max_{x, z\in \{1, \dots, K_X\}}|\widetilde{k}_{d;t}(x, z) - k_{d;t}(x, z)|\max_{x\in \{1, \dots, K_X\}}|\widetilde{p}_{X|d;t}(x)|\max_{z\in \{1, \dots, K_X\}}|\widetilde{p}_{Z|d;t}(z)|\\
        \leq& t^3 \min(C_t, L_{bounded}^{-6}L_{gap}/2) L_{bounded}^4\\
        \leq& t^3L_{bounded}^{-2}L_{gap}/2.
    \end{align*}
    We now verify the hypotheses of \cref{lemma:constant-riesz-representers} with $P \coloneqq \widetilde{P}_t$ and $p'_{X, Z, D} \coloneqq p_{X, Z, D;t}$ in that Lemma. \cref{as:dominating-measure-support} holds at $\widetilde{P}_t$ because $L_{bounded}^{-1} < \essinf \widetilde{p}_t\leq \esssup \widetilde{p}_t < L_{bounded}$ and because $\widetilde{P}_t$ and $P_0$ are supported on the same bounded rectangle. \cref{as:continuous-linear-maps} holds at $\widetilde{P}_t$ because, $X, Z, D$ being finite discrete, the spaces $L^2(\widetilde{P}_{X,D;t})$ and $L^2(\widetilde{P}_{Z,D;t})$ are finite dimensional, so that every linear functional on them is bounded, and because $\widetilde{p}_t \leq (1 + t\|s\|_{L^\infty(P_0)})p_0$ ensures that $\widetilde{m}(W, X, D;h)$ and ${m}(W, Z, D;q)$ are $\widetilde{P}_t$-integrable for every $h \in L^2(\widetilde{P}_{X,D;t})$ and $q \in L^2(\widetilde{P}_{Z,D;t})$, as they are $P_0$-integrable. We showed above that ${p}_{X|d;t} = \widetilde{p}_{X|d;t}$, ${p}_{Z|d;t} = \widetilde{p}_{Z|d;t}$ and ${p}_t(d) = \widetilde{p}_t(d)$, so that the $(X, D)$ and $(Z, D)$ marginals required by that Lemma agree. Finally, because $\widetilde{p}_t(d) < L_{bounded}$ by \cref{remark:marginal-conditional-density-bounds}, the above display gives
    \begin{align*}
        \|\widetilde{p}_{X, Z, D;t} - p_{X, Z, D;t}\|_{L^\infty(\lambda_X\otimes\lambda_Z\otimes\lambda_D)} < t^3L_{bounded}^{-1}L_{gap}/2 < L_{bounded}^{-1}L_{gap}/2,
    \end{align*}
    where the final inequality holds because $t < 1$. With $P \coloneqq \widetilde{P}_t$, and because $\lambda_X(\Omega_X) = \lambda_Z(\Omega_Z) = \lambda_D(\Omega_D) = 1$ by \cref{remark:marginal-conditional-density-bounds}, the threshold in \cref{lemma:constant-riesz-representers} is $L_{bounded}^{-1}\min(\essinf \widetilde{p}_t - L_{bounded}^{-1}, L_{bounded} - \esssup \widetilde{p}_t)$, which exceeds $L_{bounded}^{-1}L_{gap}/2$ because $t < t_s$. Hence the bound required by that Lemma holds.

    Then, for $t < t_s$, by \cref{lemma:constant-riesz-representers} there exists a law $P_t$ for $(W, X, Z, D)$ satisfying \cref{as:continuous-linear-maps} and \cref{as:dominating-measure-support} with
    \begin{enumerate}
        \item $\alpha_{P_t} = \alpha_{\widetilde{P}_t}$, $r_{P_t} = r_{\widetilde{P}_t}$ and
        \item $H^2(P_t, \widetilde{P}_t) \leq \frac{1}{2}\|\widetilde{p}_{X, Z, D;t} - p_{X, Z, D;t}\|_{L^1(\lambda_X\otimes\lambda_Z\otimes\lambda_D)}$.
    \end{enumerate}
    Hence, because $\lambda_X(\Omega_X) = \lambda_Z(\Omega_Z) = \lambda_D(\Omega_D) = 1$,
    \begin{align*}
        H^2(P_t, \widetilde{P}_t) \leq& \frac{1}{2}\|\widetilde{p}_{X, Z, D;t} - p_{X, Z, D;t}\|_{L^1(\lambda_X\otimes\lambda_Z\otimes\lambda_D)}\\
        \leq& \frac{1}{2}\|\widetilde{p}_{X, Z, D;t} - p_{X, Z, D;t}\|_{L^\infty(\lambda_X\otimes\lambda_Z\otimes\lambda_D)}\\
        \leq& \frac{t^3}{2} L_{bounded}^{-1}L_{gap}/2.
    \end{align*}

    We now claim $P_t \in \mc{P}$ for $t < t_s$. \cref{as:continuous-linear-maps} and \cref{as:dominating-measure-support} hold by \cref{lemma:constant-riesz-representers}. \cref{as:solutions-exist} holds because $\mc{E}_{P_t}$ has singular system $(\sigma_{i;d;t}, u_{i;d;t}, v_{i;d;t})$, which has all positive singular values, and so $\mathscr{R}(\mc{E}_{P_t}) = L^2(P_{Z,D;t})$, $\mathscr{R}(\mc{E}^*_{P_t}) = L^2(P_{X,D;t})$. Because $\mc{E}_{P_t}$ is a finite dimensional operator, \cref{as:beta-source} holds for any choice of $(\beta_h, \beta_q)$. \cref{as:mutual-absolute-continuity} holds because $P_t$ is mutually absolutely continuous with $\widetilde{P}_t$, and because $\widetilde{P}_t$ is mutually absolutely continuous with $P_0$.

    We conclude by arguing that $\{P_t:t \in [0, t_s)\}$ is QMD at $t=0$ with score $s$. To see this,
    \begin{align*}
        &\int \left\{\frac{\sqrt{p}_t(o) - \sqrt{p}_0(o)}{t} - \frac{1}{2}s(o)\sqrt{p}_0(o)\right\}^2d\lambda(o)\\
        =&\int \left\{\frac{\sqrt{p}_t(o) - \sqrt{\widetilde{p}}_t(o)   +  \sqrt{\widetilde{p}}_t(o)- \sqrt{p}_0(o)}{t} - \frac{1}{2}s(o)\sqrt{p}_0(o)\right\}^2d\lambda(o)\\
        \leq&4H^2(P_t, \widetilde{P}_t) /t^2  + o(1)\\
        \leq & t L_{bounded}^{-1}L_{gap} + o(1) \\
        =& o(1).
    \end{align*}
    where first inequality follows because $\{\widetilde{P}_t: t\in [0, \|s\|^{-1}_{L^\infty({P}_0)})\}$ is QMD at $t=0$ with score $s$. The final equality holds because $L_{bounded}^{-1}L_{gap}$ does not depend on $t$.

    Because $s$ was arbitrary beyond belonging to a dense subset of $L^2_0(P_0)$, we have that $\mc{T}(P_0, \mc{P}) = L^2_0(P_0)$, completing the proof.
\end{proof}

\begin{proof}[Proof of \cref{thrm:continuous-tangent-space}]

    { We begin by giving an informal sketch of the proof. The proof is similar to that of \cref{thrm:continuous-discontinuity}. As before, the main technical difficulty arises from ensuring that our proposed path $P_t$ is a valid joint distribution for $O$. We begin by taking a linear tilt QMD submodel $\{\widetilde{P}_t:t \in [0, \delta)\}$ with uniformly bounded score, and perturb this submodel so it belongs to $\mc{P}$ but such that the perturbation does not change the score. We construct this perturbation over several steps. First, by \cref{lemma:dense-continuous-densities}, we have for every $t$ that the conditional density of $(W_C, X, Z)$ given the discrete variables $(W_D, D)$ under $\widetilde{P}_t$ is arbitrarily close in $L^q(\lambda)$ to a conditional density that is infinitely differentiable on the compact support; recombining with the law of $(W_D, D)$ gives a law $P_t^\dagger$ satisfying the density bounds of \cref{as:dominating-measure-support}. In light of \cref{lemma:uniform-convergence-of-kernel}, we may truncate the singular system of $\mc{E}_{P_t^\dagger}$ sufficiently far into the tail to construct a new finite dimensional singular system that corresponds to the conditional expectation operator of a valid joint density for $X, Z$ given $D=d$. We may further perturb this finite dimensional singular system so that $\alpha_{\widetilde{P}^\dagger_t}, r_{\widetilde{P}_t^\dagger}$ are in the range of the operator and its adjoint. The $\beta$-source conditions are then satisfied for this operator, as it is finite dimensional. \cref{lemma:constant-riesz-representers} then gives a full joint distribution $P_t$ for $O$ belonging to the model $\mc{P}$. We may make these perturbations sufficiently small so that for each $t$, $H(P_t, \widetilde{P}_t) = o(t)$. Consequently, $\{P_t:t \in [0, \delta)\}$ has the same score as $\{\widetilde{P}_t:t \in [0, \delta)\}$. As the score of $\{\widetilde{P}_t:t \in [0, \delta)\}$ was arbitrary beyond belonging to a dense subspace of $L^2_0(P_0)$, this proves the result. }

    Let $s \in L^\infty(P_0) \cap L^2_0(P_0)$, and define a new density $\widetilde{p}_t(o) \coloneqq (1 + ts)p_0(o)$. Let $\widetilde{P}_t$ be the corresponding distribution. It follows that for $0 < t < \|s\|_{L^\infty(P_0)}^{-1}$, $\widetilde{p}_t$ is a density and $\widetilde{P}_t$ is mutually absolutely continuous with $P_0$.
    Let
    \begin{align*}
        L_{gap} \coloneqq \min(\essinf p_0 - L_{bounded}^{-1} , L_{bounded} - \esssup p_0),
    \end{align*}
    which is strictly positive by \cref{as:dominating-measure-support}.
    Noting that $\widetilde{p}_t \rightarrow p_0$ in $L^\infty(P_0)$, we may take $t_s < \min(1, \|s\|_{L^\infty(P_0)}^{-1})$ so that for all $t < t_s$,
    \begin{align*}
        L_{gap}/2 < \min(\essinf \widetilde{p}_t - L_{bounded}^{-1} , L_{bounded} - \esssup \widetilde{p}_t).
    \end{align*}

    Now, because \cref{as:dominating-measure-support} holds at $P_0$, and $r_{P_0} \in \mathscr{R}(\mc{E}_{P_0}), \alpha_{P_0} \in \mathscr{R}(\mc{E}_{P_0}^*)$ by \cref{as:solutions-exist}, we have by \cref{lemma:in-range-uniformly-bounded} that
    \begin{align*}
        \|\alpha_{P_0}\|_{L^\infty(P_{X,D;0})} < \infty,  &&  \|r_{P_0}\|_{L^\infty(P_{Z,D;0})} < \infty.
    \end{align*}
    By \cref{lemma:l1-bounded-functional}, \cref{as:continuous-linear-maps} holds at $\widetilde{P}_t$ and $\alpha_{\widetilde{P}_t} \in L^\infty(\widetilde{P}_{X,D;t})$, $ r_{\widetilde{P}_t} \in L^\infty(\widetilde{P}_{Z,D;t})$.

    In what follows, we write $W = (W_C, W_D)$, where $W_C$ collects the components of $W$ whose reference measure in \cref{as:dominating-measure-support} is the Lebesgue measure and $W_D$ collects those whose reference measure is the counting measure; either may be empty. Correspondingly, $\lambda_W = \lambda_{W_C}\otimes\lambda_{W_D}$, where $\lambda_{W_C}$, $\lambda_X$, $\lambda_Z$ are Lebesgue product measures and $\lambda_{W_D}$, $\lambda_D$ are counting measures. All densities are taken with respect to $\lambda$, or the appropriate product measure for functions of $W, X, Z, D$ and combinations thereof. We normalize these reference measures as in \cref{remark:marginal-conditional-density-bounds}, so that $\lambda_W(\textsf{Supp}[W])= \lambda_X(\textsf{Supp}[X])= \lambda_Z(\textsf{Supp}[Z])= \lambda_D(\textsf{Supp}[D])=1$.

    Now, for each $t$, we construct a sequence  $\{P_{j,t}^\dagger\}_{j=1}^\infty$ of distributions whose densities in $(W_C, X, Z)$ are infinitely differentiable and whose full joint densities fall within the density bounds required by \cref{as:dominating-measure-support}. Because $W_D, D$ are finite discrete, the set $\mc{S} \coloneqq \textsf{Supp}[(W_D, D);P_0]$ is finite. Let $\Omega_C$ denote the open bounded rectangle supporting $(W_C, X, Z)$ under $P_0$. Because the support of $P_0$ is a rectangle, $\Omega_C$ also supports the conditional law of $(W_C, X, Z)$ given $(W_D, D) = (w_D, d)$ for each $(w_D, d) \in \mc{S}$. Fix $t < t_s$. By the choice of $t_s$, we may fix constants $c_{lo}, c_{hi}$, depending on $t$, with
    \begin{align*}
        L_{bounded}^{-1} < c_{lo} < \essinf \widetilde{p}_t \leq \esssup \widetilde{p}_t < c_{hi} < L_{bounded}.
    \end{align*}
    For $(w_D, d) \in \mc{S}$, let $q_t(w_D, d) \coloneqq \widetilde{p}_{W_D, D;t}(w_D, d) > 0$, and let $\widetilde{p}_{W_C, X, Z|w_D, d;t} \coloneqq \widetilde{p}_t(\cdot, w_D, \cdot, \cdot, d)/q_t(w_D, d)$ denote the conditional density of $(W_C, X, Z)$ given $(W_D, D) = (w_D, d)$ with respect to $\lambda_{W_C}\otimes\lambda_X\otimes\lambda_Z$. Then
    \begin{align*}
        \frac{c_{lo}}{q_t(w_D, d)} < \widetilde{p}_{W_C, X, Z|w_D, d;t} < \frac{c_{hi}}{q_t(w_D, d)} \qquad \text{a.e. on }\Omega_C.
    \end{align*}
    We apply \cref{lemma:dense-continuous-densities} with $Y = (W_C, X, Z)$, $\Omega = \Omega_C$, $p = \widetilde{p}_{W_C, X, Z|w_D, d;t}$, $a = c_{lo}/q_t(w_D, d)$, and $b = c_{hi}/q_t(w_D, d)$, separately for each of the finitely many $(w_D, d) \in \mc{S}$. Because $\gamma(\varepsilon) \rightarrow 1$ as $\varepsilon \downarrow 0$ in that Lemma and $\mc{S}$ is finite, for each $j$ we may choose $\varepsilon = \varepsilon_j$ small enough that the resulting densities $p^\dagger_{W_C, X, Z|w_D, d;j, t}$ supported on $\Omega_C$, with corresponding constants $\gamma_{j, t}(w_D, d) > 0$, satisfy for every $(w_D, d) \in \mc{S}$:
    \begin{enumerate}
        \item $p^\dagger_{W_C, X, Z|w_D, d;j, t} \in C^\infty(\overline{\Omega}_C)$;
        \item $c_{lo}/\{q_t(w_D, d)\gamma_{j, t}(w_D, d)\} < p^\dagger_{W_C, X, Z|w_D, d;j, t} < c_{hi}/\{q_t(w_D, d)\gamma_{j, t}(w_D, d)\}$ a.e. on $\Omega_C$, where $L_{bounded}^{-1} < c_{lo}/\gamma_{j, t}(w_D, d)$ and $c_{hi}/\gamma_{j, t}(w_D, d) < L_{bounded}$;
        \item $\|p^\dagger_{W_C, X, Z|w_D, d;j, t} - \widetilde{p}_{W_C, X, Z|w_D, d;t}\|_{L^1(\lambda_{W_C}\otimes\lambda_X\otimes\lambda_Z)} < 1/j$.
    \end{enumerate}
    We then let $P^\dagger_{j, t}$ be the law with density
    \begin{align*}
        p_{j, t}^\dagger(w_C, w_D, x, z, d) \coloneqq q_t(w_D, d)\,p^\dagger_{W_C, X, Z|w_D, d;j, t}(w_C, x, z)
    \end{align*}
    for $(w_D, d) \in \mc{S}$, and $p_{j, t}^\dagger \coloneqq 0$ otherwise. Then $p_{j, t}^\dagger$ is a probability density with $p_{W_D, D;j, t}^\dagger = \widetilde{p}_{W_D, D;t}$, and in particular $P^\dagger_{j, t}(D=d) = \widetilde{P}_t(D=d)$ for each $d$. Multiplying the second item through by $q_t(w_D, d)$, and taking the minimum and maximum over the finitely many $(w_D, d) \in \mc{S}$, gives
    \begin{align*}
        L_{bounded}^{-1} < \essinf p_{j, t}^\dagger \leq \esssup p_{j, t}^\dagger < L_{bounded},
    \end{align*}
    where the essential infimum and supremum are taken over the support of $P_0$, which is also the support of $P_{j, t}^\dagger$. Hence $P_{j, t}^\dagger$ and $P_0$ are mutually absolutely continuous, and $P_{j, t}^\dagger$ satisfies \cref{as:dominating-measure-support}. By the third item,
    \begin{align*}
        \|p_{j, t}^\dagger - \widetilde{p}_t\|_{L^1(\lambda)} = \sum_{(w_D, d) \in \mc{S}} \widetilde{P}_t(W_D = w_D, D=d)\|p^\dagger_{W_C, X, Z|w_D, d;j, t} - \widetilde{p}_{W_C, X, Z|w_D, d;t}\|_{L^1(\lambda_{W_C}\otimes\lambda_X\otimes\lambda_Z)} < 1/j,
    \end{align*}
    so that $H^2(P_{j, t}^\dagger, \widetilde{P}_t) \leq \frac{1}{2}\|p_{j, t}^\dagger - \widetilde{p}_t\|_{L^1(\lambda)} < 1/j$. Because $p_{j, t}^\dagger$ and $\widetilde{p}_t$ lie between $L_{bounded}^{-1}$ and $L_{bounded}$ on their common support, $L^1(\lambda)$ convergence of $p_{j, t}^\dagger$ to $\widetilde{p}_t$ implies $L^q$ convergence for every $1 \leq q < \infty$, both of $p_{j, t}^\dagger$ and of each of its marginal and conditional densities. Finally, for each $d$, $p^\dagger_{X, Z|d;j, t} \in C^\infty(\overline{\Omega}_{X, Z})$: it is a finite positive linear combination of the functions $(x, z) \mapsto \int p^\dagger_{W_C, X, Z|w_D, d;j, t}(w_C, x, z)\,d\lambda_{W_C}(w_C)$ over $w_D$ with $(w_D, d) \in \mc{S}$ (with the integral omitted if $W_C$ is empty), and differentiation under the integral sign is justified because every derivative of $p^\dagger_{W_C, X, Z|w_D, d;j, t}$ is continuous on the compact set $\overline{\Omega}_C$, hence bounded.

    Moreover, by \cref{lemma:l1-bounded-functional}, $P_{j, t}^\dagger$ satisfies \cref{as:continuous-linear-maps} and $\alpha_{P_{j, t}^\dagger} \in L^\infty(P^\dagger_{X, D;j, t})$, $r_{P_{j, t}^\dagger} \in L^\infty(P^\dagger_{Z, D;j, t})$. Then, let $j^* = \lceil 1/t^3 \rceil +1$, and define $P_{t}^\dagger \coloneqq P_{j^*, t}^\dagger$
    so that
    \begin{align}
        \label{eq:tangent-space-dagger-tilde-distance}
        2H^2(P_{t}^\dagger, \widetilde{P}_t) \leq \|p_{t}^\dagger - \widetilde{p}_{t}\|_{L^1(\lambda)} < t^3.
    \end{align}
    Now, let
    \begin{align*}
        L_{gap;P_t^\dagger} \coloneqq \min(\essinf p_t^\dagger - L_{bounded}^{-1}, L_{bounded} - \esssup p_t^\dagger)
    \end{align*}
    which is strictly positive by construction.

    Let $(\sigma_{i;d;t}, u_{i;d;t}, v_{i;d;t})_{i=1}^\infty$ denote the singular system of $\mc{E}_{{P}^\dagger_t;{d}}$ for each $d, t$ with singular values decreasing in magnitude as $i \rightarrow \infty$ and repeated according to multiplicity. This singular system exists because $L_{bounded}^{-1} < {p}^\dagger_t(O) < L_{bounded}$ $P_t^\dagger$-a.e. implies that $\mc{E}_{P_t^\dagger}$ is compact (\cite{carrasco_chapter_2007}, discussion following Assumption A.1), and compact operators have countable singular systems with singular values accumulating at $0$ (\cite{kress_linear_2014-2}, Section 15.4). Because $\mc{E}_{P_t^\dagger;{d}}$ is a contraction, we have $\sigma_{1;d;t} \leq 1$. Furthermore, because $\mc{E}_{P_t^\dagger;{d}}^*\mc{E}_{P_t^\dagger;{d}}1 = 1$ and $\mc{E}_{P_t^\dagger;{d}}\mc{E}_{P_t^\dagger;{d}}^*1 = 1$, $\mc{E}_{P_t^\dagger;{d}}$ has $1$ as a singular value and $1$ as both left and right singular functions. We order $(\sigma_{i;d;t}, u_{i;d;t}, v_{i;d;t})_{i=1}^\infty$ such that $\sigma_{1;d;t}=1$, $u_{1;d;t}=1$, $v_{1;d;t}=1$ for each $t$.

    By the construction of $P^\dagger_{j^*, t}$ above, the conditional density $p^\dagger_{X, Z|d;t}$ defined below belongs to $C^{\infty}(\overline{\Omega}_{X, Z})$ for each $d$. Consequently, the functions
    \begin{align*}
        p^\dagger_{X, Z|{d};t}(X, Z|{d}) &\coloneqq \int p^\dagger_{W, X, Z|{d};t}(w, X, Z|{d})d\lambda_W(w), \\
        p^\dagger_{X|{d};t}(X|{d}) &\coloneqq \int p^\dagger_{X, Z|{d};t}(X, z|{d})d\lambda_Z(z), \\
        p^\dagger_{Z|{d};t}(Z|{d}) &\coloneqq \int p^\dagger_{X, Z|{d};t}(x, Z|{d})d\lambda_X(x)
    \end{align*}
    belong to $C^{\infty}(\overline{\Omega}_{X, Z})$, $C^{\infty}(\overline{\Omega}_{X})$, and $C^{\infty}(\overline{\Omega}_{Z})$, respectively.  Furthermore, because $p^\dagger_{X|{d};t}, p^\dagger_{Z|{d};t} > \{L_{bounded} \}^{-2} > 0$, we have that $k_{d;t}^\dagger \coloneqq p^\dagger_{X,Z|{d};t} / (p^\dagger_{X|{d};t}p^\dagger_{Z|{d};t}) \in C^{\infty}(\overline{\Omega}_{X, Z})$. Because $\overline{\Omega}_{X, Z}$ is compact, every derivative of $k^\dagger_{d;t}$ is uniformly bounded over $\overline{\Omega}_{X, Z}$. It then follows that $k^\dagger_{d;t} \in W_{s,2}(\Omega_{X, Z})$, the $(s,2)$-Sobolev space, for any $s \in \mathbb{N}$. In particular, $p^\dagger_{X, Z|d;t}(\cdot|{d}), k^\dagger_{d;t}$ satisfy the conditions of \cref{lemma:uniform-convergence-of-kernel} for each $t$. Consequently,
    \begin{align*}
        k_{d;t}^l \coloneqq \sum_{i=1}^l \sigma_{i;d;t}u_{i;d;t}v_{i;d;t}
    \end{align*}
    converges essentially uniformly to $k_{d;t}^\dagger$ as $l \rightarrow \infty$. Because $L_{bounded}^{-2} < p^\dagger_{X, Z|{d};t}, p^\dagger_{X|{d};t}, p^\dagger_{Z|{d};t} < L_{bounded}^2$ by \cref{remark:marginal-conditional-density-bounds}, we have that $L_{bounded}^{-6} < k^{\dagger}_{d;t}  < L_{bounded}^{6}$. Consequently, there exists some $l'$ such that for all $l > l'$, $k_{d;t}^l> 0$  $P_t^\dagger$-a.e for all $d$. For each $l > l'$, let $L_{t}^{l} \geq 1$ be such that $\{L_{t}^{l}\}^{-1} < k_{d;t}^l(X, Z) < L_t^l$ $P_t^\dagger$-a.e. for all $d$. $k_{d;t}^l$ is upper bounded because each $u_{i;d;t}$, $v_{i;d;t}$ are essentially uniformly bounded by \cref{lemma:in-range-uniformly-bounded}.

    We let
    \begin{align*}
        p^l_{X,Z|{d};t} \coloneqq k_{d;t}^lp^\dagger_{X|{d};t}p^\dagger_{Z|{d};t},
    \end{align*}
    which we claim is a density with respect to $\lambda_X\otimes\lambda_Z$. As $(x, z) \mapsto p^l_{X,Z|{d};t}(x, z)$ is a non-negative measurable function, it suffices to show that $\int\int p^l_{X,Z|{d};t}(x, z|{d})d\lambda(x)d\lambda(z) = 1$. Now,
    \begin{align*}
        &\int\int k_{d;t}^l(x, z)p^\dagger_{X|{d};t}(x)p^\dagger_{Z|{d};t}(z)d\lambda_X(x)d\lambda_Z(z)\\
        =&\int\int \sum_{i=1}^l\sigma_{i;d;t}u_{i;d;t}(x)v_{i;d;t}(z)p^\dagger_{X|{d};t}(x)p^\dagger_{Z|{d};t}(z)d\lambda_X(x)d\lambda_Z(z)\\
        =& \sum_{i=1}^l\sigma_{i;d;t}\langle u_{i;d;t}, 1\rangle_{L^2(P^\dagger_{X|d;t})}\langle v_{i;d;t}, 1\rangle_{L^2(P^\dagger_{Z|d;t})}\\
        =& 1
    \end{align*}
    where the final equality holds in light of $\sigma_{1;d;t}=1$, $u_{1;d;t}=1$, $v_{1;d;t}=1$. Additionally, $\{L_{t}^l\}^{-1}L_{bounded}^{-4} < p^l_{X,Z|{d};t} <\{L_{t}^l\}L_{bounded}^4 $ $P_t^\dagger$-a.e. by construction.

    Because $k_{d;t}^l$ converges uniformly to $k_{d;t}^\dagger$, we may choose $l $ large enough, such that
    \begin{align}
        \label{eq:p-l_t-t-infinity-bound}
        \|p_{X,Z|{d};t}^\dagger -p^l_{X,Z|{d};t}\|_{L^\infty(P^\dagger_{X,Z;t})} < \min\left(t^3, \frac{1}{2}L_{bounded}^{-4}L_{gap;P_t^\dagger}\right).
    \end{align}
    Let $l_t$ be any such $l$.

    For $i \leq l_t$ we let $\widetilde{u}_{i;d;t} \coloneqq {u}_{i;d;t}$, $\widetilde{v}_{i;d;t} \coloneqq {v}_{i;d;t}$, $\widetilde{\sigma}_{i;d;t} \coloneqq {\sigma}_{i;d;t}$. If
    \begin{align*}
        \alpha_{P_t^\dagger}(\cdot, {d}) - \sum_{i\leq l_t: \widetilde{\sigma}_{i;d;t}>0}\langle \widetilde{u}_{i;d;t}, \alpha_{P_t^\dagger}(\cdot, {d})\rangle_{L^2(P^\dagger_{X|d;t})}\widetilde{u}_{i;d;t} \not = 0,
    \end{align*}
    we let
    \begin{align}
        \label{eq:residual-projection-generic-d-tangent}
        \widetilde{u}_{l_t+1;d;t} \coloneqq& \frac{\alpha_{P_t^\dagger}(\cdot, {d}) - \sum_{i\leq l_t: \widetilde{\sigma}_{i;d;t}>0}\langle \widetilde{u}_{i;d;t}, \alpha_{P_t^\dagger}(\cdot, {d})\rangle_{L^2(P^\dagger_{X|d;t})}\widetilde{u}_{i;d;t}}{\|\alpha_{P^\dagger_t}(\cdot, {d}) - \sum_{i\leq l_t: \widetilde{\sigma}_{i;d;t}>0}^{l_t}\langle \widetilde{u}_{i;d;t}, \alpha_{P_t^\dagger}(\cdot, {d})\rangle_{L^2(P^\dagger_{X|d;t})}\widetilde{u}_{i;d;t}\|_{L^2(P^\dagger_{X|d;t})}},
    \end{align}
    and set $\widetilde{u}_{l_t+1;d;t}$ equal to an arbitrary unit length essentially uniformly bounded vector orthogonal to $\widetilde{u}_{i;d;t}$ for $i \leq l_t$, otherwise. We do similarly for $\widetilde{v}_{l_t+1;d;t}$. Recall that $\alpha_{P_t^\dagger}(\cdot, {d})$, $r_{P_t^\dagger}(\cdot,{d})$ are essentially uniformly bounded by \cref{lemma:l1-bounded-functional}, and hence $\|\widetilde{u}_{l_t+1;d;t} \|_{L^\infty(P^\dagger_{X;t})}, \|\widetilde{v}_{l_t+1;d;t} \|_{L^\infty(P^\dagger_{Z;t})} < \infty$.

    We then choose $\widetilde{\sigma}_{l_t+1;d;t}  > 0$ sufficiently small such that
    \begin{align}
        \label{eq:sigma-1-d-bound-tangent}
        \widetilde{\sigma}_{l_t + 1;d;t}\|\widetilde{u}_{l_t+1;d;t}\widetilde{v}_{l_t+1;d;t}\|_{L^\infty(P^\dagger_{X,Z;t})} < \min\left(L_{bounded}^{-4}\min\left(t^3, \frac{1}{2}L_{bounded}^{-4}L_{gap;P_t^\dagger}\right), \{L_{t}^{l_t}\}^{-1}\right).
    \end{align}
    Then, letting
    \begin{align*}
        {k}_{d;t} \coloneqq \sum_{i=1}^{l_t+1} \widetilde{\sigma}_{i;d;t}\widetilde{u}_{i;d;t}\widetilde{v}_{i;d;t}
    \end{align*}
    we have that
    \begin{align*}
        \nonumber \|k_{d;t}^{l_t} - {k}_{d;t}\|_{L^{\infty}(P^\dagger_{X,Z;t})} <  \min\left(L_{bounded}^{-4}\min\left(t^3, \frac{1}{2}L_{bounded}^{-4}L_{gap;P_t^\dagger}\right), \{L_{t}^{l_t}\}^{-1}\right).
    \end{align*}
    Because ${k}^{l_t}_{d;t} >\{L_{t}^{l_t}\}^{-1} > 0$, $\|k_{d;t}^{l_t} - {k}_{d;t}\|_{L^{\infty}(P^\dagger_{X,Z;t})} < \{L_{t}^{l_t}\}^{-1}$ implies that ${k}_{d;t} > 0$.
    Defining
    \begin{align*}
        {p}_{X,Z|{d};t}(x, z) \coloneqq {k}_{d;t}(x, z) p^\dagger_{X|{d};t}(x)p^\dagger_{Z|{d};t}(z)
    \end{align*}
    then gives that
    \begin{align}
        \label{eq:p-l_d-t-tilde-infinity-bound}
        \begin{split}
        \|p_{X, Z|{d};t}^{l_t} - {p}_{X, Z|{d};t}\|_{L^{\infty}(P^\dagger_{X,Z;t})}
        &\leq\|k_{d;t}^{l_t} - {k}_{d;t}\|_{L^{\infty}(P^\dagger_{X,Z;t})}\|p^\dagger_{X|{d};t}\|_{L^\infty(P^\dagger_{X;t})}\|p_{Z|{d};t}^\dagger\|_{L^\infty(P^\dagger_{Z;t})} \\
        &< \min\left(t^3, \frac{1}{2}L_{bounded}^{-4}L_{gap;P_t^\dagger}\right).
        \end{split}
    \end{align}
    We also have that ${p}_{X,Z|{d};t} > 0$. We further argue that ${p}_{X,Z|{d};t}$ is a density for $(X, Z)$. Indeed,
    \begin{align*}
        &\int\int {k}_{d;t}(x, z)p^\dagger_{X|{d};t}(x)p^\dagger_{Z|{d};t}(z)d\lambda_X(x)d\lambda_Z(z)\\
        =&\int\int \sum_{i=1}^{l_t+1}\widetilde{\sigma}_{i;d;t}\widetilde{u}_{i;d;t}(x)\widetilde{v}_{i;d;t}(z)p^\dagger_{X|{d};t}(x)p^\dagger_{Z|{d};t}(z)d\lambda_X(x)d\lambda_Z(z)\\
        =& \sum_{i=1}^{l_t+1}\widetilde{\sigma}_{i;d;t}\langle \widetilde{u}_{i;d;t}, 1\rangle_{L^2(P^\dagger_{X|d;t})}\langle \widetilde{v}_{i;d;t}, 1\rangle_{L^2(P^\dagger_{Z|d;t})}\\
        =& 1
    \end{align*}
    which again follows because $\widetilde{\sigma}_{1;d;t} = 1$, $\widetilde{u}_{1;d;t} = 1$, $\widetilde{v}_{1;d;t}=1$. Furthermore,
    \begin{align*}
        {p}_{X|{d};t}(x) =& \int p_{X, Z|{d};t}(x, z)d\lambda_Z(z)\\
        =&\int {k}_{d;t}(x, z)p^\dagger_{X|{d};t}(x)p^\dagger_{Z|{d};t}(z)d\lambda_Z(z)\\
        =&\int \sum_{i=1}^{l_t+1}\widetilde{\sigma}_{i;d;t}\widetilde{u}_{i;d;t}(x)\widetilde{v}_{i;d;t}(z)p^\dagger_{X|{d};t}(x)p^\dagger_{Z|{d};t}(z)d\lambda_Z(z)\\
        =& p^\dagger_{X|{d};t}(x)\sum_{i=1}^{l_t+1}\widetilde{\sigma}_{i;d;t}\widetilde{u}_{i;d;t}(x)\langle \widetilde{v}_{i;d;t}, 1\rangle_{L^2(P^\dagger_{Z|d;t})}\\
        =& p^\dagger_{X|{d};t}(x),
    \end{align*}
    and similarly ${p}_{Z|{d};t} = p^\dagger_{Z|{d};t}$.

    Now, for each $d$,
    \begin{align*}
        \|{p}_{X,Z|{d};t} - p^\dagger_{X, Z|{d};t}\|_{L^\infty(P^\dagger_{X,Z;t})}
        &\leq \|{p}^\dagger_{X,Z|{d};t} - p^{l_t}_{X, Z|{d};t}\|_{L^\infty(P^\dagger_{X,Z;t})} + \|p^{l_t}_{X, Z|{d};t} - p_{X, Z|{d};t}\|_{L^\infty(P^\dagger_{X,Z;t})} \\
        &< L_{bounded}^{-4}L_{gap;P_t^\dagger}
    \end{align*}
    Let ${p}_{X, Z, D;t}(x, z, d) \coloneqq {p}_{X,Z|{d};t}(x, z)p^\dagger_{D;t}(d)$ denote the corresponding joint density of $(X, Z, D)$. Because $p^\dagger_{D;t}(d) < L_{bounded}$ for each $d$ by \cref{remark:marginal-conditional-density-bounds}, and because ${p}_{X, Z, D;t}$ and $p^\dagger_{X, Z, D;t}$ both vanish off the support of $P_t^\dagger$, the above display gives
    \begin{align*}
        \|{p}_{X, Z, D;t} - p^\dagger_{X, Z, D;t}\|_{L^\infty(\lambda_X\otimes\lambda_Z\otimes\lambda_D)} = \max_{d \in \textsf{Supp}[D;P_t^\dagger]} p^\dagger_{D;t}(d)\|{p}_{X,Z|{d};t} - p^\dagger_{X, Z|{d};t}\|_{L^\infty(P^\dagger_{X,Z;t})} < L_{bounded}^{-3}L_{gap;P_t^\dagger} \leq L_{bounded}^{-1}L_{gap;P_t^\dagger}.
    \end{align*}
    With $P \coloneqq P_t^\dagger$, the threshold in \cref{lemma:constant-riesz-representers} is $L_{bounded}^{-1}L_{gap;P_t^\dagger}\{\lambda_X(\mathcal X)\lambda_Z(\mathcal Z)\lambda_D(\mc{D})\}^{-1} = L_{bounded}^{-1}L_{gap;P_t^\dagger}$, because the reference measures are normalized. Hence the conditions of \cref{lemma:constant-riesz-representers} are satisfied with ${p}_{X, Z, D;t}$ the density $p'$ and $P_t^\dagger$ the reference distribution $P$ in that Lemma. We let ${P}_t$ be $P^\dagger$ of that Lemma. By that Lemma,
    \begin{enumerate}
        \item $P_t^\dagger$ and ${P}_t$ are mutually absolutely continuous
        \item $L_{bounded}^{-1} < \essinf {p}_t \leq \esssup p_t < L_{bounded}$
        \item $P_t$ satisfies \cref{as:continuous-linear-maps} and $\alpha_{P_t}  = \alpha_{{P}_t^\dagger} $, $r_{{P}_t}  = r_{{P}_t^\dagger} $
        \item $H^2({P}_t, P^\dagger_t) \leq \frac{1}{2}\|{p}_{X, Z, D;t} - p^\dagger_{X, Z, D;t}\|_{L^1(\lambda_X\otimes\lambda_Z\otimes\lambda_D)}$.
    \end{enumerate}

    Now,
    \begin{align*}
        H^2(P_t, P_t^\dagger)\leq&\frac{1}{2}\|{p}_{X, Z, D;t} - p^\dagger_{X, Z, D;t}\|_{L^1(\lambda_X\otimes\lambda_Z\otimes\lambda_D)} \\
        \leq& L'\|{p}_{X, Z, D;t} - p^\dagger_{X, Z, D;t}\|_{L^\infty(\lambda_X\otimes\lambda_Z\otimes\lambda_D)}\\
        =&L'\|{p}_{X, Z, D;t} - p^\dagger_{X, Z, D;t}\|_{L^\infty(P^\dagger_{X,Z,D;t})}\\
        \leq&\sum_{d \in \textsf{Supp}[D;P_t^\dagger]}L'\|{p}_{X, Z|{d};t} - p_{X, Z|{d};t}^\dagger\|_{L^\infty(P^\dagger_{X,Z,D;t})}P_t^\dagger(D=d)\\
        \leq & \sum_{d \in \textsf{Supp}[D;P_t^\dagger]}P_t^\dagger(D=d)L' \Big\{\|{p}_{X, Z|{d};t} - p^{l_t}_{X, Z|{d};t}\|_{L^\infty(P^\dagger_{X,Z;t})}\\
        &\qquad + \|{p}^{l_t}_{X, Z|{d};t} - p^\dagger_{X, Z|{d};t}\|_{L^\infty(P^\dagger_{X,Z;t})}\Big\}\\
        \leq& 2L't^3
    \end{align*}
    for some constant $L'$ not depending on $t$.

    The final inequality follows from \eqref{eq:p-l_t-t-infinity-bound} and \eqref{eq:p-l_d-t-tilde-infinity-bound}. Consequently,
    \begin{align}
        \label{eq:P_t-widetilde-P_t-distance}
        H^2({P}_t, \widetilde{P}_t) \leq H^2({P}_t, P_t^\dagger) +2  H({P}_t, P_t^\dagger)H({P}_t^\dagger, \widetilde{P}_t) + H^2({P}_t^\dagger, \widetilde{P}_t) <  2L't^3 + 2\sqrt{2L'}t^3 + t^3
    \end{align}
    where the final inequality follows from \eqref{eq:tangent-space-dagger-tilde-distance}.

    We now argue $P_t$ belongs to $\mc{P}$ for each $t$, i.e. Assumptions \ref{as:continuous-linear-maps}, \ref{as:solutions-exist}, \ref{as:beta-source} (for all $\beta_h, \beta_q \in [0, \infty)$), \ref{as:mutual-absolute-continuity} (relative to $P_0$), and \ref{as:dominating-measure-support} hold.  \cref{as:continuous-linear-maps} holds by the construction in \cref{lemma:constant-riesz-representers} and because it holds at $P_t^\dagger$. \cref{as:solutions-exist} holds by the choice of $\widetilde{u}_{l_t + 1;d;t}$ and $\widetilde{v}_{l_t + 1;{d};t}$ in \eqref{eq:residual-projection-generic-d-tangent}. Furthermore, because $\mc{E}_{{P}_t;d}$, $\mc{E}_{{P}_t;d}^*$ only have a finite number of singular values greater than 0, \cref{as:solutions-exist} implies \cref{as:beta-source} for any choice of $\beta_h, \beta_q \in [0, \infty)$. ${P}_t$ is mutually absolutely continuous with $P_t^\dagger$ by \cref{lemma:constant-riesz-representers}, and $P_t^\dagger$ is mutually absolutely continuous with $P_0$, so \cref{as:mutual-absolute-continuity} holds. Finally, \cref{as:dominating-measure-support} also holds by the construction of ${P}_t$ in \cref{lemma:constant-riesz-representers}.

    We conclude by arguing that the submodel $\{P_t: t\in [0, t_s)\}$, with $P_{t=0} \coloneqq P_0$, is quadratic mean differentiable at $t=0$ with score $s$. By the above arguments,
    \begin{align*}
        &\int \left\{\frac{\sqrt{p_t} - \sqrt{p_0}}{t} - \frac{1}{2}s\sqrt{p_0}\right\}^2 d\lambda(o)\\
        =& \int \left\{\frac{\sqrt{p_t} - \sqrt{\widetilde{p}_t} + \sqrt{\widetilde{p}_t} - \sqrt{p_0}}{t} - \frac{1}{2}s\sqrt{p}_0\right\}^2 d\lambda(o)\\
        \leq& \frac{4}{t^2}H^2(P_t, \widetilde{P}_t)+ o(1)\\
        =& o(1).
    \end{align*}
    The inequality is a consequence of the triangle inequality applied to the squared Hellinger distance and that $\widetilde{P}_t$ is QMD at $t=0$ with score $s$. The final equality follows from \eqref{eq:P_t-widetilde-P_t-distance}. Hence, $\{P_t: t\in [0, t_s)\} \subseteq \mc{P}$ is a QMD at $t=0$ submodel with score $s$. As $s$ was arbitrary beyond belonging to the dense subset $L^\infty(P_0)\cap L^2_0(P_0)$ of $L^2_0(P_0)$, $\mc{T}(P_0, \mc{P}) = L^2_0(P_0)$.
\end{proof}

\begin{proof}[Proof of \cref{lemma:constant-rank-tangent-space}]
    We first provide a sketch of the proofs. \cref{lemma:discrete-tangent-space:part:constant-rank} relies on the fact that a differentiable (at $t=0$) matrix valued mapping $t \mapsto A(t)$ with constant rank for $t$ near $0$ has differentiable (at $t=0$) Moore-Penrose pseudo-inverse. Then, both sides of the equation $r_{P_t} = \mc{E}_{P_t}h_{P_t}^\dagger$, with $h_{P_t}^\dagger \coloneqq \mc{E}_{P_t}^\dagger r_{P_t}$, are differentiable in $t$, giving a restriction on the scores $s$. We extend slightly the argument of \cite{severiniEfficiencyBoundsEstimating2012}, noting that the above moment restriction holds for any choice of solution, not just the minimum norm solution, leading to an infinite collection of restrictions on the scores $s$. The same argument holds for the dual inverse problem, giving yet more restrictions.

    \cref{lemma:discrete-tangent-space:part:full-rank} uses that when $\mc{E}_{P_0;d}$ are full rank for each $d$, there exists an $\epsilon > 0$ such that the local model $\mc{P}_{\epsilon} \coloneqq \{P \in \mc{P}: H^2(P, P_0) < \epsilon\}$ satisfies the constant rank assumption because the minimum singular value of $\mc{E}_{P_0;d}$ is separated from zero. The maximal tangent space of $\mc{P}_{\epsilon} $ at $P_0$ is equal to the maximal tangent space of $\mc{P} $ at $P_0$ for every $\epsilon >0$, giving the result.

    \textbf{Proof of \cref{lemma:discrete-tangent-space:part:constant-rank}:}

    Let $\{P_{t}: t\in [0, \delta)\}\subseteq \mc{P}_{\text{rank}}$ be any one-dimensional QMD submodel with $P_{t=0} = P_0$. Let $s \in L^2_0(P_0)$ denote its score. We first show that $s$ satisfies that
    \begin{align*}
        \dot r_{P_0}[s](Z, D) - E_{P_0}[h_{P_0}(X, D)s_{X|Z, D}(X, Z, D)|Z, D]  \in \mathscr{R}(\mc{E}_{P_0})
    \end{align*}
    for any $h_{P_0}$ satisfying $\mc{E}_{P_0}h_{P_0} = r_{P_0}$. Because $X, Z$ are supported on $\{1, \dots, K_X\}$, $\{1, \dots, K_Z\}$, respectively, we may equivalently view $\mc{E}_{P_0;d}$ as a $K_Z \times K_X$-matrix operating on vectors of length $K_X$, after associating $L^2(P_{X|d;0})$ with $\mathbb{R}^{K_X}$ and $L^2(P_{Z|d;0})$ with $\mathbb{R}^{K_Z}$. I.e., we write $f \in L^2(P_{X|d;0})$ in vector form as $f \coloneqq (f(1), \dots, f(K_X))^T$, and write $\mc{E}_{P_0;d}$ in matrix form as
    \begin{align*}
        \mc{E}_{P_0;d} = \begin{pmatrix}
            P(X=1|Z=1, D=d) & \cdots & P(X=K_X|Z=1, D=d)\\
            \vdots \\
            P(X=1|Z=K_Z, D=d) & \cdots & P(X=K_X|Z=K_Z, D=d)
        \end{pmatrix}.
    \end{align*}
    We similarly define $\mc{E}_{P_t;d}$ at any $P_t$. We let $\mc{E}_{P_t;d}^\dagger, \mc{E}_{P_0;d}^\dagger$ denote the Moore-Penrose pseudo-inverses of these operators, taken, as elsewhere in this work, with respect to the inner products of $L^2(P_{X|d;t})$ and $L^2(P_{Z|d;t})$. All projections appearing below are orthogonal with respect to these same inner products. In particular, $\mc{E}_{P_t;d}^\dagger r_{P_t}(\cdot, d)$ is the minimum $L^2(P_{X|d;t})$-norm solution to the primal inverse problem at $P_t$.

    Now, because $\{P_t:t \in [0,\delta)\}$ is QMD at $t=0$, we have that for each $(x, z) \in \{1, \dots, K_X\}\times \{1, \dots, K_Z\}$, the derivative
    \begin{align*}
        \left.\frac{d}{dt}P_t(X=x|Z=z, D=d)\right\vert_{t=0}
    \end{align*}
    exists and is equal to $\{E_{P_0}[s(O)|X=x, Z=z, D=d] - E_{P_0}[s(O)|Z=z, D=d]\}P_0(X=x|Z=z, D=d) = s_{X|Z, D}(x, z, d)P_0(X=x|Z=z, D=d)$. Then, the $K_Z \times K_X$-matrix valued mapping $t \mapsto \mc{E}_{P_t;d}$ is differentiable at $t=0$. Consequently, because $\text{rank}(\mc{E}_{P_t;d}) = \text{rank}(\mc{E}_{P_0;d})$ for each $t$, the mappings
    \begin{align*}
        t \mapsto \mc{E}^\dagger_{P_t;d}, && t \mapsto  I - \mc{E}^\dagger_{P_t;d}\mc{E}_{P_t;d}\eqqcolon \Pi_{\mathscr{N}(\mc{E}_{P_t;d})}
    \end{align*}
    are also differentiable at $t=0$. To see this, we note that the pseudo-inverse taken with respect to the inner products of $L^2(P_{X|d;t})$ and $L^2(P_{Z|d;t})$ is the pseudo-inverse of the matrix representing $\mc{E}_{P_t;d}$ in orthonormal bases of those spaces, transported back to the canonical bases. That matrix has $(z, x)$ entry $P_t(X=x, Z=z|D=d)\{P_t(X=x|D=d)P_t(Z=z|D=d)\}^{-1/2}$, which is differentiable at $t=0$ because $\{P_t:t\in[0,\delta)\}$ is QMD at $t=0$ and, by \cref{as:dominating-measure-support}, these probabilities are bounded away from zero, and its rank is $\text{rank}(\mc{E}_{P_t;d})$, which is constant in $t$ because $\{P_t:t \in [0, \delta)\} \subseteq \mc{P}_{\text{rank}}$. \cref{lemma:pseudo-inverse-differentiability} then gives that the pseudo-inverse of that matrix is differentiable at $t=0$, and the changes of basis are diagonal with entries differentiable at $t=0$ and bounded away from zero. We note that $\Pi_{\mathscr{N}(\mc{E}_{P_t;d})} = I - \mc{E}^\dagger_{P_t;d}\mc{E}_{P_t;d}$ is the projection onto $\mathscr{N}(\mc{E}_{P_t;d})$, orthogonal in $L^2(P_{X|d;t})$.

    Now, because $P_t \in \mc{P}$ for each $t$,  by \cref{as:continuous-linear-maps} the Riesz representer $r_{P_t}$ exists. Furthermore, by \cref{as:solutions-exist}, $\mc{E}_{P_t;d}\mc{E}_{P_t;d}^\dagger r_{P_t}(\cdot, d) = r_{P_t}(\cdot, d)$. Additionally, for any $f \in L^2(P_{X|d;t})$, we have that $\mc{E}_{P_t;d}\{\mc{E}_{P_t;d}^\dagger r_{P_t}(\cdot, d) + \Pi_{\mathscr{N}(\mc{E}_{P_t;d})}f\} = r_{P_t}(\cdot, d)$. In particular, let $\eta(\cdot, d;h_{P_0}) \coloneqq \Pi_{\mathscr{N}(\mc{E}_{P_0;d})}h_{P_0}(\cdot, d)$. It follows that $h_{P_0}(\cdot, d) = \mc{E}_{P_0;d}^\dagger r_{P_0}(\cdot, d) + \eta(\cdot, d;h_{P_0})$. We then have that
    \begin{align*}
        \left.\frac{d}{dt}r_{P_t}(\cdot, d)\right\vert_{t=0} = \left.\frac{d}{dt}\mc{E}_{P_{t;d}}\{\mc{E}_{P_{t;d}}^\dagger r_{P_t}(\cdot, d) + \Pi_{\mathscr{N}(\mc{E}_{P_t;d})}\eta(\cdot, d;h_{P_0})\}\right \vert_{t=0}.
    \end{align*}
    Applying the product rule then gives that
    \begin{align*}
        \left.\frac{d}{dt}r_{P_t}(\cdot, d)\right\vert_{t=0} =& \left.\frac{d}{dt}\mc{E}_{P_t;d}\{\mc{E}_{P_{0;d}}^\dagger r_{P_0}(\cdot, d) + \Pi_{\mathscr{N}(\mc{E}_{P_0;d})}\eta(\cdot, d;h_{P_0})\}\right\vert_{t=0}\\
        &+\mc{E}_{P_0;d}\left.\frac{d}{dt}\{\mc{E}_{P_t;d}^\dagger r_{P_t}(\cdot, d) + \Pi_{\mathscr{N}(\mc{E}_{P_t;d})}\eta(\cdot, d;h_{P_0})\}\right\vert_{t=0}.
    \end{align*}
    The second term on the right hand side of the above display is an element of $\mathscr{R}(\mc{E}_{P_0;d})$. Now, because $\eta(\cdot, d;h_{P_0}) \in \mathscr{N}(\mc{E}_{P_0;d})$, $\Pi_{\mathscr{N}(\mc{E}_{P_0;d})}\eta(\cdot, d;h_{P_0}) = \eta(\cdot, d;h_{P_0})$, and so $\mc{E}_{P_{0;d}}^\dagger r_{P_0}(\cdot, d) + \Pi_{\mathscr{N}(\mc{E}_{P_0;d})}\eta(\cdot, d;h_{P_0}) = h_{P_0}(\cdot, d)$. Consequently, the first term on the right hand side written as a function of $Z$ instead of a vector, is equal to
    \begin{align*}
        E_{P_0}[h_{P_0}(X, D)s_{X|Z, D}(X, Z, D)|Z, D=d].
    \end{align*}
    Then, recalling that $\left.\frac{d}{dt}r_{P_t}(\cdot, d)\right\vert_{t=0} = \dot r_{P_0}[s](\cdot, d) $ gives, for each $d$, that $s$ must satisfy the restraint that
    \begin{align*}
        \dot r_{P_0}[s](\cdot, d) - E_{P_0}[h_{P_0}(X, D)s_{X|Z, D}(X, Z, D)|Z, D=d] \in \mathscr{R}(\mc{E}_{P_0;d})
    \end{align*}
    and hence, because $g \in \mathscr{R}(\mc{E}_{P_0})$ if and only if $g(\cdot, d) \in \mathscr{R}(\mc{E}_{P_0;d})$ for each $d$, that $s$ must satisfy
    \begin{align*}
        \dot r_{P_0}[s] - E_{P_0}[h_{P_0}(X, D)s_{X|Z, D}(X, Z, D)|Z, D] \in \mathscr{R}(\mc{E}_{P_0})
    \end{align*}
    for every $h_{P_0}$ solving $\mc{E}_{P_0}h_{P_0} = r_{P_0}$.

    The same argument, applied to the matrix representation of $\mc{E}^*_{P_0;d}$ with $(x, z)$ entry $P(Z=z|X=x, D=d)$, and using that $\text{rank}(\mc{E}^*_{P;d}) = \text{rank}(\mc{E}_{P;d})$ for every $P$, so that $\text{rank}(\mc{E}^*_{P_t;d})$ is also constant in $t$, shows that $s$ must also satisfy
    \begin{align*}
        \dot \alpha_{P_0}[s] - E_{P_0}[q_{P_0}(Z, D)s_{Z|X, D}(X, Z, D)|X, D] \in \mathscr{R}(\mc{E}^*_{P_0}).
    \end{align*}
    for every $q_{P_0}$ solving $\mc{E}_{P_0}^*q_{P_0} = \alpha_{P_0}$.

    Because $X, Z, D$ are finite discrete, $\mathscr{R}(\mc{E}_{P_0}) = \mathscr{N}(\mc{E}_{P_0}^*)^\perp$ and $\mathscr{R}(\mc{E}_{P_0}^*) = \mathscr{N}(\mc{E}_{P_0})^\perp$, so the above restrictions are exactly those defining $\mc{C}_{rank}(P_0)$. Finally, we note that $\mc{C}_{rank}(P_0)$ is a closed linear subspace of $L^2_0(P_0)$. Indeed, $s \mapsto \dot r_{P_0}[s]$ and $s \mapsto \dot \alpha_{P_0}[s]$ are linear and bounded because each of their components is an inner product against a fixed element of $L^2_0(P_0)$; the maps $s \mapsto E_{P_0}[h_{P_0}(X, D)s_{X|Z, D}(X, Z, D)|Z, D]$ and $s \mapsto E_{P_0}[q_{P_0}(Z, D)s_{Z|X, D}(X, Z, D)|X, D]$ are linear and bounded; and $\mathscr{R}(\mc{E}_{P_0})$, $\mathscr{R}(\mc{E}^*_{P_0})$ are closed linear subspaces because $X, Z, D$ are finite discrete. As the score of every submodel in $\mathscr{P}(P_0, \mc{P}_{\text{rank}})$ belongs to this set, so does the closed linear span of those scores, which is $\mc{T}(P_0, \mc{P}_{\text{rank}})$.

    Now we argue we can apply \cref{lemma:re-expressing-tangent-space} for each $h_{P_0}$. It suffices to show that $(1 + t s)p_0 \in \mc{P}_{\text{Riesz}}$ for each $s \in L^\infty(P_0) \cap L^2_0(P_0)$ for $t$ small enough. \cref{as:mutual-absolute-continuity,as:dominating-measure-support} hold for sufficiently small $t$ because $s \in L^\infty(P_0)$. \cref{as:continuous-linear-maps} then holds by \cref{lemma:l1-bounded-functional}, using \cref{as:m-m-tilde-bounded}. Applying by symmetry the analog of \cref{lemma:re-expressing-tangent-space} for $P \mapsto \alpha_{P}$ and $q_{P_0}$, and taking the intersection over all $q_{P_0}, h_{P_0}$ completes this part of the proof.

    \textbf{Proof of \cref{lemma:discrete-tangent-space:part:full-rank}:}

    By definition, $P_0 \in \mc{P}_{\text{rank}}$ and $\mc{P}_{\text{rank}} \subseteq \mc{P}$. Consequently, $\mc{T}(P_0, \mc{P}_{\text{rank}}) \subseteq \mc{T}(P_0, \mc{P})$. We now argue $ \mc{T}(P_0, \mc{P}) \subseteq\mc{T}(P_0, \mc{P}_{\text{rank}})$. Let $\{P_t:t\in [0, \delta)\}\subseteq \mc{P}$. Because this submodel is QMD, $\lim_{t\rightarrow 0 }H^2(P_t, P_0)=0$, and because $X, Z$ are supported on a rectangle for each $d$,  $P_t(X=x|Z=z, D=d) \rightarrow P_0(X=x|Z=z, D=d)$ for each $x, z, d$. Consequently, because $K_X, K_Z < \infty$, the matrix representing $\mc{E}_{P_t;d}$ in the canonical bases, with $(z, x)$ entry $P_t(X=x|Z=z, D=d)$, converges entrywise, and hence in Frobenius norm, to the matrix representing $\mc{E}_{P_0;d}$ as $t \rightarrow 0$. Consequently, the singular values of these matrix representations converge (\cite{stewart_matrix_1990} Chapter IV Theorem 1.3). On the other hand, because either $\mathscr{N}(\mc{E}_{P_0}) = \{0\}$ or $\mathscr{N}(\mc{E}^*_{P_0}) = \{0\}$, for each $d$ either $\mathscr{N}(\mc{E}_{P_0;d}) = \{0\}$ or $\mathscr{N}(\mc{E}^*_{P_0;d}) = \{0\}$, so that $\mc{E}_{P_0;d}$ has rank $\min(K_X, K_Z)$ and the smallest of the $\min(K_X, K_Z)$ singular values of its matrix representation is strictly positive.

    Then, for $T$ small enough, these singular values are strictly positive for all $t < T$. Hence, for $t < T$, $\text{rank}(\mc{E}_{P_{t};d}) = \min(K_X, K_Z) = \text{rank}(\mc{E}_{P_{0};d})$ for each $d$, the rank of an operator being independent of the bases in which it is represented. Then, $\{P_t:t\in [0, T)\} \subseteq \mc{P}_{\text{rank}}$. Consequently, the score of this submodel belongs to $\mc{T}(P_0, \mc{P}_{\text{rank}})$, and so $ \mc{T}(P_0, \mc{P}) \subseteq\mc{T}(P_0, \mc{P}_{\text{rank}})$. This completes the proof.
\end{proof}

\subsection{Proofs for \cref{subsec:recover-continuity}}

{

In the proof of \cref{thrm:bounded-means-continuous}
we use the following definitions. If $T:\mathcal{H} \mapsto \mathcal{H}$ is a compact, non-negative, and self-adjoint operator on the Hilbert space $\mathcal{H}$, with $\lambda_1,\lambda_2,\cdots$ being its eigenvalues  sorted in descending order and $u_1,u_2,\cdots$ the corresponding eigenfunctions, then $\operatorname{spec}(T) := \left\{\lambda_i\right\}_{i\ge 1} \cup\left\{0\right\}$ denotes the spectrum of $T$ and for any $a> 0$
\[
\Pi_{T}(a)h := \sum_{i:\lambda_i \ge a} \langle h,u_i\rangle_\mathcal{H} u_i, \qquad h \in \mathcal{H},
\]
denotes the orthogonal projection onto the span of the eigenfunctions with eigenvalues at least $a$. We write $\Pi[h|\mathcal{V}]$ for the orthogonal projection of $h$ onto a closed subspace $\mathcal{V}$ of $\mathcal{H}$.

Additionally, the Hausdorff distance between two sets $A$ and $B$ included in some metric space with metric $d$ is defined as $d_H(A,B):=\max\left\{\sup_{x\in B}d(x,A),\sup_{x\in A}d(x,B)\right\}$ where $d(x,A) := \inf_{y\in A}d(x,y)$.

\begin{proof} [Proof of \cref{thrm:bounded-means-continuous}]
We prove the result under \cref{as:uniform-norm-bound} \cref{as:uniform-norm-bound:part:h}. The proof under
\cref{as:uniform-norm-bound:part:q} is symmetric, using the dual inverse problem.

We begin with an informal sketch of the argument. The role of \cref{as:uniform-norm-bound} is not to make the
inverse problem well posed.  Even under the uniform bound on the minimum-norm
solution, the map
\[
    P \mapsto h_P^\dagger
\]
need not be continuous.  Small changes in the law can still produce noticeable
changes in the solution along directions that are only very weakly revealed by
the conditional expectation operator.  Thus, the proof cannot, and does not,
proceed by showing that \(h_{P_j}^\dagger\) converges to \(h_{P_0}^\dagger\).

The key observation is that continuity of the parameter requires less.  The
parameter does not depend on the whole function \(h_P^\dagger\) in isolation.
Rather, it depends on \(h_P^\dagger\) only through its weighted average against
the Riesz representer:
\(
    \varphi(P)
    =
    E_P\{\alpha_P(X,D)h_P^\dagger(X,D)\}.
\)
Thus, directions in which \(h_P^\dagger\) changes can matter for the parameter
only to the extent that they are also seen by \(\alpha_P\).  The proof exploits
this distinction.

The useful picture is the following.  For a fixed value of \(D=d\), decompose
directions of \(h(X,d)\) according to how strongly they are revealed by the
conditional distribution of \(X\) given \(Z\) and \(D=d\).  Some directions are
well revealed: changing \(h\) in those directions produces a visible change in
the conditional mean \(E\{h(X,d)\mid Z,D=d\}\).  Other directions are weakly
revealed: changing \(h\) in those directions produces only a very small change
in this conditional mean.  These weak directions are the source of ill-posedness,
because solving the inverse problem requires amplifying what is only faintly
seen in the observed law.

The proof separates these two types of directions.  First fix a threshold
\(a>0\).  Directions whose signal strength is at least \(a\) are treated as
stable.  On this stable part, the inverse problem behaves like an ordinary
finite-dimensional inverse problem with its singular values bounded away from
zero.  Combined with the uniform density bounds of \cref{as:dominating-measure-support} and the Hellinger continuity of the Riesz representer maps assumed in \cref{thrm:bounded-means-continuous}, along any sequence \(P_j\to P_0\), the stable part of the
minimum-norm solution varies continuously.  This is the portion of the argument
that shows that the truncated, well-identified component of
\(h_{P_j}^\dagger\) converges to the corresponding component of
\(h_{P_0}^\dagger\).

The unstable part consists of the remaining directions, whose signal strength is
below \(a\).  We do not try to show that the solution converges on this part.
Indeed, in an ill-posed problem there is no reason to expect such convergence.
The important point is instead that this unstable part is uniformly bounded in
norm by~\cref{as:uniform-norm-bound}.  Therefore, although it may move around as \(j\) changes,
it cannot become arbitrarily large.

It remains to see why this bounded but potentially nonconvergent part does not
affect the parameter.  This is where the fixed function \(\alpha_{P_0}\) enters.
Because \(\alpha_{P_0}\) belongs to the range of the adjoint inverse problem, it
has no component in directions that are completely invisible.  Moreover, since
\(\alpha_{P_0}\) is square integrable, its total mass in sufficiently weak
directions can be made arbitrarily small by choosing the threshold \(a\) small
enough.  Thus, even if the weak part of \(h_{P_j}^\dagger\) does not converge,
its inner product with \(\alpha_{P_0}\) is small: it is a bounded vector tested
against a part of \(\alpha_{P_0}\) whose norm is small.

The proof therefore follows a two-step logic.  First choose the cutoff \(a\)
small enough that the component of \(\alpha_{P_0}\) lying in directions weaker
than \(a\) is negligible.  Then, with this cutoff fixed, take \(j\) large enough
that the stable parts of the inverse problems under \(P_j\) and \(P_0\) are
close.  The stable part is controlled by continuity of the well-conditioned
finite-dimensional problem; the unstable part is controlled by the small tail
of \(\alpha_{P_0}\) together with the uniform bound on \(h_{P_j}^\dagger\).

Finally, the difference between integrating under \(P_j\) and integrating under
\(P_0\) is handled separately.  The uniform upper and lower density bounds in \cref{as:dominating-measure-support},
together with Hellinger convergence, imply that the relevant marginal density
ratios converge to one.  The Hellinger continuity of the Riesz representer maps assumed in \cref{thrm:bounded-means-continuous} gives convergence of the Riesz
representers.  Therefore, after rewriting the parameter under the fixed
\(P_0\)-inner product, the only nontrivial term is precisely the inner product
between \(\alpha_{P_0}\) and \(h_{P_j}^\dagger-h_{P_0}^\dagger\), which is
controlled by the stable/weak decomposition described above.

The eight steps of the proof implement this strategy.  Steps 1--3 put the
problem on fixed \(P_0\)-Hilbert spaces and show that, on these fixed spaces,
the conditional expectation operators under \(P_j\) approach the operator under
\(P_0\).  Step 4 transfers the assumed uniform bound on the minimum-norm
solutions to these fixed spaces.  Step 5 proves that, after discarding the
weakly identified directions below a fixed threshold, the stable parts of the
solutions converge.  Step 6 then uses the small spectral tail of
\(\alpha_{P_0}\), together with the uniform bound from Step 4, to show that the
remaining weak directions do not contribute to the relevant inner product.
Step 7 handles the change from integration under \(P_j\) to integration under
\(P_0\), using the convergence of the density ratios and of the Riesz
representers.  Finally, Step 8 combines these ingredients to show that
\(\varphi(P_j)-\varphi(P_0)\to0\).  Thus the proof follows exactly the heuristic
logic above: control the stable part by continuity, control the unstable part
by boundedness and a small tail of \(\alpha_{P_0}\), and then translate this
inner-product control into continuity of the parameter. We are now ready to start the proof.

Let $\{P_j\}_{j\geq 1}\subseteq\mc{P}_{cont}$ be any
sequence such that
$H(P_j,P_0)\longrightarrow 0.$
We must show that
$\varphi(P_j)\longrightarrow \varphi(P_0).$

Throughout the proof, let $p_j$ denote the density of $P_j$ with
respect to the common dominating product measure $\lambda$ from
~\cref{as:dominating-measure-support}. By \cref{as:dominating-measure-support}, there exists a constant
$L_{bounded}>1$, independent of $j$, such that
\begin{align*}
    L_{bounded}^{-1}
    \leq p_j(o)
    \leq L_{bounded}
\end{align*}
for $\lambda$-almost every $o$ and every $j\geq 0$. Here and below,
$p_0$ corresponds to $P_0$.
Since Hellinger convergence implies convergence in total variation, we have that
\begin{align}
    \|p_j-p_0\|_{L^1(\lambda)}
    \longrightarrow 0.
    \label{eq:hellinger-implies-L1}
\end{align}

We divide the proof into several steps.

\medskip

\noindent
\textbf{Step 1: Reduction to fixed Hilbert spaces.}

Fix $d\in\textsf{Supp}[D;P_0]$. Since $D$ has finite support and
all joint densities are uniformly bounded above and away from zero,
there exist constants $0<c_d<C_d<\infty$, independent of $j$, such
that
\begin{align}
    c_d
    \leq p_{X\mid d;j}(x)
    \leq C_d,
    \qquad
    c_d
    \leq p_{Z\mid d;j}(z)
    \leq C_d
    \label{eq:conditional-marginal-bounds}
\end{align}
for $\lambda_X$-almost every $x$, $\lambda_Z$-almost every $z$, and
every $j\geq 0$, where
$
    p_{X\mid d;j}(x)
    :=
    p_{X\mid D;j}(x\mid d)
    $ and $
    p_{Z\mid d;j}(z)
    :=
    p_{Z\mid D;j}(z\mid d).
$
Consequently, the spaces
   \( L^2(P_{X|d;j})\) and
    \(L^2(P_{X|d;0})\)
contain the same equivalence classes of measurable functions and have
uniformly equivalent norms. Specifically, there exists a constant
$K_d<\infty$, independent of $j$, such that
\begin{align}
    K_d^{-1}
    \|f\|_{L^2(P_{X|d;0})}
    \leq
    \|f\|_{L^2(P_{X|d;j})}
    \leq
    K_d
    \|f\|_{L^2(P_{X|d;0})}.
    \label{eq:X-norm-equivalence}
\end{align}
An analogous statement holds for
$L^2(P_{Z|d;j})$ and $L^2(P_{Z|d;0})$.
Define the fixed Hilbert spaces
\begin{align*}
    \mathbb H_{X,d}
    :=
    L^2(P_{X|d;0}) \qquad \text{and} \qquad
    \mathbb H_{Z,d}
    :=
    L^2(P_{Z|d;0}).
\end{align*}

For every $j\geq 0$, define \(
    \mathcal E_{j,d}:
    \mathbb H_{X,d}
    \longrightarrow
    \mathbb H_{Z,d}
    \)
by
\begin{align}
    (\mathcal E_{j,d}h)(z)
    &:=
    E_{P_j}\{h(X)\mid Z=z,D=d\}\nonumber\\
    &=
    \int h(x)
    p_{X\mid Z,d;j}(x\mid z)
    \,d\lambda_X(x),
    \label{eq:conditional-mean-kernel}
\end{align}
where
  \(  p_{X\mid Z,d;j}(x\mid z)
    =
    p_{X,Z\mid d;j}(x,z)/
         p_{Z\mid d;j}(z).
    \)
The operator $\mathcal E_{j,d}$ is therefore the conditional mean
operator under $P_j$, but it is regarded as an operator between the
fixed $P_0$-Hilbert spaces $\mathbb H_{X,d}$ and $\mathbb H_{Z,d}$.
Let $\mathcal E_{j,d}^{*0}$ denote its adjoint relative to the
$P_0$-inner products on these fixed spaces. The superscript $0$ is
included to emphasize that $\mathcal E_{j,d}^{*0}$ need not coincide
with the adjoint of the conditional mean operator when the domain and
codomain are equipped with their $P_j$-inner products.

\medskip

\noindent
\textbf{Step 2: Convergence of the conditional densities.}

We first show that
\begin{align}
    \left\|
    p_{X\mid Z,d;j}
    -
    p_{X\mid Z,d;0}
    \right\|_
    {L^1(\lambda_X\otimes\lambda_Z)}
    \longrightarrow 0.
    \label{eq:conditional-density-L1}
\end{align}

Write
\begin{align*}
    f_{j,d}(x,z)
    :=
    p_{X,Z\mid d;j}(x,z) \qquad \text{and} \qquad
    g_{j,d}(z)
    :=
    p_{Z\mid d;j}(z)
    =
    \int f_{j,d}(x,z)\,d\lambda_X(x).
\end{align*}
Then
\begin{align*}
    p_{X\mid Z,d;j}(x\mid z)
    =
    \frac{f_{j,d}(x,z)}{g_{j,d}(z)}.
\end{align*}

We first claim that
\begin{align}
    \|f_{j,d}-f_{0,d}\|_
    {L^1(\lambda_X\otimes\lambda_Z)}
    \longrightarrow 0.
    \label{eq:conditional-joint-L1}
\end{align}
Indeed, marginalizing \eqref{eq:hellinger-implies-L1} over $W$ gives
$
    \|p_{X,Z,D;j}-p_{X,Z,D;0}\|_
    {L^1(\lambda_X\otimes\lambda_Z\otimes\lambda_D)}
    \longrightarrow 0.
$
Because $D$ has finite support, it follows for each fixed $d$ that
$
    \left\|
    p_{X,Z,D;j}(\cdot,\cdot,d)
    -
    p_{X,Z,D;0}(\cdot,\cdot,d)
    \right\|_
    {L^1(\lambda_X\otimes\lambda_Z)}
    \longrightarrow 0.
$
In particular,
\begin{align*}
    P_j(D=d)\longrightarrow P_0(D=d).
\end{align*}
Since $P_0(D=d)>0$, and since the probabilities $P_j(D=d)$ are
uniformly bounded away from zero, division by $P_j(D=d)$ yields
\eqref{eq:conditional-joint-L1}.

Marginalization is an $L^1$-contraction, so
\begin{align}
    \|g_{j,d}-g_{0,d}\|_{L^1(\lambda_Z)}
    \leq
    \|f_{j,d}-f_{0,d}\|_
    {L^1(\lambda_X\otimes\lambda_Z)}
    \longrightarrow 0.
    \label{eq:conditional-marginal-L1}
\end{align}

Using the identity
\begin{align*}
    \frac{f_{j,d}}{g_{j,d}}
    -
    \frac{f_{0,d}}{g_{0,d}}
    =
    \frac{f_{j,d}-f_{0,d}}{g_{j,d}}
    +
    f_{0,d}
    \frac{g_{0,d}-g_{j,d}}
         {g_{j,d}g_{0,d}},
\end{align*}
the uniform lower bounds in
\eqref{eq:conditional-marginal-bounds} imply that
\begin{align*}
    \left|
    \frac{f_{j,d}}{g_{j,d}}
    -
    \frac{f_{0,d}}{g_{0,d}}
    \right|
    \leq
    c_d^{-1}|f_{j,d}-f_{0,d}|
    +
    c_d^{-2}f_{0,d}|g_{j,d}-g_{0,d}|.
\end{align*}
Integrating with respect to
$\lambda_X\otimes\lambda_Z$ gives
\begin{align*}
    \left\|
    p_{X\mid Z,d;j}
    -
    p_{X\mid Z,d;0}
    \right\|_
    {L^1(\lambda_X\otimes\lambda_Z)}
    &\leq
    c_d^{-1}
    \|f_{j,d}-f_{0,d}\|_
    {L^1(\lambda_X\otimes\lambda_Z)}
    \\
    &\quad+
    c_d^{-2}
    \int
    g_{0,d}(z)
    |g_{j,d}(z)-g_{0,d}(z)|
    \,d\lambda_Z(z).
\end{align*}
By the uniform upper bound for $g_{0,d}$,
\begin{align*}
    \left\|
    p_{X\mid Z,d;j}
    -
    p_{X\mid Z,d;0}
    \right\|_
    {L^1(\lambda_X\otimes\lambda_Z)}
    \leq
    c_d^{-1}\|f_{j,d}-f_{0,d}\|_{L^1}
    +
    c_d^{-2}C_d
    \|g_{j,d}-g_{0,d}\|_{L^1},
\end{align*}
and the right hand side converges to zero by
\eqref{eq:conditional-joint-L1} and
\eqref{eq:conditional-marginal-L1}. This proves
\eqref{eq:conditional-density-L1}.

The uniform upper and lower density bounds also imply that, for some
constant $M_d<\infty$,
\begin{align}
    \sup_{j\geq 0}
    \|p_{X\mid Z,d;j}\|_
    {L^\infty(\lambda_X\otimes\lambda_Z)}
    \leq M_d.
    \label{L_inf_XgivenZ_bound}
\end{align}
Consequently, $\|p_{X\mid Z,d;j}-p_{X\mid Z,d;0}\|_
    {L^\infty(\lambda_X\otimes\lambda_Z)}
    \leq 2M_d.$
Therefore,
\begin{align*}
    \|p_{X\mid Z,d;j}-p_{X\mid Z,d;0}\|_
    {L^2(\lambda_X\otimes\lambda_Z)}^2
    &\leq
    \|p_{X\mid Z,d;j}-p_{X\mid Z,d;0}\|_
    {L^1(\lambda_X\otimes\lambda_Z)}
    \|p_{X\mid Z,d;j}-p_{X\mid Z,d;0}\|_
    {L^\infty(\lambda_X\otimes\lambda_Z)}
    \\
    &\leq
    2M_d
    \|p_{X\mid Z,d;j}-p_{X\mid Z,d;0}\|_
    {L^1(\lambda_X\otimes\lambda_Z)}
    \longrightarrow0.
\end{align*}

Thus,
\begin{align}
    \left\|
    p_{X\mid Z,d;j}
    -
    p_{X\mid Z,d;0}
    \right\|_
    {L^2(\lambda_X\otimes\lambda_Z)}
    \longrightarrow 0.
    \label{eq:conditional-density-L2}
\end{align}

\medskip

\noindent
\textbf{Step 3: Operator-norm convergence of the conditional mean
operators.}

For $h\in\mathbb H_{X,d}$,
\begin{align*}
    \{(\mathcal E_{j,d}-\mathcal E_{0,d})h\}(z)
    =
    \int h(x)
    \Delta_{j,d}(x,z)
    \,d\lambda_X(x),
\end{align*}
where
\begin{align*}
    \Delta_{j,d}(x,z)
    :=
    p_{X\mid Z,d;j}(x\mid z)
    -
    p_{X\mid Z,d;0}(x\mid z).
\end{align*}

Viewed as an integral operator from
$L^2(P_{X|d;0})$ into $L^2(P_{Z|d;0})$, its kernel relative
to the domain measure $P_{X|d;0}$ is
\begin{align*}
    K_{j,d}(x,z)
    :=
    \frac{\Delta_{j,d}(x,z)}
         {p_{X\mid d;0}(x)}.
\end{align*}
Therefore,
\begin{align*}
    \|\mathcal E_{j,d}-\mathcal E_{0,d}\|_{\mathrm{op}}^2
    &\leq
    \|\mathcal E_{j,d}-\mathcal E_{0,d}\|_{\mathrm{HS}}^2\\
    &=
    \int
    |K_{j,d}(x,z)|^2
    \,dP_{X|d;0}(x)
    \,dP_{Z|d;0}(z)\\
    &=
    \int
    |\Delta_{j,d}(x,z)|^2
    \frac{p_{Z\mid d;0}(z)}
         {p_{X\mid d;0}(x)}
    \,d\lambda_X(x)\,d\lambda_Z(z).
\end{align*}
By \eqref{eq:conditional-marginal-bounds},
\begin{align*}
    \|\mathcal E_{j,d}-\mathcal E_{0,d}\|_{\mathrm{op}}^2
    \leq
    \frac{C_d}{c_d}
    \|\Delta_{j,d}\|_
    {L^2(\lambda_X\otimes\lambda_Z)}^2.
\end{align*}
It follows from \eqref{eq:conditional-density-L2} that
\begin{align}
    \|\mathcal E_{j,d}-\mathcal E_{0,d}\|_{\mathrm{op}}
    \longrightarrow 0.
    \label{eq:operator-norm-convergence}
\end{align}
Consequently, the adjoints satisfy
\begin{align}
    \|\mathcal E_{j,d}^{*0}
      -\mathcal E_{0,d}^{*0}\|_{\mathrm{op}}
    \longrightarrow 0.
    \label{eq:adjoint-operator-convergence}
\end{align}

\medskip

\noindent
\textbf{Step 4: Uniformly bounded solutions in the fixed Hilbert
spaces.}

Let
$
    \widetilde h_j^\dagger(\cdot,d)
$
denote the minimum-$L^2(P_{X|d;j})$-norm solution of
\begin{align}
    \mathcal E_{j,d}h=r_{P_j}(\cdot,d).
    \label{eq:fixed-space-inverse-problem}
\end{align}
By \cref{as:uniform-norm-bound} \cref{as:uniform-norm-bound:part:h},
 \(   \|\widetilde h_j^\dagger(\cdot,d)\|_
    {L^2(P_{X|d;j})}
    \leq C
\)
for every $j$ and every $d$.

Let
$
    \bar h_j^\dagger(\cdot,d)
$
denote the minimum-$\mathbb H_{X,d}$-norm solution to equation ~\eqref{eq:fixed-space-inverse-problem}.
The functions $\widetilde h_j^\dagger(\cdot,d)$ and
$\bar h_j^\dagger(\cdot,d)$ solve the same pointwise conditional
moment equation; they differ only because minimum norm is defined
using different, although equivalent, norms.

By minimum-norm optimality and \eqref{eq:X-norm-equivalence},
\begin{align*}
    \|\bar h_j^\dagger(\cdot,d)\|_{\mathbb H_{X,d}}
    \leq
    \|\widetilde h_j^\dagger(\cdot,d)\|_{\mathbb H_{X,d}}
    \leq
    K_d
    \|\widetilde h_j^\dagger(\cdot,d)\|_
    {L^2(P_{X|d;j})}
    \leq K_dC.
\end{align*}
Since $\textsf{Supp}[D;P_0]$ is finite, there exists
$B<\infty$ such that
\begin{align}
    \sup_{j\geq 0}
    \sup_{d\in\textsf{Supp}[D;P_0]}
    \|\bar h_j^\dagger(\cdot,d)\|_{\mathbb H_{X,d}}
    \leq B.
    \label{eq:fixed-space-uniform-bound}
\end{align}

\medskip

\textbf{Step 5: Convergence of spectrally truncated solutions.}

Fix $d\in\textsf{Supp}[D;P_0]$ and define
\[
    T_{j,d}
    :=
    \mathcal E_{j,d}^{*0}\mathcal E_{j,d}
    :
    \mathbb H_{X,d}
    \longrightarrow
    \mathbb H_{X,d}.
\]
Each $T_{j,d}$ is compact, nonnegative, and self-adjoint. Moreover,
by ~\eqref{eq:operator-norm-convergence} and ~\eqref{eq:adjoint-operator-convergence},
\begin{align}
    \|T_{j,d}-T_{0,d}\|_{\mathrm{op}}
    &=
    \|
        \mathcal E_{j,d}^{*0}\mathcal E_{j,d}
        -
        \mathcal E_{0,d}^{*0}\mathcal E_{0,d}
      \|_{\mathrm{op}}
      \notag\\
    &\leq
      \|\mathcal E_{j,d}^{*0}-\mathcal E_{0,d}^{*0}\|_{\mathrm{op}}
      \|\mathcal E_{j,d}\|_{\mathrm{op}}
      +
      \|\mathcal E_{0,d}^{*0}\|_{\mathrm{op}}
      \|\mathcal E_{j,d}-\mathcal E_{0,d}\|_{\mathrm{op}}
      \longrightarrow 0.
    \label{eq:squared_operator_convergence}
\end{align}

 Let $a>0$ satisfy
\(
    a\notin \operatorname{spec}(T_{0,d}).
\)
 Choose
$\eta_a>0$ such that
\begin{equation}
    \operatorname{spec}(T_{0,d})
    \cap
    [a-2\eta_a,a+2\eta_a]
    =
    \varnothing.
    \label{eq:double_neighborhood}
\end{equation}
By~\eqref{eq:squared_operator_convergence}, this implies that there exists $j_0$ such that for $j \ge j_0$,
\begin{equation}
    \operatorname{spec}(T_{j,d})
    \cap
    [a-\eta_a,a+\eta_a]
    =
    \varnothing.
    \label{eq:single_neighborhood}
\end{equation}
To see this, we first note that because $T_{j,d}$ and $T_{0,d}$ are compact self-adjoint operators it holds that (see \cite{kato_perturbation_1995} Chapter V, Section 4.3, Theorem 4.10):
\[
d_{\mathrm H}\!\left(
    \operatorname{spec}(T_{j,d}),
    \operatorname{spec}(T_{0,d})
\right)
\leq
\|T_{j,d}-T_{0,d}\|_{\mathrm{op}}.
\]
Now, by~\eqref{eq:squared_operator_convergence}, there exists $j_0$ such that for $j \ge j_0$, the right-hand side is smaller than
$\eta_a$. Hence, if
$\lambda\in\operatorname{spec}(T_{j,d})\cap
[a-\eta_a,a+\eta_a]$ for $j\ge j_0$, there would exist
$\lambda_0\in\operatorname{spec}(T_{0,d})$ such that
$|\lambda-\lambda_0|<\eta_a$, implying
$\lambda_0\in[a-2\eta_a,a+2\eta_a]$, in contradiction with ~\eqref{eq:double_neighborhood}. This proves ~\eqref{eq:single_neighborhood}.

Next, for any $j\ge0$ let $\Pi_{T_{j,d}}(a) := \Pi_{j,d}(a)$. We will show that
\begin{equation}
    \|\Pi_{j,d}(a)-\Pi_{0,d}(a)\|_{\mathrm{op}}
    \longrightarrow 0.
    \label{eq:truncated_convergence}
\end{equation}
To see this, choose a bounded continuous function
$\chi_a:[0,\infty)\to[0,1]$ satisfying
\[
    \chi_a(t)=0
    \quad\text{for }t\leq a-\eta_a,
    \qquad
    \chi_a(t)=1
    \quad\text{for }t\geq a+\eta_a.
\]
By the Spectral Theorem, $\chi_a(T_{j,d})h=\sum_i\chi_a(\lambda_{j,i})\langle h,u_{j,i}\rangle u_{j,i} + \chi_a(0) \Pi[h|\mathscr{N}(T_{j,d})]$ where $\lambda_{j,i}, i = 1,2,\cdots$ and $u_{j,i}, i = 1,2,\cdots$ are the eigenvalues and eigenfunctions of $T_{j,d}$. An analogous identity holds for $\chi_a(T_{0,d})$. By ~\eqref{eq:double_neighborhood}-- ~\eqref{eq:single_neighborhood}, it then follows that there exists $j_0$ such that if  $j\ge j_0$, it holds that
\[
\chi_a(T_{j,d})=\Pi_{j,d}(a),
    \qquad
    \chi_a(T_{0,d})=\Pi_{0,d}(a).
\]
The continuous functional calculus (\cite{conway_course_2007-1} Theorem VIII.2.6), applied through the polynomial approximation argument given in the proof of \cref{lemma:weak-convergence-of-solutions} on the common interval $[0, \sup_{j \geq 0}\|T_{j,d}\|_{\mathrm{op}}]$, which is bounded by Step 3, and~\eqref{eq:squared_operator_convergence} therefore imply~\eqref{eq:truncated_convergence}.

Now choose a bounded continuous function
$g_a:[0,\infty)\to\mathbb R$ satisfying
\begin{equation}
    g_a(t)=0
    \quad\text{for }t\leq a-\eta_a,
    \qquad
    g_a(t)=\frac{1}{t}
    \quad\text{for }t\geq a+\eta_a.
    \label{eq:g_definition}
\end{equation}
Again, because the spectra of $T_{0,d}$ and, for all sufficiently
large $j$, of $T_{j,d}$ contain no points in the corresponding
transition region, we have
\begin{equation}
    g_a(t)
    =
    \frac{1}{t}
    \qquad
    \text{for every }
    t\in\operatorname{spec}(T_{j,d})\cap[a,\infty),
    \label{eq:g_spectra}
\end{equation}
with the same assertion holding for $j=0$.

Recall that $\bar h_j^\dagger(\cdot,d)$ is the
minimum-$\mathbb H_{X,d}$-norm solution of
\(
    \mathcal E_{j,d}h=r_{P_j}(\cdot,d).
\)
It therefore belongs to
$\mathscr{N}(\mathcal E_{j,d})^\perp
 =\mathscr{N}(T_{j,d})^\perp$
and satisfies the normal equation
\begin{equation}
    T_{j,d}\bar h_j^\dagger(\cdot,d)
    =
    \mathcal E_{j,d}^{*0}r_{P_j}(\cdot,d).
    \label{eq:normal_equation}
\end{equation}
Applying the spectral decomposition of $T_{j,d}$ to~\eqref{eq:normal_equation}, and
retaining only the eigenspaces with eigenvalues at least $a$, gives
\begin{equation}
    \Pi_{j,d}(a)\bar h_j^\dagger(\cdot,d)
    =
    g_a(T_{j,d})
    \mathcal E_{j,d}^{*0}r_{P_j}(\cdot,d).
    \label{eq:finite_dim_eq}
\end{equation}

By the Hellinger continuity of $P \mapsto r_P$ at $P_0$ over $\mc{P}_{cont}$ assumed in \cref{thrm:bounded-means-continuous},
\begin{equation}
    \|r_{P_j}(\cdot,d)-r_{P_0}(\cdot,d)\|_{\mathbb H_{Z,d}}
    \longrightarrow 0.
    \label{eq:r_convergence}
\end{equation}
Furthermore, by the continuous functional calculus (\cite{conway_course_2007-1} Theorem VIII.2.6), applied through the polynomial approximation argument given in the proof of \cref{lemma:weak-convergence-of-solutions} on the same common interval, and ~\eqref{eq:squared_operator_convergence},
\begin{equation}
    \|g_a(T_{j,d})-g_a(T_{0,d})\|_{\mathrm{op}}
    \longrightarrow 0.
    \label{eq:g_convergence}
\end{equation}
Equations ~\eqref{eq:adjoint-operator-convergence}, ~\eqref{eq:r_convergence}, and ~\eqref{eq:g_convergence}, together with the uniform
boundedness of the relevant operators and of
$r_{P_j}(\cdot,d)$, now imply
\begin{align*}
&
\left\|
    g_a(T_{j,d})
    \mathcal E_{j,d}^{*0}r_{P_j}(\cdot,d)
    -
    g_a(T_{0,d})
    \mathcal E_{0,d}^{*0}r_{P_0}(\cdot,d)
\right\|_{\mathbb H_{X,d}}
\\
&\quad\leq
\|g_a(T_{j,d})-g_a(T_{0,d})\|_{\mathrm{op}}
\,
\|\mathcal E_{j,d}^{*0}r_{P_j}(\cdot,d)\|_{\mathbb H_{X,d}}
\\
&\qquad
+
\|g_a(T_{0,d})\|_{\mathrm{op}}
\,
\|(\mathcal E_{j,d}^{*0}-\mathcal E_{0,d}^{*0})
       r_{P_j}(\cdot,d)\|_{\mathbb H_{X,d}}
\\
&\qquad
+
\|g_a(T_{0,d})\|_{\mathrm{op}}
\,
\|\mathcal E_{0,d}^{*0}
       \{r_{P_j}(\cdot,d)-r_{P_0}(\cdot,d)\}
 \|_{\mathbb H_{X,d}}
\longrightarrow 0.
\end{align*}
Using~\eqref{eq:finite_dim_eq}, we conclude that, for every fixed
$a>0$ outside $\operatorname{spec}(T_{0,d})$,
\begin{equation}
    \left\|
        \Pi_{j,d}(a)\bar h_j^\dagger(\cdot,d)
        -
        \Pi_{0,d}(a)\bar h_0^\dagger(\cdot,d)
    \right\|_{\mathbb H_{X,d}}
    \longrightarrow 0.
    \label{eq:pointwise_truncated_convergence}
\end{equation}

\textbf{Step 6: Convergence of the inner products with
$\alpha_{P_0}$.}

Because ~\cref{as:solutions-exist} holds at $P_0$,
\begin{equation}
    \alpha_{P_0}(\cdot,d)
    \in
    \mathscr{R}(\mathcal E_{0,d}^{*0})
    \subseteq
    \mathscr{N}(\mathcal E_{0,d})^\perp
    =
    \mathscr{N}(T_{0,d})^\perp.
    \label{eq:alpha0}
\end{equation}
On the other hand, for any $h \in \mathbb H_{X,d}$ it holds that $\|
       \Pi_{0,d}(a)h-\Pi[h|\mathscr{N}(T_{0,d})^{\perp}]
    \|_{\mathbb H_{X,d}}
    \longrightarrow 0
    \qquad\text{as }a\downarrow0$
through values outside $\operatorname{spec}(T_{0,d})$. Therefore,
\begin{equation}
    \|
       (I-\Pi_{0,d}(a))\alpha_{P_0}(\cdot,d)
    \|_{\mathbb H_{X,d}}
    \longrightarrow 0
    \qquad\text{as }a\downarrow0.
    \label{eq:alpha0_convergence}
\end{equation}
For every fixed $a>0$ outside $\operatorname{spec}(T_{0,d})$, ~\eqref{eq:truncated_convergence} gives
\begin{equation}
\|
   (I-\Pi_{j,d}(a))\alpha_{P_0}(\cdot,d)
 \|_{\mathbb H_{X,d}}
 \longrightarrow
\|
   (I-\Pi_{0,d}(a))\alpha_{P_0}(\cdot,d)
 \|_{\mathbb H_{X,d}}.
\label{eq:Pij_alpha_convergence}
\end{equation}

Now decompose
\begin{align*}
\left\langle
    \alpha_{P_0}(\cdot,d),
    \bar h_j^\dagger(\cdot,d)-\bar h_0^\dagger(\cdot,d)
\right\rangle_{\mathbb H_{X,d}}
&=
\left\langle
    \alpha_{P_0}(\cdot,d),
    \Pi_{j,d}(a)\bar h_j^\dagger(\cdot,d)
    -
    \Pi_{0,d}(a)\bar h_0^\dagger(\cdot,d)
\right\rangle_{\mathbb H_{X,d}}
\\
&\quad+
\left\langle
    (I-\Pi_{j,d}(a))\alpha_{P_0}(\cdot,d),
    \bar h_j^\dagger(\cdot,d)
\right\rangle_{\mathbb H_{X,d}}
\\
&\quad-
\left\langle
    (I-\Pi_{0,d}(a))\alpha_{P_0}(\cdot,d),
    \bar h_0^\dagger(\cdot,d)
\right\rangle_{\mathbb H_{X,d}}.
\end{align*}
For every fixed $a$, the first term converges to zero by ~\eqref{eq:pointwise_truncated_convergence}.
By Cauchy--Schwarz and ~\eqref{eq:fixed-space-uniform-bound}, the absolute value of the sum of the
last two terms is bounded by
\begin{align*}
    B
    \|
       (I-\Pi_{j,d}(a))\alpha_{P_0}(\cdot,d)
    \|_{\mathbb H_{X,d}}
    +
    B
    \|
       (I-\Pi_{0,d}(a))\alpha_{P_0}(\cdot,d)
    \|_{\mathbb H_{X,d}}.
\end{align*}
Taking the limit superior as $j\to\infty$ and using ~\eqref{eq:Pij_alpha_convergence} yields
\begin{equation}
\limsup_{j\to\infty}
\left|
\left\langle
    \alpha_{P_0}(\cdot,d),
    \bar h_j^\dagger(\cdot,d)-\bar h_0^\dagger(\cdot,d)
\right\rangle_{\mathbb H_{X,d}}
\right|
\leq
2B
\|
   (I-\Pi_{0,d}(a))\alpha_{P_0}(\cdot,d)
\|_{\mathbb H_{X,d}}.
\label{eq:limsup}
\end{equation}
Letting $a\downarrow0$ through values outside
$\operatorname{spec}(T_{0,d})$ and applying ~\eqref{eq:limsup}, we conclude that
\begin{equation}
\left\langle
    \alpha_{P_0}(\cdot,d),
    \bar h_j^\dagger(\cdot,d)-\bar h_0^\dagger(\cdot,d)
\right\rangle_{\mathbb H_{X,d}}
\longrightarrow 0.
\label{eq:cond_weak_limit}
\end{equation}

Because $\textsf{Supp}[D;P_0]$ is finite, summing over $d$
gives
\begin{equation}
    \left\langle
        \alpha_{P_0},
        \bar h_j^\dagger-\bar h_0^\dagger
    \right\rangle_{L^2(P_{X,D;0})}
    \longrightarrow 0,
    \label{eq:weak_limit}
\end{equation}
where
\(
    \bar h_j^\dagger(x,d)
    :=
    \bar h_j^\dagger(\cdot,d)(x).
\)

\textbf{Step 7: Convergence of the change-of-measure term.}

Let
\[
    w_j(x,d)
    :=
    \frac{p_{X,D;j}(x,d)}{p_{X,D;0}(x,d)}.
\]
The uniform density bounds imply
\begin{equation}
    \sup_j
\|w_j\|_{L^\infty(P_{X,D;0})}
    <\infty.
    \label{eq:unif_lik_ratio_bound}
\end{equation}
Marginalizing ~\eqref{eq:hellinger-implies-L1} gives
\[
    \|p_{X,D;j}-p_{X,D;0}\|_{L^1(\lambda_X\otimes\lambda_D)}
    \longrightarrow0.
\]
Consequently,
\[
    \|w_j-1\|_{L^1(P_{X,D;0})}
    =
    \|p_{X,D;j}-p_{X,D;0}\|_{L^1(\lambda_X\otimes\lambda_D)}
    \longrightarrow0.
\]
Together with the uniform $L^\infty$ bound ~\eqref{eq:unif_lik_ratio_bound}, this implies
\begin{equation}
    \|w_j-1\|_{L^2(P_{X,D;0})}
    \longrightarrow0.
    \label{eq:wj_convergence}
\end{equation}

Furthermore, since \(P_0\in\mathcal P\), \cref{as:continuous-linear-maps} implies
\(\alpha_{P_0}\in L^2(P_{X,D;0})\).  By ~\eqref{eq:unif_lik_ratio_bound}, the sequence
\(\{w_j-1\}_{j\ge1}\) is uniformly bounded in \(L^\infty(P_{X,D;0})\), and by
~\eqref{eq:wj_convergence}, \(w_j-1\to0\) in \(L^2(P_{X,D;0})\), hence in
\(P_{X,D;0}\)-probability.  Therefore
\[
    |(w_j-1)\alpha_{P_0}|^2
    \le
    \sup_j\|w_j-1\|_{L^\infty(P_{X,D;0})}^2
    \alpha_{P_0}^2 ,
\]
where the right-hand side is \(P_{X,D;0}\)-integrable.  Since
\(|w_j-1|^2\alpha_{P_0}^2\to0\) in probability, by dominated convergence in probability we have
\begin{equation}
    \|(w_j-1)\alpha_{P_0}\|_{L^2(P_{X,D;0})}\to0 .
    \label{eq:wjtimesalpha_convergence}
\end{equation}

On the other hand, by ~\eqref{eq:unif_lik_ratio_bound} and the Hellinger continuity of $P \mapsto \alpha_P$ at $P_0$ over $\mc{P}_{cont}$ assumed in \cref{thrm:bounded-means-continuous}, we have
\begin{equation}
    \|w_j(\alpha_{P_j}-\alpha_{P_0})\|_{L^2(P_{X,D;0})} \le \sup_j\|w_j\|_{L^\infty(P_{X,D;0})} \|\alpha_{P_j}-\alpha_{P_0}\|_{L^2(P_{X,D;0})}
    \longrightarrow0.
    \label{eq:alphaj_convergence}
\end{equation}

Therefore, by ~\eqref{eq:wjtimesalpha_convergence} and ~\eqref{eq:alphaj_convergence} we arrive at
\begin{align*}
\|w_j\alpha_{P_j}-\alpha_{P_0}\|_{L^2(P_{X,D;0})}
&\leq
\|w_j(\alpha_{P_j}-\alpha_{P_0})\|_{L^2(P_{X,D;0})}
+
\|(w_j-1)\alpha_{P_0}\|_{L^2(P_{X,D;0})}
\longrightarrow0.
\end{align*}
Thus,
\(
    \|w_j\alpha_{P_j}-\alpha_{P_0}\|_{L^2(P_{X,D;0})}
    \longrightarrow0.
\)

The fixed-space solutions are uniformly bounded in
$L^2(P_{X,D;0})$ by ~\eqref{eq:fixed-space-uniform-bound} and the finiteness of
$\textsf{Supp}[D;P_0]$. Hence, by Cauchy--Schwarz,
\begin{equation}
    \left\langle
        w_j\alpha_{P_j}-\alpha_{P_0},
        \bar h_j^\dagger
    \right\rangle_{L^2(P_{X,D;0})}
    \longrightarrow0.
    \label{eq:weak_convergence}
\end{equation}

\textbf{Step 8: Convergence of the parameter.}

Because $\bar h_j^\dagger$ solves the primal inverse problem under
$P_j$, the defining identities for the parameter give
\[
    \varphi(P_j)
    =
    E_{P_j}
    \left[
        \alpha_{P_j}(X,D)\bar h_j^\dagger(X,D)
    \right].
\]
Therefore,
\begin{align*}
    \varphi(P_j)-\varphi(P_0)
    &=
    \left\langle
        w_j\alpha_{P_j},
        \bar h_j^\dagger
    \right\rangle_{L^2(P_{X,D;0})}
    -
    \left\langle
        \alpha_{P_0},
        \bar h_0^\dagger
    \right\rangle_{L^2(P_{X,D;0})}
\\
    &=
    \left\langle
        w_j\alpha_{P_j}-\alpha_{P_0},
        \bar h_j^\dagger
    \right\rangle_{L^2(P_{X,D;0})}
\\
    &\quad+
    \left\langle
        \alpha_{P_0},
        \bar h_j^\dagger-\bar h_0^\dagger
    \right\rangle_{L^2(P_{X,D;0})}.
\end{align*}
The first term converges to zero by ~\eqref{eq:weak_convergence}, and the second converges
to zero by~\eqref{eq:weak_limit}. Hence
\[
    \varphi(P_j)\longrightarrow\varphi(P_0).
\]
Since the sequence $\{P_j\}$ was arbitrary, the mapping
$P\mapsto\varphi(P)$ is Hellinger continuous at $P_0$ over
$\mc{P}_{cont}$.

The proof under \cref{as:uniform-norm-bound} \cref{as:uniform-norm-bound:part:q} is obtained by interchanging the
primal and dual inverse problems, replacing $\bar h_j^\dagger$ by
the corresponding fixed-space minimum-norm dual solution, and
replacing $\alpha_{P_j}$ by $r_{P_j}$.
\end{proof}
}

\subsection{Proofs for \cref{subsec:recover-differentiability}}

\begin{proof}[Proof of \cref{thrm:injective-tangent-space}]
    We first give an informal proof sketch. Under the conditions of the theorem, in particular either injectivity of $\mc{E}_{P_0}$ together with a uniform bound on the norm of the minimum norm solution to the primal inverse problem, or a uniform bound on the norm of its $\beta$-source representer, $h_{P_t}^\dagger$ converges weakly to $h_{P_0}^\dagger$ as $t \to 0$. We show that, along with other regularity conditions assumed by the theorem, this suffices for $L^2(P_0)$ differentiability of $\mc{E}_{P_0}h_{P_t}^\dagger$ at $t=0$, even though the derivative of $h_{P_t}^\dagger$ may not exist. Furthermore, this derivative is an element of $\overline{\mathscr{R}(\mc{E}_{P_0})} = \mathscr{N}(\mc{E}_{P_0}^*)^\perp$. Consequently, differentiating both sides of the moment restriction $\mc{E}_{P_t}h_{P_t} = r_{P_t}$ gives the restriction on scores derived in \cite{severiniEfficiencyBoundsEstimating2012,chen_overidentification_2018-1}, even without assuming differentiability of the solution paths $h_{P_t}^\dagger$. We then show that under the assumptions of the second part of the theorem, $\dot{r}_{P_0}[s]$ admits an alternative expression, which facilitates verifying that the resulting tangent space is a proper subspace of $L^2_0(P_0)$.

    \textbf{Preliminaries.} Let $\{P_t:t \in [0, \delta)\} \subseteq \mc{P}_{diff}$ be QMD at $P_0$. We begin by establishing some preliminaries. First, we verify the conditions of \cref{lemma:differentiability-of-operators} and \cref{lemma:weak-convergence-of-solutions}. Both require $\{P_t:t \in [0, \delta)\} \subseteq \mc{P}$, which holds because $\mc{P}_{diff} \subseteq \mc{P}_{h_{P}^\dagger \text{ bound}} \subseteq \mc{P}$, and \cref{as:dominating-measure-support} at each $P_t$, which holds because each $P_t \in \mc{P}$. Condition \cref{lemma:weak-convergence-of-solutions:densities} holds because \cref{lemma:differentiability-of-operators} gives that $t \mapsto p_t$ is differentiable in $L^2(P_0)$ at $t=0$, whence $\|p_t - p_0\|_{L^2(P_0)} = t\left\|\frac{p_t - p_0}{t}\right\|_{L^2(P_0)} \rightarrow 0$. Condition \cref{lemma:weak-convergence-of-solutions:bounded} holds because $\{P_t:t \in [0, \delta)\} \subseteq \mc{P}_{diff} \subseteq \mc{P}_{h_{P}^\dagger \text{ bound}}$. Under the first condition of \cref{thrm:sufficient-pathwise-diff}, $\mc{E}_{P_0}$ is injective, and \cref{lemma:weak-convergence-of-solutions:injective-result} applies. Under the second condition, each $P_t$ belongs to $\mc{P}_{\beta_h \text{ bound}}$, and \cref{lemma:weak-convergence-of-solutions:beta-bound-result} applies, its conclusion of convergence in norm implying weak convergence. Condition \cref{lemma:weak-convergence-of-solutions:representers} holds because pathwise differentiability of $P \mapsto r_{P}$ at $P_0$ gives $\|r_{P_t} - r_{P_0}\|_{L^2(P_{Z, D;0})} = t\left\|\frac{r_{P_t} - r_{P_0}}{t}\right\|_{L^2(P_{Z, D;0})} \rightarrow 0$, and strong convergence implies weak convergence.

    Consequently, $h_{P_t}^\dagger \rightharpoonup h_{P_0}^\dagger$ and $\frac{\{\mc{E}_{P_t} - \mc{E}_{P_0}\}}{t}$ converges to $\dot{\mc{E}}_{P_0}[s]$ in $\|\cdot\|_{L^2(P_{X, D;0}) \rightarrow L^2(P_{Z, D;0})}$ norm, where $\dot{\mc{E}}_{P_0}[s]$ is compact by \cref{lemma:differentiability-of-operators}. Here $\dot{\mc{E}}_{P_0}[s]$ depends on the submodel $\{P_t:t \in [0, \delta)\}$ only through its score $s$, as the notation records. Throughout the proof we also write $s_{X, Z, D} \coloneqq E_{P_0}[s(O)|X = \cdot, Z = \cdot, D = \cdot]$ and $s_{Z, D} \coloneqq E_{P_0}[s(O)|Z = \cdot, D = \cdot]$.

    \textbf{Differentiating the primal inverse problem.}
    We now demonstrate that the derivative of the primal inverse problem exists and can be characterized as follows. By \cref{as:solutions-exist},
    \begin{align*}
        &\mc{E}_{P_t}h_{P_t}^\dagger = r_{P_t}\\
        \Rightarrow &\{\mc{E}_{P_t} - \mc{E}_{P_0}\}h_{P_t}^\dagger + \mc{E}_{P_0}\{h^\dagger_{P_t} - h^\dagger_{P_0}\} = r_{P_t} - r_{P_0}
    \end{align*}
    for each $t \in [0, \delta)$. Now,
    \begin{align*}
        \frac{\{\mc{E}_{P_t} - \mc{E}_{P_0}\}}{t}h_{P_t}^\dagger - \dot{\mc{E}}_{P_0}[s]h_{P_0}^\dagger =\left\{\frac{\{\mc{E}_{P_t} - \mc{E}_{P_0}\}}{t} - \dot{\mc{E}}_{P_0}[s]\right\}h_{P_t}^\dagger + \dot{\mc{E}}_{P_0}[s]\{h_{P_t}^\dagger - h_{P_0}^\dagger\}.
    \end{align*}
    Now,
    \begin{align*}
        \lim_{t\downarrow 0 }\left\{\frac{\{\mc{E}_{P_t} - \mc{E}_{P_0}\}}{t} - \dot{\mc{E}}_{P_0}[s]\right\}h_{P_t}^\dagger = 0 \text{ in } L^2(P_{Z, D;0})
    \end{align*}
     because the convergence of $\frac{\{\mc{E}_{P_t} - \mc{E}_{P_0}\}}{t}$ to $\dot{\mc{E}}_{P_0}[s]$ is in operator norm and $\sup_{t \in [0, \delta)}\|h^\dagger_{P_t}\|_{L^2(P_{X, D;0}) } < \infty$ by \cref{as:uniform-norm-bound} \cref{as:uniform-norm-bound:part:h}, which holds at each $P_t \in \mc{P}_{h_{P}^\dagger \text{ bound}}$, together with \cref{as:dominating-measure-support}. On the other hand, because $\dot{\mc{E}}_{P_0}[s]$ is compact and $h^\dagger_{P_t} \rightharpoonup h^\dagger_{P_0}$, $\dot{\mc{E}}_{P_0}[s]\{h_{P_t}^\dagger - h_{P_0}^\dagger\}$ converges in $L^2(P_{Z, D;0})$ to 0. Consequently,
     \begin{align*}
         \lim_{t\downarrow 0} \frac{\{\mc{E}_{P_t} - \mc{E}_{P_0}\}}{t}h_{P_t}^\dagger = \dot{\mc{E}}_{P_0}[s]h_{P_0}^\dagger \text{ in }L^2(P_{Z, D;0}).
     \end{align*}

     By the assumption that $P\mapsto r_{P}$ is pathwise differentiable at $P_0$,
     \begin{align*}
         \lim_{t\downarrow 0}\frac{r_{P_t} - r_{P_0}}{t} = \dot r_{P_0}[s]
     \end{align*}
    where $\dot r_{P_0}: L^2_0(P_0) \rightarrow L^2(P_{Z, D;0})$ is the local parameter of $P\mapsto r_{P}$ at $P_0$ relative to $\mc{P}_{\text{Riesz}}$, which contains $\mc{P}_{diff}$. Then,
    \begin{align*}
        \lim_{t\downarrow 0}\frac{\mc{E}_{P_0}\{h^\dagger_{P_t} - h^\dagger_{P_0}\}}{t} = \dot r_{P_0}[s] - \dot{\mc{E}}_{P_0}[s]h_{P_0}^\dagger \text{ in } L^2(P_{Z, D;0}).
    \end{align*}
    \textbf{The implied restriction on the scores.}
    We now argue that the differentiability of the primal inverse problem imposes a restriction on the scores of the submodels. Notice that $\frac{\mc{E}_{P_0}\{h^\dagger_{P_t} - h^\dagger_{P_0}\}}{t} \in \mathscr{R}(\mc{E}_{P_0})$ for each $t$, so that its limit belongs to the closure of $\mathscr{R}(\mc{E}_{P_0})$. Combined with the previous display, this gives
    \begin{align*}
        \dot r_{P_0}[s] - \dot{\mc{E}}_{P_0}[s]h_{P_0}^\dagger  \in \overline{\mathscr{R}(\mc{E}_{P_0}) } = \mathscr{N}(\mc{E}_{P_0}^*)^\perp.
    \end{align*}
    Recalling from \cref{lemma:differentiability-of-operators} that $\dot{\mc{E}}_{P_0}[s]h = E_{P_0}[h(X, D)s_{X|Z, D}(X, Z, D)|Z, D]$, we obtain that for any QMD at $P_0$ submodel $\{P_t:t\in[0, \delta)\}\subseteq \mc{P}_{diff}$, its score $s$ satisfies that
    \begin{align*}
        \dot r_{P_0}[s](Z, D) - E_{P_0}[h_{P_0}^\dagger(X, D)s_{X|Z, D}(X, Z, D)|Z, D]  \in \mathscr{N}(\mc{E}_{P_0}^*)^\perp.
    \end{align*}
    As such, $\dot{\mc{T}}(P_0, \mc{P}_{diff}) \subseteq \mc{C}_{diff}(P_0)$, which we abbreviate as $\mc{C}$ in what follows. The set $\mc{C}$ is a linear subspace, being the preimage of the subspace $\mathscr{N}(\mc{E}_{P_0}^*)^\perp$ under a map that is linear in $s$, but it need not be closed. Because $\mc{T}(P_0, \mc{P}_{diff})$ is the closed linear span of $\dot{\mc{T}}(P_0, \mc{P}_{diff})$, we arrive at
    \begin{align*}
        \mc{T}(P_0, \mc{P}_{diff}) \subseteq \overline{\mc{C}}.
    \end{align*}

    Now, under the additional hypothesis of \cref{thrm:injective-tangent-space} that ${m}(W, Z, D;\eta) \in L^2(P_0)$ for every $\eta \in \mathscr{N}(\mc{E}_{P_0}^*)$, we argue we can apply \cref{lemma:re-expressing-tangent-space} with $h_{P_0} = h_{P_0}^\dagger$, which then completes the proof. It suffices to show that $(1 + t s)p_0 \in \mc{P}_{\text{Riesz}}$ for each $s \in L^\infty(P_0) \cap L^2_0(P_0)$ for $t$ small enough. \cref{as:mutual-absolute-continuity,as:dominating-measure-support} hold for sufficiently small $t$ because $s \in L^\infty(P_0)$. \cref{as:continuous-linear-maps} then holds by \cref{lemma:l1-bounded-functional}, using \cref{as:m-m-tilde-bounded}. Hence, \cref{lemma:re-expressing-tangent-space} may be applied, completing this part of the proof.
\end{proof}

\begin{proof}[Proof of \cref{thrm:sufficient-pathwise-diff}]
    The proofs under the two conditions of \cref{thrm:sufficient-pathwise-diff} are identical once we have demonstrated that in either case, $h_{P_t}^\dagger$ converges weakly to $h_{P_0}^\dagger$.

    We first provide an informal proof sketch. Recall that we must establish that for any QMD at $P_0$ submodel the derivative $\frac{d}{dt}\varphi(P_t)$ exists, and that the local parameter is bounded and linear in the score of the submodels. We begin by showing that weak convergence of $h_{P_t}^\dagger$ to $h_{P_0}^\dagger$ is implied by the conditions of the Theorem. We then express $\{\varphi(P_t) - \varphi(P_0)\}/t$ as the sum of several inner products, each of which pairs a strongly convergent sequence with either a fixed function or the weakly convergent $h_{P_t}^\dagger$. We then show that this limit is bounded and linear in the score of the submodel, which establishes pathwise differentiability. Finally, we use the hypotheses that $\widetilde{m}(W, X, D;h^\dagger_{P_0})$ and ${m}(W, Z, D;q_{P_0})$ are square integrable to rewrite this limit as an inner product against the hypothesized gradient, which identifies it as a gradient. These last two hypotheses are used only in this final step: they bear on the form of the gradient, not on whether $\varphi$ is pathwise differentiable.

    \textbf{Preliminaries.} We begin by establishing some preliminary results regarding convergences of certain functions and operators. Let $\{P_t: t\in [0, \delta)\} \subseteq \mc{P}_{diff}$ be any one-dimensional QMD at $P_0$ submodel of $\mc{P}_{diff}$. We first note that the pathwise differentiability of the Hilbert-valued parameters $P \mapsto r_{P}$ and $P \mapsto \alpha_{P}$ and \cref{as:dominating-measure-support} implies that $r_{P_t} \rightarrow r_{P_0}$ and $\alpha_{P_t} \rightarrow \alpha_{P_0}$ in $L^2(P_0)$.

    We next record the consequences of \cref{lemma:differentiability-of-operators} that we will use. Its hypotheses hold because $\{P_t: t\in [0, \delta)\} \subseteq \mc{P}_{diff} \subseteq \mc{P}$ is QMD at $P_0$ and every $P \in \mc{P}_{diff}$ satisfies \cref{as:dominating-measure-support}. Writing $s \in L^2_0(P_0)$ for the score of $\{P_t: t\in [0, \delta)\}$, and $s_{Z|X, D}$ for the conditional score defined in that lemma, it gives that $t \mapsto p_{t}$ is differentiable in $L^2(P_0)$ at $t=0$ with derivative $sp_0$, that
    \begin{align*}
        \left.\frac{d}{dt}\frac{p_{Z, D;t}}{p_{Z, D;0}}\right\vert_{t=0} = s_{Z, D} \coloneqq E_{P_0}[s(O)|Z = \cdot, D = \cdot]
    \end{align*}
    in $L^2(P_0)$, and that $t \mapsto \mc{E}^*_{P_t}$ is differentiable at $t=0$ in the operator norm on $L^2(P_{Z, D;0}) \rightarrow L^2(P_{X,  D;0})$, with derivative $\dot{\mc{E}}^*_{P_0}[s]$ Hilbert-Schmidt and hence bounded.

    That lemma further gives that $t \mapsto p_{V;t}$ and $t \mapsto \frac{p_{V;t}}{p_{V;0}}$ are differentiable in $L^2(P_0)$ at $t=0$ for every subvector $V$ of $(W, X, Z, D)$. Because $\|p_{V;t} - p_{V;0}\|_{L^2(P_0)} = t\left\|\frac{p_{V;t} - p_{V;0}}{t}\right\|_{L^2(P_0)}$ and the norm on the right converges as $t \downarrow 0$, it follows that $p_{V;t} \rightarrow p_{V;0}$ in $L^2(P_0)$, and by the same argument that $\frac{p_{V;t}}{p_{V;0}} \rightarrow 1$ in $L^2(P_0)$. In particular, $p_t$, $p_{X, Z, D;t}$, $p_{X, D;t}$ and $p_{Z, D;t}$ converge to $p_0$, $p_{X, Z, D;0}$, $p_{X, D;0}$ and $p_{Z, D;0}$, respectively, in $L^2(P_0)$.

    Next, we argue that $h_{P_t}^\dagger \rightharpoonup h_{P_0}^\dagger$ in $L^2(P_{X, D;0})$. This is the conclusion of \cref{lemma:weak-convergence-of-solutions}, applied to the family $\{P_t: t\in [0, \delta)\}$ with $h_{P_t} = h^\dagger_{P_t}$, whose hypotheses we now verify. Condition \cref{lemma:weak-convergence-of-solutions:densities} is the convergence of $p_t$ to $p_0$ in $L^2(P_0)$ established above. Condition \cref{lemma:weak-convergence-of-solutions:bounded} holds because $\{P_t: t\in [0, \delta)\} \subseteq \mc{P}_{diff}$ and, under either condition of \cref{thrm:sufficient-pathwise-diff}, $\mc{P}_{diff} \subseteq \mc{P}_{h_{P}^\dagger \text{ bound}}$, so that \cref{as:uniform-norm-bound} \cref{as:uniform-norm-bound:part:h} provides a constant $C>0$ with $\|h^\dagger_{P_t}\|_{L^2(P_{X, D;t})} < C$ for all $t \in [0, \delta)$. Condition \cref{lemma:weak-convergence-of-solutions:representers} holds because $r_{P_t} \rightarrow r_{P_0}$ in $L^2(P_{Z, D;0})$, and strong convergence implies weak convergence. Under the first condition of \cref{thrm:sufficient-pathwise-diff}, $\mc{E}_{P_0}$ is injective, so that \cref{lemma:weak-convergence-of-solutions:injective-result} gives $h_{P_t}^\dagger \rightharpoonup h_{P_0}^\dagger$. Under the second condition, each $P_t$ belongs to $\mc{P}_{\beta_h \text{ bound}}$, so that \cref{lemma:weak-convergence-of-solutions:beta-bound-result} gives $h_{P_t}^\dagger \rightarrow h_{P_0}^\dagger$ in $L^2(P_{X, D;0})$, which in particular implies $h_{P_t}^\dagger \rightharpoonup h_{P_0}^\dagger$.

    \textbf{Decomposing $\varphi(P_t) - \varphi(P_0)$.} Having established the above convergences, we now argue that $\lim_{t \downarrow 0} \frac{\varphi(P_t) - \varphi(P_0)}{t}$ exists. We begin with the following decomposition, for $q_{P_t}$ and $q_{P_0}$ any solution to the dual inverse problem at $P_t$ and $P_0$, respectively:
    \begin{align}
        \label{eq:difference-in-param-derivative-identity}
        \nonumber \varphi(P_t) - \varphi(P_0) =& \langle q_{P_t}, r_{P_t}\rangle_{L^2(P_{Z, D;t})} - \langle q_{P_0}, r_{P_0}\rangle_{L^2(P_{Z, D;0})}\\
        \nonumber=& \left\langle q_{P_t} - q_{P_0}, r_{P_t}\right\rangle_{L^2(P_{Z, D;t})} + \left\langle q_{P_0}, \frac{p_{Z, D;t}}{p_{Z, D;0}}r_{P_t} - r_{P_0}\right\rangle_{L^2(P_{Z, D;0})}\\
        \nonumber=& \left\langle q_{P_t} - q_{P_0}, \mc{E}_{P_t} h_{P_t}^\dagger\right\rangle_{L^2(P_{Z, D;t})} + \left\langle q_{P_0}, \frac{p_{Z, D;t}}{p_{Z, D;0}}r_{P_t} - r_{P_0}\right\rangle_{L^2(P_{Z, D;0})}\\
        \nonumber=& \left\langle \mc{E}_{P_t}^*\{q_{P_t} - q_{P_0}\},  h_{P_t}^\dagger\right\rangle_{L^2(P_{X, D;t})} + \left\langle q_{P_0}, \frac{p_{Z, D;t}}{p_{Z, D;0}}r_{P_t} - r_{P_0}\right\rangle_{L^2(P_{Z, D;0})}\\
        \nonumber=& \left\langle \alpha_{P_t} - \mc{E}_{P_t}^*q_{P_0},  h_{P_t}^\dagger\right\rangle_{L^2(P_{X, D;t})} + \left\langle q_{P_0}, \frac{p_{Z, D;t}}{p_{Z, D;0}}r_{P_t} - r_{P_0}\right\rangle_{L^2(P_{Z, D;0})}\\
        \nonumber=& \left\langle \alpha_{P_t} - \alpha_{P_0} + \alpha_{P_0} - \mc{E}_{P_t}^*q_{P_0},  h_{P_t}^\dagger\right\rangle_{L^2(P_{X, D;t})} + \left\langle q_{P_0}, \frac{p_{Z, D;t}}{p_{Z, D;0}}r_{P_t} - r_{P_0}\right\rangle_{L^2(P_{Z, D;0})}\\
        =& \underset{(1)}{\left\langle \alpha_{P_t} - \alpha_{P_0}, h^\dagger_{P_t}\right\rangle_{L^2(P_{X, D;t})}}  + \underset{(2)}{\left\langle \{\mc{E}_{P_0}^*- \mc{E}_{P_t}^*\}q_{P_0},  h_{P_t}^\dagger\right\rangle_{L^2(P_{X, D;t})}} + \underset{(3)}{\left\langle q_{P_0}, \frac{p_{Z, D;t}}{p_{Z, D;0}}r_{P_t} - r_{P_0}\right\rangle_{L^2(P_{Z, D;0})}}
    \end{align}
    where the third, fifth, and seventh equalities follow from \cref{as:solutions-exist}. We now analyze each of the three terms in \eqref{eq:difference-in-param-derivative-identity} separately, when standardized by $t$.

    \textbf{Analysis of term $(1)$.} For term $(1)$ of \eqref{eq:difference-in-param-derivative-identity}, pathwise differentiability of $P \mapsto \alpha_{P}$ at $P_0$ implies that $\lim_{t\downarrow 0}\frac{\alpha_{P_t} - \alpha_{P_0}}{t} = \dot\alpha_{P_0}[s]$, where $s$ is the score of the submodel $\{P_t: t\in [0, \delta)\} $ and the convergence is in $L^2(P_{X, D;0})$. Prior arguments and \cref{as:dominating-measure-support} show that $\frac{p_{X, D;t}}{p_{X, D;0}} \rightarrow 1$ in $L^2(P_{X, D;0})$ and that this ratio is bounded uniformly over $t$. Consequently, because $h_{P_t}^\dagger \rightharpoonup h_{P_0}^\dagger$, $\frac{p_{X, D;t}}{p_{X, D;0}}h_{P_t}^\dagger \rightharpoonup h_{P_0}^\dagger$ and so
    \begin{align}
        \label{eq:derivative-derivation-1}
       \nonumber\lim_{t\downarrow 0}\left\langle \frac{\alpha_{P_t} - \alpha_{P_0}}{t}, h^\dagger_{P_t}\right\rangle_{L^2(P_{X, D;t})} =& \lim_{t\downarrow 0}\left\langle \frac{\alpha_{P_t} - \alpha_{P_0}}{t}, \frac{p_{X, D;t}}{p_{X, D;0}}h^\dagger_{P_t}\right\rangle_{L^2(P_{X, D;0})}\\
       \nonumber=&\lim_{t\downarrow 0}\left\langle \frac{\alpha_{P_t} - \alpha_{P_0}}{t} - \dot\alpha_{P_0}[s], \frac{p_{X, D;t}}{p_{X, D;0}}h^\dagger_{P_t}\right\rangle_{L^2(P_{X, D;0})} +  \lim_{t\downarrow 0}\left\langle \dot\alpha_{P_0}[s], \frac{p_{X, D;t}}{p_{X, D;0}}h^\dagger_{P_t}\right\rangle_{L^2(P_{X, D;0})}\\
       =& \left\langle \dot\alpha_{P_0}[s], h^\dagger_{P_0}\right\rangle_{L^2(P_{X, D;0})}.
    \end{align}

    \textbf{Analysis of term $(2)$.} We now turn to term $(2)$ of \eqref{eq:difference-in-param-derivative-identity}. Evaluating the operator norm derivative recorded in the preliminaries at the fixed function $q_{P_0} \in L^2(P_{Z, D;0})$ gives
    \begin{align}
        \label{eq:adjoint-operator-limit}
        \lim_{t\downarrow 0} \frac{\{\mc{E}^*_{P_0} - \mc{E}^*_{P_t}\}}{t}q_{P_0} = -\dot{\mc{E}}^*_{P_0}[s]q_{P_0} = -E_{P_0}[s_{Z|X, D}(X, Z, D)q_{P_0}(Z, D)|X = \cdot,D = \cdot],
    \end{align}
    where the limit holds in $L^2(P_{X, D;0})$ because convergence in operator norm implies convergence of the image of the fixed function $q_{P_0}$, and where the limit is itself an element of $L^2(P_{X, D;0})$ because $\dot{\mc{E}}^*_{P_0}[s]$ is bounded.

    Finally, we argue that the limit $\lim_{t\downarrow 0}\left \langle\frac{\{\mc{E}_{P_0}^*- \mc{E}_{P_t}^*\}}{t}q_{P_0},h_{P_t}^\dagger \right\rangle_{L^2(P_{X, D;t})}$ exists. Recall that $\frac{p_{X, D;t}}{p_{X, D;0}}$ converges to $1$ and is bounded uniformly over $t$. Consequently, $\frac{p_{X, D;t}}{p_{X, D;0}}h_{P_t}^\dagger \rightharpoonup h_{P_0}^\dagger$. Then,
    \begin{align}
        \label{eq:derivative-derivation-2}
        \nonumber\lim_{t\downarrow 0}\left \langle\frac{\{\mc{E}_{P_0}^*- \mc{E}_{P_t}^*\}}{t}q_{P_0},h_{P_t}^\dagger \right\rangle_{L^2(P_{X, D;t})} =&\lim_{t\downarrow 0}\left \langle\frac{\{\mc{E}_{P_0}^*- \mc{E}_{P_t}^*\}}{t}q_{P_0}, \frac{p_{X, D;t}}{p_{X, D;0}}h_{P_t}^\dagger\right\rangle_{L^2(P_{X, D;0})}\\
        \nonumber=&-E_{P_0}[s_{Z|X, D}(X, Z, D)q_{P_0}(Z, D)h^\dagger_{P_0}(X, D)]\\
        =&- E_{P_0}[s(O)\{q_{P_0}(Z, D)h^\dagger_{P_0}(X, D) - \alpha_{P_0}(X, D)h^\dagger_{P_0}(X, D)\}].
    \end{align}

    \textbf{Analysis of term $(3)$.} We now turn to term $(3)$ of \eqref{eq:difference-in-param-derivative-identity}. We first argue that the limit $\lim_{t\downarrow 0}\frac{\frac{p_{Z, D;t}}{p_{Z, D;0}}r_{P_t} - r_{P_0}}{t}$ exists:
    \begin{align*}
        \lim_{t\downarrow 0}\frac{\frac{p_{Z, D;t}}{p_{Z, D;0}}r_{P_t} - r_{P_0}}{t} =& \lim_{t\downarrow 0}\frac{\left\{\frac{p_{Z, D;t}}{p_{Z, D;0}} - 1\right\}r_{P_t}}{t} + \lim_{t\downarrow 0}\frac{r_{P_t} - r_{P_0}}{t}.
    \end{align*}
    By pathwise differentiability of $P \mapsto r_{P}$ at $P_0$,
    \begin{align*}
        \lim_{t\downarrow 0}\frac{r_{P_t} - r_{P_0}}{t} = \dot r_{P_0}[s].
    \end{align*}
    Recall from the preliminaries that
    \begin{align*}
        \lim_{t\downarrow 0}\frac{\left\{\frac{p_{Z, D;t}}{p_{Z, D;0}} - 1\right\}}{t} = s_{Z, D}
    \end{align*}
    in $L^2(P_{Z, D;0})$.
    Then, because $r_{P_t} \rightarrow r_{P_0}$, $\lim_{t\downarrow 0}\frac{\left\{\frac{p_{Z, D;t}}{p_{Z, D;0}} - 1\right\}r_{P_t}}{t} =s_{Z, D}r_{P_0}$ in $L^2(P_{Z, D;0})$ if $r_{P_t}$ is bounded in $L^\infty(P_{Z, D;0})$ uniformly over $t$, say by $M$. Indeed, writing the difference as $\left[\left\{\frac{p_{Z, D;t}}{p_{Z, D;0}} - 1\right\}/t - s_{Z, D}\right]r_{P_t} + s_{Z, D}\{r_{P_t} - r_{P_0}\}$, the first term converges to zero in $L^2(P_{Z, D;0})$ by the $L^\infty$ bound, and the second by dominated convergence along subsequences, with dominating function $4M^2s_{Z, D}^2$, because $r_{P_t} \rightarrow r_{P_0}$ in $P_{Z, D;0}$-probability. Recall that by \cref{as:solutions-exist},
    \begin{align*}
        r_{P_t} = \mc{E}_{P_t}h^\dagger_{P_t}.
    \end{align*}
    The same argument as in the proof of \cref{lemma:in-range-uniformly-bounded} combined with \cref{as:dominating-measure-support} shows that $\mc{E}_{P_t}$ is bounded uniformly over $t$ as a map from $L^2(P_{X, D;0})$ to $L^\infty(P_{Z, D;0})$. Then, because $h^\dagger_{P_t}$ is $L^2(P_{X, D;0})$ bounded uniformly over $t$ by \cref{as:uniform-norm-bound} \cref{as:uniform-norm-bound:part:h}, $r_{P_t}$ is bounded in $L^\infty(P_{Z, D;0})$ uniformly over $t$.

    The above discussion shows that
    \begin{align*}
        \lim_{t\downarrow 0}\frac{\frac{p_{Z, D;t}}{p_{Z, D;0}}r_{P_t} - r_{P_0}}{t} = s_{Z, D}r_{P_0} + \dot r_{P_0}[s].
    \end{align*}
    Because $q_{P_0} \in L^2(P_{Z, D;0})$, we then have that
    \begin{align}
        \label{eq:derivative-derivation-3}
        \lim_{t\downarrow 0}\left\langle q_{P_0}, \frac{\frac{p_{Z, D;t}}{p_{Z, D;0}}r_{P_t} - r_{P_0}}{t}\right\rangle_{L^2(P_{Z, D;0})} =& \langle q_{P_0}, s_{Z, D}r_{P_0} + \dot r_{P_0}[s]\rangle_{L^2(P_{Z, D;0})}.
    \end{align}

    \textbf{Pathwise differentiability.} Combining \eqref{eq:derivative-derivation-1}, \eqref{eq:derivative-derivation-2}, and \eqref{eq:derivative-derivation-3} gives that
    \begin{align}
        \label{eq:derivative-expression-not-simplified}
        \nonumber \left.\frac{d}{dt}\varphi(P_t)\right\vert_{t=0} =& \left\langle \dot\alpha_{P_0}[s], h^\dagger_{P_0}\right\rangle_{L^2(P_{X, D;0})} +  \langle \dot r_{P_0}[s],q_{P_0} \rangle_{L^2(P_{Z, D;0})} + \langle q_{P_0}r_{P_0}, s_{Z, D}\rangle_{L^2(P_{Z, D;0})}\\
        &- E_{P_0}[s(O)\{q_{P_0}(Z, D)h^\dagger_{P_0}(X, D) - \alpha_{P_0}(X, D)h^\dagger_{P_0}(X, D)\}].\\
        \nonumber=& \left\langle \dot\alpha_{P_0}[s], h^\dagger_{P_0}\right\rangle_{L^2(P_{X, D;0})} +  \langle \dot r_{P_0}[s],q_{P_0} \rangle_{L^2(P_{Z, D;0})} + \langle q_{P_0}r_{P_0}, s\rangle_{L^2(P_{0})} + \langle h^\dagger_{P_0}\alpha_{P_0}, s\rangle_{L^2(P_{0})} - \langle h^\dagger_{P_0}q_{P_0}, s\rangle_{L^2(P_{0})}
    \end{align}
    Note that $q_{P_0}r_{P_0} \in L^2(P_0)$ because $q_{P_0} \in L^2(P_{Z, D;0})$ and $r_{P_0} \in L^\infty(P_{Z, D;0})$ by \cref{lemma:in-range-uniformly-bounded}, and by similar reasoning $h^\dagger_{P_0}\alpha_{P_0} \in L^2(P_0)$. Furthermore, $h^\dagger_{P_0}q_{P_0} \in L^2(P_0)$ by \cref{lemma:product-in-l2} applied at $P_0$, whose hypotheses hold by \cref{as:dominating-measure-support}. Hence, each of the products appearing in \eqref{eq:derivative-expression-not-simplified} is an element of $L^2(P_0)$. Furthermore, the above expression is bounded and linear in the score $s$, and holds for arbitrary QMD at $P_0$ submodels $\{P_t: t\in [0, \delta)\} \subseteq \mc{P}_{diff}$. Consequently, $P \mapsto \varphi(P)$ is pathwise differentiable at $P_0$ in the model $\mc{P}_{diff}$.

    \textbf{Expressing the gradient.} We now identify the gradient, which is where we use the hypotheses that $\widetilde{m}(W, X, D;h^\dagger_{P_0}) \in L^2(P_0)$ and ${m}(W, Z, D;q_{P_0}) \in L^2(P_0)$. We first record that, for any fixed $g \in L^2(P_0)$,
    \begin{align*}
        \left|\frac{E_{P_t}[g(O)] - E_{P_0}[g(O)]}{t} - E_{P_0}[g(O)s(O)]\right| =& \left|\left\langle \frac{g}{p_0}, \frac{p_t - p_0}{t} - sp_0\right\rangle_{L^2(P_0)}\right|\\
        \leq& L_{bounded}\|g\|_{L^2(P_0)}\left\|\frac{p_t - p_0}{t} - sp_0\right\|_{L^2(P_0)},
    \end{align*}
    which converges to $0$ by the differentiability of $t \mapsto p_t$ in $L^2(P_0)$ established above. Hence $t \mapsto E_{P_t}[g(O)]$ is differentiable at $t=0$ with derivative $E_{P_0}[g(O)s(O)]$ for every $g \in L^2(P_0)$. Now, writing $\rho_t \coloneqq p_{X, D;t}/p_{X, D;0}$,
    \begin{align*}
        \frac{\langle \alpha_{P_t}, h_{P_0}^\dagger\rangle_{L^2(P_{X, D;t})} - \langle \alpha_{P_0}, h_{P_0}^\dagger\rangle_{L^2(P_{X, D;0})}}{t} = \left\langle \frac{\alpha_{P_t} - \alpha_{P_0}}{t}, \rho_t h_{P_0}^\dagger\right\rangle_{L^2(P_{X, D;0})} + \frac{E_{P_t}[\alpha_{P_0}(X, D)h_{P_0}^\dagger(X, D)] - E_{P_0}[\alpha_{P_0}(X, D)h_{P_0}^\dagger(X, D)]}{t}.
    \end{align*}
    The first term converges to $\langle \dot\alpha_{P_0}[s],h_{P_0}^\dagger\rangle_{L^2(P_{X, D;0})}$, because $(\alpha_{P_t} - \alpha_{P_0})/t \rightarrow \dot\alpha_{P_0}[s]$ in $L^2(P_{X, D;0})$ and $\rho_t h_{P_0}^\dagger \rightarrow h_{P_0}^\dagger$ in $L^2(P_{X, D;0})$ by dominated convergence along subsequences, as $\rho_t$ is bounded uniformly over $t$ and converges to $1$ in $L^2(P_{X, D;0})$. The second term converges to $\langle \alpha_{P_0}h_{P_0}^\dagger, s\rangle_{L^2(P_{0})}$ by the display above with $g \coloneqq \alpha_{P_0}h_{P_0}^\dagger$, which belongs to $L^2(P_0)$ because $\alpha_{P_0}$ is essentially bounded by \cref{lemma:in-range-uniformly-bounded}. On the other hand, by the definition of the Riesz representer $\alpha_{P_t}$,
    \begin{align*}
        &\left.\frac{d}{dt}\langle \alpha_{P_t}, h_{P_0}^\dagger\rangle_{L^2(P_{X, D;t})} \right\vert_{t=0} = \left.\frac{d}{dt}E_{P_t}[\widetilde{m}(W, X, D;h_{P_0}^\dagger)]\right\vert_{t=0} = E_{P_0}[\widetilde{m}(W, X, D;h_{P_0}^\dagger)s(O)],
    \end{align*}
    where the final equality uses $\widetilde{m}(W, X, D;h_{P_0}^\dagger) \in L^2(P_0)$. Equating the two expressions for the derivative gives
    \begin{align*}
        E_{P_0}[\widetilde{m}(W, X, D;h_{P_0}^\dagger)s(O)] = \langle \alpha_{P_0}h_{P_0}^\dagger, s\rangle_{L^2(P_{0})} + \langle \dot\alpha_{P_0}[s],h_{P_0}^\dagger\rangle_{L^2(P_{X, D;0})}.
    \end{align*}
    The same argument, using that $r_{P_0}$ is essentially bounded by \cref{lemma:in-range-uniformly-bounded} and that ${m}(W, Z, D;q_{P_0}) \in L^2(P_0)$, shows that
    \begin{align*}
        E_{P_0}[{m}(W, Z, D;q_{P_0})s(O)] = \langle r_{P_0}q_{P_0}, s\rangle_{L^2(P_{0})} + \langle \dot r_{P_0}[s],q_{P_0}\rangle_{L^2(P_{Z, D;0})}.
    \end{align*}
    Then, \eqref{eq:derivative-expression-not-simplified} simplifies to
    \begin{align*}
        \left.\frac{d}{dt}\varphi(P_t)\right\vert_{t=0} =& E_{P_0}[\{\widetilde{m}(W, X, D;h^\dagger_{P_0}) + {m}(W, Z, D;q_{P_0})- h^\dagger_{P_0}(X, D)q_{P_0}(Z, D)\}s(O)]\\
        =& E_{P_0}[\{\widetilde{m}(W, X, D;h^\dagger_{P_0}) + {m}(W, Z, D;q_{P_0})- h^\dagger_{P_0}(X, D)q_{P_0}(Z, D) - \varphi(P_0)\}s(O)]\\
        =&E_{P_0}[\varphi^1(W, X, Z, D;h^\dagger_{P_0}, q_{P_0}, \varphi(P_0)) s(O)].
    \end{align*}
    Hence,
    \begin{align*}
        \varphi^1(W, X, Z, D;h^\dagger_{P_0}, q_{P_0}, \varphi(P_0))
    \end{align*}
    is a gradient of $P \mapsto \varphi(P)$ at $P_0$ for any $q_{P_0} \in L^2(P_{Z, D;0})$ satisfying ${m}(W, Z, D;q_{P_0}) \in L^2(P_0)$ and solving
    \begin{align*}
        \mc{E}^*_{P_0}q = \alpha_{P_0}.
    \end{align*}
\end{proof}

\subsection{Proofs for \cref{subsec:failure-differentiability}}

\begin{proof}[Proof of \cref{thrm:mismatched-injective-not-differentiable}]
    Recall that for a parameter to be pathwise differentiable relative to the set of all one-dimensional QMD at $P_0$ submodels of $\mc{P}_{h_{P}^\dagger \text{ bound}}$, the derivative
    \begin{align*}
        \frac{d}{dt}\varphi(P_t)\vert_{t=0}
    \end{align*}
    must exist for every such submodel, and there must exist a bounded linear operator $\dot\varphi: \mc{T}(P_0, \mc{P}_{h_{P}^\dagger \text{ bound}}) \rightarrow \mathbb{R}$, called the local parameter, satisfying
    \begin{align*}
        \dot\varphi[s] = \left.\frac{d}{dt}\varphi(P_t)\right\vert_{t=0}
    \end{align*}
    for $s$ the score of $\{P_t:t\in [0, \delta)\}$.

    We show failure of pathwise differentiability by exhibiting a family of submodels $\{\{P_{t, \gamma}: t\in [0, \delta(\gamma))\}, \gamma \in (a_1, a_2)\}$ for which either the derivative $\left.\frac{d}{dt}\varphi(P_t)\right\vert_{t=0}$ does not exist for some $\gamma$, or for which there does not exist an operator satisfying $\dot\varphi[s_\gamma] = \left.\frac{d}{dt}\varphi(P_{t, \gamma})\right\vert_{t=0}$.

    Our construction relies on exhibiting a family of QMD at $P_0$ submodels $\mc{P}_{counter}$ whose scores $s_\gamma$ do not depend on $\gamma$, but for which the derivatives $\left.\frac{d}{dt}\varphi(P_t)\right\vert_{t=0}$, if they exist, do depend on $\gamma$. The construction is similar in spirit to a similar construction in the proof of \cref{thrm:not-injective-not-differentiable}. We perturb the singular system of $P_0$ to construct a new distribution $P_{t, \gamma}$, for which $h_{P_{t, \gamma}}^\dagger$ does not limit to $h_{P_0}^\dagger$, and for which the norm of the dual solution $q_{P_{t, \gamma}}$ diverges. We then demonstrate that this $P_{t, \gamma} \in \mc{P}_{counter}$ belongs to $\mc{P}_{h_{P}^\dagger \text{ bound}}$ for each $\gamma$ and for $t$ small enough. We show as claimed that the score of the submodel $\{P_{t, \gamma};t\in [0, \delta(\gamma))\} \subseteq \mc{P}_{counter}$ is not a function of $\gamma$, but that the derivative $\left.\frac{d}{dt}\varphi(P_t)\right\vert_{t=0}$ is, if it exists, thus proving the claimed failure of differentiability. As in related proofs, the majority of the argument consists of showing that $P_{t, \gamma} $ is a distribution and belongs to $\mc{P}_{{h_{P}^\dagger \text{ bound}}}$.

    \textbf{Preliminaries:} In the following arguments, we let $\lambda_W, \lambda_X, \lambda_Z, \lambda_D$ be normalized so that $\lambda_W(\textsf{Supp}[P_{W;0}]) = \lambda_X(\textsf{Supp}[P_{X;0}]) = \lambda_Z(\textsf{Supp}[P_{Z;0}]) = \lambda_D(\textsf{Supp}[P_{D;0}]) = 1$ to avoid tracking unnecessary constants, as in \cref{remark:marginal-conditional-density-bounds}.

    We now define some useful objects and record some of their properties for later use. Let $d_e \in \textsf{Supp}[P_{D;0}]$ be such that $\alpha_{P_0}(\cdot, d_e), r_{P_0}(\cdot, d_e)$ are not $P_0$-a.e. constant functions, which exists by assumption. Let
    \begin{align*}
        u_{\perp;d_e} \in \mathscr{N}(\mc{E}_{P_0;d_e})\cap L^\infty(P_{X|d_e;0}) \text{ with } \|u_{\perp;d_e}\|_{L^2(P_{X|d_e;0})} = 1.
    \end{align*}
    We let $(\sigma_{i;d_e}, u_{i;d_e}, v_{i;d_e})_{i=1}^\infty$ denote a singular system for $\mc{E}_{P_0;d_e}$, with singular values arranged in non-increasing order. We take this singular system to include only singular values greater than 0, and to satisfy \eqref{eq:mismatched-injective-not-differentiable:uniform-summability}, which is possible by assumption. Because $P_{X|d_e;0}$ is a probability measure, $\|u_{i;d_e}\|_{L^\infty(P_{X|d_e;0})} \geq \|u_{i;d_e}\|_{L^2(P_{X|d_e;0})} = 1$ for every $i$, so that the first sum in \eqref{eq:mismatched-injective-not-differentiable:uniform-summability} additionally gives $\sigma_{i;d_e}\|v_{i;d_e}\|_{L^\infty(P_{Z|d_e;0})} \rightarrow 0$ as $i \rightarrow \infty$.

    We now list some properties of this singular system. First, there are infinitely many such non-zero singular values. To see this, recall that $\mc{E}_{P_0;d_e}^*$ is injective. As such, $v_{i;d_e}$ must form an orthonormal basis for $L^2(P_{Z|d_e})$. On the other hand, $L^2(P_{Z|d_e})$ is infinite dimensional because $Z$ is continuously valued. As such, there are infinitely many non-zero singular values. On the other hand, the operator $\mc{E}_{P_0}$ is compact by \cref{as:dominating-measure-support} (see the discussion in \cite{carrasco_chapter_2007} after Assumption A.1), and so $\lim_{i \rightarrow \infty}\sigma_{i;d_e} = 0$.

    Second, $\sigma_{1;d_e}$, $u_{1;d_e}$, and $v_{1;d_e}$ may all be taken equal to 1. To see this, recall that $\mc{E}_{P_0;d_e}$ is a contraction, so that $\sigma_{i;d_e} \leq 1$. On the other hand, $\mc{E}_{P_0;d_e}^*\mc{E}_{P_0;d_e}1 = 1$ and $\mc{E}_{P_0;d_e}\mc{E}_{P_0;d_e}^*1 = 1$, so it follows that $\mc{E}_{P_0;d_e}$ has 1 as its largest singular value and $1$ is both a left and right singular function associated with this singular value. We order our singular system so $1$ is the first singular value and vector. This singular value is simple: $\|\mc{E}_{P_0;d_e}u\|_{L^2(P_{Z|d_e;0})} = \|u\|_{L^2(P_{X|d_e;0})}$ forces $u(X) = g(Z)$ $P_0$-a.s. given $D = d_e$ for some $g$, which forces $u$ to be a.e. constant because $p_{X, Z|d_e;0}$ is bounded below on a rectangle by \cref{as:dominating-measure-support}. The above choice is therefore unique up to a simultaneous sign flip of $u_{1;d_e}$ and $v_{1;d_e}$.

    Finally, we let $i_{\alpha}, i_r > 1$  be the first indices greater than 1 for which
    \begin{align}
        \label{eq:mismatched-injective-not-differentiable:i_alpha-not-orthog}
        \langle \alpha_{P_0}(\cdot, d_e), u_{i_\alpha;d_e}\rangle_{L^2(P_{X|d_e;0})} \not=&0,&  \langle r_{P_0}(\cdot, d_e), v_{i_r;d_e}\rangle_{L^2(P_{Z|d_e;0})} \not=&0.
    \end{align}
    Such indices exist because neither $\alpha_{P_0}(\cdot, d_e)$ nor $r_{P_0}(\cdot;d_e)$ is $P_0$-a.e. constant, and \cref{as:solutions-exist} holding at $P_0$ implies that $\alpha_{P_0}(\cdot, d_e)$ and $r_{P_0}(\cdot, d_e)$ admit expansions purely in $(u_{i;d_e})_{i=1}^\infty$ and $ (v_{i;d_e})_{i=1}^\infty$, respectively.

    \textbf{Constructing $\mc{P}_{counter}$:} We now construct our family of submodels $\mc{P}_{counter}$. In the following, we take $\gamma \in (a_1, a_2)$ as fixed, and define our bounds $a_1, a_2$ later.

    We begin by defining a perturbation to the singular system of $\mc{E}_{P_0}$, which we use as a basis to construct $P_{t, \gamma}$. We begin with our perturbed singular functions. First, note that there exists $I \in \mathbb{N}$ such that
    \begin{align*}
        \langle v_{i;d_e}, r_{P_0}(\cdot,d_e) \rangle_{L^2(P_{Z|d_e;0})}
    \end{align*}
    takes the same sign at $i = I$ as it does for infinitely many $i$. Taking the convention that $\text{sign}(0) = 0$, such an $I$ exists because $\text{sign}$ takes only the three values $-1, 0, 1$, so at least one of these three signs is attained for infinitely many $i$; we let $I$ be any index at which such a sign is attained. We define the function $\beta: (0, \delta(\gamma)) \rightarrow \mathbb{N}$ by
    \begin{align}
        \label{eq:mismatched-injective-not-differentiable:beta-def}
        \beta(t) \coloneqq \text{The first } i > \max\{i_{\alpha}, i_{r}\} \text{ satisfying } \sigma_{i;d_e} < t^2, \text{sign}\{\langle v_{i;d_e}, r_{P_0}(\cdot,d_e) \rangle_{L^2(P_{Z|d_e;0})}\} =  \text{sign}\{\langle v_{I;d_e}, r_{P_0}(\cdot,d_e) \rangle_{L^2(P_{Z|d_e;0})}\}.
    \end{align}
    Such an $i$ exists because $\sigma_{i;d_e}$ is non-increasing with $\lim_{i\rightarrow \infty}\sigma_{i;d_e} = 0$, so that $\sigma_{i;d_e} < t^2$ for all $i$ beyond some index, while by the choice of $I$ the sign condition holds for infinitely many $i$. Furthermore, because $\sigma_{i;d_e} > 0$ for each $i$, $\lim_{t\downarrow 0}\beta(t) = \infty$.

    Now, for each $t \in (0, \delta(\gamma))$, we define
    \begin{align}
        \label{eq:mismatched-injective-not-differentiable:rotation}
        u_{i_{\alpha};d_e}(t) \coloneqq& \cos (t) u_{i_\alpha;d_e} - \sin(t) u_{\perp;d_e}, & v_{i_r;d_e}(t) \coloneqq& \cos(\epsilon(t)) v_{i_r;d_e} - \sin(\epsilon(t))v_{\beta(t);d_e}\\
        u_{\beta(t);d_e}(t) \coloneqq& \sin (t) u_{i_\alpha;d_e} + \cos(t)  u_{\perp;d_e}, & v_{\beta(t);d_e}(t) \coloneqq& \sin(\epsilon(t)) v_{i_r;d_e} + \cos(\epsilon(t))v_{\beta(t);d_e}
    \end{align}
    where $\epsilon:(0, \delta(\gamma)) \rightarrow (0, \pi/2)$ is any function satisfying that, for $t$ small enough, both
    \begin{align}
        \label{eq:mismatched-injective-not-differentiable:epsilon-def}
        \epsilon(t) \leq \sigma_{\beta(t);d_e}, \quad \left|\langle r_{P_0}(\cdot, d_e), v_{\beta(t);d_e}(t)\rangle_{L^2(P_{Z|d_e;0})}\right| \in  (0, \sigma_{\beta(t);d_e}),
    \end{align}
    and further satisfying that if $\langle r_{P_0}(\cdot, d_e), v_{I;d_e}\rangle_{L^2(P_{Z|d_e;0})}\not=0$, $\text{sign}\{\langle r_{P_0}(\cdot, d_e), v_{\beta(t);d_e}(t)\rangle_{L^2(P_{Z|d_e;0})}\} = \text{sign}\{\langle r_{P_0}(\cdot, d_e), v_{I;d_e}\rangle_{L^2(P_{Z|d_e;0})}\}$. We include this condition to ensure that $\lim_{t\downarrow 0}\text{sign}\{\langle r_{P_0}(\cdot, d_e), v_{\beta(t);d_e}(t)\rangle_{L^2(P_{Z|d_e;0})}\}$ exists and is non-zero.

    To see that such $\epsilon$ exists, note that
    \begin{align*}
        \langle r_{P_0}(\cdot, d_e), \sin(\epsilon(t)) v_{i_r;d_e} + \cos(\epsilon(t))v_{\beta(t);d_e}\rangle_{L^2(P_{Z|d_e;0})} =& \sin(\epsilon(t))\langle r_{P_0}(\cdot, d_e), v_{i_r;d_e}\rangle + \cos(\epsilon(t))\langle v_{\beta(t);d_e}, r_{P_0}(\cdot, d_e)\rangle.
    \end{align*}
    The first term in the RHS of the above display is non-zero and tends to 0 as $\epsilon(t) \rightarrow 0$. The second term is bounded in absolute value by $|\langle v_{\beta(t);d_e}, r_{P_0}(\cdot, d_e)\rangle|$. By Picard's Theorem (\cite{carrasco_chapter_2007}),
    \begin{align*}
        \sum_{i=1}^\infty\frac{1}{\sigma_{i;d_e}^2}|\langle v_{i;d_e}, r_{P_0}(\cdot, d_e)\rangle|^2 < \infty,
    \end{align*}
    and consequently
    \begin{align*}
        \frac{1}{\sigma_{i;d_e}}|\langle v_{i;d_e}, r_{P_0}(\cdot, d_e)\rangle| \rightarrow 0, \text{ equivalently } |\langle v_{i;d_e}, r_{P_0}(\cdot, d_e)\rangle| = o(\sigma_{i;d_e}).
    \end{align*}
    Then, shrinking $\delta(\gamma)$ as necessary, for all  $t < \delta(\gamma)$, $|\langle v_{\beta(t);d_e}, r_{P_0}(\cdot, d_e)\rangle| < \sigma_{\beta(t);d_e}/2$. We may then take $\epsilon(t) <\frac{\sigma_{\beta(t);d_e}}{2|\langle r_{P_0}(\cdot, d_e), v_{i_r;d_e}\rangle|}$ so that
    \begin{align*}
        |\sin(\epsilon(t))\langle r_{P_0}(\cdot, d_e), v_{i_r;d_e}\rangle| < \sigma_{\beta(t); d_e}/2
    \end{align*}
    which implies
    \begin{align*}
        |\langle r_{P_0}(\cdot, d_e), v_{\beta(t);d_e}(t)\rangle_{L^2(P_{Z|d_e;0})}| < \sigma_{\beta(t);d_e}.
    \end{align*}
    Furthermore, noting that $\text{sign}\{\langle r_{P_0}(\cdot, d_e), v_{\beta(t);d_e}\rangle_{L^2(P_{Z|d_e;0})}\} = \text{sign}\{\langle r_{P_0}(\cdot, d_e), v_{I;d_e}\rangle_{L^2(P_{Z|d_e;0})}\}$, we may further shrink $\epsilon(t)$ as necessary to ensure that $|\cos(\epsilon(t))\langle r_{P_0}(\cdot, d_e), v_{\beta(t);d_e}\rangle_{L^2(P_{Z|d_e;0})}| > |\sin(\epsilon(t))\langle r_{P_0}(\cdot, d_e), v_{i_r;d_e}\rangle_{L^2(P_{Z|d_e;0})}|$ whenever $\langle r_{P_0}(\cdot, d_e), v_{I;d_e}\rangle_{L^2(P_{Z|d_e;0})}\not=0$. Whenever $\langle r_{P_0}(\cdot, d_e), v_{I;d_e}\rangle_{L^2(P_{Z|d_e;0})}=0$, $\text{sign}\{\langle r_{P_0}(\cdot, d_e), v_{\beta(t);d_e}(t)\rangle_{L^2(P_{Z|d_e;0})}\} = \text{sign}\{\langle r_{P_0}(\cdot, d_e), v_{i_r;d_e}\rangle_{L^2(P_{Z|d_e;0})}\}$. Consequently, in either case,
    \begin{align*}
        \lim_{t\downarrow 0}\text{sign}\{\langle r_{P_0}(\cdot, d_e), v_{\beta(t);d_e}(t)\rangle_{L^2(P_{Z|d_e;0})}\}
    \end{align*}
    exists and is non-zero.

    If the initial choice of function $\epsilon$ causes $\sin(\epsilon(t))\langle r_{P_0}(\cdot, d_e), v_{i_r;d_e}\rangle$ and $\cos(\epsilon(t))\langle v_{\beta(t);d_e}, r_{P_0}(\cdot, d_e)\rangle$ to exactly cancel, shrink $\epsilon(t)$ again, which guarantees this exact cancellation no longer holds by the properties of $\sin, \cos$.

    Now, define $\zeta(t):(0, \delta(\gamma)) \rightarrow \mathbb{N}$ as
    \begin{align}
        \label{eq:mismatched-injective-not-differentiable:zeta-def}
        \zeta(t) = \text{ the first } i > \beta(t) \text{ such that } \max \left(\sqrt{\sum_{j \geq i} \sigma_{j;d_e}^2}, \ \sum_{j \geq i}\sigma_{j;d_e}\|u_{j;d_e}\|_{L^\infty(P_0)}\left\{\|v_{j;d_e}\|_{L^\infty(P_0)} + \|v_{j+1;d_e}\|_{L^\infty(P_0)}\right\}, \ \|u_{\beta(t);d_e}\|_{L^\infty(P_0)}\sigma_{i;d_e}\|v_{i;d_e}\|_{L^\infty(P_0)}\right) < t^2.
    \end{align}
    To see that such a $\zeta$ exists, we check each of the three quantities in turn. The first two are tails of convergent series, and so tend to 0 as $i \rightarrow \infty$: $\sum_{j=1}^\infty\sigma_{j;d_e}^2 < \infty$ because $\mc{E}_{P_0;d_e}$ is Hilbert-Schmidt, its kernel
    \begin{align*}
        k_{0} \coloneqq \frac{p_{X, Z|d_e;0}}{p_{X|d_e;0}p_{Z|d_e;0}}
    \end{align*}
    belonging to $L^2(P_{X|d_e;0}\otimes P_{Z|d_e;0})$ because it is essentially bounded under \cref{as:dominating-measure-support}, and
    \begin{align*}
        \sum_{j=1}^\infty\sigma_{j;d_e}\|u_{j;d_e}\|_{L^\infty(P_0)}\left\{\|v_{j;d_e}\|_{L^\infty(P_0)} + \|v_{j+1;d_e}\|_{L^\infty(P_0)}\right\} < \infty
    \end{align*}
    by \eqref{eq:mismatched-injective-not-differentiable:uniform-summability}. The third quantity also tends to 0 as $i \rightarrow \infty$, because $\beta(t)$ is already determined at this point, so that $\|u_{\beta(t);d_e}\|_{L^\infty(P_0)}$ is a fixed finite constant, while $\sigma_{i;d_e}\|v_{i;d_e}\|_{L^\infty(P_0)} \rightarrow 0$ as noted in the preliminaries.

    Now we define our remaining perturbed singular vectors. For each $t \in (0, \delta(\gamma))$, define
    \begin{align*}
        u_{i;d_e}(t) \coloneqq u_{i;d_e} \text{ for }i\in \{1, \dots, \zeta(t)-1\}\setminus  \{i_{\alpha}, \beta(t)\} \quad \text{ and } v_{i;d_e}(t) \coloneqq v_{i;d_e} \text{ for }i \in \{1, 2, \dots\}\setminus  \{i_{r}, \beta(t)\}\\
    \end{align*}
    Finally, define
    \begin{align*}
        u_{\zeta(t);d_e}(t) = u_{\beta(t);d_e} \text{ for } i = \zeta(t) \text{ and }u_{i;d_e}(t) = u_{i-1;d_e} \text{ for }i > \zeta(t).
    \end{align*}
    Note that this construction is not symmetric, as $v_{i;d_e}(t) = v_{i;d_e}$ for $i \geq \zeta(t)$.

    Our choice of \eqref{eq:mismatched-injective-not-differentiable:rotation} defines the perturbed functions $u_{i_{\alpha};d_e}(t), u_{\beta(t);d_e}(t)$ as rotations of $u_{i_{\alpha};d_e}, u_{\perp;d_e}$ in the plane that they span with rotation angle $t$. Consequently, the perturbed functions are orthonormal and are orthogonal in $L^2(P_{X|d_e;0})$ to $u_{i;d_e}$ for any $i \not = i_{\alpha}$, in particular to $u_{\beta(t);d_e}$. The analogous statement holds for $v_{i_{r};d_e}(t), v_{\beta(t);d_e}(t)$, which are rotations of $v_{i_{r};d_e}, v_{\beta(t);d_e}$ in the plane that they span with rotation angle $\epsilon(t)$, and so are orthogonal in $L^2(P_{Z|d_e;0})$ to $v_{i;d_e}$ for any $i \not \in \{i_{r}, \beta(t)\}$. Our definition of the singular functions for $i \geq \zeta(t)$ ensures that the closure of the range of the bounded linear operator mapping $L^2(P_{Z|d_e;0})$ to $L^2(P_{X|d_e;0})$ with this new singular system is strictly larger than that of $\mc{E}_{P_0;d_e}^*$. The choice of $\zeta(t)$ ensures shifting the singular functions beyond $\zeta(t)$ does not alter either uniform convergence of the corresponding kernels or contribute to the score of the forthcoming submodel.

    Now, we define our perturbed singular values. Let
    \begin{align}
        \label{eq:mismatched-injective-not-differentiable:singular-value}
        \sigma_{\beta(t);d_e}(t,\gamma) \coloneqq \gamma |\langle r_{P_0}(\cdot, d_e), v_{\beta(t);d_e}(t)\rangle_{L^2(P_{Z|d_e;0})}|
    \end{align}
    and for each $t \in (0, \delta(\gamma))$, let
    \begin{align}
        \label{eq:mismatched-injective-not-differentiable:not-beta-singular-value}
        \sigma_{i;d_e}(t) \coloneqq \sigma_{i;d_e} \text{ for } i \in \{1, \dots, \zeta(t)\} \setminus \{\beta(t)\}, \quad \sigma_{i;d_e}(t) = \sigma_{i-1;d_e} \text{ for } i > \zeta(t).
    \end{align}
    To simplify notation, we sometimes suppress the dependence of $\sigma_{\beta(t);d_e}(t,\gamma) $ on $\gamma$. We emphasize that $\beta(t)$, $\epsilon(t)$ and $\zeta(t)$ do not depend on $\gamma$; only their common domain $(0, \delta(\gamma))$ does, and $\gamma$ enters the perturbed singular system solely through $\sigma_{\beta(t);d_e}(t, \gamma)$. Note that, unlike $(\sigma_{i;d_e})_{i=1}^\infty$, the perturbed singular values $(\sigma_{i;d_e}(t))_{i=1}^\infty$ need not be arranged in non-increasing order: $\sigma_{\beta(t);d_e}(t, \gamma)$ is smaller than $\sigma_{\beta(t);d_e}$, and $\sigma_{\zeta(t);d_e}$ appears at both of the indices $\zeta(t)$ and $\zeta(t)+1$. This is immaterial in what follows, as any non-negative sequence tending to zero, paired with orthonormal families, defines a compact operator.

    \textbf{Characterizing the perturbed density.}

    We now define a perturbed density using the perturbed singular functions and values. For each $t \in (0, \delta(\gamma))$, we let
    \begin{align*}
        p_{X, Z|d_e;t, \gamma} \coloneqq p_{X|d_e;0}p_{Z|d_e;0} k_{t, \gamma} \text{ where } k_{t, \gamma}\coloneqq \sum_{i=1}^\infty \sigma_{i;d_e}(t)u_{i;d_e}(t)v_{i;d_e}(t).
    \end{align*}
    We show below that the series defining $k_{t, \gamma}$ converges absolutely and uniformly, as does the corresponding series for $k_0$.
    We claim that for $t$ small enough, $p_{X, Z|d_e;t, \gamma}$ is a joint density for $X, Z$ with margins equal to $p_{X|d_e;0}$, and $p_{Z|d_e;0}$. It suffices to argue that $p_{X, Z|d_e;t, \gamma}$ marginalizes to those marginal densities, and that $p_{X, Z|d_e;t, \gamma} > 0$ almost everywhere $\lambda_X \times \lambda_Z$.

    We first show $\int p_{X, Z|d_e;t, \gamma}d\lambda_X = p_{Z|d_e;0}$ and $\int p_{X, Z|d_e;t, \gamma}d\lambda_Z = p_{X|d_e;0}$. For the former,
    \begin{align*}
        \int p_{X, Z|d_e;t, \gamma} (x, z)d\lambda_X(x)=&  p_{Z|d_e;0}(z)\int k_{t, \gamma} p_{X|d_e;0}(x)d\lambda_X(x)\\
        =&  p_{Z|d_e;0}(z)\int \sum_{i=1}^\infty \sigma_{i;d_e}(t)u_{i;d_e}(t)(x)v_{i;d_e}(t)(z)p_{X|d_e;0}(x)d\lambda_X(x)\\
        =& p_{Z|d_e;0}(z) \sum_{i=1}^\infty\sigma_{i;d_e}(t)\langle u_{i;d_e}(t)(x), 1\rangle_{L^2(P_{X|d_e;0})}v_{i;d_e}(t)(z)\\
        =& p_{Z|d_e;0}(z) \sigma_{1;d_e}(t)\langle u_{1;d_e}(t)(x), 1\rangle_{L^2(P_{X|d_e;0})}v_{1;d_e}(t)(z)\\
        =& p_{Z|d_e;0}(z)
    \end{align*}
    where we have used in the final two equalities that $\sigma_{1;d_e}(t), u_{1;d_e}(t), v_{1;d_e}(t) = 1$. A similar argument shows that $\int p_{X, Z|d_e;t, \gamma} d\lambda_Z = p_{X|d_e;0}$.

    We now argue that $p_{X, Z|d_e;t, \gamma} > 0$ almost everywhere. To do so, we prove that $p_{X, Z|d_e;t, \gamma}$ converges uniformly to $p_{X, Z|d_e;0}$ as $t \downarrow 0$. This fact will also be useful in later arguments.

    First, note that
    \begin{align*}
        p_{X, Z|d_e;0} = p_{X|d_e;0}p_{Z|d_e;0} k_{0} \text{ where } k_0 \coloneqq \sum_{i=1}^\infty \sigma_{i;d_e}u_{i;d_e}v_{i;d_e},
    \end{align*}
    which can be derived from the fact that $\mc{E}_{P_0;d_e}$ is an integral operator from $L^2(P_{X|d_e;0})$ to $L^2(P_{Z|d_e;0})$ with $L^2$ kernel
    \begin{align*}
        k_0  = \frac{p_{X, Z|d_e;0}}{p_{X|d_e;0}p_{Z|d_e;0}}.
    \end{align*}

    Then, because $p_{X, Z|d_e;t, \gamma} - p_{X, Z|d_e;0} = p_{X|d_e;0}p_{Z|d_e;0}\{k_{t, \gamma} - k_0\}$, it suffices to study the difference of the kernels,
    \begin{align}
        \label{eq:mismatched-injective-not-differentiable:density-expansion:1}
        k_{t, \gamma}  - k_{0} =& \sigma_{\beta(t);d_e}(t)u_{\beta(t);d_e}(t)v_{\beta(t);d_e}(t) -  \sigma_{\beta(t);d_e}u_{\beta(t);d_e}v_{\beta(t);d_e}\\
        \label{eq:mismatched-injective-not-differentiable:density-expansion:2}&+\sum_{i \in \{i_\alpha, i_r\}}\sigma_{i;d_e}\{u_{i;d_e}(t)v_{i;d_e}(t) - u_{i;d_e}v_{i;d_e}\}\\
        \label{eq:mismatched-injective-not-differentiable:density-expansion:3}&+ \sum_{i > \zeta(t)}\{\sigma_{i-1;d_e}u_{i-1;d_e} - \sigma_{i;d_e}u_{i;d_e}\}v_{i;d_e}\\
        \label{eq:mismatched-injective-not-differentiable:density-expansion:4}&+\sigma_{\zeta(t);d_e}\{u_{\beta(t);d_e} - u_{\zeta(t);d_e}\}v_{\zeta(t);d_e}.
    \end{align}
    We analyze each of the terms in the above sum.

    We begin with \eqref{eq:mismatched-injective-not-differentiable:density-expansion:1}.  Because the terms of the convergent series in \eqref{eq:mismatched-injective-not-differentiable:uniform-summability} tend to 0, we have that
    \begin{align*}
        \|\sigma_{\beta(t);d_e}u_{\beta(t);d_e}v_{\beta(t);d_e}\|_{L^\infty(P_{0})} \rightarrow 0 \text{ as } \beta(t)\rightarrow \infty \text{, equivalently as $t\downarrow 0$}.
    \end{align*}
    On the other hand, because $\sigma_{\beta(t);d_e}(t) = \gamma |\langle r_{P_0}(\cdot, d_e), v_{\beta(t);d_e}(t)\rangle_{L^2(P_{Z|d_e;0})}|  <\gamma  \sigma_{\beta(t);d_e}$, we have that
    \begin{align*}
        \|\sigma_{\beta(t);d_e}(t)u_{\beta(t);d_e}(t)v_{\beta(t);d_e}(t)\|_{L^\infty(P_{0})}\leq & \sigma_{\beta(t);d_e}\gamma\|u_{\beta(t);d_e}(t)v_{\beta(t);d_e}(t)\|_{L^\infty(P_{0})}\\
        \leq &\sigma_{\beta(t);d_e}\gamma\sum_{i_u \in \{i_{\alpha}, \perp\}}\sum_{i_v \in \{i_{r}, \beta(t)\}} \|u_{i_u;d_e}\|_{L^\infty(P_{X|d_e;0})}\|v_{i_v;d_e}\|_{L^\infty(P_{Z|d_e;0})},
    \end{align*}
    where we have used that $u_{\beta(t);d_e}(t)$ and $v_{\beta(t);d_e}(t)$ are the rotations \eqref{eq:mismatched-injective-not-differentiable:rotation} of $u_{i_\alpha;d_e}, u_{\perp;d_e}$ and of $v_{i_r;d_e}, v_{\beta(t);d_e}$, respectively. Each of the four resulting terms tends to 0: those involving $v_{\beta(t);d_e}$ because $\sigma_{i;d_e}\|v_{i;d_e}\|_{L^\infty(P_0)} \rightarrow 0$, and the remaining two because $\sigma_{\beta(t);d_e} \rightarrow 0$ while $\|u_{i_\alpha;d_e}\|_{L^\infty(P_0)}, \|u_{\perp;d_e}\|_{L^\infty(P_0)}$ and $\|v_{i_r;d_e}\|_{L^\infty(P_0)}$ are fixed and finite. Consequently,
    \begin{align*}
        \lim_{t\downarrow 0}\|\sigma_{\beta(t);d_e}(t)u_{\beta(t);d_e}(t)v_{\beta(t);d_e}(t)\|_{L^\infty(P_{0})} = 0.
    \end{align*}
    As such, \eqref{eq:mismatched-injective-not-differentiable:density-expansion:1} converges to 0 in $L^\infty(P_0)$ as $t \downarrow 0$.

    We now examine \eqref{eq:mismatched-injective-not-differentiable:density-expansion:2}. For $i \in \{i_\alpha, i_r\}$,
    \begin{align*}
        \|u_{i;d_e}(t)v_{i;d_e}(t) - u_{i;d_e}v_{i;d_e}\|_{L^\infty(P_0)} \leq& \|u_{i;d_e}(t)v_{i;d_e}(t) - u_{i;d_e}(t)v_{i;d_e}\|_{L^\infty(P_0)} + \|u_{i;d_e}(t)v_{i;d_e} - u_{i;d_e}v_{i;d_e}\|_{L^\infty(P_0)}\\
        \leq& \|u_{i;d_e}(t)\|_{L^\infty(P_{X|d_e;0})}\|v_{i;d_e}(t) -v_{i;d_e}\|_{L^\infty(P_{Z|d_e;0})} + \|v_{i;d_e}\|_{L^\infty(P_{Z|d_e;0})}\|u_{i;d_e}(t) - u_{i;d_e}\|_{L^\infty(P_{X|d_e;0})}.
    \end{align*}
    The following inequalities hold:
    \begin{align*}
        \|u_{i;d_e}(t)\|_{L^\infty(P_{X|d_e;0})} \leq& \max_{i \in \{i_{\alpha}, i_{r}\}}\|u_{i_{};d_e}\|_{L^\infty(P_{X|d_e;0})} + \|u_{\perp;d_e}\|_{L^\infty(P_{X|d_e;0})} < \infty\\
        \|v_{i;d_e}\|_{L^\infty(P_{Z|d_e;0})} \leq& \max_{i \in \{i_{\alpha}, i_{r}\}}\|v_{i_{};d_e}\|_{L^\infty(P_{Z|d_e;0})}< \infty
    \end{align*}
    On the other hand,
    \begin{align*}
        \|u_{i;d_e}(t) -u_{i;d_e}\|_{L^\infty(P_{X|d_e;0})} \leq& \{1 - \cos(t)\}\|u_{i;d_e}\|_{L^\infty(P_{X|d_e;0})} + \sin(t) \|u_{\perp;d_e}\|_{L^\infty(P_{X|d_e;0})} \rightarrow 0 \text{ as } t\downarrow 0
    \end{align*}
    and
    \begin{align*}
        \|v_{i;d_e}(t) -v_{i;d_e}\|_{L^\infty(P_{Z|d_e;0})} \leq& \{1 - \cos(\epsilon(t))\}\|v_{i;d_e}\|_{L^\infty(P_{Z|d_e;0})} + \sin(\epsilon(t)) \|v_{\beta(t);d_e}\|_{L^\infty(P_{Z|d_e;0})} \\
        \leq& \{1 - \cos(\epsilon(t))\}\|v_{i;d_e}\|_{L^\infty(P_{Z|d_e;0})} + \epsilon(t) \|v_{\beta(t);d_e}\|_{L^\infty(P_{Z|d_e;0})}\\
        \leq& \{1 - \cos(\epsilon(t))\}\|v_{i;d_e}\|_{L^\infty(P_{Z|d_e;0})} + \sigma_{\beta(t);d_e} \|v_{\beta(t);d_e}\|_{L^\infty(P_{Z|d_e;0})} \rightarrow 0 \text{ as }t \downarrow 0
    \end{align*}
    where in the final inequality we have used \eqref{eq:mismatched-injective-not-differentiable:epsilon-def}. Then, \eqref{eq:mismatched-injective-not-differentiable:density-expansion:2} converges in $L^{\infty}(P_0)$ to 0 as $t \downarrow 0$.

    Now, for \eqref{eq:mismatched-injective-not-differentiable:density-expansion:3},
    \begin{align*}
        &\left\|\sum_{i > \zeta(t)}\{\sigma_{i-1;d_e}u_{i-1;d_e} - \sigma_{i;d_e}u_{i;d_e}\}v_{i;d_e}\right\|_{L^\infty(P_0)}\\
        \leq& \sum_{i> \zeta(t)}\sigma_{i-1;d_e}\left\|u_{i-1;d_e}\right\|_{L^\infty(P_0)}\left\|v_{i;d_e}\right\|_{L^\infty(P_0)} + \sum_{i> \zeta(t)}\sigma_{i;d_e}\left\|u_{i;d_e}\right\|_{L^\infty(P_0)}\left\|v_{i;d_e}\right\|_{L^\infty(P_0)}\\
        \leq& \sum_{j \geq \zeta(t)}\sigma_{j;d_e}\|u_{j;d_e}\|_{L^\infty(P_0)}\left\{\|v_{j+1;d_e}\|_{L^\infty(P_0)} + \|v_{j;d_e}\|_{L^\infty(P_0)}\right\} < t^2 \rightarrow 0 \text{ as } t\downarrow 0,
    \end{align*}
    where the second inequality follows by re-indexing the first summation by $j = i -1$ and enlarging the range of the second from $i > \zeta(t)$ to $j \geq \zeta(t)$, and the final bound holds by our choice of $\zeta(t)$ in \eqref{eq:mismatched-injective-not-differentiable:zeta-def}.

    For \eqref{eq:mismatched-injective-not-differentiable:density-expansion:4},
    \begin{align*}
        &\sigma_{\zeta(t);d_e}\|u_{\zeta(t);d_e}(t) - u_{\zeta(t);d_e}\|_{L^\infty(P_0)}\|v_{\zeta(t);d_e}\|_{L^\infty(P_0)}\\
        \leq& \|u_{\beta(t);d_e}\|_{L^\infty(P_0)}\sigma_{\zeta(t);d_e}\|v_{\zeta(t);d_e}\|_{L^\infty(P_0)} + \sigma_{\zeta(t);d_e}\|u_{\zeta(t);d_e}\|_{L^\infty(P_0)}\|v_{\zeta(t);d_e}\|_{L^\infty(P_0)} < 2t^2 \rightarrow 0 \text{ as } t\downarrow 0,
    \end{align*}
    where both summands are bounded by $t^2$ by our choice of $\zeta(t)$ in \eqref{eq:mismatched-injective-not-differentiable:zeta-def}: the first is the third quantity appearing there, evaluated at $i = \zeta(t)$, and the second is a single term of the second.

    Combining the above arguments, $\|k_{t, \gamma} - k_0\|_{L^\infty(P_0)} \rightarrow 0$ as $t \downarrow 0$. Because $p_{X|d_e;0}$ and $p_{Z|d_e;0}$ are essentially bounded above by $L_{bounded}^2$ by \cref{remark:marginal-conditional-density-bounds}, we have shown that
    \begin{align}
        \label{eq:mismatched-injective-not-differentiable:conditional-density-convergence}
        \|p_{X, Z|d_e;t, \gamma} - p_{X, Z|d_e;0}\|_{L^\infty(P_0)} \rightarrow 0 \text{ as }t\downarrow 0 \text{ for any fixed }\gamma.
    \end{align}
    Application of the $L^\infty(P_0)$ triangle inequality combined with $\essinf p_{X, Z|d_e;0} > L_{bounded}^{-2}$ shows that $\essinf p_{X, Z|d_e;t, \gamma}  > L_{bounded}^{-2}$ for $t$ small enough, so that $p_{X, Z|d_e;t, \gamma}$ is a density function. We note that $L^\infty(P_0)$ and $L^\infty(\lambda)$ are the same on the support of $X, Z$ under $P_0$, which we recall is rectangular, and note that $p_{X, Z|d_e;t, \gamma} = k_{t, \gamma}p_{X|d_e;0}p_{Z|d_e;0}$ vanishes off this support.

    \textbf{Defining the distribution $P_{t, \gamma}$.}

    Having defined a perturbed conditional density for $(X, Z)|D=d_e$, we now finish our construction of $P_{t, \gamma}$ by defining the full joint density $p_{t, \gamma}$. For $d \not=d_e$, we let $p_{X, Z|d;t, \gamma} = p_{X, Z|d;0}$ and for all $d$ we let $P_{t, \gamma}(D=d) = P_{0}(D=d)$. Then, as in the proof of \cref{lemma:constant-riesz-representers}, we let
    \begin{align*}
        p_{t, \gamma}(w, x, z, d) \coloneqq p_{0}(w, x, z, d) + p_{W;0}(w)\{p_{X, Z, D;t, \gamma}(x, z, d) - p_{X, Z, D;0}(x, z, d)\}.
    \end{align*}
    Following the proof of that lemma and uniform convergence of the conditional densities \eqref{eq:mismatched-injective-not-differentiable:conditional-density-convergence}, we may choose $\delta(\gamma)$ sufficiently small, so that the above is a valid joint density for $O$, mutually absolutely continuous with $P_0$, with the same rectangular support, and with density of $(X, Z, D)$ under $P_{t, \gamma}$ equal to $p_{X, Z, D;t, \gamma}$, satisfying \cref{as:continuous-linear-maps}, \cref{as:mutual-absolute-continuity}, and \cref{as:dominating-measure-support}. Furthermore, $p_{W, X, D;t, \gamma} = p_{W, X, D;0}$ and $p_{W, Z, D;t, \gamma} = p_{W, Z, D;0}$, so that $\alpha_{P_{t, \gamma}} = \alpha_{P_0}$ and $r_{P_{t, \gamma}} = r_{P_0}$. Furthermore, in light of the uniform convergence of $p_{X, Z|d_e;t, \gamma}$ to $p_{X, Z|d_e;0}$ in \eqref{eq:mismatched-injective-not-differentiable:conditional-density-convergence}, and because $p_{0}$ is essentially bounded above by $L_{bounded}$ and below by $L_{bounded}^{-1}$, for $\delta(\gamma)$ sufficiently small, for all $t < \delta(\gamma)$,
    \begin{align*}
        L_{bounded}^{-1} < \essinf p_{t, \gamma}, \esssup p_{t, \gamma} < L_{bounded},
    \end{align*}
    so that $P_{t, \gamma}$ satisfies \cref{as:dominating-measure-support}. Note that $\mc{E}_{P_{t,\gamma};d_e}$ has singular system equal to $(\sigma_{i;d_e}(t), u_{i;d_e}(t), v_{i;d_e}(t))_{i=1}^\infty$.

    \textbf{Verifying $P_{t, \gamma} \in \mc{P}_{h_{P}^\dagger \text{ bound}}$:}
    Having shown $P_{t, \gamma}$ satisfies \cref{as:continuous-linear-maps}, \cref{as:mutual-absolute-continuity}, and \cref{as:dominating-measure-support}, we now argue $P_{t, \gamma}$ satisfies the remaining conditions required for inclusion in $\mc{P}_{h_{P}^\dagger \text{ bound}}$.

    We first show that \cref{as:solutions-exist} and \cref{as:beta-source} are satisfied. Because $\mc{E}_{P_{t, \gamma}}$ acts element-wise in $D$, $h^\dagger_{P_0}(\cdot, d)$ solves $\mc{E}_{P_{t, \gamma};d} h(\cdot, d) = r_{P_{t, \gamma}}(\cdot, d) = r_{P_{0}}(\cdot, d)$ for each $d \not=d_e$. At $d = d_e$, the perturbed singular functions $(u_{i;d_e}(t))_{i=1}^\infty$ consist of $(u_{i;d_e})_{i=1}^\infty$ together with $u_{\perp;d_e}$, up to the rotation \eqref{eq:mismatched-injective-not-differentiable:rotation}, so that
    \begin{align*}
        \overline{\mathscr{R}(\mc{E}_{P_{t, \gamma};d_e}^*)} \supsetneq \overline{\mathscr{R}(\mc{E}_{P_{0};d_e}^*)} \ni \alpha_{P_0}(\cdot, d_e),
    \end{align*}
    while $(v_{i;d_e}(t))_{i=1}^\infty$ remains an orthonormal basis of $L^2(P_{Z|d_e;0})$, so that $\overline{\mathscr{R}(\mc{E}_{P_{t, \gamma};d_e})} = L^2(P_{Z|d_e;0}) \ni r_{P_0}(\cdot, d_e)$. It therefore remains to verify the corresponding summability conditions. Then, \cref{as:beta-source} (and \cref{as:solutions-exist}) holds if the magnitudes of $h_{P_{t, \gamma}}^\dagger$, $q_{P_{t, \gamma}}^\dagger$ are finite with respect to the norm on the RKHS's associated with $(\mc{E}_{P_{t, \gamma}}^*\mc{E}_{P_{t, \gamma}})^{\beta_h}$ and $(\mc{E}_{P_{t, \gamma}}\mc{E}_{P_{t, \gamma}}^*)^{\beta_q}$, respectively:
    \begin{align*}
        \|h_{P_{t, \gamma}}\|^2_{(\mc{E}_{P_{t, \gamma}}^*\mc{E}_{P_{t, \gamma}})^{\beta_h}} \coloneqq \sum_{i=1}^\infty \frac{1}{\sigma_{i;d_e}(t)^{2(\beta_h + 1)}}\langle v_{i;d_e}(t),r_{P_0}(\cdot, d_e) \rangle_{L^2(P_{Z|d_e;0})}^2\\
        \|q_{P_{t, \gamma}}\|^2_{(\mc{E}_{P_{t, \gamma}}\mc{E}_{P_{t, \gamma}}^*)^{\beta_q}} \coloneqq\sum_{i=1}^\infty \frac{1}{\sigma_{i;d_e}(t)^{2(\beta_q+ 1)}}\langle u_{i;d_e}(t),\alpha_{P_0}(\cdot, d_e) \rangle_{L^2(P_{X|d_e;0})}^2
    \end{align*}
    (See e.g. Definition 3.4. \cite{carrasco_chapter_2007}, noting that our definitions of $\beta_h, \beta_q$ are offset by one from those used in that work).

    We first prove $\|h_{P_{t, \gamma}}\|^2_{(\mc{E}_{P_{t, \gamma}}^*\mc{E}_{P_{t, \gamma}})^{\beta_h}} < \infty$.
    \begin{align*}
        \|h_{P_{t, \gamma}}\|^2_{(\mc{E}_{P_{t, \gamma}}^*\mc{E}_{P_{t, \gamma}})^{\beta_h}}=&\sum_{i=1}^\infty \frac{1}{\sigma_{i;d_e}(t)^{2(\beta_h + 1)}}\langle v_{i;d_e}(t),r_{P_0}(\cdot, d_e) \rangle_{L^2(P_{Z|d_e;0})}^2\\
        \leq&\frac{1}{\sigma_{\beta(t);d_e}(t)^{2(\beta_h + 1)}}\langle v_{\beta(t);d_e}(t),r_{P_0}(\cdot, d_e) \rangle_{L^2(P_{Z|d_e;0})}^2 + \sum_{i\not=\beta(t)} \frac{1}{\sigma_{i;d_e}^{2(\beta_h + 1)}}\langle v_{i;d_e}(t),r_{P_0}(\cdot, d_e) \rangle_{L^2(P_{Z|d_e;0})}^2\\
        =&\sum_{i=1}^\infty \frac{1}{\sigma_{i;d_e}^{2(\beta_h + 1)}}\langle v_{i;d_e},r_{P_0}(\cdot, d_e) \rangle_{L^2(P_{Z|d_e;0})}^2\\
        &+ \sum_{i \in \{i_\alpha, i_r, \beta(t)\}} \frac{1}{\sigma_{i;d_e}(t)^{2(\beta_h + 1)}}\langle v_{i;d_e}(t),r_{P_0}(\cdot, d_e) \rangle_{L^2(P_{Z|d_e;0})}^2 -  \frac{1}{\sigma_{i;d_e}^{2(\beta_h + 1)}}\langle v_{i;d_e},r_{P_0}(\cdot, d_e) \rangle_{L^2(P_{Z|d_e;0})}^2
    \end{align*}
    where in the inequality we have used that $\sigma_{i;d_e}(t) \geq \sigma_{i;d_e}$ for each $i\not= \beta(t)$ by \eqref{eq:mismatched-injective-not-differentiable:not-beta-singular-value}.
    The first term on the RHS of the above display is finite because $P_{0}$ satisfies \cref{as:beta-source}. The second term is bounded because it is a finite sum of finite elements.

    We now show $\|q_{P_{t, \gamma}}\|^2_{(\mc{E}_{P_{t, \gamma}}\mc{E}_{P_{t, \gamma}}^*)^{\beta_q}} < \infty$:
     \begin{align*}
         \|q_{P_{t, \gamma}}\|^2_{(\mc{E}_{P_{t, \gamma}}\mc{E}_{P_{t, \gamma}}^*)^{\beta_q}} =&\sum_{i=1}^\infty \frac{1}{\sigma_{i;d_e}(t)^{2(\beta_q + 1)}}\langle u_{i;d_e}(t),\alpha_{P_0}(\cdot, d_e) \rangle_{L^2(P_{X|d_e;0})}^2\\
         =&\sum_{i \leq \zeta(t)} \frac{1}{\sigma_{i;d_e}(t)^{2(\beta_q + 1)}}\langle u_{i;d_e}(t),\alpha_{P_0}(\cdot, d_e) \rangle_{L^2(P_{X|d_e;0})}^2\\
         &+ \sum_{i > \zeta(t)} \frac{1}{\sigma_{i;d_e}(t)^{2(\beta_q + 1)}}\langle u_{i;d_e}(t),\alpha_{P_0}(\cdot, d_e) \rangle_{L^2(P_{X|d_e;0})}^2\\
         =&\sum_{i \leq \zeta(t)} \frac{1}{\sigma_{i;d_e}(t)^{2(\beta_q + 1)}}\langle u_{i;d_e}(t),\alpha_{P_0}(\cdot, d_e) \rangle_{L^2(P_{X|d_e;0})}^2\\
         &+ \sum_{i > \zeta(t)} \frac{1}{\sigma_{i-1;d_e}^{2(\beta_q + 1)}}\langle u_{i-1;d_e},\alpha_{P_0}(\cdot, d_e) \rangle_{L^2(P_{X|d_e;0})}^2\\
         =& \sum_{i \leq \zeta(t)} \frac{1}{\sigma_{i;d_e}(t)^{2(\beta_q + 1)}}\langle u_{i;d_e}(t),\alpha_{P_0}(\cdot, d_e) \rangle_{L^2(P_{X|d_e;0})}^2\\
         &-\sum_{i < \zeta(t)} \frac{1}{\sigma_{i;d_e}^{2(\beta_q + 1)}}\langle u_{i;d_e},\alpha_{P_0}(\cdot, d_e) \rangle_{L^2(P_{X|d_e;0})}^2\\
         &+ \sum_{i =1}^\infty \frac{1}{\sigma_{i;d_e}^{2(\beta_q + 1)}}\langle u_{i;d_e},\alpha_{P_0}(\cdot, d_e) \rangle_{L^2(P_{X|d_e;0})}^2.
     \end{align*}
     As the infinite summation in the final line of the display is finite by \cref{as:beta-source}, and the finite sums are automatically finite, $\|q_{P_{t, \gamma}}\|^2_{(\mc{E}_{P_{t, \gamma}}\mc{E}_{P_{t, \gamma}}^*)^{\beta_q}}  < \infty$. Consequently, $P_{t, \gamma}$ satisfies \cref{as:beta-source}.

     We now define the bounds $(a_1, a_2)$ on $\gamma$ and demonstrate that $P_{t, \gamma}$ satisfies \cref{as:uniform-norm-bound} \cref{as:uniform-norm-bound:part:h}, i.e. that $\|h_{P_{t, \gamma}}^\dagger(\cdot, d)\|_{L^2(P_{X|d;t, \gamma})} < C$ for each $d \in \textsf{Supp}[P_{D;0}]$, for each $\gamma \in (a_1, a_2)$ and $t < \delta(\gamma)$ for $\delta(\gamma)$ small enough. For $d \not= d_e$ we have $h^\dagger_{P_{t, \gamma}}(\cdot, d) = h^\dagger_{P_0}(\cdot, d)$ and $P_{X|d;t, \gamma} = P_{X|d;0}$, so that the bound is inherited from $P_0 \in \mc{P}_{h_{P}^\dagger \text{ bound}}$, and it remains to treat $d = d_e$. Let
     \begin{align*}
         0 < C'_h \coloneqq C - \|h_{P_{0}}^\dagger(\cdot, d_e)\|_{L^2(P_{X|d_e;0})},
     \end{align*}
     where $C'_h > 0$ holds because $\|h_{P_{0}}^\dagger(\cdot, d_e)\|_{L^2(P_{X|d_e;0})}$ is strictly less than $C$. Note that
     \begin{align*}
         h^\dagger_{P_{t, \gamma}}(\cdot, d_e) - h^\dagger_{P_0}(\cdot, d_e) =&  \sum_{i=1}^\infty \frac{1}{\sigma_{i;d_e}(t)}\langle v_{i;d_e}(t),r_{P_0}(\cdot, d_e) \rangle_{L^2(P_{Z|d_e;0})}u_{i;d_e}(t) - \sum_{i=1}^\infty \frac{1}{\sigma_{i;d_e}}\langle v_{i;d_e},r_{P_0}(\cdot, d_e) \rangle_{L^2(P_{Z|d_e;0})}u_{i;d_e}\\
         =& \frac{1}{\sigma_{\beta(t);d_e}(t, \gamma)}\langle v_{\beta(t);d_e}(t), r_{P_0}(\cdot, d_e)\rangle_{L^2(P_{Z|d_e;0})}u_{\beta(t);d_e}(t) \\
         &- \frac{1}{\sigma_{\beta(t);d_e}}\langle v_{\beta(t);d_e}, r_{P_0}(\cdot, d_e)\rangle_{L^2(P_{Z|d_e;0})}u_{\beta(t);d_e}\\
         &+ \sum_{i\in \{i_\alpha, i_r\}}\frac{1}{\sigma_{i;d_e}}\langle v_{i;d_e}(t), r_{P_0}(\cdot, d_e)\rangle_{L^2(P_{Z|d_e;0})}u_{i;d_e}(t)\\
         &-\sum_{i\in \{i_\alpha, i_r\}}\frac{1}{\sigma_{i;d_e}}\langle v_{i;d_e}, r_{P_0}(\cdot, d_e)\rangle_{L^2(P_{Z|d_e;0})}u_{i;d_e}\\
         &+\sum_{i > \zeta(t)} \left\{\frac{1}{\sigma_{i-1;d_e}}u_{i-1;d_e} - \frac{1}{\sigma_{i;d_e}}u_{i;d_e}\right\}\langle v_{i;d_e}, r_{P_0}  \rangle_{L^2(P_{Z|d_e;0})}\\
         &+ \frac{1}{\sigma_{\zeta(t);d_e}} \langle v_{\zeta(t);d_e}, r_{P_0}  \rangle_{L^2(P_{Z|d_e;0})}\{u_{\beta(t);d_e} - u_{\zeta(t);d_e}\}.
     \end{align*}
    We now analyze each of the terms in the above display.

    Now,
    \begin{align*}
        \left\|\frac{1}{\sigma_{\beta(t);d_e}(t, \gamma)}\langle v_{\beta(t);d_e}(t), r_{P_0}(\cdot, d_e)\rangle_{L^2(P_{Z|d_e;0})}u_{\beta(t);d_e}(t)\right\|_{L^2(P_{X|d_e;0})} &= \frac{1}{\gamma},\\
        \lim_{t\downarrow 0}\left\|\frac{1}{\sigma_{\beta(t);d_e}}\langle v_{\beta(t);d_e}, r_{P_0}(\cdot, d_e)\rangle_{L^2(P_{Z|d_e;0})}u_{\beta(t);d_e}\right\|_{L^2(P_{X|d_e;0})} &= 0\\
        \lim_{t\downarrow 0}\left\|\sum_{i\in \{i_\alpha, i_r\}}\frac{1}{\sigma_{i;d_e}}\langle v_{i;d_e}(t), r_{P_0}(\cdot, d_e)\rangle_{L^2(P_{Z|d_e;0})}u_{i;d_e}(t)\right.&\\
        \left.-\sum_{i\in \{i_\alpha, i_r\}}\frac{1}{\sigma_{i;d_e}}\langle v_{i;d_e}, r_{P_0}(\cdot, d_e)\rangle_{L^2(P_{Z|d_e;0})}u_{i;d_e}\right\|_{L^2(P_{X|d_e;0})} &= 0
    \end{align*}
    where the first equality follows from $\sigma_{\beta(t);d_e}(t, \gamma) = |\langle v_{\beta(t);d_e}(t), r_{P_0}(\cdot, d_e)\rangle_{L^2(P_{Z|d_e;0})}|$, the second from the fact that $\beta(t) \rightarrow \infty$ as $t\rightarrow 0$ and $\frac{1}{\sigma_{i;d_e}}\langle v_{i;d_e}, r_{P_0}(\cdot, d_e)\rangle_{L^2(P_{Z|d_e;0})} \rightarrow 0$ as $i \rightarrow \infty$, and the third by the construction of $v_{i;d_e}(t), u_{i;d_e}(t)$ for $i \in \{i_\alpha, i_r\}$.

    On the other hand,
    \begin{align*}
        &\left\|\sum_{i > \zeta(t)}\left\{\frac{1}{\sigma_{i-1;d_e}}u_{i-1;d_e} - \frac{1}{\sigma_{i;d_e}}u_{i;d_e}\right\}\langle v_{i;d_e}, r_{P_0}  \rangle_{L^2(P_{Z|d_e;0})}\right\|_{L^2(P_{X|d_e;0})} \\
        \leq& \left\|\sum_{i > \zeta(t)}\frac{1}{\sigma_{i-1;d_e}}u_{i-1;d_e}\langle v_{i;d_e}, r_{P_0}  \rangle_{L^2(P_{Z|d_e;0})}\right\|_{L^2(P_{X|d_e;0})} + \left\|\sum_{i > \zeta(t)}\frac{1}{\sigma_{i;d_e}}u_{i;d_e}\langle v_{i;d_e}, r_{P_0}  \rangle_{L^2(P_{Z|d_e;0})}\right\|_{L^2(P_{X|d_e;0})}\\
        \leq& \sqrt{\sum_{i > \zeta(t)}\frac{1}{\sigma_{i-1;d_e}^2} \langle v_{i;d_e}, r_{P_0}  \rangle_{L^2(P_{Z|d_e;0})}^2} + \sqrt{\sum_{i > \zeta(t)}\frac{1}{\sigma_{i;d_e}^2} \langle v_{i;d_e}, r_{P_0}  \rangle_{L^2(P_{Z|d_e;0})}^2}\\
        \leq &2\sqrt{\sum_{i > \zeta(t)}\frac{1}{\sigma_{i;d_e}^2} \langle v_{i;d_e}, r_{P_0}  \rangle_{L^2(P_{Z|d_e;0})}^2} \rightarrow 0 \text{  as $t\downarrow 0$}.
    \end{align*}
    The first inequality follows from the triangle inequality, and the second from orthonormality of $(u_{i;d_e})_{i=1}^\infty$. The third inequality follows from $\sigma_{i-1;d_e} \geq \sigma_{i;d_e}$, as the singular values of $\mc{E}_{P_0;d_e}$ are arranged in non-increasing order. The limit follows
    in light of $\zeta(t) \rightarrow \infty$ as $t\downarrow 0$ and convergence of the infinite sum
    \begin{align*}
        \sum_{i =1}^\infty\frac{1}{\sigma_{i;d_e}} \langle v_{i;d_e}, r_{P_0}  \rangle_{L^2(P_{Z|d_e;0})}u_{i;d_e}
    \end{align*}
    to $h_{P_0}^\dagger(\cdot, d_e) \in L^2(P_{X|d_e;0})$.
    Finally,
    \begin{align*}
        \left\|\frac{1}{\sigma_{\zeta(t);d_e}} \langle v_{\zeta(t);d_e}, r_{P_0}  \rangle_{L^2(P_{Z|d_e;0})}\{u_{\beta(t);d_e} - u_{\zeta(t);d_e}\}\right\|_{L^2(P_{X|d_e;0})} \rightarrow 0
    \end{align*}
    as $t\downarrow 0$ by the same reasoning.

    Consequently,
    \begin{align*}
        \|h^\dagger_{P_{t, \gamma}}(\cdot, d_e) - h^\dagger_{P_0}(\cdot, d_e)\|_{L^2(P_{X|d_e;0})} \leq \frac{1}{\gamma} + \nu(t)
    \end{align*}
    for $\nu(t)$ some function satisfying $\lim_{t\downarrow 0}\nu(t) = 0$. Then, choosing
    \begin{align*}
        \gamma \in \left(\frac{2}{C'_h}, \frac{4}{C'_h}\right)
    \end{align*}
    gives that
    \begin{align*}
        \|h_{P_{t, \gamma}}^\dagger(\cdot, d_e)\|_{L^2(P_{X|d_e;t, \gamma})} - C =& \|h_{P_{t, \gamma}}^\dagger(\cdot, d_e)\|_{L^2(P_{X|d_e;0})} - C\\
        \leq& \|h_{P_{t, \gamma}}^\dagger(\cdot, d_e) - h_{P_{0}}^\dagger(\cdot, d_e)\|_{L^2(P_{X|d_e;0})} + \|h_{P_{0}}^\dagger(\cdot, d_e)\|_{L^2(P_{X|d_e;0})} - C\\
        =& \|h_{P_{t, \gamma}}^\dagger(\cdot, d_e) - h_{P_{0}}^\dagger(\cdot, d_e)\|_{L^2(P_{X|d_e;0})} -C'_h\\
        <& -C'_h/2 + \nu(t) < 0 \text{ for $t$ small enough.}
    \end{align*}
   Hence, shrinking $\delta(\gamma)$ as necessary, \cref{as:uniform-norm-bound} \cref{as:uniform-norm-bound:part:h} is satisfied for all $t < \delta(\gamma)$.

    Combining the above arguments, we have shown that for all $t < \delta(\gamma)$, $\gamma \in (2/C_{h}', 4/C_{h}')$,  $P_{t, \gamma} \in \mc{P}_{h_{P}^\dagger \text{ bound}}$.

    \textbf{The score of $\{P_{t, \gamma}:t \in [0, \delta(\gamma))\}$.}

    We now argue that $\{P_{t, \gamma}:t \in [0, \delta(\gamma))\}$ is QMD at $P_0$ for each $\gamma$, and compute its score. We first argue that the map $t \mapsto p_{t, \gamma}$ is differentiable at $t=0$ in $L^2(P_0)$, then proceed by continuous differentiability of the map $x \mapsto \sqrt{x}$ in $L^2(P_0)$ for $x$ uniformly bounded away from zero to show Hellinger differentiability.

    First, notice we can express the difference quotient as
    \begin{align*}
        \frac{1}{t}\{p_{t, \gamma}(w, x, z, d) - p_{0}(w, x, z, d)\} =& \frac{1}{t}p_{W;0}(w)\{p_{X, Z|D;t, \gamma}(x, z|d) - p_{X, Z|D;0}(x, z|d)\}p_{D;0}(d)\\
        =&\frac{1}{t}I(d=d_e)p_{W;0}(w)\{p_{X, Z|d_e;t, \gamma}(x, z) - p_{X, Z|d_e;0}(x, z)\}p_{D;0}(d)\\
        =&\frac{1}{t}I(d=d_e)p_{W;0}(w)p_{X|d_e;0}(x)p_{Z|d_e;0}(z)p_{D;0}(d)\{k_{t, \gamma}(x, z) - k_{0}(x, z)\}.
    \end{align*}
    The multiplier $I(d=d_e)p_{W;0}(w)p_{X|d_e;0}(x)p_{Z|d_e;0}(z)p_{D;0}(d)$ is essentially bounded, so it suffices to study the difference quotient of the kernels.
    Then, from \eqref{eq:mismatched-injective-not-differentiable:density-expansion:1}, it suffices to argue that each of the terms in the right hand side of
    \begin{align}
        \label{eq:mismatched-injective-not-differentiable:scaled-diff-one}
        \frac{1}{t}\{k_{t, \gamma}(x, z) - k_{0}(x, z)\} =&\frac{1}{t}\left\{\sigma_{\beta(t);d_e}(t)u_{\beta(t);d_e}(t)v_{\beta(t);d_e}(t) -  \sigma_{\beta(t);d_e}u_{\beta(t);d_e}v_{\beta(t);d_e}\right\}\\
        \label{eq:mismatched-injective-not-differentiable:scaled-diff-two}
        &+\frac{1}{t}\sum_{i \in \{i_\alpha, i_r\}}\sigma_{i;d_e}\{u_{i;d_e}(t)v_{i;d_e}(t) - u_{i;d_e}v_{i;d_e}\}\\
        \label{eq:mismatched-injective-not-differentiable:scaled-diff-three}
        &+\frac{1}{t}\sum_{i > \zeta(t)}\{\sigma_{i-1;d_e}u_{i-1;d_e} - \sigma_{i;d_e}u_{i;d_e}\}v_{i;d_e}\\
        \label{eq:mismatched-injective-not-differentiable:scaled-diff-four}&+\frac{1}{t}\sigma_{\zeta(t);d_e}\{u_{\beta(t);d_e} - u_{\zeta(t);d_e}\}v_{\zeta(t);d_e}
    \end{align}
    have $L^2(P_{X, Z|d_e;0})$ limits as $t\downarrow 0$.

    For \eqref{eq:mismatched-injective-not-differentiable:scaled-diff-one},
    \begin{align*}
        &\frac{1}{t}\left\|\sigma_{\beta(t);d_e}(t)u_{\beta(t);d_e}(t)v_{\beta(t);d_e}(t) -  \sigma_{\beta(t);d_e}u_{\beta(t);d_e}v_{\beta(t);d_e}\right\|_{L^2(P_{X, Z|d_e;0})}\\
        \leq&\frac{1}{t}\left\{\sigma_{\beta(t);d_e}(t)\left\|u_{\beta(t);d_e}(t)v_{\beta(t);d_e}(t)\right\|_{L^2(P_{X, Z|d_e;0})}+ \sigma_{\beta(t);d_e}\left\|u_{\beta(t);d_e}v_{\beta(t);d_e}\right\|_{L^2(P_{X, Z|d_e;0})}\right\}\\
        \leq& (\gamma + 1)t\left\{\left\|u_{\beta(t);d_e}(t)v_{\beta(t);d_e}(t)\right\|_{L^2(P_{X, Z|d_e;0})}+\left\|u_{\beta(t);d_e}v_{\beta(t);d_e}\right\|_{L^2(P_{X, Z|d_e;0})}\right\}\\
        \leq& 2(\gamma + 1)tL_{bounded}^{7/2}\left|\textsf{Supp}[P_{D;0}]\right|^{1/2},
    \end{align*}
    where in the second inequality we have used that $\sigma_{\beta(t);d_e} < t^2$ by \eqref{eq:mismatched-injective-not-differentiable:beta-def} and that $\sigma_{\beta(t);d_e}(t, \gamma) < \gamma\sigma_{\beta(t);d_e} < \gamma t^2$ by \eqref{eq:mismatched-injective-not-differentiable:epsilon-def} and \eqref{eq:mismatched-injective-not-differentiable:singular-value}, and in the third we have used \cref{lemma:product-in-l2}. From the above display, \eqref{eq:mismatched-injective-not-differentiable:scaled-diff-one} limits to 0 as $t\downarrow 0$.

    For \eqref{eq:mismatched-injective-not-differentiable:scaled-diff-three},
    \begin{align*}
        \frac{1}{t}\left\|\sum_{i > \zeta(t)}\{\sigma_{i-1;d_e}u_{i-1;d_e} - \sigma_{i;d_e}u_{i;d_e}\}v_{i;d_e}\right\|_{{L^2(P_{X, Z|d_e;0})}} \leq& \frac{1}{t}2L_{bounded}^{7/2}\left|\textsf{Supp}[P_{D;0}]\right|^{1/2}\sqrt{\sum_{i \geq \zeta(t)}\sigma_{i;d_e}^2}\\
        \leq& t2L_{bounded}^{7/2}\left|\textsf{Supp}[P_{D;0}]\right|^{1/2}
    \end{align*}
    where the first inequality uses that $p_{X, Z|d_e;0} \leq L_{bounded}^6p_{X|d_e;0}p_{Z|d_e;0}$, as in the proof of \cref{lemma:product-in-l2}, together with orthonormality of $\{u_{i-1;d_e}\otimes v_{i;d_e}\}_{i > \zeta(t)}$ and of $\{u_{i;d_e}\otimes v_{i;d_e}\}_{i > \zeta(t)}$ in $L^2(P_{X|d_e;0}\otimes P_{Z|d_e;0})$, so that $\|\sum_i c_iu_{i-1;d_e}v_{i;d_e}\|^2_{L^2(P_{X, Z|d_e;0})} \leq L_{bounded}^6\sum_i c_i^2$, and similarly for the second sum; and the second inequality uses that $\sqrt{\sum_{i \geq \zeta(t)}\sigma_{i;d_e}^2} < t^2$ by our choice of $\zeta(t)$ in \eqref{eq:mismatched-injective-not-differentiable:zeta-def}. From the above display, \eqref{eq:mismatched-injective-not-differentiable:scaled-diff-three} limits to 0 as $t\downarrow 0$. A similar argument shows that \eqref{eq:mismatched-injective-not-differentiable:scaled-diff-four} limits to 0 as $t\downarrow 0$.

    For \eqref{eq:mismatched-injective-not-differentiable:scaled-diff-two},
    \begin{align*}
        &\frac{1}{t}\sum_{i \in \{i_\alpha, i_r\}}\sigma_{i;d_e}\{u_{i;d_e}(t)v_{i;d_e}(t) - u_{i;d_e}v_{i;d_e}\}=\frac{1}{t}\sum_{i \in \{i_\alpha, i_r\}}\sigma_{i;d_e}\{u_{i;d_e}(t)\{v_{i;d_e}(t) - v_{i;d_e}\} + v_{i;d_e} \{u_{i;d_e}(t) - u_{i;d_e}\}\}.
    \end{align*}
    Now,
    \begin{align*}
        \|u_{i;d_e}(t)\{v_{i;d_e}(t) - v_{i;d_e}\}\|_{L^2(P_{X, Z|d_e;0})} \leq& L_{bounded}^{7/2}\left|\textsf{Supp}[P_{D;0}]\right|^{1/2}\|v_{i;d_e}(t) - v_{i;d_e}\|_{L^2(P_{Z|d_e;0})}\\
        \leq& L_{bounded}^{7/2}\left|\textsf{Supp}[P_{D;0}]\right|^{1/2}\{1-\cos(\epsilon(t)) + \sin(\epsilon(t))\}\\
        \leq& L_{bounded}^{7/2}\left|\textsf{Supp}[P_{D;0}]\right|^{1/2}2\epsilon(t)\\
        \leq& L_{bounded}^{7/2}\left|\textsf{Supp}[P_{D;0}]\right|^{1/2}2t^2
    \end{align*}
    for $i \in \{i_{\alpha}, i_r\}$. On the other hand,  for $i=i_{\alpha}$, a first order Taylor expansion of $\cos, \sin$ at $t=0$ gives
    \begin{align*}
        \|u_{i_{\alpha};d_e}(t) - u_{i_{\alpha};d_e} + t u_{\perp;d_e}\|_{L^2(P_0)}= o(t).
    \end{align*}
    Note that $u_{i_{r};d_e}(t) - u_{i_{r};d_e} = 0$ if $i_r \not=i_\alpha$. Consequently,
    \begin{align*}
        \lim_{t\downarrow 0}\frac{1}{t}\sum_{i \in \{i_\alpha, i_r\}}\sigma_{i;d_e}v_{i;d_e} \{u_{i;d_e}(t) - u_{i;d_e}\} = -\sigma_{i_{\alpha};d_e}u_{\perp;d_e}v_{i_{\alpha};d_e},
    \end{align*}
    and so \eqref{eq:mismatched-injective-not-differentiable:scaled-diff-two} has limit $-\sigma_{i_{\alpha};d_e}u_{\perp;d_e}v_{i_{\alpha};d_e}$ as $t\downarrow 0$. Hence,
    \begin{align*}
        \left.\frac{d}{dt}p_{t, \gamma}\right\vert_{t=0} =& -p_{W;0}I(\cdot=d_e)p_{X|d_e;0}p_{Z|d_e;0}\sigma_{i_{\alpha};d_e}u_{\perp;d_e}v_{i_{\alpha};d_e}p_{D;0}(d_e)
    \end{align*}
    where the derivative is in $L^2(P_0)$. On the other hand, because $p_{t, \gamma}$ is uniformly bounded away from 0, $L^2(P_0)$ differentiability of the map $t \rightarrow p_{t, \gamma}$ at $t=0$ implies $L^2(P_0)$ differentiability of the map $t \mapsto \sqrt{p_{t, \gamma}}$ at $t=0$. Because $P_{t, \gamma}$ and $P_0$ are mutually absolutely continuous, this is equivalent to the statement that $\{P_{t, \gamma}:t \in [0, \delta(\gamma))\}$ is QMD at $P_0$ for each $\gamma$, and in particular that the score
    \begin{align*}
        s_{\gamma} \coloneqq -\frac{p_{W;0}I(\cdot=d_e)p_{X|d_e;0}p_{Z|d_e;0}\sigma_{i_{\alpha};d_e}u_{\perp;d_e}v_{i_{\alpha};d_e}p_{D;0}(d_e)}{p_{0}}
    \end{align*}
    does not depend on $\gamma$.

    \textbf{The derivative of $t \mapsto \varphi(P_t)$.}

    Now, we argue that if the map $t \mapsto \varphi(P_t)$ is differentiable, it depends on $\gamma$. By our definition of the singular system for $P_{t, \gamma}$, we can expand $\varphi(P_t) - \varphi(P_0)$ as
    \begin{align*}
        \varphi(P_t) - \varphi(P_0)=& \langle h^\dagger_{P_{t, \gamma}}, \alpha_{P_{t, \gamma}}\rangle_{L^2(P_{X, D;t, \gamma})} - \langle h^\dagger_{P_{0}}, \alpha_{P_{0}}\rangle_{L^2(P_{X, D;0})}\\
        =& \langle h^\dagger_{P_{t, \gamma}} - h^\dagger_{P_{0}}, \alpha_{P_{0}}\rangle_{L^2(P_{X, D;0})}\\
        =& \sum_{d \in \textsf{Supp}[P_{D;0}]}P_0(D=d)\langle h^\dagger_{P_{t, \gamma}}(\cdot, d) - h^\dagger_{P_{0}}(\cdot, d), \alpha_{P_{0}}(\cdot, d)\rangle_{L^2(P_{X|d;0})}\\
        =&P_0(D=d_e)\langle h^\dagger_{P_{t, \gamma}}(\cdot, d_e) - h^\dagger_{P_{0}}(\cdot, d_e), \alpha_{P_{0}}(\cdot, d_e)\rangle_{L^2(P_{X|d_e;0})}\\
        =& P_0(D=d_e)\sum_{i=1}^\infty \sigma_{i;d_e}^{-1}(t)\langle u_{i;d_e}(t), \alpha_{P_{0}}(\cdot, d_e)\rangle_{L^2(P_{X|d_e;0})}\langle v_{i;d_e}(t), r_{P_{0}}(\cdot, d_e)\rangle_{L^2(P_{Z|d_e;0})} \\
        &-P_0(D=d_e)\sum_{i=1}^\infty  \sigma_{i;d_e}^{-1}\langle u_{i;d_e}, \alpha_{P_{0}}(\cdot, d_e)\rangle_{L^2(P_{X|d_e;0})}\langle v_{i;d_e}, r_{P_{0}}(\cdot, d_e)\rangle_{L^2(P_{Z|d_e;0})}\\
        =& P_0(D=d_e)\left\{\frac{1}{\sigma_{\beta(t);d_e}(t, \gamma)}\langle u_{\beta(t);d_e}(t), \alpha_{P_0}(\cdot, d_e)\rangle_{L^2(P_{X|d_e;0})}\langle v_{\beta(t);d_e}(t), r_{P_0}(\cdot, d_e)\rangle_{L^2(P_{Z|d_e;0})}\right\}\\
        &-P_0(D=d_e)\left\{\frac{1}{\sigma_{\beta(t);d_e}}\langle u_{\beta(t);d_e}, \alpha_{P_0}(\cdot, d_e)\rangle_{L^2(P_{X|d_e;0})}\langle v_{\beta(t);d_e}, r_{P_0}(\cdot, d_e)\rangle_{L^2(P_{Z|d_e;0})}\right\}\\
        &+P_0(D=d_e)\sum_{i \in \{i_{r}, i_{\alpha}\}\cup\{\zeta(t), \dots\}}\frac{1}{\sigma_{i;d_e}(t)}\{\langle u_{i;d_e}(t), \alpha_{P_0}(\cdot, d_e)\rangle_{L^2(P_{X|d_e;0})}\langle v_{i;d_e}(t), r_{P_0}(\cdot, d_e)\rangle_{L^2(P_{Z|d_e;0})}\}\\
        &-P_0(D=d_e)\sum_{i \in \{i_{r}, i_{\alpha}\}\cup\{\zeta(t), \dots\}}\frac{1}{\sigma_{i;d_e}}\{\langle u_{i;d_e}, \alpha_{P_0}(\cdot, d_e)\rangle_{L^2(P_{X|d_e;0})}\langle v_{i;d_e}, r_{P_0}(\cdot, d_e)\rangle_{L^2(P_{Z|d_e;0})}\}
    \end{align*}
    We define
    \begin{align*}
        \kappa_1(t, \gamma) \coloneqq& P_0(D=d_e)\left\{\frac{1}{\sigma_{\beta(t);d_e}(t, \gamma)}\langle u_{\beta(t);d_e}(t), \alpha_{P_0}(\cdot, d_e)\rangle_{L^2(P_{X|d_e;0})}\langle v_{\beta(t);d_e}(t), r_{P_0}(\cdot, d_e)\rangle_{L^2(P_{Z|d_e;0})}\right\}\\
        \kappa_2(t) \coloneqq&-P_0(D=d_e)\left\{\frac{1}{\sigma_{\beta(t);d_e}}\langle u_{\beta(t);d_e}, \alpha_{P_0}(\cdot, d_e)\rangle_{L^2(P_{X|d_e;0})}\langle v_{\beta(t);d_e}, r_{P_0}(\cdot, d_e)\rangle_{L^2(P_{Z|d_e;0})}\right\}\\
        &+P_0(D=d_e)\sum_{i \in \{i_{r}, i_{\alpha}\}\cup\{\zeta(t), \dots\}}\frac{1}{\sigma_{i;d_e}(t)}\{\langle u_{i;d_e}(t), \alpha_{P_0}(\cdot, d_e)\rangle_{L^2(P_{X|d_e;0})}\langle v_{i;d_e}(t), r_{P_0}(\cdot, d_e)\rangle_{L^2(P_{Z|d_e;0})}\}\\
        &-P_0(D=d_e)\sum_{i \in \{i_{r}, i_{\alpha}\}\cup\{\zeta(t), \dots\}}\frac{1}{\sigma_{i;d_e}}\{\langle u_{i;d_e}, \alpha_{P_0}(\cdot, d_e)\rangle_{L^2(P_{X|d_e;0})}\langle v_{i;d_e}, r_{P_0}(\cdot, d_e)\rangle_{L^2(P_{Z|d_e;0})}\}
    \end{align*}
    Note that $\kappa_1$ depends on $\gamma$ through $\sigma_{\beta(t);d_e}(t, \gamma)$, while $\kappa_2$ does not depend on $\gamma$. Consequently, if $\lim_{t\downarrow 0}\kappa_1(t, \gamma)/t$ exists and depends on $\gamma$, then
    \begin{align*}
        \lim_{t\downarrow 0}\frac{\varphi(P_{t, \gamma}) - \varphi(P_0)}{t}
    \end{align*}
    either does not exist, or if it does, the limit depends on $\gamma$.

    Now,
    \begin{align*}
        \frac{\kappa_{1}(t, \gamma)}{t}=&\frac{P_0(D=d_e)}{t}\left\{\frac{1}{\sigma_{\beta(t);d_e}(t, \gamma)}\langle u_{\beta(t);d_e}(t), \alpha_{P_0}(\cdot, d_e)\rangle_{L^2(P_{X|d_e;0})}\langle v_{\beta(t);d_e}(t), r_{P_0}(\cdot, d_e)\rangle_{L^2(P_{Z|d_e;0})}\right\}\\
        =& \frac{P_0(D=d_e)}{\gamma t}\langle u_{\beta(t);d_e}(t), \alpha_{P_0}(\cdot, d_e)\rangle_{L^2(P_{X|d_e;0})}\text{sign}\{\langle v_{\beta(t);d_e}(t), r_{P_0}(\cdot, d_e)\rangle_{L^2(P_{Z|d_e;0})}\}\\
        =& \frac{P_0(D=d_e)}{\gamma t}\langle \sin(t)u_{i_\alpha;d_e} + \cos(t)u_{\perp;d_e}, \alpha_{P_0}(\cdot, d_e)\rangle_{L^2(P_{X|d_e;0})}\text{sign}\{\langle v_{\beta(t);d_e}(t), r_{P_0}(\cdot, d_e)\rangle_{L^2(P_{Z|d_e;0})}\}\\
        =& \frac{\sin(t)}{\gamma t}P_0(D=d_e)\langle u_{i_\alpha;d_e}, \alpha_{P_0}(\cdot, d_e)\rangle_{L^2(P_{X|d_e;0})}\text{sign}\{\langle v_{\beta(t);d_e}(t), r_{P_0}(\cdot, d_e)\rangle_{L^2(P_{Z|d_e;0})}\}
    \end{align*}
    where in the second equality we have used that by \eqref{eq:mismatched-injective-not-differentiable:singular-value}, $\sigma_{\beta(t);d_e}(t, \gamma) = |\langle v_{\beta(t);d_e}(t), r_{P_0}(\cdot, d_e)\rangle_{L^2(P_{Z|d_e;0})}|$, the third by choice of $u_{\beta(t);d_e}(t)$ in \eqref{eq:mismatched-injective-not-differentiable:rotation}, and the fourth because $\alpha_{P_0}(\cdot, d_e) \in \mathscr{R}(\mc{E}_{P_0;d_e}^*) \subseteq \mathscr{N}(\mc{E}_{P_0;d_e})^\perp$ implies $\alpha_{P_0}(\cdot, d_e) \perp_{L^2(P_{X|d_e;0})} u_{\perp;d_e}$.

    Consequently,
    \begin{align*}
        &\lim_{t\downarrow 0}\frac{\kappa_1(t, \gamma)}{t} = \frac{1}{\gamma}P_0(D=d_e)\langle u_{i_\alpha;d_e}, \alpha_{P_0}(\cdot, d_e)\rangle_{L^2(P_{X|d_e;0})}\lim_{t\downarrow 0}\text{sign}\{\langle v_{\beta(t);d_e}(t), r_{P_0}(\cdot, d_e)\rangle_{L^2(P_{Z|d_e;0})}\}.
    \end{align*}
    The limit $\lim_{t\downarrow 0}\text{sign}\{\langle v_{\beta(t);d_e}(t), r_{P_0}(\cdot, d_e)\rangle_{L^2(P_{Z|d_e;0})}\}$ exists and is non-zero by our choice of $\beta(t), \epsilon(t)$ in \eqref{eq:mismatched-injective-not-differentiable:beta-def}, \eqref{eq:mismatched-injective-not-differentiable:epsilon-def}.
    Recalling that by definition of $i_{\alpha}$,
    $\langle u_{i_\alpha;d_e}, \alpha_{P_0}(\cdot, d_e)\rangle_{L^2(P_{X|d_e;0})} \not=0$, we have that $\lim_{t\downarrow 0}\frac{\kappa_1(t, \gamma)}{t}$ exists and depends non-trivially on $\gamma$.

    Hence, we have shown that either
    \begin{align*}
        \left.\frac{d}{dt}\varphi(P_{t, \gamma})\right\vert_{t=0} = \lim_{t\downarrow 0} \frac{\kappa_1(t, \gamma)}{t} + \lim_{t\downarrow 0}\frac{\kappa_2(t)}{t}
    \end{align*}
    exists and is a function of $\gamma$, or does not exist.

    In the case where the limit does not exist for some $\gamma \in (2/C'_h, 4/C'_h)$, say $\gamma_0$, $P \mapsto \varphi(P)$ fails to be pathwise differentiable at $P_0$ over $\mc{P}_{h_{P}^\dagger \text{ bound}}$ because there exists a QMD at $P_0$ submodel $\{P_{t, \gamma_0}:t\in [0, \delta(\gamma_0))\}$ for which $\frac{d}{dt}\varphi(P_t)\vert_{t=0}$ does not exist. In the case where the limit does exist for every $\gamma$,
    \begin{align*}
        \left.\frac{d}{dt}\varphi(P_{t, \gamma})\right\vert_{t=0}
    \end{align*}
    is a non-trivial function of $\gamma$ for each submodel $\{P_{t, \gamma}: t\in [0, \delta(\gamma))\}$. However, the score $s_\gamma$ of $\{P_{t, \gamma}: t\in [0, \delta(\gamma))\}$ does not depend on $\gamma$. Consequently, there does not exist a bounded linear functional $\dot\varphi_{P_0}:\mc{T}(P_0, \mc{P}_{counter}) \rightarrow \mathbb{R}$ satisfying
    \begin{align*}
        \dot\varphi_{P_0}[s_\gamma] = \left.\frac{d}{dt}\varphi(P_{t, \gamma})\right\vert_{t=0}
    \end{align*}
    for every QMD submodel in $\mc{P}_{counter}$. Because $\mc{P}_{counter} \subseteq \mc{P}_{h_{P}^\dagger \text{ bound}}$, any bounded linear $\dot\varphi: \mc{T}(P_0, \mc{P}_{h_{P}^\dagger \text{ bound}}) \rightarrow \mathbb{R}$ satisfying the display above for every QMD submodel of $\mc{P}_{h_{P}^\dagger \text{ bound}}$ would restrict to such a functional, so no local parameter exists over $\mc{P}_{h_{P}^\dagger \text{ bound}}$ either. Because the local parameter does not exist, $P \mapsto \varphi(P)$ is not pathwise differentiable at $P_0$. In either case, the parameter is not pathwise differentiable, proving the claimed result.
\end{proof}

\begin{proof}[Proof of \cref{thrm:not-injective-not-differentiable}]
    Our proof proceeds by constructing a family of QMD at $P_0$ submodels $\mathscr{P}_{counter} \coloneqq \{\{P_{t,\gamma}: t\in[0, \delta)\}\subseteq \mc{P}_{h_{P}^\dagger, q_{P}^\dagger \text{ bound}}, \gamma \in (a_1, a_2)\}$ for which $\frac{d}{dt}\varphi(P_{t,\gamma})\vert_{t=0}$ exists for every $\gamma$, but is not a linear functional of the score of $\{P_{t,\gamma}: t\in[0, \delta)\}$; exhibiting such a family suffices to show that $\varphi$ is not pathwise differentiable.

    Concretely, as in related proofs, we mix a singular pair of $\mc{E}_{P_0}$ along which $\alpha_{P_0}, r_{P_0}$ are supported with a pair of null vectors of $\mc{E}_{P_0}, \mc{E}_{P_0}^*$, rotating each pair into the other by angle $t$ so that both converge to their originals as $t \rightarrow 0$. We then let the singular value paired with the (perturbed) null vectors vanish as $t \rightarrow 0$ at the same order, indexed by $\gamma$, as the vectors themselves rotate back toward the original pair; matching these orders is what keeps $h^\dagger_{P_{t,\gamma}}, q^\dagger_{P_{t,\gamma}}$ bounded despite dividing by a vanishing singular value. The majority of the proof is dedicated to verifying that the resulting submodels $\{P_{t,\gamma}:t\in[0,\delta)\}$ are valid probability distributions belonging to $\mc{P}_{h_{P}^\dagger, q_{P}^\dagger \text{ bound}}$; the derivative and score computations that complete the argument appear at the end.

    In the following, we let $\lambda_W, \lambda_X, \lambda_Z, \lambda_D$ be normalized so that $\lambda_W(\textsf{Supp}[P_{W;0}]) = \lambda_X(\textsf{Supp}[P_{X;0}]) = \lambda_Z(\textsf{Supp}[P_{Z;0}]) = \lambda_D(\textsf{Supp}[P_{D;0}]) = 1$ to avoid tracking unnecessary constants, as in \cref{remark:marginal-conditional-density-bounds}.

    We now define our family of submodels. We begin with some preliminary definitions. Let ${d}_e \in \textsf{Supp}[P_{D;0}]$ such that $\alpha_{P_0}(\cdot, {d}_e), r_{P_0}(\cdot, {d}_e)$ are not $P_0$-a.e. constant functions, which exists by assumption. Let
    \begin{align*}
        u_{\perp;d_e} \in \mathscr{N}(\mc{E}_{P_0;d_e})\cap L^\infty(P_{X|d_e;0}), && v_{\perp;d_e} \in \mathscr{N}(\mc{E}_{P_0;d_e}^*) \cap L^\infty(P_{Z|d_e;0}),
    \end{align*}
    with $\|u_{\perp;d_e}\|_{L^2(P_{X|d_e;0})} = \|v_{\perp;d_e}\|_{L^2(P_{Z|d_e;0})} = 1$, which exist by assumption after rescaling.
    Let $(\sigma_{i;d_e}, u_{i;d_e}, v_{i;d_e})_{i=1}^{N}$ denote a singular system for $\mc{E}_{P_0;d_e}$, with singular values arranged in non-increasing order. We take this singular system to include only those singular values greater than 0, and note its length $N \in \mathbb{N}\cup \{\infty\}$ may be either finite or infinite. Because $\mc{E}_{P_0;d_e}$ is a contraction, we have that $\sigma_{1;d_e} \leq 1$. Furthermore, because $\mc{E}_{P_0;d_e}^*\mc{E}_{P_0;d_e}1 = 1$ and $\mc{E}_{P_0;d_e}\mc{E}_{P_0;d_e}^*1 = 1$, $\mc{E}_{P_0;d_e}$ has $1$ as a singular value and as both left and right singular functions. We order our singular system so that $v_{1;d_e} = 1, u_{1;d_e} = 1$.

    Then, define $i_{\alpha}, i_r > 1$ as the first indices greater than 1 for which
    \begin{align}
        \label{eq:not-injective-not-diff:i-r-i-alpha-def}
        \langle \alpha_{P_0}(\cdot, d_e), u_{i_\alpha;d_e}\rangle_{L^2(P_{X|d_e;0})} \not= 0, &&  \langle r_{P_0}(\cdot, d_e), v_{i_r;d_e}\rangle_{L^2(P_{Z|d_e;0})} \not= 0.
    \end{align}
    Such indices exist because neither $\alpha_{P_0}(\cdot, d_e)$ nor $r_{P_0}(\cdot, d_e)$ is $P_0$-a.e. constant, and \cref{as:solutions-exist} holding at $P_0$ implies that $\alpha_{P_0}(\cdot, d_e)$ and $r_{P_0}(\cdot, d_e)$ admit expansions purely in $(u_{i;d_e})_{i=1}^N$ and $(v_{i;d_e})_{i=1}^N$ respectively.

    Now, we construct our perturbed singular functions. Define
    \begin{align}
        \label{eq:not-injective-not-diff:rotation}
        \nonumber u_{i_{\alpha};d_e}(t) \coloneqq& u_{i_{\alpha};d_e}\cos(t)  - u_{\perp;d_e}\sin(t) & v_{i_{r};d_e}(t) \coloneqq& v_{i_{r};d_e}\cos(t)  - v_{\perp;d_e}\sin(t) \\
        u_{\perp;d_e}(t)\coloneqq & u_{i_{\alpha};d_e}\sin(t)  + u_{\perp;d_e}\cos(t) & v_{\perp;d_e}(t)\coloneqq & v_{i_{r};d_e}\sin(t)  + v_{\perp;d_e}\cos(t).
    \end{align}
    Hence, we perturb $u_{i_{\alpha};d_e}$ and $u_{\perp;d_e}$ by rotating them in the plane they span, with rotation angle $t$. Consequently, the perturbed singular functions $u_{i_{\alpha};d_e}(t)$ and $u_{\perp;d_e}(t)$ remain mutually orthogonal in $L^2(P_{X|d_e;0})$, are orthogonal to $u_{i;d_e}$ for any $i \not= i_{\alpha}$, and are orthogonal to $\mathscr{N}(\mc{E}_{P_0;d_e})\ominus \text{Span}(u_{\perp;d_e})$. The same holds true for $v_{i_{r};d_e}(t), v_{\perp;d_e}(t)$.

    We now define the singular value we pair with the perturbed null singular functions. Our choice of constants will ensure that for $t$ small enough, and for any choice of $\gamma$, $P_{t, \gamma}$ remains a valid distribution and belongs to $\mc{P}_{h_{P}^\dagger, q_{P}^\dagger \text{ bound}}$. In particular, they are chosen to ensure that $\|h_{P_{t, \gamma}}^\dagger\|_{L^2(P_{X|d_e;t, \gamma})}, \|q_{P_{t, \gamma}}^\dagger\|_{L^2(P_{Z|d_e;t, \gamma})} < C$.

    Let
    \begin{align*}
        C'_{h, q}  \in \left(0, C - \max(\|h_{P_0}^\dagger(\cdot, d_e)\|_{L^2(P_{X|d_e;0})}, \|q_{P_0}^\dagger(\cdot, d_e)\|_{L^2(P_{Z|d_e;0})}) \right),
    \end{align*}
    which is non-empty by definition of $\mc{P}_{h_{P}^\dagger, q_{P}^\dagger \text{ bound}}$, and set
    \begin{align}
        \label{eq:not-injective-not-diff:gamma-bounds}
        a_1 \coloneqq 2/ C'_{h, q}, a_2 \coloneqq  4/C'_{h, q}
    \end{align}
    so that for our family of QMD at $P_0$ submodels indexed by $\gamma$, $\gamma$ ranges in $(2/ C'_{h, q}, 4/ C'_{h, q})$. Taking an arbitrary $\gamma \in (2/C'_{h, q},4/C'_{h, q})$, we define
    \begin{align*}
        \sigma_{\perp;d_e}(t, \gamma) \coloneqq& \gamma \sin(t) \eta_{P_0}\\
        \eta_{P_0} \coloneqq&\max(|\langle \alpha_{P_0}(\cdot, d_e), u_{i_\alpha, d_e}\rangle_{L^2(P_{X|d_e;0})}|, |\langle r_{P_0}(\cdot, d_e), v_{i_r, d_e}\rangle_{L^2(P_{Z|d_e;0})}|).
    \end{align*}
    Note that $\sigma_{\perp;d_e}(t, \gamma) \geq 0$ whenever $\sin(t) \geq 0$. Our choice of $\sigma_{\perp;d_e}(t, \gamma)$ ensures that $\sigma_{\perp;d_e}(t, \gamma)$ tends to zero as $t\rightarrow 0$ at the same rate as $\langle \alpha_{P_0}(\cdot, d_e), u_{\perp;d_e}(t)\rangle$ tends to zero, so that $q_{P_{t, \gamma}}^\dagger$ stays bounded in norm, and similarly for $\langle r_{P_0}(\cdot, d_e), v_{\perp;d_e}(t)\rangle$ and $h_{P_{t, \gamma}}^\dagger$.

    We now construct our perturbed operator $\mc{E}_{P_{t,\gamma};d_e}$. Let
    \begin{align*}
        u_{i;d_e}(t) \coloneqq u_{i;d_e}, i \not= i_{\alpha}, v_{i;d_e}(t) \coloneqq v_{i;d_e}, i \not= i_{r}, \sigma_{i;d_e}(t) \coloneqq \sigma_{i;d_e}
    \end{align*}
    for $t \geq 0$. Then, we define the operator $\mc{E}_{P_{t,\gamma};d_e}:L^2(P_{X|d_e;0}) \rightarrow L^2(P_{Z|d_e;0})$ by
    \begin{align*}
        (\mc{E}_{t, \gamma;d_e}h)(z) \coloneqq& \sum_{i=1}^N \sigma_{i;d_e}(t)\langle h, u_{i;d_e}(t)\rangle_{L^2(P_{X|d_e;0})}v_{{i;d_e}}(t)(z)\\
        &+ \sigma_{\perp;d_e}(t, \gamma)\langle h, u_{\perp;d_e}(t)\rangle_{L^2(P_{X|d_e;0})}v_{{\perp;d_e}}(t)(z).
    \end{align*}
    This operator is bounded and linear, and has singular system when viewed as mapping $L^2(P_{X|d_e;0}) \rightarrow  L^2(P_{Z|d_e;0})$ given by $(\sigma_{i;d_e}(t), u_{i;d_e}(t), v_{{i;d_e}}(t))_{i=1}^N \cup (\sigma_{\perp;d_e}(t, \gamma), u_{\perp;{d_e}}(t), v_{\perp;{d_e}}(t))$.

    We claim the operator $\mc{E}_{t, \gamma;d_e}$ is a conditional expectation operator for some (conditional on $D=d_e$) distribution $P_{\cdot|d_e;t, \gamma}$ for $t$, $\sigma_{\perp;d_e}(\cdot, \gamma)$ small enough. We do so by arguing that this operator is an integral operator with kernel equal to the kernel of that conditional expectation operator. To see this fact, define
    \begin{align*}
        k_{d_e;0} \coloneqq& \sum_{i=1}^N \sigma_{i;d_e}u_{i;d_e}v_{i;d_e}, &k_{d_e;t, \gamma} \coloneqq& \sum_{i=1}^N \sigma_{i;d_e}(t)u_{i;d_e}(t)v_{i;d_e}(t)\\
        &&&+ \sigma_{\perp;d_e}(t, \gamma)u_{\perp;d_e}(t)v_{{\perp;d_e}}(t).
    \end{align*}
    One can show that $\mc{E}_{P_0;d_e}$ and $\mc{E}_{t, \gamma;d_e}$ are integral operators with respect to the above kernels, and that $p_{X, Z|d_{e};0} = k_{d_e;0}p_{X|d_{e};0}p_{Z|d_{e};0}$. Define
    \begin{align*}
        p_{X, Z|d_{e};t, \gamma} = k_{d_e;t, \gamma}p_{X|d_{e};0}p_{Z|d_{e};0}.
    \end{align*}
    We claim $p_{X, Z|d_{e};t, \gamma}$ is a valid joint density for $t, \gamma$ sufficiently small with appropriate choice of $t \mapsto \sigma_{\perp;d_e}(t, \gamma)$. First,
    \begin{align}
        \label{eq:not-injective-not-diff:margins-constant}
        \nonumber\int p_{X, Z|d_{e};t, \gamma}(x, z)d\lambda_X(x) =& \int\sum_{i=1}^N \sigma_{i;d_e}(t)u_{i;d_e}(t)(x)v_{i;d_e}(t)(z)p_{X|d_{e};0}(x)p_{Z|d_{e};0}(z)d\lambda_X(x)\\
        \nonumber&+ \sigma_{\perp;d_e}(t, \gamma)u_{\perp;d_e}(t)(x)v_{{\perp;d_e}}(t)(z) p_{X|d_{e};0}(x)p_{Z|d_{e};0}(z)d\lambda_X(x)\\
        \nonumber=& \int\sum_{i=1}^N \sigma_{i;d_e}(t)\langle u_{i;d_e}(t), 1 \rangle_{L^2(P_{X|{d_e};0}} v_{i;d_e}(t)(z)p_{Z|d_{e};0}(z)d\lambda_X(x)\\
        \nonumber&+ \sigma_{\perp;d_e}(t, \gamma)\langle u_{\perp;d_e}(t) , 1\rangle_{L^2(P_{X|{d_e};0}}  v_{{\perp;d_e}}(t)(z)p_{Z|d_{e};0}(z)d\lambda_X(x)\\
        \nonumber=& \sigma_{1;d_e}(t)\langle u_{1;d_e}(t), 1 \rangle_{L^2(P_{X|{d_e};0}} v_{1;d_e}(t)(z)p_{Z|d_{e};0}(z)\\
        =&p_{Z|d_{e};0}(z).
    \end{align}
    A similar argument shows $\int p_{X, Z|d_{e};t, \gamma}(x, z)d\lambda_Z(z) = p_{X|d_{e};0}(x)$. Furthermore, by the above $\int \int p_{X, Z|d_{e};t, \gamma}(x, z)d\lambda_X(x)d\lambda_Z(z) = 1$. Then, $p_{X, Z|d_{e};t, \gamma}$ is a valid version of a joint density for $(X, Z)$ if $p_{X, Z|d_{e};t, \gamma} \geq 0$ a.e. $P_0$. Recall by remark \cref{remark:marginal-conditional-density-bounds} and \cref{as:dominating-measure-support} that
    \begin{align*}
        \|p_{X, Z|d_e;0}^{-1}\|_{L^\infty(P_{X, Z|d_e;0})}, \|p_{X|d_e;0}\|_{L^\infty(P_{X|d_e;0})}, \|p_{ Z|d_e;0}\|_{L^\infty(P_{Z|d_e;0})} < L_{bounded}^2.
    \end{align*}
    It follows that
    \begin{align*}
        \essinf k_{d_e;0} \geq \frac{\essinf p_{X, Z|d_e;0}}{\|p_{X|d_e;0}\|_{L^\infty(P_{X|d_e;0})}\|p_{ Z|d_e;0}\|_{L^\infty(P_{ Z|d_e;0})}} > L_{bounded}^{-6}.
    \end{align*}
    On the other hand,
    \begin{align}
        \label{eq:not-injective-not-diff:kernel-bound}
        \zeta(t, \gamma) \coloneqq \nonumber \|k_{d_e;0} - k_{d_e;t, \gamma}\|_{L^\infty(P_{X, Z|d_e;0})} \leq& \sigma_{i_\alpha;d_e}\|u_{i_\alpha;d_e}(t)v_{i_{\alpha};d_e}(t) - u_{i_\alpha;d_e}v_{i_{\alpha};d_e}\|_{L^\infty(P_{X, Z|d_e;0})}\\
        \nonumber &+\sigma_{i_r;d_e}\|u_{i_r;d_e}(t)v_{i_{r};d_e}(t) - u_{i_r;d_e}v_{i_{r};d_e}\|_{L^\infty(P_{X, Z|d_e;0})}\\
        &+\sigma_{\perp;d_e}(t, \gamma)\|u_{\perp;d_e}(t)v_{\perp;d_e}(t) \|_{L^\infty(P_{X, Z|d_e;0})},
    \end{align}
    where this bound holds whether or not $i_\alpha = i_r$. Now, by \cref{lemma:in-range-uniformly-bounded}, and by choice of $u_{\perp;d_e}, v_{\perp;d_e}$,
    \begin{align*}
        M \coloneqq \max(\|u_{i_{\alpha};d_e}\|_{L^\infty(P_{X;0})}, \|u_{i_{r};d_e}\|_{L^\infty(P_{X;0})}, \|u_{\perp;d_e}\|_{L^\infty(P_{X;0})}, \|v_{i_{\alpha};d_e}\|_{L^\infty(P_{Z;0})}, \|v_{i_{r};d_e}\|_{L^\infty(P_{Z;0})}, \|v_{\perp;d_e}\|_{L^\infty(P_{Z;0})}) < \infty.
    \end{align*}
    On the other hand, by \eqref{eq:not-injective-not-diff:rotation},
    \begin{align*}
        \|u_{i_\alpha;d_e}(t) - u_{i_\alpha;d_e}\|_{L^\infty(P_{X;0})} \leq& |1 - \cos (t) + \sin(t)|M, &\|v_{i_r;d_e}(t) - v_{i_r;d_e}\| _{L^\infty(P_{Z;0})}\leq& |1 - \cos (t) + \sin(t)|M,\\
        \|u_{i_\alpha;d_e}(t)\|_{L^\infty(P_{X;0})} \leq &2M & \|v_{i_r;d_e}(t)\|_{L^\infty(P_{Z;0})} \leq &2M.
    \end{align*}
    Note that if $i_r = i_\alpha$,
    \begin{align*}
        \|u_{i_r;d_e}(t) - u_{i_r;d_e}\|_{L^\infty(P_{X;0})} \leq& |1 - \cos (t) + \sin(t)|M & \|v_{i_\alpha;d_e}(t) - v_{i_\alpha;d_e}\| _{L^\infty(P_{Z;0})}\leq& |1 - \cos (t) + \sin(t)|M
    \end{align*}
    and the above differences equal 0 when $i_r \not= i_\alpha$. As such,
    \begin{align*}
        \|u_{i_\alpha;d_e}(t)v_{i_{\alpha};d_e}(t) - u_{i_\alpha;d_e}v_{i_{\alpha};d_e}\|_{L^\infty(P_{X, Z|d_e;0})} \leq& \|\{u_{i_\alpha;d_e}(t) -u_{i_\alpha;d_e}\}v_{i_{\alpha};d_e}(t)\|_{L^\infty(P_{X, Z|d_e;0})}  \\
        &+\|u_{i_\alpha;d_e}\{v_{i_{\alpha};d_e}(t)  - v_{i_{\alpha};d_e}\}\|_{L^\infty(P_{X, Z|d_e;0})}\\
        \leq&  3M^2|1 - \cos (t) + \sin(t)|
    \end{align*}
    and similarly $\|u_{i_r;d_e}(t)v_{i_{r};d_e}(t) - u_{i_r;d_e}v_{i_{r};d_e}\|_{L^\infty(P_{X, Z|d_e;0})} \leq 3M^2|1 - \cos (t) + \sin(t)|$. The third term on the right hand side of \eqref{eq:not-injective-not-diff:kernel-bound} is uniformly bounded by $4M^2\sigma_{\perp;d_e}(t, \gamma) = 4M^2 \gamma \sin(t) \eta_{P_0}$. Then, we have that
    \begin{align*}
        \nonumber \|k_{d_e;0} - k_{d_e;t, \gamma}\|_{L^\infty(P_{X, Z|d_e;0})} \leq& 4M^2 \gamma \sin(t) \eta_{P_0} + 6M^2|1 - \cos (t) + \sin(t)|\\
        \leq& \{4M^2 4/C'_{h, q}\}\sin(t) \eta_{P_0} + 6M^2|1 - \cos (t) + \sin(t)|
    \end{align*}
    which tends to zero as $t \rightarrow 0$ and is non-negative for $t < \pi$, where the second inequality follows from \eqref{eq:not-injective-not-diff:gamma-bounds}.

    The above discussion shows the kernels converge uniformly as $t \rightarrow 0$. Note that
    \begin{align*}
        \|p_{X, Z|d_e;t, \gamma} - p_{X, Z|d_e;0}\|_{L^\infty(P_{X, Z|d_e;0})} \leq& \|k_{d_e;t, \gamma} - k_{d_e;0}\|_{L^\infty(P_{X, Z|d_e;0})}\| p_{X|d_e;0}\|_{L^\infty(P_{X|d_e;0})}\| p_{Z|d_e;0}\|_{L^\infty(P_{Z|d_e;0})}\\
        &\leq \|k_{d_e;t, \gamma} - k_{d_e;0}\|_{L^\infty(P_{X, Z|d_e;0})}L_{bounded}^4,
    \end{align*}
    and so the densities $p_{X, Z|d_e;t, \gamma}, p_{X, Z|d_e;0}$, also converge uniformly as $t\rightarrow 0$, uniformly over $\gamma \in (2 / C'_{h, q}, 4 / C'_{h, q})$. Recalling that $p_{X, Z|d_e;0} > L_{bounded}^{-2}$, then shows that for $t$ small enough, $p_{X, Z|d_e;t, \gamma} > 0$ a.e. $P_{X, Z|d_e;0}$. As such, it is a valid (version of a) joint density for $(X, Z)$, also with support equal the support of $(X, Z)$ under $P_0$.

    Now, set $P_{{t,\gamma}}(D=d) = P_0(D=d)$ for each $d$ and
    \begin{align*}
        p_{X, Z|d;t,\gamma} \coloneqq p_{X, Z|d;0} \text{ for } d\not= d_e.
    \end{align*}
    By \eqref{eq:not-injective-not-diff:margins-constant}, we have that $p_{X, D;t, \gamma} = p_{X, D;0}$ and similarly that $p_{Z, D;t, \gamma} = p_{Z, D;0}$.
    As in the proof of \cref{lemma:constant-riesz-representers}, we let
    \begin{align*}
        p_{t, \gamma}(w, x, z, d) \coloneqq p_0(w, x, z, d) + p_{W;0}(w)\{p_{X, Z, D;t,\gamma}(x, z, d) - p_{X, Z, D;0}(x, z, d)\}.
    \end{align*}
    Following the proof of that lemma, by choosing $\delta$ sufficiently small, for all $t < \delta$, for all $\gamma \in (2/ C'_{h, q},4/ C'_{h, q} )$, the above is a valid joint density for $O$, mutually absolutely continuous with $P_0$, with the same support, with density $(X, Z, D)$ under $P_{t, \gamma}$ equal to $p_{X, Z, D;t, \gamma}$, satisfying \cref{as:continuous-linear-maps}, and \cref{as:dominating-measure-support} with $\alpha_{P_{t,\gamma}} = \alpha_{P_0}, r_{P_{t,\gamma}} = r_{P_0}$. Furthermore, in light of the above discussion and because $p_{0}$ is bounded away from $L_{bounded}^{-1}$ below and $L_{bounded}$ above, we may choose $\delta$ sufficiently small such that for all $t < \delta$, $\gamma \in (2/ C'_{h, q},4/ C'_{h, q} )$
    \begin{align*}
        L_{bounded}^{-1} < \essinf p_{t, \gamma} , \esssup p_{t, \gamma} < L_{bounded},
    \end{align*}
    so that $P_{t, \gamma}$ satisfies \cref{as:dominating-measure-support}.

    We now show that $\mc{E}_{t, \gamma;d_e} = \mc{E}_{P_{t, \gamma};d_e}$, where we note this equality holds trivially for $d \not = d_e$. Let $h \in L^2(P_{X|d_e;0})$. Then,
    \begin{align*}
        (\mc{E}_{P_{t, \gamma};d_e} h)(z) =& \int h(x)p_{X|Z, d_e;t, \gamma}(x, z)d\lambda_X(x)\\
        =& \int h(x)p_{X, Z| d_e;t, \gamma}(x, z)p_{X|d_e;0}(x)^{-1}p_{Z|d_e;0}(z)^{-1}p_{X|d_e;0}(x)d\lambda_X(x)\\
        =& \int h(x) k_{d_e;t, \gamma}(x, z)p_{X|d_e;0}(x)d\lambda_X(x)\\
        =& (\mc{E}_{t, \gamma;d_e}h)(z).
    \end{align*}

    Now, \cref{as:solutions-exist} holds at $P_{t, \gamma}$ because $\mathscr{R}(\mc{E}_{P_{0};d_e}) \subset\mathscr{R}(\mc{E}_{P_{t, \gamma};d_e})$, $\mathscr{R}(\mc{E}_{P_{0};d_e}^*) \subset\mathscr{R}(\mc{E}_{P_{t, \gamma};d_e}^*)$ and $\mathscr{R}(\mc{E}_{P_{0};d_e}) = \mathscr{R}(\mc{E}_{P_{t, \gamma};d}), \mathscr{R}(\mc{E}^*_{P_{0};d_e}) = \mathscr{R}(\mc{E}^*_{P_{t, \gamma};d})$ for $d \not=d_e$. The inclusions at $d = d_e$ hold because the rotated pairs $u_{i_\alpha;d_e}(t), u_{\perp;d_e}(t)$ and $v_{i_r;d_e}(t), v_{\perp;d_e}(t)$ span the same planes as $u_{i_\alpha;d_e}, u_{\perp;d_e}$ and $v_{i_r;d_e}, v_{\perp;d_e}$, respectively, with strictly positive singular values, so that $\mathscr{N}(\mc{E}^*_{P_{t,\gamma};d_e}) = \mathscr{N}(\mc{E}^*_{P_0;d_e}) \ominus \text{span}\{v_{\perp;d_e}\}$ and $\mathscr{N}(\mc{E}_{P_{t,\gamma};d_e}) = \mathscr{N}(\mc{E}_{P_0;d_e}) \ominus \text{span}\{u_{\perp;d_e}\}$, while Picard's condition is unaffected because only finitely many singular values and functions change. \cref{as:beta-source} holds at $P_{t, \gamma}$ for any choice of $\beta_h, \beta_q$ for which it holds at $P_0$, as we have perturbed only a finite number of singular vectors and values. \cref{as:mutual-absolute-continuity} holds at $P_{t, \gamma}$ by construction.

    It remains to show that both parts of \cref{as:uniform-norm-bound} hold. We begin with \cref{as:uniform-norm-bound:part:h}. Because $r_{P_{t, \gamma}}(\cdot, d_e) = r_{P_0}(\cdot, d_e) \in \mathscr{R}(\mc{E}_{P_{t, \gamma}})$, the minimum norm solution to $\mc{E}_{P_{t, \gamma};d_e}h = r_{P_0}(\cdot, d_e)$ is given by
    \begin{align*}
        h_{P_{t, \gamma}}^\dagger(\cdot, d_e) =& \sum_{i=1}^N \frac{1}{\sigma_{i;d_e}} u_{i;d_e}(t) \langle r_{P_0}(\cdot, d_e), v_{i;d_e}(t)\rangle _{L^2(P_{Z|d_e;0})}\\
        &+ \frac{1}{\sigma_{\perp;d_e}(t, \gamma)} u_{\perp;d_e}(t) \langle r_{P_0}(\cdot, d_e), v_{\perp;d_e}(t)\rangle _{L^2(P_{Z|d_e;0})}.
    \end{align*}
    Then,
    \begin{align*}
        h_{P_{t, \gamma}}^\dagger(\cdot, d_e) - h_{P_0}^\dagger(\cdot, d_e)=& \nu(t)+ \frac{1}{\sigma_{\perp;d_e}(t, \gamma)} u_{\perp;d_e}(t) \langle r_{P_0}(\cdot, d_e), v_{\perp;d_e}(t)\rangle _{L^2(P_{Z|d_e;0})}.
    \end{align*}
    where
    \begin{align*}
        \nu(t) \coloneqq& \frac{1}{\sigma_{i_\alpha;d_e}} \{u_{i_{\alpha};d_e}(t) \langle r_{P_0}(\cdot, d_e), v_{i_{\alpha};d_e}(t)\rangle _{L^2(P_{Z|d_e;0})} - u_{i_{\alpha};d_e} \langle r_{P_0}(\cdot, d_e), v_{i_{\alpha};d_e}\rangle _{L^2(P_{Z|d_e;0})}\}\\
        &+ I(i_{\alpha}\not=i_r)\frac{1}{\sigma_{i_r;d_e}} u_{i_r;d_e} \{\langle r_{P_0}(\cdot, d_e), v_{i_r;d_e}(t) - v_{i_r;d_e}\rangle _{L^2(P_{Z|d_e;0})}\}
    \end{align*}
    Now, from \eqref{eq:not-injective-not-diff:rotation},
    \begin{align*}
        \|\nu(t)\|_{L^2(P_{X|d_e;0})} \rightarrow 0 \text{ as }t \rightarrow 0.
    \end{align*}
    On the other hand,
    \begin{align*}
        &\left\|\frac{1}{\sigma_{\perp;d_e}(t, \gamma)} u_{\perp;d_e}(t) \langle r_{P_0}(\cdot, d_e), v_{\perp;d_e}(t)\rangle _{L^2(P_{Z|d_e;0})}\right\|_{L^2(P_{X|d_e;0})}\\
        =& \frac{1}{\gamma \sin(t) \eta }\langle r_{P_0}(\cdot, d_e), v_{\perp;d_e}(t)\rangle _{L^2(P_{Z|d_e;0})} \left\|u_{\perp;d_e}(t) \right\|_{L^2(P_{X|d_e;0})}\\
        =& \frac{1}{\gamma \sin(t) \eta_{P_0} }\sin(t) \langle r_{P_0}(\cdot, d_e), v_{i_r;d_e}\rangle _{L^2(P_{Z|d_e;0})}\\
        \leq& \frac{1}{\gamma}
    \end{align*}
    where we have used that $r_{P_0}(\cdot, d_e) \perp_{L^2(P_{Z|d_e;0})}v_{\perp;d_e}$ and that $\eta_{P_0} \geq |\langle r_{P_0}(\cdot, d_e), v_{i_r;d_e}\rangle _{L^2(P_{Z|d_e;0})}|$. As such,
    \begin{align*}
        \|h^\dagger_{P_{t, \gamma}}(\cdot, d_e)\|_{L^2(P_{X|d_e;0})} \leq&  \|h^\dagger_{P_{0}}(\cdot, d_e)\|_{L^2(P_{X|d_e;0})} + \|h_{P_{t, \gamma}}^\dagger(\cdot, d_e) - h_{P_0}^\dagger(\cdot, d_e)\|_{L^2(P_{X|d_e;0})}\\
        \leq& \|h^\dagger_{P_{0}}(\cdot, d_e)\|_{L^2(P_{X|d_e;0})}  +\frac{1}{\gamma} +  \|\nu(t)\|_{L^2(P_{X|d_e;0})}  \\
        \leq&  \|h^\dagger_{P_{0}}(\cdot, d_e)\|_{L^2(P_{X|d_e;0})}  +\frac{C'_{h, q}}{2} +  \|\nu(t)\|_{L^2(P_{X|d_e;0})} \\
        <& C \text{ for } t \text{ s.t. } \|\nu(t)\|_{L^2(P_{X|d_e;0})}  < C'_{h, q}/2.
    \end{align*}
    Then, for $t$ small enough, not depending on $\gamma$, $\|h^\dagger_{P_{t, \gamma}}(\cdot, d_e)\|_{L^2(P_{X|d_e;0})} < C$. A similar argument gives the same bound for $\|q^\dagger_{P_{t, \gamma}}(\cdot, d_e)\|_{L^2(P_{Z|d_e;0})}$. For $d \not = d_e$, $h^\dagger_{P_{t, \gamma}}(\cdot, d) = h^\dagger_{P_{0}}(\cdot, d)$ and similarly $q^\dagger_{P_{t, \gamma}}(\cdot, d) = q^\dagger_{P_0}(\cdot, d)$. Consequently, both parts  of \cref{as:uniform-norm-bound} hold for $\gamma \in (2 / C'_{h, q}, 4/C'_{h, q})$,  $t < \delta$, for $\delta$ small enough. Combining the above discussion, we have that $P_{t, \gamma} \in \mc{P}_{h_{P}^\dagger, q_{P}^\dagger \text{ bound}}$ for $t$ small enough.

    We now compute the derivative
    \begin{align*}
        \frac{d}{dt}\varphi(P_{t,\gamma})\vert_{t=0}
    \end{align*}
    for $\gamma\in (2 / C'_{h, q}, 4/C'_{h, q})$ and show that it depends inversely on $\gamma$.
    \begin{align*}
        \varphi(P_{t, \gamma}) - \varphi(P_0) =& E_{P_{t, \gamma}}[\alpha_{P_{t, \gamma}}(X, D)h_{P_{t, \gamma}}^\dagger(X, D)] - E_{P_0}[\alpha_{P_0}(X, D)h_{P_0}^\dagger(X, D)]\\
        =& \sum_{d \in \textsf{Supp}[P_{D;0}]}p_0(D=d)\langle \alpha_{P_0}(\cdot, d), h_{P_{t, \gamma}}^\dagger(\cdot, d) - h_{P_0}^\dagger(\cdot, d)\rangle_{L^2(P_{X|d;0})}\\
        =& p_0(D=d_e)\langle \alpha_{P_0}(\cdot, d_e), h_{P_{t, \gamma}}^\dagger(\cdot, d_e) - h_{P_0}^\dagger(\cdot, d_e)\rangle_{L^2(P_{X|d;0})}.
    \end{align*}
    But, by construction,
    \begin{align*}
        \langle \alpha_{P_0}(\cdot, d_e), h_{P_{t, \gamma}}^\dagger(\cdot, d_e) - h_{P_0}^\dagger(\cdot, d_e)\rangle_{L^2(P_{X|d;0})}
        =& \langle \alpha_{P_0}(\cdot, d_e), \nu(t)\rangle_{L^2(P_{X|d_e;0})}\\
        &+ \frac{1}{\sigma_{\perp;d_e}(t, \gamma)} \langle \alpha_{P_0}(\cdot, d_e), u_{\perp;d_e}(t)\rangle_{L^2(P_{X|d;0})}  \langle r_{P_0}(\cdot, d_e), v_{\perp;d_e}(t)\rangle _{L^2(P_{Z|d_e;0})}.
    \end{align*}
    Now, by \eqref{eq:not-injective-not-diff:rotation}, if $i_{r} \not = i_{\alpha}$,
    \begin{align*}
        \langle \alpha_{P_0}(\cdot, d_e), \nu(t)\rangle_{L^2(P_{X|d_e;0})} =&\frac{1}{\sigma_{i_\alpha;d_e}} \langle \alpha_{P_0}(\cdot, d_e),u_{i_{\alpha};d_e}(t) -  u_{i_{\alpha};d_e}\rangle_{L^2(P_{X|d_e;0})}\langle r_{P_0}(\cdot, d_e),v_{i_{\alpha};d_e}\rangle_{L^2(P_{Z|d_e;0})}\\
        &+ \frac{1}{\sigma_{i_r;d_e}}\langle \alpha_{P_0}(\cdot, d_e),u_{i_r;d_e}\rangle_{L^2(P_{X|d_e;0})}\langle r_{P_0}(\cdot, d_e),v_{i_{r};d_e}(t) - v_{i_{r};d_e}\rangle_{L^2(P_{Z|d_e;0})}\\
        =& \frac{1}{\sigma_{i_\alpha;d_e}} \{\cos(t) - 1\}\langle \alpha_{P_0}(\cdot, d_e), u_{i_{\alpha};d_e}\rangle_{L^2(P_{X|d_e;0})}\langle r_{P_0}(\cdot, d_e),v_{i_{\alpha};d_e}\rangle_{L^2(P_{Z|d_e;0})}\\
        &+ \frac{1}{\sigma_{i_r;d_e}}\{\cos(t) - 1\}\langle \alpha_{P_0}(\cdot, d_e),u_{i_r;d_e}\rangle_{L^2(P_{X|d_e;0})}\langle r_{P_0}(\cdot, d_e),v_{i_{r};d_e}\rangle_{L^2(P_{Z|d_e;0})}.
    \end{align*}
    When $i_r = i_\alpha$,
    \begin{align*}
        \langle \alpha_{P_0}(\cdot, d_e), \nu(t)\rangle_{L^2(P_{X|d_e;0})} =&\frac{1}{\sigma_{i_\alpha;d_e}} \{\langle \alpha_{P_0}(\cdot, d_e), u_{i_{\alpha};d_e}(t)\rangle_{L^2(P_{X|d_e;0})} \langle r_{P_0}(\cdot, d_e), v_{i_{\alpha};d_e}(t)\rangle _{L^2(P_{Z|d_e;0})} \\
        &- \langle \alpha_{P_0}(\cdot, d_e), u_{i_{\alpha};d_e}\rangle_{L^2(P_{X|d_e;0})} \langle r_{P_0}(\cdot, d_e), v_{i_{\alpha};d_e}\rangle _{L^2(P_{Z|d_e;0})}\}\\
        =& \frac{1}{\sigma_{i_\alpha;d_e}} \{\cos(t)^2\langle \alpha_{P_0}(\cdot, d_e), u_{i_{\alpha};d_e}\rangle_{L^2(P_{X|d_e;0})} \langle r_{P_0}(\cdot, d_e), v_{i_{\alpha};d_e}\rangle _{L^2(P_{Z|d_e;0})} \\
        &- \langle \alpha_{P_0}(\cdot, d_e), u_{i_{\alpha};d_e}\rangle_{L^2(P_{X|d_e;0})} \langle r_{P_0}(\cdot, d_e), v_{i_{\alpha};d_e}\rangle _{L^2(P_{Z|d_e;0})}\}.
    \end{align*}
    In either case, because $\alpha_{P_0}(\cdot, d_e) \perp u_{\perp;d_e}$ and $r_{P_0}(\cdot, d_e) \perp v_{\perp;d_e}$, the right hand side above reduces to a multiple of $\cos(t) - 1$ (respectively $\cos(t)^2 - 1$), which is $O(t^2)$ as $t \rightarrow 0$. Hence
    \begin{align*}
        \lim_{t\downarrow 0}\frac{\langle \alpha_{P_0}(\cdot, d_e), \nu(t)\rangle_{L^2(P_{X|d_e;0})}}{t} = 0,
    \end{align*}
    independently of $\gamma$.

    On the other hand,
    \begin{align*}
        &\frac{1}{\sigma_{\perp;d_e}(t, \gamma)} \langle \alpha_{P_0}(\cdot, d_e), u_{\perp;d_e}(t)\rangle_{L^2(P_{X|d;0})}  \langle r_{P_0}(\cdot, d_e), v_{\perp;d_e}(t)\rangle _{L^2(P_{Z|d_e;0})}\\
        =&\frac{1}{\gamma \sin(t) \eta_{P_0}} \sin(t)^2\langle \alpha_{P_0}(\cdot, d_e), u_{i_\alpha;d_e}\rangle_{L^2(P_{X|d;0})}  \langle r_{P_0}(\cdot, d_e), v_{i_r;d_e}\rangle _{L^2(P_{Z|d_e;0})}\\
        =&\frac{1}{\gamma \eta_{P_0}} \sin(t)\langle \alpha_{P_0}(\cdot, d_e), u_{i_\alpha;d_e}\rangle_{L^2(P_{X|d;0})}  \langle r_{P_0}(\cdot, d_e), v_{i_r;d_e}\rangle _{L^2(P_{Z|d_e;0})}.
    \end{align*}
    Now,
    \begin{align*}
        &\lim_{t\downarrow 0}\frac{1}{t}\frac{1}{\gamma \eta_{P_0}} \sin(t)\langle \alpha_{P_0}(\cdot, d_e), u_{i_\alpha;d_e}\rangle_{L^2(P_{X|d;0})}  \langle r_{P_0}(\cdot, d_e), v_{i_r;d_e}\rangle _{L^2(P_{Z|d_e;0})}\\
        =& \frac{\langle \alpha_{P_0}(\cdot, d_e), u_{i_\alpha;d_e}\rangle_{L^2(P_{X|d;0})}  \langle r_{P_0}(\cdot, d_e), v_{i_r;d_e}\rangle _{L^2(P_{Z|d_e;0})}}{\gamma\eta_{P_0}}.
    \end{align*}
    To simplify notation, we define
    \begin{align*}
        \mu \coloneqq P_0(D=d_e)\frac{\langle \alpha_{P_0}(\cdot, d_e), u_{i_\alpha;d_e}\rangle_{L^2(P_{X|d_e;0})}  \langle r_{P_0}(\cdot, d_e), v_{i_r;d_e}\rangle _{L^2(P_{Z|d_e;0})}}{\eta_{P_0}},
    \end{align*}
    and note that $\mu\not=0$ by our choice of $i_{\alpha}, i_r$ in \cref{eq:not-injective-not-diff:i-r-i-alpha-def}.
    Then, for each $\gamma \in (2 / C'_{h, q}, 4/C'_{h, q})$,
    \begin{align}
        \label{eq:not-injective-not-diff:non-linear-derivative}
        \left.\frac{d}{dt}\varphi(P_{t, \gamma})\right\vert_{t=0} = \frac{\mu}{\gamma}.
    \end{align}

    It remains to show that for each $\gamma \in (2 / C'_{h, q}, 4/C'_{h, q})$, $\{P_{t, \gamma}:t \in [0, \delta)\}$ is QMD at $P_0$ with score that depends affinely on $\gamma$. Recall that
    \begin{align*}
        p_{t, \gamma}(w, x, z, d) = p_0(w, x, z, d) + p_{W;0}(w)\{p_{X, Z, D;t, \gamma}(x, z, d) - p_{X, Z, D;0}(x, z, d)\}.
    \end{align*}
    And notice that $p_{X, Z|d_e;t, \gamma}$ is, for each $x, z$, a polynomial in $\cos(t), \sin(t)$ whose coefficients are uniformly bounded, by $M$ and by $\gamma < 4/C'_{h, q}$. As such, it can be extended to a differentiable function in a neighborhood of $t=0$, and this derivative is uniformly bounded over $t$ in that neighborhood and over $(x, z)$ in the support of $P_{X, Z|d_e;0}$. Then,
    \begin{align*}
        \left.\frac{d}{dt}p_{t, \gamma}(w, x, z, d)\right\vert_{t=0} =&I(d=d_e)p_0(D=d_e)p_{W;0}(w)\left.\frac{d}{dt}p_{X, Z|d_e;t, \gamma}(x, z)\right\vert_{t=0}\\
        =& I(d=d_e)p_0(D=d_e)p_{W;0}(w)p_{X|d_e;0}(x)p_{Z|d_e;0}(z)\left.\frac{d}{dt}k_{d_e;t, \gamma}(x, z)\right\vert_{t=0}.
    \end{align*}
    We have that
    \begin{align*}
        \left.\frac{d}{dt}k_{d_e;t, \gamma}(x, z)\right\vert_{t=0} =& \left.\frac{d}{dt}\sigma_{i_{\alpha};d_e}u_{i_\alpha;d_e}(t)v_{i_\alpha;d_e}(t)\right\vert_{t=0} + I(i_{\alpha}\not=i_r)\left.\frac{d}{dt}\sigma_{i_{r};d_e}u_{i_r;d_e}(t)v_{i_r;d_e}(t)\right\vert_{t=0}\\
        &+ \left.\frac{d}{dt}\sigma_{\perp;d_e}(t, \gamma)u_{\perp;d_e}(t)v_{\perp;d_e}(t)\right\vert_{t=0}\\
        =& -\sigma_{i_\alpha;d_e}u_{\perp;d_e}v_{i_{\alpha};d_e} -\sigma_{i_r;d_e}u_{i_r;d_e}v_{\perp;d_e} + \gamma \eta_{P_0}u_{\perp;d_e}v_{\perp;d_e}
    \end{align*}
    where the expression holds whether or not $i_\alpha = i_r$.
    Letting
    \begin{align*}
        \dot\Delta(x, z, d) \coloneqq I(D=d_e)p_{X|d_e;0}(x)p_{Z|d_e;0}(z)p_0(D=d_e)\{-\sigma_{i_\alpha;d_e}u_{\perp;d_e}v_{i_{\alpha};d_e} -\sigma_{i_r;d_e}u_{i_r;d_e}v_{\perp;d_e} + \gamma \eta_{P_0}u_{\perp;d_e}v_{\perp;d_e}\},
    \end{align*}
    the log-density $\log p_{t, \gamma}(w, x, z, d)$ is pointwise in $w, x, z, d$ differentiable at $t=0$ with derivative
    \begin{align*}
        s_\gamma(w, x, z, d) \coloneqq \left.\frac{d}{dt}\log p_{t, \gamma}(w, x, z, d)\right\vert_{t=0} =& \frac{\left.\frac{d}{dt} p_{t, \gamma}(w, x, z, d)\right\vert_{t=0}}{p_0(w, x, z, d)} = \frac{p_{W;0}(w)\dot\Delta(x, z, d)}{p_0(w, x, z, d)}I((w, x, z, d) \in \textsf{Supp}(P_0))
    \end{align*}
    We now argue formally that the above is the score of $P_{t, \gamma}$ at $t=0$ for each $\gamma$. First, because $u_{\perp;d_e}, v_{i_{\alpha};d_e}, u_{i_r;d_e}, v_{\perp;d_e}, u_{\perp;d_e}v_{\perp;d_e}$ are each uniformly bounded above, $p_{W;0}$ is uniformly bounded above, and $p_0(w, x, z, d)$ is uniformly bounded below on the support of $O = (W, X, Z, D)$ in $P_0$, $s_\gamma(X, W, Z, D)$ has finite second moment. Next,  we show the proposed score has mean zero at $P_0$.
    \begin{align*}
        E_{P_0}[s_\gamma(O)] =& \int_{\textsf{Supp}(P_0)} p_{W;0}(w)\dot\Delta(x, z, d)d\lambda(w, x, z, d) \\
        =&\int_{\textsf{Supp}(P_{X, Z, D;0})} \frac{d}{dt}p_{X, Z, D;t, \gamma}d\lambda_X(x)d\lambda_Z(z)d\lambda_D(d)\\
        =& \frac{d}{dt}\int_{\textsf{Supp}(P_{X, Z, D;0})} p_{X, Z, D;t, \gamma}d\lambda_X(x)d\lambda_Z(z)d\lambda_D(d) \\
        =&0
    \end{align*}
    where interchanging the derivative and integral is justified because $\frac{d}{dt}p_{X, Z, D;t, \gamma}(x, z, d)$ exists for every $t$ and is uniformly bounded over $t$ in a neighborhood of $t=0$ and over $(x, z, d)$ in the bounded rectangle $\textsf{Supp}(P_{X, Z, D;0})$.

    Now,
    \begin{align*}
        \lim_{t\downarrow 0}\int \left\{\frac{\sqrt{p_{t, \gamma}(o)} - \sqrt{p_0}(o)}{t} - \frac{1}{2}s_{\gamma}(o)\sqrt{p_0(o)} \right\}^2d\lambda(o) = \int \lim_{t\downarrow 0}\left\{\frac{\sqrt{p_{t, \gamma}(o)} - \sqrt{p_0}(o)}{t} - \frac{1}{2}s_{\gamma}(o) \sqrt{p_0(o)}\right\}^2d\lambda(o) = 0
    \end{align*}
    where again interchanging the derivative and integral is justified because $\frac{d}{dt}p_{t, \gamma}(w, x, z, d)$ is, by the same argument as above, uniformly bounded for $t$ in a neighborhood of $t=0$, and because $p_{t, \gamma}$ is bounded away from $0$ by $L_{bounded}^{-1}$ for $t < \delta$; consequently, the pointwise derivative of the log-density $\frac{d}{dt}\log p_{t, \gamma}(w, x, z, d)$ exists and is uniformly bounded for every $t$ in a neighborhood of $t=0$ over $(w, x, z, d) \in \textsf{Supp}(P_0)$. Finally, the score $s_\gamma(x, w, z, d)$ is affine in $\gamma$ with
    \begin{align*}
        s_\gamma(o) =& s^*(o) + \gamma \bar{s}(o), \\
        s^*(o) =& \frac{p_{W;0}(w)}{p_0(o)}I(d=d_e)p_{X|d_e;0}(x)p_{Z|d_e;0}(z)p_0(d=d_e)\{-\sigma_{i_\alpha;d_e}u_{\perp;d_e}v_{i_{\alpha};d_e} -\sigma_{i_r;d_e}u_{i_r;d_e}v_{\perp;d_e}\}\\
        \bar{s}(o)=&\frac{p_{W;0}(w)}{p_0(o)}I(d=d_e)p_{X|d_e;0}(x)p_{Z|d_e;0}(z)p_0(d=d_e)\eta_{P_0}u_{\perp;d_e}v_{\perp;d_e}.
    \end{align*}

    Hence, we have exhibited a family of submodels $\{\{P_{t, \gamma}:t\in [0, \delta)\}, \gamma \in (2/C'_{h, q}, 4/C'_{h, q})\}\subseteq \mc{P}_{h_{P}^\dagger, q_{P}^\dagger \text{ bound}}$ for which
    \begin{align*}
        \frac{d}{dt}\varphi(P_{t, \gamma})\vert_{t=0}
    \end{align*}
    exists for each $\gamma$ but where
    \begin{align*}
        \gamma \mapsto \frac{d}{dt}\varphi(P_{t, \gamma})\vert_{t=0} = \frac{\mu}{\gamma}
    \end{align*}
    is not affine even though
    \begin{align*}
        \gamma \mapsto s_{\gamma} = s^*(o) + \gamma \bar{s}(o)
    \end{align*}
    is affine. Hence, because $\mu \not=0$, the local parameter is not a bounded linear operator of the score of the submodel $P_{t, \gamma}$, and so $P \mapsto \varphi(P)$ is not pathwise differentiable at $P_0$ over the model $\mc{P}_{h_{P}^\dagger, q_{P}^\dagger \text{ bound}}$.
\end{proof}

\printbibliography[heading=subbibintoc, title={References}]

@book{adams_sobolev_2003-1,
  title = {Sobolev Spaces},
  author = {Adams, Robert A.},
  namea = {Fournier, John J. F.},
  nameatype = {collaborator},
  date = {2003},
  series = {Pure and Applied Mathematics ; v. 140},
  edition = {2nd ed.},
  publisher = {Academic Press},
  location = {Amsterdam},
  isbn = {978-1-281-07246-7},
  pagetotal = {321}
}

@article{andrews_weak_2019,
  title = {Weak {{Instruments}} in {{Instrumental Variables Regression}}: {{Theory}} and {{Practice}}},
  shorttitle = {Weak {{Instruments}} in {{Instrumental Variables Regression}}},
  author = {Andrews, Isaiah and Stock, James H. and Sun, Liyang},
  date = {2019-08-02},
  journaltitle = {Annual Review of Economics},
  volume = {11},
  pages = {727--753},
  publisher = {Annual Reviews},
  issn = {1941-1383, 1941-1391},
  doi = {10.1146/annurev-economics-080218-025643},
  url = {https://www.annualreviews.org/content/journals/10.1146/annurev-economics-080218-025643},
  urldate = {2026-07-27},
  issue = {Volume 11, 2019},
  langid = {english}
}

@article{babiiHonestConfidenceSets2020,
  title = {Honest {{Confidence Sets}} in {{Nonparametric Iv Regression}} and {{Other Ill-Posed Models}}},
  author = {Babii, Andrii},
  date = {2020},
  journaltitle = {Econometric Theory},
  volume = {36},
  number = {4},
  eprint = {27034668},
  eprinttype = {jstor},
  pages = {658--706},
  publisher = {Cambridge University Press},
  issn = {0266-4666},
  url = {https://www.jstor.org/stable/27034668},
  urldate = {2026-10-12}
}

@article{bahadur_nonexistence_1956,
  title = {The {{Nonexistence}} of {{Certain Statistical Procedures}} in {{Nonparametric Problems}}},
  author = {Bahadur, R. R. and Savage, Leonard J.},
  date = {1956-12},
  journaltitle = {The Annals of Mathematical Statistics},
  volume = {27},
  number = {4},
  pages = {1115--1122},
  publisher = {Institute of Mathematical Statistics},
  issn = {0003-4851, 2168-8990},
  doi = {10.1214/aoms/1177728077},
  url = {https://projecteuclid.org/journals/annals-of-mathematical-statistics/volume-27/issue-4/The-Nonexistence-of-Certain-Statistical-Procedures-in-Nonparametric-Problems/10.1214/aoms/1177728077.full},
  urldate = {2026-07-27},
  langid = {english}
}

@online{bennett_inference_2023,
  title = {Inference on {{Strongly Identified Functionals}} of {{Weakly Identified Functions}}},
  author = {Bennett, Andrew and Kallus, Nathan and Mao, Xiaojie and Newey, Whitney and Syrgkanis, Vasilis and Uehara, Masatoshi},
  date = {2023-07-01},
  eprint = {2208.08291},
  eprinttype = {arXiv},
  eprintclass = {stat},
  doi = {10.48550/arXiv.2208.08291},
  url = {http://arxiv.org/abs/2208.08291},
  urldate = {2026-04-24},
  pubstate = {prepublished}
}

@online{bennettSourceConditionDouble2023b,
  title = {Source {{Condition Double Robust Inference}} on {{Functionals}} of {{Inverse Problems}}},
  author = {Bennett, Andrew and Kallus, Nathan and Mao, Xiaojie and Newey, Whitney and Syrgkanis, Vasilis and Uehara, Masatoshi},
  date = {2023-07-25},
  eprint = {2307.13793},
  eprinttype = {arXiv},
  eprintclass = {stat},
  doi = {10.48550/arXiv.2307.13793},
  url = {http://arxiv.org/abs/2307.13793},
  urldate = {2025-01-21},
  pubstate = {prepublished}
}

@book{bickelEfficientAdaptiveEstimation1998a,
  title = {Efficient and Adaptive Estimation for Semiparametric Models},
  author = {Bickel, Peter J. and Klaasen, Chris and Ritov, Ya'Acov and Wellner, Jon},
  date = {1998},
  publisher = {Springer},
  location = {New York},
  isbn = {978-0-387-98473-5},
  langid = {english}
}

@article{birman_estimates_1977,
  title = {{{ESTIMATES OF SINGULAR NUMBERS OF INTEGRAL OPERATORS}}},
  author = {Birman, M Sh and Solomyak, M Z},
  date = {1977-02-28},
  journaltitle = {Russian Mathematical Surveys},
  shortjournal = {Russ. Math. Surv.},
  volume = {32},
  number = {1},
  pages = {15--89},
  issn = {0036-0279, 1468-4829},
  doi = {10.1070/RM1977v032n01ABEH001592},
  url = {https://www.mathnet.ru/eng/rm3017},
  urldate = {2026-07-27}
}

@article{blundell_semi-nonparametric_2007,
  title = {Semi-{{Nonparametric IV Estimation}} of {{Shape-Invariant Engel Curves}}},
  author = {Blundell, Richard and Chen, Xiaohong and Kristensen, Dennis},
  date = {2007},
  journaltitle = {Econometrica},
  volume = {75},
  number = {6},
  eprint = {4502045},
  eprinttype = {jstor},
  pages = {1613--1669},
  publisher = {[Wiley, Econometric Society]},
  issn = {0012-9682},
  url = {https://www.jstor.org/stable/4502045},
  urldate = {2026-06-19}
}

@article{cai_adaptation_2004,
  title = {An Adaptation Theory for Nonparametric Confidence Intervals},
  author = {Cai, T. Tony and Low, Mark G.},
  date = {2004-10},
  journaltitle = {The Annals of Statistics},
  volume = {32},
  number = {5},
  pages = {1805--1840},
  publisher = {Institute of Mathematical Statistics},
  issn = {0090-5364, 2168-8966},
  doi = {10.1214/009053604000000049},
  url = {https://projecteuclid.org/journals/annals-of-statistics/volume-32/issue-5/An-adaptation-theory-for-nonparametric-confidence-intervals/10.1214/009053604000000049.full},
  urldate = {2026-07-27},
  langid = {english}
}

@incollection{cam_basic_1986,
  title = {Some {{Basic Inequalities}}},
  booktitle = {Asymptotic {{Methods}} in {{Statistical Decision Theory}}},
  author = {Cam, Lucien Le},
  editor = {Cam, Lucien Le},
  date = {1986},
  pages = {46--56},
  publisher = {Springer},
  location = {New York, NY},
  doi = {10.1007/978-1-4612-4946-7_4},
  url = {https://doi.org/10.1007/978-1-4612-4946-7_4},
  urldate = {2026-07-27},
  isbn = {978-1-4612-4946-7},
  langid = {english}
}

@article{canayTestabilityIdentificationNonparametric2013a,
  title = {On the {{Testability}} of {{Identification}} in {{Some Nonparametric Models With Endogeneity}}},
  author = {Canay, Ivan A. and Santos, Andres and Shaikh, Azeem M.},
  date = {2013},
  journaltitle = {Econometrica},
  volume = {81},
  number = {6},
  pages = {2535--2559},
  issn = {1468-0262},
  doi = {10.3982/ECTA10851},
  url = {https://onlinelibrary.wiley.com/doi/abs/10.3982/ECTA10851},
  urldate = {2026-09-22},
  langid = {english}
}

@incollection{carrasco_chapter_2007,
  title = {Chapter 77 {{Linear Inverse Problems}} in {{Structural Econometrics Estimation Based}} on {{Spectral Decomposition}} and {{Regularization}}},
  booktitle = {Handbook of {{Econometrics}}},
  author = {Carrasco, Marine and Florens, Jean-Pierre and Renault, Eric},
  editor = {Heckman, James J. and Leamer, Edward E.},
  date = {2007-01-01},
  volume = {6},
  pages = {5633--5751},
  publisher = {Elsevier},
  doi = {10.1016/S1573-4412(07)06077-1},
  url = {https://www.sciencedirect.com/science/article/pii/S1573441207060771},
  urldate = {2026-07-27}
}

@online{chen_equivalence_2026,
  title = {On the {{Equivalence}} between {{Neyman Orthogonality}} and {{Pathwise Differentiability}}},
  author = {Chen, Yuxi and Kennedy, Edward H. and Balakrishnan, Sivaraman},
  date = {2026-07-06},
  eprint = {2603.15817},
  eprinttype = {arXiv},
  eprintclass = {stat.ME},
  doi = {10.48550/arXiv.2603.15817},
  url = {http://arxiv.org/abs/2603.15817},
  urldate = {2026-07-27},
  pubstate = {prepublished}
}

@article{chen_overidentification_2018-1,
  title = {Overidentification in {{Regular Models}}},
  author = {Chen, Xiaohong and Santos, Andres},
  date = {2018},
  journaltitle = {Econometrica},
  volume = {86},
  number = {5},
  eprint = {44955258},
  eprinttype = {jstor},
  pages = {1771--1817},
  publisher = {[Wiley, The Econometric Society]},
  issn = {0012-9682},
  url = {https://www.jstor.org/stable/44955258},
  urldate = {2026-06-11}
}

@article{chernozhukovDoubleDebiasedMachine2018b,
  title = {Double/Debiased Machine Learning for Treatment and Structural Parameters},
  author = {Chernozhukov, Victor and Chetverikov, Denis and Demirer, Mert and Duflo, Esther and Hansen, Christian and Newey, Whitney and Robins, James},
  date = {2018-02-01},
  journaltitle = {The Econometrics Journal},
  shortjournal = {The Econometrics Journal},
  volume = {21},
  number = {1},
  pages = {C1-C68},
  issn = {1368-4221},
  doi = {10.1111/ectj.12097},
  url = {https://doi.org/10.1111/ectj.12097},
  urldate = {2025-03-27}
}

@book{conway_course_2007-1,
  title = {A {{Course}} in {{Functional Analysis}}},
  author = {Conway, John B.},
  date = {2007},
  series = {Graduate {{Texts}} in {{Mathematics}}},
  volume = {96},
  publisher = {Springer},
  location = {New York, NY},
  doi = {10.1007/978-1-4757-4383-8},
  url = {http://link.springer.com/10.1007/978-1-4757-4383-8},
  urldate = {2026-07-27},
  isbn = {978-1-4419-3092-7},
  langid = {english}
}

@article{cui_semiparametric_2024,
  title = {Semiparametric {{Proximal Causal Inference}}},
  author = {Cui, Yifan and Pu, Hongming and Shi, Xu and Miao, Wang and Tchetgen Tchetgen, Eric},
  date = {2024-04-02},
  journaltitle = {Journal of the American Statistical Association},
  volume = {119},
  number = {546},
  pages = {1348--1359},
  publisher = {Taylor \& Francis},
  issn = {0162-1459},
  doi = {10.1080/01621459.2023.2191817},
  url = {https://doi.org/10.1080/01621459.2023.2191817},
  urldate = {2026-04-02}
}

@article{darolles_nonparametric_2011-1,
  title = {Nonparametric {{Instrumental Regression}}},
  author = {Darolles, S. and Fan, Y. and Florens, J. P. and Renault, E.},
  date = {2011},
  journaltitle = {Econometrica},
  volume = {79},
  number = {5},
  eprint = {41237784},
  eprinttype = {jstor},
  pages = {1541--1565},
  publisher = {[Wiley, Econometric Society]},
  issn = {0012-9682},
  url = {https://www.jstor.org/stable/41237784},
  urldate = {2026-06-19}
}

@online{deaner_nonparametric_2022,
  title = {Nonparametric {{Instrumental Variables Estimation Under Misspecification}}},
  author = {Deaner, Ben},
  date = {2022-12-11},
  eprint = {1901.01241},
  eprinttype = {arXiv},
  eprintclass = {econ.EM},
  doi = {10.48550/arXiv.1901.01241},
  url = {http://arxiv.org/abs/1901.01241},
  urldate = {2026-07-16},
  pubstate = {prepublished}
}

@article{donoho_geometrizing_1991,
  title = {Geometrizing {{Rates}} of {{Convergence}}, {{II}}},
  author = {Donoho, David L. and Liu, Richard C.},
  date = {1991-06},
  journaltitle = {The Annals of Statistics},
  volume = {19},
  number = {2},
  pages = {633--667},
  publisher = {Institute of Mathematical Statistics},
  issn = {0090-5364, 2168-8966},
  doi = {10.1214/aos/1176348114},
  url = {https://projecteuclid.org/journals/annals-of-statistics/volume-19/issue-2/Geometrizing-Rates-of-Convergence-II/10.1214/aos/1176348114.full},
  urldate = {2026-07-27},
  langid = {english}
}

@article{donoho_one-sided_1988,
  title = {One-{{Sided Inference}} about {{Functionals}} of a {{Density}}},
  author = {Donoho, David L.},
  date = {1988},
  journaltitle = {The Annals of Statistics},
  volume = {16},
  number = {4},
  eprint = {2241771},
  eprinttype = {jstor},
  pages = {1390--1420},
  publisher = {Institute of Mathematical Statistics},
  issn = {0090-5364},
  url = {https://www.jstor.org/stable/2241771},
  urldate = {2026-07-27}
}

@article{donoho_statistical_1994,
  title = {Statistical {{Estimation}} and {{Optimal Recovery}}},
  author = {Donoho, David L.},
  date = {1994-03},
  journaltitle = {The Annals of Statistics},
  volume = {22},
  number = {1},
  pages = {238--270},
  publisher = {Institute of Mathematical Statistics},
  issn = {0090-5364, 2168-8966},
  doi = {10.1214/aos/1176325367},
  url = {https://projecteuclid.org/journals/annals-of-statistics/volume-22/issue-1/Statistical-Estimation-and-Optimal-Recovery/10.1214/aos/1176325367.full},
  urldate = {2026-07-27},
  langid = {english}
}

@article{dufour_impossibility_1997,
  title = {Some {{Impossibility Theorems}} in {{Econometrics With Applications}} to {{Structural}} and {{Dynamic Models}}},
  author = {Dufour, Jean-Marie},
  date = {1997},
  journaltitle = {Econometrica},
  volume = {65},
  number = {6},
  eprint = {2171740},
  eprinttype = {jstor},
  pages = {1365--1387},
  publisher = {[Wiley, Econometric Society]},
  issn = {0012-9682},
  doi = {10.2307/2171740},
  url = {https://www.jstor.org/stable/2171740},
  urldate = {2026-07-27}
}

@article{fiellerProblemsIntervalEstimation1954,
  title = {Some {{Problems}} in {{Interval Estimation}}},
  author = {Fieller, E. C.},
  date = {1954},
  journaltitle = {Journal of the Royal Statistical Society. Series B (Methodological)},
  volume = {16},
  number = {2},
  eprint = {2984043},
  eprinttype = {jstor},
  pages = {175--185},
  publisher = {[Royal Statistical Society, Oxford University Press]},
  issn = {0035-9246},
  url = {https://www.jstor.org/stable/2984043},
  urldate = {2026-10-02}
}

@book{folland_real_2013,
  title = {Real {{Analysis}}: {{Modern Techniques}} and {{Their Applications}}},
  shorttitle = {Real {{Analysis}}},
  author = {Folland, Gerald B.},
  date = {2013-06-11},
  eprint = {wI4fAwAAQBAJ},
  eprinttype = {googlebooks},
  publisher = {John Wiley \& Sons},
  isbn = {978-1-118-62639-9},
  langid = {english},
  pagetotal = {368}
}

@article{gleser_nonexistence_1987,
  title = {The {{Nonexistence}} of 100\$(1 - \textbackslash alpha)\$\% {{Confidence Sets}} of {{Finite Expected Diameter}} in {{Errors}}-in-{{Variables}} and {{Related Models}}},
  author = {Gleser, Leon Jay and Hwang, Jiunn T.},
  date = {1987-12},
  journaltitle = {The Annals of Statistics},
  volume = {15},
  number = {4},
  pages = {1351--1362},
  publisher = {Institute of Mathematical Statistics},
  issn = {0090-5364, 2168-8966},
  doi = {10.1214/aos/1176350597},
  url = {https://projecteuclid.org/journals/annals-of-statistics/volume-15/issue-4/The-Nonexistence-of-1001---alpha-Confidence-Sets-of/10.1214/aos/1176350597.full},
  urldate = {2026-07-27},
  langid = {english}
}

@article{griebel_decay_2018,
  title = {On the {{Decay Rate}} of the {{Singular Values}} of {{Bivariate Functions}}},
  author = {Griebel, Michael and Li, Guanglian},
  date = {2018-01},
  journaltitle = {SIAM Journal on Numerical Analysis},
  shortjournal = {SIAM J. Numer. Anal.},
  volume = {56},
  number = {2},
  pages = {974--993},
  publisher = {{Society for Industrial and Applied Mathematics}},
  issn = {0036-1429},
  doi = {10.1137/17M1117550},
  url = {https://epubs.siam.org/doi/10.1137/17M1117550},
  urldate = {2026-06-30}
}

@book{hadamard_lectures_1952,
  title = {Lectures on {{Cauchy}}'s {{Problem}} in {{Linear Partial Differential Equations}}},
  author = {Hadamard, Jacques},
  date = {1952},
  eprint = {wvRQAAAAMAAJ},
  eprinttype = {googlebooks},
  publisher = {Dover},
  langid = {english},
  pagetotal = {332}
}

@article{hirano_impossibility_2012,
  title = {Impossibility {{Results}} for {{Nondifferentiable Functionals}}},
  author = {Hirano, Keisuke and Porter, Jack R.},
  date = {2012},
  journaltitle = {Econometrica},
  volume = {80},
  number = {4},
  eprint = {23271416},
  eprinttype = {jstor},
  pages = {1769--1790},
  publisher = {[Wiley, Econometric Society]},
  issn = {0012-9682},
  url = {https://www.jstor.org/stable/23271416},
  urldate = {2026-07-24}
}

@article{horowitz_applied_2011,
  title = {Applied {{Nonparametric Instrumental Variables Estimation}}},
  author = {Horowitz, Joel L.},
  date = {2011},
  journaltitle = {Econometrica},
  volume = {79},
  number = {2},
  eprint = {41057460},
  eprinttype = {jstor},
  pages = {347--394},
  publisher = {[Wiley, Econometric Society]},
  issn = {0012-9682},
  url = {https://www.jstor.org/stable/41057460},
  urldate = {2026-06-19}
}

@article{huNONPARAMETRICIDENTIFICATIONUSING2018,
  title = {{{NONPARAMETRIC IDENTIFICATION USING INSTRUMENTAL VARIABLES}}: {{SUFFICIENT CONDITIONS FOR COMPLETENESS}}},
  shorttitle = {{{NONPARAMETRIC IDENTIFICATION USING INSTRUMENTAL VARIABLES}}},
  author = {Hu, Yingyao and Shiu, Ji-Liang},
  date = {2018-06},
  journaltitle = {Econometric Theory},
  volume = {34},
  number = {3},
  pages = {659--693},
  issn = {0266-4666, 1469-4360},
  doi = {10.1017/S0266466617000251},
  url = {https://www.cambridge.org/core/journals/econometric-theory/article/nonparametric-identification-using-instrumental-variables-sufficient-conditions-for-completeness/855B1E76EC5CB04C6BA90418CEA7AAED},
  urldate = {2026-09-22},
  langid = {english}
}

@article{kaji_theory_2021,
  title = {Theory of {{Weak Identification}} in {{Semiparametric Models}}},
  author = {Kaji, Tetsuya},
  date = {2021},
  journaltitle = {Econometrica},
  volume = {89},
  number = {2},
  pages = {733--763},
  issn = {1468-0262},
  doi = {10.3982/ECTA16413},
  url = {https://onlinelibrary.wiley.com/doi/abs/10.3982/ECTA16413},
  urldate = {2026-06-10},
  langid = {english}
}

@online{kallusCausalInferenceUnmeasured2022,
  title = {Causal {{Inference Under Unmeasured Confounding With Negative Controls}}: {{A Minimax Learning Approach}}},
  shorttitle = {Causal {{Inference Under Unmeasured Confounding With Negative Controls}}},
  author = {Kallus, Nathan and Mao, Xiaojie and Uehara, Masatoshi},
  date = {2022-10-09},
  eprint = {2103.14029},
  eprinttype = {arXiv},
  eprintclass = {stat},
  doi = {10.48550/arXiv.2103.14029},
  url = {http://arxiv.org/abs/2103.14029},
  urldate = {2025-10-09},
  pubstate = {prepublished}
}

@book{kato_perturbation_1995,
  title = {Perturbation {{Theory}} for {{Linear Operators}}},
  author = {Kato, Tosio},
  date = {1995},
  series = {Classics in {{Mathematics}}},
  volume = {132},
  publisher = {Springer},
  location = {Berlin, Heidelberg},
  doi = {10.1007/978-3-642-66282-9},
  url = {http://link.springer.com/10.1007/978-3-642-66282-9},
  urldate = {2026-07-27},
  isbn = {978-3-540-58661-6}
}

@software{kenny_simengine_2024,
  title = {{{SimEngine}}: {{A Modular Framework}} for {{Statistical Simulations}} in {{R}}},
  shorttitle = {{{SimEngine}}},
  author = {Kenny, Avi and Wolock, Charles},
  date = {2024-04-04},
  url = {https://cloud.r-project.org/web/packages/SimEngine/index.html},
  urldate = {2026-07-27},
  version = {1.4.0}
}

@article{khanIrregularIdentificationSupport2010,
  title = {Irregular {{Identification}}, {{Support Conditions}}, and {{Inverse Weight Estimation}}},
  author = {Khan, Shakeeb and Tamer, Elie},
  date = {2010},
  journaltitle = {Econometrica},
  volume = {78},
  number = {6},
  eprint = {40928469},
  eprinttype = {jstor},
  pages = {2021--2042},
  publisher = {[Wiley, Econometric Society]},
  issn = {0012-9682},
  url = {https://www.jstor.org/stable/40928469},
  urldate = {2026-10-02}
}

@article{klaassen_consistent_1987-2,
  title = {Consistent {{Estimation}} of the {{Influence Function}} of {{Locally Asymptotically Linear Estimators}}},
  author = {Klaassen, Chris A. J.},
  date = {1987-12},
  journaltitle = {The Annals of Statistics},
  volume = {15},
  number = {4},
  pages = {1548--1562},
  publisher = {Institute of Mathematical Statistics},
  issn = {0090-5364, 2168-8966},
  doi = {10.1214/aos/1176350609},
  url = {https://projecteuclid.org/journals/annals-of-statistics/volume-15/issue-4/Consistent-Estimation-of-the-Influence-Function-of-Locally-Asymptotically-Linear/10.1214/aos/1176350609.full},
  urldate = {2026-07-27},
  langid = {english}
}

@book{kress_linear_2014-2,
  title = {Linear Integral Equations},
  author = {Kress, Rainer},
  date = {2014},
  series = {Applied Mathematical Sciences ; Volume 82},
  edition = {Third edition.},
  publisher = {Springer},
  location = {New York},
  doi = {10.1007/978-1-4614-9593-2},
  isbn = {978-1-4614-9593-2}
}

@book{lecamAsymptoticsStatistics2000,
  title = {Asymptotics in {{Statistics}}},
  author = {Le Cam, Lucien and Lo Yang, Grace},
  date = {2000},
  series = {Springer {{Series}} in {{Statistics}}},
  publisher = {Springer},
  location = {New York, NY},
  doi = {10.1007/978-1-4612-1166-2},
  url = {http://link.springer.com/10.1007/978-1-4612-1166-2},
  urldate = {2026-10-02},
  isbn = {978-1-4612-7030-0}
}

@online{liNonparametricInferenceMean2022,
  title = {Nonparametric Inference about Mean Functionals of Nonignorable Nonresponse Data without Identifying the Joint Distribution},
  author = {Li, Wei and Miao, Wang and Tchetgen, Eric Tchetgen},
  date = {2022-04-06},
  eprint = {2110.05776},
  eprinttype = {arXiv},
  eprintclass = {math},
  doi = {10.48550/arXiv.2110.05776},
  url = {http://arxiv.org/abs/2110.05776},
  urldate = {2025-10-15},
  pubstate = {prepublished}
}

@article{low_nonparametric_1997,
  title = {On Nonparametric Confidence Intervals},
  author = {Low, Mark G.},
  date = {1997-12},
  journaltitle = {The Annals of Statistics},
  volume = {25},
  number = {6},
  pages = {2547--2554},
  publisher = {Institute of Mathematical Statistics},
  issn = {0090-5364, 2168-8966},
  doi = {10.1214/aos/1030741084},
  url = {https://projecteuclid.org/journals/annals-of-statistics/volume-25/issue-6/On-nonparametric-confidence-intervals/10.1214/aos/1030741084.full},
  urldate = {2026-07-27},
  langid = {english}
}

@article{luedtke_one-step_2024,
  title = {One-Step Estimation of Differentiable {{Hilbert}}-Valued Parameters},
  author = {Luedtke, Alex and Chung, Incheoul},
  date = {2024-08},
  journaltitle = {The Annals of Statistics},
  volume = {52},
  number = {4},
  pages = {1534--1563},
  publisher = {Institute of Mathematical Statistics},
  issn = {0090-5364, 2168-8966},
  doi = {10.1214/24-AOS2403},
  url = {https://projecteuclid.org/journals/annals-of-statistics/volume-52/issue-4/One-step-estimation-of-differentiable-Hilbert-valued-parameters/10.1214/24-AOS2403.full},
  urldate = {2026-07-27},
  langid = {english}
}

@article{miao_identification_2018,
  title = {Identification and Inference with Nonignorable Missing Covariate Data},
  author = {Miao, Wang and Tchetgen Tchetgen, Eric},
  date = {2018-10},
  journaltitle = {Statistica Sinica},
  shortjournal = {Stat Sin},
  volume = {28},
  number = {4},
  eprint = {33343174},
  eprinttype = {pubmed},
  pages = {2049--2067},
  issn = {1017-0405},
  doi = {10.5705/ss.202016.0322},
  url = {https://pmc.ncbi.nlm.nih.gov/articles/PMC7746016/},
  urldate = {2025-09-29},
  pmcid = {PMC7746016}
}

@article{newey_instrumental_2003-2,
  title = {Instrumental {{Variable Estimation}} of {{Nonparametric Models}}},
  author = {Newey, Whitney K. and Powell, James L.},
  date = {2003},
  journaltitle = {Econometrica},
  volume = {71},
  number = {5},
  eprint = {1555512},
  eprinttype = {jstor},
  pages = {1565--1578},
  publisher = {[Wiley, Econometric Society]},
  issn = {0012-9682},
  url = {https://www.jstor.org/stable/1555512},
  urldate = {2026-06-19}
}

@online{olivas-martinez_source-condition_2025,
  title = {Source-{{Condition Analysis}} of {{Kernel Adversarial Estimators}}},
  author = {Olivas-Martinez, Antonio and Rotnitzky, Andrea},
  date = {2025-08-24},
  eprint = {2508.17181},
  eprinttype = {arXiv},
  eprintclass = {math.ST},
  doi = {10.48550/arXiv.2508.17181},
  url = {http://arxiv.org/abs/2508.17181},
  urldate = {2026-05-20},
  pubstate = {prepublished}
}

@article{pfanzagl_nonexistence_1998,
  title = {The Nonexistence of Confidence Sets for Discontinuous Functionals},
  author = {Pfanzagl, J.},
  date = {1998-11-15},
  journaltitle = {Journal of Statistical Planning and Inference},
  shortjournal = {Journal of Statistical Planning and Inference},
  volume = {75},
  number = {1},
  pages = {9--20},
  issn = {0378-3758},
  doi = {10.1016/S0378-3758(98)00121-9},
  url = {https://www.sciencedirect.com/science/article/pii/S0378375898001219},
  urldate = {2026-07-25}
}

@article{potscher_lower_2002,
  title = {Lower {{Risk Bounds}} and {{Properties}} of {{Confidence Sets}} for {{Ill-Posed Estimation Problems}} with {{Applications}} to {{Spectral Density}} and {{Persistence Estimation}}, {{Unit Roots}}, and {{Estimation}} of {{Long Memory Parameters}}},
  author = {Pötscher, Benedikt M.},
  date = {2002},
  journaltitle = {Econometrica},
  volume = {70},
  number = {3},
  pages = {1035--1065},
  issn = {1468-0262},
  doi = {10.1111/1468-0262.00318},
  url = {https://onlinelibrary.wiley.com/doi/abs/10.1111/1468-0262.00318},
  urldate = {2026-07-10},
  langid = {english}
}

@manual{rcoreteamLanguageEnvironmentStatistical2026,
  type = {manual},
  title = {R: A Language and Environment for Statistical Computing},
  author = {{R Core Team}},
  date = {2026},
  publisher = {R Foundation for Statistical Computing},
  location = {Vienna, Austria},
  doi = {10.32614/R.manuals},
  url = {https://www.R-project.org/}
}

@article{santos_inference_2012,
  title = {Inference in {{Nonparametric Instrumental Variables}} with {{Partial Identification}}},
  author = {Santos, Andres},
  date = {2012},
  journaltitle = {Econometrica},
  volume = {80},
  number = {1},
  eprint = {41336585},
  eprinttype = {jstor},
  pages = {213--275},
  publisher = {[Wiley, Econometric Society]},
  issn = {0012-9682},
  url = {https://www.jstor.org/stable/41336585},
  urldate = {2026-06-19}
}

@article{severiniEfficiencyBoundsEstimating2012,
  title = {Efficiency Bounds for Estimating Linear Functionals of Nonparametric Regression Models with Endogenous Regressors},
  author = {Severini, Thomas A. and Tripathi, Gautam},
  date = {2012-10-01},
  journaltitle = {Journal of Econometrics},
  shortjournal = {Journal of Econometrics},
  series = {Thirtieth {{Anniversary}} of {{Generalized Method}} of {{Moments}}},
  volume = {170},
  number = {2},
  pages = {491--498},
  issn = {0304-4076},
  doi = {10.1016/j.jeconom.2012.05.018},
  url = {https://www.sciencedirect.com/science/article/pii/S0304407612001303},
  urldate = {2025-10-15}
}

@article{smuclerAsymptoticValidityConfidence2025,
  title = {On the Asymptotic Validity of Confidence Sets for Linear Functionals of Solutions to Integral Equations},
  author = {Smucler, E and Robins, J M and Rotnitzky, A},
  date = {2025-10-01},
  journaltitle = {Biometrika},
  shortjournal = {Biometrika},
  volume = {112},
  number = {4},
  pages = {asaf067},
  issn = {1464-3510},
  doi = {10.1093/biomet/asaf067},
  url = {https://doi.org/10.1093/biomet/asaf067},
  urldate = {2026-10-01}
}

@article{smuclerNotePropertiesConfidence2026,
  title = {A Note on the Properties of the Confidence Set for the Local Average Treatment Effect Obtained by Inverting the Score Test},
  author = {Smucler, Ezequiel and Lanni, Ludovico and Bonet, David Masip},
  date = {2026-01-01},
  journaltitle = {Journal of Causal Inference},
  volume = {14},
  number = {1},
  publisher = {De Gruyter},
  issn = {2193-3685},
  doi = {10.1515/jci-2025-0053},
  url = {https://www.degruyterbrill.com/document/doi/10.1515/jci-2025-0053/html?recommended=sidebar&lang=en},
  urldate = {2026-09-21},
  langid = {english}
}

@article{staiger_instrumental_1997,
  title = {Instrumental {{Variables Regression}} with {{Weak Instruments}}},
  author = {Staiger, Douglas and Stock, James H.},
  date = {1997},
  journaltitle = {Econometrica},
  volume = {65},
  number = {3},
  eprint = {2171753},
  eprinttype = {jstor},
  pages = {557--586},
  publisher = {[Wiley, Econometric Society]},
  issn = {0012-9682},
  doi = {10.2307/2171753},
  url = {https://www.jstor.org/stable/2171753},
  urldate = {2026-06-10}
}

@book{stewart_matrix_1990,
  title = {Matrix {{Perturbation Theory}}},
  author = {Stewart, G. W. and Sun, Ji-guang},
  date = {1990-06-28},
  eprint = {l78PAQAAMAAJ},
  eprinttype = {googlebooks},
  publisher = {Elsevier Science},
  isbn = {978-0-12-670230-9},
  langid = {english},
  pagetotal = {392}
}

@book{vaart_asymptotic_1998-2,
  title = {Asymptotic {{Statistics}}},
  author = {family=Vaart, given=A. W., prefix=van der, useprefix=false},
  date = {1998},
  series = {Cambridge {{Series}} in {{Statistical}} and {{Probabilistic Mathematics}}},
  publisher = {Cambridge University Press},
  location = {Cambridge},
  doi = {10.1017/CBO9780511802256},
  url = {https://www.cambridge.org/core/books/asymptotic-statistics/A3C7DAD3F7E66A1FA60E9C8FE132EE1D},
  urldate = {2025-04-15},
  isbn = {978-0-521-78450-4}
}

@article{vaart_differentiable_1991,
  title = {On {{Differentiable Functionals}}},
  author = {Vaart, Aad Van Der},
  date = {1991-03},
  journaltitle = {The Annals of Statistics},
  volume = {19},
  number = {1},
  pages = {178--204},
  publisher = {Institute of Mathematical Statistics},
  issn = {0090-5364, 2168-8966},
  doi = {10.1214/aos/1176347976},
  url = {https://projecteuclid.org/journals/annals-of-statistics/volume-19/issue-1/On-Differentiable-Functionals/10.1214/aos/1176347976.full},
  urldate = {2024-10-14}
}

\end{refsection}

\end{document}